\documentclass[a4paper,10pt, openany]{book}

\usepackage[french]{babel}
\usepackage[cp1252]{inputenc}
\usepackage{amsthm}
\usepackage{amsmath}
\usepackage{amsfonts}
\usepackage{amssymb}
\newcommand{Ž}{\'e}
\newcommand{}{\`{e}}
\newcommand{ˆ}{\`{a}}
\newcommand{}{\`{u}}
\newcommand{}{\c{c}}
\newcommand{}{\^e}
\newcommand{"}{\^{\i}}
\newcommand{•}{\"{\i}}
\newcommand{ž}{\^u}
\newcommand{‰}{\^a}
\newcommand{™}{\^{o}}

\newcommand{\stackunder}[2]{\underset{#1}{#2}}

\newcommand{\R}{\mathbb{R}}

\newcommand{\ve}{\varepsilon}

\newcommand{\g}{\mathbf{g}}
\newcommand{\E}{(E_{h,f,\mathbf{g}})}

\newtheorem{definition}{D\'efinition}
\newtheorem{theo}{Theorem}
\newtheorem*{theor}{Theorem}
\newtheorem{defi}{Definition}

\newtheorem{theorem}{Th\'eor\`eme}

\newtheorem{proposition}{Proposition}

\begin{document}

\title{Fonctions critiques et Žquations aux d\'{e}riv\'{e}es partielles elliptiques
sur les vari\'{e}t\'{e}s Riemanniennes compactes}
\author{\textbf{Stephane Collion}}

\date{Thse de Doctorat de l'UniversitŽ Pierre et Marie Curie, Paris VI \\ SpŽcialitŽ \textbf{MathŽmatiques}\\Thse soutenue le 4 DŽcembre 2004}
\maketitle

\chapter*{ }

\begin{flushright}
ˆ Alice,\\
pour chaque minute que je ne lui ai pas consacrŽe \\
ˆ cause de cette thse;\\
\end{flushright}

\begin{flushright}
ˆ Michel,\\
pour chaque minute qu'il m'a consacrŽe \\
ˆ cause de cette thse.

\end{flushright}

\tableofcontents

\chapter*{Merci...}
	Au-delˆ du respect d'une tradition, certes trs agrŽable, j'aimerais que ces quelques lignes soient vues comme l'expression sincre du respect, de l'affection et de la reconnaissance que j'ai pour tous ceux qui m'entourent et me soutiennent dans mes choix et mes passions.

	Je tiens tout d'abord ˆ remercier Messieurs Xavier CabrŽ et Franck Pacard d'avoir acceptŽ d'tre les rapporteurs de mon travail de thse. Je remercie Žgalement Messieurs Emmanuel Hebey, FrŽdŽric HŽlein et Henri Skoda de faire partie, aujourd'hui de mon jury. 
	
	Je remercie tout particulirement Henri Skoda et Guennadi Henkin qui ont dirigŽ mes premiers travaux de recherche. Ils ont tous les deux supportŽs avec beaucoup de gentillesse et de tolŽrance mes allers-retours entre mes deux passions, les mathŽmatiques et l'aviation, et donc mon inconstance. Je remercie d'ailleurs toute l'Žquipe d'analyse complexe de Paris 6, Pierre Mazet, Pierre Dolbeault, Andre• Iordan, Pascal Dingoyan et Vincent Michel  pour leur accueil toujours bienveillant ˆ mon Žgard.

	Un immense merci ˆ toute ma famille pour la grande ouverture d'esprit qu'ils ont toujours su entretenir chez nous, et pour le soutien constant apportŽ ˆ chacun d'entre nous dans ses choix et ses passions. Merci ˆ Elodie, Evelyne, David et Jean-Alain pour leur patience lorsque je travaillais ˆ Moislains. J'ai une pensŽe toute particulire pour ma grand-mre Gaby qui a eu la trs mauvaise idŽe de ne pas pouvoir tre prŽsente aujourd'hui.

	Un merci particulier ˆ ma mre qui a eu le courage de lire mon manuscrit pour essayer d'y repŽrer, au milieu de formules sans doute bien obscures pour elle, les fautes d'orthographe. Merci ˆ Pierre Leycuras d'en avoir fait aussi une relecture.

	A mes vieux amis Arnaud, CŽsar, Eric, Jean-David, Lionel, Romain 1, Romain 2; vous avez souvent cru que je n'y arriverais jamais, mais aujourd'hui c'est fait... Merci ˆ vous et ˆ Delphine, Marina, Anna, AurŽlie, Loane, ainsi qu'aux amis plus "rŽcents" mais tout aussi prŽcieux, Anne, Camille, Charles, Dominique, Bertrand, Michel, Pascal, Patrice, Pierre, Thierry, Thierry, merci ˆ tous pour votre amitiŽ.

	Merci ˆ Mermoz qui est celui qui a passŽ le plus de temps (couchŽ) sur ma thse.

	Merci ˆ Alice ... pour tout.

	Quant ˆ Michel... Je crois lui avoir beaucoup dit ce que je lui devais, sinon je le lui  rŽpterais. Pour qu'il y en ait une petite trace Žcrite, disons qu'en plus d'enseigner des mathŽmatiques, il apprend ˆ en faire et ˆ y prendre toujours plus de plaisir. Mais, au-delˆ de ses talents de mathŽmaticien et de professeur,  ce pour quoi je lui dois le plus, c'est son amitiŽ. 

\chapter*{ }

\chapter*{Introduction}
\pagestyle{myheadings}\markboth{\textbf{Introduction.}}{\textbf{Introduction}}
\textit{Pr\'{e}ambule: }

Le propos de cette introduction est de d\'{e}finir l'esprit dans lequel nous
avons r\'{e}dig\'{e} cette th\`{e}se. Nous avons voulu profiter de l'absence
de contraintes de longueur g\'{e}n\'{e}ralement impos\'{e}es aux articles
publi\'{e}s par les diff\'{e}rentes revues pour d\'{e}tailler autant que
possible nos d\'{e}monstrations. Ainsi, avant de nous lancer dans les
d\'{e}tails souvent tr\`{e}s ``techniques'' de ces d\'{e}monstrations, nous
avons tent\'{e} d'en expliciter les id\'{e}es, les principes et les
difficult\'{e}s essentielles; ce que, la plupart du temps, on n'a
malheureusement pas la place de faire dans les articles soumis aux revues.
Les id\'{e}es, en particulier \`{a} la fin du chapitre 2, sont parfois
pr\'{e}sent\'{e}es de mani\`{e}re heuristique. Le lecteur est \'{e}galement
invit\'{e} \`{a} se reporter fr\'{e}quemment au dernier appendice qui
reprend nos principales notations et conventions. Nous prions le lecteur de
nous excuser de rallonger ainsi le texte, mais notre but est de rendre notre
travail aussi clair et lisible que possible.

Nous reprenons \'{e}galement quelques d\'{e}monstrations qui ne nous
appartiennent pas, bien que l\`{a} encore, notre th\`{e}se s'en trouve
rallong\'{e}e. Notre but est triple. Tout d'abord, nous avons voulu
\'{e}viter au lecteur d'avoir \`{a} se reporter trop fr\'{e}quemment \`{a}
diff\'{e}rents articles pour suivre nos d\'{e}monstrations et faire en sorte
que notre expos\'{e} soit aussi ``complet'' que possible. Ensuite, comme
nous venons de le dire, soumises aux contraintes impos\'{e}es par les
revues, ces d\'{e}monstrations sont souvent pr\'{e}sent\'{e}es avec tr\`{e}s
peu de d\'{e}tails. Nous avons cherch\'{e} ici \`{a} les pr\'{e}senter de
notre point de vue et avec plus de d\'{e}tails, \`{a} en expliquer les
principes, dans le but de les mettre en valeur. Enfin, bien que peu de
modifications soient n\'{e}cessaires, les r\'{e}sultats que nous utilisons
n'apparaissent pas, dans ces articles, dans le m\^{e}me cadre et avec
exactement nos hypoth\`{e}ses. Nous n'avons pas voulu utiliser abruptement
la formule ``\c{c}a marche pareil'' et risquer ainsi de donner l'impression
d'occulter quelques difficult\'{e}s. Nous nous sommes efforc\'{e}s
d'indiquer clairement quelles \'{e}taient ces d\'{e}monstrations, qui
\'{e}taient leurs auteurs et quelles modifications avaient \'{e}t\'{e}
n\'{e}cessaires \`{a} notre travail.

Nous esp\'{e}rons que, r\'{e}dig\'{e}e dans cet esprit, notre th\`{e}se sera
agr\'{e}able \`{a} lire.

Voici, dans le d\'{e}tail, les diff\'{e}rents chapitres:

C\textit{hapitre 1}: \textbf{D\'{e}finitions et \'{e}nonc\'{e} des
r\'{e}sultats: }Nous donnons les d\'{e}finitions et les th\'{e}or\`{e}mes
d\'{e}montr\'{e}s dans cette th\`{e}se. Ici encore nous avons choisi un mode
de pr\'{e}sentation un peu particulier. Plut\^{o}t que de donner directement
la d\'{e}finition la plus g\'{e}n\'{e}rale possible des fonctions critiques,
nous montrons comment, au fur et \`{a} mesure de l'obtention de nos
r\'{e}sultats, nous avons fait \'{e}voluer la d\'{e}finition initialement
donn\'{e}e par E. Hebey et M. Vaugon pour l'\'{e}tude des meilleures
constantes dans les in\'{e}galit\'{e}s de Sobolev, vers celle de ''triplet
critique'' pour mettre en valeur la puissance de ce concept dans l'\'{e}tude
des EDP qui nous int\'{e}ressent. La pr\'{e}sentation de ce chapitre est
donc plus chronologique que synth\'{e}tique.

\textit{Chapitre 2}: \textbf{Trois outils fondamentaux: le point de
concentration, le changement d'echelle, le processus d'it\'{e}ration.
Principe des d\'{e}monstrations.} Ces trois outils sont \`{a} la base des
m\'{e}thodes utilis\'{e}es pour l'\'{e}tude des ph\'{e}nom\`{e}nes de
concentration (comme le nom du premier l'indique). Leur pr\'{e}sentation,
dans les articles sur le sujet, est malheureusement souvent noy\'{e}e au
milieu des d\'{e}monstrations. Nous avons voulu isoler la d\'{e}finition de
ces outils maintenant classiques pour les mettre en valeur et pour rendre
les d\'{e}monstrations plus lisibles. Nous pr\'{e}sentons \'{e}galement
\`{a} la fin de ce chapitre le principe des principales démonstrations.

\textit{Chapitre 3: }\textbf{Existence de Fonctions Extr\'{e}males, Seconde
In\'{e}galit\'{e} Fondamentale.} C'est le coeur de notre travail. Dans une
premi\`{e}re partie, nous exposons les r\'{e}sultats sur les
ph\'{e}nom\`{e}nes de concentration valables pour une famille
g\'{e}nérale d'EDP. Ils sont issus du travail de plusieurs auteurs, et
pr\'{e}sent\'{e}s selon le principe que nous avons expos\'{e} en
pr\'{e}ambule. Vient ensuite une m\'{e}thode que nous d\'{e}veloppons pour
obtenir la ``seconde in\'{e}galit\'{e} fondamentale'' et qui permet
d'aboutir au principal th\'{e}or\`{e}me.

\textit{Chapitre 4}: \textbf{Triplet Critique 1: }\textit{Existence de
Fonctions Critiques.} Notre travail portant sur les fonctions critiques, il
est utile de montrer qu'il en existe !

\textit{Chapitre 5:} \textbf{Triplet Critique 2: }Suivant la d\'{e}finition
donn\'{e}e d'un triplet critique, nous regardons ce que l'on peut dire
lorsque la m\'{e}trique varie dans une classe conforme.

\textit{Chapitre 6}: \textbf{Triplet Critique 3: }Le dernier probl\`{e}me
li\'{e} \`{a} l'existence des triplets critiques est abord\'{e}. C'est une
partie tr\'{e}s importante de notre travail, et bien que plus courte, les
r\'{e}sultats n'y sont obtenus que gr\^{a}ce \`{a} la m\'{e}thode
d\'{e}velopp\'{e}e dans le chapitre 3.

\textit{Chapitre 7:} \textbf{La dimension 3.} Les pr\'{e}c\'{e}dents
chapitres portent sur l'\'{e}tude d'EDP sur des vari\'{e}t\'{e}s de
dimension au moins 4. La dimension 3 est \`{a} part, elle est
\'{e}tudi\'{e}e dans ce chapitre.

\textit{Chapitre 8: }\textbf{Remarques sur le cas limite et le cas
d\'{e}g\'{e}n\'{e}r\'{e}. }Quelques remarques (et regrets ?) sont
propos\'{e}es sur le possible affaiblissement des hypoth\`{e}ses faites dans
les chapitres pr\'{e}c\'{e}dents. Quelques questions li\'{e}es \`{a} ce
travail auxquelles nous aimerions r\'{e}pondre dans l'avenir sont
expos\'{e}es.

\textit{Chapitre 9:} \textbf{Version Anglaise abrŽgŽe. }

\textit{Appendice A: }\textbf{D\'{e}monstration des propri\'{e}t\'{e}s
\'{e}l\'{e}mentaires des fonctions critiques cit\'{e}es au chapitre 1.}
Report\'{e}es ici pour une plus grande lisibilit\'{e} du premier chapitre.

\textit{Appendice B:} \textbf{Construction d'une fonction de Green.} Nous
nous servons de cette notion, mais bien que souvent cit\'{e}e, n'ayant pas
pu trouver de r\'{e}f\'{e}rences pr\'{e}cises, nous en donnons une
construction rapide.

\textit{Appendice C}: \textbf{Limite de }$C_{t}/\int u_{t}^{2}$. Il s'agit
d'une limite apparaissant au chapitre 3. Son calcul est issu d'un article de
Z. Djadli et O.\ Druet. Conform\'{e}ment \`{a} notre principe nous avons
voulu l'inclure car, dans l'article de Djadli et Druet, elle n'est pas
calcul\'{e}e avec exactement nos hypoth\`{e}ses. N\'{e}anmoins le calcul
\'{e}tant assez long, nous l'avons report\'{e} en appendice pour ne pas
interrompre le cours de notre d\'{e}monstration.

\textit{Appendice D:} \textbf{Notations et Conventions}. En esp\'{e}rant que
cela aide le lecteur \`{a} s'y retrouver.

\chapter{D\'{e}finitions et \'{e}nonc\'{e}s des r\'{e}sultats.}
\pagestyle{myheadings}\markboth{\textbf{Définitions et résultats.}}
{\textbf{Définitions et résultats.}}
Au commencement \'{e}tait le probl\`{e}me de Yamabe:

\textbf{Probl\`{e}me de Yamabe}\textit{: Etant donn\'{e}e une
vari\'{e}t\'{e} Riemannienne compacte }$(M,\mathbf{g})$\textit{\ de
dimension }$n\geq 3$,\textit{\ existe-t-il une m\'{e}trique conforme
\`{a} }$\mathbf{g}$\textit{\ dont la courbure scalaire est constante.}

R\'{e}soudre ce probl\`{e}me revient \`{a} prouver l'existence d'une
solution $u>0$ de l'\'{e}quation aux d\'{e}riv\'{e}es partielles: 
\[
\triangle _{\mathbf{g}}u+\frac{n-2}{4(n-1)}S_{\mathbf{g}}.u=\lambda .u^{%
\frac{n+2}{n-2}}
\]
o\`{u} $S_{\mathbf{g}}$ est la courbure scalaire de $\mathbf{g}$. Nous
proposons dans cette introduction un aper\c{c}u (pas forc\'{e}ment
chronologique et subjectivement orient\'{e} vers nos probl\`{e}mes) des d%
\'{e}veloppements qui suivirent l'\'{e}nonc\'{e} de ce probl\`{e}me, dans le
but de d\'{e}finir la notion de fonction critique introduite par E. Hebey et
M. Vaugon [20].

Le probl\`{e}me de Yamabe lan\c{c}a l'\'{e}tude d'EDP non lin\'{e}aires sur
des vari\'{e}t\'{e}s Riemanniennes compactes de la forme: 
$$(E_{h,f,\mathbf{g}}):\triangle _{\mathbf{g}}u+hu=fu^{\frac{n+2}{n-2}}$$
o\`{u} $\Delta _{\mathbf{g}}u=-\nabla ^{i}\nabla _{i}u$ est le laplacien
riemannien de la vari\'{e}t\'{e} riemannienne compacte $(M,\mathbf{g})$%
\textit{\ }de dimension \textit{\ }$n\geq 3$ et de m\'{e}trique $%
\mathbf{g}$ , et o\`{u} $h,f\in C^{\infty }(M)$ sont des fonctions
donn\'{e}es, l'inconnue \'{e}tant $u$, que l'on cherche en g\'{e}n\'{e}ral
strictement positive. Par exemple, le cas $h=\frac{n-2}{4(n-1)}S_{\mathbf{g}%
} $ correspond aux probl\`{e}mes de courbure scalaire prescrite: en effet,
si $\mathbf{g}^{\prime }=u^{\frac{4}{n-2}}.\mathbf{g}$ est une m\'{e}trique
conforme \`{a} $\mathbf{g}$, les coubures scalaires sont reli\'{e}es par
l'\'{e}quation: 
\[
\triangle _{\mathbf{g}}u+\frac{n-2}{4(n-1)}S_{\mathbf{g}}.u=\frac{n-2}{4(n-1)%
}S_{\mathbf{g}^{\prime }}.u^{\frac{n+2}{n-2}} 
\]

Chercher une m\'{e}trique conforme \`{a} $\mathbf{g}$ dont $f$ soit la
courbure scalaire revient donc \`{a} chercher une solution $u>0$ de ($E_{h,f,%
\mathbf{g}}$) avec $h=\frac{n-2}{4(n-1)}S_{\mathbf{g}}$. Le probl\`{e}me de
Yamabe est un cas particulier o\`{u} $f$ est une constante.

Dans l'\'{e}tude de ces \'{e}quations, on utilise bien s\^{u}r les espaces
de Sobolev, et T.Aubin [2] mit en \'{e}vidence un lien fondamental entre l'%
\'{e}quation ($E_{h,f,\mathbf{g}}$) et la notion de meilleure constante dans
les inclusions de Sobolev. Notons $H_{1}^{2}(M)$ l'espace de Sobolev des
fonctions $L^{2}(M)$ dont le gradient est \'{e}galement dans $L^{2}(M)$.
L'inclusion continue de $H_{1}^{2}(M)$ dans $L^{2^{*}}(M)$, o\`{u} $2^{*}=%
\frac{2n}{n-2}$ est l'exposant critique pour les inclusions de $H_{1}^{2}(M)$
dans $L^{p}(M)$ (compacte pour $p<2^{*}$ et seulement continue pour $p=2^{*}$%
), se traduit par l'existence de deux constantes A et B telles que pour
toute $u\in H_{1}^{2}(M):$%
\begin{equation}
\left( \int_{M}\left| u\right| ^{\frac{2n}{n-2}}dv_{\mathbf{g}}\right) ^{%
\frac{n-2}{n}}\leq A\int_{M}\left| \nabla u\right| _{\mathbf{g}}^{2}dv_{%
\mathbf{g}}+B\int_{M}u{{}^{2}}dv_{\mathbf{g}}
\end{equation}
On d\'{e}finit la meilleure premi\`{e}re constante comme l'inf. des A telle
qu'il existe $B>0$ avec (1.1) vraie. Sa valeur est connue: 
\[
A=K(n,2){{}^{2}}=\frac{4}{n(n-2)\omega _{n}^{\frac{2}{n}}}
\]
o\`{u} $\omega _{n}$ est le volume de la sph\`{e}re unit\'{e} de dimension
n, et cet inf. est atteint [19]. On prend alors $B_{0}(\mathbf{g})$ l'inf.
des B tel que (1.1) soit vraie avec cette valeur de A; on montre que $B_{0}(%
\mathbf{g})<+\infty $ $\left[ 17\right] .$ L'in\'{e}galit\'{e} : 
\begin{equation}
\left( \int_{M}\left| u\right| ^{\frac{2n}{n-2}}dv_{\mathbf{g}}\right) ^{%
\frac{n-2}{n}}\leq K(n,2){{}^{2}}\int_{M}\left| \nabla u\right| _{\mathbf{g}%
}^{2}dv_{g}+B_{0}(\mathbf{g})\int_{M}u{{}^{2}}dv_{\mathbf{g}} 
\end{equation}
est alors optimale en ce sens que les deux constantes ne peuvent plus \^{e}%
tre diminu\'{e}es. Si la meilleure premi\`{e}re constante A est connue et
est ind\'{e}pendante de la vari\'{e}t\'{e} $(M,\mathbf{g})$, en revanche $%
B_{0}(\mathbf{g})$, comme la notation l'indique, d\'{e}pend de la g\'{e}om%
\'{e}trie et sa recherche est difficile; cela fait l'objet de plusieurs
articles et c'est dans ce but qu'ont \'{e}t\'{e} introduites les fonctions
critiques par E.Hebey et M.Vaugon [20].(Quand il n'y aura pas d'ambigu\"{i}t%
\'{e}, nous noterons ces deux constantes $K$ et $B_{0}.$)

\textit{L'objectif de notre travail est l'\'{e}tude de l'existence de
solutions strictement positives aux \'{e}quations (}$E_{h,f,\mathbf{g}}$%
\textit{) }$\triangle _{\mathbf{g}}u+hu=fu^{\frac{n+2}{n-2}}$\textit{\ dans
les cas limites normalement non r\'{e}solus par les m\'{e}thodes
variationnelles. Les fonctions critiques (}$h$ et $f$) \textit{apparaitront
dans notre travail comme provenant de ces cas ``limites'', et c'est donc ces
fonctions que nous \'{e}tudierons.}

Nous allons rappeler les bases de ces m\'{e}thodes, ce qui mettra en \'{e}%
vidence le lien entre l'\'{e}quation ($E_{h,f,\mathbf{g}}$) et l'inclusion
de Sobolev (1.2), et nous am\`{e}nera naturellement \`{a}
la d\'{e}finition des fonctions critiques. Notons que gr\^{a}ce \`{a} la
compacit\'{e} de l'inclusion de $H_{1}^{2}(M)$ dans $L^{p}(M)$ pour $p<2^{*}$%
, les m\'{e}thodes variationnelles et les th\'{e}ories elliptiques donnent
rapidement l'existence de solutions $u>0$ \`{a} l'\'{e}quation $\triangle _{%
\mathbf{g}}u+hu=fu^{p-1}$; le cas $p=2^{*}$ est donc d\'{e}j\`{a} en lui-m\^{e}%
me un cas limite (tr\'{e}s peu de choses sont connues pour $p>2^{*}$ sans
hypoth\`{e}ses suppl\'{e}mentaires [20]).

L'\'{e}tude des \'{e}quations du type ($E_{h,f,\mathbf{g}}$) par les
m\'{e}thodes variationnelles am\`{e}ne \`{a} consid\'{e}rer la fonctionnelle
d\'{e}finie sur $H_{1}^{2}(M):$%
\[
I_{h,\mathbf{g}}(w)=\int_{M}\left| \nabla w\right| _{\mathbf{g}%
}^{2}dv_{\mathbf{g}}+\int_{M}h.w{{}^{2}}dv_{\mathbf{g}} 
\]
et le minimum de cette fonctionnelle:
\[
\lambda _{h,f,\mathbf{g}}=\stackunder{w\in \mathcal{H}_{f}}{\inf }I_{h,%
\mathbf{g}}(w) 
\]
sur l'ensemble 
\[
\mathcal{H}_{f}=\{w\in H_{1}^{2}(M)/\int_{M}f\left| w\right| ^{\frac{2n}{n-2}%
}dv_{\mathbf{g}}=1\}. 
\]
En effet, l'\'{e}quation d'Euler associ\'{e}e au probl\`{e}me de
minimisation de cette fonctionnelle \`{a} l'aide d'une fonction $u$ telle
que 
\[
I_{h,\mathbf{g}}(u)=\stackunder{w\in \mathcal{H}_{f}}{\inf }I_{h,\mathbf{g}%
}(w) 
\]
est exactement 
\[
(E_{h,f,\mathbf{g}}):\triangle _{\mathbf{g}}u+hu=\lambda _{h,f,\mathbf{g}%
}.f.u^{\frac{n+2}{n-2}} 
\]
o\`{u} $\lambda _{h,f,\mathbf{g}}\,$appara\^{i}t comme une constante de
normalisation li\'{e}e \`{a} la condition 
\[
\int_{M}f\left| u\right| ^{\frac{2n}{n-2}}dv_{\mathbf{g}}=1. 
\]
Il est parfois utile d'utiliser la fonctionnelle 
\[
J_{h,f,\mathbf{g}}(w)=\frac{\int_{M}\left| \nabla w\right| _{\mathbf{g}%
}^{2}dv_{\mathbf{g}}+\int_{M}h.w{{}^{2}}dv_{\mathbf{g}}}{\left(
\int_{M}f\left| w\right| ^{\frac{2n}{n-2}}dv_{\mathbf{g}}\right) ^{\frac{n-2%
}{n}}} 
\]
et la partie de $H_{1}^{2}(M)$ pour laquelle elle est d\'{e}finie 
\[
\mathcal{H}_{f}^{+}=\{w\in H_{1}^{2}(M)/\int_{M}f\left| w\right| ^{\frac{2n}{%
n-2}}dv_{\mathbf{g}}>0\}. 
\]
On consid\`{e}re alors le probl\`{e}me de minimisation par une fonction $u$
telle que 
\[
J_{h,f,\mathbf{g}}(u)=\stackunder{w\in \mathcal{H}_{f}^{+}}{\inf }J_{h,f,%
\mathbf{g}}(w), 
\]
l'\'{e}quation d'Euler associ\'{e}e \'{e}tant identique, mais sans constante
de normalisation. La fonctionnelle $J$ pr\'{e}sente (parfois) l'avantage
d'\^{e}tre homog\`{e}ne au sens o\`{u} $J_{h,f,\mathbf{g}}(c.w)=J_{h,f,%
\mathbf{g}}(w)$ pour toute constante $c.$ On voit donc que 
\[
\stackunder{w\in \mathcal{H}_{f}}{\inf }I_{h,\mathbf{g}}(w)=\stackunder{w\in 
\mathcal{H}_{f}^{+}}{\inf }J_{h,f,\mathbf{g}}(w)=\lambda _{h,f,\mathbf{g}} 
\]
Cette fonctionnelle $J$ a la particularit\'{e}, lorsque $h=\frac{n-2}{4(n-1)}%
S_{\mathbf{g}},$ d'\^{e}tre invariante par changement de m\'{e}trique
conforme; elle est donc particuli\`{e}rement utile dans l'\'{e}tude des
probl\`{e}mes de courbure scalaire prescrite. Nous utiliserons la plupart du
temps $I_{h,\mathbf{g}}$ et $\mathcal{H}_{f}$, mais pour certains
probl\`{e}mes, $J_{h,f,\mathbf{g}}$ s'av\`{e}rera commode lorsqu'on voudra
s'affranchir de la contrainte $\int_{M}f\left| u\right| ^{\frac{2n}{n-2}}dv_{%
\mathbf{g}}=1.$

Les meilleures constantes des inclusions de Sobolev apparurent lorsque
T.Aubin $\left[ 2\right] $ montra dans l'\'{e}tude du probl\`{e}me de
Yamabe, o\`{u} $h=\frac{n-2}{4(n-1)}S_{\mathbf{g}}$ et $f$ est une
constante, que l'\'{e}quation admettait une solution $u>0,$ ce qui r\'{e}%
solvait le probl\`{e}me, si on avait: 
\[
\lambda _{\frac{n-2}{4(n-1)}S_{\mathbf{g}},1,\mathbf{g}}<K(n,2)^{-2}
\]
Plus pr\'{e}cis\'{e}ment, cette condition permet d'\'{e}viter d'aboutir par
les m\'{e}thodes variationnelles \`{a} une solution identiquement nulle. T.
Aubin montra ensuite plus g\'{e}n\'{e}ralement que pour toutes fonctions $%
h,f\in C^{\infty }(M)$ on a : 
\[
\lambda _{h,f,\mathbf{g}}\leq \frac{1}{K(n,2){{}^{2}}(\stackunder{M}{Sup}f)^{%
\frac{n-2}{n}}}
\]
et que, si l'in\'{e}galit\'{e} est stricte, alors il existe une solution $u>0
$ \`{a} ($E_{h,f,\mathbf{g}}$) qui de plus r\'{e}alise l'inf. de la
fonctionnelle $I_{h,f,\mathbf{g}}$ sur $\mathcal{H}_{f}$. Ceci montre d\'{e}%
ja l'importance de $K(n,2)$ dans l'\'{e}tude des \'{e}quations ($E_{h,f,%
\mathbf{g}}$).

\textbf{Ce r\'{e}sultat de T. Aubin est le point de d\'{e}part de tout ce
qui suit. Il permet \`{a} partir des m\'{e}thodes variationnelles de montrer
l'existence de solutions sous l'hypoth\`{e}se } 
\[
\lambda _{h,f,\mathbf{g}}<K(n,2){{}^{-2}}(\stackunder{M}{Sup}f)^{-\frac{n-2}{%
n}}\mathit{.\ } 
\]
\textbf{Notre travail porte essentiellement sur le probl\`{e}me de
l'existence de solutions dans le cas limite: } 
\[
\lambda _{h,f,\mathbf{g}}=K(n,2){{}^{-2}}(\stackunder{M}{Sup}f)^{-\frac{n-2}{%
n}}, 
\]
\textbf{\ probl\`{e}me normalement non r\'{e}solu par les m\'{e}thodes
variationnelles; c'est pour cette \'{e}tude que nous allons introduire les
fonctions critiques.}

Une question naturelle concernant (1.2) est l'existence de
fonctions extr\'{e}males $u$, c'est \`{a} dire r\'{e}alisant
l'\'{e}galit\'{e} dans (1.2). On montre alors par les
th\'{e}ories elliptiques standard que s'il en existe, $u$ est $C^{\infty }$
et que $u<0$ ou $u>0$; attention, ce cas ne correspond plus au
th\'{e}or\`{e}me de T. Aubin car ici $\lambda _{h,f,\mathbf{g}}=\frac{1}{%
K(n,2){{}^{2}}(\stackunder{M}{Sup}f)^{\frac{n-2}{n}}}$. Quitte \`{a}
remplacer $u$ par $-u$ et \`{a} une constante multiplicative pr\`{e}s on a
alors: 
\[
\triangle _{\mathbf{g}}u+\frac{B_{0}(\mathbf{g})}{K(n,2){{}^{2}}}u=\frac{1}{%
K(n,2){{}^{2}}}u^{\frac{n+2}{n-2}}\text{ sur }M\text{, et }\int_{M}u^{\frac{%
2n}{n-2}}dv_{\mathbf{g}}=1. 
\]
Ceci montre dans l'autre sens le lien entre (1.2) et ($%
E_{h,f,\mathbf{g}}$) . (Une autre fa\c{c}on d'apercevoir ce lien qui vaut
d'\^{e}tre mentionn\'{e} est de multiplier l'\'{e}quation ci-dessus par $u$
et d'int\'{e}grer sur $M$, on obtient alors l'\'{e}galit\'{e} dans (1.2) ). On notera souvent $2^{*}=\frac{2n}{n-2}$ , remarquons
alors que $2^{*}-1=\frac{n+2}{n-2}$. Cela sugg\`{e}re \'{e}galement que la
recherche d'informations sur $B_{0}(\mathbf{g})$ am\`{e}ne \`{a} \'{e}tudier
ces \'{e}quations.

Ces rappels ont pour but de justifier maintenant l'introduction de la notion
de fonction critique. Reprenons d'abord pr\'{e}cis\'{e}ment les donn\'{e}es
intervenant dans l'\'{e}tude des \'{e}quations ($E_{h,f,\mathbf{g}}$):

\textbf{Donn\'{e}es : }On consid\`{e}re une vari\'{e}t\'{e} riemannienne
compacte $(M,\mathbf{g})$ de dimension $n\geq 3.$ Soit $f:M\rightarrow \Bbb{R%
}$ une fonction $C^{\infty }$ \textit{fix\'{e}e} telle que $\stackunder{M}{%
Max}$ $f>0$. Soit aussi $h\in C{{}^{\infty }}(M)$ avec l'hypoth\'{e}se
suppl\'{e}mentaire que l'op\'{e}rateur $\bigtriangleup _{\mathbf{g}}+h$ est
coercif si $f$ change de signe sur $M$. Notons que nous ne cherchons pas ici
les hypoth\`{e}ses de r\'{e}gularit\'{e} minimales pour $h$ et $f$, la
continuit\'{e} \'{e}tant g\'{e}n\'{e}ralement suffisante dans la plupart des
r\'{e}sultats que nous pr\'{e}sentons; dans tout notre travail nous
supposerons que $h$ et $f$ sont $C^{\infty }$.

On consid\`{e}re l'\'{e}quation $(E_{h,f,\mathbf{g}}^{\prime }):\,\triangle
_{\mathbf{g}}u+h.u=f.u^{\frac{n+2}{n-2}}$ .

\textbf{Convention:} Nous garderons ces ``notations'': $\mathbf{g,g}^{\prime
},\mathbf{g}_{t},\widetilde{\mathbf{g}},etc...$ pour les m\'{e}triques (en
gras $\mathbf{g}$ pour les distinguer plus clairement des fonctions); $%
h,h^{\prime },h_{t},etc...$ pour les fonctions du premier membre
d\'{e}finissant l'op\'{e}rateur $\bigtriangleup _{\mathbf{g}}+h$; \thinspace 
$f,f^{\prime },f_{t},etc...$ pour celles du second membre; et $u,u_{t},etc$
pour les fonctions inconnues et les solutions de $(E_{h,f,\mathbf{g}%
}^{\prime })$.

On s'int\'{e}resse aux solutions minimisantes de $E_{h,f,\mathbf{g}}^{\prime
}$: on dira que $u\in H_{1}^{2}(M)$ est minimisante (ou extr\'{e}male) pour $%
(E_{h,f,\mathbf{g}}^{\prime })$ (ou par abus de langage, minimisante pour $h$%
) si pour la fonctionnelle 
\[
I_{h,\mathbf{g}}(w)=I_{h}(w)=\int_{M}\left| \nabla w\right| _{g}^{2}dv_{%
\mathbf{g}}+\int_{M}h.w{{}^{2}}dv_{\mathbf{g}} 
\]
on a 
\[
I_{h,\mathbf{g}}(u)=\stackunder{w\in \mathcal{H}_{f}}{\inf }I_{h,\mathbf{g}%
}(w):=\lambda _{h,f,\mathbf{g}}:=\lambda _{h} 
\]
o\`{u} 
\[
\mathcal{H}_{f}=\{w\in H_{1}^{2}(M)/\int_{M}f\left| w\right| ^{\frac{2n}{n-2}%
}dv_{\mathbf{g}}=1\}. 
\]
Alors, quitte \`{a} la multiplier par une constante, $u$ est $C^{\infty }$,
strictement positive, et est solution de: 
\[
(E_{h})=(E_{h,f})=(E_{h,f,\mathbf{g}}):\,\triangle _{\mathbf{g}%
}u+h.u=\lambda _{h}.f.u^{\frac{n+2}{n-2}}\text{ avec en plus }\int_{M}fu^{%
\frac{2n}{n-2}}dv_{\mathbf{g}}=1. 
\]
On sait, d'apr\`{e}s le th\'{e}or\`{e}me de Th. Aubin dont nous avons
parl\'{e} plus haut, qu'on a toujours 
\[
\lambda _{h,f,\mathbf{g}}\leq \frac{1}{K(n,2){{}^{2}}(\stackunder{M}{Sup}f)^{%
\frac{n-2}{n}}}\,. 
\]
A partir de ces rappels, nous allons maintenant donner trois d\'{e}finitions
\'{e}quivalentes des fonctions critiques:

\begin{definition}
\textbf{\ }\textit{Avec les notations pr\'{e}c\'{e}dentes:}

\begin{itemize}
\item  $h$\textit{\ est faiblement critique pour }$f$\textit{\ et }$\mathbf{g%
}$ \textit{si }$\lambda _{h,f,\mathbf{g}}=\frac{1}{K(n,2){{}^{2}}(%
\stackunder{M}{Sup}f)^{\frac{n-2}{n}}}$

\item  $h$\textit{\ est sous-critique pour }$f$\textit{\ et }$\mathbf{g}$ 
\textit{si }$\lambda _{h,f,\mathbf{g}}<\frac{1}{K(n,2){{}^{2}}(\stackunder{M%
}{Sup}f)^{\frac{n-2}{n}}}$

\item  $h$\textit{\ est critique pour }$f$\textit{\ et }$\mathbf{g}$ \textit{%
si }$h$\textit{\ est faiblement critique et si pour toute fonction continue }%
$k\leq h,\,k\neq h$\textit{, telle que }$\bigtriangleup _{\mathbf{g}}+k$ est
coercif,\textit{\ }$k$\textit{\ est sous-critique.}
\end{itemize}
\end{definition}
\medskip

Le th\'{e}or\`{e}me de Th. Aubin rappel\'{e} plus haut permet d'obtenir une
d\'{e}finition \'{e}quivalente: d'apr\`{e}s ce th\'{e}or\`{e}me, si $h$ est
sous-critique, $(E_{h,f,\mathbf{g}})$ a une solution minimisante. On
remarque alors que si $h$ est faiblement critique et a une solution
minimisante $u$ , $h$ est critique. En effet, dans ce cas, puisque $u>0$,
pour toute fonction continue $k\leq h,\,k\neq h$, on a 
\[
I_{k,\mathbf{g}}(u)<I_{h,\mathbf{g}}(u)=\frac{1}{K(n,2){{}^{2}}(\stackunder{M%
}{Sup}f)^{\frac{n-2}{n}}} 
\]
donc $k$ est sous-critique. Par cons\'{e}quent, si $h$ est critique, pour
toute fonction $k^{\prime }\geq h,\,k^{\prime }\neq h,\,(E_{k^{\prime }})$
ne peut pas avoir de solution minimisante, sinon $k^{\prime }$ serait
faiblement critique avec une solution minimisante, donc critique et alors $h$
serait sous-critique.

D'o\`{u} la d\'{e}finition \'{e}quivalente \`{a} la pr\'{e}c\'{e}dente: 

\begin{definition}
\textbf{\ }\textit{\ }$h$\textit{\ est une fonction critique pour }$f$%
\textit{\ et }$\mathbf{g}$ \textit{si :}

\begin{itemize}
\item  \textit{pour toute fonction continue }$k\leq h,\,k\neq h$\textit{,
telle que }$\bigtriangleup _{\mathbf{g}}+k$ \textit{est coercif} (\textit{ce
qui est le cas d\`{e}s que }$k$ \textit{est assez proche de }$h$\textit{\
dans} $C^{0}$),\textit{\ }$\left( E_{k}\right) $\textit{\ a une solution
minimisante,}

\item  \textit{pour toute fonction continue }$k^{\prime }\geq h,\,k^{\prime
}\neq h$\textit{\ }$,\left( E_{k^{\prime }}\right) $\textit{\ n'a pas de
solution minimisante.}
\end{itemize}
\end{definition}
\textbf{Les fonctions critiques sont donc introduites comme s\'{e}parant les
fonctions donnant une \'{e}quation qui admet des solutions minimisantes et
les fonctions donnant une \'{e}quation qui ne peut en avoir. L'objectif de
ce travail est donc d'utiliser la notion de fonction critique pour
\'{e}tudier l'existence de solutions minimisantes aux \'{e}quations (}$%
E_{h,f,\mathbf{g}}$\textbf{). En particulier, un probl\`{e}me central sera
l'existence de solutions minimisantes pour les \'{e}quations ``critiques'',
c'est-\`{a}-dire d\'{e}finies par une fonction critique.}

\smallskip Faisons quelques remarques simples sur les fonctions critiques
(voir appendice A):

L'op\'{e}rateur $\bigtriangleup _{\mathbf{g}}+h$ est forc\'{e}ment coercif
pour $h\in C{{}^{0}}(M)$ si $\lambda _{h,f,\mathbf{g}}=\frac{1}{K(n,2){{}^{2}%
}(\stackunder{M}{Sup}f)^{\frac{n-2}{n}}}$ et si $f>0$ sur $M$.
Ind\'{e}pendamment, si $\bigtriangleup _{\mathbf{g}}+h$ est coercif, toute
fonction $k$ continue assez proche de $h$ dans $C^{0}(M)$ est telle que $%
\bigtriangleup _{\mathbf{g}}+k$ est coercif. (voir appendice A).

Si $h$ est critique, n\'{e}cessairement il existe $x\in M$ tel que $h(x)>0$
(il suffit de tester la fonction 1).

On doit dire critique pour $f$ \textit{et\ }pour la m\'{e}trique $\mathbf{g}$%
, car $\mathbf{g}$ intervient fondamentalement dans le laplacien $%
\bigtriangleup _{\mathbf{g}}.$ Nous sous-entendrons $f$ et/ou $\mathbf{g\,}$%
quand il n'y aura pas d'ambigu\"{i}t\'{e}.

On voit sur cette d\'{e}finition que se pose naturellement et
fondamentalement le probl\`{e}me de l'existence de solutions minimisantes
pour les fonctions critiques; ce sera l'objet de notre premier
th\'{e}or\`{e}me.

Enfin, on peut pr\'{e}senter la premi\`{e}re d\'{e}finition d'une autre
mani\`{e}re:

\begin{definition}

\textbf{\ }On \textit{dira que }$h$\textit{\ est critique pour }$f$\textit{\
et }$\mathbf{g}$ \textit{si on a } 
\[
\forall u\in H_{1}^{2}(M)\text{ : }\left( \int_{M}f\left| u\right| ^{\frac{2n%
}{n-2}}dv_{\mathbf{g}}\right) ^{\frac{n-2}{n}}\leq K(n,2){{}^{2}}(%
\stackunder{M}{Sup}f)^{\frac{n-2}{n}}(\int_{M}\left| \nabla u\right| _{%
\mathbf{g}}^{2}dv_{\mathbf{g}}+\int_{M}h.u{{}^{2}}dv_{\mathbf{g}}) 
\]
\textit{et si cette proposition n'est plus vraie pour toute fonction
continue }$k\leq h,\,k\neq h$\textit{\ mise \`{a} la place de }$h$.
\end{definition}
C'est cette vision qui conduisit E. Hebey et M. Vaugon \`{a}
introduire la notion de fonction critique dans l'\'{e}tude de la meilleure
seconde constante; en effet, d'apr\'{e}s cette d\'{e}finition, on peut voir
les fonctions critiques comme les ``meilleures fonctions'', au lieu des
meilleures constantes, dans l'in\'{e}galit\'{e} ci-dessus. On voit l\`{a}
aussi que la recherche d'information sur $B_{0}(\mathbf{g})$ am\`{e}ne \`{a}
\'{e}tudier les \'{e}quations $\left( E_{h,f,\mathbf{g}}\right) $.

D'ailleurs, $B_{0}(\mathbf{g})K(n,2){{}^{-2}}$ est un exemple fondamental.
En effet, par d\'{e}finition $B_{0}(\mathbf{g})K(n,2){{}^{-2}}$ est toujours
une fonction (constante) faiblement critique pour toute fonction $f$ et
toute m\'{e}trique $\mathbf{g}$. Si de plus $f\equiv 1$, si $(M,\mathbf{g})$
n'est pas conform\'{e}ment diff\'{e}omorphe \`{a} la sph\`{e}re standard et
est telle que $S_{\mathbf{g}}=cste$, alors $B_{0}(\mathbf{g})K(n,2)^{-2}$
est une fonction critique, c'est la plus petite fonction critique constante;
on le montre simplement (Voir Appendice A), mais en utilisant deux tr\`{e}s
gros th\'{e}or\`{e}mes, le th\'{e}or\`{e}me de Yamabe (!), et le th\'{e}or%
\`{e}me suivant de Z.Djadli et O.\ Druet, issu d'un article qui sera
fondamental pour nous dans la suite $\left[ 9\right] $. Notons d'abord qu'il
est connu $\left[ 7\right] $ que: 
\begin{equation}
B_{0}(\mathbf{g})\geq \max (\frac{n-2}{4(n-1)}K(n,2){{}^{2}}\stackunder{M}{%
\max }S_{\mathbf{g}},Vol_{\mathbf{g}}(M)^{-\frac{2}{n}})  
\end{equation}

Z.Djadli et O.Druet ont montr\'{e} alors que l'une des deux assertions
suivantes devait \^{e}tre vraie si $\dim M=n\geq 4$:

a/ $B_{0}(\mathbf{g})=\frac{n-2}{4(n-1)}K(n,2){{}^{2}}\stackunder{M}{\max }%
S_{\mathbf{g}}$

b/ (1.2) poss\`{e}de des fonctions extr\'{e}males.

D'autres questions se sont alors pos\'{e}es naturellement:

\begin{itemize}
\item  peut-on avoir a/ sans b/ ?

\item  peut-on avoir b/ sans a/ ?

\item  peut-on avoir a/ et b/ simultan\'{e}ment ?
\end{itemize}

E. Hebey et M. Vaugon [20] ont introduit les fonctions critiques pour r\'{e}%
pondre \`{a} ces questions. Notre but est d'utiliser cette notion, en
l'adaptant, pour l'\'{e}tude des \'{e}quations $E_{h,f,\mathbf{g}}.$

Il est tr\`{e}s important de remarquer la propri\'{e}t\'{e} suivante des
fonctions critiques: elles se transforment dans les changements de
m\'{e}trique conformes exactement comme la courbure scalaire. En effet, soit 
$u\in C^{\infty }(M),\,u>0$ et $\mathbf{g}^{\prime }=u^{\frac{4}{n-2}}%
\mathbf{g}$ une m\'{e}trique conforme \`{a} $\mathbf{g}$\textbf{.} Un calcul
montre que si on pose 
\[
h^{\prime }=\frac{\triangle _{\mathbf{g}}u+h.u}{u^{\frac{n+2}{n-2}}}\text{
c'est \`{a} dire }\triangle _{\mathbf{g}}u+h.u=h^{\prime }.u^{\frac{n+2}{n-2}%
} 
\]
on a pour toute fonction $w\in H_{1}^{2}(M):$%
\[
I_{h,\mathbf{g}}(w)=I_{h^{\prime },\mathbf{g}^{\prime }}(u^{-1}.w) 
\]
et 
\[
\int_{M}f\left| w\right| ^{2^{*}}dv_{\mathbf{g}}=\int_{M}f\left| \frac{w}{u}%
\right| ^{2^{*}}dv_{\mathbf{g}^{\prime }}\,. 
\]
De plus $\bigtriangleup _{\mathbf{g}}+h$ est coercif si, et seulement si, $%
\bigtriangleup _{\mathbf{g}^{\prime }}+h^{\prime }$ est coercif. Ceci
implique que $h$ est critique pour $f$ et $\mathbf{g}$ si, et seulement si,
la fonction 
\[
h^{\prime }=\frac{\triangle _{\mathbf{g}}u+h.u}{u^{\frac{n+2}{n-2}}} 
\]
est critique pour $f$ et $\mathbf{g}^{\prime }=u^{\frac{4}{n-2}}\mathbf{g}.$
Ou d'une autre mani\`{e}re, $h^{\prime }$ est critique pour $f$ et $\mathbf{g%
}^{\prime }=u^{\frac{4}{n-2}}\mathbf{g}$ si et seulement si la fonction 
\[
h=h^{\prime }u^{\frac{4}{n-2}}-\frac{\triangle _{\mathbf{g}}u}{u} 
\]
est critique pour $f$ et $\mathbf{g}$. Enfin, $w$ est une solution
minimisante pour $E_{h,f,\mathbf{g}}$ si, et seulement si, $\frac{w}{u}$ est
une solution minimisante pour $E_{h^{\prime },f,\mathbf{g}^{\prime }}$. Voir
l'Appendice A pour les d\'{e}tails.

Revenons \`{a} l'\'{e}valuation de $\stackunder{w\in \mathcal{H}_{f}^{+}}{%
\inf }J_{h}(w):=\lambda _{h,f,\mathbf{g}}\,$. T. Aubin introduisit dans la
fonctionnelle $J_{h,f,\mathbf{g}}$ les fonctions tests $\psi _{k}$
suivantes: 
$$\psi _{k}(Q)=\left\lbrace 
\begin{array}{c}
(\frac{1}{k}+r{{}^{2}})^{-\frac{n-2}{2}}-(\frac{1}{k}+\delta {{}^{2}})^{-%
\frac{n-2}{2}}\,\,\,\,si\,r<\delta \\ 
0\,\,\,\,\,\,\,\,si\,\,\,r\geq \delta
\end{array} \right. $$
o\`{u}: $\delta <injM$ (le rayon d'injectivit\'{e} de $M$), $P\,\in M$ est
un point fix\'{e}, $k\in \Bbb{N}^{*}$, et o\`{u} l'on note $r=d_{\mathbf{g}%
}(P,Q).$ Lorsque $P$ est un point o\`{u} $f$ est maximum sur $M$, le calcul
donne:

Si $n>4$:

\begin{center}
$J_{h,f,\mathbf{g}}(\psi _{k})=\frac{1}{K(n,2){{}^{2}}(\stackunder{M}{Sup}%
f)^{\frac{n-2}{n}}}\left\{ 1+\frac{1}{n(n-4)}\left( \frac{4(n-1)}{n-2}%
h(P)-S_{\mathbf{g}}(P)+\frac{n-4}{2}\frac{\bigtriangleup _{\mathbf{g}}f(P)}{%
f(P)}\right) \frac{1}{k}\right\} +o(\frac{1}{k})$
\end{center}

Et si $n$=4:

\begin{center}
$J_{h,f,\mathbf{g}}(\psi _{k})=\frac{1}{K(n,2){{}^{2}}(\stackunder{M}{Sup}%
f)^{\frac{1}{2}}}\left\{ 1+\left( 6h(P)-S_{\mathbf{g}}(P)\right) \frac{\log k%
}{8k}\right\} +o(\frac{1}{k})$
\end{center}
On en d\'{e}duit:

\textbf{Proposition 1:}

\textit{\ si }$h$\textit{\ est faiblement critique (donc en particulier si }$%
h$\textit{\ est critique) pour }$f$\textit{\ et }$\mathbf{g}$\textit{, comme 
}$\lambda _{h,f,\mathbf{g}}=\frac{1}{K(n,2){{}^{2}}(\stackunder{M}{Sup}f)^{%
\frac{n-2}{n}}}$\textit{, n\'{e}cessairement, si }$P$\textit{\ est un point
o\`{u} }$f$\textit{\ est maximum sur }$M,$\textit{\ on a: } 
\[
\text{\textit{Si }}\,n>4\,\,:\,\,\frac{4(n-1)}{n-2}h(P)\geq S_{\mathbf{g}%
}(P)-\frac{n-4}{2}\frac{\bigtriangleup _{\mathbf{g}}f(P)}{f(P)} 
\]
\[
\text{\textit{Si}}\,\,n=4:\,6h(P)\geq S_{\mathbf{g}}(P) 
\]

\textit{Remarque : si }$f=cste$\textit{\ cela signifie que sur tout }$M$%
\textit{\ on a }$\frac{4(n-1)}{n-2}h\geq S_{\mathbf{g}}.$

Nous allons maintenant donner les th\'{e}or\`{e}mes prouv\'{e}s dans notre
travail.

Le premier r\'{e}sultat concerne l'existence de solutions minimisantes pour
les fonctions critiques; comme nous l'avons vu dans leur deuxi\`{e}me
d\'{e}finition, elles sont en effet ``entre'' les fonctions qui ont des
solutions minimisantes et les fonctions qui n'ont pas de telles solutions.
Cette existence arrive comme corollaire d'un r\'{e}sultat plus
g\'{e}n\'{e}ral qui nous sera extr\^{e}mement utile pour la suite, et dont
la d\'{e}monstration est au coeur de notre travail. Les outils, l'id\'{e}e
et les difficult\'{e}s de cette d\'{e}monstration seront pr\'{e}sent\'{e}s
au chapitre 2, la d\'{e}monstration elle-m\^{e}me faisant l'objet du
chapitre 3. Nous ferons dans les chapitre 3 \`{a} 6 une hypoth\`{e}se
suppl\'{e}mentaire fondamentale sur $f$ et nous considèrerons des
vari\'{e}t\'{e}s $(M,\mathbf{g})$ de dimension supèrieure ou \'{e}gale
\`{a} 4.

\textbf{Hypoth\`{e}ses (H)}\textit{: On suppose que le Hessien de la
fonction }$f:M\rightarrow \Bbb{R}$, \textit{telle que }$\stackunder{M}{Sup}%
f>0$, \textit{est non d\'{e}g\'{e}n\'{e}r\'{e} en chaque point de maximum de 
}$f.$ \textit{En outre, les fonctions }$h$ \textit{consid\'{e}r\'{e}es sont
telles que }$\bigtriangleup _{\mathbf{g}}+h$ \textit{est coercif et l'on
suppose }$\dim M\geq 4$.\textit{\ On parlera des hypoth\`{e}ses (}\textbf{H}$%
_{f}$\textit{) pour d\'{e}signer celles concernant la fonction }$f$\textit{.}

\textbf{Th\'{e}or\`{e}me 1':}

\textbf{\ }\textit{Si }$h$\textit{\ est critique pour }$f$\textit{\ et }$%
\mathbf{g}$\textit{, (}$h,\,f$ et $\mathbf{g}$ v\'{e}rifiant \textbf{(H})), 
\textit{et si en tout point }$P$\textit{\ o\`{u} }$f$\textit{\ est maximum
sur }$M$\textit{\ on a: }$\frac{4(n-1)}{n-2}h(P)>S_{\mathbf{g}}(P)-\frac{n-4%
}{2}\frac{\bigtriangleup _{\mathbf{g}}f(P)}{f(P)}$\textit{\ alors il existe
une solution minimisante pour }$h$\textit{, c'est \`{a} dire minimisante
pour }$(E_{h,f,\mathbf{g}})$.

\medskip
autrement dit, pour reprendre la formulation du th\'{e}or\`{e}me de
Djadli-Druet:

\medskip

\textbf{Th\'{e}or\`{e}me 1'':}

\textbf{\ }\textit{Si }$h$\textit{\ est critique pour }$f$ \textit{et }$%
\mathbf{g}$\textit{, (}$h,\,f$ et $\mathbf{g}$ v\'{e}rifiant \textbf{(H})), 
\textit{l'une des deux assertions suivantes est vraie:}

\begin{itemize}
\item  \textit{il existe un} \textit{point }$P$\textit{\ o\`{u} }$f$\textit{%
\ est maximum sur }$M$ \textit{tel que} $\frac{4(n-1)}{n-2}h(P)=S_{\mathbf{g}%
}(P)-\frac{n-4}{2}\frac{\bigtriangleup _{\mathbf{g}}f(P)}{f(P)}$

\item  $(E_{h}^{\prime })=(E_{h,f,\mathbf{g}}^{\prime }):\,\triangle _{%
\mathbf{g}}u+h.u=f.u^{\frac{n+2}{n-2}}$ \textit{a une solution minimisante.}
\end{itemize}
\medskip
Ce th\'{e}or\`{e}me sera en fait une cons\'{e}quence imm\'{e}diate du
r\'{e}sultat suivant, plus g\'{e}n\'{e}ral mais moins parlant. Il suffit de
prendre dans ce th\'{e}or\`{e}me la suite $h_{t}=h-t$ pour $t\stackunder{>}{%
\rightarrow }0$, ces fonctions \'{e}tant sous-critiques par d\'{e}finition,
pour obtenir le th\'{e}or\`{e}me 1' (ou 1'').

\begin{theorem}

\textbf{\ }\textit{Avec l'hypoth\`{e}se }\textbf{(H}),\textit{\ soit }$h$%
\textit{\ une fonction faiblement critique pour }$f$\textit{\ et }$\mathbf{g}%
.$\textit{\ Si $\frac{4(n-1)}{n-2}h(P)>S_{\mathbf{g}}(P)-\frac{n-4%
}{2}\frac{\bigtriangleup _{\mathbf{g}}f(P)}{f(P)}$ en tout point }$P$\textit{\ o\`{u} }$f$\textit{\ est maximum
sur }$M$\textit{\ et s' il existe
une famille de fonctions }$(h_{t})$, $h_{t}\lvertneqq  h$, $h_{t}$ \textit{\
sous-critique pour tout }$t$\textit{\ dans un voisinage de }$t_{0}\in \Bbb{R}
$\textit{, telle que }$h_{t}\stackunder{t\rightarrow t_{0}}{\rightarrow }%
h\,\,dans\,$\textit{\ }$C^{0,\alpha }$\textit{\ alors il existe une solution
minimisante pour }$h$\textit{\ (i.e. pour }$(E_{h,f,\mathbf{g}})$) \textit{%
et donc }$h$\textit{\ est critique pour }$f$\textit{\ et }$\mathbf{g}$%
\textit{.}
\end{theorem}
E.Hebey et M.Vaugon, dans le cadre de leur \'{e}tude sur $B_{0}(\mathbf{g)}$%
, ont montr\'{e} ce th\'{e}or\`{e}me dans le cas o\`{u} $f\equiv cste=1$ ,
la condition \'{e}tant alors que $\frac{4(n-1)}{n-2}h>S_{\mathbf{g}}$ sur
toute la vari\'{e}t\'{e}. La d\'{e}monstration de notre th\'{e}or\`{e}me
commence de la m\^{e}me mani\`{e}re et en particulier s'appuie sur l'article
de Z. Djadli et O. Druet [9], mais la pr\'{e}sence d'une fonction $f$ non
constante fait appara\^{i}tre de nouvelles difficult\'{e}s dans l'\'{e}tude
de ce qu'on appelle les ph\'{e}nom\`{e}nes de concentration, pr\'{e}sent\'{e}%
s au chapitre 2. De plus la m\'{e}thode d\'{e}velopp\'{e}e pour surmonter
cette difficult\'{e} apporte des \'{e}l\'{e}ments nouveaux sur ces ph\'{e}nom%
\`{e}nes de concentration et semble pouvoir s'appliquer \`{a} d'autres
questions similaires: voir Zo\'{e} Faget [15]. La d\'{e}monstration fera
l'objet du chapitre 3.

La question suivante, naturelle, consiste bien s\^{u}r \`{a} savoir s'il
existe des fonctions critiques! La r\'{e}ponse, affirmative, s'obtient comme
corollaire du th\'{e}or\`{e}me 1 ci-dessus:
\\
\\
\\

\begin{theorem}

\textit{Etant donn\'{e}es la vari\'{e}t\'{e} }$(M,\mathbf{g)}$ et la
fonction $f$, \textit{v\'{e}rifiant }\textbf{(H}), \textit{il existe une
infinit\'{e} de fonctions critiques }$h$ \textit{pour }$f$ \textit{et }$%
\mathbf{g}$ \textit{, qui v\'{e}rifient, en tout point }$P$\textit{\ o\`{u} }$f$
\textit{\ est maximum sur }$M$\textit{, }$\frac{4(n-1)}{n-2}h(P)>S_{\mathbf{g%
}}(P)-\frac{n-4}{2}\frac{\bigtriangleup _{\mathbf{g}}f(P)}{f(P)}$
\textit{. Ces fonctions critiques ont des fonctions extr\'{e}males.}
\end{theorem}
Il sera de plus montr\'{e} que si $\int_{M}fdv_{\mathbf{g}}>0$, il existe
des fonctions critiques strictement positives, ce qui aura son importance
pour la suite.

Ces th\'{e}or\`{e}mes furent les premiers de notre travail. Ils nous
amen\`{e}rent \`{a} adopter une vision un peu diff\'{e}rente des fonctions
critiques pour mettre en valeur leur int\'{e}r\^{e}t dans l'\'{e}tude des
\'{e}quations $(E_{h,f,\mathbf{g}})$. Tout d'abord, on constate que dans
l'\'{e}quation $(E_{h,f,\mathbf{g}}^{^{\prime }}):\,\triangle _{\mathbf{g}%
}u+h.u=f.u^{\frac{n+2}{n-2}}$, il y a trois donn\'{e}es fondamentales que
l'on peut faire varier: les fonctions $h$ et $f$, et la m\'{e}trique $%
\mathbf{g}$ (dans la classe des m\'{e}triques qui lui sont conformes);
remarquons que $h$ et $\mathbf{g}$ d\'{e}finissent l'op\'{e}rateur $%
\triangle _{\mathbf{g}}+h$. Notons en effet que si l'on change $\mathbf{g}$
en une m\'{e}trique $\mathbf{g}^{\prime }$ qui lui est conforme, par la
r\`{e}gle de transformation du laplacien conforme, l'\'{e}quation se
transforme en une autre exactement de m\^{e}me type: $\,\triangle _{\mathbf{g%
}^{\prime }}u+h^{\prime }.u=f^{\prime }.u^{\frac{n+2}{n-2}}.$ Il sera alors
int\'{e}ressant de parler de triplet critique en adoptant la d\'{e}finition
suivante:

\begin{center}
$(h,f,\mathbf{g)}$ \textit{est un triplet critique si} $h$\textit{\ est
critique pour }$f$\textit{\ et }$\mathbf{g}$.
\end{center}

De mani\`{e}re analogue on parlera de triplet sous-critique ou faiblement
critique. On dira par ailleurs que le triplet $(h,f,\mathbf{g)}$ a des
solutions minimisantes si l'equation $(E_{h,f,\mathbf{g}})$ a des solutions
minimisantes. La question que nous nous poserons sera alors la suivante:
\'{e}tant fix\'{e}es deux des donn\'{e}es du triplet, peut-on trouver la
troisi\`{e}me pour obtenir un triplet critique. Nous montrerons que cela est
en effet possible, ce qui montrera \textit{a posteriori} que fixer une des
trois donn\'{e}es et chercher les deux autres pour obtenir un triplet
critique est aussi possible, ce qui n'est pas \'{e}vident \`{a} priori.
Ainsi par exemple, le probl\`{e}me de l'existence de fonctions critiques
revient \`{a} se fixer la fonction $f$ et la m\'{e}trique $\mathbf{g}$ et
\`{a} chercher $h$ pour que $(h,f,\mathbf{g)}$ soit un triplet critique.
Nous nous poserons alors les deux autres questions possibles, en fixant
d'abord les fonctions $h$ et $f$ et en cherchant la m\'{e}trique $\mathbf{g}$%
, puis en fixant la fonction $h$ et la m\'{e}trique $\mathbf{g}$ et en
cherchant la fonction $f$.

Nous obtenons des r\'{e}ponses compl\`{e}tes, sous la forme des deux th\'{e}%
or\`{e}mes suivants: 

\begin{theorem}
\textit{Soient donn\'{e}es la vari\'{e}t\'{e} (}$M,\mathbf{g}$\textit{) et
deux fonctions }$h^{\prime }$ \textit{et }$f$\textit{\ v\'{e}rifiant les
hypoth\`{e}ses }\textbf{(H}).\textit{\ Alors il existe une m\'{e}trique }$%
\mathbf{g}^{\prime }$ conforme \`{a} $\mathbf{g}$ \textit{telle que }$%
(h^{\prime },f,\mathbf{g}^{\prime })$\textit{\ soit un triplet critique, ou,
pour reprendre la pr\'{e}sentation premi\`{e}re, il existe une m\'{e}trique }%
$\mathbf{g}^{\prime }$ conforme \`{a} $\mathbf{g}$ \textit{telle que }$%
h^{\prime }$\textit{\ soit critique pour }$f$\textit{\ et }$\mathbf{g}%
^{\prime }$. \textit{De plus }$(h^{\prime },f,\mathbf{g}^{\prime })$\textit{%
\ a des solutions minimisantes.}
\end{theorem}
Ce th\'{e}or\`{e}me a \'{e}t\'{e} montr\'{e} par E. Humbert et M. Vaugon
dans le cas $f=cste=1$. Leur d\'{e}monstration (chapitre 5) passe dans notre
cadre plus g\'{e}n\'{e}ral une fois prouv\'{e} que l'on peut supposer
l'existence d'une fonction $h>0$ critique pour $f$\textit{\ }et\textit{\ }$%
\mathbf{g}$, r\'{e}sultat que nous avons cit\'{e} plus haut avec l'existence
de fonctions critiques et qui sera montr\'{e} dans le chapitre correspondant.

Enfin, la derni\`{e}re question nous am\`{e}ne au r\'{e}sultat suivant. Nous
avons besoin pour ce th\'{e}or\`{e}me d'augmenter un peu la dimension de la
vari\'{e}t\'{e}:
\\

\begin{theorem}

\textit{Soient donn\'{e}es la vari\'{e}t\'{e} (}$M,\mathbf{g}$\textit{), }$%
\dim M\geqslant 5$\textit{, et une fonction }$h$\textit{\ telle que
l'op\'{e}rateur }$\triangle _{\mathbf{g}}+h$\textit{\ soit coercif. Alors,
il existe une fonction }$f$ (\textit{v\'{e}rifiant }\textbf{(H}$_{f}$))%
\textit{\ telle que }$(h,f,\mathbf{g})$\textit{\ soit critique avec des
solutions minimisantes, si et seulement si }$(h,1,\mathbf{g})$\textit{\ est
sous-critique (ou }$1$ \textit{est la fonction constante }$1$).
\end{theorem}
La d\'{e}monstration, assez difficile et s'appuyant sur le principe
d\'{e}velopp\'{e} au chapitre 3, se fait en deux parties et apporte quelques
r\'{e}sultats interm\'{e}diaires int\'{e}ressants. En particulier elle nous
amena \`{a} quelques remarques suppl\'{e}mentaires sur la notion de fonction
critique.

Tout d'abord, on voit sur les d\'{e}finitions (en utilisant par exemple la
fonctionnelle $J$) que si $h$ est critique pour une fonction $f$, $h$ est
critique pour $c.f$ pour toute constante $c>0$. Il en va de m\^{e}me dans le
cas o\`{u} $h$ est sous-critique ou faiblement critique. Il serait donc plus
naturel de dire que $h$ est critique pour la ``classe'' de $f$, not\'{e}e $%
[f]$, o\`{u} $[f]=\{c.f\,/\,c>0\}$, et de consid\'{e}rer des triplets $%
(h,[f],\mathbf{g)}$. Ainsi, par exemple, dans toute ``classe'' on peut
choisir un repr\'{e}sentant tel que $Supf=1$; et pour comparer des triplets $%
(h,f,\mathbf{g})$ et $(h,f^{\prime },\mathbf{g})$ entre eux, il faut
supposer que $Supf=Supf^{\prime }$. Par ailleurs on remarque que la valeur
de $\frac{\bigtriangleup _{\mathbf{g}}f(P)}{f(P)}$ est constante sur une
classe $[f]$.

Dans la d\'{e}monstration de ce dernier th\'{e}or\`{e}me, nous avons
\'{e}t\'{e} amen\'{e} \`{a} faire ``varier'' la fonction $f$ du second
membre des \'{e}quations $(E_{h,f,\mathbf{g}})$. Cela sugg\'{e}ra une autre
d\'{e}finition possible des fonctions critiques:

\begin{definition}

\textit{Soient donn\'{e}es la vari\'{e}t\'{e} (}$M,\mathbf{g}$\textit{), }$%
\dim M\geqslant 3$\textit{, et une fonction }$h$\textit{\ telle que
l'op\'{e}rateur }$\triangle _{\mathbf{g}}+h$\textit{\ soit coercif. On
consid\`{e}re une fonction }$f\in C^{\infty }(M)$, telle que $Supf>0$.%
\textit{\ On dira que }$f$ \textit{est critique pour }$h$\textit{\ si:}

\begin{itemize}
\item  \textit{a/: }$\lambda _{h,f,\mathbf{g}}=\frac{1}{K(n,2){{}^{2}}(%
\stackunder{M}{Sup}f)^{\frac{n-2}{n}}}$

\item  \textit{b/: pour toute fonction }$f^{\prime }$ \textit{telle que } $%
Supf=Supf^{\prime }$ \textit{et }$f^{\prime }\gneqq f$, \textit{\ }$\\
\lambda_{h,f^{\prime },\mathbf{g}}<\frac{1}{K(n,2){{}^{2}}(\stackunder{M}{Sup}%
f^{\prime })^{\frac{n-2}{n}}}$

\item  \textit{Remarque: si }$Supf=Supf^{\prime }$ \textit{et }$f^{\prime
}\lvertneqq f$, \textit{\ }$\lambda _{h,f^{\prime },\mathbf{g}}=\frac{1}{%
K(n,2){{}^{2}}(\stackunder{M}{Sup}f^{\prime })^{\frac{n-2}{n}}}$ \textit{%
puisque }$J_{h,f^{\prime },\mathbf{g}}(w)\geqslant J_{h,f,\mathbf{g}}(w)$%
\textit{\ pour toute fonction }$w$.
\end{itemize}
\end{definition}
D'après ce que nous avons dit juste avant, il faut bien comparer dans
cette d\'{e}finition des fonctions de m\^{e}me Sup; ce sont en fait les
classes $[f]$ et $[f^{\prime }]$ qui importent, et l\`{a} aussi on doit dire
que c'est $[f]$ qui est critique pour $h$. Il est alors naturel de poser la
question suivante:

\begin{center}
$f$ \textit{est-elle critique pour }$h$\textit{\ si, et seulement si, }$h$ 
\textit{est critique pour }$f$ ?
\end{center}

Cette question semble assez difficile. Rappelons que dans les deux cas on a
toujours en tout point $P$\ o\`{u} $f$\ est maximum sur $M$\ :\textit{\ }$%
\frac{4(n-1)}{n-2}h(P)\geqslant S_{\mathbf{g}}(P)-\frac{n-4}{2}\frac{%
\bigtriangleup _{\mathbf{g}}f(P)}{f(P)}$\textit{\ }

Nous obtenons le th\'{e}or\`{e}me suivant:

\begin{theorem}

\textit{Soient donn\'{e}es la vari\'{e}t\'{e} (}$M,\mathbf{g}$\textit{), }$%
\dim M\geqslant 5$\textit{, et une fonction }$h$\textit{\ telle que
l'op\'{e}rateur }$\triangle _{\mathbf{g}}+h$\textit{\ soit coercif. On
consid\`{e}re une fonction }$f\in C^{\infty }(M)$, telle que $Supf>0$\textit{%
\ et v\'{e}rifiant (}\textbf{H}$_{f}$).\textit{\ On suppose de plus qu'en
tout point }$P $\textit{\ o\`{u} }$f$\textit{\ est maximum sur }$M$\textit{\
: }$\frac{4(n-1)}{n-2}h(P)>S_{\mathbf{g}}(P)-\frac{n-4}{2}\frac{%
\bigtriangleup _{\mathbf{g}}f(P)}{f(P)}$\textit{. Alors }$f$ \textit{est
critique pour }$h$\textit{\ si, et seulement si, }$h$ \textit{est critique
pour }$f$ $.$
\end{theorem}
Ce th\'{e}or\`{e}me s'obtient comme cons\'{e}quence rapide de la
d\'{e}monstration du chapitre 6, mais il semble difficile \`{a} obtenir sans
le travail des chapitres 3 et 6.

Illustrons ces derniers r\'{e}sultats avec quelques exemples.

Une des applications très importante des \'{e}quations $(E_{h,f,\mathbf{g%
}})$ est l'\'{e}tude des courbures scalaires possibles dans une classe
conforme: \'{e}tant donn\'{e}es une vari\'{e}t\'{e} \textit{(}$M,\mathbf{g}$%
\textit{) }et une fonction $f\in C^{\infty }(M)$, $f$ est-elle la courbure
scalaire d'une m\'{e}trique conforme \`{a} $\mathbf{g}$ ? Les
th\'{e}or\`{e}mes de Th. Aubin montrent que si $f$ est sous-critique pour $%
S_{\mathbf{g}}\, $alors $f$ est une courbure scalaire. Le th\'{e}or\`{e}me 4
appliqu\'{e} \`{a} $h=\frac{n-2}{4(n-1)}S_{\mathbf{g}}$ montre que

\textit{Sur une vari\'{e}t\'{e} (}$M,\mathbf{g}$) \textit{non
conform\'{e}ment diff\'{e}omorphe \`{a} la sph\`{e}re, il existe des
courbures scalaires de m\'{e}triques conformes \`{a} }$\mathbf{g}$\textit{\
qui sont seulement faiblement critiques, puisque critiques.}

Il serait int\'{e}ressant de caract\'{e}riser ces m\'{e}triques critiques.

Une des cons\'{e}quences de l'article de E. Hebey et M. Vaugon [20] est
qu'il existe beaucoup de vari\'{e}t\'{e}s \textit{(}$M,\mathbf{g}$\textit{), }%
$\dim M\geqslant 4$ qui n'ont pas de fonction critique \textit{constante}
pour $f=1$. En appliquant le th\'{e}or\`{e}me 4 \`{a} toute $h=cste=c<B_{0}(%
\mathbf{g})K(n,2){{}^{-2}}$ on voit que sur toute vari\'{e}t\'{e} il existe
des fonctions $f$ pour lesquelles $c$ sera critique. Par contre \'{e}tant
donn\'{e}es \textit{(}$M,\mathbf{g}$\textit{) et }$f$, nous ne savons pas
affirmer l'existence de fonctions critiques constantes.

Le chapitre 6 traite de la dimension 3. Il faut en effet remarquer que toute
l'\'{e}tude pr\'{e}c\'{e}dente portait sur des vari\'{e}t\'{e}s de dimension 
$\geqslant 4$. De plus la dimension 4 elle-m\^{e}me pr\'{e}sente une
particularit\'{e} puisque le terme $\frac{n-4}{2}\frac{\bigtriangleup _{%
\mathbf{g}}f(P)}{f(P)}\ $dispara\^{i}t. D'ailleurs, bien que les th\'{e}or%
\`{e}mes restent valables pour $\dim M=4$, le r\'{e}sultat fondamental que
nous obtenons sur les ph\'{e}nom\`{e}nes de concentration, pr\'{e}sent\'{e}
au chapitre suivant, n'est valable que pour $\dim M\geqslant 5$. Le cas de
la dimension 3 est, lui, radicalement diff\'{e}rent et sera pr\'{e}sent\'{e}
au chapitre 6. O. Druet $\left[ 10\right] $ a trait\'{e} le cas $f=cste=1$.
L'introduction d'une fonction non constante n'apporte ici pas de difficult%
\'{e}s, nous reprendrons donc rapidement la d\'{e}monstration d'O. Druet
pour obtenir sa g\'{e}n\'{e}ralisation au cas $f$ non constante. Cette
dimension fait intervenir de fa\c{c}on fondamentale la fonction de Green de
l'op\'{e}rateur $\triangle _{\mathbf{g}}+h$. Si cet op\'{e}rateur est
coercif, il existe une unique fonction 
\[
G_{h}:M\times M\backslash \{(x,x),x\in M\}\rightarrow \Bbb{R}
\]
sym\'{e}trique et strictement positive telle que, au sens des distributions,
on a $\forall x\in M:$%
\[
\bigtriangleup _{\mathbf{g},y}G_{h}(x,y)+h(y)G_{h}(x,y)=\delta _{x}.
\]

En dimension 3, pour un point $x\in M$, et pour $y$ proche de $x$, $G_{h}$
peut se mettre sous la forme: 
\[
G_{h}(x,y)=\frac{1}{\omega _{2}d_{\mathbf{g}}(x,y)}+M_{h}(x)+o(1) 
\]
o\`{u} $o(1)$ est \`{a} prendre pour $y\rightarrow x$. On appelle $M_{h}(x)$
la masse de la fonction de Green au point $x$.

On obtient alors facilement \`{a} partir de la m\'{e}thode d'Olivier Druet
les trois r\'{e}sultats suivants, analogues au cas $f=cste$ trait\'{e} dans
son article $\left[ 10\right] $:

\begin{theorem}

\textit{Soient (}$M,\mathbf{g}$\textit{) une vari\'{e}t\'{e} compacte de
dimension 3 et une fonction }$f\in C^{\infty }(M)$ \textit{telle que} $%
Supf>0 $ (\textit{l'hypoth\`{e}se (}\textbf{H}$_{f}$\textit{) n'est pas
n\'{e}cessaire). Alors pour toute fonction }$h$ \textit{faiblement critique
pour }$f$\textit{\ et }$\mathbf{g}$\textit{, et pour tout }$x\in M\,$\textit{%
o\`{u} }$f$\textit{\ est maximum sur} $M$,$\,$\textit{on a }$%
M_{h}(x)\leqslant 0$\textit{.}
\end{theorem}
La condition \textit{\ }$M_{h}(x)\leqslant 0\,$appara\^{i}t comme l'analogue
de la condition $\frac{4(n-1)}{n-2}h(P)\geqslant S_{\mathbf{g}}(P)-\frac{n-4%
}{2}\frac{\bigtriangleup _{\mathbf{g}}f(P)}{f(P)}$\textit{\ }que l'on avait
en dimension $\geqslant 4$. La particularit\'{e} de la dimension 3 est alors
d'offrir des fonctions critiques de toutes les formes:

\begin{theorem}

\textit{Soient (}$M,\mathbf{g}$\textit{) une vari\'{e}t\'{e} compacte de
dimension 3 et une fonction }$f\in C^{\infty }(M)$ \textit{telle que} $%
Supf>0 $ (\textit{l'hypoth\`{e}se (}\textbf{H}$_{f}$\textit{)} \textit{\
n'est pas n\'{e}cessaire). Pour toute fonction }$h$, \textit{posons }$%
B(h)=\inf \{B/$ $h+B\,\,est\,\,faiblement\,\,critique\,\,pour\,\,f\}$. 
\textit{Alors }$h+B(h)$\textit{\ est une fonction critique pour }$f$\textit{.%
}
\end{theorem}
Enfin, en ce qui concerne l'existence de fonctions extr\'{e}males, on a le
th\'{e}or\`{e}me suivant:

\begin{theorem}

\textit{Soient (}$M,\mathbf{g}$\textit{) une vari\'{e}t\'{e} compacte de
dimension 3 et une fonction }$f\in C^{\infty }(M)$ \textit{telle que} $%
Supf>0 $ (\textit{l'hypoth\`{e}se (}\textbf{H}$_{f}$\textit{)} \textit{\
n'est pas n\'{e}cessaire). Soit }$h$ \textit{une fonction critique pour }$f$%
\textit{\ et }$\mathbf{g}$.\textit{\ Alors au moins l'une des deux
conditions suivantes est remplie:}

\begin{itemize}
\item  \textit{a/: Il existe }$x\in M\,$\textit{\thinspace o\`{u} }$f$%
\textit{\ est maximum sur} $M$ \textit{tel que }$M_{h}(x)=0$\textit{.}

\item  \textit{b/: }$(E_{h,f,\mathbf{g}})$\textit{\ a des solutions extr\'{e}%
males.}
\end{itemize}
\end{theorem}
Pour finir, le chapitre 8 traite du possible affaiblissement des
hypoth\`{e}ses des th\'{e}or\`{e}mes des chapitres 3 \`{a} 6, \`{a} savoir
l'hypoth\`{e}se (\textbf{H}$_{f}$\textit{), }et l'exigence d'une
in\'{e}galit\'{e} stricte dans la condition $\frac{4(n-1)}{n-2}h(P)>S_{%
\mathbf{g}}(P)-\frac{n-4}{2}\frac{\bigtriangleup _{\mathbf{g}}f(P)}{f(P)}$
en tout point $P$\ o\`{u} $f$\ est maximum sur $M$\textit{\ . }%
Malheureusement, seules des r\'{e}ponses partielles sont obtenues. Le
chapitre pr\'{e}sente \'{e}galement quelques questions qu'il nous semble
int\'{e}ressant de consid\'{e}rer \`{a} l'issue de ce travail.

Les appendices qui suivent regroupent quelques d\'{e}monstrations que nous
voulions inclure par souci de compl\'{e}tude mais que nous avons
pr\'{e}f\'{e}r\'{e} reporter pour ne pas interrompre les raisonnements des
chapitres principaux.

\chapter{Trois outils fondamentaux: le point de concentration, le changement
d'\'{e}chelle, le processus d'it\'{e}ration. Principe des d\'{e}monstrations}
\pagestyle{myheadings}\markboth{\textbf{Trois outils fondamentaux. Principe des 
démonstrations.}}{\textbf{Trois outils fondamentaux. Principe des démonstrations.}}
Nous voulons isoler ici la pr\'{e}sentation de trois ``outils'' que nous
utiliserons constamment dans l'\'{e}tude des \'{e}quations $E_{h,f,\mathbf{g}%
}$. Ces outils ont \'{e}t\'{e} d\'{e}velopp\'{e}s par plusieurs auteurs
depuis M. Vaugon [31] et P.L. Lions [26], en particulier E. Hebey, O. Druet,
F. Robert, sur les travaux desquels nous nous appuierons.

\section{Le point de concentration}

Pour prouver l'existence de solutions $u>0$ \`{a} nos \'{e}quations 
\[
(E_{h,f,\mathbf{g}}):\,\triangle _{\mathbf{g}}u+h.u=\lambda .f.u^{\frac{n+2}{%
n-2}}, 
\]
le principe consistera le plus souvent \`{a} construire une famille
d'\'{e}quations poss\'{e}dant des solutions minimisantes $u_{t}>0$: 
\[
E_{t}:\,\,\triangle _{\mathbf{g}}u_{t}+h_{t}.u_{t}=\lambda _{t}.f.u_{t}^{%
\frac{n+2}{n-2}} 
\]
avec 
\[
h_{t}\rightarrow h\,\,\,\,dans\,\,\,\,C^{0,\alpha }(M) 
\]
et $\lambda _{t}\rightarrow \lambda $ une ``suite'' convergente de
r\'{e}els, de telle sorte que l'on ait pour une fonction $u\in
H_{1}^{2}\,:\, $ $u_{t}\rightarrow u\,$ fortement dans des espaces $L^{p}$ , 
$p<2^{*}$, et $u_{t}\rightharpoondown u$ faiblement\thinspace dans $%
H_{1}^{2},$ avec une contrainte de la forme 
\[
\int_{M}f.u_{t}^{2^{*}}dv_{\mathbf{g}}=1. 
\]
Pour fixer un peu les choses nous supposerons le plus souvent dans la suite
que toutes ces convergences sont \`{a} prendre pour $t\rightarrow t_{0}=1$;
et pour simplifier le vocabulaire nous parlerons de suite $(u_{t})$, etc.
bien que ces familles soient index\'{e}es par des r\'{e}els. La
difficult\'{e} sera de prouver que $u$ est non identiquement nulle, car dans
ce cas, par le principe du maximum, $u>0$, et il s'ensuit que $u$ est
solution minimisante de $(E_{h,f,\mathbf{g}})$. Nous proc\`{e}derons par
contradiction en supposant $u\equiv 0$. L'id\'{e}e est alors que, \`{a}
cause de la contrainte $\int_{M}f.u_{t}^{2^{*}}=1,$ toute la ``masse'' des
fonctions $u_{t}$, qui convergent vers 0 dans $L^{p}$ , $p<2^{*}$, se
concentre autour d'un point de la vari\'{e}t\'{e}. On pose ainsi:
\begin{definition}
\textit{\ }$x_{0}\in M$\textit{\ est un point de
concentration de }$(u_{t})$\textit{\ si pour tout }$\delta >0$\textit{: } 
\[
\stackunder{t\rightarrow t_{0}}{\lim \sup }\int_{B(x_{0},\delta
)}u_{t}^{2^{*}}dv_{\mathbf{g}}>0 
\]
\end{definition}
Il est facile de voir que puisque $M$ est compacte et que l'on impose $%
\int_{M}f.u_{t}^{2^{*}}dv_{\mathbf{g}}=1$, il existe forc\'{e}ment un point
de concentration. Pour donner une id\'{e}e plus pr\'{e}cise de l'allure du
ph\'{e}nom\`{e}ne, disons tout de suite que nous montrerons sous de bonnes
hypoth\`{e}ses qu'il n'existe qu'un seul point de concentration, que c'est
un point o\`{u} $f$ est maximum sur $M$, qu'il existe une suite de points $%
(x_{t})$ convergeant vers $x_{0}$ telle que 
\[
u_{t}(x_{t})=\stackunder{M}{\max }u_{t}\rightarrow +\infty , 
\]
et 
\[
u_{t}\rightarrow 0\,\,\,dans\,\,\,C_{loc}^{0}(M-\{x_{0}\}). 
\]
En fait, l'id\'{e}e est que l'on peut faire ``comme si'' les $u_{t}$
\'{e}taient \`{a} support compact dans un voisinage de $x_{0}$ lorsque $t$
est proche de $t_{0}$.

\section{Le changement d'\'{e}chelle}

Gr\^{a}ce au point de concentration, on ram\`{e}ne l'\'{e}tude globale des
solutions $u_{t}$ \`{a} ce qui se passe autour de $x_{0}$. Un bon moyen
d'obtenir des informations (ou des contradictions !) est de faire un
changement d'\'{e}chelle, ou ``blow-up'' en anglais, autour de $x_{0}$. On
appellera changement d'\'{e}chelle de centre $x_{t}$, de coefficient de
dilatation $k_{t}\in \Bbb{R}$, la succession de ``cartes'' et de changements
de m\'{e}triques suivant: On choisit $\delta >0$ assez petit, et on
consid\`{e}re:

\[
\begin{array}{ccccc}
B(x_{t},\delta ) & \stackrel{\exp _{x_{t}}^{-1}}{\rightarrow } & \stackunder{%
}{B(0,\delta )\subset \Bbb{R}^{n}} & \stackrel{\psi _{k_{t}}}{\stackunder{}{%
\rightarrow }} & B(0,k_{t}\delta )\subset \Bbb{R}^{n} \\ 
&  & x\,\,\, & \mapsto \,\, & k_{t}x \\ 
\,\mathbf{g}\, & \rightarrow & \,\,\mathbf{g}\,_{t}=\exp _{x_{t}}^{*}\mathbf{%
g}\, & \rightarrow & \,\widetilde{\,\mathbf{g}\,}_{t}=k_{t}^{2}(\psi
_{k_{t}}^{-1})^{*}\mathbf{g}_{t}
\end{array}
\]
o\`{u} $\exp _{x_{t}}^{-1}$ est la carte exponentielle en $x_{t}$,
c'est-\`{a}-dire d\'{e}duite de l'application exponentielle.

Il est important de savoir comment se transforment l'\'{e}quation $(E_{t})$,
les int\'{e}grales de la forme $\int_{B(x_{t},r)}u_{t}^{\alpha }dv_{\mathbf{g%
}}$ et vers quoi tout cela converge quand $t\rightarrow t_{0}$. Notons 
\[
\overline{u}_{t}=u_{t}\circ \exp _{x_{t}}\, 
\]
\[
\overline{f_{t}}=f\circ \exp _{x_{t}} 
\]
et 
\[
\,\overline{h_{t}}=h_{t}\circ \exp _{x_{t}} 
\]
On a alors \'{e}videmment: 
\begin{eqnarray*}
\triangle _{\mathbf{g}_{t}}\overline{u}_{t}+\overline{h}_{t}.\overline{u}%
_{t} &=&\lambda _{t}\overline{f_{t}}.\overline{u}_{t}^{\frac{n+2}{n-2}} \\
\int_{B(0,r)}\overline{u}_{t}^{\alpha }dv_{\,\mathbf{g}\,_{t}}
&=&\int_{B(x_{t},r)}u_{t}^{\alpha }dv_{\mathbf{g}}\text{ pour tout }\alpha
\geq 1
\end{eqnarray*}

Mais le plus important est la suite: on note 
\begin{eqnarray*}
m_{t} &=&\stackunder{M}{Max}\,u_{t}\, \\
\,\widetilde{u}_{t} &=&m_{t}^{-1}\overline{u}_{t}\circ \psi _{k_{t}}^{-1}\,
\\
\,\widetilde{h}_{t} &=&h_{t}\circ \psi _{k_{t}}^{-1}\, \\
\,\widetilde{f}_{t} &=&\overline{f}_{t}\circ \psi _{k_{t}}^{-1}\, \\
\,\widetilde{\,\mathbf{g}\,}_{t} &=&k_{t}^{2}(\exp _{x_{t}}\circ \psi
_{k_{t}}^{-1})^{*}\mathbf{g},
\end{eqnarray*}
ainsi en particulier 
\begin{eqnarray*}
\,\,\widetilde{u}_{t}(x) &=&m_{t}^{-1}\overline{u}_{t}(\frac{x}{k_{t}})\, \\
\,\widetilde{\,\mathbf{g}\,}_{t}(x) &=&\exp _{x_{t}}^{*}\mathbf{g}(\frac{x}{%
k_{t}}).
\end{eqnarray*}
Alors: 
\begin{eqnarray}
(\widetilde{E}_{t})\, &:&\,\triangle _{\widetilde{\mathbf{g}}_{t}}\widetilde{%
u}_{t}+\frac{1}{k_{t}^{2}}\widetilde{h}_{t}.\widetilde{u}_{t}=\frac{m_{t}^{%
\frac{4}{n-2}}}{k_{t}^{2}}\lambda _{t}\widetilde{f}_{t}.\widetilde{u}_{t}^{%
\frac{n+2}{n-2}}   \\
et\, &:&\,\int_{B(0,k_{t}r)}\widetilde{u}_{t}^{\alpha }dv_{\widetilde{\,%
\mathbf{g}\,}_{t}}=\frac{k_{t}^{n}}{m_{t}^{\alpha }}\int_{B(x_{t},r)}u_{t}^{%
\alpha }dv_{\mathbf{g}}  \nonumber
\end{eqnarray}
Nous utiliserons tr\`{e}s souvent les param\`{e}tres suivants : on consid%
\`{e}re une suite de points ($x_{t})$ tels que 
\[
m_{t}=\stackunder{M}{Max}\,u_{t}=u_{t}(x_{t}):=\mu _{t}^{-\frac{n-2}{2}}
\]
et 
\[
\,k_{t}=\mu _{t}^{-1}.
\]
Nous verrons que $\mu _{t}$ ainsi d\'{e}fini pour simplifier quelques
exposants est un param\`{e}tre fondamental dans l'\'{e}tude du ph\'{e}nom%
\`{e}ne de concentration. En notant $(x^{i})$ les coordonn\'{e}es dans $\Bbb{%
R}^{n}$ (nous en aurons besoin dans certains d\'{e}veloppements limit\'{e}%
s), on a : 
\begin{eqnarray}
(\widetilde{E}_{t})\, &:&\,\triangle _{\widetilde{\mathbf{g}}_{t}}\widetilde{%
u}_{t}+\mu _{t}^{2}.\widetilde{h}_{t}.\widetilde{u}_{t}=\lambda _{t}%
\widetilde{f}_{t}.\widetilde{u}_{t}^{\frac{n+2}{n-2}}   \\
et\, &:&\,\int_{B(0,\mu _{t}^{-1}r)}x^{i_{1}}...x^{i_{p}}.\widetilde{u}%
_{t}^{\alpha }dv_{\widetilde{\,\mathbf{g}\,}_{t}}=\mu _{t}^{-p-n+\alpha 
\frac{n-2}{2}}\int_{B(0,r)}x^{i_{1}}...x^{i_{p}}\overline{u}_{t}^{\alpha
}dv_{\,\mathbf{g}\,_{t}}  \nonumber
\end{eqnarray}
Un r\'{e}sultat tr\`{e}s important (voir par exemple le livre de Th. Aubin
[1]) est que, quand $\mu _{t}\rightarrow 0\,,\,$donc$\,k_{t}\rightarrow
+\infty $, les composantes de $\,\widetilde{\mathbf{g}}_{t}$ convergent dans 
$C_{loc}^{2}$ vers celles de la m\'{e}trique euclidienne, et $(\widetilde{E}%
_{t})\,$``converge'' vers l'\'{e}quation: 
\[
\triangle _{e}\widetilde{u}=\lambda f(x_{0}).\widetilde{u}^{\frac{n+2}{n-2}}
\]
au sens o\`{u} 
\[
\widetilde{u}_{t}\rightarrow \widetilde{u}\,\,\,\,dans\,\,\,\,C_{loc}^{2}(%
\Bbb{R}^{n}).
\]
Il est connu $\left[ 6\right] $ qu'alors 
\[
\widetilde{u}=(1+\frac{\lambda f(x_{0})}{n(n-2)}\left| x\right| ^{2})^{-%
\frac{n-2}{2}}.
\]

\section{Le processus d'it\'{e}ration}

Lorsqu'on multiplie l'\'{e}quation $(E_{t})$ par $u_{t}$ et qu'on
int\`{e}gre sur $M$ on obtient 
\[
\int_{M}\left| \nabla u_{t}\right| _{\mathbf{g}}^{2}dv_{\mathbf{g}%
}+\int_{M}h_{t}.u{{}_{t}^{2}}dv_{\mathbf{g}}=\lambda _{t} 
\]
si l'on a pos\'{e} $\int_{M}f.u_{t}^{2^{*}}dv_{\mathbf{g}}=1$. On voit alors
que, moyennant quelques hypoth\`{e}ses sur $h,$ on va obtenir des
informations sur la norme $H_{1}^{2}$ des $u_{t}$. Poussant plus loin cette
id\'{e}e, en multipliant l'\'{e}quation par $u_{t}^{k}$ et en int\'{e}grant,
on peut esp\'{e}rer obtenir des informations sur les normes $L^{q}$ avec $q$
de plus en plus grand. C'est l'id\'{e}e du processus d'it\'{e}ration de
M\"{o}ser. En fait, pour localiser l'\'{e}tude autour du point de
concentration, qui comme nous le verrons est obligatoirement un point o\`{u} 
$f$ est maximum, nous allons multiplier l'\'{e}quation par $\eta {{}^{2}}%
u_{t}^{k}$ o\`{u} $\eta $ est une fonction cut-off \'{e}gale \`{a} 1
(resp.0) sur une boule $B(x_{0},r)$ o\`{u} $f\geq 0$, \'{e}gale \`{a} 0
(resp. 1) sur $M\backslash B(x_{0},2r)$, et o\`{u} $k\geq 1$, puis
int\'{e}grer, ce qui permettra de faire des int\'{e}grations par parties
``localement'', qui resteront donc valables par exemple apr\`{e}s changement
d'\'{e}chelle. On a ainsi apr\`{e}s quelques int\'{e}grations par parties,
quelques calculs et en utilisant l'\'{e}quation $(E_{t})$ :

\begin{equation}
\frac{4k}{(k+1)^{2}}\int_{M}\left| \nabla (\eta u_{t}^{\frac{k+1}{2}%
})\right| ^{2}=\lambda _{t}\int_{M}f\eta ^{2}u_{t}^{\frac{n+2}{n-2}%
}u_{t}^{k}+\int_{M}(\frac{2}{k+1}\left| \nabla \eta \right| ^{2}+\frac{2(k-1)%
}{(k+1)^{2}}\eta \triangle \eta -\eta ^{2}h_{t})u_{t}^{k+1}
\end{equation}

o\`{u} les int\'{e}grales sont prises par rapport \`{a} la mesure $dv_{%
\mathbf{g}}$ ce que nous sous-entendrons lorsqu'il n'y a pas
d'ambiguit\'{e}s. Ensuite en utilisant l'in\'{e}galit\'{e} de H\"{o}lder, si 
$f\geq 0$ sur $Supp\,\eta $ on obtient: 
\[
\lambda _{t}\int_{M}f\eta ^{2}u_{t}^{\frac{n+2}{n-2}}u_{t}^{k}\leq \lambda
_{t}(\stackunder{Supp\,\eta }{Sup}f)^{\frac{n-2}{n}}.(\int_{Supp\,\eta
}fu_{t}^{\frac{2n}{n-2}})^{\frac{2}{n}}.(\int_{M}(\eta u_{t}^{\frac{k+1}{2}%
})^{\frac{2n}{n-2}})^{\frac{n-2}{n}} 
\]
Puis avec l'in\'{e}galit\'{e} de Sobolev : 
\[
(\int_{M}(\eta u_{t}^{\frac{k+1}{2}})^{\frac{2n}{n-2}})^{\frac{n-2}{n}}\leq
K(n,2){{}^{2}}\int_{M}\left| \nabla (\eta u_{t}^{\frac{k+1}{2}})\right|
^{2}+B\int_{M}\eta u_{t}^{k+1} 
\]
avec $B>0$. D'o\`{u}:

\smallskip

\begin{eqnarray*}
\frac{4k}{(k+1)^{2}}(\int_{M}(\eta u_{t}^{\frac{k+1}{2}})^{\frac{2n}{n-2}})^{%
\frac{n-2}{n}} \leq &\lambda _{t}K(n,2){{}^{2}}(\stackunder{Supp\,\eta }{Sup%
}f)^{\frac{n-2}{n}}.(\int_{Supp\,\eta }fu_{t}^{\frac{2n}{n-2}})^{\frac{2}{n}%
}.(\int_{M}(\eta u_{t}^{\frac{k+1}{2}})^{\frac{2n}{n-2}})^{\frac{n-2}{n}} \\
&+\int_{M}(\frac{4k}{(k+1){{}^{2}}}B\eta +\frac{2}{k+1}\left| \nabla \eta
\right| ^{2}+\frac{2(k-1)}{(k+1)^{2}}\eta \triangle \eta -\eta
^{2}h_{t})u_{t}^{k+1}
\end{eqnarray*}
Alors:

\begin{equation}
Q(t,k,\eta ).(\int_{M}(\eta u_{t}^{\frac{k+1}{2}})^{\frac{2n}{n-2}})^{\frac{%
n-2}{n}}\leq (\frac{4k}{(k+1){{}^{2}}}B+C_{0}+C_{\eta })\int_{Supp\,\eta
}u_{t}^{k+1}\,\,\,\,\,\,\,\,\,\,\,\,\,\,\,\,\,\,\,\,\,\,\,\,\,\,\,\,\,\,\,\,%
\,\,\,\,\,\,\,\,\,\,\,\,\,\,\,\,\,\,\,\,\,\,\,  
\end{equation}
o\`{u}

\[
Q(t,k,\eta )=\frac{4k}{(k+1){{}^{2}}}-\lambda _{t}K(n,2){{}^{2}}(\stackunder{%
Supp\,\eta }{Sup}f)^{\frac{n-2}{n}}.(\int_{Supp\,\eta }f.u_{t}^{2^{*}})^{%
\frac{2}{n}} 
\]
et o\`{u} l'on rappelle que $2^{*}=\frac{2n}{n-2}$ et o\`{u} $%
C_{0}\,et\,C_{\eta }$ sont ind\'{e}pendants de $k$ et $t$ et $\,$tels que $%
\forall k\geq 1,\forall t:$

\[
\,\,\left\| \frac{2}{k+1}\left| \nabla \eta \right| ^{2}+\frac{2(k-1)}{%
(k+1)^{2}}\eta \triangle \eta \right\| _{L^{\infty }(M)}\leq \,C_{\eta
}\,\,et\,\,\left\| h_{t}\right\| _{L^{\infty }(M)}\leq C_{0}\,. 
\]

Si $f$ change de signe sur $Supp\,\eta $, on reprend l'in\'{e}galit\'{e} de
H\"{o}lder: 
\[
\lambda _{t}\int_{M}f\eta ^{2}u_{t}^{\frac{n+2}{n-2}}u_{t}^{k}\leq \lambda
_{t}(\stackunder{Supp\,\eta }{Sup}\left| f\right| ).(\int_{Supp\,\eta
}u_{t}^{\frac{2n}{n-2}})^{\frac{2}{n}}.(\int_{M}(\eta u_{t}^{\frac{k+1}{2}%
})^{\frac{2n}{n-2}})^{\frac{n-2}{n}} 
\]
et on obtient (2.4) avec: 
\begin{equation}
Q(t,k,\eta )=\frac{4k}{(k+1){{}^{2}}}-\lambda _{t}K(n,2){{}^{2}}(\stackunder{%
Supp\,\eta }{Sup}\left| f\right| ).(\int_{Supp\,\eta }u_{t}^{2^{*}})^{\frac{2%
}{n}}  
\end{equation}
Si n\'{e}cessaire, on peut aussi remplacer $\stackunder{Supp\,\eta }{Sup}%
\left| f\right| $ par $\stackunder{M}{Sup}f$.

Le but est de montrer que $(\eta u_{t})$ est bornée dans $L^{\frac{k+1}{2}2^{*}}$ et 
donc qu'on peut en extraire une sous-suite qui converge fortement dans $L^{2^{*}}$.
\section{Remarque}

Ces trois ``outils'' fonctionnent exactement de la m\^{e}me mani\`{e}re pour
des familles d'\'{e}quations un peu plus g\'{e}n\'{e}rales que l'on peut
associer \`{a} $(E_{h,f,\mathbf{g}}):\,\triangle _{\mathbf{g}}u+h.u=\mu
_{h}.f.u^{\frac{n+2}{n-2}}$. Ainsi au lieu d'\'{e}quations $%
E_{t}:\,\,\triangle _{\mathbf{g}}u_{t}+h_{t}.u_{t}=\lambda _{t}.f.u_{t}^{%
\frac{n+2}{n-2}}$ o\`{u} seules les fonctions $h_{t}$ (et $\lambda _{t})\,$%
varient, on peut \^{e}tre amen\'{e} \`{a} associer \`{a} $(E_{h,f,\mathbf{g}%
})$ une famille 
\[
E_{t}:\,\,\triangle _{\mathbf{g}}u_{t}+h_{t}.u_{t}=\lambda
_{t}.f_{t}.u_{t}^{q_{t}-1} 
\]
o\`{u} $q_{t}\rightarrow 2^{*}$ et $f_{t}\rightarrow f$ dans un certain
espace $L^{p}$, avec toujours $h_{t}\rightarrow h$\thinspace dans$%
\,C^{0,\alpha }(M)$ et $\lambda _{t}\rightarrow \lambda $, et o\`{u} l'on
demande $\int_{M}f_{t}u_{t}^{q_{t}}dv_{\mathbf{g}}=1$.

$x_{0}\in M$\textit{\ }sera alors un point de concentration de $(u_{t})$\ si
on a pour tout $\delta >0$: 
\[
\stackunder{t\rightarrow t_{0}}{\lim \sup }\int_{B(x_{0},\delta
)}u_{t}^{q_{t}}>0 
\]
Le principe du changement d'\'{e}chelle est alors analogue, de m\^{e}me que
le principe d'it\'{e}ration. Les formules (2.4) et (2.5) deviennent
ainsi 
\[
Q(t,k,\eta ).(\int_{M}(\eta u_{t}^{\frac{k+1}{2}})^{q_{t}})^{\frac{2}{q_{t}}%
}\leqslant C\int_{Supp\,\eta
}u_{t}^{k+1}\,\,\,\,\,\,\,\,\,\,\,\,\,\,\,\,\,\,\,\,\,\,\,\,\,\,\,\,\,\,\,\,%
\,\,\,\,\,\,\,\,\,\,\,\,\,\,\,\,\,\,\,\,\,\,\, 
\]

\[
\text{o\`{u} }Q(t,k,\eta )=\frac{4k}{(k+1){{}^{2}}}-\lambda
_{t}(Vol_{g}(M))^{\frac{q_{t}}{2^{*}}-1}K{{}(n,2)^{2}}(\stackunder{%
Supp\,\eta }{Sup}\left| f\right| ).(\int_{Supp\,\eta }u_{t}^{q_{t}})^{\frac{%
q_{t}-2}{q_{t}}} 
\]

\section{Principe de d\'{e}monstration du th\'{e}or\`{e}me 1:}

Nous voulons prouver l'existence d'une solution $u>0$ \`{a} l'\'{e}quation 
\[
(E_{h,f,\mathbf{g}}):\,\triangle _{\mathbf{g}}u+h.u=\lambda .f.u^{\frac{n+2}{%
n-2}} 
\]
Comme nous l'avons dit au d\'{e}but de ce chapitre, nous consid\'{e}rons la
famille d'\'{e}quations et de solutions minimisantes associ\'{e}es: 
\[
E_{t}:\,\,\triangle _{\mathbf{g}}u_{t}+h_{t}.u_{t}=\lambda _{t}.f.u_{t}^{%
\frac{n+2}{n-2}} 
\]
avec 
\[
h_{t}\stackrel{\lvertneqq }{\rightarrow }h\,\,\,\,dans\,\,\,\,C^{0,\alpha
}(M) 
\]
les $h_{t}$ \'{e}tant sous-critiques par hypoth\`{e}se.

L'id\'{e}e de base est la suivante: il s'agit d'introduire dans
l'in\'{e}galit\'{e} de Sobolev l'\'{e}quation $E_{t}$ v\'{e}rifi\'{e}e par
la fonction $u_{t}$ pour obtenir une contradiction si les $u_{t}$ convergent
vers la solution nulle. Plus pr\'{e}cis\'{e}ment :

Supposons que $\mathbf{g}$ est plate au voisinage de $x_{0}$ et pour
simplifier que $f\equiv 1$. Alors $S_{\mathbf{g}}=0$ au voisinage de $x_{0}$
et notre hypoth\`{e}se est : 
\[
h(x_{0})>0=\frac{n-2}{4(n-1)}S_{\mathbf{g}}(x_{0}) 
\]
Donc \`{a} partir d'un certain rang: $h_{t}(x_{0})>0$. Or d'une part, comme $%
u_{t}$ est minimisante, 
\[
\lambda _{h_{t},f,\mathbf{g}=}J_{h_{t}}(u_{t}):=\lambda _{t}<K(n,2)^{-2} 
\]
et donc: 
\begin{equation}
\int_{M}\left| \nabla u_{t}\right| ^{2}dv_{\mathbf{g}}+\int_{M}h_{t}.u{%
{}_{t}^{2}}dv_{\mathbf{g}}=\lambda _{t}(\int_{M}u_{t}^{2^{*}}dv_{\mathbf{g}%
})^{\frac{2}{2^{*}}}=\lambda _{t}<K(n,2)^{-2}  
\end{equation}
car $\int_{M}u_{t}^{2^{*}}dv_{\mathbf{g}}=1$; et d'autre part l'identit\'{e}
de Sobolev Euclidienne donne 
\begin{equation}
K(n,2)^{-2}(\int_{\Bbb{R}^{n}}v^{2^{*}})^{\frac{2}{2^{*}}}\leq \int_{\Bbb{R}%
^{n}}\left| \nabla v\right| ^{2}\,.  
\end{equation}
Or si $u_{t}\rightarrow 0$, nous montrerons qu'il y a un ph\'{e}nom\`{e}ne
de concentration, et, comme nous l'avons dit au premier paragraphe, cela
permet de ``faire comme si'' les $u_{t}$ \'{e}taient \`{a} support compact
dans un petit voisinage $B$ de $x_{0}$ o\`{u} $h_{t}>0.$ On aurait donc
d'une part d'apr\`{e}s (2.6) 
\[
\int_{B}\left| \nabla u_{t}\right| ^{2}<K(n,2)^{-2}(\int_{B}u_{t}^{2^{*}})^{%
\frac{2}{2^{*}}}=K(n,2)^{-2} 
\]
puisque $h_{t}>0=S_{\mathbf{g}}$ sur $B$; et d'autre part d'apr\'{e}s (2.7) 
\[
\int_{B}\left| \nabla u_{t}\right| ^{2}\geq
K(n,2)^{-2}(\int_{B}u_{t}^{2^{*}})^{\frac{2}{2^{*}}}=K(n,2)^{-2} 
\]
d'o\`{u} une contradiction.

Pour passer \`{a} l'application rigoureuse de cette id\'{e}e, il faudra
multiplier les $u_{t}$ par des fonctions cut-off et faire des d\'{e}%
veloppements limit\'{e}s de la m\'{e}trique et de $f$ au voisinage de $x_{0}$
, et utiliser les r\'{e}sultats sur les ph\'{e}nom\`{e}nes de concentration
que nous montrerons au d\'{e}but du chapitre 3. Plus pr\'{e}cis\'{e}ment,
avec la fonction $f$ au second membre, on a 
\[
\int_{M}\left| \nabla u_{t}\right| _{\mathbf{g}}^{2}dv_{\mathbf{g}%
}+\int_{M}h_{t}.u{{}_{t}^{2}}dv_{\mathbf{g}}=\lambda
_{t}(\int_{M}fu_{t}^{2^{*}}dv_{\mathbf{g}})^{\frac{2}{2^{*}}}
\]
et c'est dans cette expression que nous aurons \`{a} faire un d\'{e}%
veloppement limit\'{e} de $f$ pour en faire appara\^{i}tre le laplacien en $%
x_{0}$, la contradiction s'obtenant par opposition \`{a} la condition 
\[
\frac{4(n-1)}{n-2}h(P)>S_{\mathbf{g}}(P)-\frac{n-4}{2}\frac{\bigtriangleup _{%
\mathbf{g}}f(P)}{f(P)}
\]
au point de maximum de $f$, puisque $x_{0}$, comme nous le verrons, est un
point de maximum. C'est l\`{a} qu'intervient une nouvelle difficult\'{e} par
rapport au cas $f=cste$ $.$ En effet, comme nous l'avons laiss\'{e} entendre
au paragraphe sur le changement d'\'{e}chelle, nous serons amen\'{e}s \`{a}
consid\'{e}rer une suite de points ($x_{t})$ tels que 
\[
m_{t}=\stackunder{M}{Max}\,u_{t}=u_{t}(x_{t}).
\]
Le but de l'\'{e}tude de la concentration est d'\'{e}tudier la ``forme'' des
fonctions $u_{t}$ quand $t\rightarrow t_{0}$ autour de $x_{t}$. Comme nous
l'avons dit, nous montrerons que 
\[
u_{t}\rightarrow 0\,\,\,dans\,\,\,C_{loc}^{0}(M-\{x_{0}\}).
\]
En faisant un changement d'\'{e}chelle en $x_{t}$ nous obtiendrons des
informations tr\`{e}s pr\'{e}cises sur l'allure des $u_{t}$ autour de $x_{t}$%
. La difficult\'{e} sera alors de ``relier'' ces informations en $x_{t}$
avec celles que nous avons en $x_{0}$ sur $f$ puisque c'est en $x_{0}$ que $f
$ est maximum. Tr\`{e}s pr\'{e}cis\'{e}ment, il nous faudra obtenir une
information sur la ``vitesse'' de convergence de la suite ($x_{t})$ des
points de maximum des $u_{t}$ vers $x_{0}$, point de maximum de $f$ et point
de concentration. Cette vitesse sera mesur\'{e}e par rapport \`{a} la
croissance du maximum des $u_{t}$ \`{a} savoir 
\[
m_{t}=\stackunder{M}{Max}\,u_{t}=u_{t}(x_{t}):=\mu _{t}^{-\frac{n-2}{2}}
\]
et nous voudrons obtenir la relation 
\[
\frac{d_{\mathbf{g}}(x_{t},x_{0})}{\mu _{t}}\leqslant C
\]
o\`{u} $C$ est une constante. Cette relation a \'{e}t\'{e} nomm\'{e}e
``seconde in\'{e}galit\'{e} fondamentale'' par Zo\'{e} Faget $\left[
14\right] $ tant elle est utile dans l'\'{e}tude des EDP \`{a} l'aide des
points de concentration; elle a \'{e}t\'{e} \'{e}tudi\'{e}e par plusieurs
auteurs, en particulier ceux cit\'{e}s au d\'{e}but du chapitre. Elle est en
g\'{e}n\'{e}ral difficile \`{a} obtenir. La d\'{e}finition de $\mu _{t}^{-%
\frac{n-2}{2}}=\stackunder{M}{Max}\,u_{t}=u_{t}(x_{t})$ est faite pour
simplifier l'exposant dans la relation ci-dessus; $\mu _{t}$ sera vraiment
un param\`{e}tre fondamental dans l'\'{e}tude des ph\'{e}nom\`{e}nes de
concentration.

\chapter{Existence de fonctions extr\'{e}males, seconde in\'{e}galit\'{e}
fondamentale.}
\pagestyle{myheadings}\markboth{\textbf{Existence de fonctions extrémales. Inégalité 
fondamentale.}}{\textbf{Existence de fonctions extrémales. Inégalité fondamentale.}}
Ce chapitre se divise en cinq parties. Tout d'abord nous rappelons
l'objectif principal du chapitre, la d\'{e}monstration du th\'{e}or\`{e}me
1, et nous exposons la mise en place de la d\'{e}monstration. La
deuxi\`{e}me partie, bien que s'incluant dans cette d\'{e}monstration,
s'applique en fait au cas d'une suite quelconque de solutions
d'\'{e}quations 
\[
E_{t}:\,\,\triangle _{\mathbf{g}}u_{t}+h_{t}.u_{t}=\lambda _{t}.f.u_{t}^{%
\frac{n+2}{n-2}} 
\]
d\'{e}veloppant un ph\'{e}nom\`{e}ne de concentration, et les r\'{e}sultats
et m\'{e}thodes que l'on y trouve seront r\'{e}utilis\'{e}s dans les autres
chapitres. La troisi\`{e}me partie constitue le coeur de la
d\'{e}monstration, c'est l\`{a} qu'est d\'{e}velopp\'{e} le principe
d'obtention de la ``seconde in\'{e}galit\'{e} fondamentale'' \'{e}voqu\'{e}e
au chapitre pr\'{e}c\'{e}dent. La quatri\`{e}me partie reprend
pr\'{e}cis\'{e}ment cette in\'{e}galit\'{e}, tr\`{e}s importante pour
l'\'{e}tude des ph\'{e}nom\`{e}nes de concentration, dans le cadre
g\'{e}n\'{e}ral d'une suite quelconque de solutions d'\'{e}quations 
\[
E_{t}:\,\,\triangle _{\mathbf{g}}u_{t}+h_{t}.u_{t}=\lambda _{t}.f.u_{t}^{%
\frac{n+2}{n-2}}. 
\]
Enfin, la cinqui\`{e}me partie illustre l'utilit\'{e} de cette
in\'{e}galit\'{e} pour obtenir une autre d\'{e}monstration du
th\'{e}or\`{e}me 1.

Rappelons le th\'{e}or\`{e}me que nous voulons montrer:

\textbf{Donn\'{e}es : }On consid\`{e}re une vari\'{e}t\'{e} riemannienne
compacte $(M,\mathbf{g})$ de dimension $n\geq 4.$ Soit $f:M\rightarrow \Bbb{R%
}$ une fonction $C^{\infty }$ \textit{fix\'{e}e} telle que $\stackunder{M}{%
Sup}$ $f>0$. Soit aussi $h\in C{{}^{\infty }}(M)$ avec l'hypoth\`{e}se
suppl\'{e}mentaire que l'op\'{e}rateur $\bigtriangleup _{\mathbf{g}}+h$ est
coercif si $f$ change de signe sur $M$.

On consid\`{e}re l'\'{e}quation 
\[
(E_{h}^{\prime })=(E_{h,f}^{\prime })=(E_{h,f,\mathbf{g}}^{\prime
}):\,\triangle _{\mathbf{g}}u+h.u=f.u^{\frac{n+2}{n-2}}. 
\]

\textbf{Hypoth\`{e}ses (H)}\textit{: On suppose que le Hessien de la
fonction }$f:M\rightarrow \Bbb{R}$, \textit{telle que }$\stackunder{M}{Sup}%
f>0$, \textit{est non d\'{e}g\'{e}n\'{e}r\'{e} en chaque point de maximum de 
}$f.$ \textit{En outre, les fonctions }$h$ \textit{consid\'{e}r\'{e}es sont
telles que }$\bigtriangleup _{\mathbf{g}}+h$ \textit{est coercif et l'on
suppose }$\dim M\geq 4$.\textit{\ On parlera des hypoth\`{e}ses (}\textbf{H}$%
_{f}$\textit{) pour d\'{e}signer celles concernant la fonction }$f$\textit{.}

Notre but est de prouver le r\'{e}sultat suivant:

\textbf{Th\'{e}or\`{e}me 1:}

\textit{\ Avec l'hypoth\`{e}se }\textbf{(H}),\textit{\ soit }$h$\textit{\
une fonction faiblement critique pour }$f$\textit{\ (et }$\mathbf{g}$\textit{%
)}$.$\textit{\ Si en tout point }$P$\textit{\ o\`{u} }$f$\textit{\ est
maximum sur }$M$\textit{\ on a: }$\frac{4(n-1)}{n-2}h(P)>S_{\mathbf{g}}(P)-%
\frac{n-4}{2}\frac{\bigtriangleup _{\mathbf{g}}f(P)}{f(P)}$\textit{\ et s'il
existe une famille de fonctions }$(h_{t}),\,h_{t}\lvertneqq h,\,h_{t}$%
\textit{\ sous-critique pour tout }$t$\textit{\ dans un voisinage de }$%
t_{0}\in \Bbb{R}$\textit{, telle que }$h_{t}\stackunder{t\rightarrow t_{0}}{%
\rightarrow }h\,\,dans\,$\textit{\ }$C^{0,\alpha }$\textit{\ alors il existe
une solution minimisante pour }$h$\textit{\ et donc }$h$\textit{\ est
critique.}

\section{Mise en place}

Soit donc $h$ une fonction faiblement critique pour $f$ et $\mathbf{g}$
telle qu'en tout point $P$ o\`{u} $f$ est maximum sur $M$ on ait : 
\[
\frac{4(n-1)}{n-2}h(P)>S_{\mathbf{g}}(P)-\frac{(n-4)}{2}\frac{\bigtriangleup
_{\mathbf{g}}f(P)}{f(P)} 
\]
ce que l'on peut \'{e}crire 
\[
\text{ }h(P)>\frac{n-2}{4(n-1)}S_{\mathbf{g}}(P)-\frac{(n-2)(n-4)}{8(n-1)}%
\frac{\bigtriangleup _{\mathbf{g}}f(P)}{f(P)} 
\]
et telle qu'il existe une suite de fonctions $(h_{t}),\,h_{t}\lvertneqq
h,\,h_{t}$ sous-critique pour tout $t$, v\'{e}rifiant $h_{t}\stackunder{%
t\rightarrow t_{0}}{\rightarrow }h\,\,$dans$\,$ $C^{0}$. Pour simplifier on
suppose que $t_{0}=1$ et que $t\rightarrow 1$ ce qui ne change rien. Alors
pour tout $t,\,$%
\[
\lambda _{t}:=\lambda _{h_{t},f,\mathbf{g}}<\frac{1}{K(n,2){{}^{2}}(%
\stackunder{M}{Sup}f)^{\frac{n-2}{n}}} 
\]
et il existe une suite de fonctions $u_{t}>0$ minimisantes pour $h_{t}$ qui
sont solutions de 
\[
E_{t}:\,\,\triangle _{\mathbf{g}}u_{t}+h_{t}.u_{t}=\lambda _{t}.f.u_{t}^{%
\frac{n+2}{n-2}}\text{ avec de plus }\int_{M}fu_{t}^{2^{*}}dv_{\mathbf{g}}=1 
\]
On voit alors que, puisque $\triangle _{\mathbf{g}}+h$ est coercif, la suite $(u_{t})$ est born\'{e}e 
dans $H_{1}^{2}$ (il suffit de multiplier $E_{t}$ par $u_{t}$ et d'int\'{e}grer sur $M$). Il
existe donc une fonction $u\in H_{1}^{2}\,,\,u\geq 0$ telle que, quitte
\`{a} extraire une sous-suite, 
\begin{eqnarray*}
&&u_{t}\stackrel{H_{1}^{2}}{\rightharpoondown }u\,,\, \\
&&u_{t}\stackrel{L^{2}}{\rightarrow }u\,, \\
&&\,u_{t}\stackrel{p.p.}{\rightarrow }u\,,
\end{eqnarray*}
et on peut supposer que 
\[
\lambda _{t}\stackrel{<}{\rightarrow }\lambda \leqslant \frac{1}{K(n,2){%
{}^{2}}(\stackunder{M}{Sup}f)^{\frac{n-2}{n}}}\,\,. 
\]
En particulier 
\[
u_{t}\stackrel{L^{p}}{\rightarrow }u,\,\forall p<2^{*}=\frac{2n}{n-2} 
\]
puisque l'inclusion de $H_{1}^{2}$ dans $L^{p}$ est compacte $\forall
p<2^{*} $. Alors $u$ est solution faible de 
\[
\,\triangle _{\mathbf{g}}u+h.u=\lambda .f.u^{\frac{n+2}{n-2}} 
\]
donc par les th\'{e}ories elliptiques standard $u$ est $C^{\infty }$. Le
principe du maximum nous dit alors que soit $u>0$ soit $u\equiv 0$.

Si $u>0$ alors, comme $h$ est faiblement critique, en multipliant
l'\'{e}quation ci-dessus par $u$ et en int\'{e}grant sur $M$ on voit que (voir appendice A): 
\[
\lambda =\frac{1}{K(n,2){{}^{2}}(\stackunder{M}{Sup}f)^{\frac{n-2}{n}}} 
\]
Donc $u$ est solution de 
\[
\triangle _{\mathbf{g}}u+h.u=\frac{1}{K(n,2){{}^{2}}(\stackunder{M}{Sup}f)^{%
\frac{n-2}{n}}}.f.u^{\frac{n+2}{n-2}}\text{ et }\int_{M}fu^{2^{*}}dv_{%
\mathbf{g}}=1 
\]
et donc $u$ est une solution minimisante pour $h$ et le th\'{e}or\`{e}me est
d\'{e}montr\'{e}.

Si $u\equiv 0$ on est dans le cas o\`{u} il y a ph\'{e}nom\`{e}ne de
concentration tel que d\'{e}crit pr\'{e}c\'{e}demment. Toute l'\'{e}tude de
ce ph\'{e}nom\`{e}ne qui va suivre aura pour but d'aboutir \`{a} une
contradiction. Nous supposons donc \`{a} partir de maintenant que l'on est
dans ce cas : 
\[
u\equiv 0. 
\]

\section{Etude du ph\'{e}nom\`{e}ne de concentration}

De nombreux r\'{e}sultats sur ce ph\'{e}nom\`{e}ne de concentration, d\^{u}s
\`{a} M. Vaugon, E. Hebey, O. Druet et F. Robert entre autres, sont d\'{e}ja
connus, mais ils sont souvent publi\'{e}s dans le cas o\`{u} $f\equiv 1.$
Nous allons les reprendre dans notre cas sachant que la pr\'{e}sence d'une
fonction quelconque $f$ au second membre de 
\[
E_{t}:\,\,\triangle _{\mathbf{g}}u_{t}+h_{t}.u_{t}=\lambda _{t}.f.u_{t}^{%
\frac{n+2}{n-2}} 
\]
ne change quasiment rien \`{a} de nombreuses parties des d\'{e}monstrations,
mais pose de nouveaux probl\`{e}mes dans d'autres. De plus, fait tr\`{e}s
important, elle ``fixe'' en quelque sorte la position du point de
concentration (voir le point a/). Bien que ces r\'{e}sultats fassent partie
de la d\'{e}monstration du th\'{e}or\`{e}me 1, l'\'{e}tude faite dans cette
partie s'applique au cadre plus g\'{e}n\'{e}ral d'une vari\'{e}t\'{e}
riemannienne compacte $(M,\mathbf{g)}$ de dimension $n\geq 3$ sur laquelle
on consid\`{e}re une suite $(u_{t})$ de solutions $C^{2,\alpha }$ de
l'\'{e}quation 
\[
\bigtriangleup _{\mathbf{g}}u_{t}+h_{t}u_{t}=\lambda _{t}fu_{t}^{\frac{n+2}{%
n-2}}\text{ avec }\int_{M}fu_{t}^{\frac{2n}{n-2}}dv_{\mathbf{g}}=1 
\]
o\`{u} $f$ est une fonction dont le maximum est strictement positif. On
suppose de plus que $h_{t}\rightarrow h$ dans $C^{0,\alpha }$ o\`{u} $h$ est
telle que $\bigtriangleup _{\mathbf{g}}+h$ est co\'{e}rcif. La suite $%
(u_{t}) $ est born\'{e}e dans $H_{1}^{2}$, donc $u_{t}\rightharpoondown u$
faiblement dans $H_{1}^{2}$, et on suppose que $u\equiv 0$; ainsi $%
u_{t}\rightarrow 0\,$dans tout $L^{p}$ pour $p<2^{*}$. La suite $u_{t}$
d\'{e}veloppe alors un ph\'{e}nom\`{e}ne de concentration. Nous faisons une
hypoth\`{e}se dite ``d'\'{e}nergie minimale'' (qui est v\'{e}rifi\'{e}e dans
le cadre du th\'{e}or\`{e}me 1): 
\[
\lambda _{t}\leq \frac{1}{K(n,2)^{2}(\stackunder{M}{Sup}f)^{\frac{n-2}{n}}} 
\]

Sous ces hypoth\`{e}ses, cette partie est ainsi ind\'{e}pendante du reste de
la d\'{e}monstration du th\'{e}or\`{e}me consid\'{e}r\'{e} dans ce chapitre;
elle nous servira aussi dans les chapitres suivants. Cette \'{e}tude est
valable pour $\dim M\geq 3$, sauf la partie d/, qui \'{e}tudie la
concentration $L^{2}$, valable pour $\dim M\geq 4$. Dans toute la suite, $%
c,C $ d\'{e}signent des constantes strictement positives ind\'{e}pendantes
de $t$ et $\delta $.

\subsubsection{a/: Proposition: Il existe, \`{a} extraction pr\`{e}s d'une
sous-famille de $(u_{t})$, un unique point de concentration $x_{0}$ et c'est
un point o\`{u} $f$ est maximum sur $M.$ De plus}

\[
\forall \delta >0,\,\,\,\overline{\stackunder{t\rightarrow 1}{\lim }}%
\int_{B(x_{0},\delta )}fu_{t}^{2^{*}}dv_{\mathbf{g}}=1 
\]
C'est une application du principe d'it\'{e}ration. Tout d'abord comme $M$
est compacte, il existe au moins un point de concentration. Sinon on
pourrait recouvrir $M$ par un nombre fini de boules $B(x_{i},\delta )$
telles que $\stackunder{t\rightarrow 1}{\lim }\int_{B(x_{i},\delta
)}u_{t}^{2^{*}}=0$, et on aurait $\stackunder{t\rightarrow 1}{\lim }\int_{M}u_{t}^{2^{*}}=0,$ ce qui
contredirait 
\[
1=\int_{M}fu_{t}^{2^{*}}dv_{\mathbf{g}}\leq Sup\left| f\right|
\int_{M}u_{t}^{2^{*}}dv_{\mathbf{g}} 
\]

Si $f$ change de signe sur $M$, montrons que si $f(x)\leq 0$, alors pour $%
\delta $ assez petit 
\[
\stackunder{t\rightarrow 1}{\lim }\int_{B(x,\delta )}u_{t}^{2^{*}}=0, 
\]
ce qui montre en particulier qu'un tel $x$ n'est pas un point de
concentration. Si $f(x)<0$, on choisit $\delta $ assez petit pour que $f<0$
sur $B(x,\delta )$ et on prend $\eta $ \`{a} support dans $B(x,\delta )$.
Comme les $u_{t}$ sont born\'{e}es ind\'{e}pendamment de $t$ dans $H_{1}^{2}$
donc dans $L^{2^{*}}$, on a, d'apr\`{e}s (2.3) (chapitre
pr\'{e}c\'{e}dent, principe d'it\'{e}ration), pour tout $k$ tel que $1\leq
k\leq 2^{*}-1$: 
\[
\frac{4k}{(k+1)^{2}}\int_{M}\left| \nabla (\eta u_{t}^{\frac{k+1}{2}%
})\right| ^{2}\leq \int_{M}(\frac{2}{k+1}\left| \nabla \eta \right| ^{2}+%
\frac{2(k-1)}{(k+1)^{2}}\eta \triangle \eta -\eta ^{2}h_{t})u_{t}^{k+1}\leq
C_{1} 
\]
o\`{u} $C_{1}$ est ind\'{e}pendant de $t$. Donc pour tout $k$ tel que $1\leq
k\leq 2^{*}-1$ il existe $C_{2}$ ind\'{e}pendant de $t$ tel que: 
\[
\int_{M}\left| \nabla (\eta u_{t}^{\frac{k+1}{2}})\right| ^{2}\leq C_{2} 
\]
Donc $\eta u_{t}^{\frac{k+1}{2}}$ est born\'{e}e ind\'{e}pendamment de $t$
dans $H_{1}^{2}$ et, en utilisant l'in\'{e}galit\'{e} de Sobolev, $\eta
u_{t}^{\frac{k+1}{2}}$ est born\'{e}e dans $L^{2^{*}}$ pour tout $k$ tel que 
$1\leq k\leq 2^{*}-1$.

Si $f(x)=0$, par continuit\'{e} de $f$ et en choisissant $\delta $ assez
petit, on a dans (2.5) pour tout $k$ tel que $1\leq k\leq 2^{*}-1$, $%
Q(t,k,n)\geq Q>0$. Donc avec (2.4), on obtient l\`{a} aussi que $\eta
u_{t}^{\frac{k+1}{2}}$ est born\'{e}e ind\'{e}pendamment de $t$ dans $%
L^{2^{*}}$ pour tout $k$ tel que $1\leq k\leq 2^{*}-1$. Mais par
l'in\'{e}galit\'{e} de H\"{o}lder: 
\[
(\int_{B(x,\delta /2)}u_{t}^{2^{*}})\leq (\int_{B(x,\delta /2)}(u_{t}^{\frac{%
k+1}{2}})^{2^{*}})^{\frac{n+2}{n(k+1)}}(\int_{B(x,\delta /2)}u_{t}^{\frac{%
n(k+1)}{(nk-2)}})^{\frac{nk-2}{n(k+1)}} 
\]
et si l'on suppose que $\overline{\stackunder{t\rightarrow 1}{\lim }}%
\int_{B(x,\delta )}u_{t}^{2^{*}}>0$ on obtient, quand $f(x)\leq 0$, 
\[
\overline{\stackunder{t\rightarrow 1}{\lim }}\int_{B(x,\delta /2)}u_{t}^{%
\frac{n(k+1)}{(nk-2)}}>0 
\]
alors que $\frac{n(k+1)}{(nk-2)}<2^{*}$ si $1\leq k\leq 2^{*}-1$, ce qui
contredit le fait que $u_{t}\stackrel{L^{p}}{\rightarrow }0$, $\forall
p<2^{*}$. Donc pour $\delta $ assez petit: 
\[
\stackunder{t\rightarrow 1}{\lim }\int_{B(x,\delta )}u_{t}^{2^{*}}=0, 
\]
et donc si $f(x)\leq 0,$ $x$ ne peut pas \^{e}tre un point de concentration.

Soit alors $x$ un point de concentration: $f(x)>0$ d'apr\`{e}s ce qui
pr\'{e}c\`{e}de. Posons, pour $\delta >0$ assez petit pour que $f\geq 0$ sur 
$B(x,\delta )$, 
\[
\stackunder{t\rightarrow 1}{\lim \sup }\int_{B(x,\delta
)}fu_{t}^{2^{*}}=a_{\delta } 
\]
Alors $a_{\delta }\leq 1$ car $\int_{M}fu_{t}^{2^{*}}=1.$ Supposons qu'il
existe $\delta >0$ tel que $a_{\delta }<1$. Puisque 
\[
\lambda _{t}\stackunder{\leqslant }{\rightarrow }\lambda \leqslant \frac{1}{%
K(n,2){{}^{2}}(\stackunder{M}{Sup}f)^{\frac{n-2}{n}}} 
\]
on en d\'{e}duit 
\[
\overline{\stackunder{t\rightarrow 1}{\lim }}\,\lambda _{t}K(n,2){{}^{2}}(%
\stackunder{M}{Sup}f)^{\frac{n-2}{n}}a_{\delta }<1. 
\]
Par ailleurs $\frac{4k}{(k+1)^{2}}\stackunder{k\stackunder{>}{\rightarrow }1%
}{\rightarrow }1$. Donc pour $k$ fix\'{e} assez proche de 1 tel que 
\[
\overline{\stackunder{t\rightarrow 1}{\lim }}\,\lambda _{t}K(n,2){{}^{2}}(%
\stackunder{M}{Sup}f)^{\frac{n-2}{n}}a_{\delta }<\frac{4k}{(k+1)^{2}} 
\]
on v\'{e}rifie, en prenant $\eta $ \`{a} support dans $B(x,\delta )$, que
dans (2.4): $Q(t,k,n)\geq Q>0$ pour tout $t$, o\`{u} $Q$ est
ind\'{e}pendant de $t$. On en d\'{e}duit, toujours avec (2.4), que 
\[
(\int_{B(x,\delta /2)}u_{t}^{\frac{k+1}{2}2^{*}})^{\frac{2}{2^{*}}}\leq C%
\text{ quand }t\rightarrow 1\text{ o\`{u} }C\text{ est ind\'{e}pendant de }%
t. 
\]
Alors, par l'in\'{e}galit\'{e} de H\"{o}lder: 
\[
\int_{B(x,\delta /2)}u_{t}^{2^{*}}\leq (\int_{B(x,\delta /2)}u_{t}^{\frac{k+1%
}{2}2^{*}})^{\frac{2}{2^{*}}}(\int_{M}u_{t}^{2^{*}-\frac{(k-1)2^{*}}{2^{*}-2}%
})^{\frac{2^{*}-2}{2^{*}}}\,. 
\]
Or puisque $x$ est un point de concentration 
\[
\overline{\stackunder{t\rightarrow 1}{\lim }}\int_{B(x,\delta
/2)}u_{t}^{2^{*}}>0\,\, 
\]
alors que pour $k$ assez proche de 1: 
\[
\,(\int_{M}u_{t}^{2^{*}-\frac{(k-1)2^{*}}{2^{*}-2}})^{\frac{2^{*}-2}{2^{*}}}%
\stackunder{t\rightarrow 1}{\rightarrow }0 
\]
D'o\`{u} une contradiction et $a_{\delta }=1,\,\forall \delta >0$. Donc $x$
est l'unique point de concentration que nous noterons d\'{e}sormais $x_{0}$.
Ce m\^{e}me raisonnement montre que n\'{e}cessairement 
\[
\lambda =\frac{1}{K(n,2){{}^{2}}(\stackunder{M}{Sup}f)^{\frac{n-2}{n}}}\,\,. 
\]
De la m\^{e}me mani\`{e}re, si $f(x_{0})\neq \stackunder{M}{Sup}f$, $\exists
\delta >0\,$tel que $\stackunder{B(x_{0},\delta )}{Sup}f<\stackunder{M}{Sup}%
f $. Or $\lambda _{t}\leq \frac{1}{K(n,2){{}^{2}}(\stackunder{M}{Sup}f)^{%
\frac{n-2}{n}}}$ donc 
\[
\overline{\stackunder{t\rightarrow 1}{\lim }}\,\lambda _{t}K(n,2){{}^{2}}(%
\stackunder{B(x_{0},\delta )}{Sup}f)^{\frac{n-2}{n}}(\int_{B(x_{0},\delta
)}fu_{t}^{2^{*}})^{\frac{2^{*}-2}{2^{*}}}<1\,\,. 
\]
Alors pour $k$ assez proche de 1, en prenant $\eta $ \`{a} support dans $%
B(x_{{0}},\delta )$, on obtient dans l'in\'{e}galit\'{e} (2.4): $Q(t,k,n)\geq
Q>0 $ pour tout $t$ o\`{u} $Q$ est ind\'{e}pendant de $t$. On aboutit
ensuite de la m\^{e}me mani\`{e}re \`{a} une contradiction. Donc $f(x_{0})=%
\stackunder{M}{Sup}f>0$

\subsubsection{b/: Proposition: $u_{t}\rightarrow 0$ dans $%
C_{loc}^{0}(M-\{x_{0}\})$}

Le processus d'it\'{e}ration prend ici tout son sens.

Premi\`{e}re \'{e}tape: Soit $q>0$ fix\'{e}. On montre qu'il existe pour
tout $\delta >0$, $C=C(\delta ,q)$ ind\'{e}pendant de $t$ tel que pour $t$
assez proche de 1: 
\[
\left\| u_{t}\right\| _{L^{q}(M\backslash B(x_{0},\delta ))}\leq C\left\|
u_{t}\right\| _{L^{2}(M)}\,\,. 
\]
Pour appliquer le principe d'it\'{e}ration, on construit $\eta _{1},...,\eta
_{m}$ $m$ fonctions cut-off telles que $\eta _{j}=0$ sur $B(x_{0},\delta /2)$
et \thinspace $\eta _{j}=1$ sur $M\backslash B(x_{0},\delta )$ et de telle
sorte que 
\[
M\backslash B(x_{0},\delta )\subset ...\subset \{\eta _{j+1}=1\}\subset
Supp\,\eta _{j+1}\subset \{\eta _{j}=1\}\subset ...\subset M\backslash
B(x_{0},\delta /2) 
\]
et o\`{u} $m$ est choisi tel que $2(\frac{2^{*}}{2})^{m}>q$, et on pose $%
q_{1}=2$ et $q_{j}=(\frac{2^{*}}{2})q_{j-1}$. Le processus d'it\'{e}ration
(2.4. et 2.5) donne alors 
\[
Q(t,q_{j}-1,\eta _{j}).(\int_{M}(\eta _{j}u_{t}^{\frac{q_{j}}{2}})^{2^{*}})^{%
\frac{n-2}{n}}\leq (\frac{4(q_{j}-1)}{q_{j}{{}^{2}}}B+C_{0}+C_{\eta
_{j}})\int_{Supp\,\eta _{j}}u_{t}^{q_{j}}\,\,\,. 
\]
Or pour $j\leq m$ on a $\frac{4(q_{j}-1)}{q_{j}{{}^{2}}}\geq c>0$ et
d'apr\'{e}s a/, $\int_{Supp\,\eta _{j}}u_{t}^{2^{*}}\rightarrow 0$, donc
dans (2.5), 
\[
Q(t,q_{j}-1,\eta _{j})\geq c>0,\,\forall j. 
\]
Il existe donc un voisinage $V_{j}$ de 1 et une constante $C_{j}>0$ tels que
pour $t\in V_{j}$: 
\[
(\int_{M}(\eta _{j}u_{t}^{\frac{q_{j}}{2}})^{2^{*}})^{\frac{n-2}{n}}\leq
C_{j}\int_{Supp\,\eta _{j}}u_{t}^{q_{j}}\,\,\,\,. 
\]
Alors d'apr\`{e}s le choix des $\eta _{j}$ on a 
\[
(\int_{\{\eta _{j}=1\}}u_{t}^{q_{j}\frac{2^{*}}{2}})^{\frac{n-2}{n}}\leq
C_{j}\int_{\{\eta _{j-1}=1\}}u_{t}^{q_{j}}\, 
\]
d'o\`{u} 
\[
\left\| u_{t}\right\| _{L^{q}(M\backslash B(x_{0},\delta ))}\leq
C(\prod\limits_{j=1}^{m}C_{j})\left\| u_{t}\right\| _{L^{2}(M)}\,\,\forall
t\in V_{1}\cap ...\cap V_{m}\,. 
\]

Deuxi\`{e}me \'{e}tape: Le th\'{e}or\`{e}me (8.25) de Gilbarg-Trudinger $%
\left[ 16\right] $ nous donne: si $u$ est solution d'une \'{e}quation de la
forme $E:\,\,\triangle _{\mathbf{g}}u+h.u=F$, o\`{u} $\triangle _{\mathbf{g}%
}+h$ est coercif, et si $\omega \subset \subset \omega ^{\prime }$ sont deux
ouverts, pour $r>1,\,q>n/2\,:$%
\[
\stackunder{\omega }{Sup}\,u\leq c\left\| u\right\| _{L^{r}(\omega ^{\prime
})}+c^{\prime }\left\| F\right\| _{L^{q}(\omega ^{\prime })}\,\,.
\]
La d\'{e}monstration de ce th\'{e}or\`{e}me est en fait une application du
principe d'it\'{e}ration, nous l'utilisons ici pour gagner du temps, en
l'appliquant \`{a} $E_{t}:\,\,\triangle _{\mathbf{g}}u_{t}+h_{t}.u_{t}=%
\lambda _{t}.f.u_{t}^{\frac{n+2}{n-2}}$ et \`{a} $\omega \subset \subset
\omega ^{\prime }\subset M\backslash \{x_{0}\}$.

Alors avec la premi\`{e}re \'{e}tape appliqu\'{e}e \`{a} $q\frac{n+2}{n-2}$,
en choisissant 
\[
\omega =M\backslash B(x_{0},\delta ),\,\omega ^{\prime }=M\backslash
B(x_{0},\delta /2),\,r=2,\,q>n/2 
\]
on obtient 
\begin{eqnarray*}
\stackunder{M\backslash B(x_{0},\delta )}{Sup}u_{t} & \leq &c\left\|
u_{t}\right\| _{L^{2}(\omega ^{\prime })}+c^{\prime }\lambda _{t}^{q}\left\|
u_{t}\right\| _{L^{q\frac{n+2}{n-2}}(\omega ^{\prime })}^{\frac{n+2}{n-2}}
\\ 
& \leq &c\left\| u_{t}\right\| _{L^{2}(M)}+c^{\prime \prime }\left\|
u_{t}\right\| _{L^{2}(M)}^{\frac{n+2}{n-2}}
\end{eqnarray*}
Or $\left\| u_{t}\right\| _{L^{2}(M)}\rightarrow 0$, d'o\`{u} la conclusion.

\subsubsection{c/: Estim\'{e}es ponctuelles faibles}

Reprenons les notations du changement d'\'{e}chelle d\'{e}crit au chapitre
2: on consid\`{e}re une suite de points $(x_{t})$ tel que 
\[
m_{t}=\stackunder{M}{Max}\,u_{t}=u_{t}(x_{t}):=\mu _{t}^{-\frac{n-2}{2}}. 
\]
D'apr\`{e}s ce qui pr\'{e}c\`{e}de $x_{t}\rightarrow x_{0}$ et $\mu
_{t}\rightarrow 0$. Rappelons que $\overline{u}_{t},\overline{f}_{t},%
\overline{h}_{t},\,\mathbf{g}\,_{t}$ d\'{e}signent les fonctions et la
m\'{e}trique ``lues'' dans la carte $\exp _{x_{t}}^{-1}$, et $\,\,\widetilde{%
u}_{t}\,,\,\widetilde{h}_{t}\,,\,\widetilde{f}_{t},\widetilde{\,\mathbf{g}\,}%
_{t}$ d\'{e}signent les fonctions et la m\'{e}trique ``lues'' apr\`{e}s
blow-up. \textit{D\'{e}sormais, tous les changements d'\'{e}chelle ont pour
point de d\'{e}part des boules }$B(x_{t},\delta )$\textit{\ sur lesquelles }$%
f\geq 0$\textit{\ ce qui est possible puisque }$f(x_{0})>0$\textit{.}

\textbf{Proposition 2:}

$\forall R>0$ : $\stackunder{t\rightarrow 1}{\lim }\int_{B(x_{t},R\mu
_{t})}fu_{t}^{2^{*}}dv_{\mathbf{g}}=1-\varepsilon _{R}$ o\`{u} $\varepsilon
_{R}\stackunder{R\rightarrow +\infty }{\rightarrow }0.$

\textit{D\'{e}monstration: }C'est\textit{\ }une application directe du
changement d'\'{e}chelle en $x_{t}$ avec $k_{t}=\mu _{t}^{-1}$:
\begin{center}
$\widetilde{u}_{t}\rightarrow \widetilde{u}=(1+\frac{\lambda f(x_{0})}{n(n-2)%
}\left| x\right| ^{2})^{-\frac{n-2}{2}}\,=(1+\frac{f(x_{0})^{\frac{2}{n}}}{%
K(n,2)^{2}n(n-2)}\left| x\right| ^{2})^{-\frac{n-2}{2}}$ dans $\,C_{loc}^{2}(%
\Bbb{R}^{n})$.
\end{center}
Alors:
\begin{center}
$\int_{B(x_{t},R\mu _{t})}fu_{t}^{2^{*}}dv_{\mathbf{g}}=\int_{B(0,R)}%
\widetilde{f}_{t}.\widetilde{u}_{t}^{2^{*}}dv_{\widetilde{\,\mathbf{g}\,}%
_{t}}\stackunder{t\rightarrow 1}{\rightarrow }f(x_{0})(\int_{B(0,R)}%
\widetilde{u}^{2^{*}}dx)\stackunder{R\rightarrow \infty }{\rightarrow }1$
\end{center}

\textbf{Proposition 3:}

$\exists C>0$ tel que $\forall x\in M:\,d_{\mathbf{g}}(x,x_{t})^{\frac{n-2}{2%
}}u_{t}(x)\leq C$.

\textit{D\'{e}monstration: Posons }$w_{t}(x)=\,d_{\mathbf{g}}(x,x_{t})^{%
\frac{n-2}{2}}u_{t}(x)$. On veut montrer qu'il existe $C>0$ tel que $%
\stackunder{M}{Sup}\,w_{t}\leq C$. Par contradiction on suppose que (pour une
sous-suite) $\stackunder{M}{Sup}\,w_{t}\rightarrow +\infty $. Soit $y_{t}$ un
point o\`{u} $w_{t}$ est maximum. $M$ \'{e}tant compacte, $d_{\mathbf{g}%
}(x,x_{t})$ est born\'{e}e, donc $u_{t}(y_{t})\rightarrow \infty $, donc
d'apr\`{e}s le point b/ $y_{t}\rightarrow x_{0}$. Par ailleurs la
d\'{e}finition de $\mu _{t}$ nous donne: 
\[
\frac{d_{\mathbf{g}}(y_{t},x_{t})}{\mu _{t}}\rightarrow +\infty \,\,. 
\]
On fait un changement d'\'{e}chelle en $y_{t}$ de coefficient $%
k_{t}=u_{t}(y_{t})^{\frac{2}{n-2}}$ et avec $m_{t}=u_{t}(y_{t})$. On
obtient: 
\[
\,\triangle _{\widetilde{\mathbf{g}}_{t}}\widetilde{u}_{t}+u_{t}(y_{t})^{-%
\frac{4}{n-2}}.\widetilde{h}_{t}.\widetilde{u}_{t}=\lambda _{t}\widetilde{f}%
_{t}.\widetilde{u}_{t}^{\frac{n+2}{n-2}}\,\,. 
\]
Si $x\in B(0,2):$%
\begin{eqnarray*}
d_{\mathbf{g}}(x_{t},\exp _{y_{t}}(u_{t}(y_{t})^{-\frac{2}{n-2}}x) & \geq& d_{%
\mathbf{g}}(y_{t},x_{t})-2u_{t}(y_{t})^{-\frac{2}{n-2}} \\ 
& \geq& u_{t}(y_{t})^{-\frac{2}{n-2}}(w_{t}(y_{t})^{\frac{2}{n-2}}-2)\sim d_{%
\mathbf{g}}(y_{t},x_{t})
\end{eqnarray*}
car $w_{t}(y_{t})\rightarrow \infty $ et $u_{t}(y_{t})\rightarrow \infty $,
donc pour $t$ proche de 1: 
\[
d_{\mathbf{g}}(x_{t},\exp _{y_{t}}(u_{t}(y_{t})^{-\frac{2}{n-2}}x)\geq \frac{%
1}{2}d_{\mathbf{g}}(y_{t},x_{t})\,\,. 
\]
Par cons\'{e}quent pour tout $R>0$ et pour $t$ proche de 1: 
\[
B(y_{t},2u_{t}(y_{t})^{-\frac{2}{n-2}})\cap B(x_{t},R\mu _{t})=\emptyset 
\]
Alors, d'apr\`{e}s la proposition pr\'{e}c\'{e}dente 
\begin{eqnarray*}
\int_{B(0,2)}\widetilde{f}_{t}.\widetilde{u}_{t}^{2^{*}}dv_{\widetilde{\,%
\mathbf{g}\,}_{t}} & =\int_{B(y_{t},2u_{t}(y_{t})^{-\frac{2}{n-2}%
})}fu_{t}^{2^{*}}dv_{\mathbf{g}} & \leq \int_{M\backslash B(x_{t},R\mu
_{t})}fu_{t}^{2^{*}}dv_{\mathbf{g}} \\ 
&  & \leq \int_{M}fu_{t}^{2^{*}}dv_{\mathbf{g}}-\int_{B(x_{t},R\mu
_{t})}fu_{t}^{2^{*}}dv_{\mathbf{g}} \\ 
&  & \stackunder{t\rightarrow 1,R\rightarrow \infty }{\longrightarrow }0
\end{eqnarray*}
Mais le principe d'it\'{e}ration donne alors pour $1\leq k\leq 2^{*}-1:$%
\[
\int_{B(0,1)}\widetilde{u}_{t}^{\frac{k+1}{2}2^{*}}dv_{\widetilde{\,\mathbf{g%
}\,}_{t}}\rightarrow 0 
\]
et par r\'{e}currence on obtient $\forall p\geq 1:$%
\[
\int_{B(0,1)}\widetilde{u}_{t}^{p}dv_{\widetilde{\,\mathbf{g}\,}%
_{t}}\rightarrow 0 
\]
On en d\'{e}duit alors que $\left\| \widetilde{u}_{t}\right\| _{L^{\infty
}(B(0,1))}\rightarrow 0$ bien que $\widetilde{u}_{t}(0)=1$. D'o\`{u} une
contradiction.

\textbf{Proposition 4:}

$\forall \varepsilon >0$ , $\exists R>0$ tel que $\forall t,\,\forall x\in
M: $%
\[
\,d_{\mathbf{g}}(x,x_{t})\geq R\mu _{t}\,\Rightarrow \,\,d_{\mathbf{g}%
}(x,x_{t})^{\frac{n-2}{2}}u_{t}(x)\leq \varepsilon . 
\]

\textit{D\'{e}monstration:} On utilise le m\^{e}me principe, on suppose
qu'il existe $\varepsilon _{0}>0$ et $y_{t}\in M$ tels que 
\[
\stackunder{t\rightarrow 1}{\lim }\frac{d_{\mathbf{g}}(y_{t},x_{t})}{\mu _{t}%
}=+\infty \text{ et }w_{t}(y_{t})=\,d_{\mathbf{g}}(y_{t},x_{t})^{\frac{n-2}{2%
}}u_{t}(y_{t})\geq \varepsilon _{0} 
\]
On fait un changement d'\'{e}chelle en $y_{t}$ de coefficient $%
k_{t}=u_{t}(y_{t})^{\frac{2}{n-2}}$ et avec $m_{t}=u_{t}(y_{t}).$

Alors avec ces hypoth\`{e}ses si $x\in B(0,\frac{1}{2}\varepsilon _{0}^{%
\frac{2}{n-2}}):$%
\begin{eqnarray*}
d_{\mathbf{g}}(x_{t},\exp _{y_{t}}(u_{t}(y_{t})^{-\frac{2}{n-2}}x) && \geq d_{%
\mathbf{g}}(y_{t},x_{t})-\frac{1}{2}\varepsilon _{0}^{\frac{2}{n-2}%
}u_{t}(y_{t})^{-\frac{2}{n-2}} \\ 
& &\geq \frac{1}{2}d_{\mathbf{g}}(y_{t},x_{t})
\end{eqnarray*}
Par cons\'{e}quent pour tout $R>0$ et pour $t$ proche de 1: 
\[
B(y_{t},\frac{1}{2}\varepsilon _{0}^{\frac{2}{n-2}}u_{t}(y_{t})^{-\frac{2}{%
n-2}})\cap B(x_{t},R\mu _{t})=\emptyset 
\]
donc comme pr\'{e}c\'{e}demment: 
\[
\int_{B(0,\frac{1}{2}\varepsilon _{0}^{\frac{2}{n-2}})}\widetilde{f}_{t}.%
\widetilde{u}_{t}^{2^{*}}dv_{\widetilde{\,\mathbf{g}\,}_{t}}\rightarrow 0 
\]
et on aboutit de la m\^{e}me mani\`{e}re \`{a} une contradiction.

\subsubsection{d/: Concentration $L^{2}$}

\textbf{Proposition 5:}

Si $\dim M\geq 4,$ $\forall \delta >0\,:\,$%
\[
\stackunder{t\rightarrow 1}{\lim }\frac{\int_{B(x_{0},\delta )}u_{t}^{2}dv_{%
\mathbf{g}}}{\int_{M}u_{t}^{2}dv_{\mathbf{g}}}=1 
\]

Cette limite qui d\'{e}crit un quotient de ``normes $L^{2}$'' justifie la
terminologie ``concentration $L^{2}$''.

\textit{D\'{e}monstration:}

On reprend pour commencer les deux ``\'{e}tapes'' du point b/ pour montrer
que pour tout $\delta >0$, il existe $c>0$ tel que: 
\[
\stackunder{M\backslash B(x_{0},\delta )}{Sup}u_{t}\leq c\left\|
u_{t}\right\| _{L^{2}(M)}\,\,. 
\]
Il suffit de reprendre la fin de la deuxi\`{e}me \'{e}tape de la
d\'{e}monstration du point b/: 
\begin{eqnarray*}
\stackunder{M\backslash B(x_{0},\delta )}{Sup}u_{t} && \leq c\left\|
u_{t}\right\| _{L^{2}(\omega ^{\prime })}+c^{\prime }\lambda _{t}^{q}\left\|
u_{t}^{\frac{n+2}{n-2}}\right\| _{L^{q}(\omega ^{\prime })} \\ 
&& \leq c\left\| u_{t}\right\| _{L^{2}(M)}+c^{\prime }\lambda _{t}^{q}%
\stackunder{\omega ^{\prime }}{Sup}(u_{t}^{\frac{n+2}{n-2}-1})\left\|
u_{t}\right\| _{L^{q}(\omega ^{\prime })}\\
&&\leq c^{\prime \prime }\left\|
u_{t}\right\| _{L^{2}(M)}
\end{eqnarray*}
car on sait maintenant que $\stackunder{\omega ^{\prime }}{Sup}(u_{t}^{\frac{%
n+2}{n-2}-1})\rightarrow 0$ et que, d'autre part, la premi\`{e}re \'{e}tape
donne $\left\| u_{t}\right\| _{L^{q}(\omega ^{\prime })}\leq C\left\|
u_{t}\right\| _{L^{2}(M)}$

Troisi\`{e}me \'{e}tape: Avec ce qui pr\'{e}c\`{e}de: 
\begin{eqnarray}
\left\| u_{t}\right\| _{L^{2}(M\backslash B(x_{0},\delta ))}^{2} && \leq 
\stackunder{M\backslash B(x_{0},\delta )}{Sup}u_{t}.\int_{M\backslash
B(x_{0},\delta )}u_{t}\nonumber \\ 
&& \leq c\left\| u_{t}\right\| _{L^{2}(M)}\left\| u_{t}\right\| _{L^{1}(M)}
\end{eqnarray}
On veut maintenant montrer que 
\begin{equation}
\left\| u_{t}\right\| _{L^{1}(M)}\leq c\left\| u_{t}\right\|
_{L^{2^{*}-1}(M)}^{2^{*}-1}\,\,.  
\end{equation}
Si $h>0$, il suffit pour cela d'int\'{e}grer l'\'{e}quation $E_{t}$. Sinon
comme $\lambda _{h,f,\,\mathbf{g}\,}>0$, pour tout $q\in ]2,2^{*}[$, il
existe $\varphi >0$ solution de $\triangle _{\mathbf{g}}\varphi +h\varphi
=\lambda _{h,f,\,\mathbf{g}\,}.f.\varphi ^{q-1}$. On pose 
\[
\mathbf{g}^{\prime }=\varphi ^{\frac{4}{n-2}}\mathbf{g}\text{ et }\overline{%
h_{t}}=\frac{\triangle _{\mathbf{g}}\varphi +h_{t}\varphi }{\varphi ^{\frac{%
n+2}{n-2}}} 
\]
Alors pour $t$ assez proche de 1: 
\[
\overline{h_{t}}=\varphi ^{q-2^{*}}-(h-h_{t})\varphi ^{2-2^{*}}\geq
\varepsilon _{0}>0 
\]
Par ailleurs, par invariance conforme et d'apr\`{e}s $E_{t}$, on a: 
\[
\triangle _{\mathbf{g}^{\prime }}\overline{u_{t}}+\overline{h_{t}}.\overline{%
u_{t}}=\lambda _{t}f.\overline{u_{t}}^{\frac{n+2}{n-2}} 
\]
o\`{u} $\overline{u_{t}}=\varphi ^{-1}.u_{t}$. En int\'{e}grant, on obtient: 
\[
\varepsilon _{0}\int_{M}\overline{u_{t}}dv_{\mathbf{g}^{\prime }}\leq
\lambda _{t}Supf\int_{M}\overline{u_{t}}^{\frac{n+2}{n-2}}dv_{\mathbf{g}%
^{\prime }} 
\]
et il existe donc $C>0$ tel que pour $t$ assez proche de 1 
\[
\left\| u_{t}\right\| _{L^{1}(M)}\leq C\left\| u_{t}\right\|
_{L^{2^{*}-1}(M)}^{2^{*}-1} 
\]
o\`{u} cette fois les normes sont relatives \`{a} $dv_{\mathbf{g}}$.

Quatri\`{e}me \'{e}tape: On conclut essentiellement avec l'in\'{e}galit\'{e}
de H\"{o}lder. Si $n=\dim M\geq 6$, on a en effet : 
\[
\left\| u_{t}\right\| _{L^{2^{*}-1}(M)}^{2^{*}-1}\leq \left\| u_{t}\right\|
_{L^{2}(M)}^{\frac{n+2}{n-2}}Vol_{\mathbf{g}}(M)^{\frac{n-6}{2(n-2)}}\,. 
\]
Avec (3.1) et (3.2), on obtient 
\[
\stackunder{t\rightarrow 1}{\lim }\frac{\left\| u_{t}\right\|
_{L^{2}(M\backslash B(x_{0},\delta ))}^{2}}{\left\| u_{t}\right\|
_{L^{2}(M)}^{2}}=0 
\]
ce qui prouve le r\'{e}sultat. Si $n=5$, l'in\'{e}galit\'{e} de H\"{o}lder
nous donne: 
\[
\left\| u_{t}\right\| _{L^{2^{*}-1}(M)}^{2^{*}-1}\leq \left\| u_{t}\right\|
_{L^{2}(M)}^{\frac{3}{2}}\left\| u_{t}\right\| _{L^{2}(M)}^{\frac{5}{6}} 
\]
et on conclut de la m\^{e}me mani\`{e}re avec (3.1) et (3.2). Si maintenant $n=4$,
on utilise la proposition 4 et le changement d'\'{e}chelle associ\'{e}. On a 
\[
\frac{\left\| u_{t}\right\| _{L^{3}(M)}^{3}}{\left\| u_{t}\right\|
_{L^{2}(M)}}\leq \left\| u_{t}\right\| _{L^{\infty }(M\backslash
B(x_{0},\delta ))}\left\| u_{t}\right\| _{L^{2}(M)}+\frac{\int_{B(0,\delta
\mu _{t}^{-1})}\widetilde{u}_{t}^{3}dv_{\widetilde{\,\mathbf{g}\,}_{t}}}{%
(\int_{B(0,\delta \mu _{t}^{-1})}\widetilde{u}_{t}^{2}dv_{\widetilde{\,%
\mathbf{g}\,}_{t}})^{\frac{1}{2}}}\,\,. 
\]
Alors pour tout $R>0$, d'après l'in\'{e}galit\'{e} de H\"{o}lder et la
proposition 4 on obtient: 
\[
\int_{B(0,\delta \mu _{t}^{-1})}\widetilde{u}_{t}^{3}dv_{\widetilde{\,%
\mathbf{g}\,}_{t}}\leq \int_{B(0,R)}\widetilde{u}_{t}^{3}dv_{\widetilde{\,%
\mathbf{g}\,}_{t}}+\varepsilon _{R}(\int_{B(0,\delta \mu _{t}^{-1})}%
\widetilde{u}_{t}^{2}dv_{\widetilde{\,\mathbf{g}\,}_{t}})^{\frac{1}{2}}\,\,. 
\]
Il s'ensuit pour tout $R,R^{\prime }>0$,

\[
\stackunder{t\rightarrow 1}{\lim \sup }\frac{\left\| u_{t}\right\|
_{L^{3}(M)}^{3}}{\left\| u_{t}\right\| _{L^{2}(M)}}\leq \varepsilon _{R}+%
\frac{\int_{B(0,R)}\widetilde{u}^{3}dx}{(\int_{B(0,R^{\prime })}\widetilde{u}%
^{2}dx)^{\frac{1}{2}}}\,\,. 
\]
Comme $\widetilde{u}\in L^{3}(\Bbb{R}^{4})$ et $\stackunder{R^{\prime
}\rightarrow \infty }{\lim }\int_{B(0,R^{\prime })}\widetilde{u}%
^{2}dx=+\infty $, on a finalement 
\[
\stackunder{t\rightarrow 1}{\lim \sup }\frac{\left\| u_{t}\right\|
_{L^{3}(M)}^{3}}{\left\| u_{t}\right\| _{L^{2}(M)}}=0 
\]
et on conclut encore une fois avec (3.1) et (3.2).

\subsubsection{e/: Estim\'{e}es ponctuelles fortes et concentration $L^{p}$
forte}

\textbf{Proposition 6:}

$\exists C>0$ tel que $\forall x\in M:\,d_{\mathbf{g}}(x,x_{t})^{n-2}\mu
_{t}^{-\frac{n-2}{2}}u_{t}(x)\leq C$.

\textit{D\'{e}monstration:} Elle n\'{e}cessite l'utilisation de la fonction
de Green, ainsi que les estim\'{e}es ponctuelles faibles. L'id\'{e}e est due 
\`{a} O. Druet et \`{a} F. Robert [13]. Nous en profitons pour l'exposer
avec quelques explications suppl\'{e}mentaires (et avec notre fonction $f$,
ce qui ne change pratiquement rien).

Rappelons les propri\'{e}t\'{e}s des fonctions de Green dont nous aurons
besoin (voir appendice B).

Si $\bigtriangleup _{\mathbf{g}}+h$ est un op\'{e}rateur coercif, il existe
une unique fonction (au moins $C^{2}$ dans nos hypoth\`{e}ses) 
\[
G_{h}:M\times M\backslash \{(x,x),x\in M\}\rightarrow \Bbb{R} 
\]
sym\'{e}trique et strictement positive telle que, au sens des distributions,
on a: $\forall x\in M$%
\begin{equation}
\bigtriangleup _{\mathbf{g},y}G_{h}(x,y)+h(y)G_{h}(x,y)=\delta _{x}  
\end{equation}
De plus, il existe $c>0,\,\rho >0$ tels que $\forall (x,y)$ avec $0<d_{%
\mathbf{g}}(x,y)<\rho :$%
\begin{equation}
\frac{c}{d_{\mathbf{g}}(x,y)^{n-2}}\leq G_{h}(x,y)\leq \frac{c^{-1}}{d_{%
\mathbf{g}}(x,y)^{n-2}}  
\end{equation}

\begin{equation}
\frac{\left| \nabla _{y}G_{h}(x,y)\right| }{G_{h}(x,y)}\geq \frac{c}{d_{%
\mathbf{g}}(x,y)}  
\end{equation}

\begin{center}
$c$ et $\rho $ varient continûment avec $h$ 
\begin{equation}
G_{h}(x,y)d_{\mathbf{g}}(x,y)^{n-2}\rightarrow \frac{1}{(n-2)\omega _{n-1}}%
\text{ quand }d_{\mathbf{g}}(x,y)\rightarrow 0  
\end{equation}
\end{center}

Principe de la d\'{e}monstration de la proposition: pour montrer cette estim%
\'{e}e forte il suffit, d'apr\`{e}s (3.4), de montrer que $\mu _{t}^{-\frac{n-2%
}{2}}u_{t}(x)\leq c^{\prime }G_{h}(x,x_{t})$. On remarque alors que gr\^{a}%
ce aux estim\'{e}es faibles, on a d\'{e}ja l'estim\'{e}e forte dans toute
boule $B(x_{t},R\mu _{t})$ o\`{u} $R$ est fix\'{e}. Il suffit donc de
montrer cette estim\'{e}e dans la vari\'{e}t\'{e} \`{a} bord $M\backslash
B(x_{t},R\mu _{t})$ dont le bord est $b(M\backslash B(x_{t},R\mu
_{t}))=bB(x_{t},R\mu _{t})$. On applique pour cela le principe du maximum 
\`{a} l'op\'{e}rateur: 
\[
L_{t}\varphi =\bigtriangleup _{\mathbf{g}}\varphi +h_{t}\varphi -\lambda
_{t}fu_{t}^{2^{*}-2}\varphi 
\]
et \`{a} $x\longmapsto G_{h}(x,x_{t})-c\mu _{t}^{-\frac{n-2}{2}}u_{t}(x)$.
Comme $L_{t}u_{t}=0$ avec $u_{t}>0$, $L_{t}$ v\'{e}rifie le principe du
maximum $\left[ 5\right] $. Si on montre que 
\begin{eqnarray*}
L_{t}G_{h}(x,x_{t}) &\geq &c\mu _{t}^{-\frac{n-2}{2}}L_{t}u_{t}(x)=0\text{
sur }M\backslash B(x_{t},R\mu _{t}) \\
G_{h}(x,x_{t}) &\geq &c\mu _{t}^{-\frac{n-2}{2}}u_{t}(x)\text{ sur }%
bB(x_{t},R\mu _{t})
\end{eqnarray*}
le principe du maximum nous donnera 
\[
G_{h}(x,x_{t})\geq c\mu _{t}^{-\frac{n-2}{2}}u_{t}(x)\text{ sur }M\backslash
B(x_{t},R\mu _{t})\,\,.
\]
Pour des raisons techniques, on montre d'abord que pour tout $\nu >0$ assez
petit il existe une constante $C(\nu )>0$ pour laquelle on a ``l'estim\'{e}e
forte interm\'{e}diaire'' 
\begin{equation}
\forall x\in M:\,d_{\mathbf{g}}(x,x_{t})^{n-2-\nu }\mu _{t}^{-\frac{n-2}{2}%
+\nu }u_{t}(x)\leq C(\nu )  
\end{equation}
puis on montre ensuite que c'est encore vrai pour $\nu =0$.

Pour $\nu $ assez petit, il existe $\varepsilon _{0}>0$ tel que
l'op\'{e}rateur 
\[
\bigtriangleup _{\mathbf{g}}+\frac{h-2\varepsilon _{0}}{1-\nu } 
\]
soit encore coercif; soit $\widetilde{G}$ sa fonction de Green. On va
appliquer le principe d\'{e}crit \`{a} $\widetilde{G}^{1-\nu }$. Alors on
aura (3.7) en prenant $\nu ^{\prime }=(n-2)\nu $. Pour montrer $L_{t}%
\widetilde{G}^{1-\nu }\geq 0$ sur $M\backslash B(x_{t},R\mu _{t})$ on va
montrer que 
\[
\frac{L_{t}\widetilde{G}^{1-\nu }}{\widetilde{G}^{1-\nu }}\geq 0 
\]
sur $M\backslash B(x_{t},R\mu _{t})$. Comme $\widetilde{G}^{1-\nu }>0$ on
aura bien $L_{t}\widetilde{G}^{1-\nu }\geq 0$.

Un calcul donne en utilisant (3.3) et $\delta _{x_{t}}(x)=0$ sur $M\backslash
B(x_{t},R\mu _{t})$ que $\forall x\in M\backslash B(x_{t},R\mu _{t}):$%
\[
\frac{L_{t}\widetilde{G}^{1-\nu }}{\widetilde{G}^{1-\nu }}%
(x,x_{t})=2\varepsilon _{0}+h_{t}(x)-h(x)-\lambda
_{t}f(x)u_{t}(x)^{2^{*}-2}+\nu (1-\nu )\left| \frac{\nabla \widetilde{G}}{%
\widetilde{G}}\right| ^{2}(x,x_{t}) 
\]
Or pour $t$ proche de 1, $h_{t}-h\geq -\varepsilon _{0}$ car $%
h_{t}\rightarrow h$ dans $C^{0}$. Donc 
\begin{equation}
\frac{L_{t}\widetilde{G}^{1-\nu }}{\widetilde{G}^{1-\nu }}(x,x_{t})\geq
\varepsilon _{0}-\lambda _{t}f(x)u_{t}(x)^{2^{*}-2}+\nu (1-\nu )\left| \frac{%
\nabla \widetilde{G}}{\widetilde{G}}\right| ^{2}(x,x_{t})  
\end{equation}
On s\'{e}pare maintenant $M\backslash B(x_{t},R\mu _{t})$ en deux parties
avec une boule $B(x_{t},\rho ),\,$o\`{u} $\rho >0$ est comme dans (3.4) et (3.5).
Pour $t$ assez proche de 1, $\rho >R\mu _{t}$. $R>0$ sera fix\'{e} plus tard.

1/: Comme $u_{t}\rightarrow 0$ dans $C_{loc}^{0}(M\backslash \{x_{0}\})$,
(3.8) donne pour $t$ proche de 1: 
\[
\forall x\in M\backslash B(x_{t},\rho ):\,L_{t}\widetilde{G}^{1-\nu
}(x,x_{t})\geq 0. 
\]

2/: D'après les estim\'{e}es ponctuelles faibles, dans $B(x_{t},\rho
)\backslash B(x_{t},R\mu _{t})$ :
\[
d_{\mathbf{g}}(x,x_{t})^{2}u_{t}(x)^{2^{*}-2}\leq \varepsilon _{R} 
\]
o\`{u} $\varepsilon _{R}\stackunder{R\rightarrow \infty }{\rightarrow }0$.
Alors avec (3.5) et (3.8), pour $R$ assez grand: 
\begin{eqnarray*}
\frac{L_{t}\widetilde{G}^{1-\nu }}{\widetilde{G}^{1-\nu }}(x,x_{t}) && \geq
\varepsilon _{0}-\lambda _{t}f(x)u_{t}(x)^{2^{*}-2}+\nu (1-\nu )\frac{c}{d_{%
\mathbf{g}}(x,x_{t})^{2}} \\ 
&& \geq \varepsilon _{0}-\lambda _{t}(\stackunder{B(x_{t},\rho )}{Sup}f).%
\frac{\varepsilon _{R}}{d_{\mathbf{g}}(x,x_{t})^{2}}+\nu (1-\nu )\frac{c}{d_{%
\mathbf{g}}(x,x_{t})^{2}} \\ 
&& \geq \varepsilon _{0}+\frac{c^{\prime }}{d_{\mathbf{g}}(x,x_{t})^{2}}\geq 0
\end{eqnarray*}

Remarque: Dans 1/, $\varepsilon _{0}$ sert \`{a} compenser $-\lambda
_{t}fu_{t}{}^{2^{*}-2}$ dans $M\backslash B(x_{t},\rho )$; et dans 2/, $\nu
(1-\nu )$ sert \`{a} compenser $-\lambda _{t}fu_{t}{}^{2^{*}-2}$ dans $%
B(x_{t},\rho )\backslash B(x_{t},R\mu _{t})$.

On a donc bien montr\'{e} que dans $M\backslash B(x_{t},R\mu _{t})$ et pour
toute constante $C_{t}>0$ d\'{e}pendant \'{e}ventuellement de $%
t $:
\[
L_{t}(C_{t}.\widetilde{G}^{1-\nu }(x,x_{t}))=C_{t}.L_{t}\widetilde{G}^{1-\nu
}(x,x_{t})\geq 0=L_{t}u_{t} 
\]
Enfin sur le bord $b(M\backslash B(x_{t},R\mu _{t}))$, en utilisant (3.4), on
obtient: 
\[
\widetilde{G}^{1-\nu }(x,x_{t})\geq \frac{c}{d_{\mathbf{g}%
}(x,x_{t})^{(n-2)(1-\nu )}}=\frac{c}{(R\mu _{t})^{(n-2)(1-\nu )}}\,\,. 
\]
Alors si on pose $C_{t}=c^{-1}R^{(n-2)(1-\nu )}\mu _{t}^{(n-2)(1-\nu )-\frac{%
n-2}{2}},$ on a pour $x\in bB(x_{t},R\mu _{t})=b(M\backslash B(x_{t},R\mu
_{t})):$%
\[
C_{t}.\widetilde{G}^{1-\nu }(x,x_{t})\geq \mu _{t}^{-\frac{n-2}{2}%
}=Sup\,u_{t}\geq u_{t}(x) 
\]
Donc par le principe du maximum 
\[
C_{t}.\widetilde{G}^{1-\nu }(x,x_{t})\geq u_{t}(x)\text{ dans }M\backslash
B(x_{t},R\mu _{t}) 
\]
soit 
\[
\widetilde{G}^{1-\nu }(x,x_{t})\geq C_{t}^{-1}u_{t}(x)=c\text{ }\mu _{t}^{%
\frac{n-2}{2}-(n-2)(1-\nu )}u_{t}(x) 
\]
et donc, en utilisant (3.4) : 
\[
d_{\mathbf{g}}(x,x_{t})^{(n-2)(1-\nu )}\mu _{t}^{\frac{n-2}{2}-(n-2)(1-\nu
)}u_{t}(x)\leq c 
\]
ce qui en changeant $\nu $ en $(n-2)\nu $ donne bien (3.7) sur $M$
puisque c'est vrai dans $B(x_{t},R\mu _{t})$ d'apr\`{e}s les estim\'{e}es
faibles.

On veut passer maintenant \`{a} $\nu =0$. Soit $y_{t}$ un point de $M$
o\`{u} $d_{\mathbf{g}}(x,x_{t})^{n-2}u_{t}(x_{t})u_{t}(x)$ est maximum. On veut montrer que $d_{\mathbf{g}%
}(y_{t},x_{t})^{n-2}u_{t}(x_{t})u_{t}(y_{t})$ est born\'{e}e. Soit $%
\varepsilon _{0}>0$ tel que $\bigtriangleup _{\mathbf{g}}+(h-\varepsilon
_{0})$ soit coercif et $\widehat{G}$ sa fonction de Green. Pour $t$ assez
proche de 1, on a $h_{t}\geq h-\varepsilon _{0}$. En utilisant le fait que $%
\widehat{G}>0,\,u_{t}>0,\,Maxf>0,$ on obtient: 
\begin{eqnarray*}
u_{t}(y_{t}) && =\int_{M}\widehat{G}(y_{t},x)(\bigtriangleup _{\mathbf{g}%
}u_{t}(x)+(h-\varepsilon _{0})u_{t}(x)) \\ 
&& \leq \int_{M}\widehat{G}(y_{t},x)(\bigtriangleup _{\mathbf{g}%
}u_{t}(x)+h_{t}u_{t}(x)) \\ 
&& =\lambda _{t}\int_{M}\widehat{G}(y_{t},x)f(x)u_{t}(x)^{2^{*}-1} \\ 
&& \leq \lambda _{t}Maxf\int_{M}\widehat{G}(y_{t},x)u_{t}(x)^{2^{*}-1} \\ 
&& \leq C\int_{B(x_{t},\delta )}\widehat{G}(y_{t},x)u_{t}(x)^{2^{*}-1}+C%
\int_{M\backslash B(x_{t},\delta )}\widehat{G}(y_{t},x)u_{t}(x)^{2^{*}-1}
\end{eqnarray*}
L'estim\'{e}e (3.7) , pour $0<\nu <\frac{2}{2^{*}-1}$, nous donne : 
\begin{eqnarray}
\int_{M\backslash B(x_{t},\delta )}\widehat{G}(y_{t},x)u_{t}(x)^{2^{*}-1} && 
=O(\mu _{t}^{(\frac{n+2}{2}-\nu )(2^{*}-1)}\int_{M\backslash B(x_{t},\delta
)}\widehat{G}(y_{t},x)u_{t}(x)^{2^{*}-1}) \nonumber\\ 
&& =O(\mu _{t}^{\frac{n+2}{2}-\nu (2^{*}-1)}) \nonumber\\ 
& &=o(\mu _{t}^{\frac{n-2}{2}})
\end{eqnarray}
Donc 
\begin{equation}
u_{t}(y_{t})\leq C\int_{B(x_{t},\delta )}\widehat{G}%
(y_{t},x)u_{t}(x)^{2^{*}-1}+o(\mu _{t}^{\frac{n-2}{2}})  
\end{equation}
On distingue alors trois cas:

Premier cas: $y_{t}\rightarrow y_{0}\neq x_{0}$

Soit $\delta >0$ fix\'{e} tel que $d_{\mathbf{g}}(y_{0},x_{0})\geq 3\delta $%
. Alors dans $B(x_{t},\delta ),$ $\widehat{G}(y_{t},x)$ est born\'{e}e et 
\begin{eqnarray*}
\int_{B(x_{t},\delta )}\widehat{G}(y_{t},x)u_{t}(x)^{2^{*}-1} && \leq
C\int_{B(x_{t},\delta )}u_{t}(x)^{2^{*}-1} \\ 
&& \leq C\int_{B(x_{t},\mu _{t})}u_{t}(x)^{2^{*}-1} \\ 
&& +C\mu _{t}^{\frac{n+2}{2}-\nu (2^{*}-1)}\int_{B(x_{t},\delta )\backslash
B(x_{t},\mu _{t})}d_{\mathbf{g}}(x,x_{t})^{(2^{*}-1)(\nu +2-n)} \\ 
&& \leq C\mu _{t}^{\frac{n-2}{2}}
\end{eqnarray*}
o\`{u} l'on a utilis\'{e} (3.7) en choisissant $\nu <\frac{2}{2^{*}-1}$
et les formules du changement d'\'{e}chelle pour $\int_{B(x_{t},\mu
_{t})}u_{t}{}^{2^{*}-1}=\mu _{t}^{\frac{n-2}{2}}\int_{B(0,1)}.\widetilde{u}%
_{t}^{2^{*}-1}dv_{\widetilde{g}_{t}}\sim C\mu _{t}^{\frac{n-2}{2}}$. Comme $%
d_{\mathbf{g}}(y_{t},x_{t})$ ne tend pas vers 0 on peut \'{e}crire 
\[
\int_{B(x_{t},\delta )}\widehat{G}(y_{t},x)u_{t}(x)^{2^{*}-1}\leq C\mu _{t}^{%
\frac{n-2}{2}}d_{\mathbf{g}}(y_{t},x_{t})^{2-n} 
\]
et donc en utilisant (3.10) et puisque $\nu <\frac{2}{2^{*}-1}$: 
\[
d_{\mathbf{g}}(y_{t},x_{t})^{n-2}u_{t}(x_{t})u_{t}(y_{t})\leq C\,. 
\]

Deuxi\`{e}me cas: Si, quitte \`{a} extraire, on a : $y_{t}\rightarrow x_{0}$
et $d_{\mathbf{g}}(y_{t},x_{t})\leq C\mu _{t}$.

Alors il est clair que 
\[
d_{\mathbf{g}}(y_{t},x_{t})^{n-2}u_{t}(x_{t})u_{t}(y_{t})\leq
Cu_{t}(x_{t})^{-1}u_{t}(y_{t})\leq C 
\]

Troisi\`{e}me cas: Si, quitte \`{a} extraire, on a: $y_{t}\rightarrow x_{0}$
et $r_{t}=\frac{d_{\mathbf{g}}(y_{t},x_{t})}{\mu _{t}}\rightarrow \infty $
quand $t\rightarrow 1$.

On \'{e}crit en utilisant (3.9) et (3.10) que 
\begin{eqnarray}
d_{\mathbf{g}}(y_{t},x_{t})^{n-2}u_{t}(x_{t})u_{t}(y_{t})  \leq &&C\mu _{t}^{-%
\frac{n-2}{2}}d_{\mathbf{g}}(y_{t},x_{t})^{n-2}\int_{B_{t}^{1}}\widehat{G}%
(y_{t},x)u_{t}(x)^{2^{*}-1}\nonumber \\ 
&& +C\mu _{t}^{-\frac{n-2}{2}}d_{\mathbf{g}}(y_{t},x_{t})^{n-2}%
\int_{B_{t}^{2}}\widehat{G}(y_{t},x)u_{t}(x)^{2^{*}-1}\nonumber \\ 
& &+O(\mu _{t}^{2-\nu (2^{*}-1)})
\end{eqnarray}
o\`{u} $B_{t}^{1}=\{x\in B(x_{t},\delta ):d_{\mathbf{g}}(y_{t},x)\geq \frac{1%
}{2}d_{\mathbf{g}}(y_{t},x_{t})\}$ et $B_{t}^{2}=B(x_{t},\delta )\backslash
B_{t}^{1}$. On suppose toujours $\nu <\frac{2}{2^{*}-1}$ et on utilise la
propri\'{e}t\'{e} (3.4) de la fonction de Green $\widehat{G}$. Alors (3.11) donne 
\begin{eqnarray*}
d_{\mathbf{g}}(y_{t},x_{t})^{n-2}u_{t}(x_{t})u_{t}(y_{t})  \leq&& C\mu _{t}^{-%
\frac{n-2}{2}}\int_{B_{t}^{1}}u_{t}{}^{2^{*}-1}+o(1) \\ 
&& +C\mu _{t}^{-\frac{n-2}{2}}d_{\mathbf{g}}(y_{t},x_{t})^{n-2}%
\int_{B_{t}^{2}}d_{\mathbf{g}}(y_{t},x)^{2-n}u_{t}(x)^{2^{*}-1}
\end{eqnarray*}
Comme dans le premier cas, en utilisant (3.7), on obtient: 
\[
\mu _{t}^{-\frac{n-2}{2}}\int_{B_{t}^{1}}u_{t}{}^{2^{*}-1}\leq C 
\]
et 
\begin{eqnarray*}
\mu _{t}^{-\frac{n-2}{2}}d_{\mathbf{g}}(y_{t},x_{t})^{n-2}\int_{B_{t}^{2}}d_{%
\mathbf{g}}(y_{t},x)^{2-n}u_{t}(x)^{2^{*}-1} && \leq \frac{1}{r_{t}^{2-\nu
(2^{*}-1)}d_{\mathbf{g}}(y_{t},x_{t})^{2}}\int_{B_{t}^{2}}d_{\mathbf{g}%
}(y_{t},x)^{2-n} \\ 
&& \leq \frac{C}{r_{t}^{2-\nu (2^{*}-1)}}=o(1)
\end{eqnarray*}
d\`{e}s que $\nu <\frac{2}{2^{*}-1}$ On d\'{e}duit donc l\`{a} encore que $%
d_{\mathbf{g}}(y_{t},x_{t})^{n-2}u_{t}(x_{t})u_{t}(y_{t})$ est born\'{e}.
\\

$\mathcal{7}$\textbf{Proposition 7:}

\textit{Concentration }$L^{p}$\textit{\ forte: }$\forall R>0$\textit{, }$%
\forall \delta >0$\textit{\ et }$\forall p>\frac{n}{n-2}$\textit{\ o\`{u} }$%
n=\dim M$\textit{\ : } 
\[
\stackunder{t\rightarrow 1}{\lim }\frac{\int_{B(x_{t},R\mu
_{t})}u_{t}^{p}dv_{\mathbf{g}}}{\int_{B(x_{t},\delta )}u_{t}^{p}dv_{\mathbf{g%
}}}=1-\varepsilon _{R}\text{ o\`{u} }\varepsilon _{R}\stackunder{%
R\rightarrow +\infty }{\rightarrow }0\text{ } 
\]

\textit{D\'{e}monstration: }Il suffit d'appliquer la proposition
pr\'{e}c\'{e}dente dans le changement d'\'{e}chelle en $x_{t}.$ Par
changement d'\'{e}chelle 
\begin{eqnarray*}
\int_{M}u_{t}^{p}dv_{\mathbf{g}} && \geq \int_{B(x_{t},\mu _{t})}u_{t}^{p}dv_{%
\mathbf{g}}=\mu _{t}^{n-\frac{n-2}{2}p}\int_{B(0,1)}\widetilde{u}_{t}^{p}dv_{%
\widetilde{g}_{t}} \\ 
& &\geq C\mu _{t}^{n-\frac{n-2}{2}p}
\end{eqnarray*}
D'autre part d'apr\`{e}s la proposition 6: 
\begin{eqnarray*}
\int_{M\backslash B(x_{t},R\mu _{t})}u_{t}^{p}dv_{\mathbf{g}} && \leq C\mu
_{t}^{p\frac{n-2}{2}}\int_{M\backslash B(x_{t},R\mu _{t})}d_{\mathbf{g}%
}(x_{t},x)^{(2-n)p}dv_{\mathbf{g}} \\ 
&& \leq C\mu _{t}^{n-p\frac{n-2}{2}}R^{n+(2-n)p}
\end{eqnarray*}
dès que $p>\frac{n}{n-2}$. En quotientant on obtient le corollaire.

\section{Argument central de la d\'{e}monstration du th\'{e}or\`{e}me 1}

Comme expos\'{e} au chapitre pr\'{e}c\'{e}dent dans la description de
l'id\'{e}e de la d\'{e}monstration, on ins\`{e}re dans l'in\'{e}galit\'{e}
de Sobolev euclidienne l'\'{e}quation ($E_{t})$ ``lue'' dans la carte $\exp
_{x_{t}}^{-1}$.

On reprend les notations de 3.2: $\overline{u}_{t},\overline{f}_{t},%
\overline{h}_{t},\mathbf{g}_{t}$ d\'{e}signent les fonctions et la
m\'{e}trique ``lues'' dans la carte $\exp _{x_{t}}^{-1}$, et $\,\,\widetilde{%
u}_{t}\,,\,\widetilde{h}_{t}\,,\,\widetilde{f}_{t},\widetilde{\mathbf{g}}%
_{t} $ d\'{e}signent les fonctions et la m\'{e}trique ``lues'' apr\`{e}s
blow-up de centre $x_{t}$ et de coefficient $k_{t}=\mu _{t}^{-1}$; et on se place sur un voisinage de $x_{0}$ o\`{u} $f\geq 0$. On
consid\`{e}re de plus une fonction cut-off $\eta $ sur $\Bbb{R}^{n}$
\'{e}gale \`{a} 1 sur la boule euclidienne $B(0,\delta /2)$, \'{e}gale \`{a}
0 sur $\Bbb{R}^{n}\backslash B(0,\delta )$, $0\leq \eta \leq 1$ avec $\left|
\nabla \eta \right| \leq C.\delta ^{-1}$ o\`{u} $\delta $ est assez petit
pour que, sur les boules $B(x_{t},\delta )$, $f\geq 0$. L'identit\'{e} de
Sobolev euclidienne donne d'une part 
\begin{equation}
(\int_{B(0,\delta )}(\eta \overline{u}_{t})^{2^{*}}dx)^{\frac{2}{2^{*}}}\leq
K(n,2)^{2}\int_{B(0,\delta )}\left| \nabla (\eta \overline{u}_{t})\right|
_{e}^{2}dx\,\, 
\end{equation}
\thinspace o\`{u} $\left| .\right| _{e}$ est la m\'{e}trique euclidienne de
mesure associ\'{e}e $dx$.

Par ailleurs, une int\'{e}gration par parties donne en notant que $\left|
\nabla \eta \right| =$ $\Delta \eta =0$ sur $B(0,\delta /2)$ : 
\[
\int_{B(0,\delta )}\left| \nabla (\eta \overline{u}_{t})\right|
_{e}^{2}dx\leq \int_{B(0,\delta )}\eta ^{2}\overline{u}_{t}\bigtriangleup
_{e}\overline{u}_{t}dx+C.\delta ^{-2}\int_{B(0,\delta )\backslash B(0,\delta
/2)}\overline{u}_{t}^{2}dx 
\]
En notant $\,\mathbf{g}\,_{t}^{ij}$ les composantes de $\,\mathbf{g}\,_{t}$
et $\Gamma (\,\mathbf{g}\,_{t})_{ij}^{k}$ les symboles de Chrisoffel
associ\'{e}s, on \'{e}crit : 
\[
\bigtriangleup _{e}\overline{u}_{t}=\bigtriangleup _{\mathbf{g}_{t}}%
\overline{u}_{t}+(\,\mathbf{g}\,_{t}^{ij}-\delta ^{ij})\partial _{ij}%
\overline{u}_{t}-\,\mathbf{g}\,_{t}^{ij}\Gamma (\,\mathbf{g}%
\,_{t})_{ij}^{k}\partial _{k}\overline{u}_{t} 
\]
On obtient alors \`{a} partir de cette derni\`{e}re in\'{e}galit\'{e} (point
(2.6) de l'id\'{e}e d\'{e}crite \`{a} la fin du chapitre 2), en utilisant
cette expression du laplacien comparant $\bigtriangleup _{e}$ et $%
\bigtriangleup _{\mathbf{g}_{t}}$, en utilisant l'\'{e}quation $E_{t}:\,\,\triangle _{%
\mathbf{g}}u_{t}+h_{t}.u_{t}=\lambda _{t}.f.u_{t}^{\frac{n+2}{n-2}}$ ``lue''
dans la carte exp$_{x_{t}}^{-1}$, et en utilisant le fait que $\left| \nabla \eta \right|
=\Delta \eta =0$ sur $B(0,\delta /2)$ et en faisant quelques
int\'{e}grations par parties: 
\begin{eqnarray*}
\int_{B(0,\delta )}\left| \nabla (\eta \overline{u}_{t})\right|
_{e}^{2}dx\leq & &\lambda _{t}\int_{B(0,\delta )}\eta ^{2}\overline{f}_{t}%
\overline{u}_{t}^{2^{*}}dx-\int_{B(0,\delta )}\eta ^{2}\overline{h}_{t}%
\overline{u}_{t}^{2}dx\\
&&+C.\delta ^{-2}\int_{B(0,\delta )\backslash B(0,\delta
/2)}\overline{u}_{t}^{2}dx \\ 
&& -\int_{B(0,\delta )}\eta ^{2}(\,\mathbf{g}\,_{t}^{ij}-\delta
^{ij})\partial _{i}\overline{u}_{t}\partial _{j}\overline{u}_{t}dx\\
&&+\frac{1}{2}\int_{B(0,\delta )}(\partial _{k}(\,\mathbf{g}\,_{t}^{ij}\Gamma (\,\mathbf{g%
}\,_{t})_{ij}^{k}+\partial _{ij}\,\mathbf{g}\,_{t}^{ij})(\eta \overline{u}%
_{t}^{2})dx\,.
\end{eqnarray*}
On en d\'{e}duit gr\^{a}ce \`{a} l'in\'{e}galit\'{e} de Sobolev (3.12) et
en utilisant le fait que $\lambda _{t}\leq \frac{1}{K(n,2){{}^{2}}(%
\stackunder{M}{Sup}f)^{\frac{n-2}{n}}}$: 

\begin{equation}
\int_{B(0,\delta )}\overline{h}_{t}(\eta \overline{u}_{t})^{2}dx\leq A_{t}
+B_{t}+C_{t}+C.\delta ^{-2}\int_{B(0,\delta )\backslash B(0,\delta /2)}\overline{u}%
_{t}^{2}dx
\end{equation}
avec:

$B_{t}=\frac{1}{2}\int_{B(0,\delta )}(\partial _{k}(\,\mathbf{g}%
\,_{t}^{ij}\Gamma (\,\mathbf{g}\,_{t})_{ij}^{k}+\partial _{ij}\,\mathbf{g}%
\,_{t}^{ij})(\eta \overline{u}_{t}^{2})dx$

$C_{t}=\left| \int_{B(0,\delta )}\eta ^{2}(\,\mathbf{g}\,_{t}^{ij}-\delta
^{ij})\partial _{i}\overline{u}_{t}\partial _{j}\overline{u}_{t}dx\right| $

$A_{t}=\frac{1}{K(n,2){{}^{2}}(\stackunder{M}{Sup}f)^{\frac{n-2}{n}}}%
\int_{B(0,\delta )}\overline{f}_{t}\eta ^{2}\overline{u}_{t}^{2^{*}}dx-\frac{%
1}{K(n,2){{}^{2}}}(\int_{B(0,\delta )}(\eta \overline{u}_{t})^{2^{*}}dx)^{%
\frac{2}{2^{*}}}$

Ces calculs sont d\'{e}taill\'{e}s dans l'article de Z. Djadli et O. Druet $%
\left[ 9\right] $ sur lequel nous nous appuyons. Le but va \^{e}tre
d'utiliser la ``concentration $L{{}^{2}}$'' \'{e}tudi\'{e}e en 3.2.d/,$\,$%
pour obtenir une contradiction; nous diviserons ainsi (3.13) par $%
\int_{B(0,\delta )}\overline{u}_{t}^{2}dx$ et ferons tendre $t\rightarrow
t_{0}=1$.

La concentration $L{{}^{2}}$ nous donne d\'{e}j\`{a}: 
\[
\frac{C.\delta ^{-2}\int_{B(0,\delta )\backslash B(0,\delta /2)}\overline{u}%
_{t}^{2}dx}{\int_{B(0,\delta )}\overline{u}_{t}^{2}dx}\stackunder{%
t\rightarrow 1}{\rightarrow }0\,\,. 
\]

Z.Djadli et O.Druet $\left[ 9\right] $ ont montr\'{e} (voir l'appendice C)
que : 
\[
\stackunder{t\rightarrow 1}{\overline{\lim }}\frac{C_{t}}{\int_{B(0,\delta )}%
\overline{u}_{t}^{2}dx}\leq \varepsilon _{\delta }\text{ o\`{u} }\varepsilon
_{\delta }\rightarrow 0\text{ quand }\delta \rightarrow 0\,\,.
\]

De plus, comme $x_{t}\rightarrow x_{0}$ on a $\stackunder{t\rightarrow 1}{%
\lim }(\partial _{k}(\,\mathbf{g}\,_{t}^{ij}\Gamma (\,\mathbf{g}%
\,_{t})_{ij}^{k}+\partial _{ij}\,\mathbf{g}\,_{t}^{ij})(0)=\frac{1}{3}S_{\,%
\mathbf{g}\,}(x_{0})$, d'o\`{u}, d'apr\`{e}s la concentration $L{{}^{2}}$ : 
\[
\stackunder{t\rightarrow 1}{\overline{\lim }}\frac{B_{t}}{\int_{B(0,\delta )}%
\overline{u}_{t}^{2}dx}=\frac{1}{6}S_{\,\mathbf{g}\,}(x_{0})+\varepsilon
_{\delta }\,\,. 
\]

C'est l'expression $A_{t}$ qui va nous donner les termes $\frac{n-2}{4(n-1)}%
S_{\,\mathbf{g}\,}(x_{0})$ $-\frac{1}{6}S_{\,\mathbf{g}\,}(x_{0})$ et $\frac{%
(n-2)(n-4)}{8(n-1)}\frac{\bigtriangleup _{\mathbf{g}}f(x_{0})}{f(x_{0})}$.

L'in\'{e}galit\'{e} de H\"{o}lder nous donne: 
\[
\int_{B(0,\delta )}\overline{f}_{t}\eta ^{2}\overline{u}_{t}^{2^{*}}dx\leq
(\int_{B(0,\delta )}\overline{f}_{t}\overline{u}_{t}^{2^{*}}dx)^{\frac{2}{n}%
}(\int_{B(0,\delta )}\overline{f}_{t}(\eta \overline{u}_{t})^{2^{*}}dx)^{%
\frac{n-2}{n}}\, 
\]
D'autre part : 
\[
\,dx\leq (1+\frac{1}{6}Ric(x_{t})_{ij}x^{i}x^{j}+C\left| x\right| ^{3})dv_{\,%
\mathbf{g}\,_{t}} 
\]
Un d\'{e}veloppement limit\'{e} nous permet alors d'\'{e}crire: 
\[
(\int_{B(0,\delta )}\overline{f}_{t}\overline{u}_{t}^{2^{*}}dx)^{\frac{2}{n}%
}\leq (\int_{B(0,\delta )}\overline{f}_{t}\overline{u}_{t}^{2^{*}}dv_{\,%
\mathbf{g}\,_{t}})^{\frac{2}{n}}+\frac{1}{(\int_{B(0,\delta )}\overline{f}%
_{t}\overline{u}_{t}^{2^{*}}dv_{\,\mathbf{g}\,_{t}})^{\frac{n-2}{n}}}\frac{2%
}{n}\{S_{t}\}+C\{S_{t}\}^{2} 
\]
o\`{u} 
\[
\{S_{t}\}=\frac{1}{6}Ric(x_{t})_{ij}\int_{B(0,\delta )}x^{i}x^{j}\overline{f}%
_{t}\overline{u}_{t}^{2^{*}}dv_{\,\mathbf{g}\,_{t}}+C\int_{B(0,\delta
)}\left| x\right| ^{3}\overline{u}_{t}^{2^{*}}dv_{\,\mathbf{g}\,_{t}}\,\,. 
\]
On en d\'{e}duit 
\[
A_{t}\leq \frac{1}{K(n,2){{}^{2}}(\stackunder{M}{Sup}f)^{\frac{n-2}{n}}}%
(A_{t}^{1}+A_{t}^{2}) 
\]
o\`{u} 
\[
A_{t}^{1}=(\int_{B(0,\delta )}\overline{f}_{t}\overline{u}_{t}^{2^{*}}dv_{\,%
\mathbf{g}\,_{t}})^{\frac{2}{n}}(\int_{B(0,\delta )}\overline{f}_{t}(\eta 
\overline{u}_{t})^{2^{*}}dx)^{\frac{n-2}{n}}\,-(Supf.\int_{B(0,\delta
)}(\eta \overline{u}_{t})^{2^{*}}dx)^{\frac{n-2}{n}} 
\]
et 
$$
A_{t}^{2}=\frac{2(\int_{B(0,\delta )}\overline{f}_{t}(\eta \overline{u}%
_{t})^{2^{*}}dx)^{\frac{n-2}{n}}}{n(\int_{B(0,\delta )}\overline{f}_{t}%
\overline{u}_{t}^{2^{*}}dv_{\,\mathbf{g}\,_{t}})^{\frac{n-2}{n}}}
\{\frac{1}{6}Ric(x_{t})_{ij}\int_{B(0,\delta )}x^{i}x^{j}\overline{f}_{t}%
\overline{u}_{t}^{2^{*}}dv_{\,\mathbf{g}\,_{t}}
+C\int_{B(0,\delta )}\left|
x\right| ^{3}\overline{u}_{t}^{2^{*}}dv_{\,\mathbf{g}\,_{t}}\}(1+\varepsilon
_{\delta }) 
$$
car $\{S_{t}\}\rightarrow 0$ quand $\delta \rightarrow 0\,$uniform\'{e}ment
en $t$. Ainsi l'expression $A_{t}^{2}$ va donner, par d\'{e}veloppement
limit\'{e} \`{a} l'ordre 2 de la m\'{e}trique, le terme en $S_{\,\mathbf{g}%
\,}(x_{0})\,$et l'expression $A_{t}^{1}$ donnera, par d\'{e}veloppement
limit\'{e} \`{a} l'ordre 2 de $f$, le terme en $-\frac{(n-2)(n-4)}{8(n-1)}%
\frac{\bigtriangleup _{\mathbf{g}}f(x_{0})}{f(x_{0})}$.

\textbf{Lemme:}\textit{\ Comme }$dx=(1+\frac{1}{6}%
Ric(x_{t})_{ij}x^{i}x^{j}+O(\left| x\right| ^{3}))dv_{\,\mathbf{g}\,_{t}}$%
\textit{\ et puisque }$\exp :M\times \Bbb{R}^{n}\rightarrow \Bbb{R}^{n}$%
\textit{\ est une application de classe }$C^{\infty }$\textit{, on a pour
toute fonction }$\alpha \in H_{1}^{2}(B(x_{0},2\delta )):$%
\[
\stackunder{t\rightarrow 1}{\lim }\frac{\int_{B(x_{t},\delta )}\alpha dx}{%
\int_{B(x_{t},\delta )}\alpha dv_{\,\mathbf{g}\,_{t}}}=1+O(\delta
^{2})=1+\varepsilon _{\delta } 
\]

Etudions d'abord l'expression $A_{t}^{2}$:

1/: On a $\stackunder{t\rightarrow 1}{\lim }\frac{(\int_{B(0,\delta )}%
\overline{f}_{t}(\eta \overline{u}_{t})^{2^{*}}dx)^{\frac{n-2}{n}}}{%
(\int_{B(0,\delta )}\overline{f}_{t}\overline{u}_{t}^{2^{*}}dv_{\,\mathbf{g}%
\,_{t}})^{\frac{n-2}{n}}}=1+\varepsilon _{\delta }$

2/: D'apr\`{e}s les estim\'{e}es ponctuelles faibles, $\left| x\right| ^{2}%
\overline{u}_{t}^{2^{*}}\leq c\overline{u}_{t}^{2}$, on en d\'{e}duit: 
\[
\frac{\int_{B(0,\delta )}\left| x\right| ^{3}\overline{u}_{t}^{2^{*}}dv_{\,%
\mathbf{g}\,_{t}}}{\int_{B(0,\delta )}\overline{u}_{t}^{2}dv_{\,\mathbf{g}%
\,_{t}}}\leq C.\varepsilon _{\delta }\,\,. 
\]

3/: D'apr\`{e}s les formules de changement d'\'{e}chelle on \'{e}crit: pour
tout $R>0:$%
\begin{eqnarray*}
\int_{B(0,\delta )}x^{i}x^{j}\overline{f}_{t}\overline{u}_{t}^{2^{*}}dv_{\,%
\mathbf{g}\,_{t}} && =\int_{B(0,R\mu _{t})}x^{i}x^{j}\overline{f}_{t}%
\overline{u}_{t}^{2^{*}}dv_{\,\mathbf{g}\,_{t}}+\int_{B(0,\delta )\backslash
B(0,R\mu _{t})}x^{i}x^{j}\overline{f}_{t}\overline{u}_{t}^{2^{*}}dv_{\,%
\mathbf{g}\,_{t}} \\ 
&& =\mu _{t}^{2}\int_{B(0,R)}x^{i}x^{j}\widetilde{f_{t}}\widetilde{u}%
_{t}^{2^{*}}dv_{\widetilde{\,\mathbf{g}\,}_{t}}+\mu _{t}^{2}\int_{B(0,\delta
\mu _{t}^{-1})\backslash B(0,R)}x^{i}x^{j}\widetilde{f_{t}}\widetilde{u}%
_{t}^{2^{*}}dv_{\widetilde{\,\mathbf{g}\,}_{t}}
\end{eqnarray*}
et 
\[
\int_{B(0,\delta )}\overline{u}_{t}^{2}dv_{\mathbf{g}_{t}}=\mu
_{t}^{2}\int_{B(0,\delta \mu _{t}^{-1})}\widetilde{u}_{t}^{2}dv_{\widetilde{%
\mathbf{g}}_{t}}\,\,. 
\]
En utilisant les estim\'{e}es ponctuelles faibles, on obtient : 
\[
\int_{B(0,\delta \mu _{t}^{-1})\backslash B(0,R)}x^{i}x^{j}\widetilde{f_{t}}%
\widetilde{u}_{t}^{2^{*}}dv_{\widetilde{\mathbf{g}}_{t}}\leq \varepsilon
_{R}.\int_{B(0,\delta \mu _{t}^{-1})\backslash B(0,R)}\widetilde{u}%
_{t}^{2}dv_{\widetilde{\mathbf{g}}_{t}} 
\]
$\,$donc: 
\[
\frac{\int_{B(0,\delta \mu _{t}^{-1})\backslash B(0,R)}x^{i}x^{j}\widetilde{%
f_{t}}\widetilde{u}_{t}^{2^{*}}dv_{\widetilde{\mathbf{g}}_{t}}}{%
\int_{B(0,\delta \mu _{t}^{-1})}\widetilde{u}_{t}^{2}dv_{\widetilde{\mathbf{g%
}}_{t}}}\leq \varepsilon _{R} 
\]
o\`{u} $\varepsilon _{R}\rightarrow 0$ quand $R\rightarrow +\infty $.
Maintenant, si $i\neq j:$%
\[
\stackunder{t\rightarrow 1}{\overline{\lim }}\frac{|\int_{B(0,R)}x^{i}x^{j}%
\widetilde{f_{t}}\widetilde{u}_{t}^{2^{*}}dv_{\widetilde{\mathbf{g}}_{t}}|}{%
\int_{B(0,\delta \mu _{t}^{-1})}\widetilde{u}_{t}^{2}dv_{\widetilde{\mathbf{g%
}}_{t}}}\leq \stackunder{t\rightarrow 1}{\overline{\lim }}\frac{%
|\int_{B(0,R)}x^{i}x^{j}\widetilde{f_{t}}\widetilde{u}_{t}^{2^{*}}dv_{%
\widetilde{\mathbf{g}}_{t}}|}{\int_{B(0,R)}\widetilde{u}_{t}^{2}dv_{%
\widetilde{\mathbf{g}}_{t}}}=0 
\]
car 
\[
\widetilde{u}_{t}\rightarrow \widetilde{u}=(1+\frac{f(x_{0})^{\frac{2}{n}}}{%
K(n,2)^{2}n(n-2)}\left| x\right| ^{2})^{-\frac{n-2}{2}}\,\,\,dans\,\,%
\,C^{0}(B(0,R)) 
\]
et $\widetilde{u}$ est radiale (voir le paragraphe changement d'\'{e}chelle).

Si $i=j:$%
\[
\frac{\int_{B(0,R)}x^{i}x^{i}\widetilde{f_{t}}\widetilde{u}_{t}^{2^{*}}dv_{%
\widetilde{\mathbf{g}}_{t}}}{\int_{B(0,\delta \mu _{t}^{-1})}\widetilde{u}%
_{t}^{2}dv_{\widetilde{\mathbf{g}}_{t}}}=\frac{\int_{B(0,R)}(x^{i}){{}^{2}}%
\widetilde{f_{t}}\widetilde{u}_{t}^{2^{*}}dv_{\widetilde{\mathbf{g}}_{t}}}{%
\int_{B(0,R)}\widetilde{u}_{t}^{2}dv_{\widetilde{\mathbf{g}}_{t}}}.\frac{%
\int_{B(0,R)}\widetilde{u}_{t}^{2}dv_{\widetilde{\mathbf{g}}_{t}}}{%
\int_{B(0,\delta \mu _{t}^{-1})}\widetilde{u}_{t}^{2}dv_{\widetilde{\mathbf{g%
}}_{t}}} 
\]
Or d\`{e}s que $n>4$, en utilisant la convergence $L{{}^{2}}$ forte
(proposition 7), on obtient: 
\[
\stackunder{R\rightarrow \infty }{\lim }\,\stackunder{t\rightarrow 1}{%
\overline{\lim }}\frac{\int_{B(0,R)}\widetilde{u}_{t}^{2}dv_{\widetilde{%
\mathbf{g}}_{t}}}{\int_{B(0,\delta \mu _{t}^{-1})}\widetilde{u}_{t}^{2}dv_{%
\widetilde{\mathbf{g}}_{t}}}=1 
\]
donc

\[
\stackunder{R\rightarrow \infty }{\lim }\,\stackunder{t\rightarrow 1}{%
\overline{\lim }}\frac{\int_{B(0,R)}x^{i}x^{i}\widetilde{f_{t}}\widetilde{u}%
_{t}^{2^{*}}dv_{\widetilde{\mathbf{g}}_{t}}}{\int_{B(0,\delta \mu _{t}^{-1})}%
\widetilde{u}_{t}^{2}dv_{\widetilde{\mathbf{g}}_{t}}}=f(x_{0})\frac{\int_{%
\Bbb{R}^{n}}(x^{i})^{2}.\widetilde{u}^{2}dx}{\int_{\Bbb{R}^{n}}\widetilde{u}%
^{2}dx}=f(x_{0})^{\frac{n-2}{n}}K(n,2)^{2}\frac{n-4}{4(n-1)} 
\]
d'o\`{u} 
\[
\stackunder{t\rightarrow 1}{\overline{\lim }}\frac{1}{f(x_{0})^{\frac{n-2}{n}%
}K(n,2)^{2}}\frac{A_{t}^{2}}{\int_{B(0,\delta )}\overline{u}_{t}^{2}dv_{%
\mathbf{g}_{t}}}=\frac{n-4}{12(n-1)}S_{\mathbf{g}}(x_{0})+\varepsilon
_{\delta } 
\]
ce qui avec $\stackunder{t\rightarrow 1}{\overline{\lim }}\frac{B_{t}}{%
\int_{B(0,\delta )}\overline{u}_{t}^{2}dx}=\frac{1}{6}S_{\mathbf{g}%
}(x_{0})+\varepsilon _{\delta }$ donne 
\[
\stackunder{t\rightarrow 1}{\overline{\lim }}(\frac{1}{f(x_{0})^{\frac{n-2}{n%
}}K(n,2)^{2}}\frac{A_{t}^{2}}{\int_{B(0,\delta )}\overline{u}_{t}^{2}dv_{%
\mathbf{g}_{t}}}+\frac{B_{t}}{\int_{B(0,\delta )}\overline{u}_{t}^{2}dx})=%
\frac{n-2}{4(n-1)}S_{\mathbf{g}}(x_{0})+\varepsilon _{\delta } 
\]

Si $n=4$ on \'{e}crit: 
\[
\stackunder{R\rightarrow \infty }{\lim }\stackunder{t\rightarrow 1}{%
\overline{\lim }}\frac{\int_{B(0,R)}x^{i}x^{i}\widetilde{f_{t}}\widetilde{u}%
_{t}^{2^{*}}dv_{\widetilde{\mathbf{g}}_{t}}}{\int_{B(0,\delta \mu _{t}^{-1})}%
\widetilde{u}_{t}^{2}dv_{\widetilde{\mathbf{g}}_{t}}}\leq f(x_{0})^{\frac{n-2%
}{n}}K(n,2)^{2}\frac{n-4}{4(n-1)} 
\]
et on conclut en distinguant deux cas, $S_{\mathbf{g}}(x_{0})<0$ ou $S_{%
\mathbf{g}}(x_{0})\geq 0$, la d\'{e}monstration est alors termin\'{e}e car $%
\frac{\bigtriangleup _{\mathbf{g}}f(x_{0})}{f(x_{0})}$ n'appara\^{i}t pas
(voir l'argument \`{a} la fin de la d\'{e}monstration).

Passons \`{a} l'\'{e}tude de l'expression $A_{t}^{1}$. C'est l\`{a}
qu'appara\^{i}t la principale difficult\'{e} due \`{a} la pr\'{e}sence d'une
fonction $f$ non constante au second membre.

L'id\'{e}e est la suivante: nous voudrions reprendre la m\'{e}thode
utilis\'{e}e ci-dessus pour l'expression $A_{t}^{2},$ o\`{u} nous avons
utilis\'{e} un d\'{e}veloppement limit\'{e} de la m\'{e}trique pour obtenir
le terme $S_{\mathbf{g}}(x_{0})$, en faisant un d\'{e}veloppement limit\'{e}
de $f$ pour en faire appara\^{i}tre les d\'{e}riv\'{e}es secondes et ainsi
obtenir $\frac{\bigtriangleup _{\mathbf{g}}f(x_{0})}{f(x_{0})}$. Mais la
difficult\'{e} qui appara\^{i}t alors vient du fait que toutes les
estim\'{e}es obtenues dans l'\'{e}tude du ph\'{e}nom\`{e}ne de concentration
sont centr\'{e}es en $x_{t}$. Si l'on d\'{e}veloppe $f$ en $x_{t}$, ce sont
les d\'{e}riv\'{e}es premi\`{e}res $\partial _{i}f(x_{t})$ qui interviennent
et le changement d'\'{e}chelle nous donnera alors 
\[
\int_{B(0,\delta )}\partial _{i}f(x_{t})x^{i}\overline{u}_{t}^{2^{*}}dv_{%
\mathbf{g}_{t}}=\mu _{t}\int_{B(0,\delta \mu _{t}^{-1})}\partial
_{i}f(x_{t})x^{i}\widetilde{u}_{t}^{2^{*}}dv_{\widetilde{\mathbf{g}}_{t}} 
\]
que l'on doit diviser par 
\[
\mu _{t}^{2}\int_{B(0,\delta \mu _{t}^{-1})}\widetilde{u}_{t}^{2}dv_{%
\widetilde{\mathbf{g}}_{t}}\,\,. 
\]
Il est alors n\'{e}cessaire de contrôler le rapport $\frac{\partial
_{i}f(x_{t})}{\mu _{t}}$, ce qui semble difficile. Si l'on choisit
plut\^{o}t de d\'{e}velopper $f$ au point de maximum $x_{0}$, les
d\'{e}riv\'{e}es premi\`{e}res $\partial _{i}f(x_{0})$ s'annulent, mais il
faut alors ``transposer'' les estim\'{e}es faites en $x_{t}$ au point $x_{0}$
ce qui impose cette fois l'obtention d'une in\'{e}galit\'{e} de la forme: 
\[
\frac{d_{g}(x_{t},x_{0})}{\mu _{t}}\leq C\,\,. 
\]
Cette in\'{e}galit\'{e} a \'{e}t\'{e} \'{e}tudi\'{e}e par les auteurs d\'{e}j%
\`{a} cit\'{e}s, elle est difficile \`{a} obtenir, et elle n\'{e}cessite
certaines hypoth\`{e}ses suppl\'{e}mentaires. Elle est baptis\'{e}e par Zo%
\'{e} Faget ``seconde in\'{e}galit\'{e} fondamentale'' (la premi\`{e}re \'{e}%
tant l'estim\'{e}e ponctuelle faible). Ainsi, dans l'article de O. Druet et
F. Robert [13], o\`{u} $f=cste$, il faut faire des hypoth\`{e}ses sur la
``forme'' des $h_{t}$ et sur la g\'{e}om\'{e}trie de la vari\'{e}t\'{e} au
point de concentration \'{e}ventuel pour obtenir le r\'{e}sultat. C'est
l'hypoth\`{e}se disant que les maxima de $f$ sont localement stricts
(hessien non d\'{e}g\'{e}n\'{e}r\'{e}) qui nous permettra d'obtenir cette in%
\'{e}galit\'{e}; intuitivement cette hypoth\`{e}se ``fixe'' la position du
point de concentration $x_{0}.$ En fait, la m\'{e}thode que nous d\'{e}%
velopperons donnera directement le r\'{e}sultat cherch\'{e} sur $\frac{%
\bigtriangleup _{\mathbf{g}}f(x_{0})}{f(x_{0})}\,$, et donnera cette estim%
\'{e}e \textit{a posteriori}; mais nous aurions pu faire l'inverse.

Notons $x_{0}(t)=\exp _{x_{t}}^{-1}(x_{0})=(x_{0}^{1}(t),...,x_{0}^{n}(t))$, 
ce qui a un sens dès que $t$ est assez proche de 1 pour un rayon $\delta$ fixé. Alors $x_{0}(t)\rightarrow 0$
 quand $t\rightarrow 1$. Le point $x_{0}(t)$ est un maximum localement strict de $\overline{f}_{t}$. 
 Nous ferons tendre $\delta$ vers 0 à la fin du raisonnement après avoir pris les limites en $t$. 

Le d\'{e}veloppement limit\'{e} de $\overline{f%
}_{t}$ en $x_{0}(t)$ donne: 
\[
\overline{f}_{t}(x)\leq f(x_{0})+\frac{1}{2}\partial _{kl}\overline{f}%
_{t}(x_{0}(t)).(x^{k}-x_{0}^{k}(t))(x^{l}-x_{0}^{l}(t))+c\left|
x-x_{0}(t)\right| ^{3}:=f(x_{0})+T_{t} 
\]
($T_{t}$ comme Taylor) o\`{u} ($\partial _{kl}\overline{f}_{t}(x_{0})$) est
une matrice d\'{e}finie n\'{e}gative (on notera $<0$). On notera
toujours $c,C$ des constantes ind\'{e}pendantes de $t$ et $\delta $.
Rappelons que 
\[
A_{t}^{1}=(\int_{B(0,\delta )}\overline{f}_{t}\overline{u}_{t}^{2^{*}}dv_{%
\mathbf{g}_{t}})^{\frac{2}{n}}(\int_{B(0,\delta )}\overline{f}_{t}(\eta 
\overline{u}_{t})^{2^{*}}dx)^{\frac{n-2}{n}}\,-(Supf.\int_{B(0,\delta
)}(\eta \overline{u}_{t})^{2^{*}}dx)^{\frac{n-2}{n}}\,\,. 
\]
En introduisant le d\'{e}veloppement limit\'{e} de $\overline{f}_{t}$ en $%
x_{0}(t)$ on obtient: 
\[
(\int_{B(0,\delta )}\overline{f}_{t}(\eta \overline{u}_{t})^{2^{*}}dx)^{%
\frac{n-2}{n}}\leq (\int_{B(0,\delta )}f(x_{0})(\eta \overline{u}%
_{t})^{2^{*}}dx)^{\frac{n-2}{n}}+\frac{\frac{n-2}{n}}{(\int_{B(0,\delta
)}f(x_{0})(\eta \overline{u}_{t})^{2^{*}}dx)^{\frac{2}{n}}}%
\{F_{t}\}+C.\{F_{t}\}^{2} 
\]
o\`{u} 
\[
\{F_{t}\}=\frac{1}{2}\partial _{kl}\overline{f}_{t}(x_{0}(t))\int_{B(0,%
\delta )}(x^{k}-x_{0}^{k}(t))(x^{l}-x_{0}^{l}(t))(\eta \overline{u}%
_{t})^{2^{*}}dx+C\int_{B(0,\delta )}\left| x-x_{0}(t)\right| ^{3}(\eta 
\overline{u}_{t})^{2^{*}}dx 
\]
d'o\`{u} en rappelant que $\stackunder{M}{Sup}f=f(x_{0})$ et $%
\int_{B(0,\delta )}\overline{f}_{t}\overline{u}_{t}^{2^{*}}dv_{\mathbf{g}%
_{t}}\leq 1$: 
\begin{equation}
A_{t}^{1}\leq \frac{n-2}{n}\frac{(\int_{B(0,\delta )}\overline{f}_{t}%
\overline{u}_{t}^{2^{*}}dv_{\mathbf{g}_{t}})^{\frac{2}{n}}}{%
(\int_{B(0,\delta )}f(x_{0})(\eta \overline{u}_{t})^{2^{*}}dx)^{\frac{2}{n}}}%
\{F_{t}\}(1+\varepsilon _{\delta,t })  
\end{equation}
car $C\left| F_{t}\right| =\varepsilon _{\delta,t }\rightarrow 0$
 quand $\delta \rightarrow 0$ et $t\rightarrow 1$. Par ailleurs, on peut
\'{e}crire 
\[
\{F_{t}\}=\int_{B(0,\delta )}T_{t}.(\eta \overline{u}_{t})^{2^{*}}dx=%
\int_{B(0,\delta )}T_{t}.(\eta \overline{u}_{t})^{2^{*}}dv_{\mathbf{g}_{t}}%
\frac{\int_{B(0,\delta )}T_{t}.(\eta \overline{u}_{t})^{2^{*}}dx}{%
\int_{B(0,\delta )}T_{t}.(\eta \overline{u}_{t})^{2^{*}}dv_{\mathbf{g}_{t}}}%
\,\,. 
\]
En faisant de m\^{e}me pour le premier quotient de (3.14), en utilisant le
lemme plus haut, la continuit\'{e} de $f$ en $x_{0}$, et en notant $\varepsilon _{\delta}=
\stackunder{t\rightarrow 1}{\overline{\lim }}\varepsilon _{\delta,t }$ on obtient: 
\[
\stackunder{t\rightarrow 1}{\overline{\lim }}\frac{A_{t}^{1}}{%
\int_{B(0,\delta )}\overline{u}_{t}^{2}dv_{\mathbf{g}_{t}}}\leq 
\]
$\frac{n-2}{n}(1+\varepsilon _{\delta })\stackunder{t\rightarrow 1}{\overline{\lim }}
\frac{\frac{1}{2}\partial _{kl}\overline{f}_{t}(x_{0})\int_{B(0,\delta
)}(x^{k}-x_{0}^{k}(t))(x^{l}-x_{0}^{l}(t))(\eta \overline{u}_{t})^{2^{*}}dv_{%
\mathbf{g}_{t}}+C\int_{B(0,\delta )}\left| x-x_{0}(t)\right| ^{3}(\eta 
\overline{u}_{t})^{2^{*}}dv_{\mathbf{g}_{t}}}{\int_{B(0,\delta )}\overline{u}%
_{t}^{2}dv_{\mathbf{g}_{t}}}$
\\

o\`{u} l'on \'{e}crit $\partial _{kl}\overline{f}_{t}(x_{0})$ pour $\partial
_{kl}\overline{f}_{t}(x_{0}(t)).$ Consid\'{e}rons le d\'{e}veloppement 
\[
\overline{f}_{t}(x)\leq f(x_{0})+\frac{1}{2}\partial _{kl}\overline{f}%
_{t}(x_{0}(t)).(x^{k}-x_{0}^{k}(t))(x^{l}-x_{0}^{l}(t))+c\left|
x-x_{0}(t)\right| ^{3} 
\]
\textit{A priori} $c$ d\'{e}pend de $t$, mais par la r\'{e}gularit\'{e} de $%
\exp _{x_{t}}^{-1}\circ \exp _{x_{0}}$ par rapport \`{a} toutes les
variables, on peut supposer que $c$ est ind\'{e}pendant de $t$. De plus: 
\begin{eqnarray*}
c\left| x-x_{0}(t)\right| ^{3} & &\leq c^{\prime }\left| x-x_{0}(t)\right|
\sum_{k}(x^{k}-x_{0}^{k}(t))^{2} \\ 
& &\leq 2\delta c^{\prime }\sum_{k}(x^{k}-x_{0}^{k}(t))^{2}
\end{eqnarray*}
o\`{u} l'on rapelle que $\delta $ est le rayon de la boule sur laquelle on
int\`{e}gre. On peut alors \'{e}crire: 
\[
\overline{f}_{t}(x)\leq f(x_{0})+(\frac{1}{2}\partial _{kl}\overline{f}%
_{t}(x_{0}(t))+\delta C_{kl})(x^{k}-x_{0}^{k}(t))(x^{l}-x_{0}^{l}(t)) 
\]
o\`{u} $C_{kl}=c\delta _{kl}=c$ si $k=l$ et $C_{kl}=0$ si $k\neq l$ ($\delta
_{kl}$ est ici le symbole de Kr\"{o}necker) est ind\'{e}pendant de $t$.

Introduisons une notation de plus: 
\[
D_{kl}(t,\delta )=\frac{1}{2}\partial _{kl}\overline{f}_{t}(x_{0}(t))+\delta
C_{kl}\,\,. 
\]
Alors:

1/: $\stackunder{\delta \rightarrow 0}{\lim }\stackunder{t\rightarrow 1}{%
\lim }D_{kl}(t,\delta )=\frac{1}{2}\partial _{kl}\overline{f}_{1}(x_{0}(1))$
o\`{u} $\overline{f}_{1}=f\circ \exp _{x_{0}}^{-1}$ et $x_{0}(1)=0=\exp
_{x_{0}}^{-1}(x_{0})$.

2/: pour tout $\delta $ assez petit et pour $t$ proche de 1 $D_{kl}(t,\delta
)$ est encore (d\'{e}finie) n\'{e}gative.

$D_{kl}(t,\delta )$ est le hessien de $f$ en $x_{0}(t)$ perturb\'{e} sur sa
diagonale par les termes d'ordre 3. C'est pour obtenir le point 2/ que nous
avons besoin de l'hypoth\`{e}se de non-d\'{e}g\'{e}n\'{e}rescence du hessien
de $f$ en ses points de maximum. Ainsi 
\[
\frac{1}{2}\partial _{kl}\overline{f}_{t}(x_{0})\int_{B(0,\delta
)}(x^{k}-x_{0}^{k}(t))(x^{l}-x_{0}^{l}(t))(\eta \overline{u}_{t})^{2^{*}}dv_{%
\mathbf{g}_{t}}+C\int_{B(0,\delta )}\left| x-x_{0}(t)\right| ^{3}(\eta 
\overline{u}_{t})^{2^{*}}dv_{\mathbf{g}_{t}}\leq 
\]
\[
D_{kl}(t,\delta )\int_{B(0,\delta
)}(x^{k}-x_{0}^{k}(t))(x^{l}-x_{0}^{l}(t))(\eta \overline{u}_{t})^{2^{*}}dv_{%
\mathbf{g}_{t}}\,\,. 
\]
Si 
\[
\{F_{t}^{\prime }\}=D_{kl}(t,\delta )\int_{B(0,\delta
)}(x^{k}-x_{0}^{k}(t))(x^{l}-x_{0}^{l}(t))(\eta \overline{u}_{t})^{2^{*}}dv_{%
\mathbf{g}_{t}} 
\]
$\{F_{t}^{\prime }\}$ est analogue \`{a} $\{F_{t}\}$ mais c'est la mesure $%
dv_{\mathbf{g}_{t}}$ qui intervient au lieu de $dx$ et on utilise $%
D_{kl}(t,\delta ).$ Alors 
\[
\stackunder{t\rightarrow 1}{\overline{\lim }}\frac{A_{t}^{1}}{%
\int_{B(0,\delta )}v_{t}^{2}dv_{\mathbf{g}_{t}}}\leq \frac{n-2}{n}%
\stackunder{t\rightarrow 1}{\overline{\lim }}\frac{D_{kl}(t,\delta
)\int_{B(0,\delta )}(x^{k}-x_{0}^{k}(t))(x^{l}-x_{0}^{l}(t))(\eta \overline{u%
}_{t})^{2^{*}}dv_{\mathbf{g}_{t}}}{\int_{B(0,\delta )}\overline{u}%
_{t}^{2}dv_{\mathbf{g}_{t}}}(1+\varepsilon _{\delta }) 
\]

Dans le d\'{e}veloppement de $D_{kl}(t,\delta
)(x^{k}-x_{0}^{k}(t))(x^{l}-x_{0}^{l}(t))$, ce qui nous int\'{e}resse, c'est
le premier terme, i.e $D_{kl}(t,\delta )x^{k}x^{l}$ (voir l'obtention de $%
S_{g}(x_{0})$ dans $A_{t}^{2}$), et l'on esp\`{e}re que les autres vont
\^{e}tre n\'{e}gligeables. Le principe va \^{e}tre de r\'{e}arranger $%
\{F_{t}^{\prime }\}$ et d'utiliser le fait que $D_{kl}(t,\delta )$ est une
forme bilin\'{e}aire n\'{e}gative: 
\begin{eqnarray*}
\{F_{t}^{\prime }\}=&&D_{kl}(t,\delta )\int_{B(0,\delta
)}x^{k}x^{l}(\eta \overline{u}_{t})^{2^{*}}dv_{\mathbf{g}_{t}}+D_{kl}(t,%
\delta )x_{0}^{k}(t)x_{0}^{l}(t)\int_{B(0,\delta )}(\eta \overline{u}%
_{t})^{2^{*}}dv_{\mathbf{g}_{t}} \\ 
& &-D_{kl}(t,\delta )\int_{B(0,\delta
)}(x^{k}x_{0}^{l}(t)+x^{l}x_{0}^{k}(t))(\eta \overline{u}_{t})^{2^{*}}dv_{%
\mathbf{g}_{t}}\,\,.
\end{eqnarray*}
On refactorise les deux derniers termes (en enlevant quelques $\delta $ et $%
t $ et en sous-entendant $dv_{\mathbf{g}_{t}}\,$pour \'{e}claircir):

\[
D_{kl}.x_{0}^{k}x_{0}^{l}\int_{B(0,\delta )}(\eta \overline{u}%
_{t})^{2^{*}}dv_{\mathbf{g}_{t}}-D_{kl}\int_{B(0,\delta
)}(x^{k}x_{0}^{l}+x^{l}x_{0}^{k})(\eta \overline{u}_{t})^{2^{*}}dv_{\mathbf{g%
}_{t}}= 
\]

\[
D_{kl}[x_{0}^{k}x_{0}^{l}\int_{B(0,\delta )}(\eta \overline{u}%
_{t})^{2^{*}}-x_{0}^{l}\int_{B(0,\delta )}x^{k}(\eta \overline{u}%
_{t})^{2^{*}}-x_{0}^{k}\int_{B(0,\delta )}x^{l}(\eta \overline{u}%
_{t})^{2^{*}}]= 
\]
\[
D_{kl}[x_{0}^{k}(\int_{B(0,\delta )}(\eta \overline{u}_{t})^{2^{*}})^{\frac{1%
}{2}}.x_{0}^{l}(\int_{B(0,\delta )}(\eta \overline{u}_{t})^{2^{*}})^{\frac{1%
}{2}} 
\]
\[
-x_{0}^{l}(\int_{B(0,\delta )}(\eta \overline{u}_{t})^{2^{*}})^{\frac{1}{2}}%
\frac{\int_{B(0,\delta )}x^{k}(\eta \overline{u}_{t})^{2^{*}}}{%
(\int_{B(0,\delta )}(\eta \overline{u}_{t})^{2^{*}})^{\frac{1}{2}}}%
-x_{0}^{k}(\int_{B(0,\delta )}(\eta \overline{u}_{t})^{2^{*}})^{\frac{1}{2}}%
\frac{\int_{B(0,\delta )}x^{l}(\eta \overline{u}_{t})^{2^{*}}}{%
(\int_{B(0,\delta )}(\eta \overline{u}_{t})^{2^{*}})^{\frac{1}{2}}}]\,\,. 
\]
Ainsi, en notant (d\'{e}sol\'{e}): 
\[
\varepsilon ^{k}(t)=\int_{B(0,\delta )}x^{k}(\eta \overline{u}%
_{t})^{2^{*}}dv_{\mathbf{g}_{t}} 
\]
\[
z_{t}=(\int_{B(0,\delta )}(\eta \overline{u}_{t})^{2^{*}}dv_{\mathbf{g}%
_{t}})^{\frac{1}{2}} 
\]
l'expression ci-dessus devient:
\begin{center}
\[
D_{kl}.x_{0}^{k}x_{0}^{l}\int_{B(0,\delta )}(\eta \overline{u}%
_{t})^{2^{*}}dv_{\mathbf{g}_{t}}-D_{kl}\int_{B(0,\delta
)}(x^{k}x_{0}^{l}+x^{l}x_{0}^{k})(\eta \overline{u}_{t})^{2^{*}}dv_{\mathbf{g%
}_{t}}= 
\]
\[
=D_{kl}[x_{0}^{k}(t).z_{t}.x_{0}^{l}(t).z_{t}-x_{0}^{l}(t).z_{t}.\frac{%
\varepsilon ^{k}(t)}{z_{t}}-x_{0}^{k}(t).z_{t}.\frac{\varepsilon ^{l}(t)}{%
z_{t}}] 
\]
\end{center}
\[
=D_{kl}[(x_{0}^{k}(t).z_{t}-\frac{\varepsilon ^{k}(t)}{z_{t}}%
)(x_{0}^{l}(t).z_{t}-\frac{\varepsilon ^{l}(t)}{z_{t}})-\frac{\varepsilon
^{k}(t)\varepsilon ^{l}(t)}{z_{t}^{2}}] 
\]
Par cette m\'{e}thode de refactorisation du Hessien, on a ainsi obtenu: 
\[
\frac{1}{2}\partial _{kl}\overline{f}_{t}(x_{0})\int_{B(0,\delta
)}(x^{k}-x_{0}^{k}(t))(x^{l}-x_{0}^{l}(t))(\eta \overline{u}_{t})^{2^{*}}dv_{%
\mathbf{g}_{t}}+C\int_{B(0,\delta )}\left| x-x_{0}(t)\right| ^{3}(\eta 
\overline{u}_{t})^{2^{*}}dv_{\mathbf{g}_{t}}\leq 
\]
\begin{center}
\[
D_{kl}(t,\delta )\int_{B(0,\delta )}x^{k}x^{l}(\eta \overline{u}%
_{t})^{2^{*}}dv_{\mathbf{g}_{t}}+D_{kl}(t,\delta )(x_{0}^{k}(t).z_{t}-\frac{%
\varepsilon ^{k}(t)}{z_{t}})(x_{0}^{l}(t).z_{t}-\frac{\varepsilon ^{l}(t)}{%
z_{t}})-D_{kl}(t,\delta )\frac{\varepsilon ^{k}(t)\varepsilon ^{l}(t)}{%
z_{t}^{2}} 
\]
\end{center}
\[
\leq D_{kl}(t,\delta )\int_{B(0,\delta )}x^{k}x^{l}(\eta \overline{u}%
_{t})^{2^{*}}dv_{\mathbf{g}_{t}}-D_{kl}(t,\delta )\frac{\varepsilon
^{k}(t)\varepsilon ^{l}(t)}{z_{t}^{2}} 
\]
car, et c'est l\`{a} le point fondamental : 
\[
D_{kl}(t,\delta )\omega ^{k}\omega ^{l}\leq 0\text{ \thinspace }\forall
\omega =(\omega ^{1},...,\omega ^{n}) 
\]
ce qui permet de supprimer de l'in\'{e}galit\'{e} le terme $$D_{kl}(t,\delta
)(x_{0}^{k}(t).z_{t}-\frac{\varepsilon ^{k}(t)}{z_{t}})(x_{0}^{l}(t).z_{t}-%
\frac{\varepsilon ^{l}(t)}{z_{t}})$$ C'est ce terme en fait qui donne
l'estim\'{e}e $\frac{d_{\mathbf{g}}(x_{t},x_{0})}{\mu _{t}}\leq C$ (voir la
partie suivante).

On a donc obtenu:
\[
\stackunder{t\rightarrow 1}{\overline{\lim }}\frac{A_{t}^{1}}{%
\int_{B(0,\delta )}\overline{u}_{t}^{2}dv_{\mathbf{g}_{t}}}\leq \frac{n-2}{n}%
\stackunder{t\rightarrow 1}{\overline{\lim }}\frac{D_{kl}(t,\delta
)\int_{B(0,\delta )}x^{k}x^{l}(\eta \overline{u}_{t})^{2^{*}}dv_{\mathbf{g}%
_{t}}-D_{kl}(t,\delta )\frac{\varepsilon ^{k}(t)\varepsilon ^{l}(t)}{%
z_{t}^{2}}}{\int_{B(0,\delta )}\overline{u}_{t}^{2}dv_{\mathbf{g}_{t}}}%
(1+\varepsilon _{\delta }) 
\]
Maintenant, comme pour $A_{t}^{2}$, on \'{e}crit:
\begin{eqnarray*}
\stackunder{t\rightarrow 1}{\overline{\lim }}\frac{\int_{B(0,\delta
)}x^{k}x^{l}(\eta \overline{u}_{t})^{2^{*}}dv_{\mathbf{g}_{t}}}{%
\int_{B(0,\delta )}\overline{u}_{t}^{2}dv_{\mathbf{g}_{t}}} &=&f(x_{0})^{%
\frac{-2}{n}}K(n,2)^{2}\frac{n-4}{4(n-1)}\text{ si }k=l \\
&=&0\text{ si }k\neq l
\end{eqnarray*}
d'o\`{u} 
\[
\frac{1}{K(n,2)^{2}f(x_{0})^{\frac{n-2}{n}}}\frac{n-2}{n}\stackunder{%
t\rightarrow 1}{\overline{\lim }}\frac{D_{kl}(t,\delta )\int_{B(0,\delta
)}x^{k}x^{l}(\eta \overline{u}_{t})^{2^{*}}dv_{\mathbf{g}_{t}}}{%
\int_{B(0,\delta )}\overline{u}_{t}^{2}dv_{\mathbf{g}_{t}}}= 
\]
\[
=\frac{1}{f(x_{0})}\frac{(n-2)(n-4)}{4(n-1)}\sum_{l}(\frac{1}{2}\partial
_{ll}\overline{f}_{1}(0)+c_{ll}\delta ) 
\]
\[
=-\frac{(n-2)(n-4)}{8(n-1)}\frac{\bigtriangleup _{\mathbf{g}}f(x_{0})}{%
f(x_{0})}+\varepsilon _{\delta } 
\]
car $\bigtriangleup _{\mathbf{g}}f(x_{0})=-\sum_{l}\partial _{ll}\overline{f}%
_{1}(0)$ dans la carte exponentielle en $x_{0}$.

Enfin, montrons que le terme r\'{e}siduel est n\'{e}gligeable. 
\[
\left| \varepsilon ^{k}(t)\varepsilon ^{l}(t)\right| \leq \frac{1}{2}%
(\varepsilon ^{k}(t)^{2}+\varepsilon ^{l}(t)^{2}) 
\]
or 
\[
\varepsilon ^{k}(t)^{2}=(\int_{B(0,\delta )}x^{k}(\eta \overline{u}%
_{t})^{2^{*}}dv_{\mathbf{g}_{t}})^{2} 
\]
\[
=(\int_{B(0,R\mu _{t})}x^{k}(\eta \overline{u}_{t})^{2^{*}}dv_{\mathbf{g}%
_{t}}+\int_{B(0,\delta )\backslash B(0,R\mu _{t})}x^{k}(\eta \overline{u}%
_{t})^{2^{*}}dv_{\mathbf{g}_{t}})^{2} 
\]
\[
\leq 2(\int_{B(0,R\mu _{t})}x^{k}(\eta \overline{u}_{t})^{2^{*}}dv_{\mathbf{g%
}_{t}})^{2}+2(\int_{B(0,\delta )\backslash B(0,R\mu _{t})}x^{k}(\eta 
\overline{u}_{t})^{2^{*}}dv_{\mathbf{g}_{t}})^{2} 
\]
Les formules du changement d'\'{e}chelle nous donnent pour $R$ fix\'{e}: 
\[
\frac{(\int_{B(0,R\mu _{t})}x^{k}(\eta \overline{u}_{t})^{2^{*}}dv_{\mathbf{g%
}_{t}})^{2}}{\int_{B(0,\delta )}\overline{u}_{t}^{2}dv_{\mathbf{g}_{t}}}\leq 
\frac{(\mu _{t}\int_{B(0,R)}x^{k}\widetilde{u}_{t}^{2^{*}}dv_{\widetilde{%
\mathbf{g}}_{t}})^{2}}{\mu _{t}^{2}\int_{B(0,R)}\widetilde{u}_{t}^{2}dv_{%
\widetilde{\mathbf{g}}_{t}}.}\stackunder{t\rightarrow 1}{\rightarrow }\frac{%
(\int_{B(0,R)}x^{k}\widetilde{u}^{2^{*}}dx)^{2}}{\int_{B(0,R)}\widetilde{u}%
^{2}dx}=0 
\]
car $\widetilde{u}$ est radiale.

Enfin avec les estim\'{e}es ponctuelles: $\,d_{\mathbf{g}}(x,x_{t})^{\frac{%
n-2}{2}}u_{t}(x)\leq \varepsilon $\thinspace \thinspace si $d_{\mathbf{g}%
}(x,x_{t})\geq R\mu _{t}$,$\,$et d'apr\`{e}s l'in\'{e}galit\'{e} de
H\"{o}lder: 
\begin{eqnarray*}
(\int_{B(0,\delta )\backslash B(0,R\mu _{t})}x^{k}(\eta \overline{u}%
_{t})^{2^{*}}dv_{\mathbf{g}_{t}})^{2} && \leq \varepsilon
_{R}^{2}(\int_{B(0,\delta )\backslash B(0,R\mu _{t})}\overline{u}_{t}^{2%
\frac{n-1}{n-2}}dv_{\mathbf{g}_{t}})^{2} \\ 
&& \leq \varepsilon _{R}^{2}(\int_{B(0,\delta )\backslash B(0,R\mu _{t})}%
\overline{u}_{t}^{2}dv_{\mathbf{g}_{t}})(\int_{B(0,\delta )\backslash
B(0,R\mu _{t})}\overline{u}_{t}^{\frac{2n}{n-2}}dv_{\mathbf{g}_{t}})
\end{eqnarray*}
donc 
\[
\frac{(\int_{B(0,\delta )\backslash B(0,R\mu _{t})}x^{k}(\eta \overline{u}%
_{t})^{2^{*}}dv_{\mathbf{g}_{t}})^{2}}{\int_{B(0,\delta )}\overline{u}%
_{t}^{2}dv_{\mathbf{g}_{t}}}\leq \varepsilon _{R}^{2}(\int_{B(0,\delta
)\backslash B(0,R\mu _{t})}\overline{u}_{t}^{\frac{2n}{n-2}}dv_{\mathbf{g}%
_{t}})\leq c\varepsilon _{R}^{2} 
\]
o\`{u} $\varepsilon _{R}\rightarrow 0$ quand $R\rightarrow \infty $. En
remarquant que puisque $x_{0}$ est un point de concentration: 
\[
z_{t}^{2}=\int_{B(0,\delta )}(\eta \overline{u}_{t})^{2^{*}}dv_{\mathbf{g}%
_{t}}\geq \int_{B(x_{0},\delta /4)}u_{t}^{2^{*}}dv_{\mathbf{g}}\geq c>0\text{
} 
\]
on obtient finalement 
\[
\frac{\left| \varepsilon ^{k}(t)\varepsilon ^{l}(t)\right| }{%
z_{t}^{2}\int_{B(0,\delta )}\overline{u}_{t}^{2}dv_{\mathbf{g}_{t}}}%
\stackunder{t\rightarrow 1}{\rightarrow }0 
\]

La conclusion de tout cela est que si on divise par $\int_{B(0,\delta )}%
\overline{u}_{t}^{2}dv_{\mathbf{g}_{t}}$ la relation obtenue \`{a} partir de
l'in\'{e}galit\'{e} de Sobolev: 
\begin{eqnarray*}
\int_{B(0,\delta )}\overline{h}_{t}(\eta \overline{u}_{t})^{2}dx\leq && \frac{%
1}{K(n,2){{}^{2}}(\stackunder{M}{Sup}f)^{\frac{n-2}{n}}}\int_{B(0,\delta )}%
\overline{f}_{t}\eta ^{2}\overline{u}_{t}^{2^{*}}dx \\ 
&& -\frac{1}{K(n,2){{}^{2}}}(\int_{B(0,\delta )}(\eta \overline{u}%
_{t})^{2^{*}}dx)^{\frac{2}{2^{*}}} \\ 
& &+C.\delta ^{-2}\int_{B(0,\delta )\backslash B(0,\delta /2)}\overline{u}%
_{t}^{2}dx+B_{t}+C_{t}
\end{eqnarray*}
et que l'on fait tendre $t$ vers 1, on obtient: 
\[
h(x_{0})+\varepsilon _{\delta }\leq \frac{n-2}{4(n-1)}S_{\mathbf{g}}(x_{0})-%
\frac{(n-2)(n-4)}{8(n-1)}\frac{\bigtriangleup _{\mathbf{g}}f(x_{0})}{f(x_{0})%
}+\varepsilon _{\delta }\,\,. 
\]
On fait alors tendre $\delta $ vers 0: 
\[
h(x_{0})\leq \frac{n-2}{4(n-1)}S_{\mathbf{g}}(x_{0})-\frac{(n-2)(n-4)}{8(n-1)%
}\frac{\bigtriangleup _{\mathbf{g}}f(x_{0})}{f(x_{0})}\,\,. 
\]
Or on a suppos\'{e} au contraire que: 
\[
h(x_{0})>\frac{n-2}{4(n-1)}S_{\mathbf{g}}(x_{0})-\frac{(n-2)(n-4)}{8(n-1)}%
\frac{\bigtriangleup _{\mathbf{g}}f(x_{0})}{f(x_{0})} 
\]
si $x_{0}$ \'{e}tait un point de maximum de $f$. Donc il ne peut y avoir
ph\'{e}nom\`{e}ne de concentration, et donc la limite $u$ des $u_{t}$ est
strictement positive et c'est une fonction extr\'{e}male. La fonction
faiblement critique $h$ poss\`{e}de une solution minimisante, elle est donc
critique.

\textit{Fin de la d\'{e}monstration.}

\section{Ph\'{e}nom\`{e}nes de concentration et seconde in\'{e}galit\'{e}
fondamentale.}

Nous revenons ici sur l'in\'{e}galit\'{e} 
\[
\frac{d_{\mathbf{g}}(x_{t},x_{0})}{\mu _{t}}\leq C 
\]
en nous pla\c{c}ant dans le cadre de l'\'{e}tude g\'{e}n\'{e}rale d'une
suite pr\'{e}sentant un ph\'{e}nom\`{e}ne de concentration, dans le but de
montrer cette estim\'{e}e. Les donn\'{e}es sont les suivantes:

On consid\`{e}re sur une vari\'{e}t\'{e} riemannienne compacte $(M,%
\mathbf{g)}$ de dimension $n\geq 5$ une suite $(u_{t})$ de solutions $%
C^{2,\alpha }$ de l'\'{e}quation 
\[
\bigtriangleup _{\mathbf{g}}u_{t}+h_{t}u_{t}=\lambda _{t}fu_{t}^{\frac{n+2}{%
n-2}}\text{ avec }\int_{M}fu_{t}^{\frac{2n}{n-2}}dv_{\mathbf{g}}=1 
\]
avec $f$ une fonction dont le maximum est strictement positif, et dont le
Hessien est non d\'{e}g\'{e}n\'{e}r\'{e} aux points de maximum. On suppose
de plus que 
\[
h_{t}\rightarrow h\,\,\,dans\,\,\,C^{0,\alpha } 
\]
o\`{u} $h$ est telle que $\bigtriangleup _{\mathbf{g}}+h$ est coercif. La
suite $(u_{t})$ est born\'{e}e dans $H_{1}^{2}$, donc (\`{a} extraction
pr\`{e}s) $u_{t}\rightharpoondown u$ faiblement dans $H_{1}^{2}$, et on
suppose que $u\equiv 0$, donc que la suite d\'{e}veloppe un
ph\'{e}nom\`{e}ne de concentration. Nous faisons une hypoth\`{e}se dite
``d'\'{e}nergie minimale'': 
\[
\lambda _{t}\leq \frac{1}{K(n,2)^{2}(\stackunder{M}{Sup}f)^{\frac{n-2}{n}}} 
\]
Alors, si l'on reprend ce que nous avons fait en 3.2/, nous avons (quitte
\`{a} extraire des sous-suites):

a/: Il existe un unique point de concentration $x_{0}$ et c'est un point
o\`{u} $f$ est maximum sur $M.$%
\[
\forall \delta >0\,:\,\stackunder{t\rightarrow 1}{\lim \sup }%
\int_{B(x_{0},\delta )}fu_{t}^{2^{*}}dv_{\mathbf{g}}=1 
\]

b/: $u_{t}\rightarrow 0$ dans $C_{loc}^{0}(M-\{x_{0}\})$

c/: Estim\'{e}es ponctuelles faibles:

Reprenons les notations du changement d'\'{e}chelle: on consid\`{e}re une
suite de points ($x_{t})$ tels que 
\[
m_{t}=\stackunder{M}{Max}\,u_{t}=u_{t}(x_{t}):=\mu _{t}^{-\frac{n-2}{2}}. 
\]
Alors d'apr\`{e}s ce qui pr\'{e}c\`{e}de $x_{t}\rightarrow x_{0}$ et $\mu
_{t}\rightarrow 0$. Rappelons que $\overline{u}_{t},\overline{f}_{t},%
\overline{h}_{t},\mathbf{g}_{t}$ d\'{e}signent les fonctions et la
m\'{e}trique ``lues'' dans la carte $\exp _{x_{t}}^{-1}$, et $\,\,\widetilde{%
u}_{t}\,,\,\widetilde{h}_{t}\,,\,\widetilde{f}_{t},\widetilde{\mathbf{g}}%
_{t} $ d\'{e}signent les fonctions et la m\'{e}trique ``lues'' après
blow-up. 
\[
\widetilde{u}_{t}\rightarrow \widetilde{u}=(1+\frac{f(x_{0})^{\frac{2}{n}}}{%
K(n,2)^{2}}\frac{\left| x\right| ^{2}}{n(n-2)})^{-\frac{n-2}{2}}\text{ dans }%
C_{loc}^{2}(\Bbb{R}^{n}) 
\]

\[
\forall R>0\text{ on a:}\stackunder{t\rightarrow 1}{\lim }\int_{B(x_{t},R\mu
_{t})}fu_{t}^{2^{*}}dv_{\mathbf{g}}=1-\varepsilon _{R}\text{ o\`{u} }%
\varepsilon _{R}\stackunder{R\rightarrow +\infty }{\rightarrow }0 
\]

\[
\exists C>0\text{ tel que }\forall x\in M:\,d_{\mathbf{g}}(x,x_{t})^{\frac{%
n-2}{2}}u_{t}(x)\leq C. 
\]

\[
\forall \varepsilon >0,\exists R>0\text{ tel que }\forall t,\,\forall x\in
M:\,d_{\mathbf{g}}(x,x_{t})\geq R\mu _{t}\,\Rightarrow \,\,d_{\mathbf{g}%
}(x,x_{t})^{\frac{n-2}{2}}u_{t}(x)\leq \varepsilon . 
\]

d/: Concentration $L^{2}$

\[
\forall \delta >0\,:\,\stackunder{t\rightarrow 1}{\lim }\frac{%
\int_{B(x_{0},\delta )}u_{t}^{2}dv_{\mathbf{g}}}{\int_{M}u_{t}^{2}dv_{%
\mathbf{g}}}=1 
\]

e/: Estim\'{e}es ponctuelles fortes et concentration $L^{p}$ forte

\[
\exists C>0\text{ tel que }\forall x\in M:\,d_{\mathbf{g}}(x,x_{t})^{n-2}\mu
_{t}^{-\frac{n-2}{2}}u_{t}(x)\leq C. 
\]

\textbf{\ }

$\forall R>0$ , $\forall \delta >0$ et $\forall p>\frac{n}{n-2}$ o\`{u} $%
n=\dim M$ on a : 
\[
\stackunder{t\rightarrow 1}{\lim }\frac{\int_{B(x_{t},R\mu
_{t})}u_{t}^{p}dv_{\mathbf{g}}}{\int_{B(x_{t},\delta )}u_{t}^{p}dv_{\mathbf{g%
}}}=1-\varepsilon _{R}\text{ o\`{u} }\varepsilon _{R}\stackunder{%
R\rightarrow +\infty }{\rightarrow }0\text{ } 
\]
\newline

Signalons que l'on obtient comme cons\'{e}quence relativement rapide des
estim\'{e}es fortes la convergence suivante: 
\[
\mu _{t}^{-\frac{n-2}{2}}u_{t}(x)=u_{t}(x_{t})u_{t}(x)\rightarrow (\frac{4}{%
\omega _{n}^{\frac{2}{n}}})\frac{(n-2)}{\stackunder{M}{Sup}f}\omega
_{n-1}G_{h}(x_{0},x)\text{ dans }C_{loc}^{2}(M\backslash \{x_{0}\})\,\, 
\]
o\`{u} $G_{h}(x_{0},x)$ est la valeur de la fonction de Green en $(x_{0},x)$%
. La d\'{e}monstration reprend le sch\'{e}ma de la deuxi\`{e}me partie de
celle des estim\'{e}es fortes; la pr\'{e}sence de $f$ ne changeant
strictement rien, nous renvoyons \`{a} l'article de O. Druet et F. Robert.
Nous ne nous servirons de cette relation qu'au chapitre 7.

\subsubsection{Seconde in\'{e}galit\'{e} fondamentale:}

Le point nouveau obtenu par la pr\'{e}sence d'une fonction $f$ non constante
et v\'{e}rifiant (\textbf{H}$_{f}$) au second membre de nos \'{e}quations
est alors exprim\'{e} par le r\'{e}sultat suivant, qui lie la vitesse de
convergence du point de sup $x_{t}$ vers le point de concentration $x_{0}$
\`{a} la valeur du maximum des $u_{t}$ represent\'{e}e par $\mu _{t}$ : 
\[
\]

\textbf{Th\'{e}or\`{e}me 9:}

\textit{Il existe une constante }$C>0$\textit{\ telle que }$\forall t:$%
\[
\frac{d_{\mathbf{g}}(x_{t},x_{0})}{\mu _{t}}\leq C 
\]

\medskip

En effet, en reprenant les calculs de la d\'{e}monstration
pr\'{e}c\'{e}dente, on obtient: 
\begin{eqnarray}
h(x_{0}) &\leq &\frac{n-2}{4(n-1)}S_{\mathbf{g}}(x_{0})-\frac{(n-2)(n-4)}{%
8(n-1)}\frac{\bigtriangleup _{\mathbf{g}}f(x_{0})}{f(x_{0})}+\varepsilon
_{\delta }  \\
&&+\stackunder{t\rightarrow 1}{\overline{\lim }}\frac{n-2}{n}\frac{%
D_{kl}(t,\delta )(x_{0}^{k}(t).z_{t}-\frac{\varepsilon ^{k}(t)}{z_{t}}%
)(x_{0}^{l}(t).z_{t}-\frac{\varepsilon ^{l}(t)}{z_{t}})}{\int_{B(0,\delta )}%
\overline{u}_{t}^{2}dv_{\mathbf{g}_{t}}}  \nonumber
\end{eqnarray}
o\`{u} $D_{kl}(t,\delta )$ est d\'{e}finie n\'{e}gative pour $t$ voisin de 1
et pour tout $\delta $ assez petit et o\`{u} l'on rappelle que $%
x_{0}(t)=\exp _{x_{t}}^{-1}(x_{0})=(x_{0}^{1}(t),...,x_{0}^{n}(t))$. Il
existe donc $\lambda >0$ tel que $\forall \omega \in \Bbb{R}^{n}:$%
\[
D_{kl}(t,\delta )\omega ^{k}\omega ^{l}\leq -\lambda \sum_{k}\left| \omega
^{k}\right| ^{2} 
\]
et donc
\[
D_{kl}(t,\delta )\frac{(x_{0}^{k}(t).z_{t}-\frac{\varepsilon ^{k}(t)}{z_{t}}%
)(x_{0}^{l}(t).z_{t}-\frac{\varepsilon ^{l}(t)}{z_{t}})}{(\int_{B(0,\delta )}%
\overline{u}_{t}^{2}dv_{\mathbf{g}_{t}})^{\frac{1}{2}}(\int_{B(0,\delta )}%
\overline{u}_{t}^{2}dv_{\mathbf{g}_{t}})^{\frac{1}{2}}}\leq 
\]
\[
-\lambda \sum_{k}\left| \frac{x_{0}^{k}(t).z_{t}}{(\int_{B(0,\delta )}%
\overline{u}_{t}^{2}dv_{\mathbf{g}_{t}})^{\frac{1}{2}}}-\frac{\varepsilon
^{k}(t)}{z_{t}(\int_{B(0,\delta )}\overline{u}_{t}^{2}dv_{\mathbf{g}_{t}})^{%
\frac{1}{2}}}\right| ^{2}\,\,. 
\]
Or on a vu que: 
\[
\frac{\varepsilon ^{k}(t)^{2}}{z_{t}^{2}\int_{B(0,\delta )}\overline{u}%
_{t}^{2}dv_{\mathbf{g}_{t}}}\stackunder{t\rightarrow 1}{\rightarrow }0 
\]
et d'autre part $z_{t}=(\int_{B(0,\delta )}\overline{u}_{t}^{2^{*}}dv_{%
\mathbf{g}_{t}})^{\frac{1}{2}}$,\thinspace donc puisque $x_{0}$ est un point
de concentration: 
\[
0<c\leq \lim \inf \,z_{t}\leq \lim \sup \,z_{t}\leq c^{\prime }<+\infty . 
\]
Donc n\'{e}cessairement, \`{a} cause de (3.15), pour tout $k$, il existe
une constante $C>0$ telle que pour $t\rightarrow 1$ : 
\[
\,\frac{x_{0}^{k}(t)}{(\int_{B(0,\delta )}\overline{u}_{t}^{2}dv_{\mathbf{g}%
_{t}})^{\frac{1}{2}}}\leq C 
\]
Maintenant 
\[
\int_{B(0,\delta )}\overline{u}_{t}^{2}dv_{\mathbf{g}_{t}}=\mu
_{t}^{2}\int_{B(0,\delta \mu _{t}^{-1})}\widetilde{u}_{t}^{2}dv_{\widetilde{%
\mathbf{g}}_{t}}. 
\]
mais les estim\'{e}es ponctuelles fortes nous donnent que 
\[
\stackunder{t\rightarrow 1}{\overline{\lim }}\int_{B(0,\delta \mu _{t}^{-1})}%
\widetilde{u}_{t}^{2}dv_{\widetilde{\mathbf{g}}_{t}}<+\infty 
\]
(elles donnent en fait $\widetilde{u}_{t}\leq C.\widetilde{u}$ sur $%
B(0,\delta \mu _{t}^{-1})$) donc 
\[
\int_{B(0,\delta )}\overline{u}_{t}^{2}dv_{\mathbf{g}_{t}}\sim C\mu _{t}^{2} 
\]
d'o\`{u} 
\[
\forall k:\,\frac{x_{0}^{k}(t)}{\mu _{t}}\leq C^{\prime } 
\]
et donc 
\[
\frac{d_{\mathbf{g}}(x_{t},x_{0})}{\mu _{t}}\leq C 
\]
CQFD.

\textit{Remarque 1:} Ind\'{e}pendamment, s'il y a concentration,
n\'{e}cessairement: 
\[
h(x_{0})\leq \frac{n-2}{4(n-1)}S_{\mathbf{g}}(x_{0})-\frac{(n-2)(n-4)}{8(n-1)%
}\frac{\bigtriangleup _{\mathbf{g}}f(x_{0})}{f(x_{0})} 
\]

\textit{Remarque 2:} Si on sait de plus qu'aux points de maximum de $f$ : 
\[
h(P)=\frac{n-2}{4(n-1)}S_{\mathbf{g}}(P)-\frac{(n-2)(n-4)}{8(n-1)}\frac{%
\bigtriangleup _{\mathbf{g}}f(P)}{f(P)} 
\]

alors on voit d'apr\`{e}s ce qui pr\'{e}c\`{e}de que plus
pr\'{e}cis\'{e}ment: 
\[
\frac{d_{\mathbf{g}}(x_{t},x_{0})}{\mu _{t}}\rightarrow 0 
\]

\textit{Remarque 3:} Intuitivement, les estim\'{e}es ponctuelles en $x_{t}$
donnent la ``forme'' des $u_{t}$, le th\'{e}or\`{e}me 9 donne la position de
leur maximum. Plus pr\'{e}cis\'{e}ment, il permet de remplacer dans les
estim\'{e}es fortes ou faibles et dans la concentration $L^{p}$ forte $x_{t}$
par $x_{0}$.

\textit{Remarque 4:} La m\'{e}thode que nous avons d\'{e}velopp\'{e}e pour
obtenir cette seconde in\'{e}galit\'{e} fondamentale semble pouvoir
s'appliquer \`{a} d'autres probl\`{e}mes similaires. Voir par exemple
l'article de Zo\'{e} Faget $\left[ 15\right] $.

\section{Remarque: Illustration de l'utilisation de la seconde
in\'{e}galit\'{e} fondamentale}

Pour montrer l'importance de l'estim\'{e}e 
\[
\frac{d_{\mathbf{g}}(x_{t},x_{0})}{\mu _{t}}\leq C 
\]
nous voulons ici indiquer rapidement comment, si l'on suppose qu'on a pu
l'obtenir ind\'{e}pendament, cette estim\'{e}e permet de conclure la
d\'{e}monstration du th\'{e}or\`{e}me central de ce chapitre. Il ne s'agit
donc pas d'une nouvelle d\'{e}monstration, mais d'une illustration de
l'importance de cette estim\'{e}e.

On commence de la m\^{e}me mani\`{e}re, et on arrive \`{a} l'\'{e}tude de
l'expression $A_{t}^{1}:$%
\[
A_{t}^{1}=(\int_{B(0,\delta )}\overline{f}_{t}\overline{u}_{t}^{2^{*}}dv_{%
\mathbf{g}_{t}})^{\frac{2}{n}}(\int_{B(0,\delta )}\overline{f}_{t}(\eta 
\overline{u}_{t})^{2^{*}}dx)^{\frac{n-2}{n}}\,-(Supf.\int_{B(0,\delta
)}(\eta \overline{u}_{t})^{2^{*}}dx)^{\frac{n-2}{n}}\,\,. 
\]
En introduisant le d\'{e}veloppement limit\'{e} de $f$ on obtient 
\[
\begin{tabular}{l}
$\stackunder{t\rightarrow 1}{\overline{\lim }}\frac{A_{t}^{1}}{%
\int_{B(0,\delta )}\overline{u}_{t}^{2}dv_{\mathbf{g}_{t}}}\leq $ \\ 
$\stackunder{t\rightarrow 1}{\overline{\lim }}\frac{n-2}{n}\frac{\frac{1}{2}%
\partial _{kl}\overline{f}_{t}(x_{0})\int_{B(0,\delta
)}(x^{k}-x_{0}^{k}(t))(x^{l}-x_{0}^{l}(t))(\eta \overline{u}_{t})^{2^{*}}dv_{%
\mathbf{g}_{t}}+C\int_{B(0,\delta )}\left| x-x_{0}(t)\right| ^{3}(\eta 
\overline{u}_{t})^{2^{*}}dv_{\mathbf{g}_{t}}}{\int_{B(0,\delta )}\overline{u}%
_{t}^{2}dv_{\mathbf{g}_{t}}}(1+\varepsilon _{\delta })$%
\end{tabular}
\]
o\`{u} l'on \'{e}crit $\partial _{kl}\overline{f}_{t}(x_{0})$ pour $\partial
_{kl}\overline{f}_{t}(x_{0}(t)).$ En d\'{e}veloppant 
\begin{eqnarray*}
&\frac{1}{2}\partial _{kl}\overline{f}_{t}(x_{0})\int_{B(0,\delta
)}(x^{k}-x_{0}^{k}(t))(x^{l}-x_{0}^{l}(t))(\eta \overline{u}_{t})^{2^{*}}dv_{%
\mathbf{g}_{t}}= &\\
&\frac{1}{2}\partial _{kl}\overline{f}_{t}(x_{0})%
\int_{B(0,\delta )}x^{k}x^{l}(\eta \overline{u}_{t})^{2^{*}}dv_{\mathbf{g}%
_{t}}& \\ 
& +\frac{1}{2}\partial _{kl}\overline{f}_{t}(x_{0})x_{0}^{k}(t)x_{0}^{l}(t)%
\int_{B(0,\delta )}(\eta \overline{u}_{t})^{2^{*}}dv_{\mathbf{g}_{t}} &\\ 
& +\partial _{kl}\overline{f}_{t}(x_{0})x_{0}^{l}(t)\int_{B(0,\delta
)}x^{k}(\eta \overline{u}_{t})^{2^{*}}dv_{\mathbf{g}_{t}}&
\end{eqnarray*}
soit 
\begin{eqnarray*}
\frac{1}{2}\partial _{kl}\overline{f}_{t}(x_{0})\int_{B(0,\delta
)}(x^{k}-x_{0}^{k}(t))(x^{l}-x_{0}^{l}(t))(\eta \overline{u}_{t})^{2^{*}}dv_{%
\mathbf{g}_{t}} &\leq &\frac{1}{2}\partial _{kl}\overline{f}%
_{t}(x_{0})\int_{B(0,\delta )}x^{k}x^{l}(\eta \overline{u}_{t})^{2^{*}}dv_{%
\mathbf{g}_{t}} \\
&&+\partial _{kl}\overline{f}_{t}(x_{0})x_{0}^{l}(t)\int_{B(0,\delta
)}x^{k}(\eta \overline{u}_{t})^{2^{*}}dv_{\mathbf{g}_{t}}
\end{eqnarray*}
car 
\[
\frac{1}{2}\partial _{kl}\overline{f}_{t}(x_{0})x_{0}^{k}(t)x_{0}^{l}(t)\leq
0 
\]
puisque $x_{0}$ est un point de maximum de $f$. Quand au reste: 
\[
\int_{B(0,\delta )}\left| x-x_{0}(t)\right| ^{3}(\eta \overline{u}%
_{t})^{2^{*}}dv_{\mathbf{g}_{t}}\leq C\stackunder{p+q=3}{\sum }\left|
x_{0}(t)\right| ^{q}\int_{B(0,\delta )}\left| x\right| ^{p}(\eta \overline{u}%
_{t})^{2^{*}}dv_{\mathbf{g}_{t}}\,\,. 
\]
On a d\'{e}ja vu que 
\begin{eqnarray*}
\stackunder{t\rightarrow 1}{\overline{\lim }}\frac{\int_{B(0,\delta
)}x^{k}x^{l}(\eta \overline{u}_{t})^{2^{*}}dv_{\mathbf{g}_{t}}}{%
\int_{B(0,\delta )}\overline{u}_{t}^{2}dv_{\mathbf{g}_{t}}} &=&f(x_{0})^{%
\frac{-2}{n}}K(n,2)^{2}\frac{n-4}{4(n-1)}\text{ si }k=l \\
&=&0\text{ si }k\neq l
\end{eqnarray*}
et donc 
\[
\frac{1}{K(n,2)^{2}f(x_{0})^{\frac{n-2}{n}}}\stackunder{t\rightarrow 1}{%
\overline{\lim }}\frac{n-2}{n}\frac{\frac{1}{2}\partial _{kl}\overline{f}%
_{t}(x_{0})\int_{B(0,\delta )}x^{k}x^{l}(\eta \overline{u}_{t})^{2^{*}}dv_{%
\mathbf{g}_{t}}}{\int_{B(0,\delta )}\overline{u}_{t}^{2}dv_{\mathbf{g}_{t}}}%
= 
\]
\[
=\frac{1}{f(x_{0})}\frac{(n-2)(n-4)}{4(n-1)}\sum_{l}(\frac{1}{2}\partial
_{ll}\overline{f}_{1}(0)+c_{ll}\delta ) 
\]
\[
=-\frac{(n-2)(n-4)}{8(n-1)}\frac{\bigtriangleup _{\mathbf{g}}f(x_{0})}{%
f(x_{0})}+\varepsilon _{\delta } 
\]
Par changement d'\'{e}chelle, on peut maintenant \'{e}crire 
\begin{eqnarray*}
\partial _{kl}\overline{f}_{t}(x_{0})x_{0}^{l}(t)\int_{B(0,\delta
)}x^{k}(\eta \overline{u}_{t})^{2^{*}}dv_{\mathbf{g}_{t}}= && \partial _{kl}%
\overline{f}_{t}(x_{0})x_{0}^{l}(t)\mu _{t}\int_{B(0,R)}x^{k}\widetilde{u}%
_{t}^{2^{*}}dv_{\widetilde{\mathbf{g}}_{t}} \\ 
& &+\partial _{kl}\overline{f}_{t}(x_{0})x_{0}^{l}(t)\int_{B(0,\delta
)\backslash B(0,R\mu _{t})}x^{k}(\eta \overline{u}_{t})^{2^{*}}dv_{\mathbf{g}%
_{t}}
\end{eqnarray*}
et 
\[
\int_{B(0,\delta )}\overline{u}_{t}^{2}dv_{\mathbf{g}_{t}}=\mu
_{t}^{2}\int_{B(0,\delta \mu _{t}^{-1})}\widetilde{u}_{t}^{2}dv_{\widetilde{%
\mathbf{g}}_{t}}\,\,. 
\]
Alors 
\[
\frac{x_{0}^{l}(t)\mu _{t}\int_{B(0,R)}x^{k}\widetilde{u}_{t}^{2^{*}}dv_{%
\widetilde{\mathbf{g}}_{t}}}{\int_{B(0,\delta )}\overline{u}_{t}^{2}dv_{%
\mathbf{g}_{t}}}=\frac{x_{0}^{l}(t)\int_{B(0,R)}x^{k}\widetilde{u}%
_{t}^{2^{*}}dv_{\widetilde{\mathbf{g}}_{t}}}{\mu _{t}\int_{B(0,\delta \mu
_{t}^{-1})}\widetilde{u}_{t}^{2}dv_{\widetilde{\mathbf{g}}_{t}}} 
\]
mais 
\[
\left| \frac{\int_{B(0,R)}x^{k}\widetilde{u}_{t}^{2^{*}}dv_{\widetilde{%
\mathbf{g}}_{t}}}{\int_{B(0,\delta \mu _{t}^{-1})}\widetilde{u}_{t}^{2}dv_{%
\widetilde{\mathbf{g}}_{t}}}\right| \leq \left| \frac{\int_{B(0,R)}x^{k}%
\widetilde{u}_{t}^{2^{*}}dv_{\widetilde{\mathbf{g}}_{t}}}{\int_{B(0,R)}%
\widetilde{u}_{t}^{2}dv_{\widetilde{\mathbf{g}}_{t}}}\right| \stackunder{%
t\rightarrow 1}{\rightarrow }\left| \frac{\int_{B(0,R)}x^{k}\widetilde{u}%
^{2^{*}}dx}{\int_{B(0,R)}\widetilde{u}^{2}dx}\right| =0 
\]
car $\widetilde{u}$ est radiale et $\frac{\left| x_{0}^{l}(t)\right| }{\mu
_{t}}$ est born\'{e}e par la seconde in\'{e}galit\'{e} fondamentale.
Ensuite, par l'in\'{e}galit\'{e} de H\"{o}lder et en utilisant les
estim\'{e}es ponctuelles faibles 
\begin{eqnarray*}
\left| \int_{B(0,\delta )\backslash B(0,R\mu _{t})}x^{k}(\eta \overline{u}%
_{t})^{2^{*}}dv_{\mathbf{g}_{t}}\right| && \leq \varepsilon _{R}\left|
\int_{B(0,\delta )\backslash B(0,R\mu _{t})}\overline{u}_{t}^{2\frac{n-1}{n-2%
}}dv_{\mathbf{g}_{t}}\right| \\ 
&& \leq \varepsilon _{R}(\int_{B(0,\delta )\backslash B(0,R\mu _{t})}\overline{u%
}_{t}^{2}dv_{\mathbf{g}_{t}})^{\frac{1}{2}}(\int_{B(0,\delta )\backslash
B(0,R\mu _{t})}\overline{u}_{t}^{\frac{2n}{n-2}}dv_{\mathbf{g}_{t}})^{\frac{1%
}{2}} \\ 
&&\leq c\varepsilon _{R}(\int_{B(0,\delta )\backslash B(0,R\mu _{t})}\overline{%
u}_{t}^{2}dv_{\mathbf{g}_{t}})^{\frac{1}{2}}
\end{eqnarray*}
donc, puisque $\frac{\left| x_{0}^{l}(t)\right| }{\mu _{t}}$ est born\'{e}e: 
\begin{eqnarray*}
\left| \frac{x_{0}^{l}(t)\int_{B(0,\delta )\backslash B(0,R\mu
_{t})}x^{k}(\eta \overline{u}_{t})^{2^{*}}dv_{\mathbf{g}_{t}}}{%
\int_{B(0,\delta )}\overline{u}_{t}^{2}dv_{\mathbf{g}_{t}}}\right| &\leq &%
\frac{\varepsilon _{R}.\left| x_{0}^{l}(t)\right| (\int_{B(0,\delta
)\backslash B(0,R\mu _{t})}\overline{u}_{t}^{2}dv_{\mathbf{g}_{t}})^{\frac{1%
}{2}}}{(\int_{B(0,\delta )}\overline{u}_{t}^{2}dv_{g_{t}})^{\frac{1}{2}%
}(\int_{B(0,\delta )}\overline{u}_{t}^{2}dv_{\mathbf{g}_{t}})^{\frac{1}{2}}}
\\
&\leq &\varepsilon _{R}\frac{\left| x_{0}^{l}(t)\right| }{\mu
_{t}(\int_{B(0,R)}\widetilde{u}_{t}^{2}dv_{\widetilde{\mathbf{g}}_{t}})^{%
\frac{1}{2}}}\stackunder{R\rightarrow +\infty }{\rightarrow }0
\end{eqnarray*}
Ensuite, 
\[
\frac{\left| x_{0}(t)\right| ^{3}\int_{B(0,\delta )}(\eta \overline{u}%
_{t})^{2^{*}}dv_{\mathbf{g}_{t}}}{\int_{B(0,\delta )}\overline{u}_{t}^{2}dv_{%
\mathbf{g}_{t}}}\leq c\frac{\left| x_{0}(t)\right| ^{3}}{\mu
_{t}^{2}(\int_{B(0,R)}\widetilde{u}_{t}^{2}dv_{\widetilde{\mathbf{g}}_{t}})}%
\stackunder{t\rightarrow 1}{\rightarrow }0\,\,. 
\]
car $\frac{\left| x_{0}^{l}(t)\right| }{\mu _{t}}$ est born\'{e}e et $\left|
x_{0}(t)\right| \rightarrow 0$.

Comme ci-dessus 
$$\frac{\left| x_{0}(t)\right| ^{2}\int_{B(0,\delta )}\left| x\right| (\eta 
\overline{u}_{t})^{2^{*}}dv_{\Bbb{g}_{t}}}{\int_{B(0,\delta )}\overline{u}%
_{t}^{2}dv_{\mathbf{g}_{t}}}\leq$$
$$ \frac{\left| x_{0}(t)\right|
^{2}\varepsilon _{R}(\int_{B(0,\delta )\backslash B(0,R\mu _{t})}\overline{u}%
_{t}^{2}dv_{\mathbf{g}_{t}})^{\frac{1}{2}}(\int_{B(0,\delta )\backslash
B(0,R\mu _{t})}\overline{u}_{t}^{\frac{2n}{n-2}}dv_{\mathbf{g}_{t}})^{\frac{1%
}{2}}}{\mu _{t}^{2}(\int_{B(0,R)}\widetilde{u}_{t}^{2}dv_{\widetilde{\mathbf{%
g}}_{t}})^{\frac{1}{2}}(\int_{B(0,\delta )}\overline{u}_{t}^{2}dv_{\mathbf{g}%
_{t}})^{\frac{1}{2}}}\,\stackunder{R\rightarrow +\infty }{\longrightarrow }0 $$

Puis, en utilisant les estim\'{e}es ponctuelles faibles, on obtient 
\[
\frac{\left| x_{0}(t)\right| \int_{B(0,\delta )}\left| x\right| ^{2}(\eta 
\overline{u}_{t})^{2^{*}}dv_{\mathbf{g}_{t}}}{\int_{B(0,\delta )}\overline{u}%
_{t}^{2}dv_{\mathbf{g}_{t}}}\leq C\left| x_{0}(t)\right| 
\]
et donc ce quotient tend vers 0 puisque $\left| x_{0}(t)\right| \rightarrow
0 $.

Enfin, 
\[
\frac{\int_{B(0,\delta )}\left| x\right| ^{3}(\eta \overline{u}%
_{t})^{2^{*}}dv_{\mathbf{g}_{t}}}{\int_{B(0,\delta )}\overline{u}_{t}^{2}dv_{%
\mathbf{g}_{t}}}\leq \delta \frac{\int_{B(0,\delta )}\left| x\right|
^{2}(\eta \overline{u}_{t})^{2^{*}}dv_{\mathbf{g}_{t}}}{\int_{B(0,\delta )}%
\overline{u}_{t}^{2}dv_{\mathbf{g}_{t}}}\leq C\delta \,\,. 
\]

Avec ces limites, on peut conclure comme dans la partie 3.3 que 
\[
h(x_{0})+\varepsilon _{\delta }\leq \frac{n-2}{4(n-1)}S_{\mathbf{g}}(x_{0})-%
\frac{(n-2)(n-4)}{8(n-1)}\frac{\bigtriangleup _{\mathbf{g}}f(x_{0})}{f(x_{0})%
}+\varepsilon _{\delta } 
\]
d'o\`{u} une contradiction, puis la conclusion.

Ainsi, la seconde in\'{e}galit\'{e} fondamentale des ph\'{e}nom\`{e}nes de
concentration 
\[
\frac{d_{\mathbf{g}}(x_{t},x_{0})}{\mu _{t}}\leq C 
\]
permet d'aboutir au m\^{e}me r\'{e}sultat sans utiliser la m\'{e}thode de
refactorisation du Hessien pr\'{e}sent\'{e}e dans la troisi\`{e}me partie de
ce chapitre; mais c'est gr\^{a}ce \`{a} cette m\'{e}thode que nous avons pu
montrer cette estim\'{e}e dans la partie pr\'{e}c\'{e}dente. N\'{e}anmoins,
il est important de remarquer que si l'on suppose vraie cette
in\'{e}galit\'{e} (en imaginant qu'on ait pu l'obtenir ind\'{e}pendamment),
on n'a plus besoin de supposer que le Hessien est non
d\'{e}g\'{e}n\'{e}r\'{e} au points de maximum de $f$. En revanche, on peut
trouver une suite de solutions ($u_{t})$ d'\'{e}quations de la forme 
\[
\triangle _{\mathbf{g}}u_{t}+h_{t}.u_{t}=u_{t}^{\frac{n+2}{n-2}} 
\]
o\`{u} au second membre la fonction $f$ est constante et donc a un Hessien
d\'{e}g\'{e}n\'{e}r\'{e} (!), telles que ($u_{t})$ d\'{e}veloppe un
ph\'{e}nom\`{e}ne de concentration mais ne v\'{e}rifie pas la seconde
in\'{e}galit\'{e} fondamentale. Nous en construirons un exemple au chapitre
8.

\section{Remarque finale}

Une fa\c{c}on ``intuitive '' de voir les d\'{e}tails techniques des calculs
est de consid\'{e}rer que l'on fait un d\'{e}veloppement limit\'{e} en
fonction de $\mu _{t}$ des deux membres de l'in\'{e}galit\'{e} issue de
celle de Sobolev, au sens o\`{u}, apr\`{e}s changement d'\'{e}chelle, on
\'{e}value les int\'{e}grales en fonction de ce param\`{e}tre fondamental.

Ainsi, en gros 
\[
\int \left| x\right| ^{p}\overline{u}_{t}^{2^{*}}dv_{\mathbf{g}_{t}}\sim \mu
_{t}^{p} 
\]
et 
\[
\int \overline{u}_{t}^{2}dv_{\mathbf{g}_{t}}\sim \mu _{t}^{2}\,\,. 
\]
On obtient ainsi \`{a} partir de l'in\'{e}galit\'{e} de Sobolev dans
laquelle on introduit l'\'{e}quation v\'{e}rifi\'{e}e par $\overline{u}_{t}$%
: 
\[
h(x_{0})\mu _{t}^{2}\leq (\frac{n-2}{4(n-1)}S_{\mathbf{g}}(x_{0})-\frac{%
(n-2)(n-4)}{8(n-1)}\frac{\bigtriangleup _{\mathbf{g}}f(x_{0})}{f(x_{0})})\mu
_{t}^{2}+o(\mu _{t}^{2}) 
\]
On divise alors cette relation par $\mu _{t}^{2}$ (c'est \`{a} dire qu'on
utilise ``la concentration $L^{2}$'') pour obtenir une contradiction avec
l'hypoth\`{e}se 
\[
h(x_{0})>\frac{n-2}{4(n-1)}S_{\mathbf{g}}(x_{0})-\frac{(n-2)(n-4)}{8(n-1)}%
\frac{\bigtriangleup _{\mathbf{g}}f(x_{0})}{f(x_{0})}\,\,. 
\]
C'est cette fa\c{c}on de voir qui nous a guid\'{e} dans la mise au point de
la m\'{e}thode de refactorisation du Hessien.

\chapter{Triplet Critique 1}
\pagestyle{myheadings}\markboth{\textbf{Existence de fonctions critiques.}}
{\textbf{Existence de fonctions critiques.}}
{\LARGE Existence de fonctions critiques.}

Nous montrons dans ce chapitre le r\'{e}sultat suivant:

\textbf{Th\'{e}or\`{e}me 2:}

\textit{Soit }$(M,\mathbf{g})$\textit{\ une vari\'{e}t\'{e} de dimension }$%
n\geqslant 4$\textit{\ et }$f$\textit{\ une fonction v\'{e}rifiant (H}$%
_{f}). $\textit{\ Il existe une infinit\'{e} de fonctions critiques $h$ pour }$f$%
\textit{\ et }$\mathbf{g}$\textit{, qui v\'{e}rifient }$\frac{4(n-1)}{n-2}%
h(P)>S_{\mathbf{g}}(P)-\frac{n-4}{2}\frac{\bigtriangleup _{\mathbf{g}}f(P)}{%
f(P)}$\textit{\ en tout point }$P$\textit{\ o\`{u} }$f$\textit{\ est maximum
sur }$M$\textit{. Ces fonctions critiques ont des fonctions extr\'{e}males.}

\medskip

Commen\c{c}ons par montrer ce r\'{e}sultat quand $f>0$.

L'id\'{e}e, et c'est une fa\c{c}on de voir les fonctions critiques, est de
trouver une fonction $h_{0}$ sous-critique et une fonction $h_{1}$
faiblement critique (on pourrait dire surcritique) et de passer
contin\^{u}ment par un chemin de $h_{0}$ \`{a} $h_{1}$; le th\'{e}or\`{e}me
1 nous dit que ce chemin doit ``rencontrer'' l'ensemble des fonctions
critiques.

Soit $h$ une fonction $C^{\infty \text{ }}$ faiblement critique pour $f$
telle que en tout point $P$ o\`{u} $f$ est maximum sur $M$ on ait: 
\[
\frac{4(n-1)}{n-2}h(P)>S_{\mathbf{g}}(P)-\frac{n-4}{2}\frac{\bigtriangleup _{%
\mathbf{g}}f(P)}{f(P)}. 
\]
Il suffit de prendre $h\geq B_{0}(\mathbf{g})K(n,2)^{-2}$ . Comme $f$ est
non constante, il existe une boule ouverte 
\[
B(y,2\delta )\subset M\backslash \left\{ x/f(x)=Maxf\right\} \,\,. 
\]
Soit alors $\eta $ une fonction plateau, $\eta =0\,$sur\thinspace $%
M\backslash B(y,2\delta )$ , $\eta =1\,$sur\thinspace $B(y,\delta )$, $0\leq
\eta \leq 1.$ Soit 
\[
c=\left( \int f\right) ^{-\frac{1}{2^{*}}}\,\,. 
\]
Alors $\int fc^{2^{*}}=1$ et $I_{h}(c)=\int hc{{}^{2}}$. Pour $t\in \Bbb{R}%
^{+}$ on pose 
\[
h_{t}=h-t\eta 
\]
Ainsi: $h_{t}=h$ sur $M\backslash B(y,\delta )\supset \left\{
x/f(x)=Maxf\right\} $, et 
\[
I_{h_{t}}(c)=\int h_{t}c{{}^{2}}=\int hc{{}^{2}}-c{{}^{2}}t\int_{B(y,2\delta
)}\eta \,\,. 
\]
Donc en choisissant $t$ assez grand, on obtient 
\[
I_{h_{t}}(c)<\frac{1}{K(n,2){{}^{2}}(\stackunder{M}{Sup}f)^{\frac{n-2}{n}}}%
\,\,. 
\]
Par ailleurs, sur $M\backslash B(y,\delta )$ qui contient $\left\{
x/f(x)=Maxf\right\} $, $h_{t}=h$. Donc $\forall P\in \left\{
x/f(x)=Maxf\right\} $:

\[
\frac{4(n-1)}{n-2}h_{t}(P)>S_{\mathbf{g}}(P)-\frac{n-4}{2}\frac{%
\bigtriangleup _{\mathbf{g}}f(P)}{f(P)}. 
\]
On pose 
\[
t_{0}=Inf\{t/\lambda _{h_{t},f,\mathbf{g}}<\frac{1}{K(n,2){{}^{2}}(%
\stackunder{M}{Sup}f)^{\frac{n-2}{n}}}\}\text{ }. 
\]
Alors 
\[
\lambda _{h_{t_{0}},f,\mathbf{g}}=\frac{1}{K(n,2){{}^{2}}(\stackunder{M}{Sup}%
f)^{\frac{n-2}{n}}} 
\]
et 
\[
\lambda _{h_{t},f,\mathbf{g}}<\frac{1}{K(n,2){{}^{2}}(\stackunder{M}{Sup}f)^{%
\frac{n-2}{n}}}\,\,\,si\,\,\,t>t_{0}\,\,. 
\]
De plus

\[
\frac{4(n-1)}{n-2}h_{t_{0}}(P)>S_{\mathbf{g}}(P)-\frac{n-4}{2}\frac{%
\bigtriangleup _{\mathbf{g}}f(P)}{f(P)}\,\,\,sur\,\,\,\left\{
x/f(x)=Maxf\right\} . 
\]
Enfin $h_{t}\stackunder{t\rightarrow t_{0}}{\rightarrow }h_{t_{0}}\,\,$dans$%
\,$ $C^{0,\alpha }$, et comme $h_{t_{0}}$ est faiblement critique pour $f>0$%
, $\triangle _{\mathbf{g}}+h_{t_{0}}$ est coercif. Donc d'apr\`{e}s le
th\'{e}or\`{e}me 1, $h_{t_{0}}$ est critique avec des fonctions
extr\'{e}males.

Ceci montre le th\'{e}or\`{e}me dans le cas o\`{u} $f>0$.

La difficult\'{e} quand $f$ change de signe est que pour une fonction $h$,
la condition 
\[
\lambda _{h,f,\mathbf{g}}=\frac{1}{K(n,2){{}^{2}}(\stackunder{M}{Sup}f)^{%
\frac{n-2}{n}}} 
\]
ne suffit pas \`{a} garantir que l'op\'{e}rateur $\bigtriangleup _{\mathbf{g}%
}+h$ soit coercif, ce que nous avons exig\'{e} dans la d\'{e}finition des
fonctions critiques. Rappelons que $\bigtriangleup _{\mathbf{g}}+h$ coercif
signifie qu'il existe $C>0$ tel que pour toute $u\in H_{1}^{2}$, 
\[
\int_{M}(\left| \nabla u\right| _{\mathbf{g}}^{2}+hu^{2})dv_{\mathbf{g}}\geq
C\int_{M}u^{2}dv_{\mathbf{g}}\,\,. 
\]
Une bonne fa\c{c}on d'\^{e}tre s\^{u}r que cet op\'{e}rateur est coercif,
gr\^{a}ce aux inclusions de Sobolev, est d'avoir $h>0$. Par ailleurs, nous
avons dit que si $f\equiv 1$, que $(M,\mathbf{g})$ n'est pas
conform\'{e}ment diff\'{e}omorphe \`{a} la sph\`{e}re standard, et que $S_{%
\mathbf{g}}\,$est constante, alors $B_{0}(\mathbf{g})K(n,2)^{-2}$ est une
fonction critique, c'est la plus petite fonction critique constante et en
particulier elle est strictement positive. Or quand $f$ est quelconque, il
n'est pas \'{e}vident qu'il existe de telles fonctions critiques. Le
r\'{e}sultat interm\'{e}diaire suivant a donc son int\'{e}r\^{e}t; il sera
fondamental au chapitre suivant.

\textbf{Corollaire 1:}

\textit{Soit }$(M,\mathbf{g})$\textit{\ une vari\'{e}t\'{e} de dimension }$%
n\geqslant 4$\textit{\ et }$f$\textit{\ une fonction v\'{e}rifiant (H}$%
_{f}). $ \textit{Si }$\int_{M}fdv_{\mathbf{g}}>0$, \textit{il existe une
infinit\'{e} de fonctions }$h$, \textit{strictement positives, critiques
pour }$f$ et $\mathbf{g}$\textit{, qui v\'{e}rifient }$\frac{4(n-1)}{n-2}%
h(P)>S_{\mathbf{g}}(P)-\frac{n-4}{2}\frac{\bigtriangleup _{\mathbf{g}}f(P)}{%
f(P)}$\textit{\ en tout point }$P$\textit{\ o\`{u} }$f$\textit{\ est maximum
sur }$M$\textit{. Ces fonctions critiques ont des fonctions extr\'{e}males.}

\medskip

\textit{D\'{e}monstration:}

Reprenons ce que nous avons fait pour le cas $f>0$. On part de 
\[
h=B_{0}(\mathbf{g})K(n,2)^{-2}\,\,. 
\]
En tout point $P$ o\`{u} $f$ est maximum sur $M$, d'apr\`{e}s
l'hypoth\`{e}se (H$_{f}$): 
\[
\frac{4(n-1)}{n-2}B_{0}(\mathbf{g})K(n,2)^{-2}>S_{\mathbf{g}}(P)-\frac{n-4}{2%
}\frac{\bigtriangleup _{\mathbf{g}}f(P)}{f(P)} 
\]
car 
\[
B_{0}(\mathbf{g})K(n,2)^{-2}\geq \frac{(n-2)}{4(n-1)}MaxS_{\mathbf{g}}\,\,. 
\]
Comme $f$ est non constante, il existe une fonction $\eta $ dont le support
est inclus dans $M\backslash \left\{ x/f(x)=Maxf\right\} $ et telle que $%
0\leq \eta \leq 1$. Soit 
\[
c=\left( \int_{M}fdv_{\mathbf{g}}\right) ^{-\frac{n-2}{n}} 
\]
(c'est l\`{a} qu'on utilise $\int_{M}fdv_{\mathbf{g}}>0$), on a $\int
fc^{2^{*}}dv_{\mathbf{g}}=1$ et en rappelant que 
\[
I_{h,\mathbf{g}}(w)=I_{h}(w)=\int_{M}\left| \nabla w\right| ^{2}dv_{\mathbf{g%
}}+\int_{M}h.w{{}^{2}}dv_{\mathbf{g}} 
\]
on obtient, en notant $B_{0}K^{-2}=B_{0}(\mathbf{g})K(n,2)^{-2}$%
\[
I_{B_{0}K^{-2}}(c)=\int B_{0}K^{-2}c{{}^{2}d}v{_{\mathbf{g}}=}B_{0}K^{-2}c{%
{}^{2}Vol}_{\mathbf{g}}(M)\,\,. 
\]
Pour $t\in \Bbb{R}^{+}$ on pose 
\[
h_{t}=B_{0}K^{-2}-t\eta \,\,. 
\]
Alors $h_{t}=B_{0}K^{-2}$ sur $\left\{ x/f(x)=Maxf\right\} $, et 
\[
I_{h_{t}}(c)=\int_{M}B_{0}K^{-2}c{{}^{2}}dv_{\mathbf{g}}-c{{}^{2}}%
t\int_{M}\eta dv_{\mathbf{g}}=\left( \int_{M}fdv_{\mathbf{g}}\right) ^{-%
\frac{2}{2^{*}}}(B_{0}K^{-2}{Vol}_{\mathbf{g}}(M)-t\int_{M}\eta dv_{\mathbf{g%
}}) 
\]
Donc si $t$ est assez grand, 
\[
I_{h_{t}}(c)<\frac{1}{K(n,2){{}^{2}}(\stackunder{M}{Sup}f)^{\frac{n-2}{n}}}%
\,\,. 
\]
Maintenant, on voudrait aussi que $h_{t}$ soit strictement positive sur $M$.
Il faut d'apr\`{e}s la d\'{e}finition de $h_{t}$ et puisque $\stackunder{M}{%
Sup}\,\eta =1$ que
\begin{equation}
t<B_{0}(\mathbf{g})K(n,2)^{-2}\,\,.  
\end{equation}
Mais pour que 
\[
I_{h_{t}}(c)<\frac{1}{K(n,2){{}^{2}}(\stackunder{M}{Sup}f)^{\frac{n-2}{n}}} 
\]
il faut aussi que 
\begin{equation}
t>\frac{1}{\int_{M}\eta dv_{\mathbf{g}}}(B_{0}K^{-2}{Vol}_{\mathbf{g}}(M)-%
\frac{\left( \int_{M}fdv_{\mathbf{g}}\right) ^{\frac{n-2}{n}}}{K(n,2){{}^{2}}%
(\stackunder{M}{Sup}f)^{\frac{n-2}{n}}})\,\,. 
\end{equation}
On peut donc trouver un tel $t$ d\`{e}s que 
\[
\frac{1}{\int_{M}\eta dv_{\mathbf{g}}}(B_{0}K^{-2}{Vol}_{\mathbf{g}}(M)-%
\frac{\left( \int_{M}fdv_{\mathbf{g}}\right) ^{\frac{n-2}{n}}}{K(n,2){{}^{2}}%
(\stackunder{M}{Sup}f)^{\frac{n-2}{n}}})<B_{0}K^{-2} 
\]
ce que l'on peut \'{e}crire 
\begin{equation}
\int_{M}\eta dv_{\mathbf{g}}>{Vol}_{\mathbf{g}}(M)-\frac{\left( \int_{M}fdv_{%
\mathbf{g}}\right) ^{\frac{n-2}{n}}}{B_{0}(\mathbf{g})(\stackunder{M}{Sup}%
f)^{\frac{n-2}{n}}}\,\,.  
\end{equation}
Rappelons qu'on veut $\eta $ \`{a} support inclus dans $M\backslash \left\{
x/f(x)=Maxf\right\} $ et telle que $0\leq \eta \leq 1$. Mais l'hypoth\`{e}se
que le Hessien de $f$ est non-d\'{e}g\'{e}n\'{e}r\'{e} en ses points de
maximum implique que $\left\{ x/f(x)=Maxf\right\} $, l'ensemble des points
de maximum de $f$, est un ensemble de points isol\'{e}s. On peut donc
trouver une telle fonction $\eta $ telle que de plus $\int_{M}\eta dv_{%
\mathbf{g}}$ soit aussi proche que l'on veut de ${Vol}_{\mathbf{g}}(M)$.
Comme 
\[
\frac{\left( \int_{M}fdv_{\mathbf{g}}\right) ^{\frac{n-2}{n}}}{B_{0}(\mathbf{%
g})(\stackunder{M}{Sup}f)^{\frac{n-2}{n}}}>0 
\]
on peut trouver $\eta $ v\'{e}rifiant (4.3) et donc un $t$, not\'{e} $t_{1}$%
, v\'{e}rifiant (4.1) et (4.2).

Sur $\left\{ x/f(x)=Maxf\right\} $, $h_{t}=B_{0}K^{-2}$ donc $\forall P\in
\left\{ x/f(x)=Maxf\right\} $:

\[
\frac{4(n-1)}{n-2}h_{t}(P)>S_{\mathbf{g}}(P)-\frac{n-4}{2}\frac{%
\bigtriangleup _{\mathbf{g}}f(P)}{f(P)}. 
\]
On pose alors: 
\[
t_{0}=Inf\{t\leq t_{1}\,/\,\lambda _{h_{t}}<\frac{1}{K(n,2){{}^{2}}(%
\stackunder{M}{Sup}f)^{\frac{n-2}{n}}}\}\text{ }. 
\]
N\'{e}cessairement $t_{0}<t_{1}$. On rappelle que 
\[
\lambda _{h,f,\mathbf{g}}=\lambda _{h}=\stackunder{w\in \mathcal{H}_{f}}{%
\inf }I_{h}(w) 
\]
o\`{u} 
\[
\mathcal{H}_{f}=\{w\in H_{1}^{2}(M)/\int_{M}f\left| w\right| ^{\frac{2n}{n-2}%
}dv_{\mathbf{g}}=1\}\,\,. 
\]
Alors 
\[
\lambda _{h_{t_{0}}}=\frac{1}{K(n,2){{}^{2}}(\stackunder{M}{Sup}f)^{\frac{n-2%
}{n}}}\,\,\,et\,\,\,\lambda _{h_{t}}<\frac{1}{K(n,2){{}^{2}}(\stackunder{M}{%
Sup}f)^{\frac{n-2}{n}}}\,\,\,si\,\,\,t>t_{0}. 
\]
De plus $\forall t$, $t_{0}\leq t\leq t_{1}$, $h_{t}>0$ sur $M$ et

\[
\frac{4(n-1)}{n-2}h_{t_{0}}(P)>S_{\mathbf{g}}(P)-\frac{n-4}{2}\frac{%
\bigtriangleup _{\mathbf{g}}f(P)}{f(P)}\,\,\,pour\,\,\,P\in \left\{
x/f(x)=Maxf\right\} . 
\]
Enfin $h_{t}\stackunder{t\rightarrow t_{0}}{\rightarrow }h_{t_{0}}\,\,$dans$%
\,$ $C^{0}$, et comme $h_{t_{0}}>0$, $\triangle _{\mathbf{g}}+h_{t_{0}}$ est
coercif. Donc d'apr\`{e}s le th\'{e}or\`{e}me 1, $h_{t_{0}}$ est critique
avec des fonctions extr\'{e}males.

\medskip

Venons-en au cas g\'{e}n\'{e}ral, c'est \`{a} dire o\`{u} $f$ peut changer
de signe sur $M$ et o\`{u} l'on ne suppose plus $\int fdv_{\mathbf{g}}>0$
(mais toujours $Maxf>0$). L'id\'{e}e est de modifier la fonction faiblement
critique pour $f$, $B_{0}(\mathbf{g})K(n,2)^{-2}$ non plus par une fonction
plateau $\eta $ comme dans la d\'{e}monstration pr\'{e}c\'{e}dente, mais par
les fonctions tests d\'{e}crites au chapitre 1. On peut les pr\'{e}senter de
la mani\`{e}re suivante: Pour tout point $x\in M$ et pour tout $\delta >0$
assez petit, il existe une suite de fonctions $(\psi _{k})$ \`{a} support
inclus dans $B(x,\delta )$ telle que pour toute fonction $h$: 
\[
J_{h,1,\mathbf{g}}(\psi _{k})=\frac{\int_{M}\left| \nabla \psi _{k}\right|
^{2}dv_{\mathbf{g}}+\int_{M}h.\psi _{k}{{}^{2}}dv_{\mathbf{g}}}{\left(
\int_{M}\left| \psi _{k}\right| ^{\frac{2n}{n-2}}dv_{\mathbf{g}}\right) ^{%
\frac{n-2}{n}}}\stackunder{k\rightarrow \infty }{\rightarrow }\frac{1}{K(n,2)%
{{}^{2}}} 
\]
et 
\[
\int_{M}\psi _{k}^{\frac{2n}{n-2}}dv_{\mathbf{g}}=1 
\]
cette derni\`{e}re condition \'{e}tant obtenue en multipliant les fonctions
d\'{e}crites au chapitre 1 par des constantes. Il est ici plus commode
d'utiliser la fonctionnelle $J$ car nous aurons 
\[
\int_{M}f.\psi _{k}^{\frac{2n}{n-2}}dv_{\mathbf{g}}\neq 1\,\,. 
\]
Soit alors $\psi _{k}$ une de ces fonctions o\`{u} $k$ et $B(x,\delta )$
seront fix\'{e}s plus tard. On consid\`{e}re pour $t>0$ la famille 
\[
h_{t}=B_{0}(\mathbf{g})K(n,2)^{-2}-t.\psi _{k}^{\frac{4}{n-2}}\,\,. 
\]
Tout d'abord, cherchons \`{a} quelle condition $\triangle _{\mathbf{g}%
}+h_{t} $ est coercif. Toujours en notant $B_{0}K^{-2}=B_{0}(\mathbf{g}%
)K(n,2)^{-2}$, et en sous-entendant que les int\'{e}grales sont prises pour
la mesure $dv_{\mathbf{g}}$, pour $u\in H_{1}^{2}$%
\begin{eqnarray*}
\int_{M}(\left| \nabla u\right| _{\mathbf{g}}^{2}+h_{t}u^{2}) && 
=\int_{M}(\left| \nabla u\right| _{\mathbf{g}}^{2}+B_{0}K^{-2}.u^{2})-t%
\int_{M}\psi _{k}^{\frac{4}{n-2}}.u^{2} \\ 
&& \geq K(n,2)^{-2}(\int_{M}u^{\frac{2n}{n-2}})^{\frac{n-2}{n}}-t\int_{M}\psi
_{k}^{\frac{4}{n-2}}.u^{2}
\end{eqnarray*}
d'apr\'{e}s l'in\'{e}galit\'{e} de Sobolev. Or par l'in\'{e}galit\'{e} de
H\"{o}lder: 
\[
\int_{M}\psi _{k}^{\frac{4}{n-2}}.u^{2}\leq (\int_{M}\psi _{k}^{\frac{2n}{n-2%
}})^{\frac{n-2}{n}}(\int_{M}u^{\frac{2n}{n-2}})^{\frac{2}{n}}=(\int_{M}u^{%
\frac{2n}{n-2}})^{\frac{2}{n}} 
\]
car $\int_{M}\psi _{k}^{\frac{2n}{n-2}}=1$. Donc, en utilisant
l'in\'{e}galit\'{e} de H\"{o}lder pour dire qu'il existe une constante $C>0$
telle que 
\[
C\int_{M}u^{2}\leq (\int_{M}u^{\frac{2n}{n-2}})^{\frac{n-2}{n}} 
\]
on obtient, d\`{e}s que $K(n,2)^{-2}-t>0$%
\begin{eqnarray*}
\int_{M}(\left| \nabla u\right| _{\mathbf{g}}^{2}+h_{t}u^{2})& & \geq
(K(n,2)^{-2}-t)(\int_{M}u^{\frac{2n}{n-2}})^{\frac{n-2}{n}} \\ 
& &\geq (K(n,2)^{-2}-t)C\int_{M}u^{2}\,\,\,.
\end{eqnarray*}
Ainsi $\triangle _{\mathbf{g}}+h_{t}$ est coercif d\`{e}s que $t<K(n,2)^{-2}$%
; on fixe alors $t_{1}$, $0<t_{1}<K(n,2)^{-2}$. On voudrait maintenant fixer 
$\psi _{k}$ pour que $h_{t_{1}}$ soit sous-critique pour $f$. Or, si l'on
choisit d\'{e}j\`{a} $x$ assez proche d'un point $x_{0}$ de maximum de $f$
et $\delta $ assez petit pour que $f>0$ sur $B(x,\delta )$, on obtient : 
\begin{eqnarray*}
J_{h_{t_{1}},f,\mathbf{g}}(\psi _{k}) && =\frac{\int_{M}\left| \nabla \psi
_{k}\right| ^{2}+\int_{M}B_{0}K^{-2}.\psi _{k}{{}^{2}-t}_{1}\int_{M}\psi
_{k}^{\frac{2n}{n-2}}}{\left( \int_{M}f\left| \psi _{k}\right| ^{\frac{2n}{%
n-2}}\right) ^{\frac{n-2}{n}}} \\ 
&& \leq \frac{J_{B_{0}K^{-2},1}(\psi _{k})}{(\stackunder{B(x,\delta )}{Inf}%
f)^{\frac{n-2}{n}}}-\frac{t_{1}}{(\stackunder{B(x,\delta )}{Sup}f)^{\frac{n-2%
}{n}}} \\ 
&& \leq \frac{J_{B_{0}K^{-2},1}(\psi _{k})}{(\stackunder{B(x,\delta )}{Inf}%
f)^{\frac{n-2}{n}}}-\frac{t_{1}}{(\stackunder{M}{Sup}f)^{\frac{n-2}{n}}}%
\,\,\,.
\end{eqnarray*}
Pour tout $\varepsilon >0$, on peut, par continuit\'{e} de $f$, choisir $x$
assez proche de $x_{0}$ et $\delta $ assez petit de sorte que $B(x,\delta
)\cap \left\{ x/f(x)=Maxf\right\} =\emptyset $ et que 
\[
\frac{1}{(\stackunder{B(x,\delta )}{Inf}f)^{\frac{n-2}{n}}}\leq \frac{1}{(%
\stackunder{M}{Sup}f)^{\frac{n-2}{n}}}+\varepsilon 
\]
$x$ et $\delta $ \'{e}tant fix\'{e}s, on peut maintenant choisir $k$ assez
grand pour que 
\[
J_{B_{0}K^{-2},1}(\psi _{k})\leq K(n,2)^{-2}+\varepsilon \,\,. 
\]
Par cons\'{e}quent, en choisissant $\varepsilon $ assez petit, on voit que
puisque $\frac{t_{1}}{(\stackunder{M}{Sup}f)^{\frac{n-2}{n}}}>0$ : 
\[
J_{h_{t_{1}},f,\mathbf{g}}(\psi _{k})<\frac{1}{K(n,2)^{-2}(\stackunder{M}{Sup%
}f)^{\frac{n-2}{n}}} 
\]
et donc que $h_{t_{1}}$ est sous-critique pour $f$. On conclut alors de la
m\^{e}me mani\`{e}re en posant:

\[
t_{0}=Inf\{t\leq t_{1}/\lambda _{h_{t}}<\frac{1}{K(n,2){{}^{2}}(\stackunder{M%
}{Sup}f)^{\frac{n-2}{n}}}\}\text{ }. 
\]
Alors $t_{0}\geq 0$, et 
\[
\lambda _{h_{t_{0}}}=\frac{1}{K(n,2){{}^{2}}(\stackunder{M}{Sup}f)^{\frac{n-2%
}{n}}}\,\,\,et\,\,\,\lambda _{h_{t}}<\frac{1}{K(n,2){{}^{2}}(\stackunder{M}{%
Sup}f)^{\frac{n-2}{n}}}\,\,\,si\,\,\,t>t_{0}. 
\]
De plus $\forall t$, $t_{0}\leq t\leq t_{1}$,

\[
\frac{4(n-1)}{n-2}h_{t_{0}}(P)>S_{\mathbf{g}}(P)-\frac{n-4}{2}\frac{%
\bigtriangleup _{\mathbf{g}}f(P)}{f(P)}pourP\in \left\{ x/f(x)=Maxf\right\} 
\]
puisque $B(x,\delta )\cap \left\{ x/f(x)=Maxf\right\} =\emptyset $. Enfin $%
h_{t}\stackunder{t\rightarrow t_{0}}{\rightarrow }h_{t_{0}}\,\,$dans$\,$ $%
C^{0,\alpha }$, et $\triangle _{\mathbf{g}}+h_{t_{0}}$ est coercif. Donc
d'apr\`{e}s le th\'{e}or\`{e}me 1, $h_{t_{0}}$ est critique avec des
fonctions extr\'{e}males.

\textit{Remarque 1}: Après la premi\`{e}re d\'{e}finition donn\'{e}e des
fonctions critiques (i.e. sous la forme h est critique pour f et g) la
question de l'existence de telles fonctions est naturelle; la notion de
triplet critique ne s'impose pas vraiment... N\'{e}anmoins, comme nous
l'avons fait remarquer au chapitre 1, le probl\`{e}me de l'existence se pose
naturellement dans le cadre des questions que nous avons pos\'{e}es au sujet
de ces triplets. En effet, savoir s'il existe une fonction critique pour $f$
et $\mathbf{g}$ revient au probl\`{e}me suivant: \textit{On se fixe }$f$%
\textit{\ et la m\'{e}trique }$\mathbf{g}$\textit{\ et on cherche une
fonction }$h$\textit{\ telle que (}$h,f,\mathbf{g}$\textit{) soit un triplet
critique. }D'o\`{u} le titre du chapitre...

\textit{Remarque 2:} Les d\'{e}monstrations pr\'{e}c\'{e}dentes montrent
\'{e}galement, en rempla\c{c}ant $B_{0}K^{-2}$ par $h$ que si $(h,f,\mathbf{g%
})$ est faiblement critique, et si 
\[
\frac{4(n-1)}{n-2}h(P)>S_{\mathbf{g}}(P)-\frac{n-4}{2}\frac{\bigtriangleup _{%
\mathbf{g}}f(P)}{f(P)}\text{ \thinspace pour \thinspace }P\in \left\{
x/f(x)=Maxf\right\} 
\]
alors il existe $h^{\prime }$ $\leqslant h$ telle que $(h^{\prime },f,g)$
soit critique.

Si on a seulement 
\[
\frac{4(n-1)}{n-2}h(P)\geq S_{\mathbf{g}}(P)-\frac{n-4}{2}\frac{%
\bigtriangleup _{\mathbf{g}}f(P)}{f(P)}\text{ \thinspace pour \thinspace }%
P\in \left\{ x/f(x)=Maxf\right\} 
\]
alors pour tout $\varepsilon >0$ il existe $h^{\prime }$ $\leqslant
h+\varepsilon \,\,$telle que $(h^{\prime },f,g)$ soit critique.

\chapter{Triplet Critique 2}
\pagestyle{myheadings}\markboth{\textbf{Triplet critique 2.}}
{\textbf{Triplet critique 2.}}
On d\'{e}montre le th\'{e}or\`{e}me 3 que l'on rappelle ici:

\textbf{Th\'{e}or\`{e}me 3:}

\textit{Soient donn\'{e}es la vari\'{e}t\'{e} (}$M,\mathbf{g}$\textit{) et
deux fonctions }$h^{\prime }$ \textit{et }$f$\textit{\ v\'{e}rifiant les
hypoth\`{e}ses }\textbf{(H}).\textit{\ Alors il existe une m\'{e}trique }$%
\mathbf{g}^{\prime }$ conforme \`{a} $\mathbf{g}$ \textit{telle que }$%
(h^{\prime },f,\mathbf{g}^{\prime })$\textit{\ soit un triplet critique, ou,
pour reprendre la pr\'{e}sentation premi\`{e}re, il existe une m\'{e}trique }%
$\mathbf{g}^{\prime }$ conforme \`{a} $\mathbf{g}$ \textit{telle que }$%
h^{\prime }$\textit{\ soit critique pour }$f$\textit{\ et }$\mathbf{g}%
^{\prime }$. \textit{De plus }$(h^{\prime },f,\mathbf{g}^{\prime })$\textit{%
\ a des solutions minimisantes.}

Ce th\'{e}or\`{e}me repose fondamentalement sur la formule de transformation
d'une fonction critique dans un changement de m\'{e}trique conforme:
\begin{center}
\textit{\ (}$h,f$\textit{,}$\mathbf{g)}$\textit{\ est critique\ si et
seulement si (}$h^{\prime }=\frac{\triangle _{\textbf{g}}u+h.u}{u^{\frac{n+2}{n-2}}}$%
\textit{,}$f,\mathbf{g}^{\prime }=u^{\frac{4}{n-2}}\mathbf{g)}$ \textit{est
critique.}
\end{center}
Ou de fa\c{c}on \'{e}quivalente
\begin{center}
($h^{\prime },f$\textit{,}$\mathbf{g}^{\prime }=u^{\frac{4}{n-2}}\mathbf{g)}$
\textit{est critique\ si et seulement si (}$h=h^{\prime }u^{\frac{4}{n-2}}-%
\frac{\triangle _{\mathbf{g}}u}{u}$\textit{,}$f$\textit{,}$\mathbf{g)}$ 
\textit{est critique.}
\end{center}
Notons, pour $u\in C_{+}^{\infty }(M)=\{u\in C^{\infty }(M)\,/\,u>0\}$ : 
\[
F_{h^{\prime }}(u)=h^{\prime }u^{\frac{4}{n-2}}-\frac{\bigtriangleup _{%
\mathbf{g}}u}{u} 
\]
Alors:
\begin{center}
$(h^{\prime },f,\mathbf{g}^{\prime })$\textit{\ est critique si
et seulement si }$(F_{h^{\prime }}(u),f,\mathbf{g})$\textit{\ est critique.}
\end{center}
Pour r\'{e}soudre le probl\`{e}me, il suffit donc de trouver une fonction $h$
telle que:

1/: $\triangle _{\mathbf{g}}u+h.u=h^{\prime }u^{\frac{n+2}{n-2}}$ ait une
solution $u>0$, et

2/: $(h,f,\mathbf{g})\,$soit critique.

En effet, dans ce cas $h=F_{h^{\prime }}(u)$ et $h^{\prime }$ est critique
pour $f$ et $\mathbf{g}^{\prime }=u^{\frac{4}{n-2}}\mathbf{g}$.

Pour trouver cette fonction, une premi\`{e}re tentative naturelle est de
reprendre la m\'{e}thode de la d\'{e}monstration du th\'{e}or\`{e}me de
Yamabe. D'abord, on sait qu'il existe une fonction $h$ critique pour $f$%
\textit{\ }et $\mathbf{g}$\ . On consid\`{e}re alors une suite 
\[
q_{t}\rightarrow 2^{*},\,q_{t}<2^{*}. 
\]
Pour tout $t$ il existe une fonction $u_{t}>0$ solution de 
\[
\triangle _{\mathbf{g}}u_{t}+h.u_{t}=h^{\prime }u_{t}^{q_{t}-1} 
\]
avec 
\[
\int h^{\prime }u_{t}^{q_{t}}dv_{\mathbf{g}}\leq C\text{ ind\'{e}pendant de }%
t 
\]
Apr\`{e}s extraction, $u_{t}\stackrel{H_{1}^{2}}{\rightharpoondown }%
u\geqslant 0$. Si $u>0$, alors $u$ est solution (\`{a} une constante
multiplicative pr\`{e}s) de 
\[
\triangle _{\mathbf{g}}u+h.u=h^{\prime }u^{\frac{n+2}{n-2}} 
\]
et c'est fini. Mais le probl\`{e}me est encore une fois d'\'{e}viter la
solution nulle...

E. Humbert et M. Vaugon ont r\'{e}solu le probl\`{e}me dans le cas $f=cste$
et pour une vari\'{e}t\'{e} non conform\'{e}ment diff\'{e}omorphe \`{a} la
sph\`{e}re [21]. Leur m\'{e}thode s'appuie sur le fait que, pour une vari%
\'{e}t\'{e} non conform\'{e}ment diff\'{e}omorphe \`{a} la sph\`{e}re,
quitte \`{a} faire un premier changement conforme de m\'{e}trique, $B_{0}(%
\mathbf{g})K(n,2)^{-2}$ est une fonction critique (quand il n'y aura pas
d'ambigu\"{i}t\'{e}, nous noterons ces deux constantes $K$ et $B_{0}$). En
fait, une analyse pr\'{e}cise de leur d\'{e}monstration fait appara\^{i}tre
que l'important, c'est que $B_{0}(\mathbf{g})K(n,2)^{-2}$ soit une fonction
critique strictement positive, (non qu'elle soit constante). Or nous avons
justement montr\'{e} au chapitre pr\'{e}c\'{e}dent que si $f$ v\'{e}rifiant 
\textit{\ }\textbf{(H}$_{f}$) est telle que $\int_{M}fdv_{\mathbf{g}}>0$, il
existe des fonctions critiques strictement positives. Nous avions d'ailleurs 
\'{e}tabli ce r\'{e}sultat pour montrer le th\'{e}or\`{e}me ci-dessus, mais
il semblait plus logique de le placer dans le chapitre sur l'existence de
fonctions critiques. Notons qu'avec l'hypoth\`{e}se \textbf{(H}$_{f}$), la d%
\'{e}monstration fonctionne \'{e}galement sur la sph\`{e}re, mais que, \'{e}%
videmment, cette hypoth\`{e}se ne s'applique pas \`{a} une fonction
constante !

Voici le principe de la d\'{e}monstration d'E. Humbert et M. Vaugon. Sans
pr\'{e}tendre nous l'approprier, nous la pr\'{e}sentons d'une mani\`{e}re un
peu diff\'{e}rente et en int\'{e}grant notre fonction non constante. Leur
id\'{e}e est de partir de la m\'{e}thode naturelle expos\'{e}e ci-dessus,
mais de pr\'{e}voir une solution dans le cas o\`{u} $u_{t}\stackrel{H_{1}^{2}%
}{\rightharpoondown }0$. On sait qu'il existe une suite ($h_{t}$) de
fonctions sous-critiques pour $f$\textit{\ }et $\mathbf{g}$\ telle que $h_{t}%
\stackrel{C^{2}}{\rightarrow }h$ o\`{u} $(h,f,\mathbf{g})\,$est critique et
telle qu'en tout point $P$\ o\`{u} $f$\ est maximum sur $M$ : 
\[
\frac{4(n-1)}{n-2}h(P)>S_{\mathbf{g}}(P)-\frac{n-4}{2}\frac{\bigtriangleup _{%
\mathbf{g}}f(P)}{f(P)}. 
\]
Pour une suite $q_{t}\rightarrow 2^{*},\,q_{t}<2^{*}$ on construit une suite 
$u_{t}>0$ de solutions de 
\[
\triangle _{\mathbf{g}}u+h.u=h^{\prime }u^{q_{t}-1}\,\,\,avec\,\,\,\int
h^{\prime }u_{t}^{q_{t}}dv_{\mathbf{g}}\leq C\text{ ind\'{e}pendant de }t 
\]
telle que $u_{t}\stackrel{H_{1}^{2}}{\rightharpoondown }u\geqslant 0$. Ici
encore, si $u>0$, alors $u$ est solution (\`{a} une constante multiplicative
près) de $\triangle _{\mathbf{g}}u+h.u=h^{\prime }u^{\frac{n+2}{n-2}}$
et c'est fini.

Maintenant, si $u=0$, nous montrerons que les $u_{t}$ se concentrent et que
gr\^{a}ce \`{a} cela on peut trouver $t_{0}\,$assez proche de 1 (si on a
pris $t\rightarrow 1$) et $s$ assez grand de telle sorte que $F_{h^{\prime
}}(u_{t_{0}})$ soit sous-critique et que $F_{h^{\prime }}(u_{t_{0}}^{s})$
soit faiblement critique (o\`{u} $u_{t_{0}}^{s}$ est la fonction $u_{t_{0}}$
\'{e}lev\'{e}e \`{a} la puissance $s$), avec en plus 
\[
\frac{4(n-1)}{n-2}F_{h^{\prime }}(u_{t_{0}}^{s})(P)>S_{\mathbf{g}}(P)-\frac{n-4}{2}\frac{%
\bigtriangleup _{\mathbf{g}}f(P)}{f(P)} 
\]
en tout point $P$\ o\`{u} $f$\ est maximum sur $M$\ . Alors, en
consid\'{e}rant le chemin $t\rightarrow F_{h^{\prime }}(u_{t_{0}}^{ts})$, on
obtient \`{a} l'aide du th\'{e}or\`{e}me du chapitre 3 l'existence sur ce
chemin d'une fonction critique. L'astuce consiste en ce choix du chemin,
plut\^{o}t que l'habituel chemin ``lin\'{e}aire''. Par ailleurs, c'est pour
obtenir les conditions sur $F_{h^{\prime }}(u_{t_{0}}^{s})$ qu'intervient
fondamentalement l'existence de fonctions critiques positives.

Passons \`{a} la d\'{e}monstration.

Tout d'abord, nous avons expliqu\'{e} que nous aurons besoin d'une fonction
critique positive. Or nous avons \'{e}tabli cette existence sous
l'hypoth\`{e}se $\int_{M}fdv_{\mathbf{g}}>0$. Mais $\stackunder{M}{Sup}f>0$,
donc, quitte \`{a} faire un premier changement conforme de m\'{e}trique, on
peut toujours supposer que $\int_{M}fdv_{\mathbf{g}}>0$ et donc qu'il existe
des fonctions critiques strictement positives\textit{\ }pour $f$\ et $%
\mathbf{g}$.

Ensuite, fixons quelques notations (de plus):
\begin{center}
\[
J_{h,h^{\prime },\mathbf{g,}q}(w)=\frac{\int_{M}\left| \nabla w\right|
^{2}dv_{\mathbf{g}}+\int_{M}h.w{{}^{2}}dv_{\mathbf{g}}}{\left(
\int_{M}h^{\prime }\left| w\right| ^{q}dv_{\mathbf{g}}\right) ^{\frac{2}{q}}}
\]
\end{center}

\[
\stackunder{w\in \mathcal{H}_{h^{\prime },q}^{+}}{\inf }J_{h,h^{\prime },%
\mathbf{g,}q}(w):=\lambda _{h,h^{\prime },\mathbf{g,}q} 
\]
o\`{u} 
\[
\mathcal{H}_{h^{\prime },q}^{+}=\{w\in
H_{1}^{2}(M)\,/\,\,\,w>0\,\,et\,\,\int_{M}h^{\prime }.w^{q}dv_{\mathbf{g}%
}>0\}. 
\]
et
\[
\Omega _{h,h^{\prime },\mathbf{g,}q}=\{u\in \mathcal{H}_{h^{\prime
},q}^{+}/\,J_{h,h^{\prime },\mathbf{g,}q}(u)=\lambda _{h,h^{\prime },\mathbf{%
g,}q}\,\,et\,\,\int_{M}h^{\prime }.w^{q}dv_{\mathbf{g}}=(\lambda
_{h,h^{\prime },\mathbf{g,}q}\,\,)^{\frac{q}{q-2}}\}\,\,. 
\]

Soit une suite ($h_{t}$) de fonctions sous-critiques pour $f$\textit{\ }et $%
\mathbf{g}$\ telle que $h_{t}\stackrel{C^{2}}{\rightarrow }h$ o\`{u} $(h,f,%
\mathbf{g})\,$est critique, avec $\triangle _{\mathbf{g}}+h_{t}$ coercif. On
sait que l'on peut trouver une telle suite avec $h_{t}>0$ et $h>0$, et en
plus 
\[
\frac{4(n-1)}{n-2}h_{t}(P)>S_{\mathbf{g}}(P)-\frac{n-4}{2}\frac{\bigtriangleup _{\mathbf{g}%
}f(P)}{f(P)} 
\]
en tout point $P$\ o\`{u} $f$\ est maximum sur $M$. Mais on peut faire mieux
dans le cadre de notre probl\`{e}me, et c'est l\`{a} qu'intervient
fondamentalement l'inter\^{e}t des fonctions critiques strictement
positives. En effet, pour toute constante $c>0$, si $\mathbf{g}^{\prime }=c%
\mathbf{g}$, alors $S_{\mathbf{g}^{\prime }}=c^{-1}S_{\mathbf{g}}$ et $%
\triangle _{\mathbf{g}^{\prime }}=c^{-1}\triangle _{\mathbf{g}}$ et
d'apr\`{e}s la règle de transformation des fonctions critiques:

\begin{center}
\textit{\ }$h$\textit{\ est (sous-, faiblement) critique pour }$f$\textit{\
et }$\mathbf{g}$\textit{\ si et seulement si }$c^{-1}h$\textit{\ est (sous-,
faiblement) critique pour }$f$\textit{\ et }$\mathbf{g}^{\prime }$
\end{center}

Donc, quitte \`{a} multiplier $\mathbf{g}$ par une constante, on peut, pour
toute constante $C>0$, supposer au d\'{e}part: 
\begin{eqnarray*}
h_{t} &>&C\,\,surM \\
\frac{4(n-1)}{n-2}h_{t}(P)-S_{\mathbf{g}}(P)+\frac{n-4}{2}\frac{\bigtriangleup _{\mathbf{g}%
}f(P)}{f(P)} &>&C\,\,\text{en tout point de maximum de }f
\end{eqnarray*}
et $(h,f,\mathbf{g)}$ a des solutions minimisantes d'apr\`{e}s le chapitre 4.

Une fois ceci pos\'{e}, on peut suivre le principe expos\'{e} plus haut.

\textit{Premi\`{e}re \'{e}tape:} Il est connu que pour tout $q<2^{*}$ et
pour tout $u\in \Omega _{h,h^{\prime },\mathbf{g,}q}$, $u$ est solution de $%
\triangle _{\mathbf{g}}u+h.u=h^{\prime }u^{q-1}$, les m\'{e}thodes
variationnelles classiques fonctionnant sans probl\`{e}me car l'inclusion $%
H_{1}^{2}\subset L^{q}$ est compacte. Nous voulons montrer que:

$\forall t,\,\exists \,q_{0}<2^{*}$\textit{\ tel que }$\forall q\in
[q_{0},2^{*}[$\textit{\ et }$\forall u\in \Omega _{h_{t},h^{\prime },\mathbf{%
g,}q}$\textit{\ alors }$F_{h^{\prime }}(u)$\textit{\ est sous-critique pour }%
$f$\textit{\ et }$\mathbf{g.}$

En effet, si l'on suppose le contraire alors il existe une suite $q_{i}%
\stackrel{<}{\rightarrow }2^{*}$ et des fonctions 
\[
u_{i}\in \Omega _{h_{t},h^{\prime },\mathbf{g,}q_{i}} 
\]
telles que $F_{h^{\prime }}(u_{i})$\textit{\ }est faiblement critique pour $%
f $\ et $\mathbf{g.}$ La suite ($u_{i}$) est born\'{e}e dans $H_{1}^{2}$ et
il existe $u\in H_{1}^{2}$ tel que $u_{i}\stackrel{H_{1}^{2}}{%
\rightharpoondown }u$. Les th\'{e}ories elliptiques standard nous donnent
alors deux possibilit\'{e}s: soit $u\neq 0$, soit $u\equiv 0$.

Si $u$ est non identiquement nulle, ces m\^{e}me th\'{e}ories nous donnent
que $u>0\,$et que $u\in \mathcal{H}_{h^{\prime },2^{*}}^{+}$, puis qu'\`{a}
extraction près $u_{i}\stackrel{C^{2}}{\rightarrow }u$. Alors $%
F_{h^{\prime }}(u_{i})\stackrel{C^{0}}{\rightarrow }F_{h^{\prime }}(u)\,$et
donc $F_{h^{\prime }}(u)\,$est faiblement critique. Or $u_{i}\in \Omega
_{h_{t},h^{\prime },\mathbf{g,}q_{i}}$ donc 
\begin{eqnarray*}
F_{h^{\prime }}(u_{i})& & =h^{\prime }u_{i}^{\frac{4}{n-2}}-\frac{%
\bigtriangleup _{\mathbf{g}}u_{i}}{u_{i}}=h^{\prime }u_{i}^{\frac{4}{n-2}%
}+h_{t}-h^{\prime }u_{i}^{q_{i}-2} \\ 
&& =h_{t}+h^{\prime }(u_{i}^{\frac{4}{n-2}}-u_{i}^{q_{i}-2})
\end{eqnarray*}
et, par cons\'{e}quent, quand $i\rightarrow +\infty $, $F_{h^{\prime
}}(u_{i})\rightarrow h_{t}$, et finalement $F_{h^{\prime }}(u)=h_{t}$. Or $%
h_{t}$ est sous-critique et $F_{h^{\prime }}(u)$ est faiblement critique,
d'o\`{u} une contradiction.

Si maintenant on suppose que $u\equiv 0$, $h_{t}$ \'{e}tant sous-critique,
il existe une fonction $\phi \in C^{\infty }(M)$ telle que 
\[
J_{h_{t},f,\mathbf{g,}2^{*}}(\phi )<\frac{1}{K(n,2)^{2}(Sup\,f)^{\frac{n-2}{n%
}}}\,\,. 
\]
Alors 
\[
J_{F_{h^{\prime }}(u_{i}),f,\mathbf{g,}2^{*}}(\phi )=J_{h_{t},f,\mathbf{g,}%
2^{*}}(\phi )+\frac{\int h^{\prime }(u_{i}^{\frac{4}{n-2}}-u_{i}^{q_{i}-2})%
\phi ^{2}dv_{\mathbf{g}}}{(\int f\phi ^{2^{*}}dv_{\mathbf{g}})^{\frac{2}{%
2^{*}}}}\,\,. 
\]
Or $u_{i}\stackrel{L^{\frac{4}{n-2}}}{\rightarrow }0$ et $q_{i}\leqslant
2^{*}$ donc 
\[
J_{F_{h^{\prime }}(u_{i}),f,\mathbf{g,}2^{*}}(\phi )\rightarrow J_{h_{t},f,%
\mathbf{g,}2^{*}}(\phi )<\frac{1}{K(n,2)^{2}(Sup\,f)^{\frac{n-2}{n}}} 
\]
ce qui contredit le fait que ($F_{h^{\prime }}(u_{i}),f,\mathbf{g}$) est
faiblement critique. On peut r\'{e}\'{e}crire le r\'{e}sultat de cette
\'{e}tape sous la forme suivante:

Il existe des suites $(q_{i}),(t_{i}),$ telles que 
\begin{eqnarray*}
\,2 &<&q_{i}<2^{*} \\
\,q_{i} &\rightarrow &2^{*} \\
\,t_{i} &\rightarrow &1 \\
\,h_{t_{i}} &\rightarrow &h\,
\end{eqnarray*}
et une suite ($v_{i})\in \Omega _{h_{t_{i}},h^{\prime },\mathbf{g,}q_{i}}$
telle que ($F_{h^{\prime }}(v_{i}),f,\mathbf{g)\,}$soit sous-critique.
Notons 
\[
J_{i}=J_{h_{t_{i}},h^{\prime },\mathbf{g,}q_{i}} 
\]
$\,$et 
\[
\lambda _{i}=\lambda _{h_{t_{i}},h^{\prime },\mathbf{g,}q_{i}}\,\,\,. 
\]
Alors 
\[
J_{i}(v_{i})=\lambda _{i}\,\,\,et\,\,\int h^{\prime }v_{i}^{q_{i}}dv_{%
\mathbf{g}}=\lambda _{i}^{\frac{q_{i}}{q_{i}-2}} 
\]
et $v_{i}$ est solution strictement positive de 
\[
\triangle _{\mathbf{g}}v_{i}+h_{t_{i}}.v_{i}=h^{\prime }v_{i}^{q_{i}-1}. 
\]
La suite ($v_{i})$ est born\'{e}e dans $H_{1}^{2}\,$et donc il existe $v\in
H_{1}^{2}\,$tel que 
\[
v_{i}\stackrel{H_{1}^{2}}{\rightharpoondown }v,\,\,v_{i}\stackrel{L^{2}}{%
\rightarrow }v\,\,\,et\,\,\,\,v_{i}\stackrel{L^{2^{*}-2}}{\rightarrow }v. 
\]
Ici encore se pr\'{e}sentent deux possibilit\'{e}s: $v\equiv 0$ ou $v>0$.

\textit{Deuxi\'{e}me \'{e}tape:}

Si $v>0$, comme nous l'avons expliqu\'{e} dans le principe de la
d\'{e}monstration, c'est termin\'{e}: \`{a} extraction pr\`{e}s, $v_{i}%
\stackrel{C^{2}}{\rightarrow }v$ et donc d'une part $F_{h^{\prime
}}(v_{i})\rightarrow F_{h^{\prime }}(v)$, et d'autre part 
\[
F_{h^{\prime }}(v_{i})=h_{t_{i}}+h^{\prime }(v_{i}^{\frac{4}{n-2}%
}-v_{i}^{q_{i}-2})\rightarrow h\, 
\]
soit $F_{h^{\prime }}(v)=h$ qui est critique pour $f$ et $\mathbf{g}$ avec
des solutions minimisantes. Donc $h^{\prime }$ est critique pour $f$ et $%
\mathbf{g}^{\prime }=v^{\frac{4}{n-2}}\mathbf{g}$, avec des solutions
minimisantes.

\textit{Toute la suite de la d\'{e}monstration traite donc du cas o\`{u} }$%
v\equiv 0$\textit{.}

\textit{Troisi\`{e}me \'{e}tape:} On montre qu'il y a un ph\'{e}nom\`{e}ne
de concentration.

a/: Montrons que: 
\[
0<c\leqslant \overline{\lim }\,\lambda _{i}\leqslant
K^{-2}(Sup_{M}\,h^{\prime })^{-\frac{n-2}{n}}\,\,. 
\]
Tout d'abord 
\[
\lambda _{i}=J_{i}(v_{i})=\frac{\int_{M}\left| \nabla v_{i}\right| ^{2}dv_{%
\mathbf{g}}+\int_{M}h_{t_{i}}.v_{i}{{}^{2}}dv_{\mathbf{g}}}{\left(
\int_{M}h^{\prime }v_{i}^{q_{i}}dv_{\mathbf{g}}\right) ^{\frac{2}{q_{i}}}}%
\geqslant \frac{\int_{M}\left| \nabla v_{i}\right| ^{2}dv_{\mathbf{g}%
}+\int_{M}h_{t_{i}}.v_{i}{{}^{2}}dv_{\mathbf{g}}}{(Sup\,h^{\prime })^{\frac{2%
}{2^{*}}}\left( \int_{M}v_{i}^{2^{*}}dv_{\mathbf{g}}\right) ^{\frac{2}{2^{*}}%
}(Vol_{g}(M))^{1-\frac{q_{i}}{2^{*}}}} 
\]
Or $h_{t_{i}}\stackrel{C^{2}}{\rightarrow }h$ avec $\triangle _{\mathbf{g}%
}+h $ coercif. Donc, en utilisant l'in\'{e}galit\'{e} de Sobolev, il existe $%
c,c^{\prime }>0$ ind\'{e}pendants de $i$ tels que 
\[
\int_{M}\left| \nabla v_{i}\right| ^{2}dv_{\mathbf{g}}+%
\int_{M}h_{t_{i}}.v_{i}{{}^{2}}dv_{\mathbf{g}}\geqslant c.\left\|
v_{i}\right\| _{H_{1}^{2}}^{2}\geqslant c^{\prime }.\left(
\int_{M}v_{i}^{2^{*}}dv_{\mathbf{g}}\right) ^{\frac{2}{2^{*}}}\,\,. 
\]
D'o\`{u} 
\[
\overline{\lim }\,\lambda _{i}\geqslant c^{\prime \prime }>0. 
\]
Par ailleurs, $\forall \varepsilon >0,\,$il existe une fonction $\phi >0$
telle que 
\[
J_{h,h^{\prime },\mathbf{g,}2^{*}}(\phi )\leqslant K^{-2}(Sup\,h^{\prime
})^{-\frac{2}{2^{*}}}+\varepsilon \,\,. 
\]
Or $\lim J_{i}(\phi )=J_{h,h^{\prime },\mathbf{g,}2^{*}}(\phi )$, donc 
\[
\overline{\lim }\,\lambda _{i}\leqslant K^{-2}(Sup_{M}\,h^{\prime })^{-\frac{%
n-2}{n}} 
\]
On peut donc, quitte \`{a} extraire, supposer que $\lambda _{i}\rightarrow
\lambda >0$. Alors 
\[
\lambda \leqslant K^{-2}(Sup_{M}\,h^{\prime })^{-\frac{n-2}{n}} 
\]
qui est une hypoth\`{e}se ``d'Energie Minimale'' (voir chapitre 3).

b/: Montrons que 
\[
0<\lambda ^{\frac{n}{2}}(Sup_{M}\,h^{\prime })^{-1}\leqslant \overline{\lim }%
\,\int_{M}v_{i}^{q_{i}}dv_{\mathbf{g}}\leqslant K^{2^{*}}\lambda ^{\frac{%
n2^{*}}{4}}\leqslant K^{-n}(Sup_{M}\,h^{\prime })^{-\frac{n}{2}} 
\]
On a $\int h^{\prime }v_{i}^{q_{i}}=\lambda _{i}^{\frac{q_{i}}{q_{i}-2}}$,
donc 
\[
J_{i}(v_{i})=\lambda _{i}=\frac{\int_{M}\left| \nabla v_{i}\right| ^{2}dv_{%
\mathbf{g}}+\int_{M}h_{t_{i}}.v_{i}{{}^{2}}dv_{\mathbf{g}}}{\left( \lambda
_{i}^{\frac{q_{i}}{q_{i}-2}}\right) ^{\frac{2}{q_{i}}}} 
\]
d'o\`{u} 
\[
\int_{M}\left| \nabla v_{i}\right| ^{2}dv_{\mathbf{g}}+%
\int_{M}h_{t_{i}}.v_{i}{{}^{2}}dv_{\mathbf{g}}=\lambda _{i}^{1+\frac{2}{%
q_{i}-2}}=\lambda _{i}^{\frac{q_{i}}{q_{i}-2}}\rightarrow \lambda ^{\frac{n}{%
2}} 
\]
et puisque $\int_{M}h_{t_{i}}.v_{i}{{}^{2}}dv_{\mathbf{g}}\rightarrow 0$, 
\[
\int_{M}\left| \nabla v_{i}\right| ^{2}dv_{\mathbf{g}}\rightarrow \lambda ^{%
\frac{n}{2}}\,\,. 
\]
On \'{e}crit alors (les int\'{e}grales \'{e}tant prises pour la mesure $dv_{%
\mathbf{g}})$: 
\begin{eqnarray*}
\int_{M}\left| \nabla v_{i}\right| ^{2}+\int_{M}h_{t_{i}}.v_{i}{{}^{2}}
&&\leqslant \lambda _{i}(Sup_{M}\,h^{\prime })^{\frac{2}{q_{i}}%
}(\int_{M}v_{i}^{q_{i}})^{\frac{2}{q_{i}}}\\
&&\leqslant \lambda _{i}(Sup_{M}\,h^{\prime })^{\frac{2}{q_{i}}}(Vol_{g}(M))^{1-\frac{%
q_{i}}{2^{*}}}(K^{2}\int_{M}\left| \nabla v_{i}\right|
^{2}+B_{0}\int_{M}v_{i}{{}^{2})}
\end{eqnarray*}
d'o\`{u} quand $i\rightarrow +\infty :$%
\[
\lambda ^{\frac{n}{2}}{\leqslant }\mu (Sup_{M}\,h^{\prime })^{\frac{2}{2^{*}}%
}(\overline{\lim }\int_{M}v_{i}^{q_{i}})^{\frac{2}{2^{*}}}\leqslant \lambda
(Sup_{M}\,h^{\prime })^{\frac{2}{2^{*}}}K^{2}\lambda ^{\frac{n}{2}} 
\]
soit 
\[
0<\lambda ^{\frac{n}{2}}(Sup_{M}\,h^{\prime })^{-1}\leqslant \overline{\lim }%
\,\int_{M}v_{i}^{q_{i}}dv_{\mathbf{g}}\leqslant K^{2^{*}}\lambda ^{\frac{%
n2^{*}}{4}}\leqslant K^{-n}(Sup_{M}\,h^{\prime })^{-\frac{n}{2}}\,\,. 
\]
Remarque: $K^{2^{*}}\lambda ^{\frac{n2^{*}}{4}}\leqslant
K^{-n}(Sup_{M}\,h^{\prime })^{-\frac{n}{2}}$ et si $f=cste>0$ alors $\lambda
=K^{-2}(Sup_{M}\,h^{\prime })^{-\frac{n-2}{n}}$.

c/: On dira que $x\in M$ est un point de concentration si 
\[
\forall r>0:\,\overline{\lim }\,\int_{B(x,r)}v_{i}^{q_{i}}dv_{\mathbf{g}}>0. 
\]

Reprenons rapidement les m\'{e}thodes des chapitres 2 et 3 pour l'\'{e}tude
du ph\'{e}nom\`{e}ne de concentration.

Tout d'abord, $M$ \'{e}tant compact, il existe au moins un point de
concentration $x\in M$.

Soit $\eta $ une fonction cut-off \`{a} support dans $B(x,r)\,$pour un $r>0$%
. Le principe d'it\'{e}ration vu au chapitre 2 nous donne:

\[
Q(i,k,\eta ).(\int_{M}(\eta v_{i}^{\frac{k+1}{2}})^{q_{i}})^{\frac{2}{q_{i}}%
}\leqslant
C\int_{B(x,r)}v_{i}^{k+1}
\]
o\`{u}
\[
\text{ }Q(i,k,\eta )=\frac{4k}{(k+1){{}^{2}}}-(Vol_{g}(M))^{\frac{q_{i}}{%
2^{*}}-1}K{{}^{2}}(\stackunder{B(x,r)}{Sup}\left| h^{\prime }\right|
).(\int_{B(x,r)}v_{i}^{q_{i}})^{\frac{q_{i}-2}{q_{i}}}
\]
Si $\overline{\lim }\,\int_{B(x,r)}v_{i}^{q_{i}}dv_{\mathbf{g}%
}<K^{-n}(Sup_{M}\,h^{\prime })^{-\frac{n}{2}}$ alors pour $k$ assez proche
de 1 et $i$ assez grand $Q(i,k,\eta )>Q>0$, par cons\'{e}quent: 
\[
(\int_{M}(\eta v_{i}^{\frac{k+1}{2}})^{q_{i}})^{\frac{2}{q_{i}}}\leq
C^{\prime }\int_{B(x,r)}v_{i}^{k+1}\,\,\,\, 
\]
et $\int_{B(x,r)}v_{i}^{k+1}\rightarrow 0$, donc 
\[
\int_{M}(\eta v_{i}^{\frac{k+1}{2}})^{q_{i}}\rightarrow 0. 
\]
Or par l'in\'{e}galit\'{e} de H\"{o}lder: 
\[
\int_{B(x,r/2)}v_{i}^{q_{i}}dv_{\mathbf{g}}\leqslant
C\int_{B(x,r/2)}v_{i}^{q_{i}\frac{k+1}{2}}dv_{\mathbf{g}}\,\,. 
\]
Mais si $x$ est un point de concentration, $\overline{\lim }%
\int_{B(x,r/2)}v_{i}^{q_{i}}>0$, d'o\`{u} une contradiction, et donc 
\[
\overline{\lim }\,\int_{B(x,r)}v_{i}^{q_{i}}dv_{\mathbf{g}}\geqslant
K^{-n}(Sup_{M}\,h^{\prime })^{-\frac{n}{2}}. 
\]

d/: Par cons\'{e}quent, en reprenant les m\^{e}mes m\'{e}thodes qu'au
d\'{e}but du chapitre 3:

1/: $\overline{\lim }\,\int_{B(x,r)}v_{i}^{q_{i}}dv_{\mathbf{g}%
}=K^{-n}(Sup_{M}\,h^{\prime })^{-\frac{n}{2}}$ , $\forall r>0$

2/: $x$ est l'unique point de concentration, not\'{e} $x_{0}$

3/: $\lambda =K^{-2}(Sup_{M}\,h^{\prime })^{-\frac{n-2}{n}}$

4/: $x_{0}$ est un point de maximum de $h^{\prime }$

5/: $v_{i}\rightarrow 0$ dans $C_{loc}^{2}(M-\{x_{0}\})$

Remarque: jusqu'\`{a} pr\'{e}sent, nous avons seulement utilis\'{e} le fait
que $h_{t}\rightarrow h$ avec ($h_{t},f,\mathbf{g})$ sous-critique et $%
\left( h,f,\mathbf{g}\right) $ critique.

\textit{Quatri\`{e}me \'{e}tape:}

On sait donc que la suite $(v_{i})$ se concentre en $x_{0}$ et que pour tout 
$i$ $F_{h^{\prime }}^{\prime }(v_{i})$ est sous-critique pour $f$ et $\textbf{g}$. On
voudrait trouver un $v_{i_{0}}$ , une fonction $v>0$ et un chemin continu de 
$v_{i_{0}}$ \`{a} $v$ tel que $F_{h^{\prime }}(v)$ soit faiblement critique
pour $f$ et $\textbf{g}$ et tel que 
\[
\frac{4(n-1)}{n-2}F_{h^{\prime }}(v)(P)>S_{\mathbf{g}}(P)-\frac{n-4}{2}\frac{\bigtriangleup _{%
\mathbf{g}}f(P)}{f(P)} 
\]
en tout point $P$\ o\`{u} $f$\ est maximum sur $M$\ . Alors le
th\'{e}or\`{e}me 1 nous dira que sur le chemin $u_{t}$ de $v_{i_{0}}$ \`{a} $%
v$ il existe un $u_{t}$ tel que $F_{h^{\prime }}(u_{t})$ est critique pour $%
f $ et $\mathbf{g}$ .

C'est maintenant que nous allons utiliser l'existence de fonctions critiques
strictement positives. Ainsi, comme nous l'avons dit au d\'{e}but de la
d\'{e}monstration, on peut supposer que les fonctions $h$ et $h_{t}$ sont
strictement positives.

Soit $s\geqslant 1$ et soit $v$ une fonction strictement positive. Alors 
\[
\bigtriangleup _{\mathbf{g}}(v^{s})=sv^{s-1}\bigtriangleup _{\mathbf{g}%
}v-s(s-1)v^{s-2}\left| \nabla v\right| _{\mathbf{g}}^{2}\,\,. 
\]
D'o\`{u} 
\[
F_{h^{\prime }}(v_{i}^{s})=h^{\prime }v_{i}^{s\frac{4}{n-2}%
}+sh_{t_{i}}-sh^{\prime }v_{i}^{q_{i}-2}+s(s-1)\frac{\left| \nabla v\right|
_{\mathbf{g}}^{2}}{v_{i}^{2}} 
\]
et donc 
\[
F_{h^{\prime }}(v_{i}^{s})\geqslant sh_{t_{i}}+h^{\prime }(v_{i}^{s\frac{4}{%
n-2}}-sv_{i}^{q_{i}-2})\,\,. 
\]
Maintenant:

Sur $\{x\in M\,/\,h^{\prime }(x)\leqslant 0\}$ :

$v_{i}\rightarrow 0$ uniform\'{e}ment car $x_{0}$ est un point o\`{u} $%
h^{\prime }$ est maximum sur $M,$ donc $h^{\prime }(x_{0})>0$ car par
hypoth\`{e}se $\bigtriangleup _{\mathbf{g}}+h^{\prime }$ est coercif. De
plus si $s\geqslant 1$ alors $s\frac{4}{n-2}\geqslant q_{i}-2$. Par
cons\'{e}quent, pour $i$ assez grand 
\[
F_{h^{\prime }}(v_{i}^{s})\geqslant sh_{t_{i}}\text{ sur }\{x\in
M\,/\,h^{\prime }(x)\leqslant 0\} 
\]

Sur $\{x\in M\,/\,h^{\prime }(x)>0\}$ :

consid\'{e}rons la fonction d'une variable r\'{e}elle d\'{e}finie pour $%
x\geqslant 0$ par 
\[
\beta _{i,s}(x)=x^{s\frac{4}{n-2}}-sx^{q_{i}-2}=x^{q_{i}-2}(x^{s\frac{4}{n-2}%
-q_{i}+2}-s). 
\]
Une simple \'{e}tude de cette fonction montre que 
\[
pour\,\,\,x\geqslant 0:\beta _{i,s}(x)\geqslant -s. 
\]
Or 
\[
F_{h^{\prime }}(v_{i}^{s})\geqslant sh_{t_{i}}+h^{\prime }\beta
_{i,s}(v_{i}) 
\]
donc 
\[
F_{h^{\prime }}(v_{i}^{s})\geqslant sh_{t_{i}}-sh^{\prime }\,\,\,\text{sur}%
\,\,\,\{x\in M\,/\,h^{\prime }(x)>0\}. 
\]
On peut donc \'{e}crire : 
\[
F_{h^{\prime }}(v_{i}^{s})\geqslant s(h_{t_{i}}-\stackunder{M}{Sup}\,h^{\prime
})\text{ sur }\{x\in M\,/\,h^{\prime }(x)>0\}\text{.} 
\]
Maintenant on peut utiliser le travail effectu\'{e} au d\'{e}but de la
d\'{e}monstration, \`{a} savoir que, pour toute constante $C>0$, on peut
supposer au d\'{e}part: 
\[
h_{t}>C\,\,surM 
\]
et 
\[
\frac{4(n-1)}{n-2}h_{t}(P)-S_{\mathbf{g}}(P)+\frac{n-4}{2}\frac{\bigtriangleup _{\mathbf{g}%
}f(P)}{f(P)}>C\,\,\text{en tout point de maximum de }f\,. 
\]
Alors, premi\`{e}rement, si on suppose d\'{e}j\`{a} que $h>\stackunder{M}{Sup%
}\,h^{\prime }$ sur tout $M$, on voit que pour $i$ et $s$ assez grand : 
\begin{equation}
F_{h^{\prime }}(v_{i}^{s})\geqslant B_{0}(\mathbf{g})K(n,2)^{-2}  
\end{equation}
et donc $F_{h^{\prime }}(v_{i}^{s})$ est faiblement critique pour $f$ et $%
\mathbf{g}$ . De plus pour tout $t\in [1,s]$ on a de m\^{e}me 
\[
F_{h^{\prime }}(v_{i}^{t})\geqslant t(h_{t_{i}}-\stackunder{M}{Sup}\,h^{\prime
})\text{ }\geqslant h_{t_{i}}-\stackunder{M}{Sup}\,h^{\prime }>0 
\]
donc $\bigtriangleup _{\mathbf{g}}+F_{h^{\prime }}(v_{i}^{t})$ est coercif.

Ensuite si l'on suppose en plus que 
\[
\frac{4(n-1)}{n-2}h_{t}(P)-S_{\mathbf{g}}(P)+\frac{n-4}{2}\frac{\bigtriangleup _{\mathbf{g}%
}f(P)}{f(P)}>\frac{4(n-1)}{n-2}\stackunder{M}{Sup}\ h^{\prime }\,\,\text{en tout point de
maximum de }f 
\]
alors pour tout $t\in [1,s]$: 
\begin{equation}
\frac{4(n-1)}{n-2}F_{h^{\prime }}(v_{i}^{t})(P)>S_{\mathbf{g}}(P)-\frac{n-4}{2}\frac{%
\bigtriangleup _{\mathbf{g}}f(P)}{f(P)}\,\,\text{en tout point de maximum de 
}f 
\end{equation}
d\`{e}s que $i$ est assez grand.

On fixe alors un $i$ et un $s$ assez grands pour avoir (5.1) et (5.2) et on
consid\`{e}re 
\[
s_{0}=\inf \{t>1\,/\,F_{h^{\prime }}(v_{i}^{t})\text{ est faiblement
critique\}} 
\]

On applique alors le th\'{e}or\`{e}me 1 au chemin $t\in [1,s_{0}]\mapsto
F_{h^{\prime }}(v_{i}^{t})$ pour obtenir le fait que $F_{h^{\prime
}}(v_{i}^{s_{0}})$ est critique pour $f$ et $\mathbf{g}$, avec des solutions
minimisantes. Donc $h^{\prime }$ est critique pour $f$ et $\mathbf{g}%
^{\prime }=(v_{i}^{s_{0}})^{\frac{4}{n-2}}\mathbf{g}$ avec des solutions
minimisantes.

Ceci termine la d\'{e}monstration.

Remarque: C'est \`{a} l'issue de cette d\'{e}monstration que nous avons
introduit le concept de triplet critique pour \'{e}tendre celui de fonction
critique. M. Vaugon et E. Hebey avaient en effet introduit la notion de
fonction critique pour l'\'{e}tude de certaines questions li\'{e}es \`{a} la
meilleure seconde constante $B_{0}(\mathbf{g})$, questions li\'{e}es comme
nous l'avons vu au chapitre 1 aux \'{e}quations du type $\triangle _{\mathbf{%
g}}u+h.u=u^{\frac{n+2}{n-2}}$. Il \'{e}tait donc naturel pour eux de parler
de fonction critique ''tout court'', une fois une m\'{e}trique fix\'{e}e.
Notre premi\`{e}re g\'{e}n\'{e}ralisation pour l'\'{e}tude d'\'{e}quations
du type $\triangle _{\mathbf{g}}u+h.u=f.u^{\frac{n+2}{n-2}}$ nous amenait
donc naturellement \`{a} parler de fonction $h$ critique pour $f$. Mais
\`{a} l'issue du probl\`{e}me trait\'{e} dans ce chapitre, o\`{u} la
m\'{e}trique varie dans une classe conforme, et notre but \'{e}tant
l'\'{e}tude de l'existence de solution \`{a} ces \'{e}quations, il nous a
paru int\'{e}ressant pour mettre en valeur cette notion d'introduire la
notion de triplet critique.

\chapter{Triplet Critique 3}
\pagestyle{myheadings}\markboth{\textbf{Triplet critique 3.}}{\textbf{Triplet critique 3.}}
Soit $(M,\mathbf{g})$ une vari\'{e}t\'{e} Riemannienne compacte de dimension 
$n\geqslant 3$. Soit $h$ une fonction $C^{\infty }$ fix\'{e}e telle que $%
\triangle _{\mathbf{g}}+h$ soit coercif. La question est: \textit{peut-on
trouver }$f$\textit{\ telle que }$(h,f,\mathbf{g})$\textit{\ soit critique ?}

Faisons d'abord quelques remarques. Si $h\geqslant B_{0}(\mathbf{g})K(n,2){%
{}^{-2}}$, $h$ est faiblement critique pour toute fonction $f$, et donc il
ne peut exister de fonction $f$ telle que ($h,f,\mathbf{g})$\textit{\ }soit
sous-critique. Mais plus important est la remarque suivante:

\textit{\ Si il existe une fonction }$f$\textit{\ non constante telle que
\thinspace }$(h,f,\mathbf{g})\,$\textit{soit critique avec une solution
minimisante }$u$\textit{\ alors (}$h,1,\mathbf{g}$\textit{) est
sous-critique. }

En effet, comme nous l'avons vu au chapitre 1 dans les d\'{e}finitions, on
peut supposer que $Sup\,f=1$. Alors puisque $u>0$%
\[
J_{h,1}(u)<J_{h,f}(u)=\frac{1}{K(n,2)^{2}(Sup\,f)^{\frac{2}{2^{*}}}}=\frac{1%
}{K(n,2)^{2}} 
\]
et donc $h$ est sous-critique pour 1. Notons que d'après le chapitre 3,
si $(h,f,\mathbf{g})\,$est critique et si 
\[
\frac{4(n-1)}{n-2}h(P)>S_{\mathbf{g}}(P)-\frac{n-4}{2}\frac{\bigtriangleup _{%
\mathbf{g}}f(P)}{f(P)} 
\]
\textit{\ }en tout point $P$\ o\`{u} $f$\ est maximum sur $M$, \textit{\ }$f$
v\'{e}rifiant (\textbf{H}$_{f}$), alors $(h,f,\mathbf{g})$ a des solutions
minimisantes.

Nous allons montrer le r\'{e}sultat suivant:

\textbf{Th\'{e}or\`{e}me 4:}

\textit{Soient donn\'{e}es la vari\'{e}t\'{e} (}$M,\mathbf{g}$\textit{), }$%
\dim M\geqslant 5$\textit{, et une fonction }$h$\textit{\ telle que
l'op\'{e}rateur }$\triangle _{\mathbf{g}}+h$\textit{\ soit coercif. Alors,
il existe une fonction }$f$ (\textit{v\'{e}rifiant }\textbf{(H}$_{f}$))%
\textit{\ telle que }$(h,f,\mathbf{g})$\textit{\ soit critique avec des
solutions minimisantes, si et seulement si }$(h,1,\mathbf{g})$\textit{\ est
sous-critique (où }$1$ \textit{est la fonction constante }$1$).

Nous venons de voir une des deux implications. Nous supposons donc
maintenant que $(h,1,\mathbf{g})$\textit{\ est sous-critique.}

La d\'{e}monstration se fera en deux \'{e}tapes:

Premi\`{e}rement: on montre qu'il existe une fonction $f\in C^{\infty }(M)$
v\'{e}rifiant \textbf{(H}$_{f}$), avec $\stackunder{M}{Sup}f=1$, et avec de
plus $\triangle _{\mathbf{g}}f$ aussi grand que l'on veut aux points de
maximum de $f$, telle que $(h,f,\mathbf{g})$ est faiblement critique.

Deuxi\`{e}mement: \'{e}tant donn\'{e}e cette fonction $f$, on montre qu'il
existe sur le chemin 
\[
t\rightarrow f_{t}=t.1+(1-t)f 
\]
une fonction pour laquelle $h$ est critique; autrement dit, il existe une
constante $c$ telle que ($h,f+c,\mathbf{g})$ soit critique.

\textbf{Premi\`{e}re \'{e}tape:}

On proc\`{e}de par contradiction. Supposons que pour toute fonction $f\in
C^{\infty }(M)$ v\'{e}rifiant $\stackunder{M}{Sup}f>0$, $(h,f,\mathbf{g})$
soit sous-critique. Alors, pour toutes ces fonctions il existe une solution
strictement positive $u$ \`{a} l'\'{e}quation 
\[
\triangle _{\mathbf{g}}u+h.u=\lambda .f.u^{\frac{n+2}{n-2}} 
\]
o\`{u} 
\[
\lambda =\stackunder{w\in \mathcal{H}_{f}}{\inf }I_{h,\mathbf{g}%
}(w)\,\,\,et\,\,\,\int_{M}f.u^{\frac{2n}{n-2}}dv_{\mathbf{g}}=1\,. 
\]
On rappelle que 
\[
I_{h,\mathbf{g}}(w)=\int_{M}\left| \nabla w\right| ^{2}dv_{\mathbf{g}%
}+\int_{M}h.w{{}^{2}}dv_{\mathbf{g}} 
\]
et 
\[
\mathcal{H}_{f}=\{w\in H_{1}^{2}\,/\,\int_{M}f.w^{\frac{2n}{n-2}}dv_{\mathbf{%
g}}=1\} 
\]

La m\'{e}trique $\mathbf{g}$ \'{e}tant fix\'{e}e, on sous-entendra $dv_{%
\mathbf{g}}$ dans les int\'{e}grales.

L'id\'{e}e est de construire une famille particuli\`{e}re de fonction $f_{t}$
dont le laplacien tend vers l'infini aux points de maximum pour obtenir une
contradiction. L'une de ces fonctions devra donc donner un triplet $(h,f_{t},%
\mathbf{g})$ faiblement critique. De plus le raisonnement qui suit \'{e}tant
valable pour n'importe quelle sous-partie de la famille que nous allons
construire, cette fonction aura un laplacien aussi grand que l'on veut en
ses points de maximum.

Dans $\Bbb{R}^{n}$, on construit pour $t\rightarrow 0$ une famille ($P_{t})$
de fonctions $C^{\infty }$, analogue \`{a} une suite r\'{e}gularisante,
telles que 
\begin{eqnarray*}
0 &\leqslant &P_{t}\leqslant 1 \\
\,\,P_{t}(x) &=&P_{t}(\left| x\right| ) \\
\,\,P_{t}(0) &=&1 \\
\left\| \nabla P_{t}\right\| &\sim &\frac{c_{1}}{t}\,\,\,sur\,\,\,B(0,t) \\
\left| \bigtriangleup P_{t}(0)\right| &\sim &\frac{c_{2}}{t^{2}} \\
SuppP_{t} &=&B(0,t).
\end{eqnarray*}
Soit maintenant $x_{0}$ un point de $M$ tel que $h(x_{0})>0$; ce point
existe car $\triangle _{\mathbf{g}}+h$ est coercif. On pose 
\[
f_{t}=P_{t}\circ \exp _{x_{0}}^{-1} 
\]
On suppose donc que pour tout $t:$ $(h,f_{t},\mathbf{g})$ est sous-critique
et on cherche une contradiction. Pour tout $t$ on a une solution $u_{t}>0$
de 
\[
(E_{t}):\,\,\bigtriangleup _{\mathbf{g}}u_{t}+h.u_{t}=\lambda
_{t}.f_{t}.u_{t}^{\frac{n+2}{n-2}} 
\]
avec $\int f_{t}u_{t}^{2^{*}}dv_{\mathbf{g}}=1$ et 
\[
\lambda _{t}<K^{-2}(\stackunder{M}{Sup}f_{t})^{-\frac{n-2}{n}}=K^{-2}. 
\]
Alors, ($u_{t})$ est born\'{e}e dans $H_{1}^{2}(M)$ quand $t\rightarrow 0$.
Donc ($u_{t})$ est born\'{e}e dans $L^{2^{*}}$ et ($u_{t}^{2^{*}-1})$ est
born\'{e}e dans $L^{\frac{2^{*}}{2^{*}-1}}$. Après extraction, si $f_{t}%
\stackrel{L^{2}}{\rightarrow }f$ et $u_{t}\stackrel{L^{2}}{\rightarrow }u$,
alors 
\[
f_{t}u_{t}^{2^{*}-1}\rightharpoondown fu^{2^{*}-1}. 
\]
Or ici, $f_{t}\stackrel{L^{p}}{\rightarrow }0$, donc l'\'{e}quation ($E_{t})$
``converge'' vers 
\[
\bigtriangleup _{\mathbf{g}}u+h.u=0 
\]
au sens o\`{u} $u$ est solution de cette \'{e}quation. Mais $\bigtriangleup
_{\mathbf{g}}+h$ est coercif, donc $u=0$, i.e. $u_{t}\rightarrow 0$ dans $%
L^{p}$ pour $p<2^{*}$.

La suite ($u_{t})$ d\'{e}veloppe donc un ph\'{e}nom\`{e}ne de concentration
tel qu'expos\'{e} au chapitre 3. Au chapitre 3, la fonction $f$ du second
membre ne variait pas et c'\'{e}tait au premier membre que nous avions une
suite de fonctions ($h_{t})$; ici c'est l'inverse. N\'{e}anmoins, les
r\'{e}sultats obtenus restent valables, seul le changement d'\'{e}chelle
n\'{e}cessaire pour les estim\'{e}es faibles demande une attention
particuli\`{e}re. Nous allons reprendre rapidement les \'{e}tapes de cette
\'{e}tude, ne nous arr\^{e}tant que sur les points nouveaux.

\textit{a/: Il existe, \`{a} extraction pr\`{e}s d'une sous-famille de }$%
(u_{t})$\textit{, un unique point de concentration et
c'est le point $x_{0}$ o\`{u} les }$f_{t}$\textit{\ sont maximums sur }$M.$\textit{\ De plus }

\[
\forall \delta >0\mathit{,\ }\overline{\stackunder{t\rightarrow 1}{\lim }}%
\int_{B(x_{0},\delta )}f_{t}u_{t}^{2^{*}}=1\,. 
\]

Il suffit de reprendre exactement la m\'{e}thode utilis\'{e}e au chapitre 2,
bas\'{e}e sur le principe d'it\'{e}ration. Notons que plus
pr\'{e}cis\'{e}ment, puisque le support $Supp\,f_{t}=B(x_{0},t)$, on a pour
tout $\delta >0$ et dès que $t<\delta $:
\[
\int_{B(x_{0},\delta )}f_{t}u_{t}^{2^{*}}=1. 
\]
Par ailleurs, on peut supposer que 
\[
\lambda _{t}\rightarrow \lambda =K^{-2}(\stackunder{M}{Sup}f_{t})^{-\frac{n-2%
}{n}}=K^{-2}. 
\]

\textit{b/: }$u_{t}\rightarrow 0$\textit{\ dans }$C_{loc}^{0}(M-\{x_{0}\})$

On peut reprendre exactement la d\'{e}monstration du chapitre 3

\textit{c/: Estim\'{e}es ponctuelles faibles}

Reprenons les notations du changement d'\'{e}chelle: on consid\`{e}re une
suite de points $(x_{t})$ telle que $m_{t}=\stackunder{M}{Max}\,u_{t}=u_{t}(x_{t}):=\mu _{t}^{-\frac{n-2}{2}}$. 

D'apr\`{e}s ce qui pr\'{e}c\`{e}de $x_{t}\rightarrow x_{0}$ et $\mu
_{t}\rightarrow 0$. Rappelons que $\overline{u}_{t},\overline{f}_{t},%
\overline{h}_{t},\mathbf{g}_{t}$ d\'{e}signent les fonctions et la
m\'{e}trique ``lues'' dans la carte $\exp _{x_{t}}^{-1}$, et $\,\,\widetilde{%
u}_{t}\,,\,\widetilde{h}_{t}\,,\,\widetilde{f}_{t},\widetilde{\mathbf{g}}%
_{t} $ d\'{e}signent les fonctions et la m\'{e}trique ``lues'' après
blow-up de centre $x_{t}$ et de coefficient $k_{t}=\mu _{t}^{-1}$.

En regardant attentivement la d\'{e}monstration des estim\'{e}es ponctuelles
faibles, on constate qu'elle fonctionnera dans la situation de ce chapitre
d\`{e}s qu'on aura obtenu :

\[
\forall R>0\,\,\,:\stackunder{t\rightarrow 0}{\lim }\int_{B(x_{t},R\mu
_{t})}f_{t}u_{t}^{2^{*}}dv_{\mathbf{g}}=1-\varepsilon _{R}\,\,\,o\grave{u}%
\,\,\,\varepsilon _{R}\stackunder{R\rightarrow +\infty }{\rightarrow }0\,. 
\]
Cette relation s'obtient elle-m\^{e}me \`{a} partir du changement
d'\'{e}chelle (blow-up) une fois que l'on a montr\'{e} que $\,\widetilde{u}%
_{t}\stackrel{C_{loc}^{2}(\Bbb{R}^{n})}{\rightarrow }\,\widetilde{u}$ o\`{u} 
$\widetilde{u}$ est solution de 
\[
\bigtriangleup _{e}\widetilde{u}=K^{-2}\widetilde{u}^{\frac{n+2}{n-2}}. 
\]
C'est l\`{a} que se pr\'{e}sente une difficult\'{e} due au fait que l'on a
une famille de fonctions ($f_{t}$). En effet, dans le changement
d'\'{e}chelle, l'\'{e}quation 
\[
(E_{t}):\,\,\bigtriangleup _{\mathbf{g}}u_{t}+h.u_{t}=\lambda
_{t}.f_{t}.u_{t}^{\frac{n+2}{n-2}} 
\]
devient 
\[
(\widetilde{E}_{t})\,:\,\triangle _{\widetilde{\mathbf{g}}_{t}}\widetilde{u}%
_{t}+\mu _{t}^{2}.\widetilde{h}_{t}.\widetilde{u}_{t}=\lambda _{t}\widetilde{%
f}_{t}.\widetilde{u}_{t}^{\frac{n+2}{n-2}} 
\]
et pour obtenir que cette derni\`{e}re \'{e}quation ``converge'' vers 
\[
\bigtriangleup _{e}\widetilde{u}=K^{-2}\widetilde{u}^{\frac{n+2}{n-2}} 
\]
il suffit de montrer que $\,\widetilde{f}_{t}$ converge simplement vers 1
(ce qui est \'{e}vident lorsqu'on a une fonction constante $f$ au second
membre de ($E_{h,f,\mathbf{g}}$)). Comme la suite ($\widetilde{f}_{t}$) est
uniform\'{e}ment born\'{e}e par 1 sur $\Bbb{R}^{n}$ (on admet qu'on prolonge
les $\widetilde{f}_{t}$ par 0 sur $\Bbb{R}^{n}\backslash B(0,\delta \mu
_{t}^{-1})$), on sait d\'{e}j\`{a} par le th\'{e}or\`{e}me 8.25 de
Gilbard-Trudinger $\left[ 16\right] $ et par le th\'{e}or\`{e}me d'Ascoli
qu'il existe une fonction $\widetilde{u}\in C^{0}(\Bbb{R}^{n})$ telle que,
\`{a} extraction pr\`{e}s, $\widetilde{u}_{t}\stackrel{C_{loc}^{0}(\Bbb{R}%
^{n})}{\rightarrow }\,\widetilde{u}$, avec bien s\^{u}r $\widetilde{u}(0)=1$.

Nous allons montrer que $\widetilde{f}_{t}\stackrel{p.p}{\rightarrow }1$ sur 
$\Bbb{R}^{n}$ en deux \'{e}tapes (nous allons m\^{e}me montrer un peu plus):

1/: Il existe $\widetilde{f}\in L_{loc}^{2}(\Bbb{R}^{n})$ telle que $%
\widetilde{f}_{t}\stackrel{p.p}{\rightarrow }\widetilde{f}$ sur $\Bbb{R}^{n}$

2/: $\widetilde{f}=1$ p.p sur $\Bbb{R}^{n}$

On rappelle qu'on prolonge les $\widetilde{f}_{t}$ ,d\'{e}finies \textit{a
priori} sur $B(0,\delta \mu _{t}^{-1})$, par 0 sur $\Bbb{R}^{n}\backslash
B(0,\delta \mu _{t}^{-1}).$

Premi\`{e}rement:

On a $\widetilde{f}_{t}(x)=\overline{f}_{t}(\mu _{t}x)$ et $\left| \nabla 
\overline{f}_{t}\right| \leq \frac{c}{t}$. Donc 
\[
\left| \nabla \widetilde{f}_{t}\right| \leq c.\frac{\mu _{t}}{t}\,. 
\]

On distingue alors deux cas

a/: Si ($\frac{\mu _{t}}{t}$) est born\'{e}e: Alors pour tout compact $%
K\subset \subset \Bbb{R}^{n}$, ($\widetilde{f}_{t}$) est born\'{e}e dans $%
H_{1}^{n+1}(K)$ (o\`{u} $n=\dim M)$. Donc, par compacit\'{e} de l'inclusion $%
H_{1}^{n+1}(K)\subset C^{0,\alpha }(K)$ pour un $\alpha >0$, apr\`{e}s
extraction, il existe $\widetilde{f}_{K}\in C^{0,\alpha }(K)$ telle que 
\[
\widetilde{f}_{t}\stackrel{C^{0,\alpha }(K)}{\rightarrow }\widetilde{f}_{K} 
\]
Par extraction diagonale, on construit alors $\widetilde{f}\in C^{0,\alpha }(%
\Bbb{R}^{n})$ telle que 
\[
\widetilde{f}_{t}\stackrel{C^{0,\alpha }(K^{\prime })}{\rightarrow }%
\widetilde{f} 
\]
pour tout compact $K^{\prime }$ de $\Bbb{R}^{n}$, et de plus $\widetilde{f}%
\in H_{1,loc}^{n+1}(\Bbb{R}^{n}).$ En particulier $\widetilde{f}_{t}%
\stackrel{p.p}{\rightarrow }\widetilde{f}\,\,\,sur\,\,\,\Bbb{R}^{n}\,.$

b/: Si $\frac{\mu _{t}}{t}\rightarrow +\infty $ : le support de $\widetilde{f%
}_{t}$ est 
\[
Supp\widetilde{f}_{t}=B(\frac{x_{0}(t)}{\mu _{t}},\frac{t}{\mu _{t}}), 
\]
o\`{u} $x_{0}(t)=\exp _{x_{t}}^{-1}(x_{0})$.

Si ($\frac{\left| x_{0}(t)\right| }{\mu _{t}}$) est born\'{e}e, on peut
extraire une sous-suite pour que 
\[
\frac{x_{0}(t)}{\mu _{t}}\rightarrow P\in \Bbb{R}^{n}; 
\]
et alors 
\[
\widetilde{f}_{t}\stackrel{C_{loc}^{0}(\Bbb{R}^{n}-\{P\})}{\rightarrow }0 
\]

Si $\frac{\left| x_{0}(t)\right| }{\mu _{t}}\rightarrow \infty $, alors 
\[
\widetilde{f}_{t}\stackrel{C_{loc}^{0}(\Bbb{R}^{n})}{\rightarrow }0 
\]
Dans les deux cas, $\widetilde{f}_{t}\stackrel{p.p}{\rightarrow }%
0\,\,\,sur\,\,\,\Bbb{R}^{n}.$

Dans le cas a/, $\widetilde{u}$ est solution faible de 
\[
\bigtriangleup _{e}\widetilde{u}=K^{-2}\widetilde{f}\widetilde{u}^{\frac{n+2%
}{n-2}}\,\, 
\]
avec $\widetilde{f}\geq 0$ puisque $\widetilde{f}_{t}\geq 0$, et $\widetilde{%
f}\in H_{1,loc}^{n+1}(\Bbb{R}^{n})\subset C^{0,\alpha }(\Bbb{R}^{n}).\,$

Dans le cas b/, $\,\widetilde{u}$ est solution faible de 
\[
\bigtriangleup _{e}\widetilde{u}=0.\,\, 
\]
Dans les deux cas les th\'{e}ories elliptiques et les th\'{e}or\`{e}mes de
r\'{e}gularit\'{e} standard nous donnent la régularité $C^{2}$ de $\widetilde{u}$
, et donc $\bigtriangleup _{e}\widetilde{u}\geq 0$. Le principe du
maximum nous dit alors que soit $\widetilde{u}\equiv 0$ soit $\widetilde{u}%
>0 $. Or $\widetilde{u}(0)=1$ donc $\widetilde{u}>0$. (Rappelons que par
convention $\bigtriangleup _{e}=-\sum_{i}\partial _{ii}^{2}$ ).

Deuxi\`{e}mement:

On commence par utiliser le principe d'it\'{e}ration: pour une fonction
cut-off $\eta $ valant 1 autour de $x_{0}$, on multiplie ($E_{t}$) par $\eta
^{2}u_{t}$, on int\`{e}gre et on utilise l'in\'{e}galit\'{e} de Sobolev pour
obtenir, en se souvenant que $\lambda _{t}<K^{-2}(\stackunder{M}{Sup}%
f_{t})^{-\frac{n-2}{n}}$ et que $Sup\,f=1$:
\[
(\int_{M}(\eta u_{t})^{2^{*}})^{\frac{2}{2^{*}}}\leq \lambda
_{t}K^{2}\int_{M}\eta ^{2}f_{t}u_{t}^{2^{*}}+c\int_{Supp\,\eta
}u_{t}^{2}\,\,. 
\]
On prend $\eta =1$ sur $B(x_{0},\frac{3}{2}\delta )$ et $\eta =0$ sur $%
M\backslash B(x_{0},2\delta )$. Alors pour $t$ proche de 0 
\[
Supp\,f_{t}\subset B(x_{0},t)\subset B(x_{t},\delta )\subset B(x_{0},\frac{3%
}{2}\delta ) 
\]
D'o\`{u} 
\[
(\int_{B(x_{t},\delta )}u_{t}{}^{2^{*}})^{\frac{2}{2^{*}}}\leq
\int_{B(x_{t},\delta )}f_{t}u_{t}^{2^{*}}+c\int_{M}u_{t}^{2} 
\]
et apr\`{e}s changement d'\'{e}chelle 
\[
(\int_{B(0,\delta \mu _{t}^{-1})}\widetilde{u}_{t}^{2^{*}})^{\frac{2}{2^{*}}%
}\leq \int_{B(0,\delta \mu _{t}^{-1})}\widetilde{f}_{t}\widetilde{u}%
_{t}^{2^{*}}+c\int_{M}u_{t}^{2}=1+c\int_{M}u_{t}^{2}\,\,. 
\]
Or $\int_{M}u_{t}^{2}\rightarrow 0$ donc 
\[
\stackunder{t\rightarrow 0}{\overline{\lim }}\int_{B(0,\delta \mu _{t}^{-1})}%
\widetilde{u}_{t}^{2^{*}}\leq 1\,\,. 
\]
Par ailleurs, on sait que $\widetilde{f}_{t}\stackrel{p.p}{\rightarrow }%
\widetilde{f}$ avec $\widetilde{f}\leq 1$ et $\widetilde{u}_{t}(0)=1$.
Supposons alors qu'il existe $A\subset \Bbb{R}^{n}$ avec $mes(A)>0$ telle
que $\widetilde{f}<1$ sur $A$ et \'{e}crivons $\Bbb{R}^{n}=A\cup B$ avec $%
\widetilde{f}=1$ p.p sur $B$. Alors comme $\widetilde{f}_{t}\geq 0$ et comme 
$\widetilde{u}_{t}\stackrel{C^{2}}{\rightarrow }\widetilde{u}>0$ :
\begin{eqnarray*}
1=\int_{B(0,\delta \mu _{t}^{-1})}\widetilde{f}_{t}\widetilde{u}_{t}^{2^{*}}
& &\leqslant \stackunder{t\rightarrow 0}{\overline{\lim }}\int_{B(0,\delta
\mu _{t}^{-1})\cap A}\widetilde{f}_{t}\widetilde{u}_{t}^{2^{*}}+\stackunder{%
t\rightarrow 0}{\overline{\lim }}\int_{B(0,\delta \mu _{t}^{-1})\cap B}%
\widetilde{f}_{t}\widetilde{u}_{t}^{2^{*}} \\ 
&& <\stackunder{t\rightarrow 0}{\overline{\lim }}\int_{B(0,\delta \mu
_{t}^{-1})\cap A}\widetilde{u}_{t}^{2^{*}}+\stackunder{t\rightarrow 0}{%
\overline{\lim }}\int_{B(0,\delta \mu _{t}^{-1})\cap B}\widetilde{u}%
_{t}^{2^{*}} \\ 
&& =\stackunder{t\rightarrow 0}{\overline{\lim }}\int_{B(0,\delta \mu
_{t}^{-1})}\widetilde{u}_{t}^{2^{*}}
\end{eqnarray*}
soit 
\[
1<\stackunder{t\rightarrow 0}{\overline{\lim }}\int_{B(0,\delta \mu
_{t}^{-1})}\widetilde{u}_{t}^{2^{*}} 
\]
d'o\`{u} une contradiction, et finalement $\widetilde{f}_{t}\stackrel{p.p}{%
\rightarrow }1$ sur $\Bbb{R}^{n}$.

Ainsi, comme nous l'avons dit 
\[
(\widetilde{E}_{t})\,:\,\triangle _{\widetilde{\mathbf{g}}_{t}}\widetilde{u}%
_{t}+\mu _{t}^{2}.\widetilde{h}_{t}.\widetilde{u}_{t}=\lambda _{t}\widetilde{%
f}_{t}.\widetilde{u}_{t}^{\frac{n+2}{n-2}} 
\]
``converge'' vers 
\[
\bigtriangleup _{e}\widetilde{u}=K^{-2}\widetilde{u}^{\frac{n+2}{n-2}} 
\]
au sens o\`{u} 
\[
\,\widetilde{u}_{t}\stackrel{C_{loc}^{2}(\Bbb{R}^{n})}{\rightarrow }\,%
\widetilde{u} 
\]
o\`{u} $\widetilde{u}$ est solution de $\bigtriangleup _{e}\widetilde{u}%
=K^{-2}\widetilde{u}^{\frac{n+2}{n-2}}$. Comme $\widetilde{u}(0)=1$, 
\[
\widetilde{u}(x)=(1+\frac{K^{-2}}{n(n-2)}\left| x\right| ^{2})^{-\frac{n-2}{2%
}}\,\,. 
\]

Remarque: D'apr\`{e}s ce que nous avons vu plus haut, $\widetilde{f}_{t}%
\stackrel{p.p}{\rightarrow }1$ sur $\Bbb{R}^{n}$ entra\^{i}ne que $\frac{\mu
_{t}}{t}\rightarrow 0$. Cela peut s'interpr\'{e}ter ''intuitivement'' en
disant que les fonctions $u_{t}$ se concentrent plus vite que les fonctions $%
f_{t}$.

A partir de l\`{a}, on peut reprendre exactement les d\'{e}monstrations du
chapitre 3: ainsi

\[
\forall R>0:\stackunder{t\rightarrow 0}{\lim }\int_{B(x_{t},R\mu
_{t})}f_{t}u_{t}^{2^{*}}dv_{\mathbf{g}}=1-\varepsilon _{R}\,\,\,o\grave{u}%
\,\,\,\varepsilon _{R}\stackunder{R\rightarrow +\infty }{\rightarrow }0 
\]
puis

\[
\exists C>0\,\,\,tel\,\,\,que\,\,\,\forall x\in M:\,d_{\mathbf{g}}(x,x_{t})^{%
\frac{n-2}{2}}u_{t}(x)\leq C. 
\]
et

\[
\forall \varepsilon >0,\exists R>0\,\,\,tel\,\,\,que\,\,\,\forall
t,\,\forall x\in M:\,d_{\mathbf{g}}(x,x_{t})\geq R\mu _{t}\,\Rightarrow
\,\,d_{\mathbf{g}}(x,x_{t})^{\frac{n-2}{2}}u_{t}(x)\leq \varepsilon . 
\]

d/: On obtient sans changement ce que l'on appelle la concentration $L^{2}:$

Si $\dim M\geq 4,$ 
\[
\forall \delta >0\,:\,\stackunder{t\rightarrow 0}{\lim }\frac{%
\int_{B(x_{0},\delta )}u_{t}^{2}dv_{\mathbf{g}}}{\int_{M}u_{t}^{2}dv_{%
\mathbf{g}}}=1 
\]

e/: Il n'y a pas non plus de changements pour les estim\'{e}es fortes:

\[
\exists C>0\,\,\,tel\,\,\,que\,\,\,\forall x\in M:\,d_{\mathbf{g}%
}(x,x_{t})^{n-2}\mu _{t}^{-\frac{n-2}{2}}u_{t}(x)\leq C, 
\]
ni pour la concentration $L^{p}$\ forte:

\textbf{\ }$\forall R>0$ , $\forall \delta >0$ et $\forall p>\frac{n}{n-2}$
o\`{u} $n=\dim M$ on a : 
\[
\stackunder{t\rightarrow 0}{\lim }\frac{\int_{B(x_{t},R\mu
_{t})}u_{t}^{p}dv_{\mathbf{g}}}{\int_{B(x_{t},\delta )}u_{t}^{p}dv_{\mathbf{g%
}}}=1-\varepsilon _{R}\,\,\,o\grave{u}\,\,\,\varepsilon _{R}\stackunder{%
R\rightarrow +\infty }{\rightarrow }0 
\]

Venons-en alors \`{a} l'argument central, qui repose sur la m\'{e}thode
d\'{e}velopp\'{e}e au chapitre 3:

On ins\`{e}re dans l'in\'{e}galit\'{e} de Sobolev euclidienne l'\'{e}quation
($E_{t})$ lue dans la carte $\exp _{x_{t}}^{-1}$. Reprenant les calculs de
ce chapitre (nous invitons le lecteur \`{a} s'y reporter pour revoir les
d\'{e}tails): 
\begin{eqnarray*}
\int_{B(0,\delta )}\overline{h}_{t}(\eta \overline{u}_{t})^{2}dx\leq && \frac{%
1}{K(n,2){{}^{2}}(\stackunder{M}{Sup}f)^{\frac{n-2}{n}}}\int_{B(0,\delta )}%
\overline{f}_{t}\eta ^{2}\overline{u}_{t}^{2^{*}}dx\\
&&-\frac{1}{K(n,2){{}^{2}}}%
(\int_{B(0,\delta )}(\eta \overline{u}_{t})^{2^{*}}dx)^{\frac{2}{2^{*}}} \\ 
&&+C.\delta ^{-2}\int_{B(0,\delta )\backslash B(0,\delta /2)}\overline{u}%
_{t}^{2}dx+B_{t}+C_{t}
\end{eqnarray*}
avec:

$B_{t}=\frac{1}{2}\int_{B(0,\delta )}(\partial _{k}(\,\mathbf{g}%
\,_{t}^{ij}\Gamma (\,\mathbf{g}\,_{t})_{ij}^{k}+\partial _{ij}\,\mathbf{g}%
\,_{t}^{ij})(\eta \overline{u}_{t}^{2})dx$

$C_{t}=\left| \int_{B(0,\delta )}\eta ^{2}(\,\mathbf{g}\,_{t}^{ij}-\delta
^{ij})\partial _{i}\overline{u}_{t}\partial _{j}\overline{u}_{t}dx\right| $

$A_{t}=\frac{1}{K(n,2){{}^{2}}(\stackunder{M}{Sup}f)^{\frac{n-2}{n}}}%
\int_{B(0,\delta )}\overline{f}_{t}\eta ^{2}\overline{u}_{t}^{2^{*}}dx-\frac{%
1}{K(n,2){{}^{2}}}(\int_{B(0,\delta )}(\eta \overline{u}_{t})^{2^{*}}dx)^{%
\frac{2}{2^{*}}}$

On peut \'{e}crire 
\[
A_{t}\leq \frac{1}{K(n,2){{}^{2}}(\stackunder{M}{Sup}f_{t})^{\frac{n-2}{n}}}%
(A_{t}^{1}+A_{t}^{2}) 
\]
o\`{u} $\,A_{t}^{1}=(\int_{B(0,\delta )}\overline{f}_{t}(\eta \overline{u}%
_{t})^{2^{*}}dx)^{\frac{n-2}{n}}\,-(Supf_{t}.\int_{B(0,\delta )}(\eta 
\overline{u}_{t})^{2^{*}}dx)^{\frac{n-2}{n}}\,.$

D'apr\`{e}s les calculs du chapitre 3, 
\[
\stackunder{t\rightarrow 0}{\overline{\lim }}\frac{K(n,2){{}^{-2}}(%
\stackunder{M}{Sup}f_{t})^{-\frac{n-2}{n}}A_{t}^{2}+C.\delta
^{-2}\int_{B(0,\delta )\backslash B(0,\delta /2)}\overline{u}%
_{t}^{2}dx+B_{t}+C_{t}}{\int_{B(0,\delta )}\overline{u}_{t}^{2}dx}\leq \frac{%
n-2}{4(n-1)}S_{\mathbf{g}}(x_{0})+\varepsilon _{\delta } 
\]
o\`{u} $\varepsilon _{\delta }\rightarrow 0$ quand $\delta \rightarrow 0$.

Reste \`{a} traiter : 
\[
\stackunder{t\rightarrow 0}{\overline{\lim }}\frac{A_{t}^{1}}{%
\int_{B(0,\delta )}\overline{u}_{t}^{2}dx} 
\]
Par construction, la suite $f_{t}$ est d\'{e}croissante quand $t\rightarrow
0 $ au sens o\`{u} : 
\[
si\,\,\,t\leq t^{\prime }\,\,\,alors\,\,\,f_{t}\leq f_{t^{\prime }}\,\,. 
\]
Fixons un $t_{0}.$ Alors pour tout $t\leq t_{0}$%
\begin{eqnarray*}
\int_{B(0,\delta )}\overline{f}_{t}(\eta \overline{u}_{t})^{2^{*}}dx& & 
=\int_{B(x_{t},\delta )}f_{t}.(\eta \circ \exp
_{x_{t}}^{-1})^{2^{*}}.u_{t}^{2^{*}}.(\exp _{x_{t}}^{-1})^{*}dx \\ 
&& \leq \int_{B(x_{t},\delta )}f_{t_{0}}.(\eta \circ \exp
_{x_{t}}^{-1})^{2^{*}}.u_{t}^{2^{*}}.(\exp _{x_{t}}^{-1})^{*}dx \\ 
& &=\int_{B(0,\delta )}(f_{t_{0}}\circ \exp _{x_{t}})(\eta \overline{u}%
_{t})^{2^{*}}dx\,.
\end{eqnarray*}
Notons: 
\[
\overline{f}_{t_{0},t}=f_{t_{0}}\circ \exp _{x_{t}}\, 
\]
et 
\[
\widetilde{f}_{t_{0},t}=\overline{f}_{t_{0},t}\circ \psi _{\mu _{t}^{-1}}^{-1}\,. 
\]
Alors: 
\begin{eqnarray*}
A_{t}^{1} & \leq (\int_{B(0,\delta )}\overline{f}_{t_{0},t}(\eta \overline{u}%
_{t})^{2^{*}}dx)^{\frac{n-2}{n}}\,-(Supf_{t}.\int_{B(0,\delta )}(\eta 
\overline{u}_{t})^{2^{*}}dx)^{\frac{n-2}{n}} \\ 
& \leq (\int_{B(0,\delta )}\overline{f}_{t_{0},t}(\eta \overline{u}%
_{t})^{2^{*}}dx)^{\frac{n-2}{n}}\,-(Supf_{t_{0}}.\int_{B(0,\delta )}(\eta 
\overline{u}_{t})^{2^{*}}dx)^{\frac{n-2}{n}}
\end{eqnarray*}
puisque $Supf_{t}=Supf_{t_{0}}=1=f_{t_{0}}(x_{0})\,$pour tout $t$.

On obtient alors par la m\^{e}me m\'{e}thode qu'au chapitre 3 
\[
\stackunder{t\rightarrow 0}{\overline{\lim }}\frac{A_{t}^{1}}{%
\int_{B(0,\delta )}\overline{u}_{t}^{2}dv_{\mathbf{g}_{t}}}\leq -\frac{%
(n-2)(n-4)}{8(n-1)}\frac{\bigtriangleup _{\mathbf{g}}f_{t_{0}}(x_{0})\,}{%
f_{t_{0}}(x_{0})\,}+\varepsilon _{\delta } 
\]
et apr\`{e}s avoir fait tendre $\delta $ vers 0: 
\[
h(x_{0})\leq \frac{n-2}{4(n-1)}S_{\mathbf{g}}(x_{0})-\frac{(n-2)(n-4)}{8(n-1)%
}\frac{\bigtriangleup _{\mathbf{g}}f_{t_{0}}(x_{0})\,}{f_{t_{0}}(x_{0})\,} 
\]
Or 
\[
\bigtriangleup _{\mathbf{g}}f_{t}(x_{0})\,\sim +\frac{c}{t^{2}}\stackunder{%
t\rightarrow 0}{\rightarrow }+\infty 
\]
donc en prenant $t_{0}$ assez proche de 0 on obtient une contradiction.

En cons\'{e}quence, on peut trouver dans la famille ($f_{t}$) des fonctions
avec un laplacien en $x_{0}$: $\bigtriangleup _{\mathbf{g}}f_{t}(x_{0})\,$%
aussi grand que l'on veut telles que les \'{e}quations: $\bigtriangleup _{%
\mathbf{g}}u+h.u=f_{t}.u^{\frac{n+2}{n-2}}$ n'aient pas de solutions
minimisantes et donc telles que $h$ soit faiblement critique pour $f_{t}$ et 
$\mathbf{g}$.
\\

Remarque: cela s'applique \`{a} $h=cste<B_{0}K^{-2}$ ou \`{a} $h=S_{\mathbf{g%
}}$ si $M$ n'est pas la sph\`{e}re. (voir le chapitre 1)
\\

\textit{Deuxi\`{e}me \'{e}tape:}

Pour notre fonction $h$ telle que $(h,1,\mathbf{g})$ soit sous-critique,
nous savons maintenant qu' il existe une fonction $f$, au laplacien aussi
grand que l'on veut en ses points de maximum, telle que $(h,f,\mathbf{g})$
soit faiblement critique. Plus pr\'{e}cis\'{e}ment, on a donc trouv\'{e} une
fonction $f$ telle que:

1/: $(h,f,\mathbf{g})$ est faiblement critique,

2/: $h(x_{0})>\frac{n-2}{4(n-1)}S_{\mathbf{g}}(x_{0})-\frac{(n-2)(n-4)}{%
8(n-1)}\frac{\bigtriangleup _{\mathbf{g}}f(x_{0})\,}{f(x_{0})\,}$ o\`{u}

a/: $h(x_{0})>0$

b/: \{$x_{0}\}=\{x\,/\,f(x)=\stackunder{M}{Sup}f\}$ et $f(x_{0})=1$, $0\leq
f\leq 1$, $Supp\,f=B(x_{0},r)$

c/: $\nabla ^{2}f(x_{0})<0$ .

On consid\`{e}re alors le chemin 
\[
t\rightarrow f_{t}=(1-t).1+t.f. 
\]
Notons que pour tout $t$: $\bigtriangleup _{\mathbf{g}}f_{t}=t\bigtriangleup
_{\mathbf{g}}f$ et $f_{t}(x_{0})=1=\stackunder{M}{Sup}\,f_{t}$. Posons 
\[
\lambda _{t}=Inf\,J_{h,f_{t},\mathbf{g}}. 
\]
Alors

\[
\lambda _{0}<K(n,2){{}^{-2}}(\stackunder{M}{Sup}f_{0})^{-\frac{n-2}{n}} 
\]
car ($h,1,\mathbf{g})$ est sous-critique et

\[
\lambda _{1}=K(n,2){{}^{-2}}(\stackunder{M}{Sup}f_{1})^{-\frac{n-2}{n}} 
\]
car ($h,f,\mathbf{g})$ est faiblement critique. Remarquons que $\stackunder{M%
}{Sup}\,f_{t}\,$est toujours \'{e}gal \`{a} 1.

Soit 
\[
t_{0}=Sup\{t\,/\,\lambda _{t}<K(n,2){{}^{-2}}(\stackunder{M}{Sup}f_{t})^{-%
\frac{n-2}{n}}\} 
\]
Alors $0<t_{0}\leq 1$ et 
\[
\lambda _{t_{0}}=K(n,2){{}^{-2}}(\stackunder{M}{Sup}f_{t_{0}})^{-\frac{n-2}{n%
}} 
\]
Avant de pouvoir appliquer la m\'{e}thode du chapitre 3, il faut prendre une
pr\'{e}caution: $f_{t_{0}}$ \'{e}tant faiblement critique, on sait seulement
que, au point de maximum $x_{0}$ on a 
\[
h(x_{0})\geq \frac{n-2}{4(n-1)}S_{\mathbf{g}}(x_{0})-\frac{(n-2)(n-4)}{8(n-1)%
}\frac{\bigtriangleup _{\mathbf{g}}f_{t_{0}}(x_{0})\,}{f_{t_{0}}(x_{0})\,} 
\]
car $\frac{\bigtriangleup _{\mathbf{g}}f_{t_{0}}(x_{0})\,}{f_{t_{0}}(x_{0})\,%
}=t_{0}\frac{\bigtriangleup _{\mathbf{g}}f(x_{0})\,}{f(x_{0})\,}$ avec $%
t_{0}\leq 1$, or nous avons besoin d'une in\'{e}galit\'{e} stricte.

On consid\`{e}re alors la suite, baptis\'{e}e ($f_{n})$, que l'on peut
construire \`{a} partir de la premi\`{e}re \'{e}tape: on peut trouver une
suite de fonctions ($f_{n}$), toutes telles que ($h,f_{n},\mathbf{g}$) soit
faiblement critique, avec 
\[
f_{n}(x_{0})=1=Supf_{n}\,\,\,et\,\,\,\bigtriangleup _{\mathbf{g}%
}f_{n}(x_{0})\rightarrow +\infty . 
\]
Pour chaque $f_{n}$, on note $t_{n}$ le ''$t_{0}$'' construit ci-dessus.
Donc pour tout $n$ :\thinspace 
\[
h\text{ est faiblement critique pour }(1-t_{n}).1+t_{n}.f_{n}\text{ et }%
\mathbf{g}. 
\]
Supposons que liminf $t_{n}=0$, ou, quitte \`{a} extraire, que $%
t_{n}\rightarrow 0.$ Alors, 
\[
(1-t_{n}).1+t_{n}.f_{n}\rightarrow 1 
\]
uniform\'{e}ment sur $M$ car $0\leq f_{n}\leq 1$. Or ($h,1,\mathbf{g})\,$est
sous-critique, donc il existe $u\in H_{1}^{2}(M)$ tel que 
\[
\frac{\int \left| \nabla u\right| ^{2}+\int hu^{2}}{(\int u^{2^{*}})^{\frac{2%
}{2^{*}}}}<K(n,2){{}^{-2}\,.} 
\]
Mais alors 
\[
\frac{\int \left| \nabla u\right| ^{2}+\int hu^{2}}{(\int
((1-t_{n}).1+t_{n}.f_{n})u^{2^{*}})^{\frac{2}{2^{*}}}}\rightarrow \frac{\int
\left| \nabla u\right| ^{2}+\int hu^{2}}{(\int u^{2^{*}})^{\frac{2}{2^{*}}}}%
<K(n,2){{}^{-2}} 
\]
tandis que 
\[
K(n,2){{}^{-2}=}K(n,2){{}^{-2}}(\stackunder{M}{Sup}%
((1-t_{n}).1+t_{n}.f_{n}))^{-\frac{n-2}{n}} 
\]
ce qui contredit le fait que $(h,(1-t_{n}).1+t_{n}.f_{n},\mathbf{g})$ soit
faiblement critique.

Donc, quitte \`{a} extraire $t_{n}\rightarrow t_{1}>0$

Comme $\bigtriangleup _{\mathbf{g}}f_{n}(x_{0})\rightarrow +\infty $, on
peut trouver $n$ assez grand tel que 
\[
\frac{(n-2)(n-4)}{8(n-1)}t_{n}\frac{\bigtriangleup _{\mathbf{g}%
}f_{n}(x_{0})\,}{f_{n}(x_{0})\,}>\frac{n-2}{4(n-1)}S_{g}(x_{0})-h(x_{0})\,. 
\]
Si on baptise $f$ cette derni\`{e}re fonction $f_{n}$ et $t_{0}$ ce
\thinspace $t_{n}$, on a alors un chemin 
\[
t\rightarrow f_{t}=(1-t).1+t.f 
\]
tel que :

a/: $\forall t<t_{0}$ : $(h,f_{t},\mathbf{g})$ est sous-critique,

b/: ($h,f_{t_{0}},\mathbf{g})$ est faiblement critique avec de plus :

b1/: \{$x_{0}\}=\{x\,/\,f_{t}(x)=\stackunder{M}{Sup}f_{t}\}$ et $%
f_{t}(x_{0})=1$ pour tout $t$

b2/: $h(x_{0})>\frac{n-2}{4(n-1)}S_{\mathbf{g}}(x_{0})-\frac{(n-2)(n-4)}{%
8(n-1)}\frac{\bigtriangleup _{\mathbf{g}}f_{t_{0}}(x_{0})\,}{%
f_{t_{0}}(x_{0})\,}$

b3/: $\nabla ^{2}f_{t_{0}}(x_{0})<0$

Pour tout $t<t_{0}$ il existe une solution minimisante $u_{t}$ de
l'\'{e}quation 
\[
\bigtriangleup _{\mathbf{g}}u_{t}+h.u_{t}=\lambda _{t}.f_{t}.u_{t}^{\frac{n+2%
}{n-2}} 
\]
avec $\int f_{t}u_{t}^{2^{*}}=1$. La suite ($u_{t}$) est born\'{e}e dans $%
H_{1}^{2}$ donc 
\[
u_{t}\stackrel{H_{1}^{2}}{\stackunder{t\rightarrow t_{0}}{\rightharpoondown }%
}u 
\]
et on se retrouve dans le sch\'{e}ma classique:

- soit $u>0$ et alors $u$ est solution de $\bigtriangleup _{\mathbf{g}%
}u+h.u=\lambda _{t_{0}}f_{t_{0}}.u^{\frac{n+2}{n-2}}$ et est une solution
minimisante, donc ($h,f_{t_{0}},\mathbf{g})$ est critique.

- soit $u\equiv 0$ $\,$et comme d'habitude la suite ($u_{t}$) se concentre.
Dans ce cas, l'\'{e}tude du ph\'{e}nom\`{e}ne de concentration est
n\'{e}anmoins plus simple qu'\`{a} la premi\`{e}re \'{e}tape car le chemin $%
f_{t}$ est constitu\'{e} de fonctions convergeant uniform\'{e}ment vers $f$
quand $t\rightarrow t_{0}$ avec $Supp\,f_{t}=B(x_{0},r)$. On peut trouver $%
\delta <r $ tel que $f>0$ sur $B(x_{0},\delta )$. Alors il existe $c>0$ tel
que pour tout $t$ on ait: 
\[
0<c\leq f_{t}\leq 1\text{ sur }B(x_{0},\delta ), 
\]
De plus, les $f_{t}$ atteignent toutes leur maximum en $x_{0}$, ce maximum
valant toujours 1. On peut alors reprendre toutes les \'{e}tapes du chapitre
3 (ou de la premi\`{e}re \'{e}tape), le fait de consid\'{e}rer une suite $%
f_{t}$ au second membre n'introduisant cette fois-ci aucun changement. On
aboutit \`{a} 
\[
h(x_{0})\leq \frac{n-2}{4(n-1)}S_{g}(x_{0})-\frac{(n-2)(n-4)}{8(n-1)}\frac{%
\bigtriangleup _{\mathbf{g}}f_{t_{0}}(x_{0})\,}{f_{t_{0}}(x_{0})\,} 
\]
d'o\`{u} une contradiction.

Par cons\'{e}quent ($h,f_{t_{0}},\mathbf{g})$ est critique et a une solution
minimisante.

Cette d\'{e}monstration montre en fait le r\'{e}sultat suivant, plus fort
mais moins ``parlant'':

\textbf{Th\'{e}or\`{e}me 4':}

\textit{Si }$h$\textit{\ est faiblement critique pour une fonction }$f$%
\textit{\ et une m\'{e}trique }$\mathbf{g},$ \textit{ces donn\'{e}es\
v\'{e}rifiant:}

\textit{1/: }$h(x)>\frac{n-2}{4(n-1)}S_{\mathbf{g}}(x)-\frac{(n-2)(n-4)}{%
8(n-1)}\frac{\bigtriangleup _{\mathbf{g}}f(x)\,}{f(x)\,}$\textit{\ au points
de maximum de }$f$

\textit{2/:}$\nabla ^{2}f(x)<0$\textit{\ au points de maximum de }$f$

\textit{3/: il existe une suite }$f_{t}\stackrel{C^{2}}{\stackunder{%
t\rightarrow t_{0}}{\rightarrow }}f$\textit{\ avec }$\stackunder{M}{Sup}%
f_{t}=\stackunder{M}{Sup}f$\textit{\ telle que (}$h,f_{t},\mathbf{g})$%
\textit{\ soit sous-critique pour }$t<t_{0}$

\textit{alors (}$h,f,\mathbf{g})$\textit{\ est critique et a des solutions
minimisantes.}

Comme nous l'avons dit au chapitre 1, cette d\'{e}monstration nous a
sugg\'{e}r\'{e} une autre d\'{e}finition possible pour les fonctions
critiques:
\\

\textbf{D\'{e}finition 4:}

\textit{Soient donn\'{e}es la vari\'{e}t\'{e} (}$M,\mathbf{g}$\textit{), }$%
\dim M\geqslant 3$\textit{, et une fonction }$h$\textit{\ telle que
l'op\'{e}rateur }$\triangle _{\mathbf{g}}+h$\textit{\ soit coercif. On
consid\`{e}re une fonction }$f\in C^{\infty }(M)$, telle que $Supf>0$.%
\textit{\ On dira que }$f$ \textit{est critique pour }$h$\textit{\ (et }$%
\mathbf{g}$) \textit{si:}

\textit{a/: }$\lambda _{h,f,\mathbf{g}}=\frac{1}{K(n,2){{}^{2}}(\stackunder{M%
}{Sup}f)^{\frac{n-2}{n}}}$

\textit{b/: pour toute fonction }$f^{\prime }$ \textit{telle que } $%
Supf=Supf^{\prime }$ \textit{et }$f^{\prime }\gneqq f$, \textit{\ }$\lambda
_{h,f^{\prime },\mathbf{g}}<\frac{1}{K(n,2){{}^{2}}(\stackunder{M}{Sup}%
f^{\prime })^{\frac{n-2}{n}}}\,.$

\textit{Remarque: si }$Supf=Supf^{\prime }$ \textit{et }$f^{\prime
}\lvertneqq f$, \textit{\ }$\lambda _{h,f^{\prime },\mathbf{g}}=\frac{1}{%
K(n,2){{}^{2}}(\stackunder{M}{Sup}f^{\prime })^{\frac{n-2}{n}}}$ \textit{%
puisque }$J_{h,f^{\prime },\mathbf{g}}(w)\geqslant J_{h,f,\mathbf{g}}(w)$%
\textit{\ pour toute fonction }$w$.
\\

Rappelons qu'il faut bien comparer dans cette d\'{e}finition des fonctions
de m\^{e}me Sup; ce sont en fait les classes $[f]$ et $[f^{\prime }]$ qui
importent, et on devrait en fait dire que c'est $[f]$ qui est critique pour $%
h$, (voir le chapitre 1). Il est alors naturel de se demander si les deux
d\'{e}finitions sont \'{e}quivalentes, autrement dit:

\begin{center}
$f$ \textit{est-elle critique pour }$h$\textit{\ si et seulement si }$h$ 
\textit{est critique pour }$f$ ?
\end{center}

Il est \`{a} noter que dans les deux cas on a toujours en tout point $P$\
o\`{u} $f$\ est maximum sur $M$\ :\textit{\ }$\frac{4(n-1)}{n-2}%
h(P)\geqslant S_{\mathbf{g}}(P)-\frac{n-4}{2}\frac{\bigtriangleup _{\mathbf{g%
}}f(P)}{f(P)}$\textit{\ }
\\

Cette question semble difficile. Nous obtenons l'\'{e}quivalence moyennant
quelques hypoth\`{e}ses suppl\'{e}mentaires:
\\
\\

\textbf{Th\'{e}or\`{e}me 5:}

\textit{Soient donn\'{e}es la vari\'{e}t\'{e} (}$M,\mathbf{g}$\textit{), }$%
\dim M\geqslant 5$\textit{, et une fonction }$h$\textit{\ telle que
l'op\'{e}rateur }$\triangle _{\mathbf{g}}+h$\textit{\ soit coercif. On
consid\`{e}re une fonction }$f\in C^{\infty }(M)$, telle que $Supf>0$\textit{%
\ et v\'{e}rifiant (}\textbf{H}$_{f}$).\textit{\ On suppose de plus qu'en
tout point }$P $\textit{\ o\`{u} }$f$\textit{\ est maximum sur }$M$\textit{\
:} 
\[
\mathit{\ }\frac{4(n-1)}{n-2}h(P)>S_{\mathbf{g}}(P)-\frac{n-4}{2}\frac{%
\bigtriangleup _{\mathbf{g}}f(P)}{f(P)}\mathit{\ .\ } 
\]
\textit{Alors }$f$ \textit{est critique pour }$h$\textit{\ si et seulement
si }$h$ \textit{est critique pour }$f$ $.$
\\

La d\'{e}monstration de ce th\'{e}or\`{e}me passe par la remarque suivante:
on sait que si $h$ est faiblement critique pour $f$ et $\mathbf{g}$ et que $%
\bigtriangleup _{\mathbf{g}}u+h.u=f.u^{\frac{n+2}{n-2}}$ a une solution
minimisante, alors $h$ est critique pour $f$ et $\mathbf{g}$. De m\^{e}me,
si $f$ est faiblement critique pour $h$ (au sens o\`{u} $\lambda _{h,f,%
\mathbf{g}}=K(n,2){{}^{-2}}(\stackunder{M}{Sup}f)^{-\frac{n-2}{n}}$ ) et que 
$\bigtriangleup _{\mathbf{g}}u+h.u=f.u^{\frac{n+2}{n-2}}$ a une solution
minimisante $u>0$, alors $f$ est critique pour $h$. En effet, si $f^{\prime
} $ est une fonction telle que $Supf=Supf^{\prime }$ et\textit{\ }$f^{\prime
}\gneqq f$, on a 
\[
\int f^{\prime }u^{2^{*}}>\int fu^{2^{*}} 
\]
car $u>0$. Donc 
\[
J_{h,f^{\prime },\mathbf{g}}(u)<J_{h,f,\mathbf{g}}(u)=K(n,2){{}^{-2}}(%
\stackunder{M}{Sup}f)^{-\frac{n-2}{n}}=K(n,2){{}^{-2}}(\stackunder{M}{Sup}%
f^{\prime })^{-\frac{n-2}{n}}\,. 
\]
A partir du travail que nous avons fait au chapitre 3 et dans ce chapitre,
la d\'{e}monstration devient alors rapide:

-Si $h$ est critique pour $f,$ on applique le th\'{e}or\`{e}me du chapitre
3: $\bigtriangleup _{\mathbf{g}}u+h.u=f.u^{\frac{n+2}{n-2}}$ a une solution
minimisante, et donc $f$ est critique pour $h$.

-Si $f$ est critique pour $h$, ces deux fonctions (et la m\'{e}trique)
v\'{e}rifiant les hypoth\`{e}ses du th\'{e}or\`{e}me: on a $\lambda _{h,f,%
\mathbf{g}}=K(n,2){{}^{-2}}(\stackunder{M}{Sup}f)^{-\frac{n-2}{n}},$ donc $h$
est faiblement critique pour $f$. On consid\`{e}re alors pour $t\stackrel{<}{%
\rightarrow }1$ la suite 
\[
t\rightarrow f_{t}=(1-t).Supf+t.f\,. 
\]
On a pour tout $t:\,$ $Supf_{t}=Supf$ et si $t<1$ alors $f_{t}\gneqq f$. 
Donc puisque $f$ est critique pour $h$, par d\'{e}finition:
\[
\lambda _{h,f_{t},\mathbf{g}}<K(n,2){{}^{-2}}(\stackunder{M}{Sup}f_{t})^{-%
\frac{n-2}{n}}\,. 
\]
On applique alors le th\'{e}or\`{e}me 4' donné ci-dessus pour obtenir que $h$ est 
critique pour $f$ avec des solutions minimisantes.

\section{M\'{e}thode alternative pour conclure la premi\`{e}re \'{e}tape}

Nous proposons tr\`{e}s rapidement et sch\'{e}matiquement une m\'{e}thode
reprenant la refactorisation du Hessien et ne n\'{e}cessitant plus que la
suite ($f_{t})$ soit d\'{e}croissante, et mettant en \'{e}vidence quelques
''estim\'{e}es '' pouvant \^{e}tre utiles dans un autre cadre.

1/: Comme 
\[
Supp\widetilde{f}_{t}=B(\frac{x_{0}(t)}{\mu _{t}},\frac{t}{\mu _{t}}) 
\]
et comme $\widetilde{f}_{t}\stackrel{p.p}{\rightarrow }1$,
n\'{e}cessairement 
\[
\frac{\mu _{t}}{t}\rightarrow 0\,. 
\]

2/: Supposons qu'il existe $\varepsilon _{0}>0$ tel que 
\[
\frac{\left| x_{0}(t)\right| }{t}\geq \varepsilon _{0}>0\,. 
\]
Alors $d_{\mathbf{g}}(x_{t},x_{0})\geq c_{0}t$ avec $c_{0}>0$, et par
construction des fonctions $f_{t}$: 
\[
\forall t:\,f_{t}(x_{t})\leq 1-\varepsilon _{1} 
\]
avec $\varepsilon _{1}>0$. Donc \`{a} extraction pr\`{e}s, 
\[
\widetilde{f}_{t}(0)=f_{t}(x_{t})\rightarrow 1-\varepsilon _{2} 
\]
avec $\varepsilon _{2}>0$. Or 
\[
\triangle _{\widetilde{\mathbf{g}}_{t}}\widetilde{u}_{t}+\mu _{t}^{2}.%
\widetilde{h}_{t}.\widetilde{u}_{t}=\lambda _{t}\widetilde{f}_{t}.\widetilde{%
u}_{t}^{\frac{n+2}{n-2}}\,. 
\]
En écrivant cette \'{e}quation en 0 et en faisant tendre $t$ vers 0, on
obtient 
\[
\triangle _{\mathbf{e}}\widetilde{u}(0)=\frac{1-\varepsilon _{2}}{K^{2}}%
\widetilde{u}(0)^{\frac{n+2}{n-2}} 
\]
avec 
\[
\widetilde{u}(x)=(1+\frac{K^{-2}}{n(n-2)}\left| x\right| ^{2})^{-\frac{n-2}{2%
}} 
\]
d'o\`{u} 
\[
K^{-2}=K^{-2}(1-\varepsilon _{2}) 
\]
et donc une contradiction. Par cons\'{e}quent 
\[
\frac{\left| x_{0}(t)\right| }{t}\rightarrow 0\,. 
\]

3/: On rappelle que gr\^{a}ce aux estim\'{e}es fortes 
\[
\frac{\int_{B(0,\delta )}\left| x\right| ^{p}(\eta \overline{u}%
_{t})^{2^{*}}dv_{\mathbf{g}_{t}}}{\int_{B(0,\delta )}\overline{u}_{t}^{2}dv_{%
\mathbf{g}_{t}}}\sim c\mu _{t}^{p-2} 
\]
et 
\[
\int_{B(0,\delta )}\overline{u}_{t}^{2}dv_{\mathbf{g}_{t}}\sim c\mu
_{t}^{2}\,. 
\]

4/: On reprend alors les calculs du chapitre 3 en ce qui concerne
l'expression $A_{t}^{1}$: 
\[
A_{t}^{1}=(\int_{B(0,\delta )}\overline{f}_{t}(\eta \overline{u}%
_{t})^{2^{*}}dx)^{\frac{n-2}{n}}\,-(Supf_{t}.\int_{B(0,\delta )}(\eta 
\overline{u}_{t})^{2^{*}}dx)^{\frac{n-2}{n}} 
\]
et 
\[
\overline{f}_{t}(x)\leq f(x_{0})+\frac{1}{2}\partial _{kl}\overline{f}%
_{t}(x_{0}(t)).(x^{k}-x_{0}^{k}(t))(x^{l}-x_{0}^{l}(t))+\frac{c}{t^{3}}%
\left| x-x_{0}(t)\right| ^{3}\,\,. 
\]
On a alors comme dans ce chapitre 
\[
A_{t}^{1}\leq c\frac{n-2}{n}\{F_{t}\}+C\{F_{t}\}^{2} 
\]
o\`{u} 
\[
\{F_{t}\}=\frac{1}{2}\partial _{kl}\overline{f}_{t}(x_{0}(t))\int_{B(0,%
\delta )}(x^{k}-x_{0}^{k}(t))(x^{l}-x_{0}^{l}(t))(\eta \overline{u}%
_{t})^{2^{*}}dx+\frac{c}{t^{3}}\int_{B(0,\delta )}\left| x-x_{0}(t)\right|
^{3}(\eta \overline{u}_{t})^{2^{*}}dx\,. 
\]
En reprenant la refactorisation du Hessien mais sans absorber les termes
d'ordre 3, on obtient 
\begin{eqnarray*}
\{F_{t}\} &=&\frac{1}{2}\partial _{kl}\overline{f}_{t}(x_{0}(t))\left\{
\int_{B(0,\delta )}x^{k}x^{l}(\eta \overline{u}%
_{t})^{2^{*}}dx+(x_{0}^{k}(t)z_{t}-\frac{\varepsilon _{t}^{k}}{z_{t}}%
)(x_{0}^{l}(t)z_{t}-\frac{\varepsilon _{t}^{l}}{z_{t}})-\frac{\varepsilon
_{t}^{k}\varepsilon _{t}^{l}}{z_{t}^{2}}\right\} \\
&&+\frac{c}{t^{3}}\int_{B(0,\delta )}\left| x-x_{0}(t)\right| ^{3}(\eta 
\overline{u}_{t})^{2^{*}}dx
\end{eqnarray*}
o\`{u} 
\[
\varepsilon _{t}^{k}=\int_{B(0,\delta )}x^{k}(\eta \overline{u}%
_{t})^{2^{*}}dx 
\]
\[
z_{t}=(\int_{B(0,\delta )}(\eta \overline{u}_{t})^{2^{*}}dx)^{\frac{1}{2}} 
\]
Par ailleurs 
\[
\frac{c}{t^{3}}\int_{B(0,\delta )}\left| x-x_{0}(t)\right| ^{3}(\eta 
\overline{u}_{t})^{2^{*}}dx\leq \frac{c^{\prime }}{t^{3}}\int_{B(0,\delta
)}\{\left| x\right| ^{3}+\left| x_{0}(t)\right| \left| x\right| ^{2}+\left|
x_{0}(t)\right| ^{2}\left| x\right| +\left| x_{0}(t)\right| ^{3}\}(\eta 
\overline{u}_{t})^{2^{*}}dx 
\]
On utilise alors 
\[
\frac{1}{2}\partial _{kl}\overline{f}_{t}(x_{0}(t))w^{k}w^{l}\leq -\frac{c}{%
t^{2}}\left| w\right| ^{2} 
\]
et 
\[
\frac{\left| \varepsilon _{t}^{k}\varepsilon _{t}^{l}\right| }{%
z_{t}^{2}\int_{B(0,\delta )}\overline{u}_{t}^{2}dv_{\mathbf{g}_{t}}}%
\rightarrow 0 
\]
pour obtenir que 
\begin{eqnarray*}
\{F_{t}\} &\leq &-\frac{c\mu _{t}^{2}}{t^{2}}-\frac{c\left| x_{0}(t)\right|
^{2}}{t^{2}}+o(\frac{\mu _{t}^{2}}{t^{2}}) \\
&&+c\{\frac{\mu _{t}^{3}}{t^{3}}+\frac{\mu _{t}^{2}\left| x_{0}(t)\right| }{%
t^{3}}+\frac{\mu _{t}\left| x_{0}(t)\right| ^{2}}{t^{3}}+\frac{\left|
x_{0}(t)\right| ^{3}}{t^{3}}\}\,.
\end{eqnarray*}
Ainsi 
\[
\begin{array}{ll}
\frac{\{F_{t}\}}{\int_{B(0,\delta )}\overline{u}_{t}^{2}dv_{\mathbf{g}_{t}}}%
\leq & -\frac{c_{1}}{t^{2}}-\frac{c_{2}\left| x_{0}(t)\right| ^{2}}{\mu
_{t}^{2}t^{2}}+o(\frac{1}{t^{2}}) \\ 
& +c_{3}\{\frac{\mu _{t}}{t^{3}}+\frac{\left| x_{0}(t)\right| }{t^{3}}+\frac{%
\left| x_{0}(t)\right| ^{2}}{\mu _{t}t^{3}}+\frac{\left| x_{0}(t)\right| ^{3}%
}{\mu _{t}^{2}t^{3}}\}
\end{array}
\]
et 
\[
\frac{\{F_{t}\}^{2}}{\int_{B(0,\delta )}\overline{u}_{t}^{2}dv_{\mathbf{g}%
_{t}}}\leq c_{4}\{\frac{\mu _{t}^{2}}{t^{4}}+\frac{\left| x_{0}(t)\right|
^{4}}{\mu _{t}^{2}t^{4}}+o(\frac{\mu _{t}^{2}}{t^{4}})+\frac{\mu _{t}^{4}}{%
t^{6}}+\frac{\mu _{t}^{2}\left| x_{0}(t)\right| ^{2}}{t^{6}}+\frac{\left|
x_{0}(t)\right| ^{4}}{t^{6}}+\frac{\left| x_{0}(t)\right| ^{6}}{\mu
_{t}^{2}t^{6}}\}\,. 
\]
On constate alors en utilisant 
\[
\frac{\left| x_{0}(t)\right| }{t}\rightarrow 0 
\]
et 
\[
\frac{\mu _{t}}{t}\rightarrow 0 
\]
que tout les termes strictement positifs sont n\'{e}gligeables devant $\frac{%
1}{t^{2}}$ ou $\frac{\left| x_{0}(t)\right| ^{2}}{\mu _{t}^{2}t^{2}}$ et par
cons\'{e}quent 
\[
\frac{A_{t}^{1}}{\int_{B(0,\delta )}\overline{u}_{t}^{2}dv_{\mathbf{g}_{t}}}%
\rightarrow -\infty 
\]
d'o\`{u} \`{a} partir de 
\[
h(x_{0})\leq \frac{n-2}{4(n-1)}S_{\mathbf{g}}(x_{0})+\stackunder{%
t\rightarrow 0}{\overline{\lim }}\frac{A_{t}^{1}}{\int_{B(0,\delta )}%
\overline{u}_{t}^{2}dv_{\mathbf{g}_{t}}} 
\]
une contradiction, et l'on récupère bien le théorème 4'.

\chapter{La Dimension 3}
\pagestyle{myheadings}\markboth{\textbf{La dimension 3.}}
{\textbf{La dimension 3.}}
Ce chapitre traite de la dimension 3. Il faut en effet remarquer que toute l'%
\'{e}tude des chapitres pr\'{e}c\'{e}dents portait sur des vari\'{e}t\'{e}s
de dimension $\geqslant 4$. De plus la dimension 4 elle-m\^{e}me pr\'{e}%
sente une particularit\'{e} puisque le terme 
\[
\frac{n-4}{2}\frac{\bigtriangleup _{\mathbf{g}}f(P)}{f(P)}\ 
\]
dispara\^{i}t. D'ailleurs, bien que les th\'{e}or\`{e}mes restent valables
pour $\dim M=4$, le r\'{e}sultat fondamental que nous obtenons sur les ph%
\'{e}nom\`{e}nes de concentration, pr\'{e}sent\'{e} au chapitre 3, n'est
valable que pour $\dim M\geqslant 5$. Le cas de la dimension 3 est lui
radicalement diff\'{e}rent. O. Druet $\left[ 10\right] $ a trait\'{e} ce cas
lorsque $f=cste=1$. L'introduction d'une fonction non constante n'apporte
pas ici de r\'{e}elles difficult\'{e}s, nous reprendrons donc rapidement la d%
\'{e}monstration d'Olivier Druet pour obtenir sa g\'{e}n\'{e}ralisation au
cas $f$ non constante. Cette dimension fait intervenir de fa\c{c}on
fondamentale la fonction de Green de l'op\'{e}rateur $\triangle _{\mathbf{g}%
}+h$. Si cet op\'{e}rateur est coercif, il existe une unique fonction de
classe $C^{2}$ 
\[
G_{h}:M\times M\backslash \{(x,x),x\in M\}\rightarrow \Bbb{R}
\]
sym\'{e}trique et strictement positive telle que, au sens des distributions,
on ait $\forall x\in M:$%
\[
\bigtriangleup _{\mathbf{g},y}G_{h}(x,y)+h(y)G_{h}(x,y)=\delta _{x}\,.
\]
En dimension 3, pour un point $x\in M$, et pour $y$ proche de $x$, $G_{h}$
peut se mettre sous la forme: 
\[
G_{h}(x,y)=\frac{1}{\omega _{2}d_{\mathbf{g}}(x,y)}+M_{h}(x)+o(1)
\]
o\`{u} $o(1)$ est \`{a} prendre pour $y\rightarrow x$. On appelle $M_{h}(x)$
la masse de la fonction de Green au point $x$.

On obtient une g\'{e}n\'{e}ralisation rapide des r\'{e}sultats d'Olivier
Druet traitant le cas $f=cste$:

\textbf{Th\'{e}or\`{e}me 6:}

\textit{Soient (}$M,\mathbf{g}$\textit{) une vari\'{e}t\'{e} compacte de
dimension 3 et une fonction }$f\in C^{\infty }(M)$ \textit{telle que} $%
Supf>0 $ (\textit{l'hypoth\`{e}se (}\textbf{H}$_{f}$\textit{) n'est pas
n\'{e}cessaire). Alors pour toute fonction }$h$ \textit{faiblement critique
pour }$f$\textit{\ et }$\mathbf{g}$\textit{, et pour tout }$x\in M\,$\textit{%
o\`{u} }$f$\textit{\ est maximum sur} $M$,$\,$\textit{on a }$%
M_{h}(x)\leqslant 0$\textit{.}

La condition\textit{\ } 
\[
M_{h}(x)\leqslant 0 
\]
$\,$apparait comme l'analogue de la condition 
\[
\frac{4(n-1)}{n-2}h(P)\geqslant S_{\mathbf{g}}(P)-\frac{n-4}{2}\frac{%
\bigtriangleup _{\mathbf{g}}f(P)}{f(P)} 
\]
que l'on avait en dimension $\geqslant 4$. Dans le cas $f=cste$, cette
condition doit \^{e}tre valable sur tout $M$. La particularit\'{e} de la
dimension 3 est alors d'offrir des fonctions critiques de toutes les formes:

\textbf{Th\'{e}or\`{e}me 7:}

\textit{Soient (}$M,\mathbf{g}$\textit{) une vari\'{e}t\'{e} compacte de
dimension 3 et une fonction }$f\in C^{\infty }(M)$ \textit{telle que} $%
Supf>0 $ (\textit{l'hypoth\`{e}se (}\textbf{H}$_{f}$\textit{) n'est pas
n\'{e}cessaire). Pour toute fonction }$h\in C^{\infty }(M)$, \textit{posons }%
$B(h)=\inf \{B/$ $h+B\,\,est\,\,faiblement\,\,critique\,\,pour\,\,f\}$. 
\textit{Alors }$h+B(h)$\textit{\ est une fonction critique pour }$f$\textit{.%
}

Enfin, en ce qui concerne l'existence de fonctions extr\'{e}males, on a le
th\'{e}or\`{e}me suivant:

\textbf{Th\'{e}or\`{e}me 8:}

\textit{Soient (}$M,\mathbf{g}$\textit{) une vari\'{e}t\'{e} compacte de
dimension 3 et une fonction }$f\in C^{\infty }(M)$ \textit{telle que} $%
Supf>0 $ (\textit{l'hypoth\`{e}se (}\textbf{H}$_{f}$\textit{) n'est pas
n\'{e}cessaire). Soit }$h$ \textit{une fonction critique pour }$f$\textit{\
et }$\mathbf{g}$.\textit{\ Alors au moins l'une des deux conditions
suivantes est remplie: }

\textit{a/: Il existe }$x\in M\,$\textit{\thinspace o\`{u} }$f$\textit{\ est
maximum sur} $M$ \textit{tel que }$M_{h}(x)=0$\textit{.}

\textit{b/: }$(E_{h,f,\mathbf{g}})$\textit{\ a des solutions extr\'{e}males.}

Nous allons reprendre le sch\'{e}ma des d\'{e}monstrations de l'article d'O.
Druet, l'introduction d'une fonction $f$ non constante n'apportant pas de
changements notables; essentiellement, la seule diff\'{e}rence est qu'il
faut se placer en un point de maximum de $f$.

Notons qu'en dimension 3, l'exposant critique vaut $2^{*}=\frac{2n}{n-2}=6$.

Par ailleurs, en ce qui concerne la fonction de Green, le principe du
maximum montre que si l'on a deux fonctions $h,h^{\prime }$ telles que 
\[
h\lvertneqq h^{\prime }, 
\]
alors pour tout $x\in M$ : 
\[
M_{h}(x)>M_{h^{\prime }}(x). 
\]

\section{D\'{e}monstration du Th\'{e}or\`{e}me 6}

Elle est bas\'{e}e sur l'utilisation de fonctions tests particuli\`{e}res.
Soient donc $f\in C^{\infty }(M)$ telle que $Supf>0$ (l'hypoth\`{e}se (H$%
_{f} $) n'est pas n\'{e}cessaire) et $h$ une fonction faiblement critique
pour $f$\ et $\mathbf{g}$. Par d\'{e}finition on a $\forall u\in
H_{1}^{2}(M) $ :
\begin{equation}
\left( \int_{M}f\left| u\right| ^{\frac{2n}{n-2}}dv_{\mathbf{g}}\right) ^{%
\frac{n-2}{n}}\leq K(n,2){{}^{2}}(\stackunder{M}{Sup}f)^{\frac{n-2}{n}%
}(\int_{M}\left| \nabla u\right| ^{2}dv_{\mathbf{g}}+\int_{M}h.u{{}^{2}}dv_{%
\mathbf{g}})  
\end{equation}

Soit $x_{0}$ un point de maximum de $f$ sur$M$.

Il existe une fonction $\varphi \in C^{\infty }(M)$, $\varphi >0$ avec $%
\varphi (x_{0})=1$ et $\nabla \varphi (x_{0})=0$, telle que la m\'{e}trique $%
\mathbf{g}_{\varphi }=\varphi ^{-4}\mathbf{g}$ v\'{e}rifie 
\[
S_{\mathbf{g}_{\varphi }}(x_{0})=0
\]
et 
\begin{equation}
\det \mathbf{g}_{\varphi }=1+O(d_{\mathbf{g}_{\varphi }}(x_{0},x)^{5}) 
\end{equation}
en coordonn\'{e}es normales autour de $x_{0}$ (pour l'existence de cette
fonction $\varphi $ on pourra voir l'article de Lee et Parker $\left[
25\right] $) . Soit alors $h_{\varphi }\,$la fonction obtenue par la loi de
transformation des fonctions critiques dans un changement de m\'{e}trique
conforme: 
\[
h_{\varphi }=\varphi ^{-5}(\bigtriangleup _{\mathbf{g}}\varphi +h\varphi )\,.
\]
Alors $h_{\varphi }$ est une fonction faiblement critique pour $f$ et $%
\mathbf{g}_{\varphi }$. Par la loi de transformation du Laplacien conforme,
on v\'{e}rifie que 
\[
G_{h_{\varphi }}(x_{0},x)=\varphi (x_{0})\varphi (x)G_{h}(x_{0},x)
\]
o\`{u} $G_{h_{\varphi }}\,$et $G_{h}$ sont les fonctions de Green associ\'{e}%
es respectivement \`{a} $\bigtriangleup _{\mathbf{g}_{\varphi }}+h_{\varphi
}\,$et \`{a} $\bigtriangleup _{\mathbf{g}}+h$. Alors comme $\varphi (x_{0})=1
$ et $\nabla \varphi (x_{0})=0$, on a 
\[
M_{h}(x_{0})=M_{h_{\varphi }}(x_{0})\,.
\]
Gr\^{a}ce \`{a} ces relations, il suffit donc de prouver le th\'{e}or\`{e}me
6 pour $h_{\varphi }$ , $f$ et $\mathbf{g}_{\varphi }.\,$On peut donc
supposer sans perdre en g\'{e}n\'{e}ralit\'{e} que $\varphi \equiv 1$ et
supprimer l'indice $\varphi $ dans la suite.

Fixons maintenant une fonction cut-off $\eta \in C_{c}^{\infty
}(B(x_{0},2\delta ))$, $\eta \equiv 1\,$sur $B(x_{0},\delta )\,$avec $\delta
>0$ assez petit. On peut \'{e}crire la fonction de Green $G_{h}\,$sous la
forme: 
\[
\omega _{2}G_{h}(x_{0},x)=\frac{\eta (x)}{d_{\mathbf{g}}(x,y)}+\beta (x) 
\]
o\`{u} $\beta \in C_{loc}^{\infty }(M\backslash \{x_{0}\})$. Dans $%
M\backslash B(x_{0},\delta )$, $\beta $ v\'{e}rifie 
\begin{equation}
\bigtriangleup _{\mathbf{g}}\beta +h\beta =-\bigtriangleup _{\mathbf{g}}(%
\frac{\eta (x)}{d_{\mathbf{g}}(x_{0},x)})-h\frac{\eta (x)}{d_{\mathbf{g}%
}(x_{0},x)}\,. 
\end{equation}
Et dans $B(x_{0},\delta )$, $\beta $ v\'{e}rifie en coordonn\'{e}es normales 
\begin{equation}
\bigtriangleup _{\mathbf{g}}\beta +h\beta =-\frac{1}{2d_{\mathbf{g}%
}(x_{0},x)^{2}}\partial _{r}(\ln (\det \mathbf{g}))-\frac{h(x)}{d_{\mathbf{g}%
}(x_{0},x)}\,.  
\end{equation}
On pourra se reporter \`{a} l'appendice B pour quelques indications sur ces
propri\'{e}t\'{e}s de la fonction de Green. Par les th\'{e}ories elliptiques
standard, on sait que $\beta \in C^{0}(M)\cap H_{1}^{2}(M)$ et on a par
ailleurs $\beta (x_{0})=\omega _{2}M_{h}(x_{0})$. Le but est maintenant
d'introduire des fonctions tests d\'{e}riv\'{e}es de celles
pr\'{e}sent\'{e}es au premier chapitre dans la relation (7.1). Pour $%
\varepsilon >0$ on d\'{e}finit 
\[
v_{\varepsilon }(x)=(\varepsilon ^{2}+d_{\mathbf{g}}(x_{0},x)^{2})^{-\frac{1%
}{2}} 
\]
et 
\[
u_{\varepsilon }(x)=\eta (x)v_{\varepsilon }(x)+\beta (x)\,. 
\]
Comme ($h,f,\mathbf{g}$) est faiblement critique, nous avons pour tout $%
\varepsilon >0$ d'apr\`{e}s (7.1)
\begin{equation}
K^{-2}(\stackunder{M}{Sup}f)^{-\frac{1}{3}}\left( \int_{M}f\left|
u_{\varepsilon }\right| ^{6}dv_{\mathbf{g}}\right) ^{\frac{1}{3}}\leq
\int_{M}\left| \nabla u_{\varepsilon }\right| _{\mathbf{g}}^{2}dv_{\mathbf{g}%
}+\int_{M}h.u_{\varepsilon }^{2}{}dv_{\mathbf{g}}  
\end{equation}
o\`{u} l'on note $K^{-2}=K(3,2)^{-2}$. Nous allons calculer les
d\'{e}veloppements en $\varepsilon $ des deux membres de cette
in\'{e}galit\'{e}. On prend d\'{e}j\`{a} $\delta $ assez petit pour que $f>0$
sur $B(x_{0},\delta )$, on rappelle que $x_{0}$ est un point o\`{u} $f$ est
maximum sur $M$. Nous aurons besoin d'\'{e}valuer les int\'{e}grales de la
forme 
\[
\int_{B(x_{0},\delta )}r^{p}v_{\varepsilon }^{q}dv_{\mathbf{g}} 
\]
en fonction de $\varepsilon $, o\`{u} $r=d_{\mathbf{g}}(x_{0},x)$. Pour cela
on utilise le choix de m\'{e}trique que nous avons fait, \`{a} savoir la
relation (7.2), qui permet d'\'{e}crire 
\[
\int_{B(x_{0},\delta )}r^{p}v_{\varepsilon }^{q}dv_{\mathbf{g}%
}=\int_{B(0,\delta )}r^{p}v_{\varepsilon }^{q}dx+\int_{B(0,\delta
)}O(r^{p+5})v_{\varepsilon }^{q}dx 
\]
o\`{u} les int\'{e}grales du second membre sont \`{a} comprendre comme
\'{e}tant lues dans une carte munie de la m\'{e}trique euclidienne. A partir
de l\`{a}, on peut \'{e}crire 
\begin{eqnarray*}
\int_{B(0,\delta )}r^{p}v_{\varepsilon }^{q}dx & &=\int_{B(0,\delta )}r^{p}%
\frac{1}{(\varepsilon ^{2}+r^{2})^{\frac{q}{2}}}dx \\ 
&& =\omega _{2}\int_{0}^{\delta }r^{p}\frac{1}{(\varepsilon ^{2}+r^{2})^{%
\frac{q}{2}}}r^{2}dr \\ 
&& =\omega _{2}.\varepsilon ^{p-q+3}\int_{0}^{\delta /\varepsilon }\frac{%
s^{p+2}}{(1+s^{2})^{\frac{q}{2}}}ds\,.
\end{eqnarray*}
Alors, en \'{e}crivant $\int_{0}^{\delta /\varepsilon }=\int_{0}^{+\infty
}+\int_{+\infty }^{\delta /\varepsilon }$ on obtient:

Si $q-p>3$: 
\begin{equation}
\int_{B(0,\delta )}r^{p}v_{\varepsilon }^{q}dx=\omega
_{2}I_{p,q}.\varepsilon ^{p-q+3}+O(1) 
\end{equation}

Si $q-p\leq 3$: 
\begin{equation}
\int_{B(0,\delta )}r^{p}v_{\varepsilon }^{q}dx\sim \omega _{2}.\varepsilon
^{p-q+3}(c\varepsilon ^{q-p-3}+c^{\prime })  
\end{equation}
o\`{u} $I_{p,q}=\int_{0}^{+\infty }\frac{s^{p+2}}{(1+s^{2})^{\frac{q}{2}}}ds$%
.

Notons \`{a} propos de ces int\'{e}grales que 
\[
\omega _{2}\int_{0}^{+\infty }(1+s^{2})^{-3}s^{2}ds=\frac{\omega _{3}}{8} 
\]
et que 
\[
\omega _{2}\int_{0}^{+\infty }(1+s^{2})^{-2}s^{2}ds=\frac{\omega _{3}}{2} 
\]
o\`{u} l'on rappelle que $\omega _{n}$ est le volume de la sph\`{e}re $S^{n}$%
.

Le calcul du d\'{e}veloppement du membre de droite de (7.5). ne pr\'{e}sente
strictement aucune diff\'{e}rence avec l'article d'O. Druet puisque notre
fonction $f$ n'intervient pas. Le calcul donne: 
\begin{eqnarray}
\int_{M}\left| \nabla u_{\varepsilon }\right| ^{2}dv_{\mathbf{g}%
}+\int_{M}h.u_{\varepsilon }^{2}{}dv_{\mathbf{g}} 
&=&3\varepsilon ^{-1}\omega _{2}\int_{0}^{+\infty
}(1+s^{2})^{-3}s^{2}ds+\omega _{2}\beta (x_{0})  \\
&&-h(x_{0})\omega _{2}(\int_{0}^{+\infty }(1+s^{2})^{-1}ds)\varepsilon\nonumber\\
&&+2h(x_{0})\omega _{2}[\int_{0}^{+\infty }(1+s^{2})^{-\frac{1}{2}%
}(s+(1+s^{2})^{\frac{1}{2}})ds]\varepsilon +o(\varepsilon )\nonumber
\end{eqnarray}
Venons-en au calcul du membre de gauche, en indiquant les diff\'{e}rences
introduites dans les calculs par la pr\'{e}sence de $f$. 
\begin{eqnarray*}
\int_{M}fu_{\varepsilon }^{6}dv_{\mathbf{g}} &=&\int_{M}f(\eta
v_{\varepsilon }+\beta )^{6}dv_{\mathbf{g}} \\
&=&\int_{B(x_{0},\delta )}fv_{\varepsilon }^{6}dv_{\mathbf{g}%
}+6\int_{B(x_{0},\delta )}fv_{\varepsilon }^{5}\beta dv_{\mathbf{g}%
}+15\int_{B(x_{0},\delta )}fv_{\varepsilon }^{4}\beta ^{2}dv_{\mathbf{g}%
}+\int_{M}f\beta ^{6}dv_{\mathbf{g}}+o(\varepsilon ^{-1})
\end{eqnarray*}

- Comme $\beta \in C^{0}(M)$, $\int_{M}f\beta ^{6}=O(1)=o(\varepsilon ^{-1})$%
. Il en va de m\^{e}me pour toutes les int\'{e}grales born\'{e}es
ind\'{e}pendamment de $\varepsilon $ comme par exemple celles du type $%
\int_{B(x_{0},2\delta )\backslash B(x_{0},\delta )}f\eta ^{p}v_{\varepsilon
}^{q}\beta ^{r}$, d'o\`{u} le $o(\varepsilon ^{-1})$ \`{a} la fin du
d\'{e}veloppement ci-dessus.

- Ensuite, comme $\beta \in C^{0}(M)\cap H_{1}^{2}(M)$ , il existe un $%
0<\alpha <1$ tel que $\beta \in C^{0,\alpha }(M)$. Nous pourrons donc
\'{e}crire 
\[
\beta (x)=\beta (x_{0})+O(r^{\alpha }) 
\]
Nous utiliserons \'{e}galement deux d\'{e}veloppements de $f$ en $x_{0}$
o\`{u} $f$ est maximum 
\begin{eqnarray*}
f(x) &=&f(x_{0})+\frac{1}{2}\partial _{ij}f(x_{0})x^{i}x^{j}+O(r^{3}) \\
f(x) &=&f(x_{0})+O(r^{2})\,.
\end{eqnarray*}
En introduisant ces d\'{e}veloppements et en utilisant (7.2), (7.6) et
(7.7), on obtient

- premi\`{e}rement 
\[
15\int_{B(x_{0},\delta )}fv_{\varepsilon }^{4}\beta ^{2}dv_{\mathbf{g}%
}=15f(x_{0})\beta ^{2}(x_{0})\omega _{2}I_{0,4}\varepsilon
^{-1}+o(\varepsilon ^{-1}) 
\]

- deuxi\`{e}mement 
\begin{eqnarray*}
6\int_{B(x_{0},\delta )}fv_{\varepsilon }^{5}\beta dv_{\mathbf{g}} & 
=&6f(x_{0})\int_{B(x_{0},\delta )}v_{\varepsilon }^{5}\beta dv_{\mathbf{g}%
}+6\int_{B(x_{0},\delta )}O(r^{2})v_{\varepsilon }^{5}\beta dv_{\mathbf{g}}
\\ 
& =&6f(x_{0})\int_{B(x_{0},\delta )}v_{\varepsilon }^{5}\beta dv_{\mathbf{g}%
}+o(\varepsilon ^{-1})\,.
\end{eqnarray*}
En utilisant les \'{e}quations (7.3) et (7.4) v\'{e}rifi\'{e}es par $\beta $
et le fait que 
\[
\bigtriangleup _{e}v_{\varepsilon }=3\varepsilon ^{2}v_{\varepsilon }^{5}\, 
\]
les calculs de l'article d'O. Druet donnent alors 
\begin{eqnarray*}
6\int_{B(x_{0},\delta )}fv_{\varepsilon }^{5}\beta dv_{\mathbf{g}}=&&2\omega
_{2}\varepsilon ^{-2}f(x_{0})\beta (x_{0})\\
&&+2f(x_{0})h(x_{0})\omega
_{2}[\int_{0}^{+\infty }(1+s^{2})^{-\frac{1}{2}}(s+(1+s^{2})^{\frac{1}{2}%
})ds]\varepsilon ^{-1}+o(\varepsilon ^{-1})\,. 
\end{eqnarray*}
- Enfin 
\begin{eqnarray*}
\int_{B(x_{0},\delta )}fv_{\varepsilon }^{6}dv_{\mathbf{g}} && 
=f(x_{0})\int_{B(x_{0},\delta )}v_{\varepsilon }^{6}+\frac{1}{2}\partial
_{ij}f(x_{0})\int_{B(x_{0},\delta )}x^{i}x^{j}v_{\varepsilon
}^{6}+\int_{B(x_{0},\delta )}O(r^{3})v_{\varepsilon }^{6} \\ 
&& =f(x_{0})\omega _{2}I_{0,6}\varepsilon ^{-3}-\bigtriangleup _{\mathbf{g}%
}f(x_{0})\omega _{2}I_{2,6}\varepsilon ^{-1}+o(\varepsilon ^{-1})\,.
\end{eqnarray*}
Finalement 
\begin{eqnarray*}
\int_{M}fu_{\varepsilon }^{6}dv_{\mathbf{g}}&= & f(x_{0})\omega
_{2}I_{0,6}\varepsilon ^{-3}-\bigtriangleup _{\mathbf{g}}f(x_{0})\omega
_{2}I_{2,6}\varepsilon ^{-1}+15f(x_{0})\beta ^{2}(x_{0})\omega
_{2}I_{0,4}\varepsilon ^{-1} \\ 
&& +2\omega _{2}\varepsilon ^{-2}f(x_{0})\beta
(x_{0})+2f(x_{0})h(x_{0})\omega _{2}[\int_{0}^{+\infty }(1+s^{2})^{-\frac{1}{%
2}}(s+(1+s^{2})^{\frac{1}{2}})ds]\varepsilon ^{-1}\\
&&+o(\varepsilon ^{-1})
\end{eqnarray*}
D'o\`{u} 
\begin{eqnarray*}
\left( \int_{M}f\left| u_{\varepsilon }\right| ^{6}dv_{\mathbf{g}}\right) ^{%
\frac{1}{3}}= & f(x_{0})^{\frac{1}{3}}\frac{\omega _{3}^{\frac{1}{3}}}{2}%
[\varepsilon ^{-1}-\frac{\bigtriangleup _{\mathbf{g}}f(x_{0})}{f(x_{0})}%
\frac{8\omega _{2}}{3\omega _{3}}I_{2,6}\varepsilon +20\beta
^{2}(x_{0})\varepsilon +\frac{16\omega _{2}}{3\omega _{3}}\beta (x_{0}) \\ 
& +\frac{16}{3\omega _{3}}h(x_{0})\omega _{2}[\int_{0}^{+\infty }(1+s^{2})^{-%
\frac{1}{2}}(s+(1+s^{2})^{\frac{1}{2}})ds]\varepsilon +o(\varepsilon )]\,.
\end{eqnarray*}
En notant que 
\[
K^{-2}(\stackunder{M}{Sup}f)^{-\frac{1}{3}}=K^{-2}f(x_{0})^{-\frac{1}{3}}=%
\frac{3}{4}\omega _{3}^{\frac{2}{3}}f(x_{0})^{-\frac{1}{3}} 
\]
on obtient finalement: 
\begin{eqnarray}
K^{-2}(\stackunder{M}{Sup}f)^{-\frac{1}{3}}\left( \int_{M}f\left|
u_{\varepsilon }\right| ^{6}dv_{\mathbf{g}}\right) ^{\frac{1}{3}}= && \frac{3%
}{8}\omega _{3}\varepsilon ^{-1}+2\omega _{2}\beta (x_{0})+\frac{15}{2}%
\omega _{3}\beta ^{2}(x_{0})\varepsilon -\frac{\bigtriangleup _{\mathbf{g}%
}f(x_{0})}{f(x_{0})}I_{2,6}\varepsilon \nonumber\\ 
& &+2h(x_{0})\omega _{2}[\int_{0}^{+\infty }(1+s^{2})^{-\frac{1}{2}%
}(s+(1+s^{2})^{\frac{1}{2}})ds]\varepsilon\nonumber\\
&& +o(\varepsilon )\,.
\end{eqnarray}
Maintenant on \'{e}crit que ($h,f,\mathbf{g}$) est faiblement critique,
\`{a} savoir que 
\begin{equation}
K^{-2}(\stackunder{M}{Sup}f)^{-\frac{1}{3}}\left( \int_{M}f\left|
u_{\varepsilon }\right| ^{6}dv_{\mathbf{g}}\right) ^{\frac{1}{3}}\leq
\int_{M}\left| \nabla u_{\varepsilon }\right| ^{2}dv_{\mathbf{g}%
}+\int_{M}h.u_{\varepsilon }^{2}{}dv_{\mathbf{g}} 
\end{equation}
d'o\`{u} avec (7.8) et (7.9)
\[
\omega _{2}\beta (x_{0})+h(x_{0})\omega _{2}(\int_{0}^{+\infty
}(1+s^{2})^{-1}ds)\varepsilon +\frac{15}{2}\omega _{3}\beta
^{2}(x_{0})\varepsilon -\frac{\bigtriangleup _{\mathbf{g}}f(x_{0})}{f(x_{0})}%
I_{2,6}\varepsilon +o(\varepsilon )\leq 0\,. 
\]
Ceci \'{e}tant vrai pour tout $\varepsilon >0$, on a 
\[
\omega _{2}\beta (x_{0})\leq 0 
\]
et donc 
\[
M_{h}(x_{0})\leq 0\text{ si }x_{0}\text{ est un point de maximum de }f 
\]
Il est int\'{e}ressant (et un peu surprenant) de noter que $\frac{%
\bigtriangleup _{\mathbf{g}}f(x_{0})}{f(x_{0})}\,$appara\^{i}t en dimension
3 \`{a} un ordre tel qu'il n'intervient pas dans la condition $%
M_{h}(x_{0})\leq 0$.

\section{D\'{e}monstration des th\'{e}or\`{e}mes 7 et 8}

Soit $h\in C^{\infty }(M)$. D'apr\`{e}s les d\'{e}finitions, il existe $B>0$
tel que $h+B$ est faiblement critique; il suffit que $h+B\geq B_{0}K^{-2}$.
Si l'on pose 
\[
B(h)=\inf \{B/h+B\,\,est\,\,faiblement\,\,critique\,\,pour\,\,f\}, 
\]
alors, quitte \`{a} changer $h$ en $h+B(h)$, on peut supposer que $B(h)=0$.

Pour $t\geq 0$ posons 
\[
\lambda _{t}=\stackunder{u\in H_{1}^{2},\int u^{2^{*}}=1}{\inf }%
(\int_{M}\left| \nabla u\right| _{\mathbf{g}}^{2}dv_{\mathbf{g}%
}+\int_{M}(h-t)u^{2}dv_{\mathbf{g}})\,. 
\]
Alors par d\'{e}finition de $B(h)$, pour tout $t>0:$%
\[
\lambda _{t}<K^{-1}(\stackunder{M}{Sup}f)^{-\frac{1}{3}} 
\]
et pour $t=0$%
\[
\lambda _{1}=K^{-1}(\stackunder{M}{Sup}f)^{-\frac{1}{3}} 
\]
Il existe donc pour tout $t>0$ une fonction $u_{t}>0$ solution de 
\[
\bigtriangleup _{\mathbf{g}}u_{t}+(h-t)u_{t}=\lambda _{t}fu_{t}^{5} 
\]
avec 
\[
\int_{M}fu_{t}^{2^{*}}=\int_{M}fu_{t}^{6}=1 
\]
Clairement, ($u_{t}$) est born\'{e}e dans $H_{1}^{2}$ et donc quitte \`{a}
extraire, $u_{t}\rightharpoondown u_{0}$ faiblement dans $H_{1}^{2}$. Si $%
u_{0}\neq 0$, $u_{0}$ est minimisante pour ($h+B(h),f,\mathbf{g}$) et donc $%
h+B(h)$ est critique.

Supposons donc maintenant que $u_{0}\equiv 0$. La suite d\'{e}veloppe un
ph\'{e}nom\`{e}ne de concentration. Comme au chapitre 3, on sait qu'il
existe une suite de points ($x_{t})_{t>0}$ convergeant vers un point de
concentration $x_{0}\in M,$ qui est n\'{e}cessairement un point o\`{u} $f$
est maximum sur $M,$ tels que 
\[
u_{t}(x_{t})=Sup_{M}u_{t}:=\mu _{t}^{-\frac{1}{2}}\,. 
\]
De plus, nous avons les estim\'{e}es fortes 
\begin{equation}
\mu _{t}^{-\frac{1}{2}}d_{\mathbf{g}}(x,x_{t})u_{t}(x)\leq C  
\end{equation}
et par un r\'{e}sultat d'O. Druet et F. Robert que nous avons cit\'{e} \`{a}
la fin du chapitre 3 : 
\begin{equation}
\mu _{t}^{-\frac{1}{2}}u_{t}\rightarrow 2\frac{\omega _{2}}{\omega
_{3}f(x_{0})}G_{h}(x_{0},.)  
\end{equation}
dans $C_{loc}^{2}(M\backslash \{x_{0}\})$ quand $t\rightarrow 0$ et o\`{u}
l'on rappelle que $f(x_{0})=\stackunder{M}{Sup}f$. Enfin, \`{a} partir des
estim\'{e}es fortes, on peut obtenir une estim\'{e}e sur le gradient par les
m\'{e}thodes elliptiques standard; ce genre de m\'{e}thode sera
pr\'{e}sent\'{e}e \`{a} l'appendice B lors de la construction de la fonction
de Green.
\begin{equation}
\mu _{t}^{-\frac{1}{2}}d_{\mathbf{g}}(x,x_{t})^{2}\left| \nabla
u_{t}(x)\right| _{\mathbf{g}}\leq C 
\end{equation}
Maintenant, pour $\delta >0$ on se place dans des cartes exponentielles
autour des $x_{t}$ en posant 
\begin{eqnarray*}
\mathbf{g}_{t} &=&\exp _{x_{t}}^{*}\mathbf{g} \\
\overline{u}_{t}(x) &=&u_{t}(\exp _{x_{t}}(x))
\end{eqnarray*}
et 
\begin{eqnarray*}
\overline{f}_{t}(x) &=&f(\exp _{x_{t}}(x)) \\
\overline{h}_{t} &=&h(\exp _{x_{t}}(x))\,.
\end{eqnarray*}
A la diff\'{e}rence des dimensions $\geq 4$, il faut utiliser l'identit\'{e}
de Pohozaev plut\^{o}t que l'in\'{e}galit\'{e} de Sobolev pour obtenir le
r\'{e}sultat: 
\begin{eqnarray}
\int_{B(0,\delta )}(x^{k}\partial _{k}\overline{u}_{t}+\frac{1}{2}\overline{u%
}_{t})\bigtriangleup _{e}\overline{u}_{t}dx= && \delta \int_{\partial
B(0,\delta )}[\frac{1}{2}\left| \nabla \overline{u}_{t}\right|
_{e}^{2}-(\nabla \overline{u}_{t},\nu )_{e}^{2}]d\sigma _{e}\nonumber \\ 
&& -\frac{1}{2}\int_{\partial B(0,\delta )}\overline{u}_{t}(\nabla \overline{u%
}_{t},\nu )_{e}d\sigma _{e}
\end{eqnarray}
o\`{u} $\nu $ est le vecteur normal unitaire ext\'{e}rieur. Avec (7.13)
et (7.14), on obtient 
\begin{eqnarray*}
\frac{f(x_{0})^{2}}{4\omega _{3}^{-\frac{2}{3}}\omega _{3}^{2}}\stackunder{%
t\rightarrow 0}{\lim }\mu _{t}^{-1}\int_{B(0,\delta )}(x^{k}\partial _{k}%
\overline{u}_{t}+\frac{1}{2}\overline{u}_{t})\bigtriangleup _{e}\overline{u}%
_{t}dx &=&\delta \int_{\partial B(0,\delta )}[\frac{1}{2}\left| \nabla 
\widetilde{G}\right| _{e}^{2}-(\nabla \widetilde{G},\nu )_{e}^{2}]d\sigma
_{e} \\
&&-\frac{1}{2}\int_{\partial B(0,\delta )}\widetilde{G}(\nabla \widetilde{G}%
,\nu )_{e}d\sigma _{e}
\end{eqnarray*}
o\`{u} 
\[
\widetilde{G}(x)=G_{h}((x_{0},\exp _{x_{0}}(x))\,. 
\]
On veut calculer la limite du membre de gauche. On commence par
d\'{e}velopper le laplacien: 
\[
\bigtriangleup _{e}\overline{u}_{t}=\bigtriangleup _{\mathbf{g}_{t}}%
\overline{u}_{t}+(\mathbf{g}_{t}^{ij}-\delta ^{ij})\partial _{ij}\overline{u}%
_{t}-\mathbf{g}_{t}^{ij}\Gamma (\mathbf{g}_{t})_{ij}^{m}\partial _{m}%
\overline{u}_{t} 
\]
et on \'{e}crit 
\begin{eqnarray*}
\int_{B(0,\delta )}(x^{k}\partial _{k}\overline{u}_{t}+\frac{1}{2}\overline{u%
}_{t})\bigtriangleup _{e}\overline{u}_{t}dx=& & \int_{B(0,\delta
)}(x^{k}\partial _{k}\overline{u}_{t}+\frac{1}{2}\overline{u}%
_{t})\bigtriangleup _{\mathbf{g}_{t}}\overline{u}_{t}dx \\ 
& &+\int_{B(0,\delta )}(x^{k}\partial _{k}\overline{u}_{t}+\frac{1}{2}%
\overline{u}_{t})(\mathbf{g}_{t}^{ij}-\delta ^{ij})\partial _{ij}\overline{u}%
_{t}dx \\ 
&& -\int_{B(0,\delta )}(x^{k}\partial _{k}\overline{u}_{t}+\frac{1}{2}%
\overline{u}_{t})\mathbf{g}_{t}^{ij}\Gamma (\mathbf{g}_{t})_{ij}^{m}\partial
_{m}\overline{u}_{t}dx\,.
\end{eqnarray*}
Alors, avec l'\'{e}quation v\'{e}rifi\'{e}e par les $\overline{u}_{t}$ et
quelques int\'{e}grations par parties, cette expression devient 
\begin{eqnarray*}
\int_{B(0,\delta )}(x^{k}\partial _{k}\overline{u}_{t}+\frac{1}{2}\overline{u%
}_{t})\bigtriangleup _{e}\overline{u}_{t}dx&= & \int_{B(0,\delta )}((%
\overline{h}_{t}-t)+\frac{1}{2}x^{k}\partial _{k}\overline{h}_{t})\overline{u%
}_{t}^{2}dx-\frac{\lambda _{t}}{6}\int_{B(0,\delta )}\overline{u}%
_{t}^{6}x^{k}\partial _{k}\overline{f}_{t}dx \\ 
&& -\int_{B(0,\delta )}(x^{k}\partial _{k}\overline{u}_{t}+\frac{1}{2}%
\overline{u}_{t})\partial _{j}\mathbf{g}_{t}^{ij}\partial _{i}\overline{u}%
_{t}dx\\
&&+\frac{1}{2}\int_{B(0,\delta )}x^{k}\partial _{k}\mathbf{g}%
_{t}^{ij}\partial _{i}\overline{u}_{t}\partial _{j}\overline{u}_{t}dx \\ 
&& -\int_{B(0,\delta )}(x^{k}\partial _{k}\overline{u}_{t}+\frac{1}{2}%
\overline{u}_{t})\mathbf{g}_{t}^{ij}\Gamma (\mathbf{g}_{t})_{ij}^{m}\partial
_{m}\overline{u}_{t}dx \\ 
& &-\frac{\delta }{2}\int_{\partial B(0,\delta )}(\overline{h}_{t}-t)%
\overline{u}_{t}^{2}d\sigma _{e}+\frac{\lambda _{t}}{6}\delta \int_{\partial
B(0,\delta )}\overline{f}_{t}\overline{u}_{t}^{6}d\sigma _{e} \\
&&+\int_{\partial B(0,\delta )}(x^{k}\partial _{k}\overline{u}_{t}+\frac{1}{2}%
\overline{u}_{t})(\mathbf{g}_{t}^{ij}-\delta ^{ij})\nu _{j}\partial _{i}%
\overline{u}_{t}d\sigma _{e} \\ 
&&-\frac{\delta }{2}\int_{\partial B(0,\delta )}(\mathbf{g}_{t}^{ij}-\delta
^{ij})\partial _{i}\overline{u}_{t}\partial _{j}\overline{u}_{t}d\sigma
_{e}\,.
\end{eqnarray*}
La principale diff\'{e}rence avec le cas $f=cste$ trait\'{e} par O. Druet
est le terme 
\[
\frac{\lambda _{t}}{6}\int_{B(0,\delta )}\overline{u}_{t}^{6}x^{k}\partial
_{k}\overline{f}_{t}dx\,. 
\]
Maintenant, gr\^{a}ce aux estim\'{e}es fortes (7.11) \`{a} (7.13), on peut
calculer la limite que l'on cherchait, en notant que 
\begin{eqnarray*}
\left| \partial _{k}\mathbf{g}_{t}^{ij}\right| &\leq &C\left| x\right| \\
\text{et }\left| \mathbf{g}_{t}^{ij}\Gamma (\mathbf{g}_{t})_{ij}^{m}\right|
&\leq &C\left| x\right|
\end{eqnarray*}
ce qui permet d'appliquer le th\'{e}or\`{e}me de convergence domin\'{e}e de
Lebesgue, pour obtenir 
\begin{eqnarray*}
&&\int_{B(0,\delta )}(\widetilde{h}+\frac{1}{2}x^{k}\partial _{k}\widetilde{h%
})\widetilde{G}^{2}dx-\int_{B(0,\delta )}(x^{k}\partial _{k}\widetilde{G}+%
\frac{1}{2}\widetilde{G})\partial _{j}\mathbf{g}_{t}^{ij}\partial _{i}%
\widetilde{G}dx \\
&&+\frac{1}{2}\int_{B(0,\delta )}x^{k}\partial _{k}\mathbf{g}%
_{t}^{ij}\partial _{i}\widetilde{G}\partial _{j}\widetilde{G}%
dx-\int_{B(0,\delta )}(x^{k}\partial _{k}\widetilde{G}+\frac{1}{2}\widetilde{%
G})\mathbf{g}_{t}^{ij}\Gamma (\mathbf{g}_{t})_{ij}^{m}\partial _{m}%
\widetilde{G}dx \\
&&-\stackunder{t\rightarrow 0}{\lim }[\frac{f(x_{0})^{2}}{4\omega _{3}^{-%
\frac{2}{3}}\omega _{3}^{2}}\frac{\lambda _{t}}{6}\mu
_{t}^{-1}\int_{B(0,\delta )}\overline{u}_{t}^{6}x^{k}\partial _{k}\overline{f%
}_{t}dx] \\
&=&\delta \int_{\partial B(0,\delta )}[\frac{1}{2}\left| \nabla \widetilde{G}%
\right| _{e}^{2}-(\nabla \widetilde{G},\nu )_{e}(\nabla \widetilde{G},\nu )_{%
\widetilde{\mathbf{g}}}]d\sigma _{e} \\
&&-\frac{1}{2}\int_{\partial B(0,\delta )}\widetilde{G}(\nabla \widetilde{G}%
,\nu )_{\widetilde{\mathbf{g}}}d\sigma _{e}+\frac{\delta }{2}\int_{\partial
B(0,\delta )}\widetilde{h}\widetilde{G}^{2}d\sigma _{e}
\end{eqnarray*}
o\`{u} $\widetilde{h}=h\circ \exp _{x_{0}}$ et $\widetilde{\mathbf{g}}=\exp
_{x_{0}}^{*}\mathbf{g}$.

La diff\'{e}rence par rapport au travail d'O. Druet est donc la pr\'{e}sence
de la limite 
\[
\stackunder{t\rightarrow 0}{\lim }[\mu _{t}^{-1}\int_{B(0,\delta )}\overline{%
u}_{t}^{6}x^{k}\partial _{k}\overline{f}_{t}dx] 
\]
Mais par changement d'\'{e}chelle 
\[
\mu _{t}^{-1}\int_{B(0,\delta )}\overline{u}_{t}^{6}x^{k}\partial _{k}%
\overline{f}_{t}dx\sim \int_{B(0,R)}\widetilde{u}_{t}^{6}x^{k}\partial _{k}%
\widetilde{f}_{t}dv_{\widetilde{\mathbf{g}}_{t}}+\int_{B(0,\delta \mu
_{t}^{-1})\backslash B(0,R)}\widetilde{u}_{t}^{6}x^{k}\partial _{k}%
\widetilde{f}_{t}dv_{\widetilde{\mathbf{g}}_{t}}\,. 
\]
Or 
\[
\partial _{k}\widetilde{f}_{t}\stackrel{C_{loc}^{0}}{\rightarrow }\partial
_{k}\widetilde{f}(0) 
\]
et 
\[
\widetilde{u}_{t}\stackrel{C_{loc}^{0}}{\rightarrow }\widetilde{u} 
\]
o\`{u} $\widetilde{u}$ est radiale, telle que $\widetilde{u}_{t}\leq c%
\widetilde{u}$ et $\int_{\Bbb{R}^{n}}x^{k}\widetilde{u}<+\infty $. Par
cons\'{e}quent 
\[
\int_{B(0,R)}\widetilde{u}_{t}^{6}x^{k}\partial _{k}\widetilde{f}_{t}dv_{%
\widetilde{\mathbf{g}}_{t}}\stackunder{t\rightarrow 0}{\rightarrow }0 
\]
\`{a} $R$ fix\'{e} car $\widetilde{u}$ est radiale, et 
\[
\stackunder{R\rightarrow +\infty }{\lim }(\stackunder{t\rightarrow 0}{\lim }%
\int_{B(0,\delta \mu _{t}^{-1})\backslash B(0,R)}\widetilde{u}%
_{t}^{6}x^{k}\partial _{k}\widetilde{f}_{t}dv_{\widetilde{\mathbf{g}}%
_{t}})=0 
\]
car $\int_{\Bbb{R}^{n}}x^{k}\widetilde{u}<+\infty $ et $\mu
_{t}^{-1}\rightarrow +\infty .$

On a donc en fait 
\begin{eqnarray*}
&&\int_{B(0,\delta )}(\widetilde{h}+\frac{1}{2}x^{k}\partial _{k}\widetilde{h%
})\widetilde{G}^{2}dx-\int_{B(0,\delta )}(x^{k}\partial _{k}\widetilde{G}+%
\frac{1}{2}\widetilde{G})\partial _{j}\mathbf{g}_{t}^{ij}\partial _{i}%
\widetilde{G}dx \\
&&+\frac{1}{2}\int_{B(0,\delta )}x^{k}\partial _{k}\mathbf{g}%
_{t}^{ij}\partial _{i}\widetilde{G}\partial _{j}\widetilde{G}%
dx-\int_{B(0,\delta )}(x^{k}\partial _{k}\widetilde{G}+\frac{1}{2}\widetilde{%
G})\mathbf{g}_{t}^{ij}\Gamma (\mathbf{g}_{t})_{ij}^{m}\partial _{m}%
\widetilde{G}dx \\
&=&\delta \int_{\partial B(0,\delta )}[\frac{1}{2}\left| \nabla \widetilde{G}%
\right| _{e}^{2}-(\nabla \widetilde{G},\nu )_{e}(\nabla \widetilde{G},\nu )_{%
\widetilde{\mathbf{g}}}]d\sigma _{e} \\
&&-\frac{1}{2}\int_{\partial B(0,\delta )}\widetilde{G}(\nabla \widetilde{G}%
,\nu )_{\widetilde{\mathbf{g}}}d\sigma _{e}+\frac{\delta }{2}\int_{\partial
B(0,\delta )}\widetilde{h}\widetilde{G}^{2}d\sigma _{e}
\end{eqnarray*}
c'est \`{a} dire exactement la m\^{e}me \'{e}galit\'{e} que dans le cas $%
f=cste$. Le principe est alors de calculer un d\'{e}veloppement de chaque
terme quand $\delta \rightarrow 0$ pour trouver 
\[
\frac{\omega _{2}}{2}M_{h}(x_{0})=O(\delta ) 
\]
c'est \`{a} dire 
\[
\frac{\omega _{2}}{2}M_{h}(x_{0})=0 
\]

La fonction $f$ n'intervenant plus, il n'y a aucun changement, \`{a} partir
de l\`{a}, avec l'article d'O. Druet (rappelons quand m\^{e}me qu'il faut se
placer aux points de maximum de $f$). N\'{e}anmoins celui-ci ne
pr\'{e}sentant que peu de d\'{e}tails, et dans un souci de lisibilit\'{e} de
cette th\`{e}se, nous indiquons rapidement ces d\'{e}veloppements sous
l'hypoth\`{e}se que la m\'{e}trique $\mathbf{g}$ est plate au voisinage de $%
x_{0}$; le cas g\'{e}n\'{e}ral est identique, seulement un peu plus long.

Il reste de l'\'{e}galit\'{e} ci-dessus 
\begin{eqnarray*}
\int_{B(0,\delta )}(\widetilde{h}+\frac{1}{2}x^{k}\partial _{k}\widetilde{h})%
\widetilde{G}^{2}dx &=&\delta \int_{\partial B(0,\delta )}[\frac{1}{2}\left|
\nabla \widetilde{G}\right| _{e}^{2}-(\nabla \widetilde{G},\nu
)_{e}^{2}]d\sigma _{e} \\
&&-\frac{1}{2}\int_{\partial B(0,\delta )}\widetilde{G}(\nabla \widetilde{G}%
,\nu )_{e}d\sigma _{e}+\frac{\delta }{2}\int_{\partial B(0,\delta )}%
\widetilde{h}\widetilde{G}^{2}d\sigma _{e}\,.
\end{eqnarray*}
On rappelle alors que, en notant $r=\left| x\right| $, 
\[
\omega _{2}\widetilde{G}=\frac{1}{r}+M_{h}(x_{0})+\alpha (x) 
\]
o\`{u} $\alpha (0)=0$ et $\alpha \in C^{0,\eta }(B(0,\delta ))$ pour tout $%
0<\eta <1$, $\alpha \in C^{\infty }(B(0,\delta )\backslash \{0\})$. On a $%
\nu =\nabla r$, et on d\'{e}veloppera $\widetilde{h}$ sous la forme 
\[
\widetilde{h}=\widetilde{h}(0)+x^{k}\partial _{k}\widetilde{h}(0)+o(r)=%
\widetilde{h}(0)+O(r)\,. 
\]
On utilise alors pour $k$ entier relatif diff\'{e}rent de -3 
\[
\int_{B(0,\delta )}r^{k}=c\delta ^{k+3} 
\]
et pour $k$ entier relatif 
\[
\int_{\partial B(0,\delta )}r^{k}d\sigma _{e}=\omega _{2}\delta ^{k+2}\,. 
\]
Ainsi 
\[
\frac{\delta }{2}\int_{\partial B(0,\delta )}\widetilde{h}\widetilde{G}%
^{2}d\sigma _{e}=c_{1}\delta +o(\delta )=O(\delta )\,. 
\]
Comme $(\widetilde{h}+\frac{1}{2}x^{k}\partial _{k}\widetilde{h})\widetilde{G%
}^{2}$ est continue en 0 
\[
\int_{B(0,\delta )}(\widetilde{h}+\frac{1}{2}x^{k}\partial _{k}\widetilde{h})%
\widetilde{G}^{2}dx=O(\delta )\,. 
\]
On montre que ($\nabla r,\nabla \alpha )=O(r^{\eta -1})$ et cela permet
d'obtenir 
\[
-\frac{1}{2}\int_{\partial B(0,\delta )}\widetilde{G}(\nabla \widetilde{G}%
,\nu )_{e}d\sigma _{e}=\frac{\omega _{2}}{2}M_{h}(x_{0})+\frac{\omega _{2}}{2%
}\delta ^{-1}+O(\delta ) 
\]
et 
\[
\delta \int_{\partial B(0,\delta )}[\frac{1}{2}\left| \nabla \widetilde{G}%
\right| _{e}^{2}-(\nabla \widetilde{G},\nu )_{e}^{2}]d\sigma _{e}=-\frac{%
\omega _{2}}{2}\delta ^{-1}+O(\delta )\,. 
\]
En ajoutant ces quatre d\'{e}veloppements, on trouve 
\[
\frac{\omega _{2}}{2}M_{h}(x_{0})+O(\delta )=0 
\]
et donc en faisant tendre $\delta $ vers 0 
\[
M_{h}(x_{0})=0 
\]
Autrement dit, si les $u_{t}$ se concentrent en un point $x_{0}$,
n\'{e}cessairement $M_{h}(x_{0})=0$. Mais dans ce cas, pour toute fonction $%
h^{\prime }\leq h,\,h^{\prime }\neq h$, on a par le principe du maximum 
\[
M_{h^{\prime }}(x_{0})>0 
\]
et donc, par le premier th\'{e}or\`{e}me (th\'{e}or\`{e}me 6), $h^{\prime }$
ne peut \^{e}tre faiblement critique pour $f$ et $\mathbf{g}$; donc ($%
h^{\prime },f,\mathbf{g})$ est sous-critique. Par d\'{e}finition $(h,f,%
\mathbf{g})$ est donc critique.

Le th\'{e}or\`{e}me 8 d\'{e}coule imm\'{e}diatement de cette
d\'{e}monstration.

\chapter{Remarques sur le cas limite et le cas d\'{e}g\'{e}n\'{e}r\'{e}.
Quelques questions...}
\pagestyle{myheadings}\markboth{\textbf{Remarques et questions.}}
{\textbf{Remarques et questions.}}
\section{Cas limite et cas d\'{e}g\'{e}n\'{e}r\'{e}}

Dans le th\'{e}or\`{e}me 1 d\'{e}montr\'{e} au chapitre 3, et dans la
plupart des th\'{e}or\`{e}mes qui en d\'{e}coulent, nous faisons deux
hypoth\`{e}ses importantes:

(H1): en rappelant que si ($h,f,\mathbf{g}$) est un triplet faiblement
critique, n\'{e}cessairement, au points de maximum de $f$, $h(P)\geq \frac{%
n-2}{4(n-1)}S_{\mathbf{g}}(P)-\frac{(n-2)(n-4)}{8(n-1)}\frac{\bigtriangleup
_{\mathbf{g}}f(P)}{f(P)}$, nous avons fait dans nos th\'{e}or\`{e}mes
l'hypoth\`{e}se que cette in\'{e}galit\'{e} \'{e}tait \textit{stricte:} 
\[
h(P)>\frac{n-2}{4(n-1)}S_{\mathbf{g}}(P)-\frac{(n-2)(n-4)}{8(n-1)}\frac{%
\bigtriangleup _{\mathbf{g}}f(P)}{f(P)} 
\]

(H2) ou (H$_{f}$): Nous avons suppos\'{e} que le Hessien de la fonction $f$
\'{e}tait non d\'{e}g\'{e}n\'{e}r\'{e} en ses points de maximum.

Nous avons évidemment essay\'{e} de nous passer de ces hypoth\`{e}ses, sans
succ\'{e}s... Nous voulons proposer quelques remarques et explications
concernant ces deux hypoth\`{e}ses.

Tout d'abord, en gardant l'hypoth\`{e}se (H1), nos th\'{e}or\`{e}mes 1, 1'
et 1'' restent valables si l'on suppose que $\bigtriangleup _{\mathbf{g}%
}f(P)=0$ au points de maximum, ce qui implique \'{e}videmment que le Hessien
est d\'{e}g\'{e}n\'{e}r\'{e} (il suffit de remarquer qu'avec les notations
du chapitre 3, partie 3.3, $A_{t}^{1}\leq 0$). Nous sommes donc dans la
situation o\`{u} les th\'{e}or\`{e}mes fonctionnent dans deux cas
extr\`{e}mes, Hessien non d\'{e}g\'{e}n\'{e}r\'{e} et Laplacien nul.

En ce qui concerne (H1), E. Hebey et M. Vaugon $\left[ 20\right] $ ont montr%
\'{e} dans le cadre de leur \'{e}tude sur $B_{0}(\mathbf{g)}$, qui
correspond pour nous au cas $f=cste$, le r\'{e}sultat suivant:

Supposons que la vari\'{e}t\'{e} ($M,\mathbf{g)}$ est de dimension $\geq 7$,
et soit $h$ une fonction ctirique (pour $\mathbf{g}$ et 1). Notons $T=\{x\in
M/\,h(x)=\frac{n-2}{4(n-1)}S_{\mathbf{g}}(x)\}$. Rappelons que dans ce cadre
on a toujours $h\geq \frac{n-2}{4(n-1)}S_{\mathbf{g}}$ sur tout $M$ et que
l'hypoth\`{e}se (H1) devient $h>\frac{n-2}{4(n-1)}S_{\mathbf{g}}$. On
suppose alors que pour tout point $x$ de $T$:

1: le tenseur de Weyl est nul sur un voisinage de $x$, et

2: $\nabla ^{2}(h-\frac{n-2}{4(n-1)}S_{\mathbf{g}})$ est non
d\'{e}g\'{e}n\'{e}r\'{e} en $x$

Alors $h$ a des fonctions extr\'{e}males.

On remarque alors que l'on imm\'{e}diatement \`{a} partir de leur
d\'{e}monstration le r\'{e}sultat, un peu artificiel, suivant:

\textit{Supposons que la vari\'{e}t\'{e} (}$M,\mathbf{g)}$\textit{\ est de
dimension }$\geq 7$\textit{, et soit (}$h,f,\mathbf{g)}$\textit{\ un triplet
critique. Notons }$T_{f}=\{x\in M/\,f(x)=Maxf\,\,et\,\,h(x)=\frac{n-2}{4(n-1)%
}S_{\mathbf{g}}(x)-\frac{(n-2)(n-4)}{8(n-1)}\frac{\bigtriangleup _{\mathbf{g}%
}f(P)}{f(P)}\}$\textit{. On suppose alors que }$T_{f}$\textit{\ n'est pas
dense dans }$M$\textit{\ et que pour tout point }$x$\textit{\ de }$T_{f}$%
\textit{:}

\textit{1: le tenseur de Weyl est nul sur un voisinage de }$x$\textit{, }

\textit{2: }$\nabla ^{2}(h-\frac{n-2}{4(n-1)}S_{\mathbf{g}})$\textit{\ est
non d\'{e}g\'{e}n\'{e}r\'{e} en }$x$\textit{,}

\textit{3: }$\bigtriangleup _{\mathbf{g}}f(x)=0$ si $x\in T_{f}$, \textit{et
on suppose de plus que le hessien de }$f$\textit{\ est non d\'{e}g\'{e}n\'{e}%
r\'{e} aux points de maximum qui ne sont pas dans }$T_{f}$\textit{.}

\textit{Alors (}$h,f,\mathbf{g)}$\textit{\ a des fonctions extr\'{e}males.}

L'int\'{e}r\^{e}t de ce r\'{e}sultat est de montrer qu'il existe des triplet
critique ($h,f,\mathbf{g)}$, avec $f$ non constante, qui ne v\'{e}rifie pas
(H1) et qui ont des fonctions extr\'{e}males. Indiquons très rapidement
le sch\'{e}ma de leur d\'{e}monstration.

On choisit une fonction plateau $\theta $ dont le support est disjoint de
l'ensemble des points de maximum de $f$. Soit $h_{t}=h-t\theta $; $h_{t}$
est sous-critique pour tout $t$ tendant vers 0. On a alors une suite de
solutions minimisantes $u_{t}$ associ\'{e}es aux \'{e}quations $%
\bigtriangleup _{\mathbf{g}}u_{t}+h_{t}u_{t}=fu_{t}^{2^{*}-1}$. L\`{a}
encore la suite $u_{t\text{ }}$ converge dans $H_{1}^{2}$ vers une fonction $%
u$. Si $u>0$ on a une solution minimisante. Si $u=0$, toute l'\'{e}tude sur
les ph\'{e}nom\`{e}nes de concentration du chapitre 3 est valable. Le point
nouveau utilis\'{e} par E. Hebey et M. Vaugon est une sorte d'am\'{e}%
lioration de la concentration $L^{2}$. Notons qu'il n'existe qu'un seul
point de concentration, not\'{e} $x_{0}$, et que c'est un point de maximum
de $f$. E. Hebey et M. Vaugon montrent que si $\dim M\geq 7$ alors pour tout
rayon $\delta >0$:
\begin{equation}
\stackunder{t\rightarrow 0}{\lim }\frac{\int_{M\backslash B(x_{0},\delta
)}d_{\mathbf{g}}(x,x_{0})^{2}u_{t}^{2}dv_{\mathbf{g}}}{\int_{B(x_{0},\delta
)}d_{\mathbf{g}}(x,x_{0})^{2}u_{t}^{2}dv_{\mathbf{g}}}=0  
\end{equation}
Si le Hessien de $f$ est non d\'{e}g\'{e}n\'{e}r\'{e} en $x_{0}$, on reprend
la d\'{e}monstration du chapitre 3 pour aboutir \`{a} une contradiction. On
peut donc supposer que $x_{0}\in T_{f}$. Par l'hypoth\`{e}se 1, on peut
identifier $B(x_{0},\delta )$ munie de la m\'{e}trique $\mathbf{g}$ \`{a} la
boule euclidienne $B(0,\delta )$ en identifiant $x_{0}$ \`{a} $0$. Mais
alors en reprenant les calculs du chapitre 3, partie 3.3, c'est 
\`{a} dire l'in\'{e}galit\'{e} de Sobolev euclidienne dans laquelle on
``injecte'' l'\'{e}quation v\'{e}rifi\'{e}e par $u_{t}$, on obtient 
\[
\int_{B(x_{0},\delta )}h_{t}(\eta u_{t})^{2}\leq A_{t}+\frac{c}{\delta ^{2}}%
\int_{B(x_{0},\delta )\backslash B(x_{0},\delta /2)}u_{t}^{2}\leq \frac{c}{%
\delta ^{2}}\int_{B(x_{0},\delta )\backslash B(x_{0},\delta /2)}u_{t}^{2}
\]
car $B_{t}=C_{t}=0$ puisque la m\'{e}trique est plate au voisinage de $x_{0}$%
, et puisque, c'est l\`{a} le seul point \`{a} remarquer si $f$ est non
constante, $A_{t}\leq 0$. L'hypoth\`{e}se 2 signifie que $h_{t}$, qui est 
\'{e}gale \`{a} $h$ sur $B(x_{0},\delta )$ si on choisit $\delta $ assez
petit, a un minimum strict en $x_{0}$ et que ce minimum vaut 0=$S_{\mathbf{g}%
}$. Un d\'{e}veloppement limit\'{e} de $h_{t}$ en $x_{0}$ montre que 
\[
\int_{B(x_{0},\delta )}h_{t}(\eta u_{t})^{2}\geq \lambda
\int_{B(x_{0},\delta )}\left| x\right| ^{2}u_{t}{}^{2}
\]
pour un r\'{e}el $\lambda >0.$ Mais alors 
\[
\lambda \int_{B(x_{0},\delta )}\left| x\right| ^{2}u_{t}{}^{2}\leq \frac{c}{%
\delta ^{2}}\int_{B(x_{0},\delta )\backslash B(x_{0},\delta
/2)}u_{t}^{2}\leq \frac{c}{\delta ^{4}}\int_{B(x_{0},\delta )\backslash
B(x_{0},\delta /2)}\left| x\right| ^{2}u_{t}^{2}
\]
ce qui contredit (8.1).

Notre sentiment est alors le suivant. Reprenons la remarque finale du
chapitre 3: le principe est de faire un d\'{e}veloppement limit\'{e} en
fonction de $\mu _{t}$ des deux membres de l'in\'{e}galit\'{e} issue de
celle de Sobolev, au sens o\`{u}, apr\`{e}s changement d'\'{e}chelle, on
\'{e}value les int\'{e}grales en fonction de ce param\`{e}tre fondamental.
On obtient ainsi \`{a} partir de l'in\'{e}galit\'{e} de Sobolev dans
laquelle on introduit l'\'{e}quation v\'{e}rifi\'{e}e par $\overline{u}_{t}$%
: 
\[
h(x_{0})\mu _{t}^{2}\leq (\frac{n-2}{4(n-1)}S_{\mathbf{g}}(x_{0})-\frac{%
(n-2)(n-4)}{8(n-1)}\frac{\bigtriangleup _{\mathbf{g}}f(x_{0})}{f(x_{0})})\mu
_{t}^{2}+o(\mu _{t}^{2}) 
\]
Or en reprenant le d\'{e}tails des calculs fait au chapitre 3, on constate
que dans les d\'{e}veloppements limit\'{e}s effectu\'{e}s, les coefficients
des termes d'ordre 2 en $\mu _{t}$ sont les seuls qui soient
intrins\`{e}ques, c'est \`{a} dire invariants dans les changements de cartes
( exp$_{x_{t}}^{-1}$). En effet, au premier membre o\`{u} l'on d\'{e}veloppe 
$h$, on obtient (apr\`{e}s $h(x_{0}))$ les d\'{e}riv\'{e}es premi\`{e}res $%
\partial _{i}h(x_{0})$, tandis qu'au second membre (apr\`{e}s $%
\bigtriangleup _{\mathbf{g}}f(x_{0}))$ ce sont les d\'{e}riv\'{e}es d'ordre
3 qui apparaissent; on peut faire le m\^{e}me constat pour le
d\'{e}veloppement de la m\'{e}trique. Maintenant, si l'on regarde les
hypoth\`{e}ses du r\'{e}sultat d'E.\ Hebey et de M.\ Vaugon, on constate
qu'elles reviennent \`{a} assurer la nullit\'{e} des termes d'ordre 2 en $%
\mu _{t}$ et \`{a} faire porter le poids du d\'{e}veloppement sur les termes
d'ordre 4 (la situation technique du r\'{e}sultat permet de ne pas faire
apparaitre de termes d'ordre 3), la contradiction s'obtenant gr\^{a}ce \`{a}
la condition de non dégénérécence de $\nabla ^{2}(h-\frac{n-2}{%
4(n-1)}S_{\mathbf{g}})$, condition intrins\`{e}que, et portant sur les
termes d'ordre 4. Ces remarques \'{e}tant faites, on peut alors constater en
ce qui concerne l'hypoth\`{e}se (H2), bien que cela soit moins \'{e}vident,
que d'apr\`{e}s la th\'{e}orie de Morse \'{e}l\'{e}mentaire, seuls les
points critiques non d\'{e}g\'{e}n\'{e}r\'{e}s ont un sens intrins\`{e}que.
Il semble donc, si l'on veut se passer de H1 ou H2, qu'il faille soit
trouver un autre param\`{e}tre intrins\`{e}que \`{a} faire apparaitre, soit
trouver une m\'{e}thode radicalement diff\'{e}rente pour prouver, ou
infirmer, les th\'{e}or\`{e}mes \'{e}tablis dans notre travail. Remarquons
que l'on pourrait \^{e}tre tent\'{e} d'utiliser l'identit\'{e} de Pohozahev
comme en dimension 3 pour r\'{e}soudre le probl\`{e}me, mais en fait on
tombe sur la m\^{e}me difficult\'{e}. En effet, si l'on ins\`{e}re
l'\'{e}quation dont est solution $u_{t}$ dans cette identit\'{e}, on obtient
apr\`{e}s quelques calculs, et, pour simplifier, en supposant que la
m\'{e}trique est plate au voisinage de $x_{0}$: 
\[
\int_{B(0,\delta )}\eta ^{2}(2\overline{h}_{t}+x^{k}\partial _{k}h_{t})%
\overline{u}_{t}^{2}dx+\int_{B(0,\delta )}\alpha (\eta )\overline{u}%
_{t}^{2}dx+\lambda _{t}\int_{B(0,\delta )}\beta (\eta )\overline{f}_{t}%
\overline{u}_{t}^{2^{*}}dx\leq 0 
\]
et l'on voit que les d\'{e}veloppements limit\'{e}s de $h$ ou $f$ vont faire
apparaitre les m\^{e}mes termes au m\^{e}mes ordres.

Mais l'on peut dire plus en ce qui concerne l'hypoth\`{e}se sur le Hessien
de $f$ aux points de maximum.\ Comme nous l'avons vu au chapitre 3, cette
hypoth\`{e}se et la d\'{e}monstration du th\'{e}or\`{e}me 1 sont
\'{e}troitement li\'{e}es \`{a} la ''seconde in\'{e}galit\'{e}
fondamentale'' 
\[
\frac{d_{\mathbf{g}}(x_{t},x_{0})}{\mu _{t}}\leq C 
\]
que l'on peut obtenir pour une suite ($u_{t}$) de solutions d'\'{e}quations $%
\triangle _{\mathbf{g}}u_{t}+h_{t}.u_{t}=\lambda _{t}f.u_{t}^{\frac{n+2}{n-2}%
}$ se concentrant. Or, on peut construire une suite de fonctions solutions
d'\'{e}quation de ce type v\'{e}rifiant tous les ph\'{e}nom\`{e}nes de
concentration, mais ne v\'{e}rifiant pas cette estim\'{e}e. Consid\'{e}rons
en effet la sph\`{e}re $S^{n}$ munie de m\'{e}trique standard $\mathbf{s}$.
Si l'on reformule les r\'{e}sultats connus (c.f. par exemple [17]), il existe une unique fonction
critique pour 1 et $\mathbf{s}$, \`{a} savoir 
\[
h=\frac{n-2}{4(n-1)}S_{\mathbf{s}}=\frac{n-2}{4(n-1)} 
\]
et cette fonction critique poss\`{e}de comme fonctions extr\'{e}males d'une
part les constantes, et d'autre part les fonctions de la forme 
\[
u=a(b-\cos r)^{-\frac{n-2}{2}} 
\]
o\`{u} $a\neq 0$, $b>1$, et $r$ est la distance g\'{e}od\'{e}sique \`{a} un
point fix\'{e} de $S^{n}$. Considérons alors sur $S^{n}$ une suite de
points $x_{t}$ convergeant vers un point $x_{0}$, et posons 
\[
u_{t}=\mu _{t}^{\frac{n-2}{2}}(\mu _{t}^{2}+1-\cos r_{t})^{-\frac{n-2}{2}} 
\]
o\`{u} $r_{t}(x)=d_{\mathbf{s}}(x,x_{t})$ et $\mu _{t}$ est une suite de
r\'{e}els convergeant vers 0. Alors 
\[
\int_{M}u_{t}^{2^{*}}dv_{\mathbf{s}}=1 
\]
et on obtient ainsi une suite de solutions de l'\'{e}quation 
\[
\triangle _{\mathbf{s}}u_{t}+\frac{n-2}{4(n-1)}.u_{t}=K(n,2)^{-2}u_{t}^{%
\frac{n+2}{n-2}} 
\]
où le second membre ne vérifie pas ($H_{f}$), avec 
\[
Sup_{M}u_{t}=u_{t}(x_{t})=\mu _{t}^{-\frac{n-2}{2}}. 
\]
Cette suite se concentre et v\'{e}rifie toutes les propri\'{e}t\'{e}s a/
\`{a} d/ expos\'{e}es dans la quatri\`{e}me partie du chapitre 3, ceci
quelque soit le choix de la suite $x_{t}\rightarrow x_{0}$ et de la suite $%
\mu _{t}\rightarrow 0$. Par sym\'{e}trie sph\'{e}rique, on peut facilement
trouver deux suites ($x_{t}$) et ($\mu _{t}$) telles que 
\[
\frac{d_{\mathbf{s}}(x_{t},x_{0})}{\mu _{t}}\rightarrow +\infty 
\]
en prenant par exemple $\mu _{t}=d_{\mathbf{s}}(x_{t},x_{0})^{2}.$

Encore une fois, il semble que l'hypoth\`{e}se faite sur $f$ ''fixe'' la
position du point de concentration, et ainsi ''impose'' une vitesse de
convergence \`{a} la suite ($x_{t}$).

\section{Quelques questions et perspectives.}

Nous avons vu que l'\'{e}tude des \'{e}quations $\bigtriangleup _{\mathbf{g}%
}u+hu=fu^{2^{*}-1}$ \'{e}tait li\'{e}e \`{a} celle des meilleures constantes
dans les inclusions de Sobolev de $H_{1}^{2}$ dans $L^{\frac{2n}{n-2}}$. De
la m\^{e}me mani\`{e}re, l'\'{e}tude des inclusions de Sobolev de $H_{1}^{p}$
dans $L^{\frac{pn}{n-p}}$, o\`{u} $\frac{pn}{n-p}\,$est l'exposant critique,
et des meilleures constantes associ\'{e}es, passe par l'\'{e}tude des \'{e}%
quations de la forme 
\[
\bigtriangleup _{p}u+hu=fu^{\frac{pn}{n-p}-1}
\]
o\`{u} $\bigtriangleup _{p}u=-\nabla (\left| \nabla u\right| _{\mathbf{g}%
}^{p-2}\nabla u)$ est le p-laplacien; voir par exemple O. Druet, E. Hebey $%
\left[ 11\right] $ et Z. Faget $\left[ 14\right] $. L\`{a} aussi les m\'{e}%
thodes variationnelles sont \`{a} la base de l'\'{e}tude: la fonctionnelle
consid\'{e}r\'{e}e est 
\[
I(u)=\int \left| \nabla u\right| _{\mathbf{g}}^{p}+\int hu^{p}
\]
d'o\`{u} le lien avec l'inclusion de Sobolev 
\[
(\int u^{\frac{pn}{n-p}})^{\frac{n-p}{n}}\leq K(n,p)\int \left| \nabla
u\right| _{\mathbf{g}}^{p}+B\int u^{p}
\]
o\`{u} $K(n,p)$ est la meilleure constante associ\'{e}e. Le r\'{e}sultat de d%
\'{e}part est encore le suivant: Si
\[
\stackunder{\int u^{\frac{pn}{n-p}}=1}{Inf\,}\,I(u)<K(n,p)^{-1}(Supf)^{-%
\frac{n-p}{n}}
\]
alors l'\'{e}quation a une solution minimisante $u>0$ (sachant que l'in\'{e}%
galit\'{e} large est toujours vraie). On voit donc qu'il est facile d'\'{e}%
tendre la d\'{e}finition des fonctions critiques \`{a} cette situation. Il
serait alors int\'{e}ressant de savoir si l'on peut trouver des r\'{e}%
sultats analogues \`{a} ceux de notre travail dans ce cadre.

Rappelons \'{e}galement une question que nous avons soulev\'{e}e \`{a}
l'issue du chapitre 6:

\begin{center}
$f$\textit{\ \'{e}tant donn\'{e}e, existe-t-il des fonctions critiques
constantes?}
\end{center}

Cela donnerait en quelque sorte une ''seconde meilleure constante $B_{0}(%
\mathbf{g},f)$'' li\'{e}e \`{a} $f$.

Terminons sur une question qui s'impose \`{a} l'issue de ce travail:

\begin{center}
\textit{-Pour une fonction }$h$\textit{\ quelconque sur }$M$\textit{,
existe-t-il des solutions (non minimisantes) \`{a} }$\bigtriangleup _{%
\mathbf{g}}u+hu=fu^{2^{*}-1}$\textit{\ ?}
\end{center}

Nous avons en effet vu que cette \'{e}quation avait des solutions
(minimisantes) si $h$ est sous-critique et si $h$ est critique avec les
hypoth\`{e}ses (H1) et (H2). Par contre, les th\'{e}ories variationnelles ne
donnent aucune r\'{e}ponse si $h$ est sup\'{e}rieure (et diff\'{e}rente)
\`{a} une fonction critique, ou si $\bigtriangleup _{\mathbf{g}}+h$ n'est
pas coercif. Dans ce cas, si des solutions existent, elles ne peuvent
\^{e}tre minimisantes. Il faut donc employer d'autres m\'{e}thodes pour
\'{e}tudier ces cas. Voir A. Bahri [4] qui traite le cas $f=cste$ et $3\leq dimM\leq6$.

\chapter{Abridged English Version}
\pagestyle{myheadings}\markboth{\textbf{Abridged English Version.}}
{\textbf{Abridged English Version.}}

\section{Introduction}
In the beginning was the Yamabe problem:\\

\textbf{Yamabe problem}\textit{: Given a compact Riemannian manifold }$(M,\textbf{g})$\textit{\ of dimension }$%
n\geq 3$,\textit{\ does there exist a metric} $\textbf{g'}$ \textit{conformal to }$\textbf{g}$\textit{\ having constant scalar curvature?}\\

If we write $\textbf{g'}=u^\frac{4}{n-2}.\textbf{g}$ where $u>0$ is a smooth function on $M$, the scalar 
curvatures are linked by the partial differential equation :
$$\triangle _{\textbf{g}}u+\frac{n-2}{4(n-1)}S_{\textbf{g}}.u=\frac{n-2}{4(n-1)}S_{\textbf{g'}} .u^{\frac{n+2}{n-2}} $$
where $S_{\textbf{g}}$ is the scalar curvature of $\textbf{g}$ and where $\triangle _{\textbf{g}}=-\nabla^{i}\nabla_{i}$
 is the Riemannian laplacian of $\g$.

To solve the Yamabe problem, one therefore has to prove the existence of a solution $u>0$ to this partial differential equation
when $S_{\textbf{g'}}$ is a constant. More generaly, the prescribed curvature problems, which consist in deciding, 
given a smooth function $f$ on $M$, if $f$ is the scalar curvature of a metric conformal to $\g$, come down
to prove the existence of a positive smooth solution $u$ to the above equation when $S_{\g'}$ is replaced by $f$.

These problems launched the study of elliptic PDE on compact Riemannian manifolds of the form 
$$(E_{h,f,\mathbf{g}}): \triangle _{\textbf{g}}u+h.u=f .u^{\frac{n+2}{n-2}} $$
In all this paper $M$ will be a compact Riemannian manifold of dimension $n\geqslant 3$, we will use 
the letter $\g$ or $\g'$ to denote a Riemannian metric on $M$; $h$ and $f$ will always be smooth functions on 
$M$. We will always suppose the functions to be smooth, however in the definitions and in most of the 
theorems, continuity is in general sufficient. Beside, we will keep these notations, letter $\g$ for the metrics,
letter $h$ for the function on the left of equation $E_{h,f,\g}$, (defining the opperator $ \triangle _{\textbf{g}}+h$);
and letter $f$ for the function on the right of the equation; the unknown function will be designated by $u$.

One of the possible methods to study these equation is the use of variational methods, which have the 
advantage of giving minimizing solutions, or solution of minimal energy. If one multiply equation $\E$ 
by $u$ and integrate over $M$, one gets
$$\int_{M}\left| \nabla u\right| _{\mathbf{g}%
}^{2}dv_{\mathbf{g}}+\int_{M}h.u{{}^{2}}dv_{\mathbf{g}}= \int_{M}f\left| u\right| ^{\frac{2n}{n-2}%
}dv_{\mathbf{g}}$$
The variational methods therefore lead to consider the functional
$$I_{h,\mathbf{g}}(w)=\int_{M}\left| \nabla w\right| _{\mathbf{g}%
}^{2}dv_{\mathbf{g}}+\int_{M}h.w{{}^{2}}dv_{\mathbf{g}} $$
defined for $w\in H_{1}^{2}(M)$, the Sobolev space of $L^{2}$ functions whose gradient is also in $L^{2}$, 
and the minimum of this functional
\[
\lambda _{h,f,\mathbf{g}}=\stackunder{w\in \mathcal{H}_{f}}{\inf }I_{h,%
\mathbf{g}}(w) 
\]
on the set
\[
\mathcal{H}_{f}=\{w\in H_{1}^{2}(M)/\int_{M}f\left| w\right| ^{\frac{2n}{n-2}%
}dv_{\mathbf{g}}=1\}. 
\]
The Euler equation associated with the minimization problem of this functional by a function $u$ such that
\[
I_{h,\mathbf{g}}(u)=\stackunder{w\in \mathcal{H}_{f}}{\inf }I_{h,\mathbf{g}%
}(w) 
\]
is indeed exactly
\[
(E_{h,f,\mathbf{g}}):\triangle _{\mathbf{g}}u+hu=\lambda _{h,f,\mathbf{g}%
}.f.u^{\frac{n+2}{n-2}} 
\]
where $\lambda _{h,f,\mathbf{g}}\,$ appears as a normalizing constant due to the condition
\[
\int_{M}f\left| u\right| ^{\frac{2n}{n-2}}dv_{\mathbf{g}}=1. 
\]
It is sometimes usefull to consider the functional 
\[
J_{h,f,\mathbf{g}}(w)=\frac{\int_{M}\left| \nabla w\right| _{\mathbf{g}%
}^{2}dv_{\mathbf{g}}+\int_{M}h.w{{}^{2}}dv_{\mathbf{g}}}{\left(
\int_{M}f\left| w\right| ^{\frac{2n}{n-2}}dv_{\mathbf{g}}\right) ^{\frac{n-2%
}{n}}} 
\]
and the subset of $H_{1}^{2}(M)$ where it is defined
\[
\mathcal{H}_{f}^{+}=\{w\in H_{1}^{2}(M)/\int_{M}f\left| w\right| ^{\frac{2n}{%
n-2}}dv_{\mathbf{g}}>0\}. 
\]
One then consider the minimisation problem by a function $u$ such that
\[
J_{h,f,\mathbf{g}}(u)=\stackunder{w\in \mathcal{H}_{f}^{+}}{\inf }J_{h,f,\mathbf{g}}(w), 
\]
the Euler equation being identical but without the normalizing constant. This functional sometimes 
present the advantage of being homogeneous in the sense that $J_{h,f,\mathbf{g}}(c.w)=J_{h,f,%
\mathbf{g}}(w)$ for any constant $c.$ One therefore see that
\[
\stackunder{w\in \mathcal{H}_{f}}{\inf }I_{h,\mathbf{g}}(w)=\stackunder{w\in 
\mathcal{H}_{f}^{+}}{\inf }J_{h,f,\mathbf{g}}(w)=\lambda _{h,f,\mathbf{g}} 
\]
This functional $J$ also has the particularity, when $h=\frac{n-2}{4(n-1)}S_{\mathbf{g}},$ 
of being invariant by conformal changes of metrics; it is therefore especially usefull when studying 
problems of prescribed scalar curvatures. We shall mostly use  $I_{h,\mathbf{g}}$ and $\mathcal{H}_{f}$, 
but for some problems $J_{h,f,\mathbf{g}}$ will prove to be more convenient when we shall want
 to avoid the constraint $\int_{M}f\left| u\right| ^{\frac{2n}{n-2}}dv_{\mathbf{g}}=1.$

We will say that a function $u\in H_{1}^{2}(M)$ is a solution of minimal energy, or a minimizing 
solution, if either $I_{h,\mathbf{g}}(u)=\lambda _{h,f,\mathbf{g}}$ with $\int_{M}fu^\frac{2n}{n-2}=1$, 
or $J_{h,f,\mathbf{g}}(u)=\lambda _{h,f,\mathbf{g}} $.
Then, up to multiplying it by a constant, $u$ is stricly positive and smooth, and it is a solution of
$$(E_{h,f,\mathbf{g}}):\triangle _{\mathbf{g}}u+hu=\lambda _{h,f,\mathbf{g}}.f.u^{\frac{n+2}{n-2}} $$
with or without the normalizing constant which can always be supressed just by multipliying again $u$ 
by a constant. Please, note that we will use these notations $\E$ and $\lambda _{h,f,\mathbf{g}}$ 
throughout all this article.

 Th. Aubin discoverded a very important relation between equation $\E$ and the notion of best constant 
 in the Sobolev imbedding theorems. Remember that the inclusion of $H_{1}^{2}(M)$ in $L^p(M)$ is 
 compact for $p<\frac{2n}{n-2}$ and only continuous for $p=\frac{2n}{n-2}$ which is called the 
 critical exponant for the Sobolev imbeddings and will be noted $2^{*}=\frac{2n}{n-2}$. The continuous 
 imbedding $H_{1}^{2}(M) \subset L^{2^{*}}(M)$ is expressed by the existence of two positive constants 
 $A$ and $B$ such that :
\begin{equation}
\forall u\in H_{1}^{2}(M):\,\left( \int_{M}\left| u\right| ^{\frac{2n}{n-2}}dv_{\mathbf{g}}\right) ^{%
\frac{n-2}{n}}\leq A\int_{M}\left| \nabla u\right| _{\mathbf{g}}^{2}dv_{%
\mathbf{g}}+B\int_{M}u{{}^{2}}dv_{\mathbf{g}}
\end{equation}
The best first constant is the minimum $A$ that one can put in (1) such that there exist $B$ with (1) still true. 
It was proved by E. Hebey and M. Vaugon [19] that this minimum is attained, and its value is known to 
be the same as for the sharp euclidean Sobolev inequality,
\[
A_{min}=K(n,2){{}^{2}}=\frac{4}{n(n-2)\omega _{n}^{\frac{2}{n}}}
\]
where $\omega _{n}$ is the volume of the unit sphere of dimension $n$. One then take
 $B_{0}(\mathbf{g})$ to be the minimum $B$ such that (1) remains true with $A_{min}$; it is proved that $B_{0}(%
\mathbf{g})<+\infty $ \cite{H-V 1}. The inequality: $\forall u\in H_{1}^{2}(M)$
\begin{equation}
\left( \int_{M}\left| u\right| ^{\frac{2n}{n-2}}dv_{\mathbf{g}}\right) ^{%
\frac{n-2}{n}}\leq K(n,2){{}^{2}}\int_{M}\left| \nabla u\right| _{\mathbf{g}%
}^{2}dv_{\g}+B_{0}(\mathbf{g})\int_{M}u{{}^{2}}dv_{\mathbf{g}} 
\end{equation}
is then sharp with respect to both the first and second constants, in the sense that none of them can be 
lowered. If the value of the best constant $A_{min}=K(n,2)^{2}$ is known and independent of the 
manifold $(M,\mathbf{g})$, on the other hand, $B_{0}(\mathbf{g})$, as the notation indicates, depends 
on the geometry and its study is difficult; it is for this purpose that "critical functions" were introduced 
by E.Hebey and M.Vaugon [20]. When there shall be no risk of confusion, these constants will be 
denoted by $K$ et $B_{0}.$

As a remark, note that because of the compacity of the inclusion $H_{1}^{2}(M) \subset L^{p}(M)$ for 
$p<2^{*}$, standard variational methods and elliptic theory give rapidly existence of minimizing solutions 
of the equation $\triangle _{\mathbf{g}}u+hu=f.u^{p-1}$ when $\triangle _{\mathbf{g}}+h$ is a 
coercive operator. The case $p=2^{*}$ is therefore already a limit case. (Very little is known for $p>2^{*}$ 
without additional hypothesis, like e.g. invariance by symetry, see [15].)

The best constants in the Sobolev embedding appeared in the study of equations $\E$ when Th. Aubin 
proved the following theorem:
\begin{theor}
[Aubin] For any Riemannian manifold $(M,\g)$ of dimension $n\geqslant 3$, any function $h$ such that 
$\triangle _{\mathbf{g}}+h$ is a coercive operator, and any function $f$ such that 
$\underset{M}{Sup} f>0$, 
one always has 
$$\lambda _{h,f,\mathbf{g}}\leq \frac{1}{K(n,2){{}^{2}}(\stackunder{M}{Sup}f)^{\frac{n-2}{n}}}.$$
Furthermore, if this inequality is strict, then there exists a minimizing solution for $\E$.
\end{theor}
 
This theorem is the starting point of all this work. It proves the existence of minimizing solutions 
 to equation $\E$ under the hypothesis:
 $$\lambda _{h,f,\mathbf{g}}< \frac{1}{K(n,2){{}^{2}}(\stackunder{M}{Sup}f)^{\frac{n-2}{n}}}.$$
Our work is essentially concerned with the problem of the existence of minimizing solutions 
to these equations $\E$ in the "critical case" where
$$\lambda _{h,f,\mathbf{g}}=\frac{1}{K(n,2){{}^{2}}(\stackunder{M}{Sup}f)^{\frac{n-2}{n}}},$$
problem which is normally not solved by variational methods. It is for the study of this problem 
that we are now going to define the "critical functions".
\\

Let us first review the datas:

\textbf{Datas:} Throughout this article, $(M,\g)$ will be a compact Riemannian manifold of dimension 
$n\geqslant 3$. We let $f:M\rightarrow \R$ be a fixed smooth function such that $\underset{M}{Sup} f>0$.
 Let also $h:M\rightarrow \R$ be a smooth function with the additional hypothesis that the operator 
 $\triangle _{\mathbf{g}}+h$ is coercive if $f$ is not positive on all of $M$. (Remember that continuity 
 of $h$ and $f$ is sufficient in the definitions and in most of the theorems. Also, if $f\leqslant0$ on $M$, 
 classical variational methods already give a lot of results for the existence of solutions; therefore $Sup f>0$ 
 is the most interesting case.)
 
 \begin{defi}
 With these datas, and with the above notations, we say that:
 \begin{itemize}
\item  $h$ is weakly critical for $f$ and $\mathbf{g}$ if $\lambda _{h,f,\mathbf{g}}=\frac{1}{K(n,2){{}^{2}}(%
\stackunder{M}{Sup}f)^{\frac{n-2}{n}}}$

\item  $h$ is subcritical for $f$ and $\mathbf{g}$ if $\lambda _{h,f,\mathbf{g}}<\frac{1}{K(n,2){{}^{2}}(\stackunder{M}{Sup}f)^{\frac{n-2}{n}}}$

\item  $h$ is \textbf{critical} for $f$ and $\mathbf{g}$ if $h$ 
is weakly critical and if for any function 
$k\leq h,\,k\neq h$ such that $\bigtriangleup _{\mathbf{g}}+k$ is coercive, $k$ is subcritical.
\end{itemize}
\end{defi}

Using the theorem of Th.Aubin, we can give an equivalent definition of critical functions. Indeed, using this 
theorem, it is easy to see that if $h$ is weakly critical and $\E$ has a minimizing solution $u$, then $h$ is a 
critical function; just note that for $k\leq h,\,k\neq h$, $I_{k,\g}(u)<I_{h,\g}(u)$. Therefore, we can give the following equivalent definition:
\begin{defi}
A function $h$ is critical for $f$ and $\g$ if:
\begin{itemize}
\item for any continuous function $k\leq h,\,k\neq h$ such that $\bigtriangleup _{\mathbf{g}}+k$ is coercive, 
(which is the case as soon as $k$ is close enough to $h$ in $C^{0}$), $(E_{k,f,\g})$ has a minimizing solution,
\item for any continuous function $k'\geq h,\,k'\neq h$, $(E_{k',f,\g})$ has \textbf{no} minimizing solution.
\end{itemize}
\end{defi}
Remark: if $h$ is weakly critical for a positive function $f$, necessarily, $\bigtriangleup _{\mathbf{g}}+h$ is coercive; 
just use the Sobolev inequality.

Critical functions are thus introduced as "separating" functions giving rise to an equation having 
minimizing solutions, and functions giving rise to an equation that cannot have any such solution. We therefore 
have transformed the problem of the existence of minimizing solutions when $\lambda _{h,f,\mathbf{g}}=\frac{1}{K(n,2){{}^{2}}(%
\stackunder{M}{Sup}f)^{\frac{n-2}{n}}}$ to the problem of existence of minimizing solutions to $\E$ 
when $h$ is a critical function.

Before passing to the theorems proved in this work, we have to give two very important properties of critical functions.

First, they transform in conformal changes of metric exactly like scalar curvature: indeed, let $u\in C^{\infty }(M),\,u>0$ and
 $\mathbf{g}^{\prime }=u^{\frac{4}{n-2}}\mathbf{g}$ a metric conformal to $\mathbf{g}$. Let also $h$ be 
 a smooth function. We set 
 $$h^{\prime }=\frac{\triangle _{\mathbf{g}}u+h.u}{u^{\frac{n+2}{n-2}}}.$$
 Then, some computations show that $h$ is critical for $f$ and $\g$ iff $h'$ is critical for $f$ and $\g'$.
 
 Second, we come back to the evaluation of $\lambda _{h,f,\mathbf{g}}$. Th. Aubin introduced, in the functional 
 $J_{h,f,\g}$ the following test functions:
 $$\psi _{k}(Q)=\left\lbrace 
\begin{array}{c}
(\frac{1}{k}+r{{}^{2}})^{-\frac{n-2}{2}}-(\frac{1}{k}+\delta {{}^{2}})^{-%
\frac{n-2}{2}}\,\,\,\,if\,r<\delta \\ 
0\,\,\,\,\,\,\,\,if\,\,\,r\geq \delta
\end{array} \right. $$
where: $\delta <injM$ (the injectivity radius of $M$), $P\,\in M$ is a fixed point, 
$k\in \Bbb{N}^{*}$, and where $r=d_{\mathbf{g}}(P,Q).$ When $dimM=n\geqslant 4$, we get, if $P$ is a point 
where $f$ is maximum on $M$:
 \begin{center}
$J_{h,f,\mathbf{g}}(\psi _{k})=\frac{1}{K(n,2){{}^{2}}(\stackunder{M}{Sup}%
f)^{\frac{n-2}{n}}}\left\{ 1+\frac{1}{n(n-4)}\left( \frac{4(n-1)}{n-2}%
h(P)-S_{\mathbf{g}}(P)+\frac{n-4}{2}\frac{\bigtriangleup _{\mathbf{g}}f(P)}{%
f(P)}\right) \frac{1}{k}\right\} +o(\frac{1}{k})$
\end{center}
 We therefore get the following important proposition:
 \begin{proposition}
 If $dimM\geq4$ and if $h$ is weakly critical for $f$ and $\g$ (thus in particular if it is critical), as $\lambda _{h,f,\mathbf{g}}=\frac{1}{K(n,2){{}^{2}}(\stackunder{M}{Sup}f)^{%
\frac{n-2}{n}}}$, necessarily, if $P$ is a point of maximum of $f$:
\[
\frac{4(n-1)}{n-2}h(P)\geq S_{\mathbf{g}}(P)-\frac{n-4}{2}\frac{\bigtriangleup _{\mathbf{g}}f(P)}{f(P)} 
\]
 \end{proposition}
 Remark: if $f$ is constant on $M$, this means that $\frac{4(n-1)}{n-2}h\geq S_{\mathbf{g}}$ on all of $M$. Note 
 also that in dimension 4, the term $\frac{\bigtriangleup _{\mathbf{g}}f(P)}{f(P)}$ disappears.
 
\section{Statement of the results}
In all what follows, we will make the following hypothesis:

\textbf{Hypothesis (H):} We now suppose that $dimM=n\geqslant 4$. We suppose that all our functions $h$ are such that 
$\bigtriangleup _{\mathbf{g}}+h$ is coercive. Also, $f$ will always be a smooth function such that $\underset{M}{Sup} f>0$. 
We will denote $Max\,f=\{x\in M/ f(x)=\underset{M}{Sup} f \}$.

Our first theorem concerns the existence of minimizing solutions to $\E$ when $h$ is critical.
\begin{theor}
If $h$ is a critical function for $f$ and $\g$, ($h,f,\g$ verifying \textbf{H}), and if for all point $P$ where $f$ is 
maximum on $M$, we have $$\frac{4(n-1)}{n-2}h(P)>S_{\mathbf{g}}(P)-\frac{n-4%
}{2}\frac{\bigtriangleup _{\mathbf{g}}f(P)}{f(P)},$$ then there exist a minimizing solution for $\E$.
\end{theor}
This theorem is an immediate consequence of the following result, more general but more technical in its statement. (Just 
take $h_{t}=h-t$ to get the theorem above.)
\begin{theo}
Let $h$ be a weakly critical function for $f$ and $\g$, (assuming hypothesis \textbf{H}). If, for all point $P$ where $f$ 
is maximum, we have
 $$\frac{4(n-1)}{n-2}h(P)>S_{\mathbf{g}}(P)-\frac{n-4}{2}\frac{\bigtriangleup _{\mathbf{g}}f(P)}{f(P)},$$
and if there exists a family of functions $(h_{t})$, $h_{t}\lvertneqq  h$, $h_{t}$ being sub-critical for all 
$t$ in a neighbourhood of a real $t_{0}\in \Bbb{R}$, and such that $h_{t}\stackunder{t\rightarrow t_{0}}{\rightarrow }h$
 in $C^{0,\alpha }$, then there exists a minimizing solution for $\E$, and therefore, $h$ is critical for $f$ and $\g$.
\end{theo}
E. Hebey and M. Vaugon, in the context of their study of $B_{0}(\g)$, proved this theorem in the case where $f$ is constant, 
and as them, we base our computations on the article of Djadli and Druet [9]. 
The presence of a non-constant function $f$ on the right of equation $\E$ introduces new difficulties in the proof, 
and requires the use of very powerfull estimates concerning concentration phenomena's, called $C^0-theory$, 
due to Druet and Robert [13]; 
the use of $C^0-theory$ was kindly suggested to us by E. Hebey. 
Also, an alternate proof, not using $C^0-theory$, thus in some sense more elementary, 
but requiring the additional hypothesis that the hessian of $f$ is non-degenerate at its points of maximum on $M$, 
will, as a "byproduct", prove another very important estimate concerning these concentration phenomena's, not 
available without heavy hypothesis in the case when $f$ is a constant function; this estimate concerns the speed of 
convergence to a concentration point, (see subsection 4.2), is of independent interest, and was obtained in the
author's PHD thesis to prove theorem 1.

The next natural question is of course to know if there exist critical functions. The answer, positive, will appear to be a 
consequence of theorem 1. We will say that a set $E\subset M$ is \textit{thin} if $M-E$ contains a dense open subset.
\begin{theo}
Being given the manifold $(M,\g)$ and a non constant function $f$, there exist infinitely many 
functions $h$ critical for $f$ and $\g$, which satisfy, in each point $P$ of maximum of $f$,
 $$\frac{4(n-1)}{n-2}h(P)>S_{\mathbf{g}}(P)-\frac{n-4}{2}\frac{\bigtriangleup _{\mathbf{g}}f(P)}{f(P)}\,\,\,\,(*)$$
By theorem 1, these critical functions are such that $\E$ have minimizing solutions. 
Also, if the set of maximum points of $f$ is thin and if $\int_{M}f>0$, there exist strictly positive such critical 
functions $h$, i.e. satisfying (*).
\end{theo}

These first theorems lead us to modify slightly our vision of critical functions. Note that in equation $\E$, there are 
three datas that one can modify: the functions $h$ and $f$, of course, but also the metric $\g$ in a conformal 
class, as, by the conformal laplacian transformation formula, the equation is changed in a similar one if we change 
$\g$ in $\g'=u^{\frac{4}{n-2}}.\g$. This lead us to the following definition:
\begin{center}
\textit{$(h,f,\g)$ is a critical triple if $h$ is a critical function for $f$ and $\g$.}
\end{center}
We shall say that the triple $(h,f,\g)$ has minimizing solutions if $\E$ has; we can also speak of weakly critical 
or sub-critical triples.We then asked ourselves the following question:

\textit{Being given two of the three datas of a triple, can one find the third to obtain a critical triple?}

For example, the problem of the existence of  critical functions can be formulated in the following manner: we 
are given the function $f$ and the metric $\g$, can we complete the triple $(.,f,\g)$ by a function $h$ to obtain a 
critical triple $(h,f,\g)$?

We adress the two other questions, first fixing $h$ and $f$ and seeking a conformal metric $\g'$, and then fixing 
the function $h$ and the metric $\g$ and seeking a function $f$. We obtain answers expressed by the following two 
theorems:
\begin{theo}
On the manifold $(M,\g)$, let be given a function $h$ and a function $f$, satisfying (\textbf{H}). 
We suppose that the set of maximum points of $f$ is thin.
Then, there exist a metric $\g'$ conformal to $\g$ such that $(h,f,\g')$ is a critical triple. 
Moreover, we can find $\g'$ such that $(h,f,\g')$ has minimizing solutions.
\end{theo}

This theorem was proved by E. Humbert and M. Vaugon in the case $f=cst=1$ and $M$ not conformally diffeomorphic 
to the sphere,[21]. Their method works in the case 
of a non constant function $f$ and an arbitrary manifold once it is proved that we can suppose the existence of positive critical functions 
satisfying the strict inequality (*) in theorem 2, result we included in this theorem (note that, as $Sup f>0$, we can 
always find a metric $\textbf{g'}$ conformal to $\g$ such that $\int_{M}fdv_{\g'}>0$). 
 In fact, when $M$ is not conformally diffeomorphic to the sphere and $S_{\g}$ is constant, it can be 
proved that $B_{0}(\g)K(n,2)^{-2}$ is a critical (constant) function for $1$ and $\g$, and it is obviously 
positive. We will discuss weaker hypothesis for this theorem, as well as the problem of existence of positive 
critical functions in section 6.

The last question brings us to the following answer when the dimension of $M$ is greater than 5, requirement which is 
linked to the fact that $\frac{\bigtriangleup _{\mathbf{g}}f(P)}{f(P)}$ dissapears in dimension 4 in the 
inequality of Proposition 1.
\begin{theo}
Let be given the manifold $(M,\g)$ of dimension $n\geq5$, and a function $h$ such that $\bigtriangleup _{\mathbf{g}}+h$ 
is coercive. Then, there exists a non constant function $f$ such that $(h,f,\g)$ is critical with minimizing solutions if, and 
only if, $(h,1,\g)$ is a sub-critical triple (where 1 is the constant function 1).
\end{theo}

Note that if $(h,1,\g)$ is weakly critical, then either this triple has minimizing solutions in which case it is a critical 
triple, or there is no non-constant function $f$ such that $(h,f,\g)$ is critical with minimizing solutions (see the proof 
and what follows). 
The proof of this theorem is quite difficult, and make use of the method developped for the proof of theorem 1. 
Also, this proof brought us to make some more remarks about critical functions. First, it is easily seen, by using the 
functional $J$, that if $(h,f,\g)$ is a critical triple, then, for any constant $c>0$, $(h,c.f,\g)$ is also a critical triple. 
It would therefore be more appropriate to speak of triple $(h,[f],\g)$ where $[f]=\{c.f\,/c>0\}$ could be called 
the "class" of $f$. Note for example that we can always suppose that $Sup f=1$; also, to compare two triples $(h,f,\g)$ 
and $(h,f',\g)$, one has to suppose that $Sup f=Sup f'$. Note also that on $[f]$, the quotient  $\frac{\bigtriangleup _{\mathbf{g}}f}{f}$ 
is constant. Second, in the proof of theorem 4, we had to approximate the function $f$ by a family $(f_{t})$, unlike 
theorem 1 where we used a family $(h_{t})$ approaching $h$. This suggested another possible definition of critical 
functions, dual to the first one in the sense that we exchange the role of $h$ and $f$.
\begin{defi}
Let $(M,\g)$ be of dimension $n\geq3$ and $h$ be such that $\bigtriangleup _{\mathbf{g}}+h$ is coercive. 
We shall say that a smooth function $f$ such that $\underset{M}{Sup} f>0$ is critical for $h$ and $\g$ if:
\begin{itemize}
\item  a/: $\lambda _{h,f,\mathbf{g}}=\frac{1}{K(n,2){{}^{2}}(\stackunder{M}{Sup}f)^{\frac{n-2}{n}}}$

\item  b/: for any smooth function $f^{\prime }$ such that $Supf=Supf^{\prime }$ and $f^{\prime }\gneqq f$, \\

$\lambda_{h,f^{\prime },\mathbf{g}}<\frac{1}{K(n,2){{}^{2}}(\stackunder{M}{Sup}%
f^{\prime })^{\frac{n-2}{n}}}$

\item  Remark: if $Supf=Supf^{\prime }$ and $f^{\prime
}\lvertneqq f$, then $\lambda _{h,f^{\prime },\mathbf{g}}=\frac{1}{%
K(n,2){{}^{2}}(\stackunder{M}{Sup}f^{\prime })^{\frac{n-2}{n}}}$ as $J_{h,f^{\prime },\mathbf{g}}(w)\geqslant J_{h,f,\mathbf{g}}(w)$%
for any function $w$.
\end{itemize}
\end{defi}
It is then natural to ask if the two definitions are equivalent ($\g$ being fixed):
\begin{center}
\textit{Is $f$ critical for $h$ if, and only if, $h$ is critical for $f$ ?}
\end{center} 

This question seems quite difficult. A positive answer would justify the concept of critical triple. Remember that, 
because of proposition 1, we have in both cases, when $P$ is a point where $f$ is maximum on $M$: 
$$\frac{4(n-1)}{n-2}h(P)\geqslant S_{\mathbf{g}}(P)-\frac{n-4}{2}\frac{\bigtriangleup _{\mathbf{g}}f(P)}{f(P)}.$$ 

We obtain the following theorem:
\begin{theo}
Let $(M,\g)$ be a compact manifold of dimension $n\geq5$, and let $h$ be a function such that  $\bigtriangleup _{\mathbf{g}}+h$ is coercive. 
Let $f$ be a smooth function such that $\underset{M}{Sup} f>0$. We suppose 
that for any point $P$ where $f$ is maximum on $M$:
$$\frac{4(n-1)}{n-2}h(P)> S_{\mathbf{g}}(P)-\frac{n-4}{2}\frac{\bigtriangleup _{\mathbf{g}}f(P)}{f(P)}.$$
Then, $f$ is critical for $h$ if, and only if, $h$ is critical for $f$.
\end{theo}
Remark: if $1$ is critical for $h$, then every non constant function $f$, such that $Sup f=1$, is weakly critical for $h$ 
with \textit{no} minimizing solutions. Indeed, here again if a function $f$ is weakly critical for $h$ with a minimizing 
solution, then $f$ is critical.

There is an interesting consequence of theorems 4 and 5. We said in the introduction that an important application 
of equations $\E$ was the study of prescribed scalar curvature: being given a smooth function $f$ on the manifold 
$(M,\g)$, is $f$ the scalar curvature of a metric conformal to $\g$? The theorem of Th. Aubin shows that if $f$ is 
sub-critical for $S_{\g}$, then $f$ is a scalar curvature. Theorem 4 applied to $h=\frac{n-2}{4(n-1)}S_{\g}$ 
shows that:

\textit{On a compact manifold $(M,\g)$ not conformaly diffeomorphic to the sphere, there exist scalar curvatures 
of metric conformal to $\g$ that are only weakly critical, (more precisely critical).}
\\

Another application, remarked by E. Hebey, is the study of \textit{Sobolev inequality in the presence of a twist.} 

The previous theorems all deal with manifolds of dimension at least 4, or even 5. We will give results concerning 
the dimension 3 in the last section. They are very interesting, but they are rapid generalisations of results obtained 
by O. Druet in the case $f=constant$ [10], the introduction of a non constant $f$ introducing this time no real difficulties. 

\section{The three main tools}
We want to present here the three main tools used in the proof of our various theorems. These tools were developed by 
several persons since M. Vaugon and P.L. Lions, essentially E. Hebey, O. Druet F. Robert, 
M. Struwe, E. Humbert and Z. Faget, among others. 
\subsection{The concentration point.}
To prove the existence of a solution $u>0$ to our equation 
\[
(E_{h,f,\mathbf{g}}):\,\triangle _{\mathbf{g}}u+h.u=\lambda .f.u^{\frac{n+2}{%
n-2}}, 
\]
the idea will often be to associate a family of equations having minimizing solutions $u_{t}>0$ : 
\[
E_{t}:\,\,\triangle _{\mathbf{g}}u_{t}+h_{t}.u_{t}=\lambda _{t}.f.u_{t}^{%
\frac{n+2}{n-2}} 
\]
with
\[
h_{t}\rightarrow h\,\,\,\,in \,\,\,\,C^{0,\alpha }(M) 
\]
and $\lambda _{t}\rightarrow \lambda $ a converging sequence of real numbers, in such a way that for some $u\in H_{1}^{2}$ : 
$u_{t}\rightarrow u\,$ strongly in $L^{p}$ , 
$p<2^{*}$, and $u_{t}\rightharpoondown u$ weakly in  $H_{1}^{2}$ with a constraint 
\[
\int_{M}f.u_{t}^{2^{*}}dv_{\mathbf{g}}=1. 
\]
To simplify, we will suppose that all convergences are for $t\rightarrow t_{0}=1$. The difficulty will be to prove 
that $u$ is not the trivial zero solution, as then, by the maximum principle, we have $u>0$. 
We will proceed by contradiction, and suppose $u\equiv 0$. The idea is then that, because of the condition $\int_{M}f.u_{t}^{2^{*}}=1,$ 
all the "mass" of the functions $u_{t}$, which converge to  0 in $L^{p}$ , $p<2^{*}$, concentrates around a point 
of the manifold. We thus define:
\begin{defi}
\textit{\ }$x_{0}\in M$ is a point of concentration of the sequence $(u_{t})$ if for any $\delta >0$ :
\[
\stackunder{t\rightarrow t_{0}}{\lim \sup }\int_{B(x_{0},\delta
)}u_{t}^{2^{*}}dv_{\mathbf{g}}>0 
\]
\end{defi}
It is easy to see that because $M$ is compact and we require  $\int_{M}f.u_{t}^{2^{*}}dv_{\mathbf{g}}=1$, there exist 
at least one point of concentration. We will show that there exists only one point of concentration, that it is a point where 
$f$ is maximum, and that there exist a sequence of points $x_{t}$ converging to a point $x_{0}\in M$ such that
\[
u_{t}(x_{t})=\stackunder{M}{\max }u_{t}\rightarrow +\infty , 
\]
and 
\[
u_{t}\rightarrow 0\,\,\,in\,\,\,C_{loc}^{0}(M-\{x_{0}\}). 
\]
In fact the idea is that one can do "as if" the functions $u_{t}$ have compact support in a small neigbourhood of $x_{0}$ 
when $t$ is close to $t_{0}$.

\subsection{Blow-up analysis}
Thanks to the concentration point, one brings back the study of the family $u_{t}$ converging to 0, to what happens around 
$x_{0}$. The idea of \textit{blow-up analysis} is to do a "change of scale" around $x_{0}$: we will call \textit{blow-up} 
of center $x_{t}$ and coefficient $k_{t}$ the following sequence of charts and changes of metrics. We consider, for $\delta$ 
small enough:
\[
\begin{array}{ccccc}
B(x_{t},\delta ) & \stackrel{\exp _{x_{t}}^{-1}}{\rightarrow } & \stackunder{%
}{B(0,\delta )\subset \Bbb{R}^{n}} & \stackrel{\psi _{k_{t}}}{\stackunder{}{%
\rightarrow }} & B(0,k_{t}\delta )\subset \Bbb{R}^{n} \\ 
&  & x\,\,\, & \mapsto \,\, & k_{t}x \\ 
\,\mathbf{g}\, & \rightarrow & \,\,\mathbf{g}\,_{t}=\exp _{x_{t}}^{*}\mathbf{%
g}\, & \rightarrow & \,\widetilde{\,\mathbf{g}\,}_{t}=k_{t}^{2}(\psi
_{k_{t}}^{-1})^{*}\mathbf{g}_{t}
\end{array}
\]
where $\exp _{x_{t}}^{-1}$ is the chart deduced from the exponential map in $x_{t}$. We set 
\[
\overline{u}_{t}=u_{t}\circ \exp _{x_{t}}\, ; \,
\overline{f_{t}}=f\circ \exp _{x_{t}} \, ;
\,\overline{h_{t}}=h_{t}\circ \exp _{x_{t}} 
\]
We have
\begin{eqnarray*}
\triangle _{\mathbf{g}_{t}}\overline{u}_{t}+\overline{h}_{t}.\overline{u}%
_{t} &=&\lambda _{t}\overline{f_{t}}.\overline{u}_{t}^{\frac{n+2}{n-2}} \\
\int_{B(0,r)}\overline{u}_{t}^{\alpha }dv_{\,\mathbf{g}\,_{t}}
&=&\int_{B(x_{t},r)}u_{t}^{\alpha }dv_{\mathbf{g}}\text{ for all }\alpha
\geq 1
\end{eqnarray*}
We then set
$$
m_{t} =\stackunder{M}{Max}\,u_{t}\,  ; \,
\widetilde{u}_{t} =m_{t}^{-1}\overline{u}_{t}\circ \psi _{k_{t}}^{-1}\, ;
\,\widetilde{h}_{t} =\overline{h}_{t}\circ \psi _{k_{t}}^{-1}\, ;
\,\widetilde{f}_{t} =\overline{f}_{t}\circ \psi _{k_{t}}^{-1}\, ;
\,\widetilde{\,\mathbf{g}\,}_{t} =k_{t}^{2}(\exp _{x_{t}}\circ \psi_{k_{t}}^{-1})^{*}\mathbf{g},
$$
so in particular
$
\widetilde{u}_{t}(x) =m_{t}^{-1}\overline{u}_{t}(\frac{x}{k_{t}})\, and
\,\widetilde{\,\mathbf{g}\,}_{t}(x) =\exp _{x_{t}}^{*}\mathbf{g}(\frac{x}{k_{t}}).
$
Then: 
\begin{eqnarray}
(\widetilde{E}_{t})\, &:&\,\triangle _{\widetilde{\mathbf{g}}_{t}}\widetilde{%
u}_{t}+\frac{1}{k_{t}^{2}}\widetilde{h}_{t}.\widetilde{u}_{t}=\frac{m_{t}^{%
\frac{4}{n-2}}}{k_{t}^{2}}\lambda _{t}\widetilde{f}_{t}.\widetilde{u}_{t}^{%
\frac{n+2}{n-2}}   \\
and\, &:&\,\int_{B(0,k_{t}r)}\widetilde{u}_{t}^{\alpha }dv_{\widetilde{\,%
\mathbf{g}\,}_{t}}=\frac{k_{t}^{n}}{m_{t}^{\alpha }}\int_{B(x_{t},r)}u_{t}^{%
\alpha }dv_{\mathbf{g}}  \nonumber
\end{eqnarray}
We will mostly use the following parameters : we consider a sequence of points ($x_{t})$ such that: 
\[
m_{t}=\stackunder{M}{Max}\,u_{t}=u_{t}(x_{t}):=\mu _{t}^{-\frac{n-2}{2}}
\]
and
\[
\,k_{t}=\mu _{t}^{-1}.
\]
$\mu _{t}$ will appear to be a fundamental parameter in the study of concentration phenomena's. Noting $(x^{i})$ the coordinates in $\Bbb{%
R}^{n}$, one has : 
\begin{eqnarray}
(\widetilde{E}_{t})\, &:&\,\triangle _{\widetilde{\mathbf{g}}_{t}}\widetilde{%
u}_{t}+\mu _{t}^{2}.\widetilde{h}_{t}.\widetilde{u}_{t}=\lambda _{t}%
\widetilde{f}_{t}.\widetilde{u}_{t}^{\frac{n+2}{n-2}}   \\
and\, &:&\,\int_{B(0,\mu _{t}^{-1}r)}x^{i_{1}}...x^{i_{p}}.\widetilde{u}%
_{t}^{\alpha }dv_{\widetilde{\,\mathbf{g}\,}_{t}}=\mu _{t}^{-p-n+\alpha 
\frac{n-2}{2}}\int_{B(0,r)}x^{i_{1}}...x^{i_{p}}\overline{u}_{t}^{\alpha
}dv_{\,\mathbf{g}\,_{t}}  \nonumber
\end{eqnarray}
A very important result is that when $\mu _{t}\rightarrow 0$ and therefore $k_{t}\rightarrow
+\infty $, the components of $\,\widetilde{\mathbf{g}}_{t}$ converge in 
$C_{loc}^{2}$ to those of the euclidean metric, and $(\widetilde{E}_{t})$``converges'' to the equation: 
\[
\triangle _{e}\widetilde{u}=\lambda f(x_{0}).\widetilde{u}^{\frac{n+2}{n-2}}
\]
in the sense that
\[
\widetilde{u}_{t}\rightarrow \widetilde{u}\,\,\,\,in \,\,\,\,C_{loc}^{2}(\Bbb{R}^{n}).
\]
It is known, then, that
\[
\widetilde{u}=(1+\frac{\lambda f(x_{0})}{n(n-2)}\left| x\right| ^{2})^{-%
\frac{n-2}{2}}.
\]

\subsection{The iteration process}
The idea of the Möser iteration process is to multiply the equations $(E_{t})$ by succesive powers $u_{t}^k$ of the functions 
$u_{t}$ and to integrate over $M$ to obtain bounds on increasing $L^p$-norms of the $u_{t}$. To localize the study around 
the concentration point $x_{0}$, which is a maximum point for $f$, we shall in fact multiply the equations by $\eta {{}^{2}}u_{t}^{k}$ where $\eta $ is a 
cut-off function equal to 1 (resp.0) on a ball $B(x_{0},r)$ where $f\geq 0$, and equal to 0
(resp. 1) on $M\backslash B(x_{0},2r)$, and where $k\geq 1$, then integrate by part. We will therefore be able to study 
blow-up around $x_{0}$ using this method. We get after some integrations by parts, and using equation $(E_{t})$ :
\begin{equation}
\frac{4k}{(k+1)^{2}}\int_{M}\left| \nabla (\eta u_{t}^{\frac{k+1}{2}%
})\right| ^{2}=\lambda _{t}\int_{M}f\eta ^{2}u_{t}^{\frac{n+2}{n-2}%
}u_{t}^{k}+\int_{M}(\frac{2}{k+1}\left| \nabla \eta \right| ^{2}+\frac{2(k-1)%
}{(k+1)^{2}}\eta \triangle \eta -\eta ^{2}h_{t})u_{t}^{k+1}
\end{equation}
where the integrals are taken with the measure $dv_{\mathbf{g}}$. Then using Hölder inequality, 
if $f\geq 0$ on $Supp\,\eta $ we obtain: 
\[
\lambda _{t}\int_{M}f\eta ^{2}u_{t}^{\frac{n+2}{n-2}}u_{t}^{k}\leq \lambda
_{t}(\stackunder{Supp\,\eta }{Sup}f)^{\frac{n-2}{n}}.(\int_{Supp\,\eta
}fu_{t}^{\frac{2n}{n-2}})^{\frac{2}{n}}.(\int_{M}(\eta u_{t}^{\frac{k+1}{2}%
})^{\frac{2n}{n-2}})^{\frac{n-2}{n}} 
\]
Then using Sobolev inequality : 
\[
(\int_{M}(\eta u_{t}^{\frac{k+1}{2}})^{\frac{2n}{n-2}})^{\frac{n-2}{n}}\leq
K(n,2){{}^{2}}\int_{M}\left| \nabla (\eta u_{t}^{\frac{k+1}{2}})\right|
^{2}+B\int_{M}\eta u_{t}^{k+1} 
\]
with $B>0$. Therefore:

\smallskip

\begin{eqnarray*}
\frac{4k}{(k+1)^{2}}(\int_{M}(\eta u_{t}^{\frac{k+1}{2}})^{\frac{2n}{n-2}})^{%
\frac{n-2}{n}} \leq &\lambda _{t}K(n,2){{}^{2}}(\stackunder{Supp\,\eta }{Sup%
}f)^{\frac{n-2}{n}}.(\int_{Supp\,\eta }fu_{t}^{\frac{2n}{n-2}})^{\frac{2}{n}%
}.(\int_{M}(\eta u_{t}^{\frac{k+1}{2}})^{\frac{2n}{n-2}})^{\frac{n-2}{n}} \\
&+\int_{M}(\frac{4k}{(k+1){{}^{2}}}B\eta +\frac{2}{k+1}\left| \nabla \eta
\right| ^{2}+\frac{2(k-1)}{(k+1)^{2}}\eta \triangle \eta -\eta
^{2}h_{t})u_{t}^{k+1}
\end{eqnarray*}
Then:
\begin{equation}
Q(t,k,\eta ).(\int_{M}(\eta u_{t}^{\frac{k+1}{2}})^{\frac{2n}{n-2}})^{\frac{%
n-2}{n}}\leq (\frac{4k}{(k+1){{}^{2}}}B+C_{0}+C_{\eta })\int_{Supp\,\eta
}u_{t}^{k+1}\,\,\,\,\,\,\,\,\,\,\,\,\,\,\,\,\,\,\,\,\,\,\,\,\,\,\,\,\,\,\,\,%
\,\,\,\,\,\,\,\,\,\,\,\,\,\,\,\,\,\,\,\,\,\,\,  
\end{equation}
where
\[
Q(t,k,\eta )=\frac{4k}{(k+1){{}^{2}}}-\lambda _{t}K(n,2){{}^{2}}(\stackunder{%
Supp\,\eta }{Sup}f)^{\frac{n-2}{n}}.(\int_{Supp\,\eta }f.u_{t}^{2^{*}})^{%
\frac{2}{n}} 
\]
where we remind that $2^{*}=\frac{2n}{n-2}$ and where $C_{0}\,et\,C_{\eta }$ are constants independant of
 $k$ and $t$ and such that $\forall k\geq 1,\forall t:$
\[
\,\,\left\| \frac{2}{k+1}\left| \nabla \eta \right| ^{2}+\frac{2(k-1)}{%
(k+1)^{2}}\eta \triangle \eta \right\| _{L^{\infty }(M)}\leq \,C_{\eta
}\,\,and\,\,\left\| h_{t}\right\| _{L^{\infty }(M)}\leq C_{0}\,. 
\]
If the sign of $f$ changes on $Supp\,\eta $, we go back to Hölder's inequality:
\[
\lambda _{t}\int_{M}f\eta ^{2}u_{t}^{\frac{n+2}{n-2}}u_{t}^{k}\leq \lambda
_{t}(\stackunder{Supp\,\eta }{Sup}\left| f\right| ).(\int_{Supp\,\eta
}u_{t}^{\frac{2n}{n-2}})^{\frac{2}{n}}.(\int_{M}(\eta u_{t}^{\frac{k+1}{2}%
})^{\frac{2n}{n-2}})^{\frac{n-2}{n}} 
\]
to obtain (6) with: 
\begin{equation}
Q(t,k,\eta )=\frac{4k}{(k+1){{}^{2}}}-\lambda _{t}K(n,2){{}^{2}}(\stackunder{%
Supp\,\eta }{Sup}\left| f\right| ).(\int_{Supp\,\eta }u_{t}^{2^{*}})^{\frac{2%
}{n}}  
\end{equation}
One can also replace $\stackunder{Supp\,\eta }{Sup}%
\left| f\right| $ by $\stackunder{M}{Sup}f$.

The goal is to show that $(\eta u_{t})$ is bounded in $L^{\frac{k+1}{2}2^{*}}$ and therefore that we can 
extract a sub-sequence converging strongly in $L^{2^{*}}$.

\textbf{Remark}
Those three tools also work for more general equations that we can associate to $(E_{h,f,\mathbf{g}}):\,\triangle _{\mathbf{g}}u+h.u=\mu
_{h}.f.u^{\frac{n+2}{n-2}}$. like e.g.
$
E_{t}:\,\,\triangle _{\mathbf{g}}u_{t}+h_{t}.u_{t}=\lambda
_{t}.f_{t}.u_{t}^{q_{t}-1} 
$
where $q_{t}\rightarrow 2^{*}$ and $f_{t}\rightarrow f$ in some $L^{p}$, still with $h_{t}\rightarrow h$\thinspace in $%
\,C^{0,\alpha }(M)$ and $\lambda _{t}\rightarrow \lambda $.

\section{Proof of theorem 1}
\subsection{Setup}
Let $h$ be a weakly critical function for $f$ and $\mathbf{g}$
such that for any $P\in M$ where $f$ is maximum on $M$ we have : 
\[
\text{ }h(P)>\frac{n-2}{4(n-1)}S_{\mathbf{g}}(P)-\frac{(n-2)(n-4)}{8(n-1)}%
\frac{\bigtriangleup _{\mathbf{g}}f(P)}{f(P)} 
\]
and such that there exist a family $(h_{t}),\,h_{t}\lvertneqq
h,\,h_{t}$ sub-critical for every $t$, and satisfying $h_{t}\stackunder{%
t\rightarrow t_{0}}{\rightarrow }h\,\,$in$\,$ $C^{0,\alpha}$. To simplify, we suppose that $t_{0}=1$ and that 
$t\rightarrow 1$. Then for every $t$ :
\[
\lambda _{t}:=\lambda _{h_{t},f,\mathbf{g}}<\frac{1}{K(n,2){{}^{2}}(%
\stackunder{M}{Sup}f)^{\frac{n-2}{n}}} 
\]
and there exist a family $u_{t}$ of minimizing solutions of the equations
\[
E_{t}:\,\,\triangle _{\mathbf{g}}u_{t}+h_{t}.u_{t}=\lambda _{t}.f.u_{t}^{%
\frac{n+2}{n-2}}\text{ with }\int_{M}fu_{t}^{2^{*}}dv_{\mathbf{g}}=1 
\]
We then see, as $\triangle _{\mathbf{g}}+h$ is coercive, that the sequence $(u_{t})$ is bounded in 
$H_{1}^{2}$ (just multiply $E_{t}$ by $u_{t}$ and integrate on $M$). Thus, there exist a function $u\in H_{1}^{2}\,,\,u\geq 0$ 
such that, after extracting a subsequence, 
\begin{eqnarray*}
&&u_{t}\stackrel{H_{1}^{2}}{\rightharpoondown }u\,,\, \\
&&u_{t}\stackrel{L^{2}}{\rightarrow }u\,, \\
&&\,u_{t}\stackrel{p.p.}{\rightarrow }u\,,
\end{eqnarray*}
and we can suppose
\[
\lambda _{t}\stackrel{<}{\rightarrow }\lambda \leqslant \frac{1}{K(n,2){%
{}^{2}}(\stackunder{M}{Sup}f)^{\frac{n-2}{n}}}\,\,. 
\]
In particular
\[
u_{t}\stackrel{L^{p}}{\rightarrow }u,\,\forall p<2^{*}=\frac{2n}{n-2} 
\]
as the inclusion of $H_{1}^{2}$ in $L^{p}$ is compact $\forall
p<2^{*} $. Therefore $u$ is a weak solution of
\[
\,\triangle _{\mathbf{g}}u+h.u=\lambda .f.u^{\frac{n+2}{n-2}} 
\]
and by standard elliptic theory, $u$ is $C^{\infty }$. The maximum principle then gives us that either $u>0$ or 
$u\equiv 0$.

If $u>0$ then, using elliptic theory and iteration process, and the fact that $h$ is weakly critical, one can prove that:
\[
\lambda =\frac{1}{K(n,2){{}^{2}}(\stackunder{M}{Sup}f)^{\frac{n-2}{n}}} 
\]
and then that $u$ is a minimizing positive solution of
\[
\triangle _{\mathbf{g}}u+h.u=\frac{1}{K(n,2){{}^{2}}(\stackunder{M}{Sup}f)^{%
\frac{n-2}{n}}}.f.u^{\frac{n+2}{n-2}}\text{ with }\int_{M}fu^{2^{*}}dv_{%
\mathbf{g}}=1 
\]
and the theorem is proved.

If $u\equiv 0$, we will show that there is a concentration phenomena. All the study that follows will aim at finding a 
contradiction. From now, we suppose that we are in this case:
\[
u\equiv 0. 
\]

\subsection{Concentration phenomena}
In this section we study the behavior of a family of $C^{2,\alpha }$ solutions $(u_{t})$ of
\[
\bigtriangleup _{\mathbf{g}}u_{t}+h_{t}u_{t}=\lambda _{t}fu_{t}^{\frac{n+2}{%
n-2}}\text{ with }\int_{M}fu_{t}^{\frac{2n}{n-2}}dv_{\mathbf{g}}=1 
\]
where $f$ is a smooth function such that $\underset{M}{Sup}f>0$. We also suppose that 
$h_{t}\rightarrow h$ in $C^{0,\alpha }$ where $h$ is such that $\bigtriangleup _{\mathbf{g}}+h$ is coercive.
 The sequence $(u_{t}) $ is bounded in $H_{1}^{2}$, therefore,up to a subsequence, $u_{t}\rightharpoondown u$
weakly in $H_{1}^{2}$, and we supose that $u\equiv 0$; that is $u_{t}\rightarrow 0$ in any $L^{p}$ for $p<2^{*}$. 
We also make the following "minimal energy" hypothesis:
\[
\lambda _{t}\leq \frac{1}{K(n,2)^{2}(\stackunder{M}{Sup}f)^{\frac{n-2}{n}}} 
\]
and we can suppose that $\lambda_{t}\rightarrow \lambda$. All this hypothesis are satisfied by the $u_{t}$ of the preceding section. 
The results of this section are valid for $dimM=3$, exept 
$L^{2}$-concentration, valid for $\dim M\geq 4$. In all this text, $c,C$ are constants independant of $t$ and $\delta $.
\begin{proposition}There exist, after extraction of a subsequence, exactly one concentration point $x_{0}$, and it is 
a point where $f$ is maximum on $M$. Moreover
\[
\forall \delta >0,\,\,\,\overline{\stackunder{t\rightarrow 1}{\lim }}%
\int_{B(x_{0},\delta )}fu_{t}^{2^{*}}dv_{\mathbf{g}}=1 
\]
\end{proposition}
\textit{Proof :}
We apply the iteration process. First, as $M$
is compact, there exist at least one point of concentration. Otherwise, we could cover $M$ by a finite number of balls 
$B(x_{i},\delta )$ such that $\stackunder{t\rightarrow 1}{\lim }\int_{B(x_{i},\delta
)}u_{t}^{2^{*}}=0$, and we would have $\stackunder{t\rightarrow 1}{\lim }\int_{M}u_{t}^{2^{*}}=0,$ which 
would contradict
\[
1=\int_{M}fu_{t}^{2^{*}}dv_{\mathbf{g}}\leq Sup\left| f\right|
\int_{M}u_{t}^{2^{*}}dv_{\mathbf{g}} 
\]
The principle of iteration process is the following: if we find, for a point $x$, a cut-off function $\eta$ equal 
to 1 around $x$ such that $Q(t,k,\eta)\geq Q>0$, we get, using formula (6) or (7), that 
$(\eta u_{t}^{\frac{k+1}{2}})$ is bounded in $L^{2^{*}}$, and therefore we can extract a 
subsequence such that $(\eta u_{t})$ converges strongly to 0 in $L^{2^{*}}$; thus $x$ cannot be a concentration 
point.

Let us prove now that we can do this for a point $x$ such that $f(x)\leq 0$. If $f(x)<0$, we choose $\delta $ small enough 
such that $f<0$
on $B(x,\delta )$ and we choose $\eta $ with support in $B(x,\delta )$.
As $(u_{t})$ is bounded in $H_{1}^{2}$
and thus in $L^{2^{*}}$, we get using formula (5), that for any $k$ such that $1\leq
k\leq 2^{*}-1$: 
\[
\frac{4k}{(k+1)^{2}}\int_{M}\left| \nabla (\eta u_{t}^{\frac{k+1}{2}%
})\right| ^{2}\leq \int_{M}(\frac{2}{k+1}\left| \nabla \eta \right| ^{2}+%
\frac{2(k-1)}{(k+1)^{2}}\eta \triangle \eta -\eta ^{2}h_{t})u_{t}^{k+1}\leq
C_{1} 
\]
where $C_{1}$ is independent of $t$. Therefore for any $k$ such that $1\leq
k\leq 2^{*}-1$ there exist $C_{2}$ independent of $t$ such that: 
\[
\int_{M}\left| \nabla (\eta u_{t}^{\frac{k+1}{2}})\right| ^{2}\leq C_{2} 
\]
Therefore $(\eta u_{t}^{\frac{k+1}{2}})$ is bounded in $H_{1}^{2}$ and, using Sobolev inequality, $(\eta
u_{t}^{\frac{k+1}{2}})$ is bounded in $L^{2^{*}}$ for any $k$ such that $1\leq k\leq 2^{*}-1$.

If $f(x)=0$, by continuity of $f$ and choosing $\delta $ small 
enough, we get in (7) that for any $k$ such that $1\leq k\leq 2^{*}-1$, $%
Q(t,k,\eta)\geq Q>0$. Therefore, as we said, here again $(\eta u_{t}^{\frac{k+1}{2}})$ is bounded in $L^{2^{*}}$, 
and therefore we can extract a 
subsequence such that $(\eta u_{t})$ converges strongly to 0 in $L^{2^{*}}$. 
Thus, when $f(x)\leq0$, $x$ cannot be a concentration point.

Now, let $x$ be a concentration point: $f(x)>0$ as we just saw. 
For $\delta >0$ such that $f\geq 0$ on $B(x,\delta )$, set
\[
\stackunder{t\rightarrow 1}{\lim \sup }\int_{B(x,\delta
)}fu_{t}^{2^{*}}=a_{\delta } 
\]
Then as $\int_{M}fu_{t}^{2^{*}}=1$ and as $\overline{lim}\int_{B}fu_{t}^{2^{*}}=0$ 
if $f\leqslant 0$ on $B$ from what we saw above, necessarilly, $a_{\delta }\leq 1$. 
Suppose that there exist $\delta >0$ such that $a_{\delta }<1$. Because
\[
\lambda _{t}\stackunder{\leqslant }{\rightarrow }\lambda \leqslant \frac{1}{%
K(n,2){{}^{2}}(\stackunder{M}{Sup}f)^{\frac{n-2}{n}}} 
\]
we get 
\[
\overline{\stackunder{t\rightarrow 1}{\lim }}\,\lambda _{t}K(n,2){{}^{2}}(%
\stackunder{M}{Sup}f)^{\frac{n-2}{n}}a_{\delta }<1. 
\]
Beside, $\frac{4k}{(k+1)^{2}}\stackunder{k\stackunder{>}{\rightarrow }1%
}{\rightarrow }1$. Therefore, for $k$ close to 1 such that
\[
\overline{\stackunder{t\rightarrow 1}{\lim }}\,\lambda _{t}K(n,2){{}^{2}}(%
\stackunder{M}{Sup}f)^{\frac{n-2}{n}}a_{\delta }<\frac{4k}{(k+1)^{2}} 
\]
we get, taking $\eta $ with support in $B(x,\delta )$, that in formula (6): 
$Q(t,k,\eta)\geq Q>0$ for all $t$, where $Q$ is independent of $t$. So, as before, $x$ cannot be a concentration point, 
and we have a contradiction.
Thus $a_{\delta }=1,\,\forall \delta >0$. Therefore $x$
is the only concentration point, that we will now denote $x_{0}$.
The same reasonning shows that, necessarilly,
\[
\lambda =\frac{1}{K(n,2){{}^{2}}(\stackunder{M}{Sup}f)^{\frac{n-2}{n}}}\,\,. 
\]
In the same way, if $f(x_{0})\neq \stackunder{M}{Sup}f$, there exist $\delta >0$ such that $\stackunder{B(x_{0},\delta )}{Sup}f<\stackunder{M}{Sup}%
f $. But $\lambda _{t}\leq \frac{1}{K(n,2){{}^{2}}(\stackunder{M}{Sup}f)^{%
\frac{n-2}{n}}}$, so
\[
\overline{\stackunder{t\rightarrow 1}{\lim }}\,\lambda _{t}K(n,2){{}^{2}}(%
\stackunder{B(x_{0},\delta )}{Sup}f)^{\frac{n-2}{n}}(\int_{B(x_{0},\delta
)}fu_{t}^{2^{*}})^{\frac{2^{*}-2}{2^{*}}}<1\,\,. 
\]
Then for $k$ close enough to 1, taking $\eta $ with support in $%
B(x_{{0}},\delta )$, we get in (6): $Q(t,k,n)\geq
Q>0 $ for all $t$; and once again we have a contradiction. Therefore $f(x_{0})=\stackunder{M}{Sup}f>0$.

Note that this is the main particularity introduced by the function $f$ on the right of equation $\E$. It gives a precise 
location for the concentration point.

The next propositions concerning the concentration phenomenom are now quite standard, even though they are mostly 
published in the case $f=constant$ and often with few details. We shall therefore give possible proofs, refering to the  
books \cite{D-H} and \cite{DHR} for more information, the presence of a function $f$ introducing only slight modifications that we will indicate 
when necessary.

\begin{proposition}
$u_{t}\rightarrow 0$ in $C_{loc}^{0}(M-\{x_{0}\})$.
\end{proposition}
\textit{Proof :}
It is a typical aplication of the iteration process in standard elliptic theory.
First step: Let $q>0$ be fixed. We prove that for any $\delta >0$, there exists $C=C(\delta ,q)$ independent of 
$t$ such that for $t$ close enough to 1: 
\begin{equation}
\left\| u_{t}\right\| _{L^{q}(M\backslash B(x_{0},\delta ))}\leq C\left\|
u_{t}\right\| _{L^{2}(M)}\,\,. 
\end{equation}
To apply the iteration process, we build a sequence $\eta _{1},...,\eta
_{m}$ of $m$ cut-off functions such that $\eta _{j}=0$ on $B(x_{0},\delta /2)$
and $\eta _{j}=1$ on $M\backslash B(x_{0},\delta )$ and such that
\[
M\backslash B(x_{0},\delta )\subset ...\subset \{\eta _{j+1}=1\}\subset
Supp\,\eta _{j+1}\subset \{\eta _{j}=1\}\subset ...\subset M\backslash
B(x_{0},\delta /2) 
\]
and where $m$ is chosen such that $2(\frac{2^{*}}{2})^{m}>q$. We set $q_{1}=2$ and 
$q_{j}=(\frac{2^{*}}{2})q_{j-1}$. The iteration process (6), (7), gives that
\[
Q(t,q_{j}-1,\eta _{j}).(\int_{M}(\eta _{j}u_{t}^{\frac{q_{j}}{2}})^{2^{*}})^{%
\frac{n-2}{n}}\leq (\frac{4(q_{j}-1)}{q_{j}{{}^{2}}}B+C_{0}+C_{\eta
_{j}})\int_{Supp\,\eta _{j}}u_{t}^{q_{j}}\,\,\,. 
\]
But for $j\leq m$ we have $\frac{4(q_{j}-1)}{q_{j}{{}^{2}}}\geq c>0$ and 
from proposition 2, $\int_{Supp\,\eta _{j}}u_{t}^{2^{*}}\rightarrow 0$, therefore in (7), 
\[
Q(t,q_{j}-1,\eta _{j})\geq c>0,\,\forall j. 
\]
Thus there exists a neighborhood $V_{j}$ of 1 and a constant $C_{j}>0$ such that for $t\in V_{j}$: 
\[
(\int_{M}(\eta _{j}u_{t}^{\frac{q_{j}}{2}})^{2^{*}})^{\frac{n-2}{n}}\leq
C_{j}\int_{Supp\,\eta _{j}}u_{t}^{q_{j}}\,\,\,\,. 
\]
Then by construction of the $\eta _{j}$ we have
\[
(\int_{\{\eta _{j}=1\}}u_{t}^{q_{j}\frac{2^{*}}{2}})^{\frac{n-2}{n}}\leq
C_{j}\int_{\{\eta _{j-1}=1\}}u_{t}^{q_{j}}\, 
\]
and thus
\[
\left\| u_{t}\right\| _{L^{q}(M\backslash B(x_{0},\delta ))}\leq
C(\prod\limits_{j=1}^{m}C_{j})\left\| u_{t}\right\| _{L^{2}(M)}\,\,\forall
t\in V_{1}\cap ...\cap V_{m}\,. 
\]

Second step: By Gilbarg-Trudinger theorem (8.25) \cite{G-T}, we have : if $u$ is solution of an equation 
$E:\,\,\triangle _{\mathbf{g}}u+h.u=F$, where $\triangle _{\mathbf{g}}+h$ is coercive, and if 
$\omega \subset \subset \omega ^{\prime }$ are two open set, for $r>1,\,q>n/2\,:$%
\begin{equation}
\stackunder{\omega }{Sup}\,u\leq c\left\| u\right\| _{L^{r}(\omega ^{\prime
})}+c^{\prime }\left\| F\right\| _{L^{q}(\omega ^{\prime })}\,\,.
\end{equation}
This theorem is also an application of the iteration process. We apply it to 
$E_{t}:\,\,\triangle _{\mathbf{g}}u_{t}+h_{t}.u_{t}=\lambda _{t}.f.u_{t}^{\frac{n+2}{n-2}}$ 
and to $\omega \subset \subset
\omega ^{\prime }\subset M\backslash \{x_{0}\}$.

Then with the first step applied to $q\frac{n+2}{n-2}$,
and chosing
\[
\omega =M\backslash B(x_{0},\delta ),\,\omega ^{\prime }=M\backslash
B(x_{0},\delta /2),\,r=2,\,q>n/2 
\]
we obtain
\begin{eqnarray*}
\stackunder{M\backslash B(x_{0},\delta )}{Sup}u_{t} & \leq &c\left\|
u_{t}\right\| _{L^{2}(\omega ^{\prime })}+c^{\prime }\lambda _{t}\left\|
u_{t}\right\| _{L^{q\frac{n+2}{n-2}}(\omega ^{\prime })}^{\frac{n+2}{n-2}}
\\ 
& \leq &c\left\| u_{t}\right\| _{L^{2}(M)}+c^{\prime \prime }\left\|
u_{t}\right\| _{L^{2}(M)}^{\frac{n+2}{n-2}}
\end{eqnarray*}
But $\left\| u_{t}\right\| _{L^{2}(M)}\rightarrow 0$, thus the result.

We recall now the notations of subsection (3.2): we consider a sequence of points $(x_{t})$ such that
\[
m_{t}=\stackunder{M}{Max}\,u_{t}=u_{t}(x_{t}):=\mu _{t}^{-\frac{n-2}{2}}. 
\]
From proposition 3, $x_{t}\rightarrow x_{0}$ and $\mu
_{t}\rightarrow 0$. Remember that $\overline{u}_{t},\overline{f}_{t},%
\overline{h}_{t},\,\mathbf{g}\,_{t}$ are the functions and the metric "viewed" in the chart $\exp _{x_{t}}^{-1}$, and $\,\,\widetilde{%
u}_{t}\,,\,\widetilde{h}_{t}\,,\,\widetilde{f}_{t},\widetilde{\,\mathbf{g}\,}%
_{t}$ are the functions and the metric after
blow-up. From now, all the blow-up's will be made on balls $B(x_{t},\delta )$ where $%
f\geq 0$, which is possible as $f(x_{0})>0$.

\begin{proposition}

$\forall R>0$ : $\stackunder{t\rightarrow 1}{\lim }\int_{B(x_{t},R\mu
_{t})}fu_{t}^{2^{*}}dv_{\mathbf{g}}=1-\varepsilon _{R}$ where $\varepsilon
_{R}\stackunder{R\rightarrow +\infty }{\rightarrow }0.$
\end{proposition}

\textit{Proof:} 
This is a direct application of blow-up analysis
in $x_{t}$ with $k_{t}=\mu _{t}^{-1}$:
\begin{center}
$\widetilde{u}_{t}\rightarrow \widetilde{u}=(1+\frac{\lambda f(x_{0})}{n(n-2)%
}\left| x\right| ^{2})^{-\frac{n-2}{2}}\,=(1+\frac{f(x_{0})^{\frac{2}{n}}}{%
K(n,2)^{2}n(n-2)}\left| x\right| ^{2})^{-\frac{n-2}{2}}$ in $\,C_{loc}^{2}(%
\Bbb{R}^{n})$.
\end{center}
Then:
\begin{center}
$\int_{B(x_{t},R\mu _{t})}fu_{t}^{2^{*}}dv_{\mathbf{g}}=\int_{B(0,R)}%
\widetilde{f}_{t}.\widetilde{u}_{t}^{2^{*}}dv_{\widetilde{\,\mathbf{g}\,}%
_{t}}\stackunder{t\rightarrow 1}{\rightarrow }f(x_{0})(\int_{B(0,R)}%
\widetilde{u}^{2^{*}}dx)\stackunder{R\rightarrow \infty }{\rightarrow }1$
\end{center}

\begin{proposition}\textbf{Weak estimates, first part.}

$\exists C>0$ such that $\forall x\in M:\,d_{\mathbf{g}}(x,x_{t})^{\frac{n-2}{2}}u_{t}(x)\leq C$.
\end{proposition}

\textit{Proof :}
Define $w_{t}(x)=\,d_{\mathbf{g}}(x,x_{t})^{%
\frac{n-2}{2}}u_{t}(x)$. We want to prove that there exists $C>0$ such that $\stackunder{M}{Sup}\,w_{t}\leq C$. 
By contradiction, we suppose that (for a subsequence) $\stackunder{M}{Sup}\,w_{t}\rightarrow +\infty $. 
Let $y_{t}$ be a point where $w_{t}$ is maximum. $M$ being compact, $d_{\mathbf{g}}(x,x_{t})$ 
is bounded, therefore $u_{t}(y_{t})\rightarrow \infty $, and thus from proposition 3, $y_{t}\rightarrow x_{0}$. 
Besides, the definition of $\mu _{t}$ gives: 
$
d_{\mathbf{g}}(y_{t},x_{t})\mu _{t}^{-1}\rightarrow +\infty \,\,. 
$

We do a \textit{blow-up} of center $y_{t}$ and coefficient $k_{t}=u_{t}(y_{t})^{\frac{2}{n-2}}$. 
If $x\in B(0,2)$:
$
d_{\mathbf{g}}(x_{t},\exp _{y_{t}}(u_{t}(y_{t})^{-\frac{2}{n-2}}x)  \geq 
d_{\mathbf{g}}(y_{t},x_{t})-2u_{t}(y_{t})^{-\frac{2}{n-2}}  
 \geq u_{t}(y_{t})^{-\frac{2}{n-2}}(w_{t}(y_{t})^{\frac{2}{n-2}}-2)\sim d_{\mathbf{g}}(y_{t},x_{t})
$
as $w_{t}(y_{t})\rightarrow \infty $ and $u_{t}(y_{t})\rightarrow \infty $.
Therefore, for $t$ close to 1: 
$
d_{\mathbf{g}}(x_{t},\exp _{y_{t}}(u_{t}(y_{t})^{-\frac{2}{n-2}}x)\geq \frac{%
1}{2}d_{\mathbf{g}}(y_{t},x_{t})\,\,. 
$

By consequence, for any $R>0$ and $t$ close to 1: 
$
B(y_{t},2u_{t}(y_{t})^{-\frac{2}{n-2}})\cap B(x_{t},R\mu _{t})=\emptyset 
$
Thus, by proposition 4,
\begin{eqnarray*}
\int_{B(0,2)}\widetilde{f}_{t}.\widetilde{u}_{t}^{2^{*}}dv_{\widetilde{\,%
\mathbf{g}\,}_{t}} & =\int_{B(y_{t},2u_{t}(y_{t})^{-\frac{2}{n-2}%
})}fu_{t}^{2^{*}}dv_{\mathbf{g}} & \leq \int_{M\backslash B(x_{t},R\mu
_{t})}fu_{t}^{2^{*}}dv_{\mathbf{g}} \\ 
&  & \leq \int_{M}fu_{t}^{2^{*}}dv_{\mathbf{g}}-\int_{B(x_{t},R\mu
_{t})}fu_{t}^{2^{*}}dv_{\mathbf{g}} \\ 
&  & \stackunder{t\rightarrow 1,R\rightarrow \infty }{\longrightarrow }0
\end{eqnarray*}
But the iteration process then gives that for $1\leq k\leq 2^{*}-1:$%
\[
\int_{B(0,1)}\widetilde{u}_{t}^{\frac{k+1}{2}2^{*}}dv_{\widetilde{\,\mathbf{g%
}\,}_{t}}\rightarrow 0 
\]
and by iteration we obtain that $\forall p\geq 1:$%
\[
\int_{B(0,1)}\widetilde{u}_{t}^{p}dv_{\widetilde{\,\mathbf{g}\,}%
_{t}}\rightarrow 0 
\]
We deduce that $\left\| \widetilde{u}_{t}\right\| _{L^{\infty
}(B(0,1))}\rightarrow 0$ whereas $\widetilde{u}_{t}(0)=1$. Thus a contradiction.

\begin{proposition}\textbf{Weak estimates, second part.}

$\forall \varepsilon >0$ , $\exists R>0$ such that $\forall t,\,\forall x\in
M: $%
\[
\,d_{\mathbf{g}}(x,x_{t})\geq R\mu _{t}\,\Rightarrow \,\,d_{\mathbf{g}%
}(x,x_{t})^{\frac{n-2}{2}}u_{t}(x)\leq \varepsilon . 
\]
\end{proposition}
 
 \textit{Proof :}
 We use the same method, supposing the existence of a $\varepsilon _{0}>0$ and $y_{t}\in M$ such that 
$
\stackunder{t\rightarrow 1}{\lim }d_{\mathbf{g}}(y_{t},x_{t})\mu _{t}^{-1}=+\infty 
\text{ and }w_{t}(y_{t})=\,d_{\mathbf{g}}(y_{t},x_{t})^{\frac{n-2}{2%
}}u_{t}(y_{t})\geq \varepsilon _{0}. 
$
We do a blow-up of center $y_{t}$ and coefficient $k_{t}=u_{t}(y_{t})^{\frac{2}{n-2}}$ and with 
$m_{t}=u_{t}(y_{t}).$
Then, as in proposition 5, for any $R>0$ and $t$ close to 1: 
$
B(y_{t},\frac{1}{2}\varepsilon _{0}^{\frac{2}{n-2}}u_{t}(y_{t})^{-\frac{2}{%
n-2}})\cap B(x_{t},R\mu _{t})=\emptyset .
$
Therefore, as previously: 
$
\int_{B(0,\frac{1}{2}\varepsilon _{0}^{\frac{2}{n-2}})}\widetilde{f}_{t}.%
\widetilde{u}_{t}^{2^{*}}dv_{\widetilde{\,\mathbf{g}\,}_{t}}\rightarrow 0 
$
and we obtain in the same way a contradiction.

\begin{proposition} \textbf{$L^2$-concentration.} 

If $\dim M\geq 4,$ $\forall \delta >0\,:\,$%
\[
\stackunder{t\rightarrow 1}{\lim }\frac{\int_{B(x_{0},\delta )}u_{t}^{2}dv_{%
\mathbf{g}}}{\int_{M}u_{t}^{2}dv_{\mathbf{g}}}=1 
\]
\end{proposition}
 
 \textit{Proof :}
 We first use the two first step of the proof of proposition 3 to show that there exists $c>0$ such that: 
\[
\stackunder{M\backslash B(x_{0},\delta )}{Sup}u_{t}\leq c\left\|
u_{t}\right\| _{L^{2}(M)}\,\,. 
\]
Indeed, going over what we did there,: 
\begin{eqnarray*}
\stackunder{M\backslash B(x_{0},\delta )}{Sup}u_{t} && \leq c\left\|
u_{t}\right\| _{L^{2}(\omega ^{\prime })}+c^{\prime }\lambda _{t}\left\|
u_{t}^{\frac{n+2}{n-2}}\right\| _{L^{q}(\omega ^{\prime })} \\ 
&& \leq c\left\| u_{t}\right\| _{L^{2}(M)}+c^{\prime }\lambda _{t}^{q}%
\stackunder{\omega ^{\prime }}{Sup}(u_{t}^{\frac{n+2}{n-2}-1})\left\|
u_{t}\right\| _{L^{q}(\omega ^{\prime })}\\
&&\leq c^{\prime \prime }\left\|
u_{t}\right\| _{L^{2}(M)}
\end{eqnarray*}
as we know now that $\stackunder{\omega ^{\prime }}{Sup}(u_{t}^{\frac{%
n+2}{n-2}-1})\rightarrow 0$ and that, on the other hand, the first step of the proof of proposition 3 gives 
$\left\| u_{t}\right\| _{L^{q}(\omega ^{\prime })}\leq C\left\|u_{t}\right\| _{L^{2}(M)}$

Third step: Using this: 
\begin{eqnarray}
\left\| u_{t}\right\| _{L^{2}(M\backslash B(x_{0},\delta ))}^{2} && \leq 
\stackunder{M\backslash B(x_{0},\delta )}{Sup}u_{t}.\int_{M\backslash
B(x_{0},\delta )}u_{t}\nonumber \\ 
&& \leq c\left\| u_{t}\right\| _{L^{2}(M)}\left\| u_{t}\right\| _{L^{1}(M)}
\end{eqnarray}
We now want to prove that
\begin{equation}
\left\| u_{t}\right\| _{L^{1}(M)}\leq c\left\| u_{t}\right\|
_{L^{2^{*}-1}(M)}^{2^{*}-1}\,\,.  
\end{equation}
If $h>0$, we get the result by integrating equation $E_{t}$. Otherwise, 
for any $q\in ]2,2^{*}[$, there exists $\varphi >0$ 
solution of $\triangle _{\mathbf{g}}\varphi +h\varphi
=\lambda _{h,f,\,\mathbf{g}\,}.f.\varphi ^{q-1}$. We set
\[
\mathbf{g}^{\prime }=\varphi ^{\frac{4}{n-2}}\mathbf{g}\text{ and }\overline{%
h_{t}}=\frac{\triangle _{\mathbf{g}}\varphi +h_{t}\varphi }{\varphi ^{\frac{%
n+2}{n-2}}} 
\]
Then for $t$ close to 1
$$\overline{h_{t}}=\varphi ^{q-2^{*}}-(h-h_{t})\varphi ^{2-2^{*}}\geq
\varepsilon _{0}>0 $$

Besides, by conformal invariance, and using $E_{t}$, we have: 
\[
\triangle _{\mathbf{g}^{\prime }}\overline{u_{t}}+\overline{h_{t}}.\overline{%
u_{t}}=\lambda _{t}f.\overline{u_{t}}^{\frac{n+2}{n-2}} 
\]
where $\overline{u_{t}}=\varphi ^{-1}.u_{t}$. Integrating, we obtain: 
\[
\varepsilon _{0}\int_{M}\overline{u_{t}}dv_{\mathbf{g}^{\prime }}\leq
\lambda _{t}Supf\int_{M}\overline{u_{t}}^{\frac{n+2}{n-2}}dv_{\mathbf{g}%
^{\prime }} 
\]
and thus there exists $C>0$ such that for $t$ close to 1 
\[
\left\| u_{t}\right\| _{L^{1}(M)}\leq C\left\| u_{t}\right\|
_{L^{2^{*}-1}(M)}^{2^{*}-1} 
\]
where the norms are now relative to $dv_{\mathbf{g}}$.

Fourth step: We conclude using Hölder's inequality. If $n=\dim M\geq 6$ : 
\[
\left\| u_{t}\right\| _{L^{2^{*}-1}(M)}^{2^{*}-1}\leq \left\| u_{t}\right\|
_{L^{2}(M)}^{\frac{n+2}{n-2}}Vol_{\mathbf{g}}(M)^{\frac{n-6}{2(n-2)}}\,. 
\]
With (10) and (11), we obtain :
\[
\stackunder{t\rightarrow 1}{\lim }\frac{\left\| u_{t}\right\|
_{L^{2}(M\backslash B(x_{0},\delta ))}^{2}}{\left\| u_{t}\right\|
_{L^{2}(M)}^{2}}=0 
\]
which proves the result. If $n=5$, Hölder's inequality gives: 
\[
\left\| u_{t}\right\| _{L^{2^{*}-1}(M)}^{2^{*}-1}\leq \left\| u_{t}\right\|
_{L^{2}(M)}^{\frac{3}{2}}\left\| u_{t}\right\| _{L^{2^{*}}(M)}^{\frac{5}{6}} 
\]
and we also conclude using (10) and (11). If now $n=4$,
we have to use proposition 6 and the associated blow-up. 
\begin{proposition} \textbf{Strong estimates.}

For any $\nu$, $0< \nu < n-2$, there exists a constant $C(\nu )>0$ such that
$$
\forall x\in M:\,d_{\mathbf{g}}(x,x_{t})^{n-2-\nu }\mu _{t}^{-\frac{n-2}{2}%
+\nu }u_{t}(x)\leq C(\nu )  
$$
\end{proposition}

\textit{Proof :} The proof requires the use of the Green function and of the weak estimates. The idea is due to 
O. Druet and F. Robert and can be found in [13]. We recall first the property of the Green function.
If $\bigtriangleup _{\mathbf{g}}+h$ is a coercive operator, there exists a unique function (at least $C^{2}$ 
with our hypothesis) 
$
G_{h}:M\times M\backslash \{(x,x),x\in M\}\rightarrow \Bbb{R} 
$
symetric and positive, such that in the sense of distributions, we have: $\forall x\in M$%
\begin{equation}
\bigtriangleup _{\mathbf{g},y}G_{h}(x,y)+h(y)G_{h}(x,y)=\delta _{x}  
\end{equation}
Furthermore, there exists $c>0,\,\rho >0$ such that $\forall (x,y)$ with $0<d_{%
\mathbf{g}}(x,y)<\rho :$%
\begin{equation}
\frac{c}{d_{\mathbf{g}}(x,y)^{n-2}}\leq G_{h}(x,y)\leq \frac{c^{-1}}{d_{%
\mathbf{g}}(x,y)^{n-2}}  
\end{equation}

\begin{equation}
\frac{\left| \nabla _{y}G_{h}(x,y)\right| }{G_{h}(x,y)}\geq \frac{c}{d_{%
\mathbf{g}}(x,y)}  
\end{equation}

\begin{center}
$c$ and $\rho $ vary continuously with $h$ 
\begin{equation}
G_{h}(x,y)d_{\mathbf{g}}(x,y)^{n-2}\rightarrow \frac{1}{(n-2)\omega _{n-1}}%
\text{ when }d_{\mathbf{g}}(x,y)\rightarrow 0  
\end{equation}
\end{center}

To prove these strong estimates, it is sufficient, considering (13), to prove that 
$\mu _{t}^{\frac{n-2}{2}-(n-2)(1-\nu)}u_{t}(x)\leq c^{\prime }G_{h}^{1-\nu}(x,x_{t})$, (just change $\nu$ 
by $(n-2)\nu$). 
First, notice that, using for example the weak estimates, the strong estimates are true in any ball 
$B(x_{t},R\mu _{t})$ where $R$ is fixed. We therefore have to prove the estimates in the manifold with 
boundary $M\backslash B(x_{t},R\mu _{t})$ whose boundary is $b(M\backslash B(x_{t},R\mu
_{t}))=bB(x_{t},R\mu _{t})$. 
For $\nu $ small, there exists $\varepsilon _{0}>0$ such that he operator
\[
\bigtriangleup _{\mathbf{g}}+\frac{h-2\varepsilon _{0}}{1-\nu } 
\]
is still coercive; let $\widetilde{G}$ be its Green function.
To prove our esimate, we apply the maximum principle to : 
$
L_{t}\varphi =\bigtriangleup _{\mathbf{g}}\varphi +h_{t}\varphi -\lambda
_{t}fu_{t}^{2^{*}-2}\varphi 
$
and to $x\longmapsto \widetilde{G}^{1-\nu}(x,x_{t})-c\mu _{t}^{\frac{n-2}{2}-(n-2)(1-\nu)}u_{t}(x)$.
As $L_{t}u_{t}=0$ with $u_{t}>0$, $L_{t}$ satisfies the maximum principle (see [13]). 
Using (12), the fact that $\delta _{x_{t}}(x)=0$ on $M\backslash B(x_{t},R\mu _{t})$ and 
the fact that for $t$ close to 1, $h_{t}-h\geq -\varepsilon _{0}$ (as $h_{t}\rightarrow h$ in $C^{0}$),
some computations give that $\forall x\in M\backslash B(x_{t},R\mu _{t}):$
\begin{equation}
\frac{L_{t}\widetilde{G}^{1-\nu }}{\widetilde{G}^{1-\nu }}(x,x_{t})\geq
\varepsilon _{0}-\lambda _{t}f(x)u_{t}(x)^{2^{*}-2}+\nu (1-\nu )\left| \frac{%
\nabla \widetilde{G}}{\widetilde{G}}\right| ^{2}(x,x_{t})  
\end{equation}
We now separate $M\backslash B(x_{t},R\mu _{t})$ in two parts using a ball 
$B(x_{t},\rho )$ where $\rho >0$ is as in (13) and (14).
For $t$ close to 1, $\rho >R\mu _{t}$. $R>0$ will be fixed later.

1/:As $u_{t}\rightarrow 0$ in $C_{loc}^{0}(M\backslash \{x_{0}\})$,
(16) gives for $t$ close to 1: 
\[
\forall x\in M\backslash B(x_{t},\rho ):\,L_{t}\widetilde{G}^{1-\nu
}(x,x_{t})\geq 0. 
\]

2/: Using the weak estimates (second part), in $B(x_{t},\rho
)\backslash B(x_{t},R\mu _{t})$ :
\[
d_{\mathbf{g}}(x,x_{t})^{2}u_{t}(x)^{2^{*}-2}\leq \varepsilon _{R} 
\]
where $\varepsilon _{R}\stackunder{R\rightarrow \infty }{\rightarrow }0$.
Then, with (14) et (16), for $R$ big enough: 
\begin{eqnarray*}
\frac{L_{t}\widetilde{G}^{1-\nu }}{\widetilde{G}^{1-\nu }}(x,x_{t}) && \geq
\varepsilon _{0}-\lambda _{t}f(x)u_{t}(x)^{2^{*}-2}+\nu (1-\nu )\frac{c}{d_{%
\mathbf{g}}(x,x_{t})^{2}} \\ 
&& \geq \varepsilon _{0}-\lambda _{t}(\stackunder{B(x_{t},\rho )}{Sup}f).%
\frac{\varepsilon _{R}}{d_{\mathbf{g}}(x,x_{t})^{2}}+\nu (1-\nu )\frac{c}{d_{%
\mathbf{g}}(x,x_{t})^{2}} \\ 
&& \geq \varepsilon _{0}+\frac{c^{\prime }}{d_{\mathbf{g}}(x,x_{t})^{2}}\geq 0
\end{eqnarray*}

We have proved that in $M\backslash B(x_{t},R\mu _{t})$ and for any constant $C_{t}>0$ 
which can depend of $t$ :
\[
L_{t}(C_{t}.\widetilde{G}^{1-\nu }(x,x_{t}))=C_{t}.L_{t}\widetilde{G}^{1-\nu
}(x,x_{t})\geq 0=L_{t}u_{t} 
\]
At last, on the boundary $b(M\backslash B(x_{t},R\mu _{t}))$, using (13), we obtain :
\[
\widetilde{G}^{1-\nu }(x,x_{t})\geq \frac{c}{d_{\mathbf{g}%
}(x,x_{t})^{(n-2)(1-\nu )}}=\frac{c}{(R\mu _{t})^{(n-2)(1-\nu )}}\,\,. 
\]
So, if we let $C_{t}=c^{-1}R^{(n-2)(1-\nu )}\mu _{t}^{(n-2)(1-\nu )-\frac{%
n-2}{2}},$ we have for $x\in bB(x_{t},R\mu _{t})=b(M\backslash B(x_{t},R\mu
_{t})):$%
\[
C_{t}.\widetilde{G}^{1-\nu }(x,x_{t})\geq \mu _{t}^{-\frac{n-2}{2}%
}=Sup\,u_{t}\geq u_{t}(x) 
\]
Therefore, by the maximum principle :
\[
C_{t}.\widetilde{G}^{1-\nu }(x,x_{t})\geq u_{t}(x)\text{ in }M\backslash
B(x_{t},R\mu _{t}) 
\]
which can be rewriten
\[
\widetilde{G}^{1-\nu }(x,x_{t})\geq C_{t}^{-1}u_{t}(x)=c\text{ }\mu _{t}^{%
\frac{n-2}{2}-(n-2)(1-\nu )}u_{t}(x) 
\]
and therefore, using (13) : 
\[
d_{\mathbf{g}}(x,x_{t})^{(n-2)(1-\nu )}\mu _{t}^{\frac{n-2}{2}-(n-2)(1-\nu
)}u_{t}(x)\leq c 
\]
which gives the strong estimates by changing $\nu $ in $(n-2)\nu $.

\begin{proposition}\textbf{Corollary: Strong $L^p$-concentration.}

$\forall R>0$, $%
\forall \delta >0$\ and $\forall p>\frac{n}{n-2}$
$$
\stackunder{t\rightarrow 1}{\lim }\frac{\int_{B(x_{t},R\mu
_{t})}u_{t}^{p}dv_{\mathbf{g}}}{\int_{B(x_{t},\delta )}u_{t}^{p}dv_{\mathbf{g%
}}}=1-\varepsilon _{R}  \text{ where } \varepsilon _{R}\stackunder{%
R\rightarrow +\infty }{\rightarrow }0.
$$
\end{proposition}

\textit{Proof :}
Just apply the strong estimates to a blow-up in $x_{t}.$ By blow-up formulae
\begin{eqnarray*}
\int_{M}u_{t}^{p}dv_{\mathbf{g}} && \geq \int_{B(x_{t},\mu _{t})}u_{t}^{p}dv_{%
\mathbf{g}}=\mu _{t}^{n-\frac{n-2}{2}p}\int_{B(0,1)}\widetilde{u}_{t}^{p}dv_{%
\widetilde{g}_{t}} \\ 
& &\geq C\mu _{t}^{n-\frac{n-2}{2}p}
\end{eqnarray*}
On the other hand, by the strong estimates: 
\begin{eqnarray*}
\int_{M\backslash B(x_{t},R\mu _{t})}u_{t}^{p}dv_{\mathbf{g}} && \leq C\mu
_{t}^{p\frac{n-2}{2}}\int_{M\backslash B(x_{t},R\mu _{t})}d_{\mathbf{g}%
}(y_{t},x)^{(2-n)p}dv_{\mathbf{g}} \\ 
&& \leq C\mu _{t}^{n-p\frac{n-2}{2}}R^{n+(2-n)p}
\end{eqnarray*}
as soon as $p>\frac{n}{n-2}$. Dividing, we obtain the corollary.
\\

At this point, to carry on the proof of theorem 1, we need a powerfull extension of the strong estimates, called 
$C^0-theory$, which is in fact a complete control of the sequence 
$d_{\mathbf{g}}(x,x_{t})^{n-2}\mu _{t}^{-\frac{n-2}{2}}u_{t}(x)$; it is expressed by the next theorem of Druet 
and Robert, and proved in arbitrary energy in [O. Druet, E. Hebey, and F. Robert, "Blow-up Theory for Elliptic PDEs in Riemannian Geometry", Mathematical notes, Princeton University Press, vol. 45, 2004.]

Another approach, also accessible at this point and originally used in the author's PHD thesis, 
is to prove another very important estimate 
concerning the "speed" of convergence of $(x_{t})$ to $x_{0}$, but it requires the additional hypothesis that the Hessian of $f$ is non-degenerate 
at the points of maximum of $f$; it will be our theorem 6, whose proof is independent of the theorem of Druet-Robert, 
only requiring the results up to proposition 9, and appears as a byproduct of an alternative proof of theorem 1. It is 
however of independent interest, as it is a very important estimate concerning concentration phenomena's which has been 
studied by various authors.

We now state the theorem of Druet and Robert and refer for its proof to the reference cited above, the function $f$ introducing no difficulties. 
It says first that one can take $\nu =0$ in the strong estimates, but also that one has somehow the reverse estimate.
\\

\begin{theor} [Druet, Robert]

For any $\ve>0$, there exist $\delta_{\ve}>0$ such that, up to a subsequence, for any $t$ and any 
$x\in B(x_{0},\delta_{\ve})$ :
$$(1-\ve) B_{t}(x) \leq u_{t}(x) \leq (1+\ve) B_{t}(x)$$
where 
$$ B_{t}(x)=\mu_{t}^{-\frac{n-2}{2}}\Bigl( 1+\frac{\lambda f(x_{0})}{n(n-2)} \frac{d_{\g}(x_{t},x)^2}{\mu_{t}^2}\Bigr) ^{-\frac{n-2}{2}} $$
is the "standard bubble".
\end{theor}
Note that in the proof of theorem 1, we will need the minoration:\\
 $(1-\ve) B_{t}(x) \leq u_{t}(x)$, which is 
a stronger result than $u_{t}(x) \leq (1+\ve) B_{t}(x)$ which must first be proved to get the minoration.

Finally, we come to our main result concerning the concentration phenomenom, which is the "missing link" between 
the sequence $(x_{t})$ and $x_{0}$.
\begin{theo}\textbf{"Second fundamental estimate".}
Suppose that $dimM\geq5$ and that the hessian of the function $f$ is non-degenerate at each of its points of maximum. Then, there 
exist a constant $C$ such that for all $t$ :
$$\frac{d_{\mathbf{g}}(x_{t},x_{0})}{\mu _{t}}\leq C .$$
Moreover, if for each point $P$ of maximum of $f$ we have
$$h(P)=\frac{n-2}{4(n-1)}S_{\mathbf{g}}(P)-\frac{(n-2)(n-4)}{8(n-1)}\frac{%
\bigtriangleup _{\mathbf{g}}f(P)}{f(P)} ,$$
then more precisely
$$\frac{d_{\mathbf{g}}(x_{t},x_{0})}{\mu _{t}}\rightarrow 0 .$$
\end{theo}
To understand the significance of this theorem, note that the weak and strong estimates, the strong $L^p$-concentration 
and the estimates in the theorem of Druet-Robert, are "centered" in $x_{t}$. Theorem 6 allows one to "translate" 
these estimates in $x_{0}$ in the sense that one can now replace $x_{t}$ by $x_{0}$. This estimate, called by Zoé Faget "second fundamental 
estimate", (the "first one" being the strong estimate), joined with the estimates of $C^0-theory$ presented in 
the theorem of Druet and Robert above, 
gives a complete description of the behavior of a sequence of solutions of equations 
$\bigtriangleup _{\mathbf{g}}u_{t}+h_{t}u_{t}=\lambda _{t}fu_{t}^{\frac{n+2}{n-2}}$ 
in the spirit of the study of Palais-Smale sequences associated to these equations. It has been studied, for example, by 
Druet and Robert in the case $f=constant=1$ in [13] where they require strong hypothesis on the shape of the functions 
$h_{t}$ and on the geometry of the manifold near the concentration point, or by Hebey in the euclidean setting.
Intuitively, it seems that our hypothesis on $f$ "fixes" the position of the concentration point, and so we get a control 
on the distance between $x_{t}$ and $x_{0}$. Also, our method seems to be applicable to other settings, 
see e.g. [15].

\subsection{Proof of theorem 1}
Remember that $\overline{u}_{t},\overline{f}_{t},\overline{h}_{t},\,\mathbf{g}\,_{t}$ 
are the functions and the metric "viewed" in the chart $\exp _{x_{t}}^{-1}$, and 
$\,\,\widetilde{u}_{t}\,,\,\widetilde{h}_{t}\,,\,\widetilde{f}_{t},\widetilde{\,\mathbf{g}\,}_{t}$ 
are the functions and the metric after
blow-up with center $x_{t}$ and coefficient $k_{t}=\mu_{t}^{-1}$. From now, all the blow-up's will be made on balls $B(x_{t},\delta )$ where $%
f\geq 0$, which is possible as $f(x_{0})>0$.

Let also $\eta $ be a cut-off function on $\Bbb{R}^{n}$ equal to 1 on the euclidean ball $B(0,\delta /2)$, 
and equal to 0 on $\Bbb{R}^{n}\backslash B(0,\delta )$, $0\leq \eta \leq 1$ with $\left|
\nabla \eta \right| \leq C.\delta ^{-1}$ where $\delta $ is chosen small enough to have $f\geq 0$ on the balls 
$B(x_{t},\delta )$. 
The Sobolev inequality gives on the one hand
\begin{equation}
(\int_{B(0,\delta )}(\eta \overline{u}_{t})^{2^{*}}dx)^{\frac{2}{2^{*}}}\leq
K(n,2)^{2}\int_{B(0,\delta )}\left| \nabla (\eta \overline{u}_{t})\right|
_{e}^{2}dx\,\, 
\end{equation}
where $\left| .\right| _{e}$ is the euclidean metric of associated measure $dx$.

On the other hand, integration by part gives, noting that $\left|
\nabla \eta \right| =$ $\Delta \eta =0$ on $B(0,\delta /2)$ : 
\[
\int_{B(0,\delta )}\left| \nabla (\eta \overline{u}_{t})\right|
_{e}^{2}dx\leq \int_{B(0,\delta )}\eta ^{2}\overline{u}_{t}\bigtriangleup
_{e}\overline{u}_{t}dx+C.\delta ^{-2}\int_{B(0,\delta )\backslash B(0,\delta
/2)}\overline{u}_{t}^{2}dx 
\]
Noting $\,\mathbf{g}\,_{t}^{ij}$ the components of  $\,\mathbf{g}\,_{t}$
and $\Gamma (\,\mathbf{g}\,_{t})_{ij}^{k}$ the associated Christoffel symbols, we write : 
\[
\bigtriangleup _{e}\overline{u}_{t}=\bigtriangleup _{\mathbf{g}_{t}}%
\overline{u}_{t}+(\,\mathbf{g}\,_{t}^{ij}-\delta ^{ij})\partial _{ij}%
\overline{u}_{t}-\,\mathbf{g}\,_{t}^{ij}\Gamma (\,\mathbf{g}%
\,_{t})_{ij}^{k}\partial _{k}\overline{u}_{t} 
\]
We get from this inequallity, using 
using this expression of the laplacian, equation $E_{t}:\,\,\triangle _{%
\mathbf{g}}u_{t}+h_{t}.u_{t}=\lambda _{t}.f.u_{t}^{\frac{n+2}{n-2}}$ ``viewed''
in the chart exp$_{x_{t}}^{-1}$, and using the fact that $\left| \nabla \eta \right|
=\Delta \eta =0$ on $B(0,\delta /2)$ and with some integration by parts: 
\begin{eqnarray*}
\int_{B(0,\delta )}\left| \nabla (\eta \overline{u}_{t})\right|
_{e}^{2}dx\leq & &\lambda _{t}\int_{B(0,\delta )}\eta ^{2}\overline{f}_{t}%
\overline{u}_{t}^{2^{*}}dx-\int_{B(0,\delta )}\eta ^{2}\overline{h}_{t}%
\overline{u}_{t}^{2}dx\\
&&+C.\delta ^{-2}\int_{B(0,\delta )\backslash B(0,\delta
/2)}\overline{u}_{t}^{2}dx \\ 
&& -\int_{B(0,\delta )}\eta ^{2}(\,\mathbf{g}\,_{t}^{ij}-\delta
^{ij})\partial _{i}\overline{u}_{t}\partial _{j}\overline{u}_{t}dx\\
&&+\frac{1}{2}\int_{B(0,\delta )}(\partial _{k}(\,\mathbf{g}\,_{t}^{ij}\Gamma (\,\mathbf{g%
}\,_{t})_{ij}^{k}+\partial _{ij}\,\mathbf{g}\,_{t}^{ij})(\eta \overline{u}%
_{t}^{2})dx\,.
\end{eqnarray*}
Using the Sobolev inequality (17) and the fact that 
 $\lambda _{t}\leq \frac{1}{K(n,2){{}^{2}}(%
\stackunder{M}{Sup}f)^{\frac{n-2}{n}}}$, we obtain at last: 

\begin{equation}
\int_{B(0,\delta )}\overline{h}_{t}(\eta \overline{u}_{t})^{2}dx\leq A_{t}
+B_{t}+C_{t}+C.\delta ^{-2}\int_{B(0,\delta )\backslash B(0,\delta /2)}\overline{u}%
_{t}^{2}dx
\end{equation}
where:

$B_{t}=\frac{1}{2}\int_{B(0,\delta )}(\partial _{k}(\,\mathbf{g}%
\,_{t}^{ij}\Gamma (\,\mathbf{g}\,_{t})_{ij}^{k}+\partial _{ij}\,\mathbf{g}%
\,_{t}^{ij})(\eta \overline{u}_{t}^{2})dx$

$C_{t}=\left| \int_{B(0,\delta )}\eta ^{2}(\,\mathbf{g}\,_{t}^{ij}-\delta
^{ij})\partial _{i}\overline{u}_{t}\partial _{j}\overline{u}_{t}dx\right| $

$A_{t}=\frac{1}{K(n,2){{}^{2}}(\stackunder{M}{Sup}f)^{\frac{n-2}{n}}}%
\int_{B(0,\delta )}\overline{f}_{t}\eta ^{2}\overline{u}_{t}^{2^{*}}dx-\frac{%
1}{K(n,2){{}^{2}}}(\int_{B(0,\delta )}(\eta \overline{u}_{t})^{2^{*}}dx)^{%
\frac{2}{2^{*}}}$

These computations were developed in the article of Djadli and Druet \cite{D-D}. Our goal is to use $L{{}^{2}}$-concentration (proposition 7) 
to obtain a contradiction; we shall divide (18) by $%
\int_{B(0,\delta )}\overline{u}_{t}^{2}dx$ and take the limit when $t\rightarrow
t_{0}=1$.

$L{{}^{2}}$-concentration first gives : 
\[
\frac{C.\delta ^{-2}\int_{B(0,\delta )\backslash B(0,\delta /2)}\overline{u}%
_{t}^{2}dx}{\int_{B(0,\delta )}\overline{u}_{t}^{2}dx}\stackunder{%
t\rightarrow 1}{\rightarrow }0\,\,. 
\]

Z.Djadli and O.Druet \cite{D-D} showed (see also [10] for full details):
\[
\stackunder{t\rightarrow 1}{\overline{\lim }}\frac{C_{t}}{\int_{B(0,\delta )}%
\overline{u}_{t}^{2}dx}\leq \varepsilon _{\delta }\text{ where }\varepsilon
_{\delta }\rightarrow 0\text{ when }\delta \rightarrow 0\,\,.
\]

Furthermore, as $x_{t}\rightarrow x_{0}$ we have $\stackunder{t\rightarrow 1}{%
\lim }(\partial _{k}(\,\mathbf{g}\,_{t}^{ij}\Gamma (\,\mathbf{g}%
\,_{t})_{ij}^{k}+\partial _{ij}\,\mathbf{g}\,_{t}^{ij})(0)=\frac{1}{3}S_{\,%
\mathbf{g}\,}(x_{0})$, therefore, using $L{{}^{2}}$-concentration : 
\[
\stackunder{t\rightarrow 1}{\overline{\lim }}\frac{B_{t}}{\int_{B(0,\delta )}%
\overline{u}_{t}^{2}dx}=\frac{1}{6}S_{\,\mathbf{g}\,}(x_{0})+\varepsilon
_{\delta }\,\,. 
\]

It is the expression $A_{t}$ which will give $\frac{n-2}{4(n-1)}%
S_{\,\mathbf{g}\,}(x_{0})$ $-\frac{1}{6}S_{\,\mathbf{g}\,}(x_{0})$ and $\frac{%
(n-2)(n-4)}{8(n-1)}\frac{\bigtriangleup _{\mathbf{g}}f(x_{0})}{f(x_{0})}$.

By Hölder's inequality: 
\[
\int_{B(0,\delta )}\overline{f}_{t}\eta ^{2}\overline{u}_{t}^{2^{*}}dx\leq
(\int_{B(0,\delta )}\overline{f}_{t}\overline{u}_{t}^{2^{*}}dx)^{\frac{2}{n}%
}(\int_{B(0,\delta )}\overline{f}_{t}(\eta \overline{u}_{t})^{2^{*}}dx)^{%
\frac{n-2}{n}}\, 
\]
Beside : 
\[
\,dx\leq (1+\frac{1}{6}Ric(x_{t})_{ij}x^{i}x^{j}+C\left| x\right| ^{3})dv_{\,%
\mathbf{g}\,_{t}} 
\]
Using this development and $(1+x)^\alpha \leq 1+\alpha x$ for $0<\alpha \leq1$: 
\[
(\int_{B(0,\delta )}\overline{f}_{t}\overline{u}_{t}^{2^{*}}dx)^{\frac{2}{n}%
}\leq (\int_{B(0,\delta )}\overline{f}_{t}\overline{u}_{t}^{2^{*}}dv_{\,%
\mathbf{g}\,_{t}})^{\frac{2}{n}}+\frac{1}{(\int_{B(0,\delta )}\overline{f}%
_{t}\overline{u}_{t}^{2^{*}}dv_{\,\mathbf{g}\,_{t}})^{\frac{n-2}{n}}}\frac{2%
}{n}\{S_{t}\}+C\{S_{t}\}^{2} 
\]
where 
\[
\{S_{t}\}=\frac{1}{6}Ric(x_{t})_{ij}\int_{B(0,\delta )}x^{i}x^{j}\overline{f}%
_{t}\overline{u}_{t}^{2^{*}}dv_{\,\mathbf{g}\,_{t}}+C\int_{B(0,\delta
)}\left| x\right| ^{3}\overline{u}_{t}^{2^{*}}dv_{\,\mathbf{g}\,_{t}}\,\,. 
\]
We deduce
\[
A_{t}\leq \frac{1}{K(n,2){{}^{2}}(\stackunder{M}{Sup}f)^{\frac{n-2}{n}}}%
(A_{t}^{1}+A_{t}^{2}) 
\]
where
\[
A_{t}^{1}=(\int_{B(0,\delta )}\overline{f}_{t}\overline{u}_{t}^{2^{*}}dv_{\,%
\mathbf{g}\,_{t}})^{\frac{2}{n}}(\int_{B(0,\delta )}\overline{f}_{t}(\eta 
\overline{u}_{t})^{2^{*}}dx)^{\frac{n-2}{n}}\,-(Supf.\int_{B(0,\delta
)}(\eta \overline{u}_{t})^{2^{*}}dx)^{\frac{n-2}{n}} 
\]
and
$$
A_{t}^{2}=\frac{2(\int_{B(0,\delta )}\overline{f}_{t}(\eta \overline{u}%
_{t})^{2^{*}}dx)^{\frac{n-2}{n}}}{n(\int_{B(0,\delta )}\overline{f}_{t}%
\overline{u}_{t}^{2^{*}}dv_{\,\mathbf{g}\,_{t}})^{\frac{n-2}{n}}}
\{\frac{1}{6}Ric(x_{t})_{ij}\int_{B(0,\delta )}x^{i}x^{j}\overline{f}_{t}%
\overline{u}_{t}^{2^{*}}dv_{\,\mathbf{g}\,_{t}}
+C\int_{B(0,\delta )}\left|
x\right| ^{3}\overline{u}_{t}^{2^{*}}dv_{\,\mathbf{g}\,_{t}}\}(1+\varepsilon
_{\delta }) 
$$
as $\{S_{t}\}\rightarrow 0$ when $\delta \rightarrow 0\,$uniformly in $t$.  $A_{t}^{2}$ will give, by 
developing the metric, $S_{\,\mathbf{g}%
\,}(x_{0})\,$while $A_{t}^{1}$ will give, by developing $f$, $-\frac{(n-2)(n-4)}{8(n-1)}%
\frac{\bigtriangleup _{\mathbf{g}}f(x_{0})}{f(x_{0})}$.

Note that for any $\alpha \in H_{1}^{2}(B(x_{0},2\delta )):$%
\[
\stackunder{t\rightarrow 1}{\lim }\frac{\int_{B(x_{t},\delta )}\alpha dx}{%
\int_{B(x_{t},\delta )}\alpha dv_{\,\mathbf{g}\,_{t}}}=1+O(\delta
^{2})=1+\varepsilon _{\delta } 
\]

We start by studying $A_{t}^{2}$:

1/: We have $\stackunder{t\rightarrow 1}{\lim }\frac{(\int_{B(0,\delta )}%
\overline{f}_{t}(\eta \overline{u}_{t})^{2^{*}}dx)^{\frac{n-2}{n}}}{%
(\int_{B(0,\delta )}\overline{f}_{t}\overline{u}_{t}^{2^{*}}dv_{\,\mathbf{g}%
\,_{t}})^{\frac{n-2}{n}}}=1+\varepsilon _{\delta }$

2/: Using the weak estimates (proposition 5), $\left| x\right| ^{2}%
\overline{u}_{t}^{2^{*}}\leq c\overline{u}_{t}^{2}$, from where we get: 
\[
\frac{\int_{B(0,\delta )}\left| x\right| ^{3}\overline{u}_{t}^{2^{*}}dv_{\,%
\mathbf{g}\,_{t}}}{\int_{B(0,\delta )}\overline{u}_{t}^{2}dv_{\,\mathbf{g}%
\,_{t}}}\leq C.\varepsilon _{\delta }\,\,. 
\]

3/: Using the blow-up formula's we write: for all $R>0:$%
\begin{eqnarray*}
\int_{B(0,\delta )}x^{i}x^{j}\overline{f}_{t}\overline{u}_{t}^{2^{*}}dv_{\,%
\mathbf{g}\,_{t}} && =\int_{B(0,R\mu _{t})}x^{i}x^{j}\overline{f}_{t}%
\overline{u}_{t}^{2^{*}}dv_{\,\mathbf{g}\,_{t}}+\int_{B(0,\delta )\backslash
B(0,R\mu _{t})}x^{i}x^{j}\overline{f}_{t}\overline{u}_{t}^{2^{*}}dv_{\,%
\mathbf{g}\,_{t}} \\ 
&& =\mu _{t}^{2}\int_{B(0,R)}x^{i}x^{j}\widetilde{f_{t}}\widetilde{u}%
_{t}^{2^{*}}dv_{\widetilde{\,\mathbf{g}\,}_{t}}+\mu _{t}^{2}\int_{B(0,\delta
\mu _{t}^{-1})\backslash B(0,R)}x^{i}x^{j}\widetilde{f_{t}}\widetilde{u}%
_{t}^{2^{*}}dv_{\widetilde{\,\mathbf{g}\,}_{t}}
\end{eqnarray*}
and
\[
\int_{B(0,\delta )}\overline{u}_{t}^{2}dv_{\mathbf{g}_{t}}=\mu
_{t}^{2}\int_{B(0,\delta \mu _{t}^{-1})}\widetilde{u}_{t}^{2}dv_{\widetilde{%
\mathbf{g}}_{t}}\,\,. 
\]
Using the weak estimates again, we get : 
\[
\int_{B(0,\delta \mu _{t}^{-1})\backslash B(0,R)}x^{i}x^{j}\widetilde{f_{t}}%
\widetilde{u}_{t}^{2^{*}}dv_{\widetilde{\mathbf{g}}_{t}}\leq \varepsilon
_{R}.\int_{B(0,\delta \mu _{t}^{-1})\backslash B(0,R)}\widetilde{u}%
_{t}^{2}dv_{\widetilde{\mathbf{g}}_{t}} 
\]
thus:
\[
\frac{\int_{B(0,\delta \mu _{t}^{-1})\backslash B(0,R)}x^{i}x^{j}\widetilde{%
f_{t}}\widetilde{u}_{t}^{2^{*}}dv_{\widetilde{\mathbf{g}}_{t}}}{%
\int_{B(0,\delta \mu _{t}^{-1})}\widetilde{u}_{t}^{2}dv_{\widetilde{\mathbf{g%
}}_{t}}}\leq \varepsilon _{R} 
\]
where $\varepsilon _{R}\rightarrow 0$ when $R\rightarrow +\infty $.
Now, if $i\neq j:$%
\[
\stackunder{t\rightarrow 1}{\overline{\lim }}\frac{|\int_{B(0,R)}x^{i}x^{j}%
\widetilde{f_{t}}\widetilde{u}_{t}^{2^{*}}dv_{\widetilde{\mathbf{g}}_{t}}|}{%
\int_{B(0,\delta \mu _{t}^{-1})}\widetilde{u}_{t}^{2}dv_{\widetilde{\mathbf{g%
}}_{t}}}\leq \stackunder{t\rightarrow 1}{\overline{\lim }}\frac{%
|\int_{B(0,R)}x^{i}x^{j}\widetilde{f_{t}}\widetilde{u}_{t}^{2^{*}}dv_{%
\widetilde{\mathbf{g}}_{t}}|}{\int_{B(0,R)}\widetilde{u}_{t}^{2}dv_{%
\widetilde{\mathbf{g}}_{t}}}=0 
\]
because
\[
\widetilde{u}_{t}\rightarrow \widetilde{u}=(1+\frac{f(x_{0})^{\frac{2}{n}}}{%
K(n,2)^{2}n(n-2)}\left| x\right| ^{2})^{-\frac{n-2}{2}}\,\,\,in\,\,%
\,C^{0}(B(0,R)) 
\]
and $\widetilde{u}$ is radial (see subsection 3.2).

If $i=j:$%
\[
\frac{\int_{B(0,R)}x^{i}x^{i}\widetilde{f_{t}}\widetilde{u}_{t}^{2^{*}}dv_{%
\widetilde{\mathbf{g}}_{t}}}{\int_{B(0,\delta \mu _{t}^{-1})}\widetilde{u}%
_{t}^{2}dv_{\widetilde{\mathbf{g}}_{t}}}=\frac{\int_{B(0,R)}(x^{i}){{}^{2}}%
\widetilde{f_{t}}\widetilde{u}_{t}^{2^{*}}dv_{\widetilde{\mathbf{g}}_{t}}}{%
\int_{B(0,R)}\widetilde{u}_{t}^{2}dv_{\widetilde{\mathbf{g}}_{t}}}.\frac{%
\int_{B(0,R)}\widetilde{u}_{t}^{2}dv_{\widetilde{\mathbf{g}}_{t}}}{%
\int_{B(0,\delta \mu _{t}^{-1})}\widetilde{u}_{t}^{2}dv_{\widetilde{\mathbf{g%
}}_{t}}} 
\]
But as soon as $n>4$, using strong $L{{}^{2}}$-concentration 
(proposition 9), we obtain: 
\[
\stackunder{R\rightarrow \infty }{\lim }\,\stackunder{t\rightarrow 1}{%
\overline{\lim }}\frac{\int_{B(0,R)}\widetilde{u}_{t}^{2}dv_{\widetilde{%
\mathbf{g}}_{t}}}{\int_{B(0,\delta \mu _{t}^{-1})}\widetilde{u}_{t}^{2}dv_{%
\widetilde{\mathbf{g}}_{t}}}=1 
\]
therefore

\[
\stackunder{R\rightarrow \infty }{\lim }\,\stackunder{t\rightarrow 1}{%
\overline{\lim }}\frac{\int_{B(0,R)}x^{i}x^{i}\widetilde{f_{t}}\widetilde{u}%
_{t}^{2^{*}}dv_{\widetilde{\mathbf{g}}_{t}}}{\int_{B(0,\delta \mu _{t}^{-1})}%
\widetilde{u}_{t}^{2}dv_{\widetilde{\mathbf{g}}_{t}}}=f(x_{0})\frac{\int_{%
\Bbb{R}^{n}}(x^{i})^{2}.\widetilde{u}^{2}dx}{\int_{\Bbb{R}^{n}}\widetilde{u}%
^{2}dx}=f(x_{0})^{\frac{n-2}{n}}K(n,2)^{2}\frac{n(n-4)}{4(n-1)} 
\]
and thus
\[
\stackunder{t\rightarrow 1}{\overline{\lim }}\frac{1}{f(x_{0})^{\frac{n-2}{n}%
}K(n,2)^{2}}\frac{A_{t}^{2}}{\int_{B(0,\delta )}\overline{u}_{t}^{2}dv_{%
\mathbf{g}_{t}}}=\frac{n-4}{12(n-1)}S_{\mathbf{g}}(x_{0})+\varepsilon
_{\delta } 
\]
which, with $\stackunder{t\rightarrow 1}{\overline{\lim }}\frac{B_{t}}{%
\int_{B(0,\delta )}\overline{u}_{t}^{2}dx}=\frac{1}{6}S_{\mathbf{g}%
}(x_{0})+\varepsilon _{\delta }$ gives 
\[
\stackunder{t\rightarrow 1}{\overline{\lim }}(\frac{1}{f(x_{0})^{\frac{n-2}{n%
}}K(n,2)^{2}}\frac{A_{t}^{2}}{\int_{B(0,\delta )}\overline{u}_{t}^{2}dv_{%
\mathbf{g}_{t}}}+\frac{B_{t}}{\int_{B(0,\delta )}\overline{u}_{t}^{2}dx})=%
\frac{n-2}{4(n-1)}S_{\mathbf{g}}(x_{0})+\varepsilon _{\delta } 
\]

If $n=4$ we write: 
\[
\stackunder{R\rightarrow \infty }{\lim }\stackunder{t\rightarrow 1}{%
\overline{\lim }}\frac{\int_{B(0,R)}x^{i}x^{i}\widetilde{f_{t}}\widetilde{u}%
_{t}^{2^{*}}dv_{\widetilde{\mathbf{g}}_{t}}}{\int_{B(0,\delta \mu _{t}^{-1})}%
\widetilde{u}_{t}^{2}dv_{\widetilde{\mathbf{g}}_{t}}}\leq f(x_{0})^{\frac{n-2%
}{n}}K(n,2)^{2}\frac{n-4}{4(n-1)} 
\]
and we get the conclusion by distinguishing two cases, $S_{\mathbf{g}}(x_{0})<0$ or $S_{%
\mathbf{g}}(x_{0})\geq 0$, the proof being finished as $%
\frac{\bigtriangleup _{\mathbf{g}}f(x_{0})}{f(x_{0})}$ does not appear in dimension 4 (see the end of the proof).

Let us now consider $A^1_{t}$.
\[
A_{t}^{1}=(\int_{B(0,\delta )}\overline{f}_{t}\overline{u}_{t}^{2^{*}}dv_{%
\mathbf{g}_{t}})^{\frac{2}{n}}(\int_{B(0,\delta )}\overline{f}_{t}(\eta 
\overline{u}_{t})^{2^{*}}dx)^{\frac{n-2}{n}}\,-(Supf.\int_{B(0,\delta
)}(\eta \overline{u}_{t})^{2^{*}}dx)^{\frac{n-2}{n}}\,\,. 
\]
We write $f_{t}=f(x_{0})+g_{t}$. Remembering that $f(x_{0})=Sup f$, we have $g_{t}(x_{0})=0$ and 
$g_{t}\leq0$. Using $(1+x)^\alpha \leq 1+\alpha x$ for $0<\alpha \leq1$:
\[
(\int_{B(0,\delta )}\overline{f}_{t}(\eta \overline{u}_{t})^{2^{*}}dx)^{%
\frac{n-2}{n}}\leq (\int_{B(0,\delta )}f(x_{0})(\eta \overline{u}%
_{t})^{2^{*}}dx)^{\frac{n-2}{n}}+\frac{n-2}{n}\frac{\int_{B(0,\delta )}\overline{g}_{t}(\eta \overline{u}_{t})^{2^{*}}dx}
{(\int_{B(0,\delta)}f(x_{0})(\eta \overline{u}_{t})^{2^{*}}dx)^{\frac{2}{n}}}%
\]
where $\overline{g}_{t}$ is $g_{t}$ in the exponential chart in $x_{t}$. 
We now use the theorem of Druet and Robert to write in $B(0,\delta)$: 
$$\overline{u}_{t}\geq (1-\ve_{\delta})\overline{B}_{t},$$
where $\overline{B}_{t}$ is $B_{t}$ in the exponential chart in $x_{t}$.
Because $g_{t}\leq0$, we have:
$$\int_{B(0,\delta )}\overline{g}_{t}(\eta \overline{u}_{t})^{2^{*}}dx\leq
(1-\ve_{\delta})\int_{B(0,\delta )}\overline{g}_{t}(\eta \overline{B}_{t})^{2^{*}}dx.$$
Combining this with the expansion above and the fact that 
$\int_{B(0,\delta )}\overline{f}_{t}\overline{u}_{t}^{2^{*}}dv_{\mathbf{g}_{t}}\leq1$, we obtain:
$$A^1_{t}\leq(1-\ve_{\delta})\frac{n-2}{n}\frac{\int_{B(0,\delta )}\overline{g}_{t}(\eta \overline{B}_{t})^{2^{*}}dx}
{(\int_{B(0,\delta)}f(x_{0})(\eta \overline{u}_{t})^{2^{*}}dx)^{\frac{2}{n}}}$$
We now expand $g_{t}$ noting that 
$\partial _{i}\overline{g}_{t}=\partial _{i}\overline{f}_{t}$ and 
$\partial _{ij}^{2}\overline{g}_{t}=\partial _{ij}^{2}\overline{f}_{t}$.
\[
\overline{g}_{t}(x)\leq g_{t}(x_{t})+x^{i}\partial _{i}\overline{f}_{t}(x_{t})+\frac{1}{2}\partial _{kl}\overline{f}%
_{t}(x_{t}).x^{k}x^{l}+c\left|x\right| ^{3}
\]
Thus
\begin{eqnarray*}
\int_{B(0,\delta )}\overline{g}_{t}(\eta \overline{B}_{t})^{2^{*}}dx\leq
&& g_{t}(x_{t})\int_{B(0,\delta)}(\eta \overline{B}_{t})^{2^{*}}dx\\
&&+\partial _{i}\overline{f}_{t}(x_{t})\int_{B(0,\delta)}x^{i}(\eta \overline{B}_{t})^{2^{*}}dx\\
&&+\frac{1}{2}\partial _{kl}\overline{f}_{t}(x_{t})\int_{B(0,\delta)}x^{k}x^{l}(\eta \overline{B}_{t})^{2^{*}}dx\\
&&+C\int_{B(0,\delta )}\left| x\right| ^{3}(\eta \overline{B}_{t})^{2^{*}}dx
\end{eqnarray*}
Now, first $g_{t}(x_{t})\leq0$, and second, and this is the main point for which we need the theorem of Druet and Robert (see the reason 
at the beginning of the next section), as $\overline{B}_{t}$ is radial, we have
$$\partial _{i}\overline{f}_{t}(x_{t})\int_{B(0,\delta)}x^{i}(\eta \overline{B}_{t})^{2^{*}}dx=0.$$
Therefore, introducing all this in the last inequality for $A^1_{t}$, we have
\[
\stackunder{t\rightarrow 1}{\overline{\lim }}\frac{A_{t}^{1}}{%
\int_{B(0,\delta )}\overline{u}_{t}^{2}dv_{\mathbf{g}_{t}}}\leq 
\]
$$\frac{n-2}{n}(1-\varepsilon _{\delta })\stackunder{t\rightarrow 1}{\overline{\lim }}
\frac{\frac{1}{2}\partial _{kl}\overline{f}_{t}(x_{t})\int_{B(0,\delta
)}x^{k}x^{l}(\eta \overline{B}_{t})^{2^{*}}dv_{%
\mathbf{g}_{t}}+C\int_{B(0,\delta )}\left| x\right| ^{3}(\eta 
\overline{B}_{t})^{2^{*}}dv_{\mathbf{g}_{t}}}{\int_{B(0,\delta )}\overline{u}%
_{t}^{2}dv_{\mathbf{g}_{t}}}$$
where we have replaced $dx$ by $dv_{\mathbf{g}_{t}}$ using the remark made at the beginning of the study of $A^2_{t}$.

Now, as for $A_{t}^{2}$, we write:
\begin{eqnarray*}
\stackunder{t\rightarrow 1}{\overline{\lim }}\frac{\int_{B(0,\delta
)}x^{k}x^{l}(\eta \overline{B}_{t})^{2^{*}}dv_{\mathbf{g}_{t}}}{%
\int_{B(0,\delta )}\overline{u}_{t}^{2}dv_{\mathbf{g}_{t}}} &=&f(x_{0})^{%
\frac{-2}{n}}K(n,2)^{2}\frac{n-4}{4(n-1)}\text{ if }k=l \\
&=&0\text{ if }k\neq l
\end{eqnarray*}
and therefore
\[
\frac{1}{K(n,2)^{2}f(x_{0})^{\frac{n-2}{n}}}\frac{n-2}{n}\stackunder{%
t\rightarrow 1}{\overline{\lim }}\frac{\frac{1}{2}\partial _{kl}\overline{f}_{t}(x_{t})\int_{B(0,\delta
)}x^{k}x^{l}(\eta \overline{B}_{t})^{2^{*}}dv_{\mathbf{g}_{t}}}{%
\int_{B(0,\delta )}\overline{u}_{t}^{2}dv_{\mathbf{g}_{t}}}= 
\]
\[
=\frac{1}{f(x_{0})}\frac{(n-2)(n-4)}{4(n-1)}\sum_{l}\frac{1}{2}\partial_{ll}\overline{f}_{1}(0)
\]
\[
=-\frac{(n-2)(n-4)}{8(n-1)}\frac{\bigtriangleup _{\mathbf{g}}f(x_{0})}{f(x_{0})}
\]
as $\bigtriangleup _{\mathbf{g}}f(x_{0})=-\sum_{l}\partial _{ll}\overline{f}%
_{1}(0)$ in the exponential chart in $x_{0}$.
Also
\[
\frac{\int_{B(0,\delta )}\left| x\right| ^{3}\overline{u}_{t}^{2^{*}}dv_{\,%
\mathbf{g}\,_{t}}}{\int_{B(0,\delta )}\overline{u}_{t}^{2}dv_{\,\mathbf{g}%
\,_{t}}}\leq C.\varepsilon _{\delta }\,\,. 
\]
Thus, we have proved that dividing inequality (18) by $\int_{B(0,\delta )}\overline{u}_{t}^{2}dv_{\mathbf{g}_{t}}$ 
and letting $t$ go to 1, we get
\[
h(x_{0})+\varepsilon _{\delta }\leq \frac{n-2}{4(n-1)}S_{\mathbf{g}}(x_{0})-%
\frac{(n-2)(n-4)}{8(n-1)}\frac{\bigtriangleup _{\mathbf{g}}f(x_{0})}{f(x_{0})%
}+\varepsilon _{\delta }\,\,. 
\]
Letting $\delta $ tend to 0: 
\[
h(x_{0})\leq \frac{n-2}{4(n-1)}S_{\mathbf{g}}(x_{0})-\frac{(n-2)(n-4)}{8(n-1)%
}\frac{\bigtriangleup _{\mathbf{g}}f(x_{0})}{f(x_{0})} 
\]
which contradict our hypothesis: 
\[
h(x_{0})>\frac{n-2}{4(n-1)}S_{\mathbf{g}}(x_{0})-\frac{(n-2)(n-4)}{8(n-1)}%
\frac{\bigtriangleup _{\mathbf{g}}f(x_{0})}{f(x_{0})} 
\]
when $x_{0}$ is a point of maximum of $f$. This prove that $u\not\equiv 0$, and therefore $u_{t}\rightarrow u>0$, a minimizing solution for $\E$,
and thus the weakly critical function $h$ is in fact critical.

\subsection{Alternate proof, proof of the fundamental estimate}
As we saw in the last part of the proof, the difficulty introduced by the presence of the function $f$ is to control 
the first derivatives of $f$, $\partial _{i}f(x_{t})$, as blow-up gives 
\[
\int_{B(0,\delta )}\partial _{i}f(x_{t})x^{i}\overline{u}_{t}^{2^{*}}dv_{%
\mathbf{g}_{t}}=\mu _{t}\int_{B(0,\delta \mu _{t}^{-1})}\partial
_{i}f(x_{t})x^{i}\widetilde{u}_{t}^{2^{*}}dv_{\widetilde{\mathbf{g}}_{t}} 
\]
to be divided by
\[
\mu _{t}^{2}\int_{B(0,\delta \mu _{t}^{-1})}\widetilde{u}_{t}^{2}dv_{%
\widetilde{\mathbf{g}}_{t}}\,\,, 
\]
and it would be necessary to control $\frac{\partial
_{i}f(x_{t})}{\mu _{t}}$, which seems to be difficult. But thanks to the theorem of Druet and Robert, we can replace 
$u(t)$ by $B(t)$ near $x_{t}$, and after blow-up 
$$\mu _{t}\int_{B(0,\delta \mu _{t}^{-1})}\partial
_{i}f(x_{t})x^{i}\widetilde{B}_{t}^{2^{*}}dv_{\widetilde{\mathbf{g}}_{t}} =0$$
as $\widetilde{B}_{t}$ is radial. Of course, the proof is then short, but the proof of the theorem of Druet and Robert 
is quite involved, even though the strong estimates (proposition 8) is the first step. 

The other way to get over the problem of the first derivatives of $f$ is to expand $f$ in $x_{0}$ as then 
 $\partial _{i}f(x_{0})=0$ because $x_{0}$ is a point of maximum of $f$. 
But then, one has to transpose the weak and strong estimates from $x_{t}$ to $x_{0}$,
which, as we said in the section about concentration phenomenom, requires to prove the following estimate: 
\[
\frac{d_{g}(x_{t},x_{0})}{\mu _{t}}\leq C\,\,. 
\]
As we said, this estimate is important and of independent interest, as it gives a complete description of the sequence 
$(u_{t})$. This is why we give this alternate proof of theorem 1, even though it requires an additional hypothesis.
This proof, which gives at the same time the proof of theorem 1 and of the estimate, is, we think, interesting, and is available 
directly after proposition 9, i.e it does not require the theorem of Druet and Robert.

We now make the hypothesis that the hessian of $f$ is nondegenerate at its points of maximum. 
We also suppose now that $dim M \geq 5$, even though our proof gives theorem 1 in 
dimension 4.

Let us note $x_{0}(t)=\exp _{x_{t}}^{-1}(x_{0})=(x_{0}^{1}(t),...,x_{0}^{n}(t))$, 
which is possible as soon as $t$ is close enough to 1 for a fixed radius $\delta$. Then $x_{0}(t)\rightarrow 0$
 when $t\rightarrow 1$. The point $x_{0}(t)$ is a locally strict maximum of $\overline{f}_{t}$. 
 We will let $\delta$ go to 0 at the end of the reasoning, after having taken the limit when $t\rightarrow 1$. 

The expansion of $\overline{f}_{t}$ in $x_{0}(t)$ gives: 
\[
\overline{f}_{t}(x)\leq f(x_{0})+\frac{1}{2}\partial _{kl}\overline{f}%
_{t}(x_{0}(t)).(x^{k}-x_{0}^{k}(t))(x^{l}-x_{0}^{l}(t))+c\left|
x-x_{0}(t)\right| ^{3}:=f(x_{0})+T_{t} 
\]
($T_{t}$ like Taylor) where ($\partial _{kl}\overline{f}_{t}(x_{0})$) is a negative definite matrix 
(we shall write $<0$). $c,C$ will always be constants independent of $t$ and $\delta $.
Remember that
\[
A_{t}^{1}=(\int_{B(0,\delta )}\overline{f}_{t}\overline{u}_{t}^{2^{*}}dv_{%
\mathbf{g}_{t}})^{\frac{2}{n}}(\int_{B(0,\delta )}\overline{f}_{t}(\eta 
\overline{u}_{t})^{2^{*}}dx)^{\frac{n-2}{n}}\,-(Supf.\int_{B(0,\delta
)}(\eta \overline{u}_{t})^{2^{*}}dx)^{\frac{n-2}{n}}\,\,. 
\]
Introducing the expansion of $\overline{f}_{t}$ in $x_{0}(t)$, 
and using again the fact that $(1+x)^\alpha \leq 1+\alpha x$ for $0<\alpha \leq1$, we get: 
\[
(\int_{B(0,\delta )}\overline{f}_{t}(\eta \overline{u}_{t})^{2^{*}}dx)^{%
\frac{n-2}{n}}\leq (\int_{B(0,\delta )}f(x_{0})(\eta \overline{u}%
_{t})^{2^{*}}dx)^{\frac{n-2}{n}}+\frac{\frac{n-2}{n}}{(\int_{B(0,\delta
)}f(x_{0})(\eta \overline{u}_{t})^{2^{*}}dx)^{\frac{2}{n}}}%
\{F_{t}\} 
\]
where
\[
\{F_{t}\}=\frac{1}{2}\partial _{kl}\overline{f}_{t}(x_{0}(t))\int_{B(0,%
\delta )}(x^{k}-x_{0}^{k}(t))(x^{l}-x_{0}^{l}(t))(\eta \overline{u}%
_{t})^{2^{*}}dx+C\int_{B(0,\delta )}\left| x-x_{0}(t)\right| ^{3}(\eta 
\overline{u}_{t})^{2^{*}}dx 
\]
from where, remembering that $\stackunder{M}{Sup}f=f(x_{0})$ and that $%
\int_{B(0,\delta )}\overline{f}_{t}\overline{u}_{t}^{2^{*}}dv_{\mathbf{g}%
_{t}}\leq 1$: 
\begin{equation}
A_{t}^{1}\leq \frac{n-2}{n}\frac{(\int_{B(0,\delta )}\overline{f}_{t}%
\overline{u}_{t}^{2^{*}}dv_{\mathbf{g}_{t}})^{\frac{2}{n}}}{%
(\int_{B(0,\delta )}f(x_{0})(\eta \overline{u}_{t})^{2^{*}}dx)^{\frac{2}{n}}}%
\{F_{t}\}
\end{equation}
Therefore, we obtain:
\[
\stackunder{t\rightarrow 1}{\overline{\lim }}\frac{A_{t}^{1}}{%
\int_{B(0,\delta )}\overline{u}_{t}^{2}dv_{\mathbf{g}_{t}}}\leq 
\]
$\frac{n-2}{n}(1+\varepsilon _{\delta })\stackunder{t\rightarrow 1}{\overline{\lim }}
\frac{\frac{1}{2}\partial _{kl}\overline{f}_{t}(x_{0})\int_{B(0,\delta
)}(x^{k}-x_{0}^{k}(t))(x^{l}-x_{0}^{l}(t))(\eta \overline{u}_{t})^{2^{*}}dv_{%
\mathbf{g}_{t}}+C\int_{B(0,\delta )}\left| x-x_{0}(t)\right| ^{3}(\eta 
\overline{u}_{t})^{2^{*}}dv_{\mathbf{g}_{t}}}{\int_{B(0,\delta )}\overline{u}%
_{t}^{2}dv_{\mathbf{g}_{t}}}$
\\

where we write $\partial _{kl}\overline{f}_{t}(x_{0})$ for $\partial
_{kl}\overline{f}_{t}(x_{0}(t)).$ Considering the expansion
\[
\overline{f}_{t}(x)\leq f(x_{0})+\frac{1}{2}\partial _{kl}\overline{f}%
_{t}(x_{0}(t)).(x^{k}-x_{0}^{k}(t))(x^{l}-x_{0}^{l}(t))+c\left|
x-x_{0}(t)\right| ^{3},
\]
note that by the regularity of $\exp _{x_{t}}^{-1}\circ \exp _{x_{0}}$ with respect to all the variables, 
we can suppose that $c$ is independent of $t$. Moreover: 
\begin{eqnarray*}
c\left| x-x_{0}(t)\right| ^{3} & &\leq c^{\prime }\left| x-x_{0}(t)\right|
\sum_{k}(x^{k}-x_{0}^{k}(t))^{2} \\ 
& &\leq 2\delta c^{\prime }\sum_{k}(x^{k}-x_{0}^{k}(t))^{2}
\end{eqnarray*}
where we remind that $\delta $ is the radius of the ball of integration. 
We can then write:
\[
\overline{f}_{t}(x)\leq f(x_{0})+(\frac{1}{2}\partial _{kl}\overline{f}%
_{t}(x_{0}(t))+\delta C_{kl})(x^{k}-x_{0}^{k}(t))(x^{l}-x_{0}^{l}(t)) 
\]
where $C_{kl}=c\delta _{kl}=c$ if $k=l$ and $C_{kl}=0$ if $k\neq l$ ($\delta
_{kl}$ is the Kr\"{o}necker symbol) is independent of $t$.

We introduce one more notation:
\[
D_{kl}(t,\delta )=\frac{1}{2}\partial _{kl}\overline{f}_{t}(x_{0}(t))+\delta
C_{kl}\,\,. 
\]
Then:

1/: $\stackunder{\delta \rightarrow 0}{\lim }\stackunder{t\rightarrow 1}{%
\lim }D_{kl}(t,\delta )=\frac{1}{2}\partial _{kl}\overline{f}_{1}(x_{0}(1))$
where $\overline{f}_{1}=f\circ \exp _{x_{0}}^{-1}$ and \\
$x_{0}(1)=0=\exp_{x_{0}}^{-1}(x_{0})$.

2/: for any $\delta $ small enough and for all $t$ close to 1, $D_{kl}(t,\delta)$ is still negative definite.

$D_{kl}(t,\delta )$ is the hessian of $f$ in $x_{0}(t)$ perturbated on its diagonal by the third order terms. 
It is for the second point that we need the hypothesis that the hessian of $f$ is non degenerate. Thus
\[
\frac{1}{2}\partial _{kl}\overline{f}_{t}(x_{0})\int_{B(0,\delta
)}(x^{k}-x_{0}^{k}(t))(x^{l}-x_{0}^{l}(t))(\eta \overline{u}_{t})^{2^{*}}dv_{%
\mathbf{g}_{t}}+C\int_{B(0,\delta )}\left| x-x_{0}(t)\right| ^{3}(\eta 
\overline{u}_{t})^{2^{*}}dv_{\mathbf{g}_{t}}\leq 
\]
\[
D_{kl}(t,\delta )\int_{B(0,\delta
)}(x^{k}-x_{0}^{k}(t))(x^{l}-x_{0}^{l}(t))(\eta \overline{u}_{t})^{2^{*}}dv_{%
\mathbf{g}_{t}}\,\,. 
\]
Let
\[
\{F_{t}^{\prime }\}=D_{kl}(t,\delta )\int_{B(0,\delta
)}(x^{k}-x_{0}^{k}(t))(x^{l}-x_{0}^{l}(t))(\eta \overline{u}_{t})^{2^{*}}dv_{%
\mathbf{g}_{t}} 
\]
We have
\[
\stackunder{t\rightarrow 1}{\overline{\lim }}\frac{A_{t}^{1}}{%
\int_{B(0,\delta )}v_{t}^{2}dv_{\mathbf{g}_{t}}}\leq \frac{n-2}{n}%
\stackunder{t\rightarrow 1}{\overline{\lim }}\frac{D_{kl}(t,\delta
)\int_{B(0,\delta )}(x^{k}-x_{0}^{k}(t))(x^{l}-x_{0}^{l}(t))(\eta \overline{u%
}_{t})^{2^{*}}dv_{\mathbf{g}_{t}}}{\int_{B(0,\delta )}\overline{u}%
_{t}^{2}dv_{\mathbf{g}_{t}}}(1+\varepsilon _{\delta }) 
\]

In the expansion of $D_{kl}(t,\delta
)(x^{k}-x_{0}^{k}(t))(x^{l}-x_{0}^{l}(t))$, we are interested by the first term, 
i.e $D_{kl}(t,\delta )x^{k}x^{l}$ (look back how we obtained $%
S_{g}(x_{0})$ in $A_{t}^{2}$), and we are going to show that the other terms can be neglected. 
The idea is to reorganize the expansion of $%
\{F_{t}^{\prime }\}$ and to use the fact $D_{kl}(t,\delta )$ is a negative bilinear form: 
\begin{eqnarray*}
\{F_{t}^{\prime }\}=&&D_{kl}(t,\delta )\int_{B(0,\delta
)}x^{k}x^{l}(\eta \overline{u}_{t})^{2^{*}}dv_{\mathbf{g}_{t}}+D_{kl}(t,%
\delta )x_{0}^{k}(t)x_{0}^{l}(t)\int_{B(0,\delta )}(\eta \overline{u}%
_{t})^{2^{*}}dv_{\mathbf{g}_{t}} \\ 
& &-D_{kl}(t,\delta )\int_{B(0,\delta
)}(x^{k}x_{0}^{l}(t)+x^{l}x_{0}^{k}(t))(\eta \overline{u}_{t})^{2^{*}}dv_{%
\mathbf{g}_{t}}\,\,.
\end{eqnarray*}
We rewrite the two last terms (suppressing some $\delta $ et $%
t $ and all integral being taken with respect to $dv_{\mathbf{g}_{t}}$):

\[
D_{kl}.x_{0}^{k}x_{0}^{l}\int_{B(0,\delta )}(\eta \overline{u}%
_{t})^{2^{*}}dv_{\mathbf{g}_{t}}-D_{kl}\int_{B(0,\delta
)}(x^{k}x_{0}^{l}+x^{l}x_{0}^{k})(\eta \overline{u}_{t})^{2^{*}}dv_{\mathbf{g%
}_{t}}= 
\]

\[
D_{kl}\Big[ x_{0}^{k}x_{0}^{l}\int_{B(0,\delta )}(\eta \overline{u}%
_{t})^{2^{*}}-x_{0}^{l}\int_{B(0,\delta )}x^{k}(\eta \overline{u}%
_{t})^{2^{*}}-x_{0}^{k}\int_{B(0,\delta )}x^{l}(\eta \overline{u}%
_{t})^{2^{*}}\Big]= 
\]
\[
D_{kl}\Big[x_{0}^{k}(\int_{B(0,\delta )}(\eta \overline{u}_{t})^{2^{*}})^{\frac{1%
}{2}}.x_{0}^{l}(\int_{B(0,\delta )}(\eta \overline{u}_{t})^{2^{*}})^{\frac{1%
}{2}} 
\]
\[
-x_{0}^{l}(\int_{B(0,\delta )}(\eta \overline{u}_{t})^{2^{*}})^{\frac{1}{2}}%
\frac{\int_{B(0,\delta )}x^{k}(\eta \overline{u}_{t})^{2^{*}}}{%
(\int_{B(0,\delta )}(\eta \overline{u}_{t})^{2^{*}})^{\frac{1}{2}}}%
-x_{0}^{k}(\int_{B(0,\delta )}(\eta \overline{u}_{t})^{2^{*}})^{\frac{1}{2}}%
\frac{\int_{B(0,\delta )}x^{l}(\eta \overline{u}_{t})^{2^{*}}}{%
(\int_{B(0,\delta )}(\eta \overline{u}_{t})^{2^{*}})^{\frac{1}{2}}}\Big]\,\,. 
\]
Thus, setting (sorry): 
\[
\varepsilon ^{k}(t)=\int_{B(0,\delta )}x^{k}(\eta \overline{u}%
_{t})^{2^{*}}dv_{\mathbf{g}_{t}} 
\]
\[
z_{t}=(\int_{B(0,\delta )}(\eta \overline{u}_{t})^{2^{*}}dv_{\mathbf{g}%
_{t}})^{\frac{1}{2}} 
\]
the expression above becomes:
\begin{center}
\[
D_{kl}.x_{0}^{k}x_{0}^{l}\int_{B(0,\delta )}(\eta \overline{u}%
_{t})^{2^{*}}dv_{\mathbf{g}_{t}}-D_{kl}\int_{B(0,\delta
)}(x^{k}x_{0}^{l}+x^{l}x_{0}^{k})(\eta \overline{u}_{t})^{2^{*}}dv_{\mathbf{g%
}_{t}}= 
\]
\[
=D_{kl}\Big[x_{0}^{k}(t).z_{t}.x_{0}^{l}(t).z_{t}-x_{0}^{l}(t).z_{t}.\frac{%
\varepsilon ^{k}(t)}{z_{t}}-x_{0}^{k}(t).z_{t}.\frac{\varepsilon ^{l}(t)}{z_{t}}\Big] 
\]
\end{center}
\[
=D_{kl}\Big[(x_{0}^{k}(t).z_{t}-\frac{\varepsilon ^{k}(t)}{z_{t}}%
)(x_{0}^{l}(t).z_{t}-\frac{\varepsilon ^{l}(t)}{z_{t}})-\frac{\varepsilon
^{k}(t)\varepsilon ^{l}(t)}{z_{t}^{2}}\Big] 
\]
By this method of reorganization of the hessian, we have obtained: 
\[
\frac{1}{2}\partial _{kl}\overline{f}_{t}(x_{0})\int_{B(0,\delta
)}(x^{k}-x_{0}^{k}(t))(x^{l}-x_{0}^{l}(t))(\eta \overline{u}_{t})^{2^{*}}dv_{%
\mathbf{g}_{t}}+C\int_{B(0,\delta )}\left| x-x_{0}(t)\right| ^{3}(\eta 
\overline{u}_{t})^{2^{*}}dv_{\mathbf{g}_{t}}\leq 
\]
\begin{center}
\[
D_{kl}(t,\delta )\int_{B(0,\delta )}x^{k}x^{l}(\eta \overline{u}%
_{t})^{2^{*}}dv_{\mathbf{g}_{t}}+D_{kl}(t,\delta )(x_{0}^{k}(t).z_{t}-\frac{%
\varepsilon ^{k}(t)}{z_{t}})(x_{0}^{l}(t).z_{t}-\frac{\varepsilon ^{l}(t)}{%
z_{t}})-D_{kl}(t,\delta )\frac{\varepsilon ^{k}(t)\varepsilon ^{l}(t)}{%
z_{t}^{2}} 
\]
\end{center}
\[
\leq D_{kl}(t,\delta )\int_{B(0,\delta )}x^{k}x^{l}(\eta \overline{u}%
_{t})^{2^{*}}dv_{\mathbf{g}_{t}}-D_{kl}(t,\delta )\frac{\varepsilon
^{k}(t)\varepsilon ^{l}(t)}{z_{t}^{2}} 
\]
because, and that is the fundamental point : 
\[
D_{kl}(t,\delta )\omega ^{k}\omega ^{l}\leq 0\text{ \thinspace }\forall
\omega =(\omega ^{1},...,\omega ^{n}) 
\]
which allows to suppress from the inequality $$D_{kl}(t,\delta
)(x_{0}^{k}(t).z_{t}-\frac{\varepsilon ^{k}(t)}{z_{t}})(x_{0}^{l}(t).z_{t}-%
\frac{\varepsilon ^{l}(t)}{z_{t}})$$ It is this term that will give us the estimate $\frac{d_{\mathbf{g}}(x_{t},x_{0})}{\mu _{t}}\leq C$ 
(see below).

We have therefore obtained:
\[
\stackunder{t\rightarrow 1}{\overline{\lim }}\frac{A_{t}^{1}}{%
\int_{B(0,\delta )}\overline{u}_{t}^{2}dv_{\mathbf{g}_{t}}}\leq \frac{n-2}{n}%
\stackunder{t\rightarrow 1}{\overline{\lim }}\frac{D_{kl}(t,\delta
)\int_{B(0,\delta )}x^{k}x^{l}(\eta \overline{u}_{t})^{2^{*}}dv_{\mathbf{g}%
_{t}}-D_{kl}(t,\delta )\frac{\varepsilon ^{k}(t)\varepsilon ^{l}(t)}{%
z_{t}^{2}}}{\int_{B(0,\delta )}\overline{u}_{t}^{2}dv_{\mathbf{g}_{t}}}%
(1+\varepsilon _{\delta }) 
\]
Now, as for $A_{t}^{2}$, we write:
\begin{eqnarray*}
\stackunder{t\rightarrow 1}{\overline{\lim }}\frac{\int_{B(0,\delta
)}x^{k}x^{l}(\eta \overline{u}_{t})^{2^{*}}dv_{\mathbf{g}_{t}}}{%
\int_{B(0,\delta )}\overline{u}_{t}^{2}dv_{\mathbf{g}_{t}}} &=&f(x_{0})^{%
\frac{-2}{n}}K(n,2)^{2}\frac{n-4}{4(n-1)}\text{ if }k=l \\
&=&0\text{ if }k\neq l
\end{eqnarray*}
and therefore
\[
\frac{1}{K(n,2)^{2}f(x_{0})^{\frac{n-2}{n}}}\frac{n-2}{n}\stackunder{%
t\rightarrow 1}{\overline{\lim }}\frac{D_{kl}(t,\delta )\int_{B(0,\delta
)}x^{k}x^{l}(\eta \overline{u}_{t})^{2^{*}}dv_{\mathbf{g}_{t}}}{%
\int_{B(0,\delta )}\overline{u}_{t}^{2}dv_{\mathbf{g}_{t}}}= 
\]
\[
=\frac{1}{f(x_{0})}\frac{(n-2)(n-4)}{4(n-1)}\sum_{l}(\frac{1}{2}\partial
_{ll}\overline{f}_{1}(0)+c_{ll}\delta ) 
\]
\[
=-\frac{(n-2)(n-4)}{8(n-1)}\frac{\bigtriangleup _{\mathbf{g}}f(x_{0})}{%
f(x_{0})}+\varepsilon _{\delta } 
\]
as $\bigtriangleup _{\mathbf{g}}f(x_{0})=-\sum_{l}\partial _{ll}\overline{f}%
_{1}(0)$ in the exponential chart in $x_{0}$.

At last, let us show that the residual term can be neglected. 
\[
\left| \varepsilon ^{k}(t)\varepsilon ^{l}(t)\right| \leq \frac{1}{2}%
(\varepsilon ^{k}(t)^{2}+\varepsilon ^{l}(t)^{2}) 
\]
But 
\[
\varepsilon ^{k}(t)^{2}=(\int_{B(0,\delta )}x^{k}(\eta \overline{u}%
_{t})^{2^{*}}dv_{\mathbf{g}_{t}})^{2} 
\]
\[
=(\int_{B(0,R\mu _{t})}x^{k}(\eta \overline{u}_{t})^{2^{*}}dv_{\mathbf{g}%
_{t}}+\int_{B(0,\delta )\backslash B(0,R\mu _{t})}x^{k}(\eta \overline{u}%
_{t})^{2^{*}}dv_{\mathbf{g}_{t}})^{2} 
\]
\[
\leq 2(\int_{B(0,R\mu _{t})}x^{k}(\eta \overline{u}_{t})^{2^{*}}dv_{\mathbf{g%
}_{t}})^{2}+2(\int_{B(0,\delta )\backslash B(0,R\mu _{t})}x^{k}(\eta 
\overline{u}_{t})^{2^{*}}dv_{\mathbf{g}_{t}})^{2} 
\]
The blow-up formula's give, for a fixed $R$ : 
\[
\frac{(\int_{B(0,R\mu _{t})}x^{k}(\eta \overline{u}_{t})^{2^{*}}dv_{\mathbf{g%
}_{t}})^{2}}{\int_{B(0,\delta )}\overline{u}_{t}^{2}dv_{\mathbf{g}_{t}}}\leq 
\frac{(\mu _{t}\int_{B(0,R)}x^{k}\widetilde{u}_{t}^{2^{*}}dv_{\widetilde{%
\mathbf{g}}_{t}})^{2}}{\mu _{t}^{2}\int_{B(0,R)}\widetilde{u}_{t}^{2}dv_{%
\widetilde{\mathbf{g}}_{t}}.}\stackunder{t\rightarrow 1}{\rightarrow }\frac{%
(\int_{B(0,R)}x^{k}\widetilde{u}^{2^{*}}dx)^{2}}{\int_{B(0,R)}\widetilde{u}%
^{2}dx}=0 
\]
because $\widetilde{u}$ is radial.

At last, using the weak estimates: $\,d_{\mathbf{g}}(x,x_{t})^{\frac{%
n-2}{2}}u_{t}(x)\leq \varepsilon $ if $d_{\mathbf{g}%
}(x,x_{t})\geq R\mu _{t}$, and using the Hölder's inequality: 
\begin{eqnarray*}
(\int_{B(0,\delta )\backslash B(0,R\mu _{t})}x^{k}(\eta \overline{u}%
_{t})^{2^{*}}dv_{\mathbf{g}_{t}})^{2} && \leq \varepsilon
_{R}^{2}(\int_{B(0,\delta )\backslash B(0,R\mu _{t})}\overline{u}_{t}^{2%
\frac{n-1}{n-2}}dv_{\mathbf{g}_{t}})^{2} \\ 
&& \leq \varepsilon _{R}^{2}(\int_{B(0,\delta )\backslash B(0,R\mu _{t})}%
\overline{u}_{t}^{2}dv_{\mathbf{g}_{t}})(\int_{B(0,\delta )\backslash
B(0,R\mu _{t})}\overline{u}_{t}^{\frac{2n}{n-2}}dv_{\mathbf{g}_{t}})
\end{eqnarray*}
therefore 
\[
\frac{(\int_{B(0,\delta )\backslash B(0,R\mu _{t})}x^{k}(\eta \overline{u}%
_{t})^{2^{*}}dv_{\mathbf{g}_{t}})^{2}}{\int_{B(0,\delta )}\overline{u}%
_{t}^{2}dv_{\mathbf{g}_{t}}}\leq \varepsilon _{R}^{2}(\int_{B(0,\delta
)\backslash B(0,R\mu _{t})}\overline{u}_{t}^{\frac{2n}{n-2}}dv_{\mathbf{g}%
_{t}})\leq c\varepsilon _{R}^{2} 
\]
where $\varepsilon _{R}\rightarrow 0$ when $R\rightarrow \infty $. Remarking that because $x_{0}$ 
is a concentration point: 
\[
z_{t}^{2}=\int_{B(0,\delta )}(\eta \overline{u}_{t})^{2^{*}}dv_{\mathbf{g}%
_{t}}\geq \int_{B(x_{0},\delta /4)}u_{t}^{2^{*}}dv_{\mathbf{g}}\geq c>0\text{
} 
\]
we have obtained
\[
\frac{\left| \varepsilon ^{k}(t)\varepsilon ^{l}(t)\right| }{%
z_{t}^{2}\int_{B(0,\delta )}\overline{u}_{t}^{2}dv_{\mathbf{g}_{t}}}%
\stackunder{t\rightarrow 1}{\rightarrow }0 
\]
We have therefore obtained once again that
\[
h(x_{0})+\varepsilon _{\delta }\leq \frac{n-2}{4(n-1)}S_{\mathbf{g}}(x_{0})-%
\frac{(n-2)(n-4)}{8(n-1)}\frac{\bigtriangleup _{\mathbf{g}}f(x_{0})}{f(x_{0})%
}+\varepsilon _{\delta }\,\,. 
\]
Letting $\delta $ tend to 0: 
\[
h(x_{0})\leq \frac{n-2}{4(n-1)}S_{\mathbf{g}}(x_{0})-\frac{(n-2)(n-4)}{8(n-1)%
}\frac{\bigtriangleup _{\mathbf{g}}f(x_{0})}{f(x_{0})} 
\]
which contradict our hypothesis: 
\[
h(x_{0})>\frac{n-2}{4(n-1)}S_{\mathbf{g}}(x_{0})-\frac{(n-2)(n-4)}{8(n-1)}%
\frac{\bigtriangleup _{\mathbf{g}}f(x_{0})}{f(x_{0})} 
\]
when $x_{0}$ is a point of maximum of $f$. This prove that $u_{t}\rightarrow u>0$, a minimizing solution for $\E$,
and therefore the weakly critical function $h$ is in fact critical.

We now prove the estimate
\[
\frac{d_{\mathbf{g}}(x_{t},x_{0})}{\mu _{t}}\leq C 
\]

Going back to the computations above, we have obtained:
\begin{eqnarray}
h(x_{0}) &\leq &\frac{n-2}{4(n-1)}S_{\mathbf{g}}(x_{0})-\frac{(n-2)(n-4)}{%
8(n-1)}\frac{\bigtriangleup _{\mathbf{g}}f(x_{0})}{f(x_{0})}+\varepsilon
_{\delta }  \\
&&+\stackunder{t\rightarrow 1}{\overline{\lim }}\frac{n-2}{n}\frac{%
D_{kl}(t,\delta )(x_{0}^{k}(t).z_{t}-\frac{\varepsilon ^{k}(t)}{z_{t}}%
)(x_{0}^{l}(t).z_{t}-\frac{\varepsilon ^{l}(t)}{z_{t}})}{\int_{B(0,\delta )}%
\overline{u}_{t}^{2}dv_{\mathbf{g}_{t}}}  \nonumber
\end{eqnarray}
where $D_{kl}(t,\delta )$ is negative definite for $t$ close to 1
and for all $\delta $ small enough, and where we remind that $%
x_{0}(t)=\exp _{x_{t}}^{-1}(x_{0})=(x_{0}^{1}(t),...,x_{0}^{n}(t))$. So, there exists a $\lambda >0$ such that 
$\forall \omega \in \Bbb{R}^{n}:$%
\[
D_{kl}(t,\delta )\omega ^{k}\omega ^{l}\leq -\lambda \sum_{k}\left| \omega
^{k}\right| ^{2} 
\]
and so
\[
D_{kl}(t,\delta )\frac{(x_{0}^{k}(t).z_{t}-\frac{\varepsilon ^{k}(t)}{z_{t}}%
)(x_{0}^{l}(t).z_{t}-\frac{\varepsilon ^{l}(t)}{z_{t}})}{(\int_{B(0,\delta )}%
\overline{u}_{t}^{2}dv_{\mathbf{g}_{t}})^{\frac{1}{2}}(\int_{B(0,\delta )}%
\overline{u}_{t}^{2}dv_{\mathbf{g}_{t}})^{\frac{1}{2}}}\leq 
\]
\[
-\lambda \sum_{k}\left| \frac{x_{0}^{k}(t).z_{t}}{(\int_{B(0,\delta )}%
\overline{u}_{t}^{2}dv_{\mathbf{g}_{t}})^{\frac{1}{2}}}-\frac{\varepsilon
^{k}(t)}{z_{t}(\int_{B(0,\delta )}\overline{u}_{t}^{2}dv_{\mathbf{g}_{t}})^{%
\frac{1}{2}}}\right| ^{2}\,\,. 
\]
Moreover, we already proved that: 
\[
\frac{\varepsilon ^{k}(t)^{2}}{z_{t}^{2}\int_{B(0,\delta )}\overline{u}%
_{t}^{2}dv_{\mathbf{g}_{t}}}\stackunder{t\rightarrow 1}{\rightarrow }0 
\]
as we also have $z_{t}=(\int_{B(0,\delta )}\overline{u}_{t}^{2^{*}}dv_{%
\mathbf{g}_{t}})^{\frac{1}{2}}$, and therefore as $x_{0}$ is a concentration point: 
\[
0<c\leq \lim \inf \,z_{t}\leq \lim \sup \,z_{t}\leq c^{\prime }<+\infty . 
\]
Therefore, necessarilly, because of (20), for all $k$, there exists a constant 
$C>0$ such that for $t\rightarrow 1$ : 
\[
\,\frac{x_{0}^{k}(t)}{(\int_{B(0,\delta )}\overline{u}_{t}^{2}dv_{\mathbf{g}%
_{t}})^{\frac{1}{2}}}\leq C 
\]
Now
\[
\int_{B(0,\delta )}\overline{u}_{t}^{2}dv_{\mathbf{g}_{t}}=\mu
_{t}^{2}\int_{B(0,\delta \mu _{t}^{-1})}\widetilde{u}_{t}^{2}dv_{\widetilde{%
\mathbf{g}}_{t}}. 
\]
But the strong estimates give that 
\[
\stackunder{t\rightarrow 1}{\overline{\lim }}\int_{B(0,\delta \mu _{t}^{-1})}%
\widetilde{u}_{t}^{2}dv_{\widetilde{\mathbf{g}}_{t}}<+\infty 
\]
therefore
\[
\int_{B(0,\delta )}\overline{u}_{t}^{2}dv_{\mathbf{g}_{t}}\sim C\mu _{t}^{2} 
\]
from where we have
\[
\forall k:\,\frac{x_{0}^{k}(t)}{\mu _{t}}\leq C^{\prime } 
\]
and so 
\[
\frac{d_{\mathbf{g}}(x_{t},x_{0})}{\mu _{t}}\leq C 
\]

If we have furthermore that at the points of maximum of $f$ : 
\[
h(P)=\frac{n-2}{4(n-1)}S_{\mathbf{g}}(P)-\frac{(n-2)(n-4)}{8(n-1)}\frac{%
\bigtriangleup _{\mathbf{g}}f(P)}{f(P)} 
\]
then we have more precisely that 
\[
\frac{d_{\mathbf{g}}(x_{t},x_{0})}{\mu _{t}}\rightarrow 0 
\]

\textit{Remark:} Note that when concentration occurs we have: 
\[
h(x_{0})\leq \frac{n-2}{4(n-1)}S_{\mathbf{g}}(x_{0})-\frac{(n-2)(n-4)}{8(n-1)%
}\frac{\bigtriangleup _{\mathbf{g}}f(x_{0})}{f(x_{0})} 
\]

\section{Critical triple 1: existence of critical functions}
The idea to prove the existence of critical functions (theorem 2), is to find, being given the manifold $(M,\g)$ 
and the function $f$, a subcritical function $h_{0}$ and a weakly critical function $h_{1}$ and then to join these 
two functions by a continuous path; theorem 1 then shows that this path must "cross" the set of critical functions.

Note first that, by the sharp Sobolev inequality (2), $B_{0}(\mathbf{g})K(n,2)^{-2}$ is a weakly critical function 
for any manifold $(M,\g)$ and any function $f$. Also, it is known that 
$$B_{0}(\mathbf{g})K(n,2)^{-2}\geq \frac{n-2}{4(n-1)} \underset{M}{Sup}\, S_{\g}$$
Therefore, for any $\alpha >0$, and for any point $P$ where $f$ is maximum on $M$, we have
$$B_{0}(\mathbf{g})K(n,2)^{-2}+ \alpha \geq \frac{n-2}{4(n-1)}S_{\mathbf{g}}(P)-\frac{(n-2)(n-4)}{8(n-1)%
}\frac{\bigtriangleup _{\mathbf{g}}f(P)}{f(P)}.$$

Now, we are going to modify the weakly critical function $B_{0}(\mathbf{g})K(n,2)^{-2}+ \alpha$ by the test 
functions presented in the introduction. They can be seen under the following form: 
 for any $x\in M$ and any $\delta >0$
small enough, there exists a sequence of functions $(\psi _{k})$ with compact support in 
$B(x,\delta )$ such that for any function $h$: 
\[
J_{h,1,\mathbf{g}}(\psi _{k})=\frac{\int_{M}\left| \nabla \psi _{k}\right|
^{2}dv_{\mathbf{g}}+\int_{M}h.\psi _{k}{{}^{2}}dv_{\mathbf{g}}}{\left(
\int_{M}\left| \psi _{k}\right| ^{\frac{2n}{n-2}}dv_{\mathbf{g}}\right) ^{%
\frac{n-2}{n}}}\stackunder{k\rightarrow \infty }{\rightarrow }\frac{1}{K(n,2)%
{{}^{2}}} 
\]
and
\[
\int_{M}\psi _{k}^{\frac{2n}{n-2}}dv_{\mathbf{g}}=1 
\]
this last condition being obtained by multiplying the functions in the introduction by suitable constants. We will 
use the functional $J$ here, as 
\[
\int_{M}f.\psi _{k}^{\frac{2n}{n-2}}dv_{\mathbf{g}}\neq 1\,\,. 
\]
Let then $\psi _{k}$ be one of these functions, where $k$ and $B(x,\delta )$
will be fixed later. We consider, for $t>0$ the sequence
\[
h_{t}=B_{0}(\mathbf{g})K(n,2)^{-2}+\alpha -t.\psi _{k}^{\frac{4}{n-2}}\,\,. 
\]
First, we seek a condition for $\triangle _{\mathbf{g}%
}+h_{t} $ to be coercive. Noting $B_{0}K^{-2}=B_{0}(\mathbf{g}%
)K(n,2)^{-2}$, and taking all integrals for the measure $dv_{\mathbf{g}}$, we have for $u\in H_{1}^{2}$%
\begin{eqnarray*}
\int_{M}(\left| \nabla u\right| _{\mathbf{g}}^{2}+h_{t}u^{2}) && 
=\int_{M}(\left| \nabla u\right| _{\mathbf{g}}^{2}+B_{0}K^{-2}.u^{2})-(t-\alpha)%
\int_{M}\psi _{k}^{\frac{4}{n-2}}.u^{2} \\ 
&& \geq K(n,2)^{-2}(\int_{M}u^{\frac{2n}{n-2}})^{\frac{n-2}{n}}-(t-\alpha)\int_{M}\psi
_{k}^{\frac{4}{n-2}}.u^{2}
\end{eqnarray*}
by Sobolev inequality. But using Hölder's inequality:
\[
\int_{M}\psi _{k}^{\frac{4}{n-2}}.u^{2}\leq (\int_{M}\psi _{k}^{\frac{2n}{n-2%
}})^{\frac{n-2}{n}}(\int_{M}u^{\frac{2n}{n-2}})^{\frac{2}{n}}=(\int_{M}u^{%
\frac{2n}{n-2}})^{\frac{2}{n}} 
\]
as $\int_{M}\psi _{k}^{\frac{2n}{n-2}}=1$. Thus, using Hölder's inequality again to get the existence of a constant
 $C>0$ such that 
\[
C\int_{M}u^{2}\leq (\int_{M}u^{\frac{2n}{n-2}})^{\frac{n-2}{n}} 
\]
we have as soon as $K(n,2)^{-2}-(t-\alpha)>0$
\begin{eqnarray*}
\int_{M}(\left| \nabla u\right| _{\mathbf{g}}^{2}+h_{t}u^{2})& & \geq
(K(n,2)^{-2}-(t-\alpha))(\int_{M}u^{\frac{2n}{n-2}})^{\frac{n-2}{n}} \\ 
& &\geq (K(n,2)^{-2}-(t-\alpha))C\int_{M}u^{2}\,\,\,.
\end{eqnarray*}
So $\triangle _{\mathbf{g}}+h_{t}$ is coercive as soon as $t-\alpha <K(n,2)^{-2}$; 
we then fix $t_{1}$ such that $\alpha <t_{1}<K(n,2)^{-2}+ \alpha$. 

We now want to fix 
$\psi _{k}$ so that $h_{t_{1}}$ is subcritical for $f$. We pick first $x$ close enough to a point $x_{0}$ 
of maximum of $f$
and $\delta $ small enough such that $f>0$ on $B(x,\delta )$, to obtain : 
\begin{eqnarray*}
J_{h_{t_{1}},f,\mathbf{g}}(\psi _{k}) && =\frac{\int_{M}\left| \nabla \psi
_{k}\right| ^{2}+\int_{M}B_{0}K^{-2}.\psi _{k}{{}^{2}}-(t_{1}-\alpha)\int_{M}\psi
_{k}^{\frac{2n}{n-2}}}{\left( \int_{M}f\left| \psi _{k}\right| ^{\frac{2n}{%
n-2}}\right) ^{\frac{n-2}{n}}} \\ 
&& \leq \frac{J_{B_{0}K^{-2},1}(\psi _{k})}{(\stackunder{B(x,\delta )}{Inf}%
f)^{\frac{n-2}{n}}}-\frac{t_{1}-\alpha}{(\stackunder{B(x,\delta )}{Sup}f)^{\frac{n-2%
}{n}}} \\ 
&& \leq \frac{J_{B_{0}K^{-2},1}(\psi _{k})}{(\stackunder{B(x,\delta )}{Inf}%
f)^{\frac{n-2}{n}}}-\frac{t_{1}-\alpha}{(\stackunder{M}{Sup}f)^{\frac{n-2}{n}}}%
\,\,\,.
\end{eqnarray*}
For any $\varepsilon >0$, by continuity of $f$, we can choose $x$
close enough to a point of maximum $x_{0}$ and $\delta $ small enough such that $B(x,\delta
)\cap \left\{ x/f(x)=Maxf\right\} =\emptyset $ and
\[
\frac{1}{(\stackunder{B(x,\delta )}{Inf}f)^{\frac{n-2}{n}}}\leq \frac{1}{(%
\stackunder{M}{Sup}f)^{\frac{n-2}{n}}}+\varepsilon 
\]
$x$ and $\delta $ being fixed, we can now choose $k$ large enough to have
\[
J_{B_{0}K^{-2},1}(\psi _{k})\leq K(n,2)^{-2}+\varepsilon \,\,. 
\]
Therfore, choosing $\varepsilon $ small enough, we see that because$\frac{t_{1}-\alpha}{(\stackunder{M}{Sup}f)^{\frac{n-2}{n}}}>0$ : 
\[
J_{h_{t_{1}},f,\mathbf{g}}(\psi _{k})<\frac{1}{K(n,2)^{-2}(\stackunder{M}{Sup%
}f)^{\frac{n-2}{n}}} 
\]
and therefore $h_{t_{1}}$ is subcritical for $f$. We now set:
\[
t_{0}=Inf\{t\leq t_{1}/\lambda _{h_{t}}<\frac{1}{K(n,2){{}^{2}}(\stackunder{M%
}{Sup}f)^{\frac{n-2}{n}}}\}\text{ }. 
\]
Then $t_{0}\geq 0$, and
\[
\lambda _{h_{t_{0}}}=\frac{1}{K(n,2){{}^{2}}(\stackunder{M}{Sup}f)^{\frac{n-2%
}{n}}}\,\,\,and\,\,\,\lambda _{h_{t}}<\frac{1}{K(n,2){{}^{2}}(\stackunder{M}{%
Sup}f)^{\frac{n-2}{n}}}\,\,\,if\,\,\,t>t_{0}. 
\]
Furthermore $\forall t$, $t_{0}\leq t\leq t_{1}$,

\[
\frac{4(n-1)}{n-2}h_{t_{0}}(P)>S_{\mathbf{g}}(P)-\frac{n-4}{2}\frac{%
\bigtriangleup _{\mathbf{g}}f(P)}{f(P)}\, for \, P\in \left\{ x/f(x)=Maxf\right\} 
\]
because $B(x,\delta )\cap \left\{ x/f(x)=Maxf\right\} =\emptyset $. At last, $%
h_{t}\stackunder{t\rightarrow t_{0}}{\rightarrow }h_{t_{0}}$ in 
$C^{0,\alpha }$, and $\triangle _{\mathbf{g}}+h_{t_{0}}$ is coercive. Therefore by theorem 1, $h_{t_{0}}$ 
is critical and $\E$ has minimizing solutions.

Now, we prove that if $\left\{ x/f(x)=Maxf\right\}$ is thin and if $\int_{M}f>0$, there exist positive critical functions. 
We start again with $h=B_{0}(\mathbf{g})K(n,2)^{-2}+\alpha$, with $\alpha>0$.
For all $P$ where $f$ is maximum on $M$ : 
\[
\frac{4(n-1)}{n-2}B_{0}(\mathbf{g})K(n,2)^{-2}+\alpha>S_{\mathbf{g}}(P)-\frac{n-4}{2%
}\frac{\bigtriangleup _{\mathbf{g}}f(P)}{f(P)} 
\]
as
\[
B_{0}(\mathbf{g})K(n,2)^{-2}\geq \frac{(n-2)}{4(n-1)}MaxS_{\mathbf{g}}\,\,. 
\]
As $f$ is not constant, there exist $\eta $ with support in $M\backslash \left\{ x/f(x)=Maxf\right\} $ 
and such that $0\leq \eta \leq 1$. Let
\[
c=\left( \int_{M}fdv_{\mathbf{g}}\right) ^{-\frac{n-2}{n}} 
\]
That's where we need $\int_{M}fdv_{\mathbf{g}}>0$. We have $\int
fc^{2^{*}}dv_{\mathbf{g}}=1$.
For $t\in \Bbb{R}^{+}$ we set
\[
h_{t}=B_{0}K^{-2}+\alpha -t\eta \,\,. 
\]
Then $h_{t}=B_{0}K^{-2}+\alpha $ on $\left\{ x/f(x)=Maxf\right\} $, and
\[
I_{h_{t}}(c)=\int_{M}(B_{0}K^{-2}+\alpha)c{{}^{2}}dv_{\mathbf{g}}-c{{}^{2}}%
t\int_{M}\eta dv_{\mathbf{g}}=\left( \int_{M}fdv_{\mathbf{g}}\right) ^{-%
\frac{2}{2^{*}}}((B_{0}K^{-2}+\alpha){Vol}_{\mathbf{g}}(M)-t\int_{M}\eta dv_{\mathbf{g%
}}) 
\]
So, if $t$ is large enough, 
\[
I_{h_{t}}(c)<\frac{1}{K(n,2){{}^{2}}(\stackunder{M}{Sup}f)^{\frac{n-2}{n}}}%
\,\,. 
\]
We also want $h_{t}$ to be positive on $M$.
By the definition of $h_{t}$ and because $\stackunder{M}{Sup}\,\eta =1$, it is the case if 
\begin{equation}
t<B_{0}(\mathbf{g})K(n,2)^{-2}+\alpha.  
\end{equation}
But we also want that
\[
I_{h_{t}}(c)<\frac{1}{K(n,2){{}^{2}}(\stackunder{M}{Sup}f)^{\frac{n-2}{n}}} 
\]
which requires
\begin{equation}
t>\frac{1}{\int_{M}\eta dv_{\mathbf{g}}}\Bigl( (B_{0}K^{-2}+\alpha){Vol}_{\mathbf{g}}(M)-%
\frac{\left( \int_{M}fdv_{\mathbf{g}}\right) ^{\frac{n-2}{n}}}{K(n,2){{}^{2}}%
(\stackunder{M}{Sup}f)^{\frac{n-2}{n}}}\Bigr)\,\,. 
\end{equation}
We can find such a $t$ if : 
\[
\frac{1}{\int_{M}\eta dv_{\mathbf{g}}}\Bigl( (B_{0}K^{-2}+\alpha){Vol}_{\mathbf{g}}(M)-%
\frac{\left( \int_{M}fdv_{\mathbf{g}}\right) ^{\frac{n-2}{n}}}{K(n,2){{}^{2}}%
(\stackunder{M}{Sup}f)^{\frac{n-2}{n}}}\Bigr)<B_{0}K^{-2}+\alpha
\]
which can be writen 
\begin{equation}
\int_{M}\eta dv_{\mathbf{g}}>{Vol}_{\mathbf{g}}(M)-\frac{K^{-2}\left( \int_{M}fdv_{%
\mathbf{g}}\right) ^{\frac{n-2}{n}}}{(B_{0}K^{-2}+\alpha)(\stackunder{M}{Sup}%
f)^{\frac{n-2}{n}}}\,\,.  
\end{equation}
Remember that we want $\eta $ to have support in $M\backslash \left\{
x/f(x)=Maxf\right\} $ with $0\leq \eta \leq 1$. But we made the hypothesis that $\left\{ x/f(x)=Maxf\right\} $, 
the set of maximum points of $f$, is a thin set. We can therefore find such a function $\eta $ with
$\int_{M}\eta dv_{\mathbf{g}}$ as close as we want to ${Vol}_{\mathbf{g}}(M)$.
As
\[
\frac{\left( \int_{M}fdv_{\mathbf{g}}\right) ^{\frac{n-2}{n}}}{B_{0}(\mathbf{%
g})(\stackunder{M}{Sup}f)^{\frac{n-2}{n}}}>0 
\]
we can find $\eta $ satisfying (23) and a real $t$, denoted $t_{1}$, satisfying (21) and (22).

On the set $\left\{ x/f(x)=Maxf\right\} $, $h_{t}=B_{0}K^{-2}+\alpha$, so $\forall P\in
\left\{ x/f(x)=Maxf\right\} $:

\[
\frac{4(n-1)}{n-2}h_{t}(P)>S_{\mathbf{g}}(P)-\frac{n-4}{2}\frac{%
\bigtriangleup _{\mathbf{g}}f(P)}{f(P)}. 
\]
We then set: 
\[
t_{0}=Inf\{t\leq t_{1}\,/\,\lambda _{h_{t}}<\frac{1}{K(n,2){{}^{2}}(%
\stackunder{M}{Sup}f)^{\frac{n-2}{n}}}\}\text{ }. 
\]
Necessarilly, $t_{0}<t_{1}$. We remind that (see section 1):
\[
\lambda _{h,f,\mathbf{g}}=\lambda _{h}=\stackunder{w\in \mathcal{H}_{f}}{%
\inf }I_{h}(w) 
\]
Therefore
\[
\lambda _{h_{t_{0}}}=\frac{1}{K(n,2){{}^{2}}(\stackunder{M}{Sup}f)^{\frac{n-2%
}{n}}}\,\,\,and\,\,\,\lambda _{h_{t}}<\frac{1}{K(n,2){{}^{2}}(\stackunder{M}{%
Sup}f)^{\frac{n-2}{n}}}\,\,\,if\,\,\,t>t_{0}. 
\]
Furthermore $\forall t$, $t_{0}\leq t\leq t_{1}$, $h_{t}>0$ on $M$ and

\[
\frac{4(n-1)}{n-2}h_{t_{0}}(P)>S_{\mathbf{g}}(P)-\frac{n-4}{2}\frac{%
\bigtriangleup _{\mathbf{g}}f(P)}{f(P)}\,\,\,for\,\,\,P\in \left\{
x/f(x)=Maxf\right\} . 
\]
At last $h_{t}\stackunder{t\rightarrow t_{0}}{\rightarrow }h_{t_{0}}$ in $C^{0}$, and as 
$h_{t_{0}}>0$, $\triangle _{\mathbf{g}}+h_{t_{0}}$ is coercive. 
Therefore by theorem 1, $h_{t_{0}}$ is critical and $\E$ has minimizing solutions.

\textit{Remark:} The precceding proofs also show, by replacing $B_{0}K^{-2}$ by $h$, that if 
$(h,f,\mathbf{g})$ is weakly critical, and if
\[
\frac{4(n-1)}{n-2}h(P)>S_{\mathbf{g}}(P)-\frac{n-4}{2}\frac{\bigtriangleup _{%
\mathbf{g}}f(P)}{f(P)}\text{ \thinspace for \thinspace }P\in \left\{
x/f(x)=Maxf\right\} 
\]
then there exists $h^{\prime }$ $\leqslant h$ such that$ (h^{\prime },f,g)$ is critical.

If we only have
\[
\frac{4(n-1)}{n-2}h(P)\geq S_{\mathbf{g}}(P)-\frac{n-4}{2}\frac{%
\bigtriangleup _{\mathbf{g}}f(P)}{f(P)}\text{ \thinspace for \thinspace }%
P\in \left\{ x/f(x)=Maxf\right\} 
\]
then for any $\varepsilon >0$ there exists $h^{\prime }$ $\leq h+\varepsilon $ such that $(h^{\prime },f,g)$ 
is critical.

Weaker hypothesis are sufficient to prove the existence of positive critical functions: for example, it suffices that 
the \textit{boundary} of the set $Max f$ is a set of null measure; see \cite{C3} for full details.

\section{Critical triple 2}
We want to prove here theorem 3.
This theorem lies on the transformation formula for a critical function in a conformal change of metric (seen at 
the end of the introduction): 
\begin{center}
($h^{\prime },f$\textit{,}$\mathbf{g}^{\prime }=u^{\frac{4}{n-2}}\mathbf{g)}$
\textit{ is critical if and only if (}$h=h^{\prime }u^{\frac{4}{n-2}}-%
\frac{\triangle _{\mathbf{g}}u}{u}$\textit{,}$f$\textit{,}$\mathbf{g)}$ 
\textit{is critical.}
\end{center}
We set, for $u\in C_{+}^{\infty }(M)=\{u\in C^{\infty }(M)\,/\,u>0\}$ : 
\[
F_{h^{\prime }}(u)=h^{\prime }u^{\frac{4}{n-2}}-\frac{\bigtriangleup _{%
\mathbf{g}}u}{u} 
\]
Then:
\begin{center}
$(h^{\prime },f,\mathbf{g}^{\prime })$\textit{\ is critical if and only if }
$(F_{h^{\prime }}(u),f,\mathbf{g})$\textit{\ is critical.}
\end{center}
To prove the theorem, we therefore have to prove the existence of a function $h$ such that:

1/: $\triangle _{\mathbf{g}}u+h.u=h^{\prime }u^{\frac{n+2}{n-2}}$ has a solution $u>0$, and

2/: $(h,f,\mathbf{g})$ is critical.

Indeed, in this case $h=F_{h^{\prime }}(u)$ and $h^{\prime }$ is critical for 
$f$ and $\mathbf{g}^{\prime }=u^{\frac{4}{n-2}}\mathbf{g}$.

E. Humbert et M. Vaugon proved this theorem in the case $f=cste$
and for a manifold not conformaly diffeomorphic to the sphere [21]. Their method lies on the fact that for such a manifold ,
after a first conformal change of metric, $B_{0}(%
\mathbf{g})K(n,2)^{-2}$ is a critical function, (we will denote these two constants $K$ et $B_{0}$). In fact, 
a careful study of their proof shows that what is needed is in fact that $B_{0}K^{-2}$ is positive. But we proved 
in the previous section the existence of positive critical functions under a geometric hypothesis concerning $f$. 
Remark that our proof will work on the sphere, but only for a non-constant function $f$.

The principle of the proof of E. Humbert and M. Vaugon is the following. We know that there exists a sequence 
($h_{t}$) of sub-critical functions for $f$ and $\mathbf{g}$ such that $h_{t}%
\stackrel{C^{2}}{\rightarrow }h$ where $(h,f,\mathbf{g})$ is critical and such that for any point $P$ where $f$ is maximum on $M$
\[
\frac{4(n-1)}{n-2}h(P)>S_{\mathbf{g}}(P)-\frac{n-4}{2}\frac{\bigtriangleup _{%
\mathbf{g}}f(P)}{f(P)}. 
\]
For a sequence $q_{t}\rightarrow 2^{*},\,q_{t}<2^{*}$ we build a sequence
$u_{t}>0$ of solutions of 
\[
\triangle _{\mathbf{g}}u+h.u=h^{\prime }u^{q_{t}-1}\,\,\,with\,\,\,\int
h^{\prime }u_{t}^{q_{t}}dv_{\mathbf{g}}\leq C\text{ independant of }t 
\]
such that $u_{t}\stackrel{H_{1}^{2}}{\rightharpoondown }u\geqslant 0$. Here again, if $u>0$, 
then $u$ is solution (up to a multiplicative constant) of $\triangle _{\mathbf{g}}u+h.u=h^{\prime }u^{\frac{n+2}{n-2}}$
and we are done.

Now, if $u=0$, one shows that the $u_{t}$ concentrate and that using this phenomenom, one can find a $t_{0}$ 
close to 1 ( if e.g. $t\rightarrow 1$) and a real $s$ large, such that $F_{h^{\prime
}}(u_{t_{0}})$ is sub-critical and $F_{h^{\prime }}(u_{t_{0}}^{s})$
is weaklly critical, with furthermore
\[
\frac{4(n-1)}{n-2}F_{h^{\prime }}(u_{t_{0}}^{s})(P)>S_{\mathbf{g}}(P)-\frac{n-4}{2}\frac{%
\bigtriangleup _{\mathbf{g}}f(P)}{f(P)} 
\]
at any point $P$ where $f$ is maximum. Then, considering the path 
$t\rightarrow F_{h^{\prime }}(u_{t_{0}}^{ts})$ and using theorem 1, we get the existence of a critical 
function on this path. It is to obtain the conditions on $F_{h^{\prime }}(u_{t_{0}}^{s})$ at the maximum points 
of $f$ that we need the existence of positive critical functions.

We will now give the scheme of the proof, refering for complete details to the article of E. Humbert and M. Vaugon or to our PHD 
thesis available online, 
and we will only indicate the modifications due to our function $f$ and the necessity of positive critical functions.

First, we said that we will need positive critical functions. Their existence was proved under the hypothesis that 
$Max f$ is thin and that $\int_{M}fdv_{\mathbf{g}}>0$. But $\stackunder{M}{Sup}f>0$,
so, after making if necessary a first conformal change of metric, we can suppose that 
$\int_{M}fdv_{\mathbf{g}}>0$, and we supposed in the hypothesis of theorem 3 that $Max f$ is thin, and therefore 
we can suppose that we have positive critical function for $f$ and $\mathbf{g}$.

Then, we fix some (more) notations:
\begin{center}
\[
J_{h,h^{\prime },\mathbf{g,}q}(w)=\frac{\int_{M}\left| \nabla w\right|
^{2}dv_{\mathbf{g}}+\int_{M}h.w{{}^{2}}dv_{\mathbf{g}}}{\left(
\int_{M}h^{\prime }\left| w\right| ^{q}dv_{\mathbf{g}}\right) ^{\frac{2}{q}}}
\]
\end{center}

\[
\stackunder{w\in \mathcal{H}_{h^{\prime },q}^{+}}{\inf }J_{h,h^{\prime },%
\mathbf{g,}q}(w):=\lambda _{h,h^{\prime },\mathbf{g,}q} 
\]
where
\[
\mathcal{H}_{h^{\prime },q}^{+}=\{w\in
H_{1}^{2}(M)\,/\,\,\,w>0\,\,and\,\,\int_{M}h^{\prime }.w^{q}dv_{\mathbf{g}%
}>0\}. 
\]
and
\[
\Omega _{h,h^{\prime },\mathbf{g,}q}=\{u\in \mathcal{H}_{h^{\prime
},q}^{+}/\,J_{h,h^{\prime },\mathbf{g,}q}(u)=\lambda _{h,h^{\prime },\mathbf{%
g,}q}\,\,and\,\,\int_{M}h^{\prime }.w^{q}dv_{\mathbf{g}}=(\lambda
_{h,h^{\prime },\mathbf{g,}q}\,\,)^{\frac{q}{q-2}}\}\,\,. 
\]

Let ($h_{t}$) be a sequence of sub-critical functions for $f$ and $%
\mathbf{g}$ such that $h_{t}\stackrel{C^{2}}{\rightarrow }h$ where $(h,f,%
\mathbf{g})$ is critical, with $\triangle _{\mathbf{g}}+h_{t}$ coercive. We know that we can find such a 
sequence with $h_{t}>0$ et $h>0$, and also
\[
\frac{4(n-1)}{n-2}h_{t}(P)>S_{\mathbf{g}}(P)-\frac{n-4}{2}\frac{\bigtriangleup _{\mathbf{g}%
}f(P)}{f(P)} 
\]
for all $P\in Max f$. But here, we can say more, and that is where the existence of positive critical functions is 
crucial. Indeed, for any constant $c>0$, if $\mathbf{g}^{\prime }=c%
\mathbf{g}$, then $S_{\mathbf{g}^{\prime }}=c^{-1}S_{\mathbf{g}}$ and $%
\triangle _{\mathbf{g}^{\prime }}=c^{-1}\triangle _{\mathbf{g}}$ and by the transformation formula 
for critical functions:

\begin{center}
\textit{\ }$h$\textit{\ is (sub-, weakly) critical for }$f$\textit{\
and }$\mathbf{g}$\textit{\ if and onlu if }$c^{-1}h$\textit{\ is (sub-,
weakly) critical for }$f$\textit{\ and }$\mathbf{g}^{\prime }$.
\end{center}

Therefore, up to multiplying $\mathbf{g}$ by a constant, we can, for any constant $C>0$, suppose : 
\begin{eqnarray*}
h_{t} &>&C\,\,\,\,on\, M \\
\frac{4(n-1)}{n-2}h_{t}(P)-S_{\mathbf{g}}(P)+\frac{n-4}{2}\frac{\bigtriangleup _{\mathbf{g}%
}f(P)}{f(P)} &>&C\,\,\,\,\forall P \in Max f
\end{eqnarray*}
and $(h,f,\mathbf{g)}$ has minimizing solutions.

We can now follow the method exposed above; we only give the scheme of the proof.

\textit{First step:} Thanks to the compacity of the inclusion $H_{1}^{2}\subset L^{q}$, it is known 
that  $\forall q<2^{*}$ and $\forall u\in \Omega _{h,h^{\prime },\mathbf{g,}q}$, $u$ is solution of 
$\triangle _{\mathbf{g}}u+h.u=h^{\prime }u^{q-1}$. Using this fact, in the first step, one proves the following:

There exist sequences $(q_{i}),(t_{i}),$ such that
\begin{eqnarray*}
\,2 &<&q_{i}<2^{*} \\
\,q_{i} &\rightarrow &2^{*} \\
\,t_{i} &\rightarrow &1 \\
\,h_{t_{i}} &\rightarrow &h\,
\end{eqnarray*}
and a sequence ($v_{i})\in \Omega _{h_{t_{i}},h^{\prime },\mathbf{g,}q_{i}}$
such that ($F_{h^{\prime }}(v_{i}),f,\mathbf{g)\,}$ is sub-critical.

We note 
\[
J_{i}=J_{h_{t_{i}},h^{\prime },\mathbf{g,}q_{i}} 
\]
$\,$and
\[
\lambda _{i}=\lambda _{h_{t_{i}},h^{\prime },\mathbf{g,}q_{i}}\,\,\,. 
\]
Then
\[
J_{i}(v_{i})=\lambda _{i}\,\,\,and\,\,\int h^{\prime }v_{i}^{q_{i}}dv_{%
\mathbf{g}}=\lambda _{i}^{\frac{q_{i}}{q_{i}-2}} 
\]
and $v_{i}$ is a positive solution of 
\[
\triangle _{\mathbf{g}}v_{i}+h_{t_{i}}.v_{i}=h^{\prime }v_{i}^{q_{i}-1}. 
\]
The sequence ($v_{i})$ is bounded in $H_{1}^{2}$ and thus there exists $v\in
H_{1}^{2}$ such that
\[
v_{i}\stackrel{H_{1}^{2}}{\rightharpoondown }v,\,\,v_{i}\stackrel{L^{2}}{%
\rightarrow }v\,\,\,et\,\,\,\,v_{i}\stackrel{L^{2^{*}-2}}{\rightarrow }v. 
\]
Once again, we have two possibilities: $v\equiv 0$ or $v>0$.

\textit{Second step:}

If $v>0$, as we said above, the proof is over: up to a subsequence, $v_{i}%
\stackrel{C^{2}}{\rightarrow }v$ and so on one hand $F_{h^{\prime
}}(v_{i})\rightarrow F_{h^{\prime }}(v)$, and on the other hand
\[
F_{h^{\prime }}(v_{i})=h_{t_{i}}+h^{\prime }(v_{i}^{\frac{4}{n-2}%
}-v_{i}^{q_{i}-2})\rightarrow h\, ,
\]
that is, $F_{h^{\prime }}(v)=h$ which is critical for $f$ and $\mathbf{g}$ with minimizing solutions. 
Thus $h^{\prime }$ is critical for $f$ and $\mathbf{g}^{\prime }=v^{\frac{4}{n-2}}\mathbf{g}$, with 
minimizing solutions.

\textit{The rest of the proof is therefore concerned with the case }$v\equiv 0$\textit{.}

\textit{Third step:} One proves that there is a concentration phenomenom:

a/: One first shows that: 
\[
0<c\leqslant \overline{\lim }\,\lambda _{i}\leqslant
K^{-2}(Sup_{M}\,h^{\prime })^{-\frac{n-2}{n}}\,\,. 
\]

b/: Second, one shows that: 
\[
0<\lambda ^{\frac{n}{2}}(Sup_{M}\,h^{\prime })^{-1}\leqslant \overline{\lim }%
\,\int_{M}v_{i}^{q_{i}}dv_{\mathbf{g}}\leqslant K^{2^{*}}\lambda ^{\frac{%
n2^{*}}{4}}\leqslant K^{-n}(Sup_{M}\,h^{\prime })^{-\frac{n}{2}} 
\]
where $\lambda >0$ is such that, after extraction, $\lambda_{i} \rightarrow \lambda$.

c/: We say that $x\in M$ is a concentration point if
\[
\forall r>0:\,\overline{\lim }\,\int_{B(x,r)}v_{i}^{q_{i}}dv_{\mathbf{g}}>0. 
\]

Using a/ and b/, and method analogous to section 4.2, one gets the following:

First, as $M$ is compact, there exists at least one concentration point $x\in M$.

Then, using the iteration process, one shows that 
\[
\overline{\lim }\,\int_{B(x,r)}v_{i}^{q_{i}}dv_{\mathbf{g}}\geqslant
K^{-n}(Sup_{M}\,h^{\prime })^{-\frac{n}{2}}. 
\]

d/: Therefore using the method of section 4.2, we get:

1/: $\overline{\lim }\,\int_{B(x,r)}v_{i}^{q_{i}}dv_{\mathbf{g}%
}=K^{-n}(Sup_{M}\,h^{\prime })^{-\frac{n}{2}}$ , $\forall r>0$

2/: $x$ is the only concentration point, denoted $x_{0}$

3/: $\lambda =K^{-2}(Sup_{M}\,h^{\prime })^{-\frac{n-2}{n}}$

4/: $x_{0}$ is a point of maximum of $h^{\prime }$

5/: $v_{i}\rightarrow 0$ in $C_{loc}^{2}(M-\{x_{0}\})$
\\

\textit{Fourth step:}

We know now that the sequence $(v_{i})$ concentrates in $x_{0}$ and that for any 
$i$ $F_{h^{\prime }}^{\prime }(v_{i})$ is sub-critical for $f$ and $\textbf{g}$. We would like to find a 
$v_{i_{0}}$ , a function $v>0$ and a continuous path from
$v_{i_{0}}$ to $v$ such that $F_{h^{\prime }}(v)$ is weakly critical for $f$ and $\textbf{g}$ and such that
\[
\frac{4(n-1)}{n-2}F_{h^{\prime }}(v)(P)>S_{\mathbf{g}}(P)-\frac{n-4}{2}\frac{\bigtriangleup _{%
\mathbf{g}}f(P)}{f(P)} 
\]
for all $P \in Max f$ . Then, the theorem 1 will tell us that on the path $u_{t}$ from $v_{i_{0}}$ to $v$ 
there exists a $u_{t}$ such that $F_{h^{\prime }}(u_{t})$ is critical for $f $ and $\mathbf{g}$ .

That is where we are going to use the existence of positive critical functions.

Let $s\geqslant 1$ and let $v$ be a positive function. Then 
\[
\bigtriangleup _{\mathbf{g}}(v^{s})=sv^{s-1}\bigtriangleup _{\mathbf{g}%
}v-s(s-1)v^{s-2}\left| \nabla v\right| _{\mathbf{g}}^{2}\,\,. 
\]
Thus
\[
F_{h^{\prime }}(v_{i}^{s})=h^{\prime }v_{i}^{s\frac{4}{n-2}%
}+sh_{t_{i}}-sh^{\prime }v_{i}^{q_{i}-2}+s(s-1)\frac{\left| \nabla v\right|
_{\mathbf{g}}^{2}}{v_{i}^{2}} 
\]
and therefore
\[
F_{h^{\prime }}(v_{i}^{s})\geqslant sh_{t_{i}}+h^{\prime }(v_{i}^{s\frac{4}{%
n-2}}-sv_{i}^{q_{i}-2})\,\,. 
\]
Now:

On $\{x\in M\,/\,h^{\prime }(x)\leqslant 0\}$ :

$v_{i}\rightarrow 0$ uniformly because $x_{0} \in Max \, h^{\prime }$ and $h^{\prime }(x_{0})>0$ as we 
have supposed that $\bigtriangleup _{\mathbf{g}}+h^{\prime }$ is coercive. Furthermore if 
$s\geqslant 1$ then $s\frac{4}{n-2}\geqslant q_{i}-2$. Thus, for $i$ large enough 
\[
F_{h^{\prime }}(v_{i}^{s})\geqslant sh_{t_{i}}\text{ on }\{x\in
M\,/\,h^{\prime }(x)\leqslant 0\} 
\]

On $\{x\in M\,/\,h^{\prime }(x)>0\}$ :

We consider the function of a real variable defined for $x\geqslant 0$ by 
\[
\beta _{i,s}(x)=x^{s\frac{4}{n-2}}-sx^{q_{i}-2}=x^{q_{i}-2}(x^{s\frac{4}{n-2}%
-q_{i}+2}-s). 
\]
An easy study of this function shows that
\[
for\,\,\,x\geqslant 0:\beta _{i,s}(x)\geqslant -s. 
\]
But 
\[
F_{h^{\prime }}(v_{i}^{s})\geqslant sh_{t_{i}}+h^{\prime }\beta
_{i,s}(v_{i}) 
\]
therefore 
\[
F_{h^{\prime }}(v_{i}^{s})\geqslant sh_{t_{i}}-sh^{\prime }\,\,\,\text{on}%
\,\,\,\{x\in M\,/\,h^{\prime }(x)>0\}. 
\]
We can therefore write : 
\[
F_{h^{\prime }}(v_{i}^{s})\geqslant s(h_{t_{i}}-\stackunder{M}{Sup}\,h^{\prime
})\text{ on }\{x\in M\,/\,h^{\prime }(x)>0\}\text{.} 
\]
We now use our work from the beginning of the proof, that is that, for any $C>0$, we can suppose that: 
\[
h_{t}>C\,\,on\, M 
\]
and 
\[
\frac{4(n-1)}{n-2}h_{t}(P)-S_{\mathbf{g}}(P)+\frac{n-4}{2}\frac{\bigtriangleup _{\mathbf{g}%
}f(P)}{f(P)}>C, \,\,\,\, \forall P\in Max f\,. 
\]
Then, first, if we suppose that $h>\stackunder{M}{Sup}\,h^{\prime }$ on $M$, we see that for $i$ and $s$ 
large enough : 
\begin{equation}
F_{h^{\prime }}(v_{i}^{s})\geqslant B_{0}(\mathbf{g})K(n,2)^{-2}  
\end{equation}
and therefore $F_{h^{\prime }}(v_{i}^{s})$ is weakly critical for $f$ and $\mathbf{g}$ . 
Beside, for all $t\in [1,s]$ we also have
\[
F_{h^{\prime }}(v_{i}^{t})\geqslant t(h_{t_{i}}-\stackunder{M}{Sup}\,h^{\prime
})\text{ }\geqslant h_{t_{i}}-\stackunder{M}{Sup}\,h^{\prime }>0 
\]
so $\bigtriangleup _{\mathbf{g}}+F_{h^{\prime }}(v_{i}^{t})$ is coercive.

Secondly, if we also suppose that
\[
\frac{4(n-1)}{n-2}h_{t}(P)-S_{\mathbf{g}}(P)+\frac{n-4}{2}\frac{\bigtriangleup _{\mathbf{g}%
}f(P)}{f(P)}>\frac{4(n-1)}{n-2}\stackunder{M}{Sup}\ h^{\prime }\,\,\,\,\, \forall P \in Max f
\]
we have for all $t\in [1,s]$: 
\begin{equation}
\frac{4(n-1)}{n-2}F_{h^{\prime }}(v_{i}^{t})(P)>S_{\mathbf{g}}(P)-\frac{n-4}{2}\frac{%
\bigtriangleup _{\mathbf{g}}f(P)}{f(P)}\,\,\,\,\,\,\, \forall P \in Max f
\end{equation}
as soon as $i$ is large enough.

We therefore fix $i$ and $s$ large enough to have (24) et (25) and we consider
\[
s_{0}=\inf \{t>1\,/\,F_{h^{\prime }}(v_{i}^{t})\text{ is weakly critical\}} 
\]

We then apply theorem 1 to the path $t\in [1,s_{0}]\mapsto
F_{h^{\prime }}(v_{i}^{t})$ to obtain that $F_{h^{\prime
}}(v_{i}^{s_{0}})$ is critical for $f$ and $\mathbf{g}$, with minimizing solutions. 
Therefore $h^{\prime }$ is critical for $f$ and $\mathbf{g}^{\prime }=(v_{i}^{s_{0}})^{\frac{4}{n-2}}\mathbf{g}$ 
with minimizing solutions.

This ends the proof.

\section{Critical triple 3}
Let $(M,\mathbf{g})$ be a compact riemannian manifold of dimension 
$n\geqslant 3$. Let $h$ be a fixed $C^{\infty }$ function such that $\triangle _{\mathbf{g}}+h$ 
is coercive. The problem we want to study is the following: \textit{can we find a function $f$ such that 
$(h,f,\g)$ is a critical triple ?}

We first make a remark. If $h\geqslant B_{0}(\mathbf{g})K(n,2){%
{}^{-2}}$, then $h$ is weakly critical for any function $f$, and there cannot exist a function $f$ such that 
($h,f,\mathbf{g})$ is subcritical. But more important is the next observation:

\textit{If there exist a non constant function $f$ such that $(h,f,\g)$ is critical with a minimizing solution $u$, 
then $(h,1,\g)$ is sub-critical. }

Indeed, as we saw in section 1, we can suppose that $Sup\,f=1$. Then, as $u>0$%
\[
J_{h,1}(u)<J_{h,f}(u)=\frac{1}{K(n,2)^{2}(Sup\,f)^{\frac{2}{2^{*}}}}=\frac{1%
}{K(n,2)^{2}} 
\]
and therefore $h$ is subcritical for $1$.

We want to prove that, at least if $dim M\geq5$, this necessary condition is sufficient, i.e we want to prove theorem 4. 
We thus suppose now that $(h,1,\mathbf{g})$\textit{\ is sub-critical.}

The proof will proceed in two steps:

First step: we prove that there exist a function $f\in C^{\infty }(M)$
such that $\stackunder{M}{Sup}f=1$, with
$\triangle _{\mathbf{g}}f$ being as large as we want in its maximum points, 
and such that $(h,f,\mathbf{g})$ is weakly critical.

Second step: being given this function $f$, we prove that there exists on the path 
\[
t\rightarrow f_{t}=t.1+(1-t)f 
\]
a function for which $h$ is critical.

\textbf{First step:}

We proceed by contradiction. We suppose that for any $f\in
C^{\infty }(M)$ such that $\stackunder{M}{Sup}f>0$, $(h,f,\mathbf{g})$
is sub-critical. Then, for all such function, there exit positive solution $u$ to the equation 
\[
\triangle _{\mathbf{g}}u+h.u=\lambda .f.u^{\frac{n+2}{n-2}} 
\]
where 
\[
\lambda =\stackunder{w\in \mathcal{H}_{f}}{\inf }I_{h,\mathbf{g}%
}(w)\,\,\,and\,\,\,\int_{M}f.u^{\frac{2n}{n-2}}dv_{\mathbf{g}}=1\,. 
\]

The metric $\mathbf{g}$ being fixed, we will not write $dv_{%
\mathbf{g}}$ in the integrals.

The idea is to build a familly of functions $f_{t}$ whose laplacians tend to infinity at the maximum points. 
One of these function will then give a weakly critical triple $(h,f_{t},\mathbf{g})$. 
Furthermore, our proof holding for any subsequence of this familly. this function will have a laplacian as large 
as we want in its point of maximum.

In $\Bbb{R}^{n}$, we build for $t\rightarrow 0$ a familly ($P_{t})$
of $C^{\infty }$ functions, similar to a regularizing sequence, such that
\begin{eqnarray*}
0 &\leqslant &P_{t}\leqslant 1 \\
\,\,P_{t}(x) &=&P_{t}(\left| x\right| ) \\
\,\,P_{t}(0) &=&1 \\
\left\| \nabla P_{t}\right\| &\sim &\frac{c_{1}}{t}\,\,\,on\,\,\,B(0,t) \\
\left| \bigtriangleup P_{t}(0)\right| &\sim &\frac{c_{2}}{t^{2}} \\
SuppP_{t} &=&B(0,t).
\end{eqnarray*}
Let now $x_{0}$ be a point of $M$ such that $h(x_{0})>0$; this point exists because 
$\triangle _{\mathbf{g}}+h$ is coercive. We define 
\[
f_{t}=P_{t}\circ \exp _{x_{0}}^{-1} 
\]
We are therefore supposing that, for all $t$, $(h,f_{t},\mathbf{g})$ is sub-critical and we are looking for a 
contradiction. For all $t$ we have a solution $u_{t}>0$ of 
\[
(E_{t}):\,\,\bigtriangleup _{\mathbf{g}}u_{t}+h.u_{t}=\lambda
_{t}.f_{t}.u_{t}^{\frac{n+2}{n-2}} 
\]
with $\int f_{t}u_{t}^{2^{*}}dv_{\mathbf{g}}=1$ and
\[
\lambda _{t}<K^{-2}(\stackunder{M}{Sup}f_{t})^{-\frac{n-2}{n}}=K^{-2}. 
\]
Then, ($u_{t})$ is bounded in $H_{1}^{2}(M)$ when $t\rightarrow 0$.
So ($u_{t})$ is bounded in $L^{2^{*}}$ and ($u_{t}^{2^{*}-1})$ is bounded in 
$L^{\frac{2^{*}}{2^{*}-1}}$. After extraction of a subsequence, if 
$f_{t}\stackrel{L^{2}}{\rightarrow }f$ and $u_{t}\stackrel{L^{2}}{\rightarrow }u$,
then
\[
f_{t}u_{t}^{2^{*}-1}\rightharpoondown fu^{2^{*}-1}. 
\]
But here, $f_{t}\stackrel{L^{p}}{\rightarrow }0$, therefore the equation ($E_{t})$
``converge'' to 
\[
\bigtriangleup _{\mathbf{g}}u+h.u=0 
\]
in the sense that $u$ is solution of this equation. But $\bigtriangleup
_{\mathbf{g}}+h$ is coercive, therefore $u=0$, i.e. $u_{t}\rightarrow 0$ in $L^{p}$ for $p<2^{*}$.

The sequence ($u_{t})$ therefore concentrates in the sense we saw in subsection 4.2. 
But in subsection 4.2, the function $f$ on the right handside of the equation was constant and it was on the 
left handside that we had a sequence $(h_{t})$. However the results we saw there remain true, only the blow-up 
necessary for the weak estimates requires a new treatment. We will go over these results, only detailing the new 
difficulties.

\textit{a/: There exists, up to a subsequence of $(u_{t})$, exactly one concentration point 
and it is the point $x_{0}$ where the $f_{t}$ are maximum on $M$. Moreover}

\[
\forall \delta >0\mathit{,\ }\overline{\stackunder{t\rightarrow 1}{\lim }}%
\int_{B(x_{0},\delta )}f_{t}u_{t}^{2^{*}}=1\,. 
\]

The method of subsection 4.2 works here. More precisely, as $Supp\,f_{t}=B(x_{0},t)$, we have for all 
$\delta >0$ and as soon as $t<\delta $:
\[
\int_{B(x_{0},\delta )}f_{t}u_{t}^{2^{*}}=1. 
\]
We can also suppose that
\[
\lambda _{t}\rightarrow \lambda =K^{-2}(\stackunder{M}{Sup}f_{t})^{-\frac{n-2%
}{n}}=K^{-2}. 
\]

\textit{b/: }$u_{t}\rightarrow 0$\textit{\ in }$C_{loc}^{0}(M-\{x_{0}\})$

Same proof as in subsection 4.2.

\textit{c/: weak estimates}

We consider a sequence of points $(x_{t})$ such that $$m_{t}=\stackunder{M}{Max}\,u_{t}=u_{t}(x_{t}):=\mu _{t}^{-\frac{n-2}{2}}.$$

From the previous point, $x_{t}\rightarrow x_{0}$ and $\mu
_{t}\rightarrow 0$. Remember that $\overline{u}_{t},\overline{f}_{t},%
\overline{h}_{t},\mathbf{g}_{t}$ are the functions and the metric seen in the chart $\exp _{x_{t}}^{-1}$, 
and $\,\,\widetilde{u}_{t}\,,\,\widetilde{h}_{t}\,,\,\widetilde{f}_{t},\widetilde{\mathbf{g}}%
_{t} $ are the functions after
blow-up of center $x_{t}$ and coefficient $k_{t}=\mu _{t}^{-1}$.

Reviewing the proof of the weak estimates in section 4.2, we see that it will work here if we obtain :

\[
\forall R>0\,\,\,:\stackunder{t\rightarrow 0}{\lim }\int_{B(x_{t},R\mu
_{t})}f_{t}u_{t}^{2^{*}}dv_{\mathbf{g}}=1-\varepsilon _{R}\,\,\,where%
\,\,\,\varepsilon _{R}\stackunder{R\rightarrow +\infty }{\rightarrow }0\,. 
\]
This relation is itself proved using blow-up theory once it is proved that $\,\widetilde{u}%
_{t}\stackrel{C_{loc}^{2}(\Bbb{R}^{n})}{\rightarrow }\,\widetilde{u}$ where
$\widetilde{u}$ is solution of :
\[
\bigtriangleup _{e}\widetilde{u}=K^{-2}\widetilde{u}^{\frac{n+2}{n-2}}. 
\]
This is where we have the main difficulty due to the presence of a familly ($f_{t}$). 
Indeed, after blow-up, the equation
\[
(E_{t}):\,\,\bigtriangleup _{\mathbf{g}}u_{t}+h.u_{t}=\lambda
_{t}.f_{t}.u_{t}^{\frac{n+2}{n-2}} 
\]
becomes
\[
(\widetilde{E}_{t})\,:\,\triangle _{\widetilde{\mathbf{g}}_{t}}\widetilde{u}%
_{t}+\mu _{t}^{2}.\widetilde{h}_{t}.\widetilde{u}_{t}=\lambda _{t}\widetilde{%
f}_{t}.\widetilde{u}_{t}^{\frac{n+2}{n-2}} 
\]
and to obtain that this equation "converges" to
\[
\bigtriangleup _{e}\widetilde{u}=K^{-2}\widetilde{u}^{\frac{n+2}{n-2}} 
\]
we need to show that $(\,\widetilde{f}_{t})$ is simply convergent to 1
(which is obvious when we have a constant function $f$ on the right handside of 
($E_{h,f,\mathbf{g}}$)). As the sequence ($\widetilde{f}_{t}$) is uniformly bounded by 1 
on $\Bbb{R}^{n}$ (considering we have extended
$\widetilde{f}_{t}$ by 0 on $\Bbb{R}^{n}\backslash B(0,\delta \mu
_{t}^{-1})$), we have, using e.g. theorem 8.25 of
Gilbard-Trudinger \cite{G-T} and Ascoli's theorem,
the existence of a function $\widetilde{u}\in C^{0}(\Bbb{R}^{n})$ such that, 
after extraction, $\widetilde{u}_{t}\stackrel{C_{loc}^{0}(\Bbb{R}%
^{n})}{\rightarrow }\,\widetilde{u}$, with $\widetilde{u}(0)=1$.

We are going to prove that $\widetilde{f}_{t}\stackrel{a.e.}{\rightarrow }1$ on 
$\Bbb{R}^{n}$ in two steps (we will prove a little bit more):

1/: There exists $\widetilde{f}\in L_{loc}^{2}(\Bbb{R}^{n})$ such that $%
\widetilde{f}_{t}\stackrel{a.e.}{\rightarrow }\widetilde{f}$ on $\Bbb{R}^{n}$

2/: $\widetilde{f}=1$ a.e. on $\Bbb{R}^{n}$\\

First step:

We have $\widetilde{f}_{t}(x)=\overline{f}_{t}(\mu _{t}x)$ and $\left| \nabla 
\overline{f}_{t}\right| \leq \frac{c}{t}$. Therefore
\[
\left| \nabla \widetilde{f}_{t}\right| \leq c.\frac{\mu _{t}}{t}\,. 
\]

We consider two cases:

a/: If ($\frac{\mu _{t}}{t}$) is bounded: Then for any compact set $%
K\subset \subset \Bbb{R}^{n}$, ($\widetilde{f}_{t}$) is bounded in $%
H_{1}^{n+1}(K)$ (where $n=\dim M)$. Thus, by compacity of the inclusion $%
H_{1}^{n+1}(K)\subset C^{0,\alpha }(K)$ for some $\alpha >0$, up to a subsequence, 
there exists $\widetilde{f}_{K}\in C^{0,\alpha }(K)$ such that
\[
\widetilde{f}_{t}\stackrel{C^{0,\alpha }(K)}{\rightarrow }\widetilde{f}_{K} 
\]
By diagonal extraction, we constuct $\widetilde{f}\in C^{0,\alpha }(%
\Bbb{R}^{n})$ such that
\[
\widetilde{f}_{t}\stackrel{C^{0,\alpha }(K^{\prime })}{\rightarrow }%
\widetilde{f} 
\]
for any compact set $K^{\prime }$ of $\Bbb{R}^{n}$, and moreover $\widetilde{f}%
\in H_{1,loc}^{n+1}(\Bbb{R}^{n}).$ So $\widetilde{f}_{t}%
\stackrel{a.e.}{\rightarrow }\widetilde{f}\,\,\,on\,\,\,\Bbb{R}^{n}\,.$

b/: If $\frac{\mu _{t}}{t}\rightarrow +\infty $ : the support of $\widetilde{f%
}_{t}$ is
\[
Supp\widetilde{f}_{t}=B(\frac{x_{0}(t)}{\mu _{t}},\frac{t}{\mu _{t}}), 
\]
where $x_{0}(t)=\exp _{x_{t}}^{-1}(x_{0})$.

If ($\frac{\left| x_{0}(t)\right| }{\mu _{t}}$) is bounded, there is after extraction a subsequence
\[
\frac{x_{0}(t)}{\mu _{t}}\rightarrow P\in \Bbb{R}^{n}; 
\]
and therefore
\[
\widetilde{f}_{t}\stackrel{C_{loc}^{0}(\Bbb{R}^{n}-\{P\})}{\rightarrow }0 
\]

If $\frac{\left| x_{0}(t)\right| }{\mu _{t}}\rightarrow \infty $, then
\[
\widetilde{f}_{t}\stackrel{C_{loc}^{0}(\Bbb{R}^{n})}{\rightarrow }0 
\]
In both cases, $\widetilde{f}_{t}\stackrel{p.p}{\rightarrow }%
0\,\,\,on\,\,\,\Bbb{R}^{n}.$

In case a/, $\widetilde{u}$ is a weak solution of 
\[
\bigtriangleup _{e}\widetilde{u}=K^{-2}\widetilde{f}\widetilde{u}^{\frac{n+2%
}{n-2}}\,\, 
\]
with $\widetilde{f}\geq 0$ as $\widetilde{f}_{t}\geq 0$, and $\widetilde{%
f}\in H_{1,loc}^{n+1}(\Bbb{R}^{n})\subset C^{0,\alpha }(\Bbb{R}^{n}).\,$

In case b/, $\,\widetilde{u}$ is a weak solution of
\[
\bigtriangleup _{e}\widetilde{u}=0.\,\, 
\]
In both cases, elliptic thory and standard regularity thorems gives the $C^{2}$ regularity of $\widetilde{u}$
, and therefore $\bigtriangleup _{e}\widetilde{u}\geq 0$. The maximum principle then shows that 
either $\widetilde{u}\equiv 0$ or $\widetilde{u}%
>0 $. But $\widetilde{u}(0)=1$ thus $\widetilde{u}>0$.

Second step:

We start using the iteration process : for some cut-off function 
$\eta $ equal to 1 near $x_{0}$, we multiply ($E_{t}$) by $\eta
^{2}u_{t}$, integrate and use the Sobolev inequality to obtain, remembering that $\lambda _{t}<K^{-2}(\stackunder{M}{Sup}%
f_{t})^{-\frac{n-2}{n}}$ and that $Sup\,f_{t}=1$:
\[
(\int_{M}(\eta u_{t})^{2^{*}})^{\frac{2}{2^{*}}}\leq \lambda
_{t}K^{2}\int_{M}\eta ^{2}f_{t}u_{t}^{2^{*}}+c\int_{Supp\,\eta
}u_{t}^{2}\,\,. 
\]
We take $\eta =1$ on $B(x_{0},\frac{3}{2}\delta )$ and $\eta =0$ on $%
M\backslash B(x_{0},2\delta )$. Then for $t$ close to 0 
\[
Supp\,f_{t}\subset B(x_{0},t)\subset B(x_{t},\delta )\subset B(x_{0},\frac{3%
}{2}\delta ) 
\]
So
\[
(\int_{B(x_{t},\delta )}u_{t}{}^{2^{*}})^{\frac{2}{2^{*}}}\leq
\int_{B(x_{t},\delta )}f_{t}u_{t}^{2^{*}}+c\int_{M}u_{t}^{2} 
\]
and after blow-up
\[
(\int_{B(0,\delta \mu _{t}^{-1})}\widetilde{u}_{t}^{2^{*}})^{\frac{2}{2^{*}}%
}\leq \int_{B(0,\delta \mu _{t}^{-1})}\widetilde{f}_{t}\widetilde{u}%
_{t}^{2^{*}}+c\int_{M}u_{t}^{2}=1+c\int_{M}u_{t}^{2}\,\,. 
\]
But $\int_{M}u_{t}^{2}\rightarrow 0$ therefore
\[
\stackunder{t\rightarrow 0}{\overline{\lim }}\int_{B(0,\delta \mu _{t}^{-1})}%
\widetilde{u}_{t}^{2^{*}}\leq 1\,\,. 
\]
Beside, we know that $\widetilde{f}_{t}\stackrel{a.e.}{\rightarrow }%
\widetilde{f}$ with $\widetilde{f}\leq 1$ and $\widetilde{u}_{t}(0)=1$.
Let suppose that there exists a set $A\subset \Bbb{R}^{n}$ with $mes(A)>0$ such that $\widetilde{f}<1$ on 
$A$ and write $\Bbb{R}^{n}=A\cup B$ with $%
\widetilde{f}=1$ a.e. on $B$. Then, as $\widetilde{f}_{t}\geq 0$ and as 
$\widetilde{u}_{t}\stackrel{C^{2}}{\rightarrow }\widetilde{u}>0$ :
\begin{eqnarray*}
1=\int_{B(0,\delta \mu _{t}^{-1})}\widetilde{f}_{t}\widetilde{u}_{t}^{2^{*}}
& &\leqslant \stackunder{t\rightarrow 0}{\overline{\lim }}\int_{B(0,\delta
\mu _{t}^{-1})\cap A}\widetilde{f}_{t}\widetilde{u}_{t}^{2^{*}}+\stackunder{%
t\rightarrow 0}{\overline{\lim }}\int_{B(0,\delta \mu _{t}^{-1})\cap B}%
\widetilde{f}_{t}\widetilde{u}_{t}^{2^{*}} \\ 
&& <\stackunder{t\rightarrow 0}{\overline{\lim }}\int_{B(0,\delta \mu
_{t}^{-1})\cap A}\widetilde{u}_{t}^{2^{*}}+\stackunder{t\rightarrow 0}{%
\overline{\lim }}\int_{B(0,\delta \mu _{t}^{-1})\cap B}\widetilde{u}%
_{t}^{2^{*}} \\ 
&& =\stackunder{t\rightarrow 0}{\overline{\lim }}\int_{B(0,\delta \mu
_{t}^{-1})}\widetilde{u}_{t}^{2^{*}}
\end{eqnarray*}
so 
\[
1<\stackunder{t\rightarrow 0}{\overline{\lim }}\int_{B(0,\delta \mu
_{t}^{-1})}\widetilde{u}_{t}^{2^{*}} 
\]
which is a contradiction, and therefore $\widetilde{f}_{t}\stackrel{a.e.}{%
\rightarrow }1$ on $\Bbb{R}^{n}$.

Thus, as we said
\[
(\widetilde{E}_{t})\,:\,\triangle _{\widetilde{\mathbf{g}}_{t}}\widetilde{u}%
_{t}+\mu _{t}^{2}.\widetilde{h}_{t}.\widetilde{u}_{t}=\lambda _{t}\widetilde{%
f}_{t}.\widetilde{u}_{t}^{\frac{n+2}{n-2}} 
\]
``converges'' to 
\[
\bigtriangleup _{e}\widetilde{u}=K^{-2}\widetilde{u}^{\frac{n+2}{n-2}} 
\]
in the sense that
\[
\,\widetilde{u}_{t}\stackrel{C_{loc}^{2}(\Bbb{R}^{n})}{\rightarrow }\,%
\widetilde{u} 
\]
where $\widetilde{u}$ is a solution of $\bigtriangleup _{e}\widetilde{u}%
=K^{-2}\widetilde{u}^{\frac{n+2}{n-2}}$. As $\widetilde{u}(0)=1$, 
\[
\widetilde{u}(x)=(1+\frac{K^{-2}}{n(n-2)}\left| x\right| ^{2})^{-\frac{n-2}{2%
}}\,\,. 
\]

Now, we can proceed exactly as in subsection 4.2.. We have:

\[
\forall R>0:\stackunder{t\rightarrow 0}{\lim }\int_{B(x_{t},R\mu
_{t})}f_{t}u_{t}^{2^{*}}dv_{\mathbf{g}}=1-\varepsilon _{R}\,\,\,where%
\,\,\,\varepsilon _{R}\stackunder{R\rightarrow +\infty }{\rightarrow }0 
\]
then

\[
\exists C>0\,\,\,such\,\,\,that\,\,\,\forall x\in M:\,d_{\mathbf{g}}(x,x_{t})^{%
\frac{n-2}{2}}u_{t}(x)\leq C. 
\]
and

\[
\forall \varepsilon >0,\exists R>0\,\,\,such\,\,\,that\,\,\,\forall
t,\,\forall x\in M:\,d_{\mathbf{g}}(x,x_{t})\geq R\mu _{t}\,\Rightarrow
\,\,d_{\mathbf{g}}(x,x_{t})^{\frac{n-2}{2}}u_{t}(x)\leq \varepsilon . 
\]

d/: We have here again the $L^{2}$-concentration:

If $\dim M\geq 4,$ 
\[
\forall \delta >0\,:\,\stackunder{t\rightarrow 0}{\lim }\frac{%
\int_{B(x_{0},\delta )}u_{t}^{2}dv_{\mathbf{g}}}{\int_{M}u_{t}^{2}dv_{%
\mathbf{g}}}=1 
\]

e/: We also have the strong estimates:
For $0<\nu <\frac{n-2}{2}$
\[
\exists C(\nu)>0\,\,\,such\,\,\,that\,\,\,\forall x\in M:\,d_{\mathbf{g}%
}(x,x_{t})^{n-2-\nu}\mu _{t}^{-\frac{n-2}{2}+\nu}u_{t}(x)\leq C, 
\]
and therefore the strong $L^{p}$-concentration:

\textbf{\ }$\forall R>0$ , $\forall \delta >0$ and $\forall p>\frac{n}{n-2}$
where $n=\dim M$:
\[
\stackunder{t\rightarrow 0}{\lim }\frac{\int_{B(x_{t},R\mu
_{t})}u_{t}^{p}dv_{\mathbf{g}}}{\int_{B(x_{t},\delta )}u_{t}^{p}dv_{\mathbf{g%
}}}=1-\varepsilon _{R}\,\,\,where\,\,\,\varepsilon _{R}\stackunder{%
R\rightarrow +\infty }{\rightarrow }0 
\]

We can now proceed with the central part of the proof of theorem 4:

We consider the euclidean Sobolev inequality and equation 
($E_{t})$ viewed in the chart $\exp _{x_{t}}^{-1}$. Using the same computations as in subsection 4.3, we get:
\begin{eqnarray*}
\int_{B(0,\delta )}\overline{h}_{t}(\eta \overline{u}_{t})^{2}dx\leq && \frac{%
1}{K(n,2){{}^{2}}(\stackunder{M}{Sup}f)^{\frac{n-2}{n}}}\int_{B(0,\delta )}%
\overline{f}_{t}\eta ^{2}\overline{u}_{t}^{2^{*}}dx\\
&&-\frac{1}{K(n,2){{}^{2}}}%
(\int_{B(0,\delta )}(\eta \overline{u}_{t})^{2^{*}}dx)^{\frac{2}{2^{*}}} \\ 
&&+C.\delta ^{-2}\int_{B(0,\delta )\backslash B(0,\delta /2)}\overline{u}%
_{t}^{2}dx+B_{t}+C_{t}
\end{eqnarray*}
with

$B_{t}=\frac{1}{2}\int_{B(0,\delta )}(\partial _{k}(\,\mathbf{g}%
\,_{t}^{ij}\Gamma (\,\mathbf{g}\,_{t})_{ij}^{k}+\partial _{ij}\,\mathbf{g}%
\,_{t}^{ij})(\eta \overline{u}_{t}^{2})dx$

$C_{t}=\left| \int_{B(0,\delta )}\eta ^{2}(\,\mathbf{g}\,_{t}^{ij}-\delta
^{ij})\partial _{i}\overline{u}_{t}\partial _{j}\overline{u}_{t}dx\right| $

$A_{t}=\frac{1}{K(n,2){{}^{2}}(\stackunder{M}{Sup}f)^{\frac{n-2}{n}}}%
\int_{B(0,\delta )}\overline{f}_{t}\eta ^{2}\overline{u}_{t}^{2^{*}}dx-\frac{%
1}{K(n,2){{}^{2}}}(\int_{B(0,\delta )}(\eta \overline{u}_{t})^{2^{*}}dx)^{%
\frac{2}{2^{*}}}$

We can write
\[
A_{t}\leq \frac{1}{K(n,2){{}^{2}}(\stackunder{M}{Sup}f_{t})^{\frac{n-2}{n}}}%
(A_{t}^{1}+A_{t}^{2}) 
\]
where $\,A_{t}^{1}=(\int_{B(0,\delta )}\overline{f}_{t}(\eta \overline{u}%
_{t})^{2^{*}}dx)^{\frac{n-2}{n}}\,-(Supf_{t}.\int_{B(0,\delta )}(\eta 
\overline{u}_{t})^{2^{*}}dx)^{\frac{n-2}{n}}\,.$

from the computation of subsection 4.3, 
\[
\stackunder{t\rightarrow 0}{\overline{\lim }}\frac{K(n,2){{}^{-2}}(%
\stackunder{M}{Sup}f_{t})^{-\frac{n-2}{n}}A_{t}^{2}+C.\delta
^{-2}\int_{B(0,\delta )\backslash B(0,\delta /2)}\overline{u}%
_{t}^{2}dx+B_{t}+C_{t}}{\int_{B(0,\delta )}\overline{u}_{t}^{2}dx}\leq \frac{%
n-2}{4(n-1)}S_{\mathbf{g}}(x_{0})+\varepsilon _{\delta } 
\]
where $\varepsilon _{\delta }\rightarrow 0$ when $\delta \rightarrow 0$.

We now consider : 
\[
\stackunder{t\rightarrow 0}{\overline{\lim }}\frac{A_{t}^{1}}{%
\int_{B(0,\delta )}\overline{u}_{t}^{2}dx} 
\]
We remark that from its definition, $f_{t}$ is decreasing when $t\rightarrow
0 $ in the sense that: 
\[
if\,\,\,t\leq t^{\prime }\,\,\,then\,\,\,f_{t}\leq f_{t^{\prime }}\,\,. 
\]
We fix a $t_{0}.$ Then, for any $t\leq t_{0}$%
\begin{eqnarray*}
\int_{B(0,\delta )}\overline{f}_{t}(\eta \overline{u}_{t})^{2^{*}}dx& & 
=\int_{B(x_{t},\delta )}f_{t}.(\eta \circ \exp
_{x_{t}}^{-1})^{2^{*}}.u_{t}^{2^{*}}.(\exp _{x_{t}}^{-1})^{*}dx \\ 
&& \leq \int_{B(x_{t},\delta )}f_{t_{0}}.(\eta \circ \exp
_{x_{t}}^{-1})^{2^{*}}.u_{t}^{2^{*}}.(\exp _{x_{t}}^{-1})^{*}dx \\ 
& &=\int_{B(0,\delta )}(f_{t_{0}}\circ \exp _{x_{t}})(\eta \overline{u}%
_{t})^{2^{*}}dx\,.
\end{eqnarray*}
We note: 
\[
\overline{f}_{t_{0},t}=f_{t_{0}}\circ \exp _{x_{t}}\, 
\]
and
\[
\widetilde{f}_{t_{0},t}=\overline{f}_{t_{0},t}\circ \psi _{\mu _{t}^{-1}}^{-1}\,. 
\]
Then : 
\begin{eqnarray*}
A_{t}^{1} & \leq (\int_{B(0,\delta )}\overline{f}_{t_{0},t}(\eta \overline{u}%
_{t})^{2^{*}}dx)^{\frac{n-2}{n}}\,-(Supf_{t}.\int_{B(0,\delta )}(\eta 
\overline{u}_{t})^{2^{*}}dx)^{\frac{n-2}{n}} \\ 
& \leq (\int_{B(0,\delta )}\overline{f}_{t_{0},t}(\eta \overline{u}%
_{t})^{2^{*}}dx)^{\frac{n-2}{n}}\,-(Supf_{t_{0}}.\int_{B(0,\delta )}(\eta 
\overline{u}_{t})^{2^{*}}dx)^{\frac{n-2}{n}}
\end{eqnarray*}
as $Supf_{t}=Supf_{t_{0}}=1=f_{t_{0}}(x_{0})\,$for all $t$.

We therefore obtain by the same method than that of section 4.3: 
\[
\stackunder{t\rightarrow 0}{\overline{\lim }}\frac{A_{t}^{1}}{%
\int_{B(0,\delta )}\overline{u}_{t}^{2}dv_{\mathbf{g}_{t}}}\leq -\frac{%
(n-2)(n-4)}{8(n-1)}\frac{\bigtriangleup _{\mathbf{g}}f_{t_{0}}(x_{0})\,}{%
f_{t_{0}}(x_{0})\,}+\varepsilon _{\delta } 
\]
and thus, after letting $\delta $ tend to 0, we obtain: 
\[
h(x_{0})\leq \frac{n-2}{4(n-1)}S_{\mathbf{g}}(x_{0})-\frac{(n-2)(n-4)}{8(n-1)%
}\frac{\bigtriangleup _{\mathbf{g}}f_{t_{0}}(x_{0})\,}{f_{t_{0}}(x_{0})\,} 
\]
But
\[
\bigtriangleup _{\mathbf{g}}f_{t}(x_{0})\,\sim +\frac{c}{t^{2}}\stackunder{%
t\rightarrow 0}{\rightarrow }+\infty 
\]
so taking $t_{0}$ close to 0 we obtain a contradiction.

This proves that we can find in the sequence ($f_{t}$) functions
whith laplacian in $x_{0}$, $\bigtriangleup _{\mathbf{g}}f_{t}(x_{0})$, 
as large as we want such that the equations: $\bigtriangleup _{%
\mathbf{g}}u+h.u=f_{t}.u^{\frac{n+2}{n-2}}$ do \textit{not} have minimizing solutions 
and therefore such that $h$ is weakly critical for $f_{t}$ and $\mathbf{g}$.
\\

Remark 1: We also have in this setting the analog of theorem 6 on the speed of convergence of $(x_{t})$ 
to $x_{0}$.

Remark 2: this can be apply to $h=cste<B_{0}K^{-2}$ or to $h=S_{\mathbf{g}}$ if $M$ is not the sphere.
\\

\textit{Second step:}

For our function $h$ such that $(h,1,\mathbf{g})$ is subcritical,
we know now that there exists a function $f$, with a laplacian as large as we want at its maximum points, such that $(h,f,\mathbf{g})$
is weakly critical. More precisely, we found a function $f$ such that:

1/: $(h,f,\mathbf{g})$ is weakly critical,

2/: $h(x_{0})>\frac{n-2}{4(n-1)}S_{\mathbf{g}}(x_{0})-\frac{(n-2)(n-4)}{%
8(n-1)}\frac{\bigtriangleup _{\mathbf{g}}f(x_{0})\,}{f(x_{0})\,}$ where

a/: $h(x_{0})>0$

b/: \{$x_{0}\}=\{x\,/\,f(x)=\stackunder{M}{Sup}f\}$ and $f(x_{0})=1$, $0\leq
f\leq 1$, $Supp\,f=B(x_{0},r)$

c/: $\nabla ^{2}f(x_{0})<0$ .

We now consider the path
\[
t\rightarrow f_{t}=(1-t).1+t.f. 
\]
Remark that for all $t$: $\bigtriangleup _{\mathbf{g}}f_{t}=t\bigtriangleup
_{\mathbf{g}}f$ and $f_{t}(x_{0})=1=\stackunder{M}{Sup}\,f_{t}$. We set
\[
\lambda _{t}=Inf\,J_{h,f_{t},\mathbf{g}}. 
\]
Then

\[
\lambda _{0}<K(n,2){{}^{-2}}(\stackunder{M}{Sup}f_{0})^{-\frac{n-2}{n}} 
\]
because ($h,1,\mathbf{g})$ is sub-critical and

\[
\lambda _{1}=K(n,2){{}^{-2}}(\stackunder{M}{Sup}f_{1})^{-\frac{n-2}{n}} 
\]
as ($h,f,\mathbf{g})$ is weakly critical. Remark that $\stackunder{M%
}{Sup}\,f_{t}\,$is always equal to 1.

Let
\[
t_{0}=Sup\{t\,/\,\lambda _{t}<K(n,2){{}^{-2}}(\stackunder{M}{Sup}f_{t})^{-%
\frac{n-2}{n}}\} 
\]
Then $0<t_{0}\leq 1$ and
\[
\lambda _{t_{0}}=K(n,2){{}^{-2}}(\stackunder{M}{Sup}f_{t_{0}})^{-\frac{n-2}{n%
}} 
\]
Before applying the method of section 4.3, we need to prove one more thing : as $h$ is weakly critical for $f_{t_{0}}$, 
we know that at the maximum point $x_{0}$ we have
\[
h(x_{0})\geq \frac{n-2}{4(n-1)}S_{\mathbf{g}}(x_{0})-\frac{(n-2)(n-4)}{8(n-1)%
}\frac{\bigtriangleup _{\mathbf{g}}f_{t_{0}}(x_{0})\,}{f_{t_{0}}(x_{0})\,} 
\]
because $\frac{\bigtriangleup _{\mathbf{g}}f_{t_{0}}(x_{0})\,}{f_{t_{0}}(x_{0})\,%
}=t_{0}\frac{\bigtriangleup _{\mathbf{g}}f(x_{0})\,}{f(x_{0})\,}$ with $%
t_{0}\leq 1$, but we need a strict inequality.

We consider the sequence ($f_{i})$, that we can construct using the first step: $f_{i}$ is such that 
($h,f_{i},\mathbf{g}$) is weakly critical with
\[
f_{i}(x_{0})=1=Supf_{i}\,\,\,et\,\,\,\bigtriangleup _{\mathbf{g}%
}f_{i}(x_{0})\rightarrow +\infty . 
\]
For each $f_{i}$, we note $t_{i}$ the ''$t_{0}$'' built above.
Therefore for any $i$ :\thinspace 
\[
h\text{ is weakly critical for }(1-t_{i}).1+t_{i}.f_{i}\text{ and }%
\mathbf{g}. 
\]
Suppose that liminf $t_{i}=0$, or, after extracting, that $%
t_{i}\rightarrow 0.$ Then, 
\[
(1-t_{i}).1+t_{i}.f_{i}\rightarrow 1 
\]
uniformly on $M$ as $0\leq f_{i}\leq 1$. But ($h,1,\mathbf{g})$ is sub-critical, thus there exists $u\in H_{1}^{2}(M)$ such that
\[
\frac{\int \left| \nabla u\right| ^{2}+\int hu^{2}}{(\int u^{2^{*}})^{\frac{2%
}{2^{*}}}}<K(n,2){{}^{-2}\,.} 
\]
But then
\[
\frac{\int \left| \nabla u\right| ^{2}+\int hu^{2}}{(\int
((1-t_{i}).1+t_{i}.f_{i})u^{2^{*}})^{\frac{2}{2^{*}}}}\rightarrow \frac{\int
\left| \nabla u\right| ^{2}+\int hu^{2}}{(\int u^{2^{*}})^{\frac{2}{2^{*}}}}%
<K(n,2){{}^{-2}} 
\]
whereas
\[
K(n,2){{}^{-2}=}K(n,2){{}^{-2}}(\stackunder{M}{Sup}%
((1-t_{i}).1+t_{i}.f_{i}))^{-\frac{n-2}{n}} 
\]
which contradict the fact that $(h,(1-t_{i}).1+t_{i}.f_{i},\mathbf{g})$ is weakly critical.

Therefore, up to extraction, $t_{i}\rightarrow t_{1}>0$

As $\bigtriangleup _{\mathbf{g}}f_{i}(x_{0})\rightarrow +\infty $, we can find $i$ 
large enough so that
\[
\frac{(n-2)(n-4)}{8(n-1)}t_{i}\frac{\bigtriangleup _{\mathbf{g}%
}f_{i}(x_{0})\,}{f_{i}(x_{0})\,}>\frac{n-2}{4(n-1)}S_{g}(x_{0})-h(x_{0})\,. 
\]
If we now denote $f$ this last function $f_{i}$ and $t_{0}$ this
\thinspace $t_{i}$, we get a path
\[
t\rightarrow f_{t}=(1-t).1+t.f 
\]
such that :

a/: $\forall t<t_{0}$ : $(h,f_{t},\mathbf{g})$ is sub-critical,

b/: ($h,f_{t_{0}},\mathbf{g})$ is weakly critical with :

b1/: \{$x_{0}\}=\{x\,/\,f_{t}(x)=\stackunder{M}{Sup}f_{t}\}$ and $%
f_{t}(x_{0})=1$ for all $t$

b2/: $h(x_{0})>\frac{n-2}{4(n-1)}S_{\mathbf{g}}(x_{0})-\frac{(n-2)(n-4)}{%
8(n-1)}\frac{\bigtriangleup _{\mathbf{g}}f_{t_{0}}(x_{0})\,}{%
f_{t_{0}}(x_{0})\,}$

b3/: $\nabla ^{2}f_{t_{0}}(x_{0})<0$

For any $t<t_{0}$ there exists a minimizing solution $u_{t}$ of the equation
\[
\bigtriangleup _{\mathbf{g}}u_{t}+h.u_{t}=\lambda _{t}.f_{t}.u_{t}^{\frac{n+2%
}{n-2}} 
\]
with $\int f_{t}u_{t}^{2^{*}}=1$. The sequence ($u_{t}$) is bounded in $%
H_{1}^{2}$ therefore 
\[
u_{t}\stackrel{H_{1}^{2}}{\stackunder{t\rightarrow t_{0}}{\rightharpoondown }%
}u 
\]
and we are once again in the situation where :

- either $u>0$ and then $u$ is a minimizing solution of $\bigtriangleup _{\mathbf{g}%
}u+h.u=\lambda _{t_{0}}f_{t_{0}}.u^{\frac{n+2}{n-2}}$, and therefore 
($h,f_{t_{0}},\mathbf{g})$ is critical.

- either $u\equiv 0$ $\,$ and once again the sequence ($u_{t}$) concentrates.
In this case, the sudy of the concentration phenomenom is easier than in the first step as the family $(f_{t})$ tend uniformly to $f$
when $t\rightarrow t_{0}$ with $Supp\,f_{t}=B(x_{0},r)$. We can find $%
\delta <r $ such that $f>0$ on $B(x_{0},\delta )$. Then there exists $c>0$ such that for any $t$ we have: 
\[
0<c\leq f_{t}\leq 1\text{ on }B(x_{0},\delta ), 
\]
Furthermore, the $f_{t}$ all reach their maximum at $x_{0}$, this maximum being always 1. 
We can then go over all the results and methods of section 4.3, 
the functions $f_{t}$ bringing this time no changes. We finally obtain
\[
h(x_{0})\leq \frac{n-2}{4(n-1)}S_{g}(x_{0})-\frac{(n-2)(n-4)}{8(n-1)}\frac{%
\bigtriangleup _{\mathbf{g}}f_{t_{0}}(x_{0})\,}{f_{t_{0}}(x_{0})\,} 
\]
thus a contradiction. Therefore ($h,f_{t_{0}},\mathbf{g})$ is critical with a minimizing solution.

This proof in fact shows the following result:

\textbf{Theorem 4':}

\textit{If }$h$\textit{\ is weakly critical for a function }$f$%
\textit{\ and a metric }$\mathbf{g},$ \textit{these datas satisfying:}

\textit{1/: }$h(x)>\frac{n-2}{4(n-1)}S_{\mathbf{g}}(x)-\frac{(n-2)(n-4)}{%
8(n-1)}\frac{\bigtriangleup _{\mathbf{g}}f(x)\,}{f(x)\,}$\textit{\ at the maximum points of }$f$

\textit{2/:}$\nabla ^{2}f(x)<0$\textit{\ at the maximum points of }$f$

\textit{3/: there exists a sequence }$f_{t}\stackrel{C^{2}}{\stackunder{%
t\rightarrow t_{0}}{\rightarrow }}f$\textit{\ with }$\stackunder{M}{Sup}%
f_{t}=\stackunder{M}{Sup}f$\textit{\ such that (}$h,f_{t},\mathbf{g})$%
\textit{\ is subcritical for }$t<t_{0}$

\textit{then (}$h,f,\mathbf{g})$\textit{\ is critical and has minimizing solutions.}

As we said in the introduction, this leads to another, dual, definition of critical functions, that is definition 3.
The natural question is then

\begin{center}
\textit{Is }$f$ \textit{critical for }$h$\textit{\ if and only if }$h$ 
\textit{is critical for }$f$ ?
\end{center}

Remark that in both cases, if $P$ is a point where $f$ is maximum on $M$\ :\textit{\ }$\frac{4(n-1)}{n-2}%
h(P)\geqslant S_{\mathbf{g}}(P)-\frac{n-4}{2}\frac{\bigtriangleup _{\mathbf{g%
}}f(P)}{f(P)}$\textit{\ }
\\

This problem seeems difficult. We prove here the result we obtain, theorem 5.

The proof starts with the following remark:
We have seen that if $h$ is weakly critical for $f$ and $\mathbf{g}$ and that $%
\bigtriangleup _{\mathbf{g}}u+h.u=f.u^{\frac{n+2}{n-2}}$ has a minimizing solution, then $h$ is 
critical for $f$ and $\mathbf{g}$. In the same way,
if $f$ is weakly critical for $h$ (in the sense that $\lambda _{h,f,%
\mathbf{g}}=K(n,2){{}^{-2}}(\stackunder{M}{Sup}f)^{-\frac{n-2}{n}}$ ) and if 
$\bigtriangleup _{\mathbf{g}}u+h.u=f.u^{\frac{n+2}{n-2}}$ has a minimizing solution $u>0$, 
then $f$ is critical for $h$. Indeed, if $f^{\prime
} $ is a function such that $Supf=Supf^{\prime }$ and $f^{\prime
}\gneqq f$, we have
\[
\int f^{\prime }u^{2^{*}}>\int fu^{2^{*}} 
\]
because $u>0$. Therefore
\[
J_{h,f^{\prime },\mathbf{g}}(u)<J_{h,f,\mathbf{g}}(u)=K(n,2){{}^{-2}}(%
\stackunder{M}{Sup}f)^{-\frac{n-2}{n}}=K(n,2){{}^{-2}}(\stackunder{M}{Sup}%
f^{\prime })^{-\frac{n-2}{n}}\,. 
\]
Using our work of section 4.3 and of this section, the proof is now short:

-If $h$ is critical for $f,$ we apply theorem 1: $\bigtriangleup _{\mathbf{g}}u+h.u=f.u^{\frac{n+2}{n-2}}$ 
has a minimizing solution, and therefore $f$ is critical for $h$.

-If $f$ is critical for $h$, these two functions (and the metric)
satisfying the hypothesis of the theorem, we have $\lambda _{h,f,%
\mathbf{g}}=K(n,2){{}^{-2}}(\stackunder{M}{Sup}f)^{-\frac{n-2}{n}},$ so $h$
is weakly critical for $f$. We then consider, for $t\stackrel{<}{%
\rightarrow }1$, the sequence
\[
t\rightarrow f_{t}=(1-t).Supf+t.f\,. 
\]
For allt $t:\,$we have $Supf_{t}=Supf$ and if $t<1$ then $f_{t}\gneqq f$. 
Therefore as $f$ is critical for $h$, by definition:
\[
\lambda _{h,f_{t},\mathbf{g}}<K(n,2){{}^{-2}}(\stackunder{M}{Sup}f_{t})^{-%
\frac{n-2}{n}}\,. 
\]
We then apply theorem 4' above to obtain that $h$ is critical for $f$ with minimizing solutions.

\section{The case of the dimension 3; ending remarks}
\subsection{The case of the dimension 3.}
We just state the results in the case of dimension 3, as they are immediate generalisations of results of O. Druet 
proved in the case where $f$ is a constant; we refer to his article for the proofs [10].
The dimension 3 requires fondamentaly the use of the Green function. We refer to the proof of proposition 8 in 
section 4.2 for the definition and the property of the Green function.
In dimension 3, for any point $x\in M$, and for $y$ close to $x$, $G_{h}$
can be writen in the following way: 
\[
G_{h}(x,y)=\frac{1}{\omega _{2}d_{\mathbf{g}}(x,y)}+M_{h}(x)+o(1)
\]
where $o(1)$ is to be taken for $y\rightarrow x$. We call $M_{h}(x)$
the mass of the Green function at $x$.

The generalisation of the results of O. Druet to the case of an arbitrary function $f$ in $\E$ gives the following:

Let $(M,\g)$ be a compact manifold of dimension 3, and let $f\in C^{\infty }(M)$ be such that $Sup f >0$.
We have the following results:
\begin{itemize}
\item For any function $h$ weakly critical for $f$ and $\g$, and for any $x\in Max f$, we have $M_{h}(x)\leq 0$.
\item For any $h\in C^{\infty }(M)$, let $B(h)=\inf \{B/$ $h+B\,\,is\,\,weakly\,\,critical\,\,for\,\,f\}$. 
Then $h+B(h)$is a critical function for $f$.
\item Let $h$ be a critical function for $f$ and $\g$. Then one of the following condition is true:
\begin{enumerate}
\item There exists $x\in Max f$ such that $M_{h}(x)=0$.
\item ($\E$) has minimizing solutions.
\end{enumerate}
\end{itemize}

Remarks:

-The condition
\[
M_{h}(x)\leqslant 0 
\]
appears as the analog of the condition
\[
\frac{4(n-1)}{n-2}h(P)\geqslant S_{\mathbf{g}}(P)-\frac{n-4}{2}\frac{%
\bigtriangleup _{\mathbf{g}}f(P)}{f(P)} 
\]
we had in dimension $\geqslant 4$. In the case $f=cst$, this condition must be satisfied on all of $M$.

-The particularity of dimension 3 is to offer critical functions of any shape, that is the meaning of the second point.

-The main difference with the case $f=cst$ studied by O. Druet is that the conditions on the mass of the Green 
function are to be considered only at the point of maximum of $f$.

\subsection{Degenerate hessian at the point of maximum and fundamental estimate.}
In theorem 6, we made the hypothesis that the hessian of $f$ is non degenerate at each of its points of maximum. 
We give here a conterexample to show that this hypothesis is necessary.
Consider the $n$-dimensional sphere $S^{n}$ with its standard metric $\mathbf{s}$.
Rewriting known results (c.f. for example [17]), there exists a unique critical function for 1 et $\mathbf{s}$,
which is
\[
h=\frac{n-2}{4(n-1)}S_{\mathbf{s}}=\frac{n-2}{4(n-1)} 
\]
and this critical function has only two type of extremal functions, the constants and the functions of the form
\[
u=a(b-\cos r)^{-\frac{n-2}{2}} 
\]
where $a\neq 0$, $b>1$, and $r$ is the geodesic distance to some fixed point of $S^{n}$. 
Consider now on $S^{n}$ a sequence of points $x_{t}$ converging to a poit $x_{0}$, and let 
\[
u_{t}=\mu _{t}^{\frac{n-2}{2}}(\mu _{t}^{2}+1-\cos r_{t})^{-\frac{n-2}{2}} 
\]
where $r_{t}(x)=d_{\mathbf{s}}(x,x_{t})$ and $\mu _{t}$ is a sequence of real converging to 0. Then 
\[
\int_{M}u_{t}^{2^{*}}dv_{\mathbf{s}}=1 
\]
and we obtain in this way a sequence of solutions of the equation 
\[
\triangle _{\mathbf{s}}u_{t}+\frac{n-2}{4(n-1)}.u_{t}=K(n,2)^{-2}u_{t}^{%
\frac{n+2}{n-2}} 
\]
where obviously the function $f=K(n,2)^{-2}$ has degenerate hessian at its maximum points ! Furthermore
\[
Sup_{M}u_{t}=u_{t}(x_{t})=\mu _{t}^{-\frac{n-2}{2}}. 
\]
This sequence concentrates and satisfies all the propositions 2 to 9 seen in section 4.2, whatever the choice of 
the sequence $x_{t}\rightarrow x_{0}$ and of the sequence $\mu _{t}\rightarrow 0$. By spherical symetry, we can 
easily find two sequences ($x_{t}$) and ($\mu _{t}$) such that
\[
\frac{d_{\mathbf{s}}(x_{t},x_{0})}{\mu _{t}}\rightarrow +\infty 
\]
by taking for example $\mu _{t}=d_{\mathbf{s}}(x_{t},x_{0})^{2}.$

Once again, it seems that the hypothesis on the hessian of $f$ ''fixes'' the position of the concentration point, 
and so imposes a speed of convergence of the sequence ($x_{t}$).
\subsection{Further questions.}
First a remark concerning the requirement of a strict inequality at the point of maximum of $f$ in theorem 1.
An easy but somewhat artificial extension of a result of hebey and Vaugon is the following:

\textit{Suppose that the manifold (}$M,\mathbf{g)}$\textit{\ is of dimension }$\geq 7$\textit{, 
and let (}$h,f,\mathbf{g)}$\textit{\ be a critical triple. Let }$T_{f}=\{x\in M/\,f(x)=Maxf\,\,and\,\,h(x)=\frac{n-2}{4(n-1)%
}S_{\mathbf{g}}(x)-\frac{(n-2)(n-4)}{8(n-1)}\frac{\bigtriangleup _{\mathbf{g}%
}f(P)}{f(P)}\}$\textit{. We suppose that }$T_{f}$\textit{\ is not dense in }$M$\textit{\ 
and that for any point }$x$\textit{\ of }$T_{f}$\textit{:}

\textit{1:The Weyl tensor vanishes on a neihbourhood of }$x$\textit{, }

\textit{2: }$\nabla ^{2}(h-\frac{n-2}{4(n-1)}S_{\mathbf{g}})$\textit{\ is not degenerate in }$x$\textit{,}

\textit{3: }$\bigtriangleup _{\mathbf{g}}f(x)=0$ if $x\in T_{f}$, \textit{and we suppose that  }$f$
\textit{\ is non degeneraye at the points of maximum which are not in }$T_{f}$\textit{.}

\textit{Then (}$h,f,\mathbf{g)}$\textit{\ has minimizing solutions.}

The main interest of this result is that we can expect existence of solutions in this case. Looking to our method, 
it seems that one need to find some other intrinsec parameters, i.e. invariant by the exponential charts $exp_{x_{t}}$.
See our thesis for more precision.
\\

Another question is the following: We saw that the study of equations $\bigtriangleup _{\mathbf{g}%
}u+hu=fu^{2^{*}-1}$ is linked to the study of the best constants in the Sobolev inclusions of $H_{1}^{2}$ in 
$L^{\frac{2n}{n-2}}$. In the same way, the study of the Sobolev inclusions of $H_{1}^{p}$
in $L^{\frac{pn}{n-p}}$, where $\frac{pn}{n-p}$ is the critical exponent,
and of the associated best constants, goes through the study of equations of the form
\[
\bigtriangleup _{p}u+hu=fu^{\frac{pn}{n-p}-1}
\]
where $\bigtriangleup _{p}u=-\nabla (\left| \nabla u\right| _{\mathbf{g}%
}^{p-2}\nabla u)$ is the p-laplacian; see for example O. Druet, E. Hebey and Z. Faget [F2]. Here also variational methods are used :
 the functional used is :
\[
I(u)=\int \left| \nabla u\right| _{\mathbf{g}}^{p}+\int hu^{p}
\]
from where we see the link with the Sobolev inclusion
\[
(\int u^{\frac{pn}{n-p}})^{\frac{n-p}{n}}\leq K(n,p)\int \left| \nabla
u\right| _{\mathbf{g}}^{p}+B\int u^{p}
\]
where $K(n,p)$ is the associated best constant. The starting point is again the following : If
\[
\stackunder{\int u^{\frac{pn}{n-p}}=1}{Inf\,}\,I(u)<K(n,p)^{-1}(Supf)^{-%
\frac{n-p}{n}}
\]
then the equation has a minimizing solution $u>0$ (knowing that the large inequality is always true). 
We therefore see that it is easy to extend the definition of critical functions to this case. It would therefore be 
interesting to know if our results can be extended to this setting.

Another question that can be asked after our work is the following :

\begin{center}
$f$\textit{ being given, is there constant critical functions ?}
\end{center}

This would give some kind of  ''best second constant $B_{0}(\mathbf{g},f)$'' linked to $f$.

At last, there is a question which emerges from our work:

\begin{center}
\textit{For a given arbitrary function }$h$\textit{\ on }$M$\textit{,
does there exist solutions (not minimizing) to the equation }$\bigtriangleup _{%
\mathbf{g}}u+hu=fu^{2^{*}-1}$\textit{\ ?}
\end{center}

Indeed, we saw that this equation has (minimizing) solutions when $h$ is sub-critical and when $h$ is 
critical with some hypothesis. However, variational methods do not give any answer is $h$ larger and different than 
some critical function, or if $\bigtriangleup _{\mathbf{g}}+h$ is not coercive. 
In this cases, if solutions exist, they cannot be minimizing. One therefore needs other methods for these cases. 
See [4] who study the case $f=cst$ and $3\leq dimM\leq6$.
\\

\chapter{Appendice A: d\'{e}monstrations de quelques propri\'{e}t\'{e}s}
\pagestyle{myheadings}\markboth{\textbf{Appendice A.}}{\textbf{Appendice A.}}
Nous donnons ici quelques d\'{e}monstrations des propri\'{e}t\'{e}s
\'{e}l\'{e}mentaires des fonctions critiques que nous avons expos\'{e}es au
chapitre 1. Nous les avons report\'{e}es ici pour rendre plus lisible ce
premier chapitre, et parce qu'elles sont tr\`{e}s simples (bien
qu'importantes).
\\

\textit{L'op\'{e}rateur }$\bigtriangleup _{\mathbf{g}}+h$\textit{\ est
forc\'{e}ment coercif pour }$h\in C{{}^{0}}(M)$\textit{\ si }$\lambda _{h,f,%
\mathbf{g}}=\frac{1}{K(n,2){{}^{2}}(\stackunder{M}{Sup}f)^{\frac{n-2}{n}}}$%
\textit{\ et si }$f>0$\textit{\ sur }$M$\textit{.}
\\

En effet, par d\'{e}finition 
\begin{eqnarray*}
\int_{M}\left| \nabla u\right| _{\mathbf{g}}^{2}dv_{g}+\int_{M}h.u{{}^{2}}%
dv_{\mathbf{g}} && \geq \lambda _{h,f,\mathbf{g}}\left( \int_{M}f\left|
u\right| ^{\frac{2n}{n-2}}dv_{\mathbf{g}}\right) ^{\frac{n-2}{n}} \\ 
&&\geq \lambda _{h,f,\mathbf{g}}.\stackunder{M}{Inf}f.\left( \int_{M}\left|
u\right| ^{\frac{2n}{n-2}}dv_{\mathbf{g}}\right) ^{\frac{n-2}{n}} \\ 
&&\geq \lambda _{h,f,\mathbf{g}}.\stackunder{M}{Inf}f.c.\left(
\int_{M}\left| u\right| ^{2}dv_{\mathbf{g}}\right)
\end{eqnarray*}
la derni\`{e}re in\'{e}galit\'{e} \'{e}tant obtenue par l'in\'{e}galit\'{e}
de H\"{o}lder.
\\

\textit{Si }$\bigtriangleup _{\mathbf{g}}+h$\textit{\ est coercif, toute
fonction }$k$\textit{\ continue assez proche de }$h$\textit{\ dans }$%
C^{0}(M) $\textit{\ est telle que }$\bigtriangleup _{\mathbf{g}}+k$\textit{\
est coercif.}
\\

Simplement par continuit\'{e} de l'int\'{e}grale.
\\

\textit{Si }$h$\textit{\ est critique, n\'{e}cessairement il existe }$x\in M$%
\textit{\ tel que }$h(x)>0$\textit{\ .}
\\

Cela d\'{e}coule de l'exigence que nous avons faite que $\Delta _{\mathbf{g}%
}+h$ soit coercif. Il suffit d'appliquer la coercivit\'{e} \`{a} la fonction
1: 
\begin{eqnarray*}
\int_{M}hdv_{\mathbf{g}}=\int_{M}\left| \nabla 1\right| _{\mathbf{g}%
}^{2}dv_{g}+\int_{M}h.1{{}^{2}}dv_{\mathbf{g}} & &\geq \int_{M}1^{2}dv_{%
\mathbf{g}} \\ 
& &>0
\end{eqnarray*}
\\

\textit{Par d\'{e}finition }$B_{0}(g)K(n,2){{}^{-2}}$\textit{\ est toujours
une fonction (constante) faiblement critique pour toute fonction }$f$\textit{%
\ et toute m\'{e}trique }$\mathbf{g}$\textit{. Si de plus }$f\equiv 1$%
\textit{, que }$(M,\mathbf{g})$\textit{\ n'est pas conform\'{e}ment
diff\'{e}omorphe \`{a} la sph\`{e}re standard et est telle que }$S_{\mathbf{g%
}}=cste$\textit{, }$B_{0}(\mathbf{g})K(n,2)^{-2}$\textit{\ est une fonction
critique, c'est la plus petite fonction critique constante.}
\\

D'apr\`{e}s la r\'{e}solution du probl\`{e}me de Yamabe, on a avec ces
hypoth\`{e}ses: 
\[
\lambda _{\mathbf{g}}:=\stackunder{H_{1}^{2}-\{0\}}{Inf}\frac{\int_{M}\left|
\nabla u\right| _{\mathbf{g}}^{2}dv_{g}+\int_{M}\frac{n-2}{4(n-1)}S_{\mathbf{%
g}}.u{{}^{2}}dv_{\mathbf{g}}}{\left( \int_{M}\left| u\right| ^{\frac{2n}{n-2}%
}dv_{\mathbf{g}}\right) ^{\frac{n-2}{n}}}<K(n,2)^{-2} 
\]
($\lambda _{\mathbf{g}}$ est ''l'invariant de Yamabe''). Donc il existe $u>0$
solution de 
\[
\triangle _{\mathbf{g}}u+h.u=\lambda _{\mathbf{g}}.u^{\frac{n+2}{n-2}}\text{
avec }\int_{M}u^{2^{*}}dv_{\mathbf{g}}=1\,. 
\]
Alors 
\begin{eqnarray*}
1=(\int_{M}u^{2^{*}}dv_{\mathbf{g}})^{\frac{2}{2^{*}}} & \leq&
K(n,2)^{2}\int_{M}\left| \nabla u\right| _{\mathbf{g}}^{2}dv_{g}+B_{0}(%
\mathbf{g)}\int_{M}u^{2}dv_{\mathbf{g}} \\ 
& \leq& K(n,2)^{2}(\int_{M}\left| \nabla u\right| _{\mathbf{g}%
}^{2}dv_{g}+\int_{M}\frac{n-2}{4(n-1)}S_{\mathbf{g}}.u{{}^{2}}dv_{\mathbf{g}%
}) \\ 
&&+(B_{0}(\mathbf{g)-}\frac{n-2}{4(n-1)}S_{\mathbf{g}}.K(n,2)^{2})%
\int_{M}u^{2}dv_{\mathbf{g}}\\ 
& \leq& K(n,2)^{2}\lambda _{\mathbf{g}}+(B_{0}(\mathbf{g)-}\frac{n-2}{4(n-1)}%
S_{\mathbf{g}}.K(n,2)^{2})\int_{M}u^{2}dv_{\mathbf{g}}\,.
\end{eqnarray*}
Or 
\[
1-K(n,2)^{2}\lambda _{\mathbf{g}}>0 
\]
donc 
\[
B_{0}(\mathbf{g)>}\frac{n-2}{4(n-1)}S_{\mathbf{g}}.K(n,2)^{2} 
\]
et donc d'apr\`{e}s le th\'{e}or\`{e}me de Z. Djadli et O. Druet cit\'{e} au
chapitre 1, (1.2) a des fonctions extr\'{e}males. Mais cela
signifie exactement que la fonction faiblement critique $B_{0}(\mathbf{g}%
)K(n,2)^{-2}$ a des solutions minimisantes, elle est donc critique.
\\

\textit{Les fonctions critiques se transforment dans les changements de
m\'{e}trique conformes exactement comme la courbure scalaire.}
\\

En effet, soit $u\in C^{\infty }(M),\,u>0$ et $\mathbf{g}^{\prime }=u^{\frac{%
4}{n-2}}\mathbf{g}$ une m\'{e}trique conforme \`{a} $\mathbf{g}$\textbf{.}
Posons 
\[
h^{\prime }=\frac{\triangle _{\mathbf{g}}u+h.u}{u^{\frac{n+2}{n-2}}}\text{
c'est \`{a} dire }\triangle _{\mathbf{g}}u+h.u=h^{\prime }.u^{\frac{n+2}{n-2}%
}\,. 
\]
Pour montrer que pour toute fonction $w\in H_{1}^{2}(M)-\{0\}:$%
\[
J_{h,f,\mathbf{g}}(w)=J_{h^{\prime },f,\mathbf{g}^{\prime }}(u^{-1}.w) 
\]
il suffit de remarquer que: 
\[
dv_{\mathbf{g}^{\prime }}=u^{\frac{2n}{n-2}}.dv_{\mathbf{g}}=u^{2^{*}}.dv_{%
\mathbf{g}} 
\]
et 
\[
\left| \nabla w\right| _{\mathbf{g}^{\prime }}^{2}=u^{-\frac{4}{n-2}}\left|
\nabla w\right| _{\mathbf{g}}^{2}\,. 
\]
Car alors 
\[
\int_{M}f\left| w\right| ^{2^{*}}dv_{\mathbf{g}}=\int_{M}f\left| w\right|
^{2^{*}}u^{-2^{*}}dv_{\mathbf{g}^{\prime }}=\int_{M}f\left| \frac{w}{u}%
\right| ^{2^{*}}dv_{\mathbf{g}^{\prime }} 
\]
et 
\[
\int_{M}h^{\prime }(\frac{w}{u})^{2}dv_{\mathbf{g}^{\prime
}}=\int_{M}h^{\prime }w^{2}u^{2^{*}-2}dv_{\mathbf{g}}=\int_{M}h^{\prime
}w^{2}u^{\frac{4}{n-2}}dv_{\mathbf{g}}\,. 
\]
Enfin, par int\'{e}gration par parties 
\[
\int_{M}\left| \nabla (\frac{w}{u})\right| _{\mathbf{g}^{\prime }}^{2}dv_{%
\mathbf{g}^{\prime }}=\int_{M}\left| \nabla w\right| _{\mathbf{g}}^{2}dv_{%
\mathbf{g}}-\int_{M}\Delta _{\mathbf{g}}u.w^{2}.u^{-1}dv_{\mathbf{g}}\,. 
\]
Ainsi 
\begin{eqnarray*}
I_{h^{\prime },\mathbf{g}^{\prime }}(\frac{w}{u}) & &:=\int_{M}\left| \nabla (%
\frac{w}{u})\right| _{\mathbf{g}^{\prime }}^{2}dv_{\mathbf{g}^{\prime
}}+\int_{M}h^{\prime }(\frac{w}{u})^{2}dv_{\mathbf{g}^{\prime }} \\ 
&& =\int_{M}\left| \nabla w\right| _{\mathbf{g}}^{2}dv_{\mathbf{g}%
}+\int_{M}(h^{\prime }u^{\frac{4}{n-2}}-\frac{\Delta _{\mathbf{g}}u}{u}%
).w^{2}dv_{\mathbf{g}} \\ 
& &=I_{h,\mathbf{g}}(w)
\end{eqnarray*}
avec 
\[
h=h^{\prime }u^{\frac{4}{n-2}}-\frac{\triangle _{\mathbf{g}}u}{u} 
\]
et donc 
\[
J_{h,f,\mathbf{g}}(w)=J_{h^{\prime },f,\mathbf{g}^{\prime }}(u^{-1}.w)\,. 
\]
De plus si $\bigtriangleup _{\mathbf{g}}+h$ est coercif $\forall w\in
H_{1}^{2}(M)$%
\[
\int_{M}\left| \nabla w\right| _{\mathbf{g}}^{2}dv_{\mathbf{g}%
}+\int_{M}hw^{2}dv_{\mathbf{g}}\geq c\int_{M}w^{2}dv_{\mathbf{g}} 
\]
et alors $\forall w\in H_{1}^{2}(M)$%
\begin{eqnarray*}
\int_{M}\left| \nabla w\right| _{\mathbf{g}^{\prime }}^{2}dv_{\mathbf{g}%
^{\prime }}+\int_{M}h^{\prime }w^{2}dv_{\mathbf{g}^{\prime }}=I_{h,\mathbf{g}%
}(uw) & &\geq c\int_{M}u^{2}w^{2}dv_{\mathbf{g}} \\ 
& &\geq \frac{c}{\stackunder{M}{Sup\,}u^{2^{*}-2}}\int_{M}w^{2}.u^{2^{*}}dv_{%
\mathbf{g}} \\ 
& &\geq \frac{c}{\stackunder{M}{Sup\,}u^{2^{*}-2}}\int_{M}w^{2}dv_{\mathbf{g}%
^{\prime }}
\end{eqnarray*}
et donc $\bigtriangleup _{\mathbf{g}^{\prime }}+h^{\prime }$ est coercif. En
\'{e}changeant $h$ et $\mathbf{g}$ avec $h^{\prime }$ et $\mathbf{g}^{\prime
}$, on voit que $\bigtriangleup _{\mathbf{g}^{\prime }}+h^{\prime }$ est
coercif si et seulement si $\bigtriangleup _{\mathbf{g}}+h$ est coercif.

Ceci implique que $h$ est critique pour $f$ et $\mathbf{g}$ si et seulement
si 
\[
h^{\prime }=\frac{\triangle _{\mathbf{g}}u+h.u}{u^{\frac{n+2}{n-2}}} 
\]
est critique pour $f$ et $\mathbf{g}^{\prime }=u^{\frac{4}{n-2}}\mathbf{g}.$
Ou d'une autre mani\`{e}re, $h^{\prime }$ est critique pour $f$ et $\mathbf{g%
}^{\prime }=u^{\frac{4}{n-2}}\mathbf{g}$ si et seulement si 
\[
h=h^{\prime }u^{\frac{4}{n-2}}-\frac{\triangle _{\mathbf{g}}u}{u} 
\]
est critique pour $f$ et $\mathbf{g}$. (Ceci est aussi valable pour les
fonctions sous-critiques et faiblement critiques). De plus, $w$ est une
solution minimisante pour ($h,f,\mathbf{g}$) si et seulement si $\frac{w}{u}$
est une solution minimisante pour ($h^{\prime },f,\mathbf{g}^{\prime }$).
\\

Dans le chapitre 3 §3.2 (mise en place), nous avons affirmé que si $u_{t}\rightharpoondown u>0$, 
alors $u$ est une solution minimisante. On montre en fait que $u_{t}\rightarrow u$ dans 
 $L^{2^{*}}(M)$. On utilise pour cela le principe d'itération 2.3 et une définition légérement 
 modifiée mais équivalente du point de concentration: $x$ est un point de concentration si 
 $\exists \varepsilon >0$ tel que $\forall \delta >0$ 
$\stackunder{t\rightarrow1}{\lim \sup }\int_{B(x,\delta
)}u_{t}^{2^{*}}dv_{\mathbf{g}}\geqslant \varepsilon$. Remarquons que si $u_{t}\rightarrow u$ dans 
 $L^{2^{*}}$ sur un voisinage de $x$ alors $x$ ne peut être un point de concentration. On utilise alors 
 la formule (2.5) de la manière suivante: si pour $\eta$ valant 1 au voisinage de $x$ on a 
 $limsup\,\lambda _{t}K(n,2){{}^{2}}(\stackunder{Supp\,\eta }{Sup}\left| f\right| ).(\int_{Supp\,\eta }u_{t}^{2^{*}})^{\frac{2%
}{n}}  <1$ alors, pour $k$ et $t$ assez proche de 1, $Q(t,k,\eta )>0$ et donc, avec (2.4), $(\eta u_{t})$ 
est bornée dans $L^{\frac{k+1}{2}2^{*}}$ et on peut en extraire une sous-suite qui 
converge fortement dans $L^{2^{*}}$, donc $u_{t}\rightarrow u$ dans 
 $L^{2^{*}}$ sur un voisinage de $x$. Le principe si $u_{t}\rightharpoondown u>0$ 
 est alors le suivant. S'il n'existe pas de point de concentration, on voit avec cette méthode que pour
 tout $x$, il existe une petite boule $B(x,\delta )$ sur laquelle $u_{t}\rightarrow u$ dans $L^{2^{*}}$. 
Mais alors, on recouvre $M$ par un nombre fini de telles boules, et, par exemple en utilisant une 
partition de l'unité, on obtient que $u_{t}\rightarrow u$ dans $L^{2^{*}}(M)$. S'il existe au moins un 
point de concentration $x$, on commence par montrer, toujours en utilisant (2.3),(2.4) et (2.5),
 que si $f(x)\leqslant 0$, $u_{t}\rightarrow u$ dans $L^{2^{*}}$ au voisinage de $x$, donc 
 que $x$ n'est pas un point de concentration. On montre alors de même que si $x$ est un point de 
 concentration, nécessairement, $\forall \delta >0$, $limsup\,\int_{B(x,\delta)}f.u_{t}^{2^{*}}dv_{\mathbf{g}}=1$, 
 et donc qu'il n'existe qu'un seul point de concentration $x_{0}$. Mais alors on voit que 
 $u_{t}\rightarrow 0$ dans $L_{loc}^{2^{*}}(M-\{x_{0}\})$ ce qui contredit le fait que 
 $u_{t}\rightharpoondown u>0$ dans $L^{2^{*}}(M)$. Par conséquent 
 $u_{t}\rightarrow  u>0$ dans $L^{2^{*}}(M)$ et $\int_{M}f.u^{2^{*}}dv_{\mathbf{g}}=1$. 
 On utilise ensuite le fait que $h$ est faiblement critique pour montrer que $\lambda =\frac{1}{K(n,2){{}^{2}}(\stackunder{M}{Sup}f)^{\frac{n-2}{n}}}$
et donc que $u$ est solution de 
$\triangle _{\mathbf{g}}u+h.u=\frac{1}{K(n,2){{}^{2}}(\stackunder{M}{Sup}f)^{%
\frac{n-2}{n}}}.f.u^{\frac{n+2}{n-2}}$ et $\int_{M}fu^{2^{*}}dv_{\mathbf{g}}=1 $, ce qui 
montre que $u$ est bien une solution minimisante.
 Cette méthode marche aussi bien sûr, quand, comme au chapitre 6, on a une suite $f_{t}\rightarrow f$.

\chapter{Appendice B: construction d'une fonction de Green}
\pagestyle{myheadings}\markboth{\textbf{Appendice B.}}{\textbf{Appendice B.}}
Consid\'{e}rons une vari\'{e}t\'{e} $(M,\mathbf{g})$ compacte, sans bord, de
dimension $n\geq 3.$ Soit $h\in C^{\infty }(M)$ telle que l'op\'{e}rateur 
\[
\bigtriangleup _{\mathbf{g}}+h 
\]
soit coercif. On cherche \`{a} construire une fonction $C^{2}$ 
\[
G_{h}:M\times M\backslash \{(x,x),x\in M\}\rightarrow \Bbb{R} 
\]
strictement positive telle que, au sens des distributions, on ait $\forall
x\in M:$%
\begin{equation}
\bigtriangleup _{\mathbf{g},y}G_{h}(x,y)+h(y)G_{h}(x,y)=\delta _{x}\,. 
\end{equation}
De plus, $G_{h}$ v\'{e}rifiera les estim\'{e}es suivantes: il existe $%
c>0,\,\rho >0$ tels que $\forall (x,y)$ avec $0<d_{\mathbf{g}}(x,y)<\rho :$%
\begin{equation}
\frac{c}{d_{\mathbf{g}}(x,y)^{n-2}}\leq G_{h}(x,y)\leq \frac{c^{-1}}{d_{%
\mathbf{g}}(x,y)^{n-2}} 
\end{equation}

\begin{equation}
\frac{\left| \nabla _{y}G_{h}(x,y)\right| }{G_{h}(x,y)}\geq \frac{c}{d_{%
\mathbf{g}}(x,y)}  
\end{equation}

\begin{center}
$c$ et $\rho $ varient continuement avec $h$
\begin{equation}
G_{h}(x,y)d_{\mathbf{g}}(x,y)^{n-2}\rightarrow \frac{1}{(n-2)\omega _{n-1}}%
\text{ quand }d_{\mathbf{g}}(x,y)\rightarrow 0  
\end{equation}
\end{center}

Beaucoup d'articles utilisent l'existence et les propriétés d'une telle fonction de Green,
mais nous n'avons pas trouvé de r\'{e}f\'{e}rence qui en donnent une d\'{e}monstration prècise.
 Nous proposons donc ici un sch\'{e}ma rapide de
construction de $G_{h}$, ce qui est tr\'{e}s classique, mais surtout
l'obtention des estim\'{e}es (10.2) et (10.3) en utilisant des m\'{e}thodes
proches de celles intervenant dans notre travail.

\subsection{Premi\`{e}re \'{e}tape}

On montre:

\textit{Soit }$\Gamma \in L^{p}(M)$\textit{, }$1<p<\infty $\textit{. Alors
il existe une solution }$\beta \in H_{2}^{p}(M)$\textit{\ de l'\'{e}quation }%
$\Delta _{g}\beta +h\beta =\Gamma $\textit{.}

Si $p\geq 2,$ il suffit d'appliquer les m\'{e}thodes variationnelles
classiques; il faut faire un peu plus attention si $p<2.$ Pour cela, on
commence par utiliser une estim\'{e}e bien connue:

\textit{Si }$u$\textit{\ est solution faible de }$\bigtriangleup _{\mathbf{g}%
}u=\Gamma $\textit{, }$\Gamma \in L^{p}(M)$\textit{, }$1<p<\infty $\textit{,
alors }$u\in H_{2}^{p}(M)$\textit{\ et } 
\[
\left\| u\right\| _{H_{2}^{p}}\leq c\left\| \Gamma \right\|
_{L^{p}}+c\left\| u\right\| _{L^{p}}\,. 
\]

\textit{Or, si }$u$\textit{\ est solution de }$\Delta _{\mathbf{g}%
}u+hu=\Gamma $\textit{\ , on sait que }$u\in H_{2}^{p}$\textit{\ et par
cons\'{e}quent} 
\begin{equation}
\left\| u\right\| _{H_{2}^{p}}\leq c\left\| \Delta _{\mathbf{g}}u+hu\right\|
_{L^{p}}+c\left\| u\right\| _{L^{p}}\,.  
\end{equation}

En fait $+c\left\| u\right\| _{L^{p}}$ est l\`{a} pour tenir compte du fait
que $\Delta _{\mathbf{g}}c=0$ pour toute constante $c$. Mais quand $\Delta _{%
\mathbf{g}}+h$ est coercif, on peut obtenir mieux:

\textit{Si }$u$\textit{\ est solution faible de }$\bigtriangleup _{\mathbf{g}%
}u+hu=\Gamma $\textit{, }$\Gamma \in L^{p}(M)$\textit{, }$1<p<\infty $%
\textit{, et si }$\Delta _{\mathbf{g}}+h$\textit{\ est coercif, alors }$u\in
H_{2}^{p}$\textit{\ et } 
\[
\left\| u\right\| _{H_{2}^{p}}\leq c\left\| \Delta _{\mathbf{g}}u+hu\right\|
_{L^{p}} 
\]

\textit{On montre en fait qu'il existe une constante }$C>0$\textit{\ telle
que pour toute fonction }$u\in H_{2}^{p}$: 
\begin{equation}
\left\| u\right\| _{H_{2}^{p}}\leq C\left\| \Delta _{\mathbf{g}}u+hu\right\|
_{L^{p}}\,.  
\end{equation}

En effet, si une telle constante $C$ n'existe pas, sur $H_{2}^{p}-\{0\}$ 
\[
\frac{\left\| \Delta _{\mathbf{g}}u+hu\right\| _{L^{p}}}{\left\| u\right\|
_{H_{2}^{p}}} 
\]
n'est pas minor\'{e} par une contante $c>0$. Il existe donc une suite $%
(v_{m})\in H_{2}^{p}$ , $v_{m}\neq 0$ telle que 
\[
\frac{\left\| \Delta _{\mathbf{g}}v_{m}+hv_{m}\right\| _{L^{p}}}{\left\|
v_{m}\right\| _{H_{2}^{p}}}\rightarrow 0\,. 
\]
En prenant $\frac{v_{m}}{\left\| v_{m}\right\| _{L^{p}}}$, on peut supposer
que $\left\| v_{m}\right\| _{L^{p}}=1$. Alors $\left\| v_{m}\right\|
_{H_{2}^{p}}\geq 1$, donc 
\[
\left\| \Delta _{\mathbf{g}}v_{m}+hv_{m}\right\| _{L^{p}}\rightarrow 0 
\]
et par (10.5), $\left\| v_{m}\right\| _{H_{2}^{p}}$ est born\'{e}e. Par
reflexivit\'{e}, il existe $v\in H_{2}^{p}$ telle que 
\[
v_{m}\stackrel{H_{2}^{p}}{\rightharpoondown }v\,. 
\]
Or l'inclusion $L^{p}\subset H_{2}^{p}$ est compacte, donc 
\[
v_{m}\stackrel{L^{p}}{\rightarrow }v 
\]
et $\left\| v\right\| _{L^{p}}=1$. De plus comme $v_{m}\stackrel{H_{2}^{p}}{%
\rightharpoondown }v$, pour tout $\phi \in \mathcal{D}(M)$%
\[
\int v_{m}(\Delta _{\mathbf{g}}\phi +h\phi )\rightarrow \int v(\Delta _{%
\mathbf{g}}\phi +h\phi )\,. 
\]
Or $\left\| \Delta _{g}v_{m}+hv_{m}\right\| _{L^{p}}\rightarrow 0$, donc
pour tout $\phi \in L^{p^{\prime }}$ ($p^{\prime }$ conjugu\'{e} de $p$) 
\[
\int \phi (\Delta _{\mathbf{g}}v_{m}+hv_{m})\rightarrow 0 
\]
Donc $\forall \phi \in \mathcal{D}(M):$ $\int \phi (\Delta _{\mathbf{g}%
}v+hv)=0$ et donc $\Delta _{\mathbf{g}}v+hv=0$ faiblement et donc fortement.
Mais comme $\Delta _{\mathbf{g}}+h$ est coercif, 
\[
\int v(\Delta _{\mathbf{g}}v+hv)=\int (\left| \nabla v^{2}\right|
+hv^{2})=0\,\Rightarrow \,v=0 
\]
ce qui contredit $\left\| v\right\| _{L^{p}}=1$.

Montrons donc maintenant le r\'{e}sultat annonc\'{e}:

\textit{Soit }$\Gamma \in L^{p}(M)$\textit{, }$1<p<\infty $\textit{. Alors
il existe une solution }$\beta \in H_{2}^{p}(M)$\textit{\ de l'\'{e}quation }%
$\Delta _{\mathbf{g}}\beta +h\beta =\Gamma $\textit{.}

Il existe ($\Gamma _{k})\in C^{\infty }(M)$ telle que $\Gamma _{k}\stackrel{%
L^{p}}{\rightarrow }\Gamma $. Comme $\Gamma _{k}\in L^{2}$, les m\'{e}thodes
variationnelles classiques donnent l'existence de solutions $\beta _{k}\in
C^{\infty }$ de 
\[
\Delta _{\mathbf{g}}\beta _{k}+h\beta _{k}=\Gamma _{k} 
\]
Comme $\left\| \Gamma _{k}\right\| _{L^{p}}$ est born\'{e}e, par
l'estim\'{e}e (10.6) la suite ($\beta _{k})$ est born\'{e}e dans $H_{2}^{p}$.
Il existe donc \textit{\ }$\beta \in H_{2}^{p}$ telle que 
\[
\beta _{k}\stackrel{H_{2}^{p}}{\rightharpoondown }\beta \,. 
\]
Alors $\forall \phi \in \mathcal{D}(M):$ 
\[
\int \beta _{k}(\Delta _{\mathbf{g}}\phi +h\phi )\rightarrow \int \beta
(\Delta _{\mathbf{g}}\phi +h\phi )\,. 
\]
Mais comme $\beta _{k}\in C^{\infty }$%
\[
\int \beta _{k}(\Delta _{\mathbf{g}}\phi +h\phi )=\int \phi (\Delta _{%
\mathbf{g}}\beta _{k}+h\beta _{k})=\int \phi \Gamma _{k}\rightarrow \int
\phi \Gamma 
\]
et donc $\forall \phi \in \mathcal{D}(M):$%
\[
\int \beta (\Delta _{\mathbf{g}}\phi +h\phi )=\int \phi \Gamma \,. 
\]
Ainsi $\beta $ est solution faible de $\Delta _{g}\beta +h\beta =\Gamma $%
\textit{.}

\subsection{Deuxi\`{e}me \'{e}tape:}

Soit $y\in M$ fix\'{e}. On pose $r:=r_{y}:=d_{\mathbf{g}}(x,y)$. Soit
\'{e}galement une fonction cut-off $\eta _{y}=\eta $ valant 1 dans $%
B(y,\delta )$ et valant 0 dans $M\backslash B(0,2\delta )$ pour $\delta >0$
assez petit. On pose alors 
\[
\Gamma _{y}:=\Gamma =-\Delta _{\mathbf{g}}(\frac{\eta }{r^{n-2}})-h\frac{%
\eta }{r^{n-2}}\,. 
\]
Notons que 
\[
\left| \Delta _{\mathbf{g}}\frac{1}{r^{n-2}}\right| =\left| \frac{n-2}{%
r^{n-1}}\partial _{r}(\ln \sqrt{\det \mathbf{g}})\right| \leq \frac{c}{%
r^{n-2}} 
\]
o\`{u} cette expression est \`{a} ``lire'' dans une carte exponentielle. Donc $\Gamma
_{y}\in L^{p}(M)$ pour tout $p<\frac{n}{n-2}$.

Soit enfin $\beta _{y}:=\beta $ une solution faible de $\Delta _{\mathbf{g}%
}\beta _{y}+h\beta _{y}=\Gamma _{y}$

On construit alors la \textit{fonction de Green de l'op\'{e}rateur }$\Delta
_{\mathbf{g}}+h$ en posant 
\[
G_{y}:=\frac{1}{(n-2)\omega _{n-1}}(\beta _{y}+\frac{\eta _{y}}{r_{y}^{n-2}}%
) 
\]
$G_{y}$ est $C^{\infty }$ sur $M\backslash \{y\}$ et est solution sur $M,$
au sens des distribution, de 
\[
\bigtriangleup _{\mathbf{g}}G_{y}+hG_{y}=\delta _{y} 
\]
La d\'{e}monstration de cette derni\`{e}re identit\'{e} est classique et
tout \`{a} fait analogue au cas du laplacien dans $\Bbb{R}^{n}$; nous serons
donc rapide.

Il faut utiliser la formule de Green 
\[
-\int_{B}v\Delta _{\mathbf{g}}u+\int_{B}(\nabla u,\nabla v)_{\mathbf{g}%
}=\int_{\partial B}v\frac{\partial u}{\partial \nu } 
\]
o\`{u} $\nu $ est la normale ext\'{e}rieure au domaine $B$.

Soit $u\in C^{\infty }(M)$ et soit $\varepsilon >0$ petit. On note $\nu ^{-}$
la normale int\'{e}rieure \`{a} la boule $B_{\varepsilon }:=B(y,\varepsilon
) $. Alors 
\[
\int_{M}G_{y}(\Delta _{\mathbf{g}}u+hu)=\int_{M\backslash B_{\varepsilon
}}G_{y}(\Delta _{\mathbf{g}}u+hu)+\int_{B_{\varepsilon }}G_{y}(\Delta _{%
\mathbf{g}}u+hu)\,. 
\]
Sur $M\backslash B_{\varepsilon },$ $G_{y}$ est $C^{\infty }$ donc 
\[
\int_{M\backslash B_{\varepsilon }}G_{y}(\Delta _{\mathbf{g}%
}u+hu)=\int_{M\backslash B_{\varepsilon }}u(\Delta _{\mathbf{g}%
}G_{y}+hG_{y})+\int_{\partial (M\backslash B_{\varepsilon })}(u\frac{%
\partial G_{y}}{\partial \nu ^{-}}-G_{y}\frac{\partial u}{\partial \nu ^{-}}%
)\,. 
\]
Or par d\'{e}finition de $G_{y}$, sur $M\backslash B_{\varepsilon }$ : $%
\Delta _{\mathbf{g}}G_{y}+hG_{y}=0$. D'o\`{u} 
\[
\int_{M}G_{y}(\Delta _{\mathbf{g}}u+hu)=\int_{B_{\varepsilon }}G_{y}(\Delta
_{\mathbf{g}}u+hu)+\int_{\partial B_{\varepsilon }}(G_{y}\frac{\partial u}{%
\partial \nu }-u\frac{\partial G_{y}}{\partial \nu })\,. 
\]
- Comme $G_{y}\in L^{1}(M)$ et que $u\in C^{\infty }(M)$ 
\[
\int_{B_{\varepsilon }}G_{y}(\Delta _{\mathbf{g}}u+hu)\stackunder{%
\varepsilon \rightarrow 0}{\rightarrow }0 
\]
- Ensuite 
\begin{eqnarray*}
\int_{\partial B_{\varepsilon }}(G_{y}\frac{\partial u}{\partial \nu }-u%
\frac{\partial G_{y}}{\partial \nu })= && \frac{1}{(n-2)\omega _{n-1}}%
\int_{\partial B_{\varepsilon }}(\beta \frac{\partial u}{\partial \nu }-u%
\frac{\partial \beta }{\partial \nu }) \\ 
&& +\frac{1}{(n-2)\omega _{n-1}}\int_{\partial B_{\varepsilon }}\frac{1}{%
r^{n-2}}\frac{\partial u}{\partial \nu }-\frac{1}{(n-2)\omega _{n-1}}%
\int_{\partial B_{\varepsilon }}u\frac{\partial r^{-(n-2)}}{\partial \nu }
\end{eqnarray*}
Comme $\beta \in H_{2}^{p}$ pour $p>1$, $\beta \frac{\partial u}{\partial
\nu }-u\frac{\partial \beta }{\partial \nu }\in L^{1}(M)$ et donc 
\[
\frac{1}{(n-2)\omega _{n-1}}\int_{\partial B_{\varepsilon }}(\beta \frac{%
\partial u}{\partial \nu }-u\frac{\partial \beta }{\partial \nu })%
\stackunder{\varepsilon \rightarrow 0}{\rightarrow }0\,. 
\]
On a \'{e}galement 
\[
\left| \int_{\partial B_{\varepsilon }}\frac{1}{r^{n-2}}\frac{\partial u}{%
\partial \nu }\right| \leq \stackunder{M}{Sup}\left| \nabla u\right|
\int_{\partial B_{\varepsilon }}\frac{1}{r^{n-2}}\leq c\stackunder{M}{Sup}%
\left| \nabla u\right| .\varepsilon \stackunder{\varepsilon \rightarrow 0}{%
\rightarrow }0 
\]
- Enfin 
\[
\begin{array}{ll}
-\frac{1}{(n-2)\omega _{n-1}}\int_{\partial B_{\varepsilon }}u\frac{\partial
r^{-(n-2)}}{\partial \nu } & =\frac{1}{\omega _{n-1}}\int_{\partial
B_{\varepsilon }}u\frac{1}{r^{n-1}} \\ 
& =\frac{1}{\omega _{n-1}}.\frac{1}{\varepsilon ^{n-1}}\int_{\partial
B_{\varepsilon }}u \\ 
& \stackunder{\varepsilon \rightarrow 0}{\rightarrow }u(y)
\end{array}
\]
d'o\`{u} finalement 
\[
\int_{M}G_{y}(\Delta _{\mathbf{g}}u+hu)=u(y) 
\]

\subsection{Troisi\`{e}me \'{e}tape:}

C'est la partie qui int\'{e}resse le plus directement notre travail
puisqu'elle concerne les estim\'{e}es sur la fonction de Green. Nous allons
montrer par r\'{e}currence que pour tout $\varepsilon >0$%
\[
\left| \beta _{y}\right| \leq \frac{c}{r^{n-3+\varepsilon }}%
\,\,\,et\,\,\,\left| \nabla \beta _{y}\right| \leq \frac{c}{%
r^{n-2+\varepsilon }}\,. 
\]

Pour cela on montre que $r^{n-3+\varepsilon }\beta _{y}$ et $\left| \nabla
\beta _{y}\right| r^{n-2+\varepsilon }$ sont dans $C^{0}(M)$. Le principe
est d'\'{e}tudier l'\'{e}quation v\'{e}rifi\'{e}e par $\Delta (r^{p}\beta )$
et d'en d\'{e}duire avec les th\'{e}or\`{e}mes de r\'{e}gularit\'{e}
classiques \`{a} quels $H_{k}^{p}$ appartiennent $r^{p}\beta $ et $\nabla
(r^{p}\beta )$. Ces estim\'{e}es entrainent la propri\'{e}t\'{e} (10.2).

Nous traiterons le cas o\`{u} la dimension $n\geq 4$, le cas de la dimension
3 \'{e}tant analogue mais plus facile. Pour all\'{e}ger l'\'{e}criture nous
noterons $u\in $($H_{k}^{p})^{-}$ pour dire que $u\in $($H_{k}^{q})$ pour
tout $q<p$ aussi proche que l'on veut. Ainsi 
\[
\frac{1}{r^{n-2}}\in (L^{\frac{n}{n-2}})^{-} 
\]

et donc 
\[
\beta _{y}\in (H_{2}^{\frac{n}{n-2}})^{-}\,\,. 
\]

Rappelons que:

-si $\frac{1}{p}\geq \frac{1}{q}-\frac{k-m}{n}>0$ alors $H_{k}^{q}(M)\subset
H_{m}^{p}(M)$

-si $\frac{1}{q}<\frac{k-m}{n}\,$alors $H_{k}^{q}(M)\subset C^{m}(M)$

-si pour $0<\alpha <1:$ $\frac{1}{q}<\frac{1-\alpha }{n}\,$alors $%
H_{1}^{q}(M)\subset C^{0,\alpha }(M)$

Ainsi 
\[
\beta _{y}\in (H_{2}^{\frac{n}{n-2}})^{-}\subset (H_{1}^{\frac{n}{n-3}%
})^{-}\subset (L^{\frac{n}{n-4}})^{-} 
\]
et donc $\left| \nabla \beta _{y}\right| \in (L^{\frac{n}{n-3}})^{-}$.

Soit $p$ un r\'{e}el $1<p\leq 2$. On a en utilisant l'\'{e}quation 
\[
\Delta _{\mathbf{g}}\beta _{y}+h\beta _{y}=-\Delta _{\mathbf{g}}(\frac{\eta 
}{r^{n-2}})-h\frac{\eta }{r^{n-2}} 
\]
v\'{e}rifi\'{e}e par $\beta _{y}:$%
\[
\begin{array}{cc}
& \Delta _{\mathbf{g}}(r^{p}\beta )=r^{p}\Delta _{\mathbf{g}}\beta +\beta
\Delta _{\mathbf{g}}r^{p}-2(\nabla r^{p},\nabla \beta )_{\mathbf{g}} \\ 
& =-hr^{p}\beta -r^{p}(\Delta _{\mathbf{g}}(\frac{\eta }{r^{n-2}})+h\frac{%
\eta }{r^{n-2}})+\beta \Delta _{\mathbf{g}}r^{p}-2pr^{p-1}(\nabla r,\nabla
\beta )_{\mathbf{g}}\,.
\end{array}
\]
Or:

i/: $hr^{p}\beta \in (L^{\frac{n}{n-4}})^{-}\,$car $\beta _{y}\in (L^{\frac{n%
}{n-4}})^{-}$

ii/: $\left| r^{p}(\Delta _{\mathbf{g}}(\frac{\eta }{r^{n-2}})+h\frac{\eta }{%
r^{n-2}})\right| \leq \frac{c}{r^{n-2-p}}\in (L^{\frac{n}{n-(p+2)}})^{-}$

iii/: $\left| \beta \Delta _{\mathbf{g}}r^{p}\right| \leq cr^{p-2}\left|
\beta \right| \in (L^{\frac{n}{n-4}})^{-}$

iv/: $\left| r^{p-1}(\nabla r,\nabla \beta )_{\mathbf{g}}\right| \leq
r^{p-1}\left| \nabla \beta _{y}\right| \in (L^{\frac{n}{n-3}})^{-}$.

Si $1<p\leq 2$%
\[
L^{\frac{n}{n-4}}\subset L^{\frac{n}{n-(p+2)}}\subset L^{\frac{n}{n-3}%
}\subset L^{\frac{n}{n-(p+1)}} 
\]
Donc 
\[
Si\,\,\,1<p\leq 2:\,\,r^{p}\beta \in (H_{2}^{\frac{n}{n-(p+1)}})^{-} 
\]
et on peut commencer la r\'{e}currence \`{a} un tel $p$.

Prenons donc comme hypoth\`{e}se de r\'{e}currence que si $p<n-4:$%
\[
r^{p}\beta \in (H_{2}^{\frac{n}{n-(p+1)}})^{-}\,. 
\]
Alors tant que $p<n-3:$%
\[
r^{p}\beta \in (H_{2}^{\frac{n}{n-(p+1)}})^{-}\subset (H_{1}^{\frac{n}{%
n-(p+2)}})^{-}\subset (L^{\frac{n}{n-(p+3)}})^{-} 
\]
et 
\[
\Delta _{\mathbf{g}}(r^{p+1}\beta )=-hr^{p+1}\beta -r^{p+1}(\Delta _{\mathbf{%
g}}(\frac{\eta }{r^{n-2}})+h\frac{\eta }{r^{n-2}})+\beta \Delta _{\mathbf{g}%
}r^{p+1}-2(p+1)r^{p}(\nabla r,\nabla \beta )_{\mathbf{g}}\,. 
\]
On a comme ci-dessus:

i/: $hr^{p+1}\beta =hr(r^{p}\beta )\in (L^{\frac{n}{n-(p+3)}})^{-}\,$par
hypoth\`{e}se de r\'{e}currence.

ii/: $\left| r^{p+1}(\Delta _{\mathbf{g}}(\frac{\eta }{r^{n-2}})+h\frac{\eta 
}{r^{n-2}})\right| \leq \frac{c}{r^{n-3-p}}\in (L^{\frac{n}{n-(p+3)}})^{-}$

iii/: $\left| \beta \Delta _{\mathbf{g}}r^{p+1}\right| \leq cr^{p-1}\left|
\beta \right| \in (L^{\frac{n}{n-(p+2)}})^{-}$ par hypothèse de
r\'{e}currence.

iv/: $\left| r^{p}(\nabla r,\nabla \beta )_{\mathbf{g}}\right| \leq \left|
r^{p}\nabla \beta \right| $. Or 
\[
r^{p}\nabla \beta =\nabla (r^{p}\beta )-\beta \nabla r^{p}=\nabla
(r^{p}\beta )-p(r^{p-1}\beta )\nabla r 
\]
or par hypoth\`{e}se de r\'{e}currence, $\nabla (r^{p}\beta )\in (L^{\frac{n%
}{n-(p+2)}})^{-}\,$et $\left| (r^{p-1}\beta )\nabla r\right| =\left|
r^{p-1}\beta \right| \in (L^{\frac{n}{n-(p+2)}})^{-}$.

Comme 
\[
L^{\frac{n}{n-(p+3)}}\subset L^{\frac{n}{n-(p+2)}} 
\]
on a finalement que $\Delta _{\mathbf{g}}(r^{p+1}\beta )\in (L^{\frac{n}{%
n-(p+2)}})^{-}$ et donc $r^{p+1}\beta \in (H_{2}^{\frac{n}{n-(p+2)}})^{-}\,$%
d'o\`{u} le fonctionnement de la r\'{e}currence.

Ainsi pour tout r\'{e}el $1<p<n-3:\,\,r^{p}\beta \in (H_{2}^{\frac{n}{n-(p+1)%
}})^{-}$

Maintenant pour $0<\varepsilon <1$%
\[
\Delta _{\mathbf{g}}(r^{n-3+\varepsilon }\beta )=-hr^{n-3+\varepsilon }\beta
-r^{n-3+\varepsilon }(\Delta _{\mathbf{g}}(\frac{\eta }{r^{n-2}})+h\frac{%
\eta }{r^{n-2}})+\beta \Delta _{\mathbf{g}}r^{n-3+\varepsilon
}-cr^{n-4+\varepsilon }(\nabla r,\nabla \beta )_{\mathbf{g}} 
\]

i/: $hr^{n-3+\varepsilon }\beta =hr^{2}(r^{n-5+\varepsilon }\beta )\in (L^{%
\frac{n}{2-\varepsilon }})^{-}\,$

ii/: $\left| r^{n-3+\varepsilon }(\Delta _{\mathbf{g}}(\frac{\eta }{r^{n-2}}%
)+h\frac{\eta }{r^{n-2}})\right| \leq \frac{c}{r^{1-\varepsilon }}\in (L^{%
\frac{n}{1-\varepsilon }})^{-}\subset (L^{\frac{n}{2-\varepsilon }})^{-}$

iii/: $\left| \beta \Delta _{\mathbf{g}}r^{n-3+\varepsilon }\right| \leq
cr^{n-5+\varepsilon }\left| \beta \right| \in (L^{\frac{n}{2-\varepsilon }%
})^{-}$

iv/: $\left| r^{n-4+\varepsilon }(\nabla r,\nabla \beta )_{\mathbf{g}%
}\right| \leq \left| r^{n-4+\varepsilon }\nabla \beta \right| \leq \left|
\nabla (r^{n-4+\varepsilon }\beta )\right| +c\left| r^{n-5+\varepsilon
}\beta \right| \in (L^{\frac{n}{2-\varepsilon }})^{-}$

Donc $\Delta _{\mathbf{g}}(r^{n-3+\varepsilon }\beta )\in (L^{\frac{n}{%
2-\varepsilon }})^{-}$ et donc $r^{n-3+\varepsilon }\beta \in (H_{2}^{\frac{n%
}{2-\varepsilon }})^{-}$ . Autrement dit, pour tout $0<
\varepsilon ^{\prime }<<\varepsilon $%
\[
r^{n-3+\varepsilon }\beta \in H_{2}^{\frac{n}{2-\varepsilon +\varepsilon
^{\prime }}} 
\]
donc 
\[
r^{n-3+\varepsilon }\beta \in H_{2}^{\frac{n}{2-\varepsilon +\varepsilon
^{\prime }}}\subset H_{1}^{\frac{n}{1-(\varepsilon -\varepsilon ^{\prime })}%
}\subset C^{0,(\varepsilon -\varepsilon ^{\prime })} 
\]
et donc $r^{n-3+\varepsilon }\beta $ est born\'{e}e sur $M$.

Enfin 
\[
\Delta _{\mathbf{g}}(r^{n-2+\varepsilon }\beta )=-hr^{n-2+\varepsilon }\beta
-r^{n-2+\varepsilon }(\Delta _{\mathbf{g}}(\frac{\eta }{r^{n-2}})+h\frac{%
\eta }{r^{n-2}})+\beta \Delta _{\mathbf{g}}r^{n-2+\varepsilon
}-cr^{n-3+\varepsilon }(\nabla r,\nabla \beta )_{\mathbf{g}} 
\]

i/: $hr^{n-2+\varepsilon }\beta =hr(r^{n-3+\varepsilon }\beta )\in C^{0}$
d'après ce qui pr\'{e}cède$\,$

ii/: $\left| r^{n-2+\varepsilon }(\Delta _{\mathbf{g}}(\frac{\eta }{r^{n-2}}%
)+h\frac{\eta }{r^{n-2}})\right| \leq cr^{\varepsilon }\in C^{0}$

iii/: $\left| \beta \Delta _{\mathbf{g}}r^{n-2+\varepsilon }\right| \leq
cr^{n-4+\varepsilon }\left| \beta \right| \in (L^{\frac{n}{1-\varepsilon }%
})^{-}$

iv/: $\left| r^{n-3+\varepsilon }(\nabla r,\nabla \beta )_{\mathbf{g}%
}\right| \leq \left| r^{n-4+\varepsilon }\nabla \beta \right| \leq \left|
\nabla (r^{n-3+\varepsilon }\beta )\right| +c\left| r^{n-4+\varepsilon
}\beta \right| \in (L^{\frac{n}{1-\varepsilon }})^{-}$

donc $r^{n-2+\varepsilon }\beta \in (H_{2}^{\frac{n}{1-\varepsilon }})^{-}$.
Or $(H_{2}^{\frac{n}{1-\varepsilon }})^{-}\subset C^{1}$, donc $\nabla
(r^{n-2+\varepsilon }\beta )\in C^{0}$. Mais 
\[
\nabla (r^{n-2+\varepsilon }\beta )=cr^{n-3+\varepsilon }\beta
+r^{n-2+\varepsilon }\nabla \beta 
\]
donc $\left| r^{n-2+\varepsilon }\nabla \beta \right| \in C^{0}$.

\subsection{Quatri\`{e}me \'{e}tape:}

Il nous reste \`{a} montrer que $G_{y}>0$ et l'estim\'{e}e (10.3). Il s'agit
en fait essentiellement d'appliquer le principe du maximum. (Rappelons que
notre Laplacien $\Delta _{\mathbf{g}}$ est le laplacien des
g\'{e}om\`{e}tres, c'est-\`{a}-dire avec la convention ``du signe moins'': $%
\Delta _{\mathbf{g}}=-\nabla ^{i}\nabla _{i}$ ; le principe du maximum
devient donc un principe du minimum !)

Sur $M\backslash \{y\}$ : $\Delta _{\mathbf{g}}G_{y}+hG_{y}=0$, donc $G_{y}$
est $C^{\infty }$ sur $M\backslash \{y\}$.

Notons $U_{t}=M\backslash B(y,t)$ pour $t$ petit. Le principe du maximum
appliqu\'{e} \`{a} $G_{y}$ sur $U_{t}$ nous dit que si $G_{y}$ atteint un
minimum n\'{e}gatif ou nul \`{a} l'int\'{e}rieur de $U_{t}$ alors $%
G_{y}=cste.$ Sinon $G_{y}$ atteint son minimum sur le bord $\partial U_{t}$.
Mais dans $B(y,t)$ 
\[
\left| \beta _{y}\right| \leq \frac{c}{r^{n-3+\varepsilon }}%
=r^{1-\varepsilon }\frac{c}{r^{n-2}} 
\]
donc dans $B(y,t)$ si $t<\delta $%
\[
G_{y}(x)\sim \frac{c}{r^{n-2}}\stackunder{x\rightarrow y}{\rightarrow }%
+\infty \,. 
\]
En prenant $t$ assez petit, ceci prouve que d'une part $G_{y}$ n'est pas
constante, et d'autre part que $G_{y}>0$.

Enfin, sur $M$%
\[
\left| \nabla \beta _{y}\right| \leq \frac{c}{r^{n-2+\varepsilon }}%
=r^{1-\varepsilon }\frac{c}{r^{n-1}} 
\]
et 
\[
\left| \nabla \frac{\eta }{r^{n-2}}\right| \leq \frac{c}{r^{n-1}}\,. 
\]
Donc sur $M$%
\[
\left| \nabla G_{y}\right| \leq \frac{c}{r^{n-1}} 
\]
Sur $B(y,\delta )$%
\[
\left| \nabla \frac{\eta }{r^{n-2}}\right| =\frac{n-2}{r^{n-1}} 
\]
et donc sur $B(y,\delta )$%
\[
\frac{c^{-1}}{r^{n-1}}\leq \left| \nabla G_{y}\right| \leq \frac{c}{r^{n-1}} 
\]
et cette estim\'{e}e est donc vraie sur tout $M$.

Comme par ailleurs $G_{y}>0$ sur $M$, on a aussi sur tout $M:$%
\[
\frac{c^{-1}}{r^{n-2}}\leq \left| G_{y}\right| \leq \frac{c}{r^{n-2}}\,. 
\]
On en d\'{e}duit (10.3).

\chapter{Appendice C}
\pagestyle{myheadings}\markboth{\textbf{Appendice C.}}{\textbf{Appendice C.}}
Dans la partie centrale du chapitre 3, nous avons affirm\'{e} que 
\[
\stackunder{t\rightarrow 1}{\overline{\lim }}\frac{C_{t}}{\int_{B(0,\delta )}%
\overline{u}_{t}^{2}dx}\leq \varepsilon _{\delta } 
\]
o\`{u} $\varepsilon _{\delta }\rightarrow 0$ quand $\delta \rightarrow 0\,$%
et o\`{u}

\[
C_{t}=\left| \int_{B(0,\delta )}\eta ^{2}(\mathbf{g}_{t}^{ij}-\delta
^{ij})\partial _{i}\overline{u}_{t}\partial _{j}\overline{u}_{t}dx\right| 
\]

Ceci est montr\'{e} dans l'article de Z. Djadli et O. Druet [9], sur lequel
nous nous appuyons, dans le cas $f=cste$, et il n'y a pas de changements
pour une fonction $f$ non constante. C'est pourquoi pour rendre plus lisible
la d\'{e}monstration du chapitre 3 nous avons pr\'{e}f\'{e}r\'{e} reporter
le calcul de cette limite. N\'{e}anmoins, la d\'{e}monstration de Z. Djadli
et O. Druet se faisant dans un but et dans un cadre diff\'{e}rents, et avec $%
f=cste,$ pour que notre travail soit ``complet'', nous reprenons rapidement
pour le lecteur int\'{e}ress\'{e} la d\'{e}monstration dans cet appendice.

Le d\'{e}veloppement de Cartan de la m\'{e}trique $\mathbf{g}$ autour de $%
x_{t}$, nous donne 
\begin{eqnarray*}
C_{t} &\leq &C\int_{B(0,\delta )}\eta ^{2}\left| Rm_{\mathbf{g}%
}(x_{t})(\nabla \overline{u}_{t},x,\nabla \overline{u}_{t},x)\right| dv_{%
\mathbf{g}_{t}} \\
&&+\varepsilon _{\delta }\int_{B(0,\delta )}\eta ^{2}\left| \nabla \overline{%
u}_{t}\right| _{\mathbf{g}_{t}}^{2}dv_{\mathbf{g}_{t}}\,.
\end{eqnarray*}
En int\'{e}grant par parties et en utilisant l'\'{e}quation v\'{e}rifi\'{e}e
par les fonctions $\overline{u}_{t}$ et l'estim\'{e}e faible $\left| x^{2}%
\overline{u}_{t}^{2^{*}}\right| \leq c\overline{u}_{t}^{2}$\thinspace on
obtient 
\[
\int_{B(0,\delta )}\eta ^{2}\left| \nabla \overline{u}_{t}\right| _{\mathbf{g%
}_{t}}^{2}dv_{\mathbf{g}_{t}}\leq C\int_{B(0,\delta )}\overline{u}%
_{t}^{2}dv_{\mathbf{g}_{t}}\,. 
\]
Ensuite le changement d'\'{e}chelle faisant passer de $\overline{u}_{t}$
\`{a} $\widetilde{u}_{t}$ nous donne 
\[
\nabla \widetilde{u}_{t}=\nabla (\mu _{t}^{\frac{n-2}{2}}\overline{u}%
_{t}(\mu _{t}x))=\mu _{t}^{\frac{n}{2}}(\nabla \overline{u}_{t})(\mu _{t}x) 
\]
et par cons\'{e}quent 
\[
\int_{B(0,\delta )}\eta ^{2}\left| Rm_{\mathbf{g}}(x_{t})(\nabla \overline{u}%
_{t},x,\nabla \overline{u}_{t},x)\right| dv_{\mathbf{g}_{t}}=\mu
_{t}^{2}\int_{B(0,\delta \mu _{t}^{-1})}\eta ^{2}(\mu _{t}x)\left| Rm_{%
\mathbf{g}}(x_{t})(\nabla \widetilde{u}_{t},x,\nabla \widetilde{u}%
_{t},x)\right| dv_{\widetilde{\mathbf{g}}_{t}}\,. 
\]
On coupe l'int\'{e}grale de droite en deux 
\[
\int_{B(0,\delta \mu _{t}^{-1})}=\int_{B(0,R)}+\int_{B(0,\delta \mu
_{t}^{-1})\backslash B(0,R)}\,. 
\]
Or sur $B(0,R)$ : $\widetilde{u}_{t}\stackrel{C^{0}}{\stackunder{%
t\rightarrow t_{0}}{\rightarrow }}\widetilde{u}$ qui est radiale, donc
puisque $x$ et $\nabla \widetilde{u}$ sont colin\'{e}aires 
\[
\int_{B(0,R)}\eta ^{2}(\mu _{t}x)\left| Rm_{\mathbf{g}}(x_{t})(\nabla 
\widetilde{u}_{t},x,\nabla \widetilde{u}_{t},x)\right| dv_{\widetilde{%
\mathbf{g}}_{t}}\stackunder{t\rightarrow t_{0}}{\rightarrow }0\,. 
\]
En rappelant que 
\[
\int_{B(0,\delta )}\overline{u}_{t}^{2}dv_{\mathbf{g}_{t}}=\mu
_{t}^{2}\int_{B(0,\delta \mu _{t}^{-1})}\widetilde{u}_{t}^{2}dv_{\widetilde{%
\mathbf{g}}_{t}} 
\]
on voit que 
\[
\stackunder{t\rightarrow 1}{\overline{\lim }}\frac{C_{t}}{\int_{B(0,\delta )}%
\overline{u}_{t}^{2}dx}\leq C\stackunder{t\rightarrow 1}{\overline{\lim }}%
\frac{\int_{B(0,\delta \mu _{t}^{-1})\backslash B(0,R)}\eta ^{2}(\mu
_{t}x)\left| Rm_{\mathbf{g}}(x_{t})(\nabla \widetilde{u}_{t},x,\nabla 
\widetilde{u}_{t},x)\right| dv_{\widetilde{\mathbf{g}}_{t}}}{%
\int_{B(0,\delta \mu _{t}^{-1})}\widetilde{u}_{t}^{2}dv_{\widetilde{\mathbf{g%
}}_{t}}}+\varepsilon _{\delta }\,. 
\]
Maintenant, pour des vecteurs $x,y$, par d\'{e}finition, 
\[
Rm_{\mathbf{g}}(P)(x,y,x,y)=K(P)(\left| x\right| _{\mathbf{g}}\left|
y\right| _{\mathbf{g}}-(x,y)_{\mathbf{g}}^{2}) 
\]
o\`{u} $K(P)$ est la courbure sectionnelle au point $P$.

Donc en majorant $K(P)\,$dans un voisinage de $x_{0}$%
\[
\left| Rm_{\mathbf{g}}(x_{t})(\nabla \widetilde{u}_{t},x,\nabla \widetilde{u}%
_{t},x)\right| \leq C[\left| \nabla (\left| x\right| \widetilde{u}%
_{t})\right| _{\widetilde{\mathbf{g}}_{t}}^{2}-(\nabla (\left| x\right| 
\widetilde{u}_{t}),\nu )_{\widetilde{\mathbf{g}}_{t}}^{2}] 
\]
o\`{u} $\nu =\frac{x}{\left| x\right| }$. En int\'{e}grant par parties, on
obtient 
\begin{eqnarray*}
&&\int_{B(0,\delta \mu _{t}^{-1})\backslash B(0,R)}\eta ^{2}(\mu
_{t}x)[\left| \nabla (\left| x\right| \widetilde{u}_{t})\right| _{\widetilde{%
\mathbf{g}}_{t}}^{2}-(\nabla (\left| x\right| \widetilde{u}_{t}),\nu )_{%
\widetilde{\mathbf{g}}_{t}}^{2}]dv_{\widetilde{\mathbf{g}}_{t}} \\
&\leq &C\int_{\partial B(0,R)}\left| \nabla (\left| x\right| ^{2}\widetilde{u%
}_{t}^{2})\right| _{e}d\sigma _{e} \\
&&+\int_{B(0,\delta \mu _{t}^{-1})\backslash B(0,R)}\eta ^{2}(\mu
_{t}x)\Delta _{\widetilde{\mathbf{g}}_{t}}(\left| x\right| \widetilde{u}%
_{t})\left| x\right| \widetilde{u}_{t}dv_{\widetilde{\mathbf{g}}_{t}} \\
&&-\int_{B(0,\delta \mu _{t}^{-1})\backslash B(0,R)}\eta ^{2}(\mu
_{t}x)(\nabla (\left| x\right| \widetilde{u}_{t}),\nu )_{\widetilde{\mathbf{g%
}}_{t}}^{2}dv_{\widetilde{\mathbf{g}}_{t}} \\
&&+C\int_{B(0,\delta \mu _{t}^{-1})\backslash B(0,\frac{\delta }{2}\mu
_{t}^{-1})}\widetilde{u}_{t}^{2}dv_{\widetilde{\mathbf{g}}_{t}}
\end{eqnarray*}
o\`{u} l'on a utilis\'{e} entre autres que 
\begin{eqnarray*}
\left| \nabla \eta \right| &\leq &c/\delta \\
\left| \Delta \eta \right| &\leq &c/\delta ^{2} \\
\,\,\Delta (\eta ^{2}(\mu _{t}x)) &=&\mu _{t}^{2}(\Delta \eta ^{2})(\mu
_{t}x)
\end{eqnarray*}
$\,$et que $\mu _{t}^{2}\left| x\right| ^{2}\leq \delta ^{2}$ sur $%
B(0,\delta \mu _{t}^{-1})$.

De plus $\eta ^{2}(\mu _{t}x)\rightarrow 1$ uniform\'{e}ment sur $\partial
B(0,R)$, donc en d\'{e}veloppant l'in\'{e}galit\'{e} ci-dessus, pour $t$
assez proche de $t_{0}$: 
\begin{eqnarray*}
&&\int_{B(0,\delta \mu _{t}^{-1})\backslash B(0,R)}\eta ^{2}(\mu
_{t}x)[\left| \nabla (\left| x\right| \widetilde{u}_{t})\right| _{\widetilde{%
\mathbf{g}}_{t}}^{2}-(\nabla (\left| x\right| \widetilde{u}_{t}),\nu )_{%
\widetilde{\mathbf{g}}_{t}}^{2}]dv_{\widetilde{\mathbf{g}}_{t}} \\
&\leq &C\int_{\partial B(0,R)}\left| \nabla (\left| x\right| ^{2}\widetilde{u%
}_{t}^{2})\right| _{e}d\sigma _{e} \\
&&+\int_{B(0,\delta \mu _{t}^{-1})\backslash B(0,R)}\eta ^{2}(\mu
_{t}x)\left| x\right| ^{2}\widetilde{u}_{t}\Delta _{\widetilde{\mathbf{g}}%
_{t}}\widetilde{u}_{t}\left| x\right| \widetilde{u}_{t}dv_{\widetilde{%
\mathbf{g}}_{t}} \\
&&+\int_{B(0,\delta \mu _{t}^{-1})\backslash B(0,R)}\eta ^{2}(\mu
_{t}x)\left| x\right| \Delta _{\widetilde{\mathbf{g}}_{t}}(\left| x\right| )%
\widetilde{u}_{t}^{2}dv_{\widetilde{\mathbf{g}}_{t}} \\
&&+C\int_{B(0,\delta \mu _{t}^{-1})\backslash B(0,\frac{\delta }{2}\mu
_{t}^{-1})}\widetilde{u}_{t}^{2}dv_{\widetilde{\mathbf{g}}_{t}} \\
&&-2\int_{B(0,\delta \mu _{t}^{-1})\backslash B(0,R)}\eta ^{2}(\mu _{t}x)%
\widetilde{u}_{t}(\nabla \widetilde{u}_{t},x)_{\widetilde{\mathbf{g}}%
_{t}}^{2}dv_{\widetilde{\mathbf{g}}_{t}} \\
&&-\int_{B(0,\delta \mu _{t}^{-1})\backslash B(0,R)}\eta ^{2}(\mu
_{t}x)\left[ (\nabla \widetilde{u}_{t},x)_{\widetilde{\mathbf{g}}_{t}}+%
\widetilde{u}_{t}\right] ^{2}dv_{\widetilde{\mathbf{g}}_{t}}\,\,.
\end{eqnarray*}
On utilise maintenant la relation suivante: 
\[
\left| x\right| \Delta _{\widetilde{\mathbf{g}}_{t}}(\left| x\right| )\leq
\left| x\right| \Delta _{e}(\left| x\right| )+c\mu _{t}^{2}\left| x\right|
^{2}\leq -(n-1)+c\mu _{t}^{2}\left| x\right| ^{2} 
\]
et l'\'{e}quation v\'{e}rifi\'{e}e par $\widetilde{u}_{t}$

\[
\widetilde{u}_{t}\Delta _{\widetilde{\mathbf{g}}_{t}}\widetilde{u}%
_{t}=\lambda _{t}\widetilde{f}_{t}\widetilde{u}_{t}^{2^{*}}-\mu _{t}^{2}%
\widetilde{h}_{t}\widetilde{u}_{t}^{2} 
\]
pour obtenir: 
\begin{eqnarray*}
&&\int_{B(0,\delta \mu _{t}^{-1})\backslash B(0,R)}\eta ^{2}(\mu
_{t}x)[\left| \nabla (\left| x\right| \widetilde{u}_{t})\right| _{\widetilde{%
\mathbf{g}}_{t}}^{2}-(\nabla (\left| x\right| \widetilde{u}_{t}),\nu )_{%
\widetilde{\mathbf{g}}_{t}}^{2}]dv_{\widetilde{\mathbf{g}}_{t}} \\
&\leq &C\int_{\partial B(0,R)}\left| \nabla (\left| x\right| ^{2}\widetilde{u%
}_{t}^{2})\right| _{e}d\sigma _{e}+C\int_{B(0,\delta \mu
_{t}^{-1})\backslash B(0,\frac{\delta }{2}\mu _{t}^{-1})}\widetilde{u}%
_{t}^{2}dv_{\widetilde{\mathbf{g}}_{t}} \\
&&+\int_{B(0,\delta \mu _{t}^{-1})\backslash B(0,R)}\eta ^{2}(\mu
_{t}x)\lambda _{t}\widetilde{f}_{t}\widetilde{u}_{t}^{2^{*}}\left| x\right|
^{2}dv_{\widetilde{\mathbf{g}}_{t}}-\int_{B(0,\delta \mu
_{t}^{-1})\backslash B(0,R)}\eta ^{2}(\mu _{t}x)\mu _{t}^{2}\widetilde{h}_{t}%
\widetilde{u}_{t}^{2}\left| x\right| ^{2}dv_{\widetilde{\mathbf{g}}_{t}} \\
&&-(n-1)\int_{B(0,\delta \mu _{t}^{-1})\backslash B(0,R)}\eta ^{2}(\mu _{t}x)%
\widetilde{u}_{t}^{2}dv_{\widetilde{\mathbf{g}}_{t}}+C\int_{B(0,\delta \mu
_{t}^{-1})\backslash B(0,R)}\eta ^{2}(\mu _{t}x)\mu _{t}^{2}\widetilde{u}%
_{t}^{2}\left| x\right| ^{2}dv_{\widetilde{\mathbf{g}}_{t}} \\
&&-\int_{B(0,\delta \mu _{t}^{-1})\backslash B(0,R)}\eta ^{2}(\mu
_{t}x)\left[ (\nabla \widetilde{u}_{t},x)_{\widetilde{\mathbf{g}}_{t}}+%
\widetilde{u}_{t}\right] ^{2}dv_{\widetilde{\mathbf{g}}_{t}}\\
&&-2\int_{B(0,\delta \mu _{t}^{-1})\backslash B(0,R)}\eta ^{2}(\mu _{t}x)\widetilde{u}%
_{t}(\nabla \widetilde{u}_{t},x)_{\widetilde{\mathbf{g}}_{t}}^{2}dv_{%
\widetilde{\mathbf{g}}_{t}}
\end{eqnarray*}
Or puisque $n\geq 4$%
\[
-(n-1)\widetilde{u}_{t}^{2}-2\widetilde{u}_{t}(\nabla \widetilde{u}_{t},x)_{%
\widetilde{\mathbf{g}}_{t}}^{2}-\left[ (\nabla \widetilde{u}_{t},x)_{%
\widetilde{\mathbf{g}}_{t}}+\widetilde{u}_{t}\right] ^{2}=-(n-4)\widetilde{u}%
_{t}^{2}-\left[ (\nabla \widetilde{u}_{t},x)_{\widetilde{\mathbf{g}}_{t}}+2%
\widetilde{u}_{t}\right] ^{2}\leq 0\,\,. 
\]
De plus avec l'estim\'{e}e 
\[
\left| x^{2}\widetilde{u}_{t}^{2^{*}}\right| \leq \varepsilon _{R}\widetilde{%
u}_{t}^{2} 
\]
on obtient 
\[
\int_{B(0,\delta \mu _{t}^{-1})\backslash B(0,R)}\eta ^{2}(\mu _{t}x)\lambda
_{t}\widetilde{f}_{t}\widetilde{u}_{t}^{2^{*}}\left| x\right| ^{2}dv_{%
\widetilde{\mathbf{g}}_{t}}\leq \varepsilon _{R}\int_{B(0,\delta \mu
_{t}^{-1})}\widetilde{u}_{t}^{2}dv_{\widetilde{g}_{t}} 
\]
et enfin 
\[
\int_{B(0,\delta \mu _{t}^{-1})\backslash B(0,R)}\eta ^{2}(\mu _{t}x)\mu
_{t}^{2}(C-\widetilde{h}_{t})\widetilde{u}_{t}^{2}\left| x\right| ^{2}dv_{%
\widetilde{\mathbf{g}}_{t}}\leq C\delta ^{2}\int_{B(0,\delta \mu _{t}^{-1})}%
\widetilde{u}_{t}^{2}dv_{\widetilde{\mathbf{g}}_{t}}=\varepsilon _{\delta
}\int_{B(0,\delta \mu _{t}^{-1})}\widetilde{u}_{t}^{2}dv_{\widetilde{\mathbf{%
g}}_{t}} 
\]
o\`{u} $\varepsilon _{R}\stackunder{R\rightarrow +\infty }{\rightarrow }0$
et $\varepsilon _{\delta }\stackunder{\delta \rightarrow 0}{\rightarrow }0$.

D'o\`{u} finalement 
\begin{eqnarray*}
&&\int_{B(0,\delta \mu _{t}^{-1})\backslash B(0,R)}\eta ^{2}(\mu
_{t}x)[\left| \nabla (\left| x\right| \widetilde{u}_{t})\right| _{\widetilde{%
\mathbf{g}}_{t}}^{2}-(\nabla (\left| x\right| \widetilde{u}_{t}),\nu )_{%
\widetilde{\mathbf{g}}_{t}}^{2}]dv_{\widetilde{\mathbf{g}}_{t}} \\
&\leq &C\int_{\partial B(0,R)}\left| \nabla (\left| x\right| ^{2}\widetilde{u%
}_{t}^{2})\right| _{e}d\sigma _{e}+C\int_{B(0,\delta \mu
_{t}^{-1})\backslash B(0,\frac{\delta }{2}\mu _{t}^{-1})}\widetilde{u}%
_{t}^{2}dv_{\widetilde{\mathbf{g}}_{t}} \\
&&+\varepsilon _{R}\int_{B(0,\delta \mu _{t}^{-1})}\widetilde{u}_{t}^{2}dv_{%
\widetilde{\mathbf{g}}_{t}}+\varepsilon _{\delta }\int_{B(0,\delta \mu
_{t}^{-1})}\widetilde{u}_{t}^{2}dv_{\widetilde{\mathbf{g}}_{t}}\,\,.
\end{eqnarray*}
Mais maintenant, \`{a} $R$ fix\'{e} 
\[
\widetilde{u}_{t}\stackrel{C^{1}}{\stackunder{t\rightarrow t_{0}}{%
\rightarrow }}\widetilde{u}=(1+\frac{\lambda f(x_{0})}{n(n-2)}\left|
x\right| ^{2})^{-\frac{n-2}{2}}\,\text{sur }\partial B(0,R) 
\]
donc 
\[
\nabla (\left| x\right| ^{2}\widetilde{u}_{t}^{2})\nabla (\left| x\right|
^{2}\widetilde{u}^{2})=\left| x\right| ^{2}\nabla (\widetilde{u}^{2})+%
\widetilde{u}^{2}\nabla (\left| x\right| ^{2})\,\,. 
\]
Or 
\[
\int_{\partial B(0,R)}\left| \nabla (\left| x\right| ^{2}\widetilde{u}%
^{2})\right| _{e}d\sigma _{e}\sim C.R^{-(n-4)} 
\]
et de plus 
\[
\stackunder{t\rightarrow 1}{\overline{\lim }}\frac{1}{\int_{B(0,\delta \mu
_{t}^{-1})}\widetilde{u}_{t}^{2}dv_{\widetilde{\mathbf{g}}_{t}}}\leq 
\stackunder{t\rightarrow 1}{\overline{\lim }}\frac{1}{\int_{B(0,R^{\prime })}%
\widetilde{u}_{t}^{2}dv_{\widetilde{\mathbf{g}}_{t}}}\sim \frac{1}{%
\int_{B(0,R^{\prime })}\widetilde{u}^{2}dx} 
\]
pour tout $R^{\prime }$ fix\'{e}, avec si $n=4$ $\int_{B(0,R^{\prime })}%
\widetilde{u}^{2}dx\rightarrow +\infty $ $\,$lorsque $R^{\prime }\rightarrow
+\infty \,.$

En utilisant la concentration $L^{2}$ (point d/ de l'\'{e}tude du
ph\'{e}nom\`{e}ne de concentration du chapitre 3), on obtient donc en fin de
compte 
\[
\stackunder{t\rightarrow 1}{\overline{\lim }}\frac{C_{t}}{\int_{B(0,\delta )}%
\overline{u}_{t}^{2}dx}\leq \varepsilon _{R}+\varepsilon _{\delta
}+C.R^{-(n-4)}\stackunder{t\rightarrow 1}{\overline{\lim }}\frac{1}{%
\int_{B(0,R^{\prime })}\widetilde{u}_{t}^{2}dv_{\widetilde{\mathbf{g}}_{t}}} 
\]
d'o\`{u} en faisant tendre $R$ vers l'infini 
\[
\stackunder{t\rightarrow 1}{\overline{\lim }}\frac{C_{t}}{\int_{B(0,\delta )}%
\overline{u}_{t}^{2}dx}\leq \varepsilon _{\delta } 
\]
o\`{u} $\varepsilon _{\delta }\rightarrow 0$ quand $\delta \rightarrow 0\,$.

Comme on le voit, la pr\'{e}sence d'une fonction $f$ au second membre de
l'\'{e}quation $E_{h,f,\mathbf{g}}$ ne pose pas de difficult\'{e} dans
l'obtention de cette limite, car on peut toujours majorer $f$ par son Sup
dans les int\'{e}grales. De plus, aucune d\'{e}riv\'{e}e de $f$
n'apparaissant, ces calculs restent valables dans le cas d'une famille
d'\'{e}quations du type 
\[
\Delta _{\mathbf{g}}u_{t}+h_{t}u_{t}=f_{t}u_{t} 
\]
o\`{u} l'on a une famille de fonctions ($f_{t}$) au second membre, tant que
les $f_{t}$ sont uniform\'{e}ment born\'{e}es. Nous nous servons de cela au
chapitre 6.

\chapter{Appendice D: notations et conventions}
\pagestyle{myheadings}\markboth{\textbf{Notations et conventions.}}
{\textbf{Notations et conventions.}}
Nous nous sommes efforc\'{e}s de garder les notations et les conventions
suivantes dans notre travail.

\textbf{Donn\'{e}es : }On consid\`{e}re une vari\'{e}t\'{e} riemannienne
compacte $(M,\mathbf{g})$ de dimension $n\geq 3.$ Soit $f:M\rightarrow \Bbb{R%
}$ une fonction $C^{\infty }$ \textit{fix\'{e}e} telle que $\stackunder{M}{%
Max}$ $f>0$. Soit aussi $h\in C{{}^{\infty }}(M)$ avec l'hypoth\`{e}se
suppl\'{e}mentaire que l'op\'{e}rateur $\bigtriangleup _{\mathbf{g}}+h$ est
coercif si $f$ change de signe sur $M$.

On consid\`{e}re l'\'{e}quation 
\[
(E_{h}^{\prime })=(E_{h,f}^{\prime })=(E_{h,f,\mathbf{g}}^{\prime
}):\,\triangle _{\mathbf{g}}u+h.u=f.u^{\frac{n+2}{n-2}}. 
\]

On note souvent $2^{*}=\frac{2n}{n-2}$ , remarquons alors que $2^{*}-1=\frac{%
n+2}{n-2}$.

\textbf{Convention:} Nous utilisons toujours ces ``notations'': $\mathbf{g,g}%
^{\prime },\mathbf{g}_{t},etc$ pour les m\'{e}triques (en gras $\mathbf{g}$
pour les distinguer plus clairement des fonctions); $h,h^{\prime },h_{t},etc$
pour les fonctions du premier membre d\'{e}finissant l'op\'{e}rateur $%
\bigtriangleup _{\mathbf{g}}+h$; \thinspace et enfin $f,f^{\prime
},f_{t},etc $ pour celles du deuxi\`{e}me membre; de plus $u,u_{t},etc$
d\'{e}signent les fonctions inconnues ou les solutions de $\,\triangle _{%
\mathbf{g}}u+h.u=f.u^{\frac{n+2}{n-2}}$ .

On s'int\'{e}resse aux solutions minimisantes de $E_{h,f,\mathbf{g}}^{\prime
}$: on dit que $u\in H_{1}^{2}(M)$ est minimisante pour $(E_{h,f,\mathbf{g}%
}^{\prime })$ (ou par abus de langage, minimisante pour $h$) si pour 
\[
I_{h,\mathbf{g}}(w)=\int_{M}\left| \nabla w\right| _{\mathbf{g}}^{2}dv_{%
\mathbf{g}}+\int_{M}h.w{{}^{2}}dv_{\mathbf{g}} 
\]
on a 
\[
I_{h,\mathbf{g}}(u)=\stackunder{w\in \mathcal{H}_{f}}{\inf }I_{h,\mathbf{g}%
}(w):=\lambda _{h,f,\mathbf{g}} 
\]
o\`{u} 
\[
\mathcal{H}_{f}=\{w\in H_{1}^{2}(M)/\int_{M}f\left| w\right| ^{\frac{2n}{n-2}%
}dv_{\mathbf{g}}=1\}. 
\]
On utilise aussi la fonctionnelle 
\[
J_{h,f,\mathbf{g}}(w)=J_{h}(w)=\frac{\int_{M}\left| \nabla w\right| _{%
\mathbf{g}}^{2}dv_{\mathbf{g}}+\int_{M}h.w{{}^{2}}dv_{\mathbf{g}}}{\left(
\int_{M}f\left| w\right| ^{\frac{2n}{n-2}}dv_{\mathbf{g}}\right) ^{\frac{n-2%
}{n}}} 
\]
et la partie de $H_{1}^{2}(M)$ pour laquelle elle est d\'{e}finie 
\[
\mathcal{H}_{f}^{+}=\{w\in H_{1}^{2}(M)/\int_{M}f\left| w\right| ^{\frac{2n}{%
n-2}}dv_{\mathbf{g}}>0\}. 
\]
Alors 
\[
\stackunder{w\in \mathcal{H}_{f}}{\inf }I_{h,\mathbf{g}}(w)=\stackunder{w\in 
\mathcal{H}_{f}^{+}}{\inf }J_{h,f,\mathbf{g}}(w)=\lambda _{h,f,\mathbf{g}} 
\]
L'\'{e}quation d'Euler associ\'{e}e au probl\`{e}me de minimisation de $I_{h,%
\mathbf{g}}(w)\,$ sur $\mathcal{H}_{f}$ est 
\[
(E_{h})=(E_{h,f})=(E_{h,f,\mathbf{g}}):\,\triangle _{\mathbf{g}%
}u+h.u=\lambda _{h,f,\mathbf{g}}f.u^{\frac{n+2}{n-2}} 
\]
o\`{u} $\lambda _{h,f,\mathbf{g}}\,$apparait comme une constante de
normalisation li\'{e}e \`{a} la condition 
\[
\int_{M}f\left| w\right| ^{\frac{2n}{n-2}}dv_{\mathbf{g}}=1 
\]
L'\'{e}quation d'Euler associ\'{e}e \`{a} $J_{h,f,\mathbf{g}}(w)$ et $%
\mathcal{H}_{f}^{+}$ est identique, mais sans cette constante.

Pour \'{e}tudier ces \'{e}quations, nous utilisons des \'{e}quations
associ\'{e}es 
\[
(E_{t}):\,\triangle _{\mathbf{g}}u_{t}+h_{t}.u_{t}=\lambda _{t}f_{t}.u_{t}^{%
\frac{n+2}{n-2}} 
\]
o\`{u} en g\'{e}n\'{e}ral $\lambda _{t}=\lambda _{h_{t},f_{t},\mathbf{g}}$
et dans certains cas $\lambda _{t}\leq \lambda _{h_{t},f_{t},\mathbf{g}}$.

Les th\'{e}or\`{e}mes des chapitres 3 \`{a} 6 utilisent l'hypoth\`{e}se
suivante :

\textbf{Hypoth\`{e}ses (H)}\textit{: On suppose que le Hessien de la
fonction }$f:M\rightarrow \Bbb{R}$, \textit{telle que }$\stackunder{M}{Sup}%
f>0$, \textit{est non d\'{e}g\'{e}n\'{e}r\'{e} en chaque point de maximum de 
}$f.$ \textit{En outre, les fonctions }$h$ \textit{consid\'{e}r\'{e}es sont
telles que }$\bigtriangleup _{\mathbf{g}}+h$ \textit{est coercif et l'on
suppose }$\dim M\geq 4$.\textit{\ On parle des hypoth\`{e}ses (}\textbf{H}$%
_{f}$\textit{) pour d\'{e}signer celles concernant la fonction }$f$\textit{.}

En ce qui concerne les notations les plus g\'{e}n\'{e}rales:

Nous notons 
\[
L^{p}(M)=L^{p} 
\]
l'ensemble des fonctions de puissance $p$-i\`{e}me int\'{e}grable; l'espace $%
M$ consid\'{e}r\'{e} \'{e}tant sous-entendu s'il n'y a pas d'ambiguit\'{e}. 
\[
\left\| .\right\| _{L^{p}(M)}=\left\| .\right\| _{L^{p}}=\left\| .\right\|
_{p} 
\]
est la norme associ\'{e}e. 
\[
H_{k}^{p}(M)=H_{k}^{p} 
\]
d\'{e}signe l'espace de Sobolev des fonctions dont les d\'{e}riv\'{e}es
jusqu'\`{a} l'odre $k$ sont dans $L^{p}$ et 
\[
\left\| .\right\| _{H_{k}^{p}(M)}=\left\| .\right\| _{H_{k}^{p}} 
\]
les normes correspondantes. En ce qui concerne les convergences, la
fl\`{e}che $\rightharpoondown $ d\'{e}signe une convergence faible, et la
fl\`{e}che $\rightarrow $ d\'{e}signe une convergence forte.

Pour la vari\'{e}t\'{e} Riemannienne $(M,\mathbf{g})$ 
\[
dv_{\mathbf{g}} 
\]
d\'{e}signe l'\'{e}l\'{e}ment de volume riemannien associ\'{e} \`{a} $%
\mathbf{g.\,}$Le laplacien riemannien associ\'{e} \`{a} $\mathbf{g}$ est
not\'{e} 
\[
\Delta _{\mathbf{g}}=-\nabla ^{i}\nabla _{i} 
\]
(attention au signe moins) 
\[
\left| \nabla u\right| _{\mathbf{g}} 
\]
est la norme de $\nabla u$ pour la m\'{e}trique $\mathbf{g}$. Lorsqu'il n'y
a pas de doute sur la m\'{e}trique consid\'{e}r\'{e}e, $\mathbf{g}$ est
sous-entendue; ainsi 
\[
\int_{M}\left| \nabla u\right| =\int_{M}\left| \nabla u\right| _{\mathbf{g}%
}dv_{\mathbf{g}}\,\,. 
\]
Par ailleurs, nous utilisons la m\'{e}trique euclidienne $\xi $, les
notations correspondantes \'{e}tant 
\[
\Delta _{e}=-\sum_{i}\partial _{ii}^{2}\,,\,\,\left| \nabla u\right| _{%
\mathbf{e}},\,\,dx,\,\,d\sigma _{e} 
\]
o\`{u} $dx$ est la mesure de Lebesgue et $d\sigma _{e}$ la mesure induite
sur la sph\`{e}re de dimension $n-1$.

Enfin, les meilleures constantes $K(n,2)$ et $B_{0}(\mathbf{g})$ sont
parfois not\'{e}es $K$ et $B$ pour simplifier quelques expressions. Le
principe g\'{e}n\'{e}ral est d'ailleurs que, pour all\'{e}ger certaines
expressions assez longues apparaissant dans ce travail, certains indices
sont omis lorsqu'il n'y a pas d'ambigu\"{i}t\'{e}: 
\begin{eqnarray*}
\lambda _{h,f,\mathbf{g}} &=&\lambda _{h} \\
J_{h,f,\mathbf{g}} &=&J_{h}
\end{eqnarray*}
et ainsi de suite...

Dans le d\'{e}tail:

- $B(x,\delta )$ est la boule de centre $x$ et de rayon $\delta $ pour la
distance g\'{e}od\'{e}sique.

- $x_{0}$ est toujours l'(unique) point de concentration de la suite $u_{t}.$

- $x_{t}$ est un point de maximum de $u_{t}$, et, \`{a} extraction près, 
$x_{t}\rightarrow x_{0}$

-$\mu _{t}$ param\`{e}tre fondamental est d\'{e}fini par 
\[
\stackunder{M}{Max}\,u_{t}=u_{t}(x_{t}):=\mu _{t}^{-\frac{n-2}{2}} 
\]

- $\delta $ d\'{e}signe toujours le rayon d'une petite boule autour de $%
x_{t} $ ou de $x_{0}$.

- $\eta $ d\'{e}signe toujours une fonction cut-off.

- les familles associ\'{e}es \`{a} l'\'{e}quation $(E_{h,f,\mathbf{g}})$ par
les \'{e}quations $(E_{t})$ sont index\'{e}es par un param\`{e}tre $%
t\rightarrow t_{0}$ o\`{u} nous prenons souvent $t_{0}=1.$

- $c,C,c^{\prime },C^{\prime }$ d\'{e}signent toujours des constantes
ind\'{e}pendantes des param\`{e}tres variables tels que $t$ ou $\delta $ ou $%
R$.

-$\varepsilon _{\delta }\stackunder{\delta \rightarrow 0}{\rightarrow }0$, 
$\varepsilon _{R}\stackunder{R\rightarrow +\infty }{\rightarrow }0$, 
$\varepsilon _{\delta,t }\stackunder{\delta \rightarrow 0, t\rightarrow1}{\longrightarrow }0$
\\
Pour les équations numérotées (a.b), a est le numéro du chapitre et b est le numéro de l'équation 
dans le chapitre.

\chapter{Bibliographie}
\pagestyle{myheadings}\markboth{\textbf{Bibliographie.}}{\textbf{Bibliographie.}}

[1]: T. AUBIN: Some nonlinear problems in Riemannian geometry, Springer
monograph in mathematics, 1998

[2]: T. AUBIN: Equation diff\'{e}rentielles non lin\'{e}aires et
probl\`{e}me de Yamabe concernant la courbure scalaire, J. Math. Pures
Appl., 55, 1976.

[3]: T. AUBIN: Probl\`{e}mes isop\'{e}rim\'{e}triques et espaces de Sobolev,
J. of Diff. Geometry, 11, 1976.

[4]: T. AUBIN: Meilleures constantes dans le th\'{e}or\`{e}me d'inclusion de
Sobolev et un th\'{e}or\`{e}me de Fredholm non lin\'{e}aire pour la
transformation conforme de la courbure scalaire, Journal of Functional
Analysis, 32, 1979.

[4]: A. BAHRI: C. R. Acad. Sci. de Paris, 307, (1998), n°11.

[5]: H. BERESTYCKI, L. NIRENBERG, S. VARADHAN: The principal eigenvalue and
maximum principle for second order elliptic operators in general domains,
Comm. Pure Appl. Math., 47, 1994.

[6]: L.A. CAFFARELLI, B. GIDAS, J. SPRUCK: Asymptotic symmetry and local
behavior of semilinear elliptic equations with Sobolev growth, Comm. in Pure
and Applied Math., 42, 1989.

[7]: S. COLLION: Fonction critique et EDP elliptiques sur les vari\'{e}t\'{e}%
s riemanniennes compactes, pr\'{e}publication de l'institut Elie Cartan,
2003, n${{}^{\circ }}$ 28.

[8]: S. COLLION: Transformation d'Abel et formes diff\'{e}rentielles
alg\'{e}briques, C.R. Acad\'{e}mie des Sciences de Paris, t.323, 1996.

[9]: Z. DJADLI-O. DRUET: Extremal functions for optimal Sobolev inequalities
on compact manifolds, Calc. Var., 12, 59-84, 2001

[10]: O. DRUET: Optimal Sobolev inequalities and extremal functions. The
three dimensional case. Indiana Univ. Math. J., to appear.

[11]: O. DRUET, E; HEBEY: The AB program in geometric analysis: sharp
Sobolev inequalities and related problems, Memoirs of the AMS, 761.

[12]: O. DRUET, E; HEBEY, M. VAUGON, Pohozahev type obstructions and
solutions of bounded energy for quasilinear elliptic equations with critical
Sobolev growth. The conformally flat case. Universit\'{e} de Cergy-Pontoise,
n${{}^{\circ }}$13, avril 2000.

[13]: O. DRUET-F. ROBERT: Asymptotic profile and blow-up estimates on
compact Riemannian manifolds, preprint, disponible dans Memoirs of the AMS,
761: The AB program in geometric analysis: sharp Sobolev inequalities and
related problems, par E. HEBEY et O. DRUET

[14]: Z. FAGET: Optimal constants in critical Sobolev inequalities on
riemannian manifolds in the presence of symmetries, Annals of Global
Analysis and Geometry, 24: 161-200, 2003.

[15]: Z. FAGET: Second best constant and extremal functions in Sobolev
inequalities in the presence of symmetries, preprint, \`{a} para\^{i}tre.

[16]: GILBARG-TRUDINGER: Elliptic partial differential equations of second
order, Springer 1985.

[17]: E. HEBEY: Sobolev spaces on riemannian manifolds, Lecture notes in
Mathematics, 1635, Springer, 1996.

[18]: E. HEBEY: Introduction \`{a} l'analyse non lin\'{e}aire sur les
vari\'{e}t\'{e}s, Diderot, 1997

[19]: E. HEBEY- M. VAUGON: The best constant problem in the Sobolev
embedding theorem for complete riemannian manifolds, Duke math. Journal, vol
79, July 1995.

[20]: E. HEBEY- M. VAUGON: From best constants to Critical functions, Math.
Z., 237, 737-767,2001.

[21]: E. HUMBERT-M. VAUGON: The problem of prescribed critical functions,
Preprint de l'institut Elie Cartan, 2003

[22]: J. JOST: Partial differential equations, Springer 2002.

[23]: J. M. LEE: Riemannian Manifolds, Springer Verlag, 2002.

[24]: J. M. LEE: Introduction to smooth manifolds, Springer Verlag, 2002.

[25]: J. M. LEE - T. H. PARKER: The Yamabe Problem, Bulletin of the AMS, 17,
n${{}^{\circ }}$1, 1987.

[26]: P.L. LIONS: The concentration-compactness principle in the calculus of
variations. The locally compact case. Part I, Annales de l'institut Henri
Poincar\'{e}, 1, 1984.

[27]: S. POHOZAEV: Eigenfunctions of the equations $\Delta u+\lambda f(u)=0$%
, Soviet. Math. Dokl., vol 6, 1965.

[28]: M. STRUWE: Variational Methods. Ergebnisse der Mathematik und ihrer Grenzgebiete,
34, Springer-Verlag, 1996.

[29]: G. TALENTI: Best constant in Sobolev inequality, Annali di Matematica
pura ed Applicata, 110, 1976.

[30]: M. VAUGON: Equations diff\'{e}rentielles non lin\'{e}aires sur les
vari\'{e}t\'{e}s riemanniennes compactes, Bull. des Sciences
Math\'{e}matiques, 106, 1982.

[31]: M. VAUGON: Transformation conforme de la courbure scalaire sur la
sph\`{e}re, Ann. Inst. Henri Poincar\'{e}, Vol 3, 1986, p55-65

[32]: M. VAUGON: Transformation conforme de la courbure scalaire sur une
vari\'{e}t\'{e} Riemannienne compacte, Journal of Functional Analysis, Vol
71, March 1987.

\end{document}